\pdfoutput=1
\documentclass[nofrench,12pt]{thloria}
\usepackage[pdftex]{tlhypref}
\usepackage[latin1]{inputenc}
\usepackage[french]{babel}
\usepackage{latexsym,amssymb,amsmath,amsfonts,amsthm,exscale}
\usepackage{graphicx}
\usepackage{graphics}
\usepackage{algorithm}
\usepackage{algorithmic}
\usepackage{cite}
\usepackage{xcolor}
\usepackage{aeguill}
\usepackage{enumerate}
\def\pathPic{Illustrations_arxiv}
\usepackage{url}
\usepackage{wasysym}

\usepackage{minitoc} 

\def\PP{\rm \hbox{I\kern-.2em\hbox{P}}}
\def\RR{\rm \hbox{I\kern-.2em\hbox{R}}}
\def\NN{\rm \hbox{I\kern-.2em\hbox{N}}}
\def\ZZ{\rm {{\rm Z}\kern-.28em{\rm Z}}}
\def\CC{\rm \hbox{C\kern -.5em {\raise .32ex \hbox{$\scriptscriptstyle
|$}}\kern-.22em{\raise .6ex \hbox{$\scriptscriptstyle |$}}\kern .4em}}
\def\vp{\varphi}
\def\<{\langle}
\def\>{\rangle}
\def\e{\varepsilon}
\def\sm{\setminus}
\def\nl{\newline}
\def\o{\overline}

\def\bq{{\bf q}}
\def\cT{{\cal T}}
\def\cA{{\cal A}}
\def\cI{{\cal I}}
\def\cV{{\cal V}}
\def\cB{{\cal B}}
\def\cR{{\cal R}}
\def\cS{{\cal S}}
\def\cD{{\cal D}}
\def\cH{{\cal H}}
\def\cL{{\cal L}}
\def\cN{{\cal N}}
\def\cC{{\cal C}}
\def\cF{{\cal F}}
\def\cP{{\cal P}}
\def\cQ{{\cal Q}}

\def\cO{{\cal O}}
\def\Chi{\raise .3ex
\hbox{\large $\chi$}} \def\vp{\varphi}
\def\lsima{\hbox{\kern -.6em\raisebox{-1ex}{$~\stackrel{\textstyle<}{\sim}~$}}\kern -.4em}
\def\lsim{\hbox{\kern -.2em\raisebox{-1ex}{$~\stackrel{\textstyle<}{\sim}~$}}\kern -.2em}
\def\({\Bigl (}
\def\){\Bigr )}
\def\({\Bigl (}
\def\){\Bigr )}

\newcommand{\be}{\begin{equation}}
\newcommand{\ee}{\end{equation}}
\newcommand{\bea}{$$ \begin{array}{lll}}
\newcommand{\eea}{\end{array} $$}
\newcommand{\bi}{\begin{itemize}}
\newcommand{\ei}{\end{itemize}}
\newcommand{\iref}[1]{(\ref{#1})}
\newtheorem{theorem}{Theorem}[section]
\newtheorem{remark}[theorem]{Remark}
\newtheorem{lemma}[theorem]{Lemma}
\newtheorem{proposition}[theorem]{Proposition}
\newtheorem{corollary}[theorem]{Corollary}
\newtheorem{definition}[theorem]{Definition}

\newtheorem{prop}[theorem]{Proposition}

\newtheorem{conjecture}[theorem]{Conjecture}

\newtheorem{theoremIntro}{Theorem}
\newtheorem{FRAtheoremIntro}{Théorème}

\def\I{{\rm \hbox{I\kern-.2em\hbox{I}}}}
\def\P{{\rm \hbox{I\kern-.2em\hbox{P}}}}
\def\sP{{\rm \hbox{\scriptsize I\kern-.2em\hbox{\scriptsize P}}}}
\def\H{{\rm \hbox{I\kern-.2em\hbox{H}}}}
\def\R{{\rm \hbox{I\kern-.2em\hbox{R}}}}
\def\N{{\rm \hbox{I\kern-.2em\hbox{N}}}}
\def\Z{{\rm {{\rm Z}\kern-.28em{\rm Z}}}}
\def\C{{\rm \hbox{C\kern -.5em {\raise .32ex \hbox{$\scriptscriptstyle
|$}}\kern-.22em{\raise .6ex \hbox{$\scriptscriptstyle |$}}\kern .4em}}}
\def\sC{{\rm \hbox{\scriptsize \rm C\kern -.6em {\raise .4ex \hbox{$\scriptscriptstyle \scriptsize|$}}\kern-.22em{\raise .5ex \hbox{$\scriptscriptstyle \scriptsize |$}}\kern .4em}}}

\def\proof{{\noindent \bf Proof: }}

\def\ra{\rightarrow}
\def\Ra{\Rightarrow}
\def\ve{\varepsilon}
\def\sq{\hfill $\diamond$\\}
\def\cO{\mathcal O}
\def\R{\mathbb R}

\def\b{\mathbf}
\def\mE{\mathbb E}
\def\cE{\mathcal E}

\def\mT{\mathbb T}

\def\nsdf{\|\sqrt{\det d^2f}\|}

\newcommand\sep{\; ; \;}
\newcommand\ssep{; \,}

\newcommand\trans{\mathrm T}

\def\cM{\mathcal M}
\def\cG{\mathcal G}

\newcommand\Keq{K_{\rm eq}}
\newcommand\Kpol{\b K}

\def\TEq{{T_{\rm eq}}}
\def\n{{\rm \bf n}}

\def\bac{$$\left\{\begin{array}}
\def\eac{\end{array}\right.$$}
\def\beqa{\begin{eqnarray*}}
\def\eeqa{\end{eqnarray*}}

\def\gsim{\hbox{\kern -.2em\raisebox{-1ex}{$~\stackrel{\textstyle>}{\sim}~$}}\kern -.2em}

\DeclareMathOperator\Id{Id}

\DeclareMathOperator\diam{diam}
\DeclareMathOperator\argmax{argmax}
\DeclareMathOperator\argmin{argmin}

\DeclareMathOperator\Cvx{Cvx}
\DeclareMathOperator\bary{bary}

\DeclareMathOperator\Vect{Vect}
\DeclareMathOperator\diver{div}

\DeclareMathOperator\Tr{Tr}
\DeclareMathOperator\supp{supp}

\DeclareMathOperator\disc{disc}
\DeclareMathOperator\diag{diag}
\DeclareMathOperator\interp{{\rm I}}

\DeclareMathOperator\Span{Span}
\def\bT{\b T}
\def\bH{\b H}

\newcommand\seqT{(\cT_N)_{N\geq N_0}}
\newcommand\mhalf{{s_m}}
\newcommand\cTr{{\cT_s^\text{reg}}}
\newcommand\cTb{{\cT_s^\text{bd}}}
\newcommand\mi{{\mathbf i}}

\DeclareMathOperator{\cyc}{\mathbb S}
\newcommand\NF{N\hspace{-0.8mm}F}

\def\cPi{{\text{\large$\boldsymbol\pi$}}}
\def\Lpol{\b L}
\DeclareMathOperator\M{M}
\DeclareMathOperator\GL{GL}
\DeclareMathOperator\SL{SL}

\def\TRect{{T_0}}
\def\TAc{{T_{\rm ac}}}
\def\ti{\tilde}
\newcommand\tb{\textcolor{blue}}
\newcommand\tr{\textcolor{red}}
\def\KEq{K_*}
\def\m{{\rm \bf m}}

\newtheorem{question}{Question}
\newtheorem{programme}[question]{Programmation}

\def\bques{\begin{question}}
\def\eques{\end{question}}
\def\bprog{\begin{programme}}
\def\eprog{\end{programme}}

\usepackage{stmaryrd}

\def\Th{{\rm Th }}

\def\mI{{\mathbb I}}
\def\u{{\bf u}}
\def\v{{\bf v}}
\def\seqR{(\cR_N)_{N\geq 0}}

\def\bt{{\rm \bf t}}
\def\bn{{\rm \bf n}}
\def\ssdelta{{\scriptscriptstyle\delta}}
\def\II{\mathrm {I\kern-0.1exI}}
\DeclareMathOperator\proj{P}
\DeclareMathOperator\move{U}

\DeclareMathOperator{\TV}{TV}
\DeclareMathOperator{\PSNR}{PSNR}
\newcommand\stext[1]{\ \text{ #1 } \ }

\def\TO{{T_0}}

\def\TRef{T_{\rm ref}}
\def\cPRef{\cP_{\rm ref}}
\newtheorem{assumptions}[theorem]{Assumptions}

\newtheorem{property}[theorem]{Property}

\DeclareMathOperator\Avg{Avg}
\DeclareMathOperator\dil{dil}
\DeclareMathOperator\Lip{Lip}
\DeclareMathOperator\interior{int}

\DeclareMathOperator\Refine{Refine}
\DeclareMathOperator\Vor{Vor}

\DeclareMathOperator\wed{wedge}
\DeclareMathOperator\Inv{Inv}

\def\gT{\mathfrak T}

\def\bbS{{\mathbb S}}
\def\sR{{\rm \hbox{\scriptsize I\kern-.2em\hbox{\scriptsize R}}}}

\def\VEq{{V_{\rm eq}}}

\def\per{{\rm per}}
\def\loc{{\rm loc}}

\def\sob{{\rm sob}}
\def\proj{{\rm proj}}

\def\sN{{\rm \hbox{\scriptsize I\kern-.2em\hbox{\scriptsize N}}}}
\def\sC{{\rm \hbox{\scriptsize \rm C\kern -.6em {\raise .4ex \hbox{$\scriptscriptstyle \scriptsize|$}}\kern-.22em{\raise .5ex \hbox{$\scriptscriptstyle \scriptsize |$}}\kern .4em}}}

\def\gL{{\mathfrak L}}

\DeclareMathOperator\Ker{Ker}
\DeclareMathOperator\Image{Im}
\def\bbB{{\mathbb B}}

\ThesisTitle{Approximation adaptative et anisotrope par éléments finis\\Théorie et Algorithmes.}
\ThesisDate{6 décembre 2010}
\ThesisAuthor{Jean-Marie Mirebeau}
\ThesisNancyI

\President={}
\Rapporteurs={
Weiming {\sc Cao}\\
Gabriel {\sc Peyré}
}							
\Examinateurs={
Albert {\sc Cohen} & Directeur de th\`ese\\ 
Jean-Daniel {\sc Boissonnat}\\
Bruno {\sc Després}\\
Ronald {\sc DeVore}\\
Frédéric {\sc Hecht}\\
Yves {\sc Meyer}
} 

\SetRealMargins{30mm}{20mm}

\begin{document}
\dominitoc

\nthks
\MakeThesisTitlePage

\begin{ThesisAcknowledgments}
\bigskip

%


En tout premier lieu, je voudrais remercier Albert Cohen pour ces trois années passées sous sa direction.
En particulier, je lui suis reconnaissant de m'avoir dès le début proposé des sujets de recherche passionnants,
puis d'avoir progressivement stimulé et éduqué ma prise d'indépendance,
par son regard mathématique précis et critique, mais aussi bienveillant et constructif.
Je le remercie chaleureusement pour sa très grande disponibilité et son attention constante
qui m'ont poussé à donner le meilleur de moi-même.
Etre son étudiant fut un honneur et une chance exceptionnels.
\\



Je remercie vivement Weiming Cao et Gabriel Peyré d'avoir accepté d'être les rapporteurs de cette thèse. La précision de leurs commentaires et l'enthousiasme de certaines remarques sont pour moi un encouragement sans égal.
En dehors de ce rôle de rapporteur, je suis reconnaissant à Weiming Cao d'avoir très tôt apprécié et encouragé ma contribution à son domaine de recherche, notamment en m'invitant à un minisymposium; et
à Gabriel Peyré pour son dynamisme et son enthousiasme en mathématiques,
et pour de nombreuses et instructives discussions. \\




Je suis très honoré que Jean-Daniel Boissonnat, Bruno Després, Ronald DeVore, Frédéric Hecht et Yves Meyer témoignent de l'intérêt qu'ils portent à ma thèse en prenant part à son jury.
Je remercie par ailleurs Jean-Daniel Boissonnat et Frédéric Hecht pour nos discussions
qui m'ont beaucoup instruit sur la génération pratique de maillages;
je suis également reconnaissant à Frédéric Hecht pour m'avoir aidé à intégrer dans FreeFem++
certains outils développés dans cette thèse.
Un grand merci à Ronald DeVore pour m'avoir très tôt accepté comme collaborateur
et pour ses nombreux conseils sur le plan mathématique comme sur le plan humain.
\\

%
%
%
%





Je remercie vivement Nira Dyn, qui m'a invité à discuter avec elle, Albert Cohen et Shai Dekel une semaine
durant laquelle de nombreuses idées ont pu émerger.
J'ai beaucoup apprécié la collaboration avec Yuliya Babenko, que je remercie pour
son enthousiasme et son efficacité.
Je suis très reconnaissant à Suzanne Brenner et Li-Yen Sung,
pour leur invitation en Louisiane et l'extrême gentillesse de leur accueil.
\\

Je tiens à remercier l'ensemble du Laboratoire Jacques Louis Lions,
chercheurs, informaticiens, secrétaires et doctorants,
qui donnent à ce lieu de travail une ambiance très agréable.
Parmi les thésards, je voudrais remercier en particulier :
Evelyne pour les sorties escalade, partagées avec Frank, Frédéric
 et ma femme;
Mathieu qui m'a fait découvrir la symphonie à Gollum;
Alexis, Marianne et Pierre, ainsi que les nouveaux venus,
Marie, Malek et Ange,
avec qui j'ai eu le plaisir de partager le bureau 3D23.
\\

%
%












Je voudrais remercier mes amis, et parmi eux ceux de toujours, Emmanuel, Romain, Valentin.

Je voudrais remercier ma famille pour son soutien inconditionnel. Je remercie mes parents et mon frère de leur curiosité, parfois incrédule, de mes sujets d'étude.

Je dédie cette thèse à 
ma femme Jennifer, sans qui rien de ceci n'aurait de sens, et à mon fils Nathanaël pour le bonheur qu'il m'apporte chaque jour. 


\end{ThesisAcknowledgments}

\begin{ThesisDedication}
A Jennifer\\
et Nathanaël
\end{ThesisDedication}

\renewcommand{\contentsname}{Summary}
\renewcommand{\partname}{Part}
\WritePartLabelInToc
\renewcommand{\chaptername}{Chapter}
\renewcommand{\bibname}{Bibliography}
\WriteChapterLabelInToc
\tableofcontents
\mainmatter

\chapter*{Introduction}
\DontWriteThisInToc
\begin{quotation}
{\it
Although this may seem a paradox, all exact science is dominated by the idea of approximation. \hfill (Bertrand Russel, logicien et prix Nobel)} 
\end{quotation}

\section[Approximation adaptative et anisotrope : pourquoi et comment?]{Approximation adaptative et anisotrope par éléments finis : pourquoi et comment?} 

Cette thèse est consacrée au problème de l'approximation de fonctions par des éléments finis polynomiaux par morceaux sur des triangulations, et sur des maillages plus généraux. Nous nous intéressons tout particulièrement à la situation où le maillage est {\it construit en adaptation} avec la fonction approchée. Ce maillage peut donc comporter des éléments de taille, de rapport d'aspect et d'orientation fortement variables. \\

L'approximation par des fonctions polynomiales par morceaux est une procédure qui intervient dans de nombreuses applications. Dans certaines d'entre elles comme la compression de données de terrain, d'images ou de surfaces, la fonction $f$ approchée peut être connue exactement. Dans d'autres applications comme le débruitage de données, l'apprentissage statistique ou la discrétisation par éléments finis d'Equations aux Dérivées Partielles (EDP), la fonction approchée n'est connue que partiellement, voire totalement inconnue initialement.
Dans toutes ces applications, on établit usuellement une distinction entre l'approximation {\it uniforme} ou {\it adaptative}. Dans le cadre uniforme, le domaine de définition de la fonction est décomposé en une partition formée d'éléments de taille et de forme comparables, alors que ces attributs peuvent varier fortement dans le cadre adaptatif.
La partition peut dans ce dernier cas être adaptée aux propriétés locales de la fonction $f$, dans l'objectif d'optimiser le compromis entre la précision et la complexité de l'approximation.

\begin{figure}
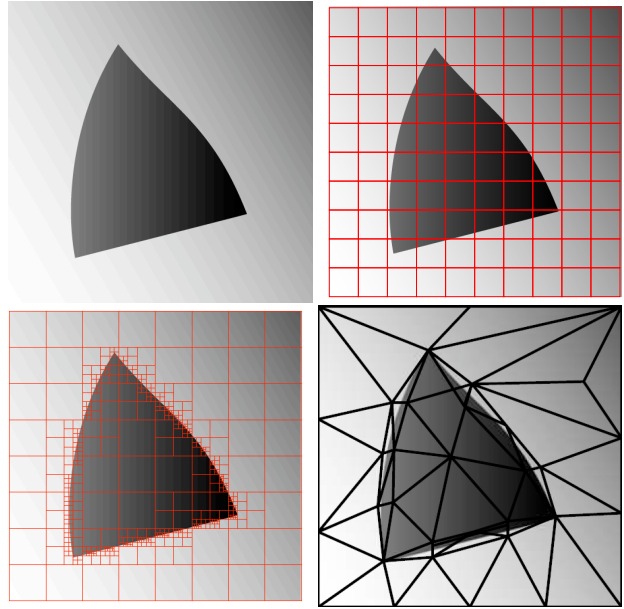

\centering
\includegraphics[height = 4cm,width = 4cm]{\pathPic/Cartoon/Triangle.pdf}
\includegraphics[height = 4cm,width = 4cm]{\pathPic/Cartoon/TriangleUnif.pdf}\\
\includegraphics[height = 4cm,width = 4cm]{\pathPic/Cartoon/Haar.pdf}
\includegraphics[height = 4cm,width = 4cm]{\pathPic/Cartoon/TriangleMeshed.pdf}
\caption{Une fonction régulière par morceaux (haut, gauche), Partition uniforme du domaine de définition (haut, droite), Partition associée à l'approximation adaptative fondée sur les plus grands coefficients du développement de la fonction sur la base de Haar (bas, gauche), Partition associée à l'approximation adaptative anisotrope de la fonction par des éléments finis (bas, droite). (source : G. Peyr\'e)\label{FRAfigEllTri0}}
\end{figure}

\begin{figure}
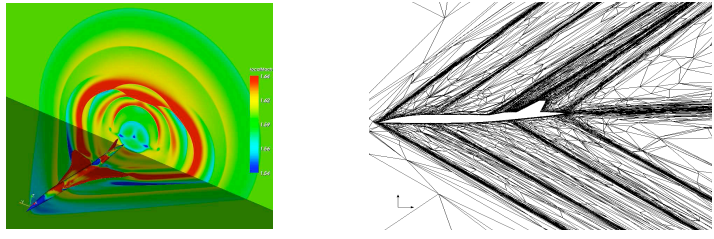

	\centering		
	\includegraphics[width=3.5cm,height=3cm]{\pathPic/SimuExt/AlauzetShock.pdf}
	\hspace{1cm}
	\includegraphics[width=4.5cm,height=3cm]{\pathPic/SimuExt/AlauzetPlane.pdf}
	\caption{Flux d'air autour d'un avion supersonique, calculé à l'aide d'un maillage tridimensionnel fortement anisotrope. (source : F. Alauzet \cite{A}.\label{FRAfigInria})}
\end{figure}

Du point de vue de la théorie de l'approximation, ce compromis entre précision et complexité est généralement lié à la régularité de la fonction: on s'attend typiquement à des vitesses de convergence plus élevées pour des fonctions plus régulières. Les fonctions qui se présentent dans les applications concrètes peuvent cependant présenter des propriétés hétérogènes de régularité, dans le sens où elles sont régulières dans certaines régions, qui séparent des discontinuités localisées.
Deux exemples typiques, illustrés Figure \ref{FRAfigEllTri0} et Figure \ref{FRAfigInria}, sont (i) les bords des objets dans les fonctions qui représentent des images, et (ii) les chocs dans les solutions d'EDP hyperboliques et non-linéaires.
Les méthodes numériques destinées au traitement de l'image, pour le débruitage ou la compression notamment, ou dédiées à la simulation des EDP bénéficient grandement de représentations économiques et précises de telles fonctions.
Une approximation polynomiale par morceaux sur une partition uniforme du domaine n'est généralement pas suffisante à cet effet. Un premier pas vers l'adaptativité consiste à faire varier la taille des éléments qui forment la partition en fonction des propriétés locales de régularité de la fonction. Un procédé très similaire, souvent utilisé en traitement de l'image et illustré Figure \ref{FRAfigEllTri0} (bas, gauche), consiste à retenir les termes d'une décomposition en ondelettes de $f$ qui correspondent aux plus grands coefficients. Un second pas dans l'adaptativité est de remarquer qu'une plus grande résolution de la partition est requise dans la direction orthogonale à la courbe de discontinuité que dans la direction tangentielle, et de tirer parti de cette propriété en employant une partition {\it anisotrope} du domaine. En deux dimensions ces partitions sont typiquement construites à l'aide de triangles de fort rapport d'aspect alignés avec les discontinuités, comme illustré Figure \ref{FRAfigEllTri0} (bas droite) et Figure \ref{FRAfigInria}.

Dans le contexte de la simulation numérique des EDP, l'adaptativité signifie également que le maillage de calcul n'est pas fixé a priori, mais dynamiquement mis à jour au cours de la simulation à mesure que la solution exacte se dévoile. D'un point de vue numérique ces méthodes requièrent des algorithmes et des structures de données plus complexes que leurs pendants non-adaptatifs. D'un point de vue théorique l'analyse de ces algorithmes est difficile, lorsqu'elle est possible. On ne sait en fait rigoureusement établir que les méthodes adaptatives améliorent la vitesse de convergence des solutions approchées vers la solution exacte, que pour un nombre réduit de systèmes d'EDP et seulement dans le cas de l'adaptation {\it isotrope} de maillage. Nous référons au survey \cite{Nochetto} pour une vue d'ensemble de ces résultats dans le cas des équations elliptiques. Nous devons mentionner que ces difficultés sont accentuées lorsque des éléments finis {\it anisotropes} sont utilisés. Plusieurs logiciels comme \cite{FreeFem, Inria, A} utilisent cependant avec succès les maillages anisotropes pour la simulation numérique des EDP, comme illustré par exemple Figure \ref{FRAfigInria}. D'un point de vue numérique l'amélioration apportée par ces méthodes semble évidente en comparaison avec les méthodes non-adaptatives ou adaptatives isotropes. Cependant de nombreux aspects de l'analyse théorique de ces méthodes restent des questions ouvertes.\\

Cette thèse étudie le problème de l'adaptation de maillage anisotrope pour l'approximation d'une fonction {\it connue}. Ceci peut être considéré comme une étape préliminaire pour l'analyse de l'adaptation anisotrope de maillage dans la simulation des EDP, mais d'autres applications peuvent en tirer parti comme le traitement de données de terrain, de surfaces ou d'images.\\

Etant donnée une triangulation $\cT$ d'un domaine borné et polygonal $\Omega\subset \R^2$, et un entier fixé $k \geq 1$, nous notons $V_k(\cT)$ l'espace des éléments finis de degré $k$ sur $\cT$. L'espace $V_k(\cT)$ est formé de toutes les fonctions qui coïncident sur chaque triangle $T\in \cT$ avec un polynôme de degré total $k$
$$
V_k( \cT) := \{g\; ; \; g_{|T} \in \P_k,\; T \in \cT\}.
$$
La dimension de $V_k(\cT)$ est de l'ordre de $\cO(k^2 \#(\cT))$. Etant donnée une fonction $f : \Omega \to \R$ et une triangulation $\cT$ de $\Omega$, l'erreur de meilleure approximation de $f$ dans $V_k(\cT)$ est définie par 
\be
\label{FRAeqApproxG}
e_\cT(f)_X := \inf_{g\in V_k(\cT)} \|f-g\|_X.
\ee
La lettre $X$ désigne la norme ou la semi-norme dans laquelle l'erreur d'approximation $\|f-g\|_X$ est mesurée. Dans cette thèse nous restreignons notre attention à la norme $L^p$ et à la semi-norme $W^{1,p}$, où $1\leq p \leq \infty$. Elles sont définies comme suit:
$$
\|h\|_{L^p(\Omega)} := \left(\int_\Omega |h|^p\right)^{\frac 1 p} \stext{ et } |h|_{W^{1,p}(\Omega)} := \left(\int_\Omega |\nabla h|^p\right)^{\frac 1 p},
$$
avec la modification usuelle lorsque $p=\infty$. Notons que l'on doit imposer la continuité globale de $g$ dans la définition ci dessus de $V_k(\cT)$ lorsqu'on utilise la semi-norme $W^{1,p}$.

La meilleure approximation $g\in V_k(\cT)$ de $f$ peut être calculée exactement dans le cas de la norme $L^2$ ou de la semi-norme $W^{1,2}$ (ou $H^1$): $g$ est la projection orthogonale de $f$ sur $V_k(\cT)$ par rapport au produit scalaire associé à la norme d'intérêt. Dans le cas $p\neq 2$ de normes non hilbertiennes la meilleure approximation de $f$ est généralement difficile à calculer, mais des approximations ``satisfaisantes'' peuvent être obtenues par différentes méthodes. Si la fonction est continue, on peut utiliser l'interpolation de Lagrange, tandis qu'un opérateur de quasi-interpolation est préféré pour les fonctions non-lisses, voir Chapitre \ref{chapApproxMet}.
Plus généralement, si $P_\cT$ est un opérateur arbitraire de projection de l'espace $X$ sur $V_k(\cT)$, il est aisé de voir que l'on a pour toute $f\in X$
\be
\label{FRAstableproj}
 \|f-P_\cT f\|_X \leq Ce_\cT(f)_X,
 \ee
où $C:=1+\|P_\cT\|_{X\to X}$. 
Le problème de l'approximation d'une fonction $f$ sur une triangulation {\it donnée} $\cT$, par des éléments finis de degré $k$, est donc en bonne partie résolu.\\

Dans le contexte de l'approximation adaptative, la triangulation $\cT$ du domaine $\Omega$ n'est pas fixée, mais peut être choisie librement en adaptation avec $f$ (par contraste nous supposons toujours dans cette thèse que l'entier $k$ est fixé, bien qu'arbitraire).
Ceci nous mène naturellement à l'objectif de caractériser et de construire un maillage optimal pour une fonction donnée $f$.  Etant donnée une norme $X$ d'intérêt et une fonction $f$ à approcher, nous formulons le problème de {\it l'adaptation optimale de maillage}, comme la minimisation de l'erreur d'approximation parmi toutes les triangulations de {\it cardinalité donnée}. Nous définissons donc l'erreur de meilleure approximation adaptative comme suit:
\be
\label{FRAeqApproxT}
e_N(f)_X := \inf_{\#(\cT) \leq N} e_\cT(f)_X=\inf_{\#(\cT) \leq N}\inf_{g\in V_k(\cT)} \|f-g\|_X.
\ee
Par contraste avec le problème de la meilleure approximation par éléments finis sur un maillage fixé, l'approximation adaptative et anisotrope n'est pas encore bien comprise. En particulier (i) {\it comment le maillage optimal dépend-il de la fonction} $f$, et (ii) {\it comment l'erreur optimale $e_N(f)_X$ décroît-elle lorsque $N$ augmente} ? Ces problèmes sont bien compris dans le cadre isotrope, où l'optimisation est restreinte aux triangulations dans lesquelles les triangles satisfont uniformément une contrainte sur leur rapport d'aspect
$$
\diam(T)^2\leq C |T|
$$
où $\diam(T)$ et $|T|$ représentent le diamètre et l'aire de $T$ respectivement, et $C>0$ est une constante fixée. Dans le cadre général des triangulations potentiellement anisotropes, ce sont des problèmes ouverts.

Heuristiquement, la simplicité de \iref{eqApproxG} en comparaison avec \iref{eqApproxT} vient ce que l'optimisation est posée sur l'espace linéaire $V_k(\cT)$, et qu'une solution presque optimale peut donc être obtenue en appliquant à $f$ un opérateur de projection stable comme indiqué dans \iref{stableproj}.
Par contraste, le problème d'optimisation \iref{eqApproxT} est posé sur la {\it réunion} d'espaces $V_k(\cT)$ pour toutes les triangulations $\cT$ satisfaisant $\#(\cT) \leq N$. Il s'agit donc d'un problème {\it d'approximation non-linéaire}. D'autres exemples de ce type de problème sont l'approximation par les $N$ meilleurs termes dans un dictionnaire de fonctions, ou l'approximation par des fonctions rationnelles. Nous référons le lecteur à \cite{De} pour un survey sur l'approximation non-linéaire.\\

L'objectif de cette thèse est de mieux comprendre le problème d'adaptation de maillage optimale posé sur {\it la classe entière des triangulations potentiellement anisotropes}. Les quatre parties de cette thèse sont consacrées respectivement aux quatre questions ci-dessous:
\begin{enumerate}[I.]
\item 
Comment l'erreur d'approximation $e_N(f)_X$ se comporte-t-elle dans le régime asymptotique où le nombre $N$ de triangles tend vers l'infini, lorsque $f$ est une fonction suffisamment régulière? Nous établissons dans ce contexte une caractérisation mathématique du maillage optimal, ainsi que des estimations précises supérieures et inférieures de $e_N(f)_X$ à l'aide de $N$ et de quantités qui dépendent {\it non linéairement} des dérivées de $f$.
\item 
Quelles classes de fonctions gouvernent la vitesse de décroissance de $e_N(f)_X$ lorsque $N$ augmente, et sont en ce sens naturellement liées au problème d'adaptation optimale de maillage?
Nous pensons en particulier à ce qu'on appelle les {\it fonctions cartoon}, qui par définition sont régulières excepté le long d'une famille de courbes de discontinuité. Il s'agit d'un modèle d'image populaire, illustré par exemple Figure \ref{FRAfigEllTri0}. Nous verrons que ce modèle s'inscrit naturellement dans une classe de fonctions plus riche qui correspond à une vitesse donnée de décroissance de $e_N(f)_X$.
\item 
Le problème d'optimisation \iref{eqApproxT}, qui porte sur les triangulations $\cT$ de cardinalité donnée $N$, peut-il être remplacé par un problème équivalent mais plus accessible? Les triangulations sont en effet des objets combinatoires et discrets, décrits par leurs sommets et arêtes, ce qui est peu commode lorsque l'on résout des problèmes d'optimisation de la forme de \iref{eqApproxT}. Nous étudions la correspondance entre certaines classes de triangulations et de {\it métriques riemanniennes} qui sont par contraste des objets continus. Ceci nous permet de reformuler le problème original d'optimisation par un problème plus facilement soluble posé sur l'ensemble des métriques riemanniennes.
\item 
Est-il possible de construire une suite quasi-optimale de triangulations $(\cT_N)_{N \geq 0}$, où $\#(\cT_N) = N$, en utilisant une {\it procédure hiérarchique de raffinement}? 
La propriété de hiérarchie garantit l'inclusion des espaces d'éléments finis associés aux triangulations : $V_k(\cT_N) \subset V_k(\cT_{N+1})$. Elle est requise dans des applications comme le codage progressif ou le traitement de données en ligne (i.e. à mesure qu'elles sont transmises). Nous proposons un algorithme simple et explicite qui donne une réponse positive à cette question sous certaines conditions.
\end{enumerate}

\DontWriteThisInToc
\section{Etat de l'art}

Avant d'entrer dans le détail du contenu de la thèse, nous donnons dans cette section un aperçu rapide de l'état de l'art dans les deux problèmes de l'étude de $e_N(f)_X$ lorsque $N$ augmente et de la construction d'une triangulation proche de l'optimale.
Nous devons pour cela introduire certaines notations. Nous supposons dans la suite que $f\in C^0(\overline \Omega)$, où $\Omega \subset \R^2$ est un domaine polygonal borné et $\overline \Omega$ désigne l'adhérence de $\Omega$. Pour chaque triangulation $\cT$ de $\Omega$ nous notons $\interp^k_\cT$ l'opérateur usuel d'interpolation de Lagrange sur les éléments finis de degré $k$ sur $\cT$: sur chaque triangle $T\in \cT$ l'interpolation $\interp_\cT^k f$ est l'unique élément de $\P_k$ qui coïncide avec $f$ aux points de coordonnées barycentriques $\{0,\frac 1 k, \cdots , \frac {k-1} k, 1\}$. Lorsque $k=1$ ces points sont simplement les trois sommets de $T$.
Si la triangulation $\cT$ est conforme (chaque arête de chaque triangle est soit sur le bord de $\Omega$, soit coïncide avec l'arête entière d'un autre triangle), alors $\interp_\cT^k f$ est continu.
Pour être consistant avec le reste de cette thèse nous définissons $m := k+1 \geq 2$, et nous avons donc 
$$
e_\cT(f)_{L^p} \leq \|f-\interp_\cT^{m-1} f\|_{L^p(\Omega)} \stext{ et } e_\cT(f)_{W^{1,p}} \leq  \|\nabla f-\nabla \interp_\cT^{m-1} f\|_{L^p(\Omega)}.
$$
Toute estimation sur l'erreur d'interpolation donne donc une estimation supérieure sur l'erreur de meilleure approximation $e_\cT(f)_X$. De plus si $f$ est suffisamment régulière et $\cT$ suffisamment fine, alors ces quantités sont généralement comparables.\\

L'un des résultats fondateurs en approximation adaptative et anisotrope par des éléments finis porte sur les éléments finis affines par morceaux ($m=2$) lorsque l'erreur est mesurée en norme $L^p$, voir aussi \cite{CSX} ou Chapitre \ref{chapOptAniso}. Ce résultat peut être formulé comme suit: pour tout domaine polygonal borné $\Omega\subset \R^2$, pour tout $1\leq p < \infty$ et pour toute fonction $f\in C^2(\overline \Omega)$, il existe une suite $(\cT_N)_{N \geq N_0}$ de triangulations de $\Omega$ satisfaisant $\#(\cT_N) \leq N$, et telles que 
\be
\label{FRAeqCSX}
\limsup_{N\to \infty} N \|f-\interp^1_{\cT_N} f\|_{L^p(\Omega)} \leq C \left\| \sqrt{|\det (d^2f)|}\right\|_{L^\tau(\Omega)},
\ee
où l'exposant $\tau \in (0, \infty)$ est défini par  
$$
\frac 1 \tau := 1 + \frac 1 p\, ,
$$ 
et $C$ est une constante universelle ($C$ est indépendante de $p$, $\Omega$ et $f$).
Nous rappelons que la limite supérieure d'une suite $(u_N)_{N \geq N_0}$ est définie par 
\be
\label{FRAdefLimSup}
\limsup_{N \to \infty} u_N := \lim_{N \to \infty} \sup_{n \geq N} u_n,
\ee
et est en général strictement inférieure au supremum $\sup_{N \geq N_0} u_N$. Trouver une majoration pertinente de $\sup_{N \geq N_0} \|f-\interp_{\cT_N}^1 f\|_{L^p(\Omega)}$ reste aujourd'hui un problème ouvert lorsque des triangulations anisotropes adaptées de manière optimale sont utilisées. Le résultat \iref{eqCSX} s'étend à l'exposant $p=\infty$, et aux maillages simpliciaux de domaines $\Omega$ de plus grande dimension, mais les maillages $\cT_N$ peuvent ne pas être conformes.

Le résultat \iref{eqCSX} révèle que la précision de l'approximation de $f$ est gouvernée par la quantité $\sqrt{|\det (d^2 f)|}$ qui dépend non linéairement de la matrice hessienne $d^2 f$. Cette dépendance non linéaire est fortement liée au fait que nous autorisons des triangles de formes potentiellement fortement anisotropes. L'estimation d'erreur \iref{eqCSX} est obtenue en produisant des triangulations suffisamment fines qui combinent les deux propriétés heuristiques suivantes:
\begin{enumerate}[a)]
\item {\it Equidistribution des erreurs:} la contribution $\|f-\interp_T^1 f\|_{L^p(T)}$ de chaque triangle $T\in \cT_N$ à l'erreur d'interpolation globale  $\|f-\interp_{\cT_N}^1 f\|_{L^p(T)}$ est du même ordre. Cette condition se traduit par une contrainte locale sur l'aire des triangles, qui est dictée par le comportement local de $f$ et en particulier par $\det (d^2 f(z))$.
\item {\it Forme optimale des triangles:} le rapport d'aspect et l'orientation d'un triangle $T\in \cT_N$ est dicté par le rapport des valeurs propres et par la direction des vecteurs propres de la matrice hessienne $d^2 f(z)$ pour $z\in T$.
\end{enumerate}

La méthode la plus simple pour construire une suite $(\cT_N)_{N \geq N_0}$ de triangulations satisfaisant \iref{eqCSX} est d'utiliser une stratégie de ``patchs locaux'' que l'on peut décrire intuitivement comme suit. Dans une première étape le domaine $\Omega$ est découpé en régions $\Omega_i$, $1\leq i \leq k$, suffisamment petites pour que la matrice hessienne $d^2 f(z)$ varie peu sur chaque $\Omega_i$ autour d'une valeur moyenne $M_i$. En d'autres termes $f$ est bien approchée par un polynôme quadratique sur chaque $\Omega_i$. Chaque région $\Omega_i$ est ensuite pavée par une triangulation uniforme $\cT_N^i$ dont les mailles sont de taille, de rapport d'aspect et d'orientation dictés par $M_i$. Les triangulations $\cT_N^i$ de $\Omega$ sont ensuite recollées de manière conforme pour former une triangulation $\cT_N$ de $\Omega$, au prix de quelques triangles supplémentaires aux interfaces entre les $\Omega_i$.

\begin{figure}
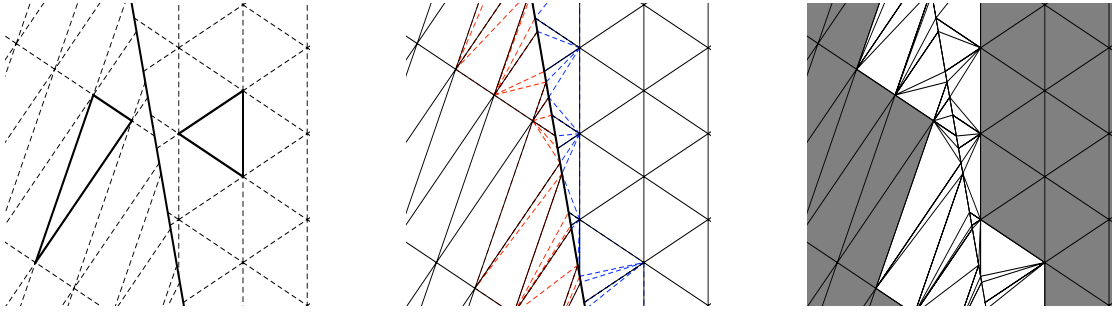

	\centering
		\includegraphics[width=4cm,height=4cm]{\pathPic/PaperOptAniso/NonConformTiling.pdf}
		\hspace{1cm}
		\includegraphics[width=4cm,height=4cm]{\pathPic/PaperOptAniso/ConformTiling.pdf}
		\hspace{1cm}
		\includegraphics[width=4cm,height=4cm]{\pathPic/PaperOptAniso/TriangleClasses.pdf}
	\caption{\label{FRAfigLocPatch} Illustration de la stratégie des patchs locaux : les régions $\Omega_i$ sont couvertes par des pavages périodiques, qui sont ensuite recollés.}
\end{figure}

\begin{figure}
\centering
\includegraphics[height = 4cm,width = 8cm]{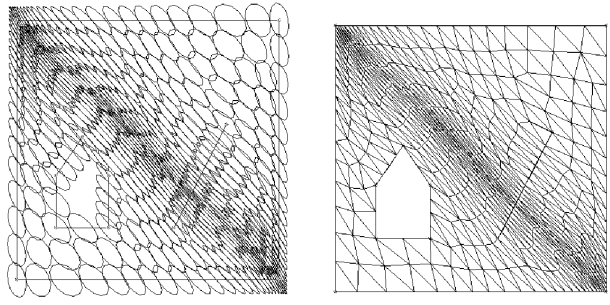}
\caption{Une métrique et une triangulation adaptée. (source : J. Schoen \cite{Schoen})\label{FRAfigEllTri}}
\end{figure}

Cette construction, illustrée Figure \ref{FRAfigLocPatch} est suffisante pour établir le résultat asymptotique \iref{eqCSX}, mais pas pour des applications pratiques car elle ne devient efficace que pour un grand nombre de triangles. L'approche suivante, fondée sur les métriques riemanniennes, est souvent préférée dans les applications. Pour simplifier l'exposition nous supposons que la matrice hessienne $M(z) := d^2 f(z)$ est définie positive en chaque point $z\in \Omega$, et nous définissons 
\be
\label{FRAdefHHessian}
H(z) := \lambda^2 (\det M(z))^{-\frac 1 {2p+2}} M(z)
\ee
où $\lambda>0$ est une constante dont le rôle est de contrôler la résolution de la triangulation. Une telle fonction $H$, qui associe continûment à chaque point $z\in \Omega$ une matrice symétrique définie positive $H(z)\in S_2^+$, est appelée une métrique riemannienne. Notez que pour chaque $z\in \Omega$, la matrice $H(z)$ définit une ellipse
$$
\cE_z := \{u \in \R^2 \sep u^\trans H(z) u\leq 1\}.
$$ 
Comme illustré Figure \ref{FRAfigEllTri}, partie gauche, une métrique encode à chaque point $z\in \Omega$ une information d'aire, de rapport d'aspect et d'orientation, sous la forme d'une matrice symétrique définie positive $H(z)$ ou de manière équivalente d'une ellipse $\cE_z$.
Plusieurs algorithmes de génération de maillage, comme \cite{FreeFem, Inria}, sont capables de produire une triangulation adaptée à la métrique $z\mapsto H(z)$, dans le sens où pour chaque $z\in \Omega$ le triangle $T\in \cT$ contenant $z$ a une forme ``similaire'' à l'ellipse $\cE_z$ comme illustré Figure \ref{FRAfigEllTri}, partie droite.
En termes mathématiques, cela signifie que pour chaque triangle $T\in \cT$ et chaque $z \in T$, l'on a 
\be
\label{FRAtriellipse}
b_T+c_1 \cE_z\subset T \subset b_T+c_2 \cE_z,
\ee
où $0<c_1 <c_2$ sont des constantes fixées et $b_T$ désigne le barycentre de $T$, ce qui signifie aussi que $T$ est ``proche'' d'être un triangle équilatéral dans la métrique $H(z)$.
Si $\cT$ est une telle triangulation, et si $\lambda$ est suffisamment grand, un argument heuristique (qui sera rappelé Chapitre \ref{chapOptAniso}) montre que  
$$
\#(\cT) \|f-\interp_{\cT}^1 f\|_{L^p(\Omega)} \leq C \|\sqrt{|\det (d^2f)|}\|_{L^\tau(\Omega)},
$$
où la constante $C$ dépend de $c_1$ et $c_2$.

En faisant varier le paramètre $\lambda$ nous obtenons différentes triangulations $\cT_\lambda$, de cardinalité proportionnelle à $\lambda^2$, ce qui mène à l'estimation d'erreur \iref{eqCSX}.
Mentionnons que la construction du maillage $\cT$ à partir de la métrique $H$ n'est pas évidente.
De plus il n'existe pas de preuve rigoureuse que la condition de similarité \iref{triellipse} est vérifiée par les algorithmes de génération de maillage les plus courants, à l'exception notable de \cite{Shew} et \cite{Bois}.\\

La borne d'approximation \iref{eqCSX} est 
optimale si l'on se restreint aux triangulations 
qui satisfont une condition technique définie comme suit.
Nous disons qu'une suite $(\cT_N)_{N\geq N_0}$ de triangulations d'un domaine polygonal $\Omega \subset \R^2$ est admissible si $\#(\cT_N) \leq N$ et si 
\be
\label{FRAdefAdmiIntro}
\sup_{N \geq N_0} \left(\sqrt N \max_{T \in \cT_N} \diam(T)\right) <\infty.
\ee
Pour toute suite admissible  $(\cT_N)_{N \geq N_0}$ de triangulations de $\Omega$, on peut établir la minoration 
\be
\label{FRAeqCSXOpt}
\liminf_{N\to \infty} N \|f-\interp^1_{\cT_N} f\|_{L^p(\Omega)} \geq c \left\| \sqrt{|\det (d^2f)|}\right\|_{L^\tau(\Omega)},
\ee
où la constante $c>0$ est universelle (voir aussi Théorème \ref{FRAthLowerLPIntro} dans le plan de la thèse ci dessous). 
De plus la condition d'admissibilité n'est pas trop restrictive : pour chaque $\ve>0$ il existe une suite admissible de triangulations $(\cT_N)_{N \geq N_0}$ qui satisfait la majoration d'erreur \iref{eqCSX} à la constante $\ve$ près ajoutée au terme de droite.\\

Des résultats similaires à \iref{eqCSX} et \iref{eqCSXOpt} peuvent être développés pour les maillages isotropes, dans lesquels la taille des triangles peut varier mais pas leur forme, en ce sens que la mesure de dégénérescence $\rho(T) := \diam(T)^2/|T|$ est uniformément bornée par une constante $\rho_0$, voir par exemple \cite{CMi2}.
Dans ce cas l'estimation \iref{eqCSX} doit être remplacée par 
\be
\label{FRAeqCSXiso}
\limsup_{N\to \infty} N \|f-\interp^1_{\cT_N} f\|_{L^p(\Omega)} \leq C \| d^2f\|_{L^\tau(\Omega)},
\ee
avec la même valeur de $\tau$, et l'estimation
 \iref{eqCSXOpt} a un pendant similaire. Les constantes $C$ et $c$ apparaissant dans ces estimations dépendent désormais de la borne $\rho_0$ sur la mesure de dégénérescence. Ainsi la quantité non linéaire $\sqrt{|\det (d^2 f)|}$ est remplacée par le terme linéaire $d^2f$ dans la norme $L^\tau$, et ces résultats sont désormais très similaires à ceux de meilleure approximation par $N$ ondelettes \cite{Co}.

En termes des valeurs propres $\lambda_1(z), \lambda_2(z)$ de la matrice symétrique $d^2 f(z)$ nous remplaçons donc la moyenne géométrique $\sqrt{|\lambda_1(z) \lambda_2(z)|}$ par $\max \{|\lambda_1(z)|, \, |\lambda_2(z)|\}$, qui peut être significativement plus grand quand ces valeurs propres sont d'ordres de grandeur différents. C'est le cas typiquement si la fonction $f$ approchée présente des caractéristiques fortement anisotropes, et nous pouvons donc nous attendre à une amélioration substantielle des propriétés d'approximation
lorsque des maillages anisotropes sont utilisés pour de telles fonctions.\\

Le résultat \iref{eqCSX} donne un compte rendu précis de l'amélioration que peuvent apporter des triangulations anisotropes en comparaison avec les triangulations isotropes, mais malheureusement seulement dans un cadre restreint, ce qui a motivé notre travail:
\begin{enumerate}[I.]
\item Le résultat original ne s'applique qu'à l'erreur d'interpolation linéaire mesurée en norme $L^p$, alors que les éléments finis de plus haut degré et les normes de Sobolev $W^{1,p}$ sont aussi pertinents. En particulier la norme $W^{1,2}$ (ou $H^1$) apparaît très naturellement dans le contexte des EDP elliptiques.
\item La fonction approchée $f$ doit être $C^2$, alors que les applications les plus intéressantes de l'approximation adaptative font intervenir des fonctions non lisses voire discontinues. Quel sens peut-on 
donner à la quantité $\|\sqrt{|\det(d^2 f)|}\|_{L^\tau(\Omega)}$ lorsque $f$ est une fonction discontinue?
\item La métrique riemannienne $z\mapsto H(z)$ est utilisée comme objet intermédiaire pour la génération de maillages dans les applications numériques. Cette approche manque cependant d'un résultat précis d'équivalence entre ces objets continus et les différentes classes de triangulations anisotropes. Sous quelles conditions peut-on associer à une métrique $z\mapsto H(z)$ une triangulation $\cT$ qui lui est adaptée au sens de \iref{triellipse}?
\item Les algorithmes précédemment mentionnés de génération de maillages anisotropes ne sont pas hiérarchiques, dans le sens où une meilleure précision n'est pas atteinte par un 
raffinement local mais par la re-génération globale d'un nouveau maillage. Peut-on proposer un algorithme de raffinement hiérarchique qui permet d'obtenir l'erreur d'approximation optimale \iref{eqCSX}?
\end{enumerate}

Nous devons aussi mentionner deux problèmes fondamentaux qui sont discutés dans cette thèse, mais qui sont restés des problèmes ouverts et feront l'objet de travaux futurs. En premier lieu l'estimation d'erreur \iref{eqCSX} ne donne qu'une information asymptotique, lorsque le nombre de triangles tend vers $\infty$, alors qu'une estimation d'erreur portant sur toutes les valeurs de $N$ est fortement souhaitée. En second lieu l'extension de ce résultat à des fonctions définies sur des domaines de dimension supérieure à deux est entravée par un problème difficile de géométrie algorithmique : la production de {\it maillages conformes et anisotropes} en dimension $3$ ou plus élevée. Des méthodes numériques comme \cite{Inria} s'attaquent à ce problème, mais il n'est pas résolu d'un point de vue théorique.

\DontWriteThisInToc
\section{Plan de la thèse}

La thèse est constituée de quatre parties qui visent à résoudre les quatre problèmes clefs, numérotés I à IV, que nous avons rencontrés dans les résultats antérieurs sur l'approximation adaptative et anisotrope par éléments finis.
Ces quatre parties sont essentiellement auto-consistantes et peuvent donc être lus indépendamment.
Les chapitres constituant chaque partie gagnent à être lus dans l'ordre, à l'exception du chapitre 1 dont la 
lecture peut être ignorée par le lecteur souhaitant aller plus rapidement au coeur du sujet. Les chapitres
1, 2, 3, 4, 7, 8, ainsi que la troisième partie du chapitre 9, sont respectivement issus des articles 
\cite{BLM}, \cite{Mi}, \cite{Mi2}, \cite{CM2}, \cite{CDHM}, \cite{CM} 
et \cite{CMi2}. Les notations utilisées dans 
ces articles ont été unifiés pour la clarté de l'ensemble. Les chapitres 5 et 6 
sont l'objet d'articles en préparation.

\subsection*{Partie I. \ Eléments finis de degré arbitraire, et normes de Sobolev} 

Nous généralisons dans cette partie l'estimation d'erreur asymptotique \iref{eqCSX} aux éléments finis de degré arbitraire et aux semi-normes de Sobolev $W^{1,p}$ pour la mesure de l'erreur. Cette analyse nous amène à introduire des concepts clefs pour l'adaptation optimale de maillage.\\

\begin{figure}
\centerline{
\includegraphics[width=4cm,height=4cm]{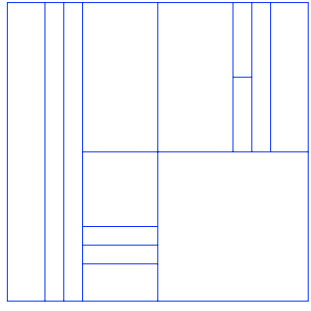}
}
\caption{Partition anisotrope d'un domaine en rectangles.\label{FRAfigRectPart}}
\end{figure}

Pour commencer nous considérons dans le Chapitre \ref{chapBlocks} des partitions d'un domaine rectangulaire en rectangles alignés avec les axes de coordonnées, comme illustré Figure \ref{FRAfigRectPart}, à la place de triangles de directions arbitraires. De telles partitions sont pertinentes lorsque les axes de coordonnées jouent un rôle privilégié, de sorte que les traits anisotropes de la fonction $f$ sont alignés avec les axes de coordonnées. Nous obtenons une estimation d'erreur asymptotique optimale dans ce contexte. Nos résultats s'appliquent à des polynômes par morceaux de degré arbitraire, en dimension quelconque $d>1$, lorsque l'erreur d'approximation est mesurée en norme $L^p$. Nous ne considérons pas ici les normes $W^{1,p}$ car ces approximations polynomiales par morceaux sont généralement discontinues.

Le principal avantage de ce cadre est que les détails techniques requis pour la construction d'une partition anisotrope du domaine, ainsi que l'analyse d'erreur, sont simplifiés par la présence de directions privilégiées.
Nous tirons avantage de ce contexte simple pour introduire et étudier un concept clef appelé la {\it fonction de forme}, ou {\it shape function} en anglais, qui gouverne l'erreur d'approximation locale après une adaptation optimale des éléments de la partition aux propriétés locales de la fonction approchée. Cette fonction de forme est aussi définie et utilisée dans les Chapitres 2 et 3 pour des triangulations anisotropes. Nous donnons ci-dessous sa définition précise dans ce cadre. \\

Le Chapitre \ref{chapOptAniso} est consacré aux éléments finis triangulaires de degré arbitraire $m-1$, où $m\geq 2$, lorsque l'erreur est mesurée en norme $L^p$.
Pour présenter nos résultats, nous devons introduire quelques notations.
Nous notons $\H_m$ l'espace vectoriel des polynômes homogènes de degré $m$:$$
\H_m := \Span\{x^k y^l \sep k+l = m\}.
$$
Un ingrédient clef de notre approche est la fonction de forme $K_{m,p} : \H_m \to \R_+$, où $1\leq p \leq \infty$ est un exposant donné. Cette fonction est définie par une optimisation 
de l'erreur $L^p$ d'interpolation parmi les triangles d'aire $1$ 
de toutes les formes possibles: pour tout $\pi\in \H_m$,
\be
\label{FRAdefKMPIntro}
K_{m,p}(\pi) := \inf_{|T| = 1} \|\pi-\interp_T^{m-1} \pi\|_{L^p(T)},
\ee
où $\interp_T^{m-1}$ désigne l'opérateur d'interpolation locale sur $T$.
Notre résultat principal est l'estimation d'erreur asymptotique suivante.
\begin{FRAtheoremIntro}
\label{FRAthOptiLPIntro}
Soit $\Omega\subset \R^2$ un domaine polygonal borné, soit $f\in C^m(\overline \Omega)$ et soit $1\leq p < \infty$. Il existe une suite $(\cT_N)_{N \geq N_0}$, $\#(\cT_N)\leq N$, de triangulations conformes de $\Omega$ telles que 
\be
\label{FRAeqLPIntro}
\limsup_{N \to \infty} N^{\frac m 2}\|f-\interp_{\cT_N}^{m-1} f\|_{L^p(\Omega)} \leq \left\|K_{m,p}\left(\frac{d^m f}{m!}\right) \right\|_{L^\tau(\Omega)}
\ee
où l'exposant $\tau\in (0, \infty)$ est défini par $\frac 1 \tau := \frac m 2 + \frac 1 p$.
\end{FRAtheoremIntro}

Dans l'estimation d'erreur \iref{eqLPIntro}, nous identifions la donnée $d^m f(z)$ des dérivées d'ordre $m$ de $f$ au point $z$ au polynôme homogène qui lui correspond dans le développement de Taylor de $f$ en $z$. Sous forme mathématique 
$$
\frac {d^m f(z)}{m!} \sim \sum_{k+l = m} \frac{\partial^m f}{\partial^k x\, \partial^l y}(z) \, \frac {x^k}{k!} \frac {y^l}{l!}.
$$

Le Théorème \ref{FRAthOptiLPIntro} étend le résultat connu \iref{eqCSX} aux éléments finis de degré arbitraire $m-1$. De manière similaire, la qualité de l'approximation adaptative et anisotrope de $f$ est déterminée par une expression non-linéaire des dérivées de $f$: la fonction de forme $K_{m,p}(d^m f(z))$ est la ``généralisation'' aux dérivées d'ordre supérieur du déterminant $\sqrt{|\det (d^2f(z))|}$ apparaissant dans \iref{eqCSX}.

Ces résultats motivent une étude approfondie de la fonction de forme $K_{m,p}$. Nous avons montré que $K_{2,p}$, qui correspond au cas $m=2$ de l'approximation affine par morceaux, est proportionnel à la racine carrée du déterminant du polynôme quadratique $\pi = a x^2 + 2 b xy + cy^2\in \H^2$:
$$
K_{2,p}(\pi) = c_{2,p} \sqrt{|\det \pi|} = c_{2,p}\sqrt {|ac-b^2|},
$$
où la constante $c_{2,p}>0$ dépend seulement du signe de $\det \pi$. Nous retrouvons donc le résultat antérieur \iref{eqCSX}. Dans le cas $m=3$ de l'approximation quadratique par morceaux, nous montrons que 
la fonction de forme $K_{3,p}$ est la racine quatrième du discriminant du polynôme cubique homogène $\pi = a x^3+3bx^2y+3cxy^2+dy^3\in \H_3$:
$$
K_{3,p}(\pi) = c_{3,p} \sqrt[4]{|\disc \pi|} = c_{3,p}\sqrt[4]{|4(ac-b^2)(bd-c^2) - (ad-bc)^2|},
$$
où la constante $c_{3,p}>0$ dépend seulement du signe de $\disc \pi$. 
Pour les plus grandes valeurs de $m$, $m\geq 4$, nous n'avons pas obtenu d'expression explicite de la fonction de forme, mais une quantité explicite qui lui est uniformément équivalente.
Cet équivalent a la forme suivante : il existe un polynôme $Q_m(a_0, \cdots, a_m)$ des $m+1$ variables $a_0, \cdots, a_m$, et une constante $C_m\geq 1$ telle que pour tout $\pi= a_0 x^m+ a_1 x^{m-1}y+ \cdots a_m y^m\in \H_m$ on ait en notant $r_m := \deg Q_m$,
\be
\label{FRAequivKQ}
C_m^{-1} \sqrt[r_m]{Q_m(a_0, \cdots , a_m)} \leq K_{m,p}(\pi) \leq  C_m\sqrt[r_m]{Q_m(a_0, \cdots , a_m)}.
\ee
Le polynôme $Q_m$ s'obtient à l'aide de la théorie des polynômes invariants développée par Hilbert dans \cite{Hilbert}. Nous caractérisons aussi les zéros de la fonction de forme, et donc les cas possibles de ``super-convergence'': $K_{m,p}(\pi) = 0$ si et seulement si $\pi$ se factorise par un facteur linéaire $ax+by$ de multiplicité $s>m/2$, en d'autres termes si le polynôme homogène $\pi$ est suffisamment dégénéré.\\

La preuve du Theorème \ref{FRAthOptiLPIntro} est fondée sur la ``stratégie des patchs locaux''
qui a été évoquée précédemment: on considère en premier lieu une ``macro-triangulation'' $\cR$ du domaine initial $\Omega$, qui est suffisamment fine pour que les dérivées d'ordre $m$ de $f$ varient peu sur chaque triangle $R\in \cR$ autour d'une valeur moyenne $\pi_R\in \H_m$. A chaque polynôme $\pi_R$, $R\in \cR$, on associe ensuite un triangle $T_R$ qui minimise, ou presque, le problème d'optimisation définissant $K_{m,p}(\pi_R)$. On pave ensuite chaque ``macro-triangle'' $R\in \cR$ de manière périodique en utilisant le triangle $T_R$ convenablement mis à l'échelle et son symétrique par rapport à l'origine. Finalement, comme illustré Figure \ref{FRAfigLocPatch}, la triangulation $\cT_N$ est obtenue en recollant ensemble les pavages périodiques définis sur chaque $R\in \cR$, à l'aide de quelques triangles supplémentaires aux interfaces pour obtenir un maillage globalement conforme.\\

Le théorème suivant établit que l'estimation asymptotique \iref{eqLPIntro} est optimale, 
si l'on se restreint aux suites admissibles de triangulations qui sont définies par
la condition \iref{defAdmiIntro}. Ce théorème 
montre de plus que la condition d'admissibilité n'est pas trop restrictive.

\begin{FRAtheoremIntro}
\label{FRAthLowerLPIntro}
Soit $\Omega\subset \R^2$ un domaine polygonal borné, soit $f\in C^m(\overline \Omega)$ et soit $1\leq p \leq \infty$. Soit $(\cT_N)_{N \geq N_0}$, $\#(\cT_N) \leq N$, une suite admissible de triangulations de $\Omega$. Alors 
\be
\label{FRAeqLPLowerIntro}
\liminf_{N \to \infty} N^{\frac m 2} \|f-\interp_{\cT_N}^{m-1} f\|_{L^p(\Omega)} \geq \left\|K_{m,p}\left(\frac{d^m f}{m!}\right) \right\|_{L^\tau(\Omega)},
\ee
où $\frac 1 \tau := \frac m 2+ \frac 1 p$. De plus pour chaque $\ve>0$ il existe une suite admissible $(\cT_N^\ve)_{N \geq N_0}$ de triangulations de $\Omega$, $\#(\cT_N) \leq N$, telle que:
\be
\label{FRAeqLPEpsIntro}
\limsup_{N \to \infty} N^{\frac m 2} \|f-\interp_{\cT^\ve_N}^{m-1} f\|_{L^p(\Omega)} \leq \left\|K_{m,p}\left(\frac{d^m f}{m!}\right) \right\|_{L^\tau(\Omega)}+\ve.
\ee
\end{FRAtheoremIntro}

Le Chapitre \ref{chapW1P} est consacré aux versions des Théorèmes \ref{FRAthOptiLPIntro} et \ref{FRAthLowerLPIntro} lorsque l'erreur d'interpolation est mesurée dans la semi-norme de Sobolev $W^{1,p}$, $1\leq p < \infty$. Ces estimations font intervenir l'analogue $L_{m,p}$ de la fonction de forme $K_{m,p}$ qui est défini comme suit: pour tout $\pi\in \H_m$
$$
L_{m,p}(\pi) := \inf_{|T|=1} \|\nabla (\pi-\interp_T^{m-1} \pi)\|_{L^p(T)}.
$$
Nous donnons de nouveau des équivalents explicites de $L_{m,p}(\pi)$, définis de la même manière que \iref{equivKQ} par une expression algébrique en les coefficients de $\pi\in \H_m$. Nos résultats pour les semi-normes $W^{1,p}$ sont donc extrêmement similaires aux résultats obtenus pour les normes $L^p$.

\begin{figure}
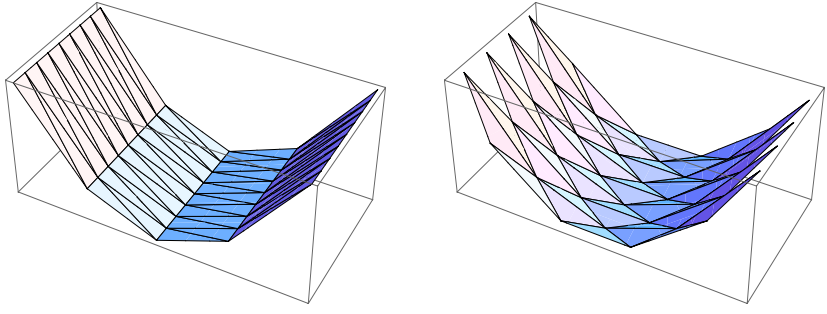

\centerline{
\includegraphics[width=5cm,height=4cm]{\pathPic/Triangles/ParaboleAigue.pdf}
\hspace{0.5cm}
\includegraphics[width=5cm,height=4cm]{\pathPic/Triangles/ParaboleObtuse.pdf}
}
\caption{Interpolation d'une fonction parabolique bi-dimensionelle sur un maillage composé de triangles aigus (gauche), ou bien ``plats'' et fortement obtus (droite).\label{FRAfigParaAO}}
\end{figure}

L'adaptation des preuves n'est en revanche pas évidente, à cause du phénomène suivant: la présence de triangles {\it fortement obtus} (avec un angle proche de $\cPi$) dans un maillage peut causer des oscillations du gradient de l'interpolation d'une fonction, comme illustré Figure \ref{FRAfigParaAO}. Ce phénomène détériore l'erreur d'interpolation dans la semi-norme $W^{1,p}$, mais pas dans la norme $L^p$. Ces triangles ``plats'' doivent donc être évités avec précaution. En résumé, les triangles longs et fins peuvent être souhaitables mais il ne doivent pas être trop fortement obtus.

Avant de continuer la description de cette thèse, nous rappelons au lecteur que les trois chapitres qui composent la Partie \ref{partAsymptApprox} sont auto-consistants et peuvent donc être lus indépendamment.

\subsection*{Partie II. \ Classes d'approximation anisotropes et modèles d'images}

Dans cette partie, formée de l'unique Chapitre 4, nous discutons de l'extension de nos résultats d'approximation aux fonctions non lisses.\\

Il existe des moyens variés de mesurer la régularité de fonctions définies sur un domaine $\Omega\subset \R^2$, le plus souvent au moyen d'un {\it espace de régularité} approprié et d'une norme associée. Des exemples classiques sont les espaces de Sobolev et de Besov. Ces espaces sont souvent utilisés pour décrire la régularité de solutions d'EDP. D'un point de vue numérique ils caractérisent précisément la vitesse à laquelle une fonction $f$ peut être approchée par des fonctions plus simples comme les séries de Fourier, les éléments finis (sur des triangulations isotropes), les fonctions splines ou les ondelettes.

Le résultat d'approximation adaptative anisotrope \iref{FRAeqCSX} et sa généralisation
par le Théorème \ref{FRAthOptiLPIntro} indiquent que la qualité de l'approximation d'une fonction $f$ par des éléments finis sur des triangulations anisotropes est gouvernée par une quantité non-linéaire de ses dérivées, d'un point de vue asymptotique du moins. Dans le cas des éléments finis de degré $1$, et de l'approximation en norme $L^2$, la quantité pertinente est la suivante
$$
A(f) := \left\|\sqrt {|\det(d^2f)|}\right\|_{L^{2/3}(\Omega)}.
$$
La fonctionnelle $A$ diffère fortement des normes de Sobolev, Hölder ou Besov, car elle est fortement non linéaire: $A$ ne satisfait pas l'inégalité triangulaire, ni aucune quasi-inégalité triangulaire. En d'autres termes pour chaque constante $C$ il existe $f,g\in C^2(\overline \Omega)$ telles que 
\be
\label{FRAeqNoTri}
A(f+g) > C(A(f)+A(g)).
\ee
L'absence d'une inégalité triangulaire interdit d'utiliser les techniques classiques de l'analyse linéaire pour définir un espace de régularité attaché à la fonctionnelle $A$.
L'extension du résultat d'approximation \iref{eqCSX} aux fonctions qui ne sont pas $C^2$ n'est donc pas évidente. \\

Les fonctions qui se présentent dans les applications concrètes, comme par exemple en traitement de l'image ou comme solutions d'EDP hyperboliques, présentent souvent des zones de régularité séparées par des discontinuités localisées.
Un modèle mathématique simple pour ce type de comportement est donné par la classe des fonctions cartoon, qui sont régulières excepté le long d'une famille de courbes elles mêmes régulières, à travers lesquelles elles sont discontinues.
Une analyse heuristique présentée dans le Chapitre \ref{chapCartoon} suggère que pour toute fonction cartoon $f$ définie sur un domaine polygonal borné $\Omega$, il existe une suite $(\cT_N)_{N \geq N_0}$ de triangulations anisotropes de $\Omega$, $\#(\cT_N)\leq N$, telles que 
\be
\label{FRAeqApproxCartoon}
\sup_{N\geq N_0} N \|f-\interp_{\cT_N}^1 f\|_{L^2(\Omega)} < \infty. 
\ee
Comme illustré Figure \ref{FRAfigEllTri0}, les triangulations $(\cT_N)_{N \geq N_0}$ se composent de triangles fortement anisotropes alignés avec les discontinuités de $f$, et de grands triangles dans les régions où $f$ est régulière.
Le résultat d'approximation \iref{eqApproxCartoon} fait espérer qu'il existe une estimation d'erreur asymptotique précise, pour l'approximation anisotrope des fonctions cartoon, qui étende le résultat \iref{eqCSX} connu lorsque $f$ est $C^2$, à savoir
\be
\label{FRAeqApproxA}
\limsup_{N \to \infty} N \|f-\interp_{\cT_N}^1 f\|_{L^2(\Omega)} \leq C A(f).
\ee

Nous n'avons rempli pour l'instant qu'une partie de ce programme: l'extension de la fonctionnelle $A$ aux fonctions cartoon. Plus précisément, considérons une fonction $\vp\in C^\infty(\R^2)$ radiale, à support compact et d'intégrale unité. Définissons $\vp_\delta(z) := \frac 1 {\delta^2} \vp\left(\frac z \delta\right)$ pour chaque $\delta>0$ et notons $f_\delta := f*\vp_\delta$ la convolution de $f$ avec $\vp_\delta$.
Nous prouvons que si $f$ est une fonction cartoon, alors $A(f_\delta)$ converge lorsque $\delta\to 0$ vers une expression explicite:
\begin{eqnarray*}
\lim_{\delta \to 0} A(f_\delta)^{2/3} &=& \left\|\sqrt{|\det(d^2 f)|}\right\|_{L^{2/3}(\Omega\sm\Gamma)}^{2/3} + C(\vp) \left\|\sqrt{|\kappa|} \gamma\right\|_{L^{2/3}(\Gamma)}^{2/3}\\ 
&=& \int_{\Omega\sm\Gamma} \sqrt[3]{|\det (d^2 f(z))|}dz+ C(\vp)\int_\Gamma |\kappa(s)|^{1/3} |\gamma(s)|^{2/3}ds.
\end{eqnarray*} 
Nous avons noté ici $\Gamma$ la famille de courbes le long desquelles la fonction cartoon $f$ est discontinue, $\gamma(s)$ le saut de $f$ en un point $s\in \Gamma$, et $\kappa(s)$ la courbure de $\Gamma$ en $s$. La constante $C(\vp)$ est strictement positive et s'exprime explicitement en fonction de $\vp$.\\

L'extension de la fonctionnelle non linéaire $A$ aux fonctions cartoon met en lumière un lien, exploré en \S\ref{secCartoon4}, entre l'approximation par éléments finis sur des triangulations anisotropes et plusieurs autres domaines mathématiques. Nous pensons en particulier aux méthodes présentées dans \cite{Ca} de traitement de l'image invariant par transformation affine, et à la définition donnée dans \cite{DPW} d'espaces de fonctions définis par la régularité des lignes de niveau de leurs éléments.\\

\begin{figure}
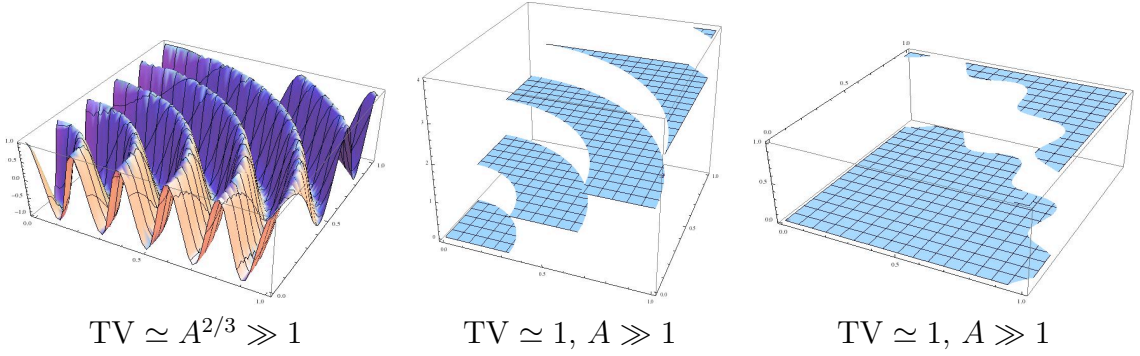

\centering
\begin{tabular}[c]{ccc}
\includegraphics[height=3.5cm,width=5cm]{\pathPic/Bayesien/SinaNoText.pdf}
&\includegraphics[height=4cm,width=4cm]{\pathPic/Bayesien/StairaNoText.pdf}
&\includegraphics[height=3.5cm,width=5cm]{\pathPic/Bayesien/FloweraNoText.pdf}\\
 $\TV \simeq A^{2/3} \gg 1$
& $\TV \simeq 1$, $A \gg 1$
& $\TV \simeq 1$, $A \gg 1$
\end{tabular}
\caption{Comportement des fonctionnelles $\TV$ et $A$ sur différents types d'images.\label{FRAfigTVA}}
\end{figure}
On peut aussi penser la quantité $A$ comme le
pendant ``d'ordre deux'' de la semi-norme $\TV$ de variation totale, une mesure de régularité définie en termes des dérivées d'ordre un de $f$ et qui est aussi finie lorsque $f$ est une fonction cartoon. La variation totale joue un rôle central en traitement de l'image et dans l'analyse des équations de transport, deux domaines dans lesquels les fonctions régulières par morceaux apparaissent naturellement. Lorsque $f$ est une fonction cartoon, sa variation totale est donnée par la formule suivante:
$$
\TV(f) = \int_{\Omega\sm\Gamma} |\nabla f(z)|dz+\int_\Gamma |\gamma(s)| ds.
$$
Nous comparons le comportement des quantités $A(f)$ et $\TV(f)$ pour différentes familles de fonctions cartoon $f$. Les quantités $\TV(f)$ et $A(f)^{2/3}$ se révèlent équivalentes lorsque $f$ est la fonction oscillant de manière lisse $f(z) := \cos(\omega |z|)$, illustrée Figure \ref{FRAfigTVA} (i), où $\omega$ est un grand paramètre.
En revanche, les discontinuités sont pénalisées de manière différente par ces deux fonctionnelles. Pour une fonction ``en escalier'', comme illustré Figure \ref{FRAfigTVA} (ii), la variation totale $\TV$ reste bornée alors que $A$ tend vers l'infini à mesure que le nombre de marches croît, à cause du terme de saut $|\gamma(s)|^{2/3}$ dans l'intégrale sur l'ensemble $\Gamma$ des discontinuités. Par ailleurs, à cause du terme de courbure $|\kappa(s)|^\frac 1 3$, la fonctionnelle $A$ est bien plus grande que $\TV$ pour les fonctions caractéristiques d'ensembles ayant un bord complexe ou oscillant comme illustré Figure \ref{FRAfigTVA}.

La fonctionnelle $A$ peut donc être considérée comme un modèle d'image quantitatif : une image monochrome, décrite par sa luminosité $f:[0,1]^2 \to [0,1]$, est plausible si $A(f)$ est suffisamment petit.
A partir de l'approche introduite dans \cite{LMo}, 
nous proposons un algorithme de débruitage d'images utilisant un a-priori bayésien fondé sur ce modèle.
Dans sa version actuelle, cet algorithme n'est pas satisfaisant en terme
de vitesse de convergence, et pour cette raison nous présentons uniquement 
des illustrations numériques dans un cadre simplifié unidimensionel. \\

Enfin nous étudions l'extension aux fonctions cartoon des autres quantités non-linéaires qui apparaissent en approximation par éléments finis sur des triangulations anisotropes, comme la norme $\|K_{m,p}(d^m f)\|_{L^\tau}$ de la fonction de forme pour $m\geq 2$, ou l'analogue de cette quantité lorsque $f$ est une fonction de plus de deux variables.

\subsection*{Partie III. \ Génération de maillage anisotrope et métriques riemanniennes}

Les triangulations sont des objets discrets de nature combinatoire : elles peuvent être décrites par une famille de sommets et d'arêtes les joignant. Cette description est fructueuse pour la démonstration de résultats algébriques comme la formule d'Euler, ou pour le traitement informatisé.
Par contraste, comme expliqué précédemment, de nombreuses approches en adaptation anisotrope de maillage \cite{Shew,Bois,A} se fondent sur un objet continu équivalent aux triangulations, à savoir une métrique riemannienne $z\mapsto H(z)$. Il s'agit en d'autres termes d'une fonction continue $H$ de $\Omega$ dans l'ensemble $S_2^+$ des matrices symétriques définies positives. Une fois cette métrique conçue convenablement, un algorithme de génération de maillages a la charge de fabriquer une triangulation qui lui correspond dans le sens de \iref{triellipse}.
L'objectif de la Partie \ref{partRiemann} est de formuler des résultats précis d'équivalence entre certaines classes de triangulations et de métriques riemanniennes. Cette équivalence traduit les contraintes géométriques que satisfont les triangulations sous la formes de propriétés de régularité des métriques riemanniennes qui leur sont équivalentes.\\

Pour énoncer nos résultats nous devons introduire certaines notations. Nous associons à chaque triangle $T$ son barycentre $z_T\in \R^2$ et la matrice symétrique définie positive $\cH_T\in S_2^+$ telle que l'ellipse $\cE_T$ définie par 
$$
\cE_T := \{z\in \R^2 \sep (z-z_T)^\trans \cH_T (z-z_T)\},
$$
est l'ellipse d'aire minimale contenant $T$.
Le point $z_T$ indique donc la position de $T$, tandis que la matrice $\cH_T\in S_2^+$ décrit
son aire, son rapport d'aspect et son orientation.

\begin{figure}
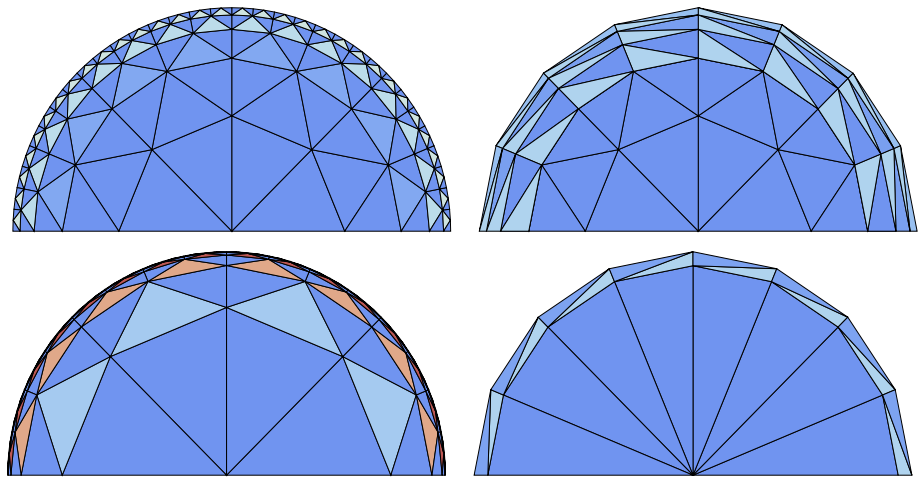

	\centering
		\includegraphics[width=6cm,height=3.2cm]{\pathPic/CoucheLimite/DiskIso.pdf}
		\includegraphics[width=6cm,height=3.2cm]{\pathPic/CoucheLimite/DiskAigu.pdf}
		\includegraphics[width=6cm,height=3.2cm]{\pathPic/CoucheLimite/DiskGen.pdf}
		\includegraphics[width=6cm,height=3.2cm]{\pathPic/CoucheLimite/DiskNoReg.pdf}
	\caption{\label{FRAfigDisk} Génération d'une couche de triangles de largeur $\delta$ le long du cercle unité en utilisant $\cO(\delta^{-1})$ triangles isotropes (haut gauche), $\cO(\delta^{-\frac 1 2} |\ln \delta|)$ 
triangles ``quasi-aigus'' (haut droit), $\cO(\delta^{-\frac 1 2})$ triangles satisfaisant, ou non, la condition d'étagement (bas gauche et bas droit).}
\end{figure}

Nous notons $\bT$ la famille de toutes les triangulations conformes du domaine infini $\R^2$. Le choix de considérer des triangulations infinies est guidé par la simplicité, et un travail futur sera consacré aux triangulations de domaines polygonaux bornés.
Nous notons $\bH := C^0(\R^2, S_2^+)$ la famille des métriques riemanniennes sur $\R^2$. Une métrique $H\in \bH$ associe continûment à chaque point $z\in \R^2$ une matrice symétrique définie positive $H(z)\in S_2^+$. Soit $C\geq 1$, nous disons qu'une triangulation $\cT\in \bT$ est $C$-équivalente à une métrique $H\in \bH$ si pour tout $T\in \cT$ et tout $z\in T$ nous avons au sens des matrices symétriques
$$
C^{-2} H(z) \leq \cH_T \leq C^2 H(z).
$$
Nous disons qu'une famille $\bT_*\subset \bT$ de triangulations est équivalente à une famille $\bH_* \subset \bH$ de métriques s'il existe une constante uniforme $C\geq 1$ telle que 
\begin{itemize}
\item 
Pour chaque triangulation $\cT\in \bT_*$ il existe une métrique $H\in \bH_*$ telle que $\cT$ et $H$ soient $C$-équivalents.
\item 
Pour chaque métrique $H\in \bH_*$ il existe une triangulation $\cT\in \bT_*$ telle que $\cT$ et $H$ soient $C$-équivalents.\\
\end{itemize}
Nous considérons trois familles pertinentes de triangulations de $\R^2$
$$
\bT_{i,C} \subset \bT_{a,C} \subset \bT_{g,C}
$$
qui dépendent de manière mineure d'un paramètre $C\geq 1$.
La famille $\bT_{g,C}$ des triangulation {\it étagées}, en anglais {\it graded}, est définie par la condition suivante qui impose un minimum de consistance dans les formes des triangles voisins.
Une triangulation $\cT$ appartient à $\bT_{g,C}$ si pour tous $T,T'\in \cT$ nous avons 
$$
T \cap T' \neq \emptyset \quad \Ra \quad C^{-2} \cH_T \leq \cH_{T'} \leq C^2 \cH_T.
$$
Une triangulation appartient à la classe $\bT_{i,C}$ des triangulations {\it isotropes} si $\cT$ est étagée, $\cT\in \bT_{g,C}$, et si les éléments de $\cT$ sont suffisamment proches du triangle équilatéral, ce qu'exprime la condition suivante: pour tout $T\in \cT$
$$
\|\cH_T\| \|\cH_T^{-1}\|\leq C^2.
$$
La classe $\bT_{a,C}$ des triangulations {\it quasi-aigues} est définie par une condition légèrement plus technique sur les angles maximaux des triangles, voir Chapitre \ref{chapMeshMet}.
Comme expliqué dans la description ci-dessus du chapitre 3, éviter les angles trop fortement obtus est nécessaire pour garantir la stabilité du gradient lorsqu'on applique l'opérateur d'interpolation.
L'un de nos résultats clefs est la reformulation de cette condition sous la forme d'une hypothèse de régularité de la métrique riemannienne équivalente.
D'un point de vue pratique, cette condition n'est malheureusement pas garantie par les programmes existants de génération de maillage. Les contraintes satisfaites par ces familles de triangulations de $\R^2$ sont illustrées 
pour des triangulations de domaines bornés sur la Figure \ref{FRAfigDisk}.\\

Les résultats du Chapitre \ref{chapMeshMet} établissent que lorsque $C$ est suffisamment grand, ces trois classes de triangulations sont équivalentes à trois familles de métriques, respectivement
$$
\bH_i \subset \bH_a \subset \bH_g,
$$
qui sont définies par des conditions de régularités précises sur la fonction $z\mapsto H(z)$.
En particulier la famille $\bH_i$ est formée des métriques $H$ qui sont proportionnelles à l'identité
$
H(z) = \Id /s(z)^2 ,
$
et telles que le facteur de proportionnalité $s$ satisfait l'une des conditions suivantes qui de manière surprenante sont {\it équivalentes}
\begin{enumerate}[$\bullet$]
\item Propriété de Lipschitz Euclidienne: $|s(z) - s(z')| \leq |z-z'|$ pour tous $z,z'\in \R^2.$
\item Propriété de Lipschitz Riemannienne: $|\ln s(z) - \ln s(z')| \leq d_H(z,z')$ pour tous $z,z'\in \R^2.$
\end{enumerate}
Rappelons que la distance riemannienne $d_H(z,z')$ mesure la longueur, au sens de la métrique $H$, du plus court chemin joignant $z$ à $z'$:
$$
d_H(z,z') := \inf_{\substack{\gamma(0) = z\\ \gamma(1) = z'}} \int_0^1 \|\gamma'(t) \|_{H(\gamma(t))} dt
$$
où $\|u\|_M := \sqrt{u^\trans M u}$ et où l'infimum est pris parmi tous les chemins $\gamma\in C^1([0,1], \R^d)$ joignant $z$ à $z'$.

Les deux propriétés de Lipschitz présentées ci-dessus s'étendent naturellement aux métriques riemanniennes générales, {\it mais ne sont alors plus équivalentes}. Pour $H\in \bH$,
on pose $S(z) := H(z)^{-\frac 1 2}$, et on introduit les deux propriétés distinctes
\be
\label{FRAeqDPlusIntro}
d_+(H(z), H(z') ) \leq  |z-z'| \stext{ pour tout} z,z'\in \R^2,
\ee
où $d_+(H(z), H(z') ) := \|S(z) - S(z')\|$,
et 
\be
\label{FRAeqDTimesIntro}
d_\times(H(z), H(z') ) \leq d_H(z,z') \stext{ pour tout } z,z'\in \R^2.
\ee
où 
$d_\times(H(z), H(z') ) := \ln \max \{\|S(z) S(z')^{-1}\| , \|S(z') S(z)^{-1}\| \}$.
Notons que $d_+$ et $d_\times$ sont des distances sur $S_2^+$.

Le résultat principal du Chapitre \ref{chapMeshMet} caractérise les familles de métriques $\bH_g$ et $\bH_a$ (associées aux familles $\bT_{g,C}$ des triangulations étagées, et $\bT_{a,C}$ des triangulations quasi-aigues) sous la forme des propriétés de Lipschitz ci-dessus. 
\begin{FRAtheoremIntro}
\label{FRAthMeshMetIntro}
Si la constante $C$ est suffisamment grande, alors la famille $\bT_{g,C}$ des triangulations étagées est équivalente à la famille $\bH_g$ des métriques qui satisfont \iref{eqDTimesIntro}, et la famille $\bT_{a,C}$ des triangulations quasi-aigues est équivalente à la famille $\bH_a$ des métriques qui satisfont simultanément \iref{eqDTimesIntro} et \iref{eqDPlusIntro}.
\end{FRAtheoremIntro}

Nous donnons Chapitre \ref{chapApproxMet} des applications de ces résultats dans le contexte de la théorie de l'approximation et de la génération contrainte de maillage. Décrivons ce dernier exemple.
Pour chaque ensemble fermé $E\subset \R^d$ et chaque triangulation $\cT\in \bT$ nous notons $V_\cT(E)$ le voisinage de $E$ dans la triangulation $\cT$, qui est défini comme suit
$$
V_\cT(E) := \bigcup_{\substack{T\in \cT\\ T \cap E \neq \emptyset}} T.
$$
Nous disons qu'une triangulation $\cT\in \bT$ sépare deux ensembles fermés et disjoints $X,Y\subset \R^2$ si $V_\cT(X) \cap V_\cT(Y) = \emptyset$.
La propriété analogue pour les métriques est la suivante : une métrique $H\in \bH$ sépare $X$ et $Y$ si $d_H(x,y)\geq 1$ pour tous $x\in X$ et $y\in Y$. Nous montrons que ces deux propriétés sont rigoureusement équivalentes, et nous utilisons cette reformulation pour calculer le plus petit nombre de triangles (à périodicité près) requis pour séparer des ensembles (périodiques) en utilisant une triangulation (périodique) isotrope, quasi-aigue ou étagée.
Comme illustré sur la Figure \ref{FRAfigDisk}, imposer davantage de contraintes sur la triangulation augmente typiquement le nombre de triangles requis pour réaliser la même tâche.\\

La suite de ce chapitre est consacrée au contrôle de l'erreur d'approximation par éléments finis d'une fonction sur une triangulation $\cT$, en norme $L^p$ ou $W^{1,p}$, par une quantité $e_H(f)_p$, $e_H^a(\nabla f)_p$ ou $e_H^g(\nabla f)_p$ attachée à une métrique $H$ équivalente à $\cT$. Cette analyse fait apparaître, dans le cas des normes de Sobolev, le rôle particulier joué par les conditions d'angle et de régularité qui définissent les triangulations et métriques appartenant à $\bT_a$ et $\bH_a$ respectivement. 
Finalement nous étendons aux métriques les estimations d'erreur asymptotiques développés pour des triangulations Chapitres 2 et 3. 

\subsection*{Partie IV. \ Algorithmes de raffinement hiérarchique}

\begin{figure}
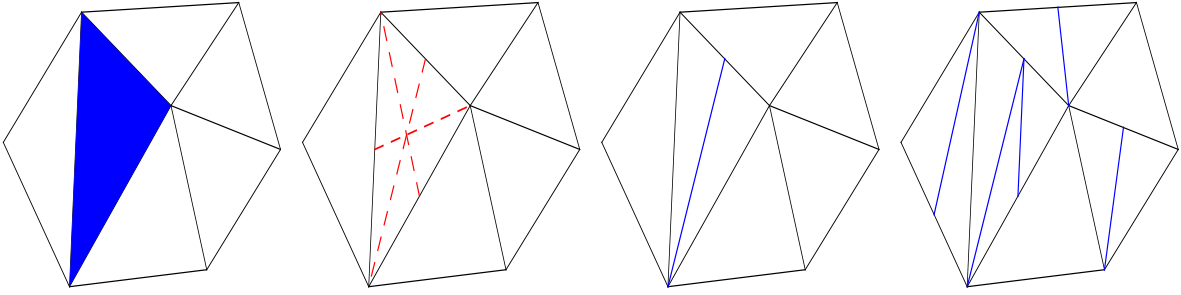

	\centering
		\includegraphics[width=3.8cm,height=4cm]{\pathPic/Subdivision/Algo1.pdf}
		\includegraphics[width=3.8cm,height=4cm]{\pathPic/Subdivision/Algo2.pdf}
		\includegraphics[width=3.8cm,height=4cm]{\pathPic/Subdivision/Algo2p.pdf}
		\includegraphics[width=3.8cm,height=4cm]{\pathPic/Subdivision/Algo3.pdf}
	\caption{\label{FRAfigRefinement}
	L'algorithme de raffinement : choix du triangle qui maximise l'erreur d'approximation (i, en sombre), choix d'une bissection parmi les trois possibilités (ii, iii), itération de ces deux étapes (iv).}
\end{figure}

La dernière partie de cette thèse est consacrée à l'étude d'un algorithme proposé par Cohen, Dyn et Hecht qui produit des suites hiérarchiques $(\cT_N)_{N \geq N_0}$ de triangulations anisotropes (non-conformes) adaptées à une fonction donnée $f$.
Etant donnée une triangulation $\cT$ d'un domaine $\Omega$ et une fonction $f\in L^p(\Omega)$, cet algorithme crée en une étape une triangulation $\cT'$ de $\Omega$, de cardinalité $\#(\cT') = \#(\cT)+1$, de la façon suivante:
\begin{enumerate}
\item (Sélection ``greedy'' du triangle à raffiner) 
On sélectionne un triangle $T\in \cT$ dont la contribution à l'erreur d'approximation est maximale
$$
T := \underset {T'\in \cT} \argmax \, \|f-\cA_{T'} f\|_{L^p(T')},
$$
où $\cA_T : L^p(T) \to \P_{m-1}$ est un opérateur de projection, par exemple la projection $L^2(T)$ orthogonale sur $\P_{m-1}$.
\item 
(Choix d'une bissection) Une arête $e\in \{a,b,c\}$ de $T$ est choisie en minimisant une \emph{fonction de décision} donnée $e\mapsto d_T(e,f)$ parmi les trois arêtes.
Le triangle est découpé le long du segment joignant le point milieu de l'arête choisie au sommet opposé, ce qui crée les sous-triangles $T_e^1$ et $T_e^2$. La nouvelle triangulation est donc
$$
\cT' := \cT-\{T\} + \{T_e^1, T_e^2\}.
$$
\end{enumerate}
Partant d'une triangulation $\cT_0$ de cardinalité $N_0$ du domaine $\Omega$, l'algorithme produit pas après pas une suite $(\cT_N)_{N \geq 0}$, voir Figure \ref{FRAfigRefinement}, de triangulations ``adaptées'' à une fonction donnée $f\in L^p(\Omega)$. Les propriétés de ces triangulations dépendent fortement de la fonction approchée $f$ et du choix de la fonction de décision $d_T(e,f)$, qui guide la création de l'anisotropie. Par contraste l'opérateur de projection $\cA_T$ joue un rôle plutôt mineur. Un choix typique de la fonction de décision $e\mapsto d_T(e,f)$ est l'erreur locale après bissection, en d'autres termes l'algorithme choisit la bissection qui réduit le plus possible l'erreur d'approximation.\\

Nous décrivons cet algorithme plus en détail dans le
Chapitre \ref{chapCDHM}, et nous établissons sa convergence au sens où les approximations polynomiales par morceaux $\cA_{\cT_N} f$ définies par 
$$
\cA_{\cT_N} f(z)=\cA_T(z),\; \; z\in T,\;\; T\in \cT_N
$$
convergent vers $f$ dans $L^p$ lorsque $N \to \infty$ pour n'importe quelle $f\in L^p$, sous certaines hypothèses sur la fonction de décision $e\mapsto d_T(e,f)$.
Nous discutons aussi de la possibilité d'utiliser la structure hiérarchique multi-échelle 
pour définir des approximations multi-résolution, des ondelettes et un algorithme de type CART.
Nous illustrons l'adaptation anisotrope donnée par l'algorithme sur plusieurs types de fonctions et d'images qui présentent des transitions rapides le long de lignes courbes.\\

\begin{figure}
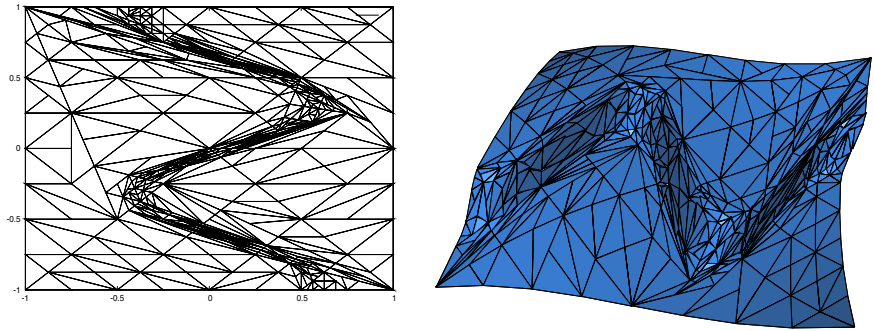

	\centering
		\includegraphics[width=6cm,height=5cm]{\pathPic/Triangles/f1-512.pdf}
		\includegraphics[width=6cm,height=4cm]{\pathPic/Triangles/f1-512_2.pdf}
	\caption{\label{FRAfigSinGreedy}
	Triangulation produite par l'algorithme de raffinement hiérarchique, et interpolation, pour une fonction ayant une variation brusque le long d'une courbe sinusoïdale.}
\end{figure}

Nous faisons une analyse plus approfondie de la convergence de l'algorithme Chapitre \ref{chapBisecOpt}, dans le cas $m=2$ de l'approximation linéaire par morceaux.
Notre résultat principal montre que lorsque 
$f$ est $C^2$ et strictement convexe, pour un choix particulier de la fonction de décision fondé sur l'erreur d'interpolation $L^1$, la suite $(\cT_N)_{N \geq N_0}$ de triangulations générée par cet algorithme satisfait l'estimation asymptotiquement optimale de convergence
$$
\limsup_{N \to \infty} N \| f-\cA_{\cT_N} f\|_{L^p(\Omega)} \leq C \|\sqrt{|\det (d^2f)|}\|_{L^\tau(\Omega)},
$$
où $\frac 1 \tau := 1+\frac 1 p$. L'observation clef qui mène à ce résultat est que lorsque $f$ est un polynôme quadratique convexe, la minimisation de la fonction de décision choisit {\it l'arête la plus longue} dans la métrique associée à la partie homogène quadratique de ce polynôme.
Cette propriété permet de montrer que les triangles générés par l'algorithme tendent en majorité à adopter un rapport d'aspect optimal. Le bon comportement de l'algorithme est aussi observé pour des fonctions générales non-convexes, comme illustré Figure \ref{FRAfigSinGreedy}.
Cependant prouver le résultat de convergence optimale ci-dessus reste un problème ouvert dans ce cadre.\\

Nous étudions Chapitre \ref{chapVariations} des variantes de l'algorithme de raffinement hiérarchique présenté ci-dessus. Nous considérons en premier lieu une fonction de décision fondée sur l'erreur locale de projection $L^2$, qui est particulièrement bien adaptée à l'implémentation numérique car elle peut être évaluée en un temps machine réduit.

Nous nous concentrons ensuite sur le comportement de l'algorithme lorsqu'il est appliqué à des fonctions cartoon. L'algorithme original ne satisfait pas la meilleure estimation possible de convergence pour de telles fonctions. Nous montrons que la vitesse optimale de convergence pourrait être rétablie en remplaçant la bissection du point milieu d'une arête vers le sommet opposé par d'autres choix de découpages géométriques.

Finalement nous considérons une autre variante de l'algorithme, fondée sur des rectangles alignés avec les axes de coordonnées au lieu de triangles de direction arbitraire, dans l'esprit des partitions rectangulaires étudiées Chapitre 1. Cette simplification mène à un résultat qui garantit la meilleure estimation possible de convergence pour toutes les fonctions $C^1$, dans le cas d'approximations constantes par morceaux.

\chapter*{Introduction (English Version)}

\begin{quotation}
{\it
Although this may seem a paradox, all exact science is dominated by the idea of approximation. \hfill (Bertrand Russel, logician and Nobel prize)} 
\end{quotation}

\DontWriteThisInToc
\section{Anisotropic finite element approximation : why and how?} 

This thesis is devoted to the problem of approximating functions by piecewise polynomial
finite elements over triangulations, and more general meshes. We are particularly
concerned with the setting where the mesh is {\it adaptively designed} depending 
on the function to be approximated. This mesh may therefore include elements 
of strongly varying size, aspect ratio and orientation.\\

Approximation by piecewise polynomial functions
is a procedure that occurs in numerous applications. 
In some of them such as 
terrain data simplification, surface or image compression,
the function $f$ to be approximated might be fully known.
In other applications such as denoising, statistical learning or 
in the finite element discretization of PDE's,
it might be only partially known or fully
unknown.
In all these applications, one usually makes the 
distinction between {\it uniform} and {\it adaptive} approximation.
In the uniform case, the domain of interest is decomposed into
a partition where all elements have comparable shape and size,
while these attributes are allowed to vary strongly in the adaptive case.
The partition may therefore be adapted to the local properties
of $f$, with the objective of 
optimizing the trade-off between accuracy and complexity of
the approximation.

\begin{figure}
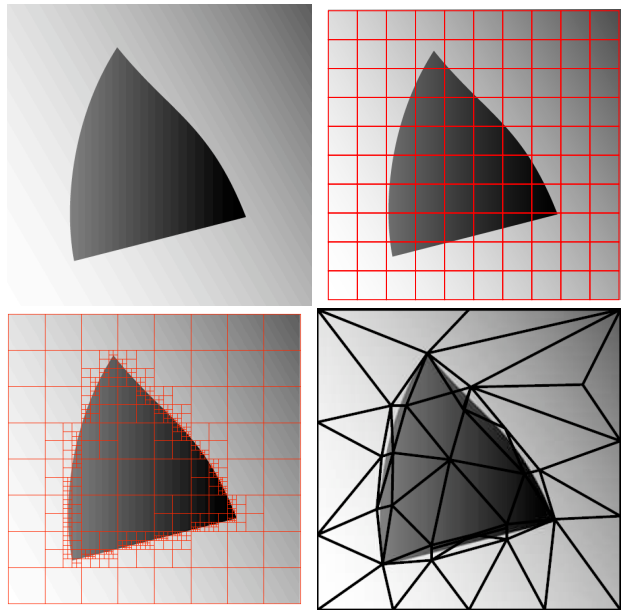

\centering
\includegraphics[height = 4cm,width = 4cm]{\pathPic/Cartoon/Triangle.pdf}
\includegraphics[height = 4cm,width = 4cm]{\pathPic/Cartoon/TriangleUnif.pdf}\\
\includegraphics[height = 4cm,width = 4cm]{\pathPic/Cartoon/Haar.pdf}
\includegraphics[height = 4cm,width = 4cm]{\pathPic/Cartoon/TriangleMeshed.pdf}
\caption{A piecewise smooth function (top left), Uniform partition of the domain (top right), Partition associated to the adaptive approximation based on the largest coefficients
in the Haar wavelet basis (bottom left), Partition associated to an adaptive anisotropic finite element approximation (bottom right). (credit : G. Peyr\'e)\label{figEllTri0}}
\end{figure}

\begin{figure}
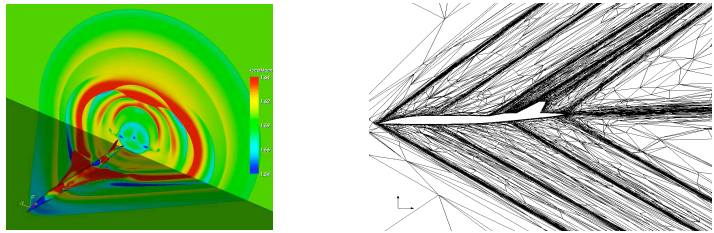

	\centering		
	\includegraphics[width=3.5cm,height=3cm]{\pathPic/SimuExt/AlauzetShock.pdf}
	\hspace{1cm}
	\includegraphics[width=4.5cm,height=3cm]{\pathPic/SimuExt/AlauzetPlane.pdf}
	\caption{Airflow around a supersonic plane, computed using a three dimensional highly anisotropic mesh. (credit : F. Alauzet \cite{A}.\label{figInria})}
\end{figure}

From an approximation theoretic point of view, the trade-off between
accuracy and complexity if usually tied to the smoothness properties of the
function: typically one expects higher convergence rates for smoother functions.
Functions arising in concrete applications may however have inhomogeneous smoothness properties, in the sense that they exhibit areas of smoothness separated by localized discontinuities.
Two typical instances displayed in Figure \ref{figEllTri0} 
and Figure \ref{figInria} are (i) edges in functions representing images, and (ii) shock profiles in the solutions to non-linear hyperbolic PDE's. 
Numerical procedures for image processing, such as image denoising or compression, or for the simulation of PDE's  greatly benefit from economical and faithful approximations of such functions.
A piecewise polynomial approximation on a uniform partition of the domain is generally not sufficient for these purposes. A first step toward adaptivity is to vary the size of the elements forming the partition
according to the local smoothness properties of the function. A very similar procedure,
often used in image processing, consists in retaining the largest terms in a wavelet
decomposition, such as displayed on Figure \ref{figEllTri0} (bottom left). A second step is to observe that a higher resolution is needed in the direction orthogonal to the curve of discontinuity than in the tangential direction, and to take advantage of this property by using an {\it anisotropic} partition of the domain. In the two dimensional cases, such partitions
are typically built from triangles of high aspect ratio aligned with the discontinuities, 
as displayed Figure \ref{figEllTri0} (bottom right) and Figure \ref{figInria}.

In the context of numerical PDE's, adaptivity also refers to the fact that
the computational mesh is not fixed in advance, but instead is dynamically
updated based on the available information on the exact solution
gained through the solution process. From a numerical point of view, 
such methods require more complex algorithms and more intricate
 data structures than their non-adaptive counterparts. 
From a theoretical point of view the analysis of these adaptive
algorithms, when it is possible, is generally involved. As a matter of fact, 
the improvement brought by adaptivity in terms of convergence
rate is rigorously established only for few systems of PDE's,
and in the case of {\it isotropic} finite element meshes. 
We refer to the survey paper \cite{Nochetto} for a complete overview
on these aspects in the case of elliptic equations.
We need to mention that these difficulties are 
exacerbated when anisotropic elements are used.
Several numerical mesh generation software such as \cite{FreeFem, Inria, A}, nevertheless
successfully use anisotropic adaptativity for the numerical simulation of PDE's,
as displayed for instance on Figure \ref{figInria}. From a numerical point of view, the improvement brought by these methods seems obvious compared to non-adaptive or adaptive isotropic methods. However, many aspects of the theoretical analysis of anisotropic methods remain open questions.\\

This thesis studies the problem of anisotropic mesh adaptation for the approximation of a {\it known} function, which may be regarded as a preliminary step for the analysis of anisotropic mesh adaptation for the numerical simulation of PDE's, but may also serve in other applications such as terrain data, surface and image processing.\\

Given a triangulation $\cT$ of a bounded polygonal domain $\Omega\subset\R^2$, and a fixed integer $k\geq 1$, we denote by $V_k(\cT)$ the space of finite elements of degree $k$ on $\cT$. The space $V_k(\cT)$ consists of all
functions which coincide on each triangle $T\in \cT$ with a polynomial of total degree $k$
$$
V_k( \cT) := \{g\; ; \; g_{|T} \in \P_k,\; T \in \cT\}.
$$
The dimension of $V_k(\cT)$ is of the order $\cO(k^2 \#(\cT))$. 
Given a function $f : \Omega \to \R$ and a triangulation $\cT$ of $\Omega$, 
the best approximation error of $f$ in $V_k(\cT)$ is defined by
\be
\label{eqApproxG}
e_\cT(f)_X := \inf_{g\in V_k(\cT)} \|f-g\|_X.
\ee
The letter $X$ indicates the norm or semi-norm in which the approximation error $\|f-g\|_X$ is measured. In this thesis we restrict our attention to the $L^p$ norm and the $W^{1,p}$ semi-norm, where $1\leq p \leq \infty$. They are defined as follows: 
$$
\|h\|_{L^p(\Omega)} := \left(\int_\Omega |h|^p\right)^{\frac 1 p} \stext{ and } |h|_{W^{1,p}(\Omega)} := \left(\int_\Omega |\nabla h|^p\right)^{\frac 1 p},
$$
with the standard modification when $p=\infty$. Note that when using the $W^{1,p}$ semi-norm,
we need to impose global continuity of $g$ in the above definition of $V_k(\cT)$.

The best approximation $g\in V_k(\cT)$ of $f$ can be exactly computed in the case of the $L^2$ norm or the $W^{1,2}$ (or $H^1$) semi-norm: $g$ is the orthogonal projection of $f$ onto $V_k(\cT)$ with respect to the scalar product associated to the norm of interest. In the case $p\neq 2$ of non hilbertian norms the best approximation of $f$ is generally hard to compute, but ``satisfactory'' approximations can be obtained by different methods. If the function is smooth (at least continuous), one may
use the standard Lagrange interpolant, while for non-smooth functions a 
quasi-interpolant operator is preferred, see Chapter \ref{chapApproxMet}.
More generally, if $P_{\cT}$ is any continuous projection operator from the space $X$ to $V_k(\cT)$, 
it is easily seen that
for any $f\in X$, one has
\be
\label{stableproj}
 \|f-P_\cT f\|_X \leq Ce_\cT(f)_X,
 \ee
where $C:=1+\|P_\cT\|_{X\to X}$. The problem of approximating a function $f$ on a {\it given} triangulation $\cT$, using finite elements of degree $k$, is thus solved in good part.\\
 
In the context of adaptive approximation, the triangulation $\cT$ of the domain $\Omega$ is not fixed, but can be freely chosen depending on the function $f$ (in contrast we always assume in this thesis 
that the integer $k$ is fixed although arbitrary). This naturally raises the objective of characterizing and constructing an optimal mesh for a given function $f$. Given a norm $X$ of interest
and a function $f$ to be approximated, we formulate the problem of {\it optimal mesh adaptation}, as 
minimizing the approximation error over all triangulations of {\it prescribed cardinality}. We therefore
define the adaptive best approximation error by
\be
\label{eqApproxT}
e_N(f)_X := \inf_{\#(\cT) \leq N} e_\cT(f)_X=\inf_{\#(\cT) \leq N}\inf_{g\in V_k(\cT)} \|f-g\|_X.
\ee
In contrast to the procedure of best finite element approximation on a fixed mesh, adaptive and anisotropic approximation is not yet well understood. In particular (i) {\it how does the optimal
mesh depend on the function} $f$ and (ii) {\it how does the optimal error $e_N(f)_X$ decay
as $N$ grows} ? These problems are well understood is the isotropic setting, 
for which the optimization is restricted to triangulations for which all triangles satisfy a uniform shape constraint
$$
\diam(T)^2\leq C |T|
$$
where $\diam(T)$ and $|T|$ stand for the diameter and area of $T$ respectively, and $C>0$ is a fixed constant. In the general setting of potentially anisotropic triangulations, 
they are open problems.

Heuristically, the simplicity of \iref{eqApproxG} comparatively to \iref{eqApproxT} is due to the fact that the optimization is posed on the linear space $V_k(\cT)$,
and that a near best solution may therefore be obtained by applying to $f$ a 
stable projection operator as expressed by \iref{stableproj}.
In contrast the optimization problem \iref{eqApproxT} is posed on the {\it union} of
spaces $V_k(\cT)$ for all triangulations $\cT$ satisfying $\#(\cT)\leq N$, which is certainly not a linear space. This problem is therefore an instance of {\it nonlinear approximation}. Other instances
include best $N$-terms approximations in a dictionary of functions, or best approximation by rational function. We refer to \cite{De} for a survey on nonlinear approximation.\\

The purpose of this thesis is to better understand optimal
mesh adaptation posed on the {\it full class of potentially anisotropic triangulations}. 
The four parts of the thesis are respectively devoted to the four questions below: 
\begin{enumerate}[I.]
\item 
How does the approximation error $e_N(f)_X$ behaves in the asymptotic regime
when the number of triangles $N$ grows to $+\infty$, when $f$ is
a smooth function ? In that context, we establish a mathematical characterization
of the optimal mesh, as well as sharp estimates of $e_N(f)_X$ by above and below
in terms of $N$ and quantities that depend {\it nonlinearly} on the derivatives of $f$.
\item 
Which classes of functions govern the rate of decay of $e_N(f)_X$ as $N$ grows,
and are in that sense naturally tied to the problem of optimal mesh adaptation? 
In particular, we have in mind the model of the so-called {\it cartoon functions},
which by definition are smooth except along a collection of smooth curves of discontinuity.
This is a popular image model in the image processing community (see for instance
Figure \ref{figEllTri0} for an instance of a cartoon image). We shall see that such a model 
naturally fits in a richer function class corresponding to a given rate of decay of $e_N(f)_X$.
\item 
Could the optimization problem \iref{eqApproxT} posed on triangulations 
$\cT$ of a given cardinality $N$, be replaced by an equivalent more tractable problem ?
Triangulations are indeed discrete combinatorial objects, described in terms of points and edges, 
which is not handy when solving optimization problems of the form \iref{eqApproxT}. 
We study the correspondence between certain classes of triangulations and of 
{\it riemanninan metrics} which in contrast are continuous objects. This allows us to reformulate 
and to solve the original optimization problem as a more tractable problem posed on the
set of riemannian metrics. 
\item 
Is it possible to produce a near-optimal sequence of triangulations $(\cT_N)_{N \geq _0}$ with
$\#(\cT_N)=N$, using a {\it hierarchical refinement procedure}? The property
of hierarchy guarantees the inclusion of the associated finite element spaces $V_k(\cT_N) \subset V_k(\cT_{N+1})$. It is required in applications such as progressive encoding or online
data processing. We provide with a simple and explicit algorithm 
which gives a positive answer to this question under some conditions.
\end{enumerate}

\DontWriteThisInToc
\section{State of the art}

Before detailing the content of the thesis, we
give in this section a short overview of the state of the art on both problems
of the study of $e_N(f)_X$ as $N$ grows,
and the construction of a near-optimal triangulation. For that purpose we need to introduce some notations.
We assume in the following that $f\in C^0(\overline \Omega)$, where $\Omega \subset \R^2$ is a bounded polygonal domain and $\overline \Omega$ denotes the closure of $\Omega$. For each triangulation $\cT$ of $\Omega$ we denote by $\interp^k_\cT$ the standard interpolation operator on Lagrange finite elements of degree $k$ on $\cT$: on each triangle $T\in\cT$
the interpolation $\interp^k_\cT f$ is the unique element of $\P_k$ which agrees with
$f$ on the points of barycentric coordinates in $\{0,\frac 1 k,\cdots,\frac {k-1} k,1\}$.
In the case $k=1$, these points are simply the three vertices of $T$.
If the triangulation $\cT$ is conforming (each edge of a triangle is either on the boundary
of $\Omega$ or coincides with the entire edge of another triangle), then $\interp_\cT^k f$ is continuous. In order to be consistent with the rest of this thesis we define $m := k+1\geq 2$, and we thus have 
$$
e_\cT(f)_{L^p} \leq \|f-\interp_\cT^{m-1} f\|_{L^p(\Omega)} \stext{ and } e_\cT(f)_{W^{1,p}} \leq  \|\nabla f-\nabla \interp_\cT^{m-1} f\|_{L^p(\Omega)}.
$$
Any estimate on the interpolation error thus automatically yields an upper estimate on the best approximation error $e_\cT(f)_X$. Furthermore if $f$ is sufficiently smooth and if $\cT$ is sufficiently fine, then these quantities are generally comparable.\\

One of the founding results of adaptive anisotropic finite element approximation deals with
the case of piecewise linear elements ($m=2$) with the error measured in the $L^p$ norm, see \cite{CSX} or Chapter \ref{chapOptAniso}. This result may be stated as follows: for any bounded polygonal domain $\Omega\subset \R^2$, for any $1 \leq p < \infty$, and for any function $f\in C^2(\overline \Omega)$ there exists a sequence $(\cT_N)_{N \geq N_0}$ of triangulations of $\Omega$, satisfying $\#(\cT_N) \leq N$, and such that 
\be
\label{eqCSX}
\limsup_{N\to \infty} N \|f-\interp^1_{\cT_N} f\|_{L^p(\Omega)} \leq C \left\| \sqrt{|\det (d^2f)|}\right\|_{L^\tau(\Omega)},
\ee
where the exponent $\tau \in (0, \infty)$ is defined by 
$$
\frac 1 \tau := 1 + \frac 1 p\, ,
$$ 
and $C$ is a universal constant ($C$ is independent of $p$, $\Omega$ and $f$).
We recall that the upper limit of a sequence $(u_N)_{N \geq N_0}$ is defined by 
\be
\label{defLimSup}
\limsup_{N \to \infty} u_N := \lim_{N \to \infty} \sup_{n \geq N} u_n,
\ee
and is in general strictly smaller than the supremum $\sup_{N\geq N_0} u_N$. It is still an open question to find an appropriate upper bound for $\sup_{N\geq N_0} \|f-\interp^1_{\cT_N} f\|_{L^p(\Omega)}$ when optimally adapted anisotropic partitions are used.
The result \iref{eqCSX} can be extended to the exponent $p=\infty$, and to simplicial partitions
of domains $\Omega$ of higher dimension, but the meshes $\cT_N$ may then not be conforming.

The result \iref{eqCSX} reveals that the accuracy of the approximation to $f$ is governed by the quantity $\sqrt{|\det(d^2 f)|}$, which depends non-linearly on the hessian matrix $d^2 f$. This non-linear dependency is heavily tied to the fact that we authorize triangles of potentially highly anisotropic shape. The estimate \iref{eqCSX} is obtained by producing sufficiently fine triangulations which combine the two following heuristical properties:
\begin{enumerate}[a)]
\item {\it Error equidistribution:} the contribution 
$\|f-\interp_T^1 f\|_{L^p(T)}$ of each triangle $T \in \cT_N$ to the global interpolation error $\|f-\interp_{\cT_N}^1 f\|_{L^p(\Omega)}$ is approximately the same. This condition can be translated into a local constraint on the area of the triangles, which is dictated by the local behavior of $f$ and in particular by $\det(d^2f(z))$.
\item {\it Optimal shape of the triangles:} the aspect ratio and the orientation of a triangle $T\in \cT_N$ 
is dictated by the ratio of eigenvalues and by the eigenvectors 
of the hessian matrix $d^2f(z)$ for $z\in T$.
\end{enumerate}

The simplest method for producing a sequence $(\cT_N)_{N \geq N_0}$ 
of triangulations satisfying \iref{eqCSX} is to use a ``local patching strategy'' 
that may be intuitively described as follows. In a first step the domain $\Omega$ is split into regions $\Omega_i$, $1\leq i \leq k$, sufficiently small so that the hessian matrix $d^2f(z)$ varies little
on each $\Omega_i$ around an average value $M_i$.
In other words $f$ is well approximated by a quadratic polynomial on each $\Omega_i$. 
Each region $\Omega_i$ is then tiled using a uniform triangulation $\cT_N^i$ which cells have an area,
aspect ratio and orientation based on $M_i$, and the triangulations $\cT_N^i$ of $\Omega_i$ are glued together into a triangulation $\cT_N$ of $\Omega$, at the price of a few additional triangles at the interfaces between 
the $\Omega_i$. 

\begin{figure}
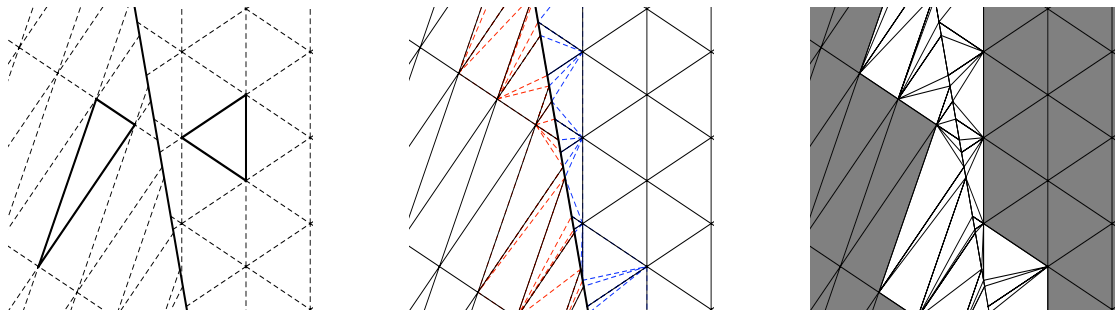

	\centering
		\includegraphics[width=4cm,height=4cm]{\pathPic/PaperOptAniso/NonConformTiling.pdf}
		\hspace{1cm}
		\includegraphics[width=4cm,height=4cm]{\pathPic/PaperOptAniso/ConformTiling.pdf}
		\hspace{1cm}
		\includegraphics[width=4cm,height=4cm]{\pathPic/PaperOptAniso/TriangleClasses.pdf}
	\caption{\label{figLocPatch} Illustration of the local patching strategy : the regions $\Omega_i$ are covered by periodic tilings, which are then glued together.}
\end{figure}

\begin{figure}
\centering
\includegraphics[height = 4cm,width = 8cm]{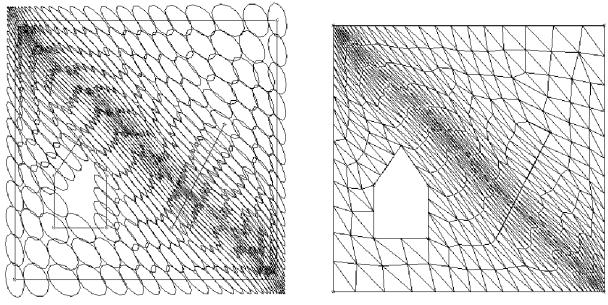}
\caption{A metric and an adapted triangulation. (credit : J. Schoen \cite{Schoen})\label{figEllTri}}
\end{figure}

This construction, which is illustrated on Figure \ref{figLocPatch} is sufficient for the purpose of
proving the asymptotical result \iref{eqCSX} but not for practical applications since it only becomes efficient for a large number of triangles. The following approach, based on riemannian metrics, is often preferred in applications.
We assume for simplicity that the hessian matrix $M(z) := d^2f(z)$ is positive definite for each $z\in \Omega$, and we define 
\be
\label{defHHessian}
H(z) := \lambda^2 (\det M(z))^{-\frac 1 {2p+2}} M(z)
\ee
where $\lambda>0$ is a constant which role is to control the resolution of the
triangulation. Such a map $H$ is called a riemannian metric, and continuously associates to each point $z\in \Omega$ a symmetric positive definite matrix $H(z)$. Note that for each $z\in \Omega$ the matrix $H(z)$ defines an ellipse
$$
\cE_z := \{u\in \R^2 \sep u^\trans H(z) u\leq 1\}.
$$ 
As illustrated on the left part of Figure \ref{figEllTri}, a metric encodes at each point $z\in \Omega$ an information of area, aspect ratio and orientation, under the form of a symmetric positive definite matrix $H(z)$, or equivalently of an ellipse $\cE_z$. 
Several mesh generation algorithms, such as \cite{FreeFem,Inria} are able to produce a triangulation adapted to a given metric $z\mapsto H(z)$ in the sense that for each $z\in \Omega$ the triangle $T\in \cT$ containing $z$ has a shape ``similar'' to the ellipse $\cE_z$  as illustrated on the right part of Figure \ref{figEllTri}.
In mathematical terms, this means that one has for each triangle $T$ and $z\in T$,
\be
\label{triellipse}
b_T+c_1 \cE_z\subset T \subset b_T+c_2 \cE_z,
\ee
where $0<c_1<c_2$ are fixed constants and $b_T$ is the barycenter of $T$,
or equivalently that $T$ is ``close'' to an equilateral 
triangle of unit area in the metric $H(z)$.
If $\cT$ is such a triangulation and if $\lambda$ is sufficiently large, a heuristic argument 
(which will be recalled in Chapter \ref{chapOptAniso}) shows that 
$$
\#(\cT) \|f-\interp_{\cT}^1 f\|_{L^p(\Omega)} \leq C \|\sqrt{|\det (d^2f)|}\|_{L^\tau(\Omega)},
$$
where the constant $C$ depends on $c_1$ and $c_2$.

Varying the parameter $\lambda$ we obtain different triangulations $\cT_\lambda$, of cardinality proportional to $\lambda^2$, which leads to the error estimate \iref{eqCSX}.
The production of the mesh $\cT$ from the metric $H$ is not straightforward. In
addition,  there exists no rigorous proof that
the similarity condition \iref{triellipse} is achieved for most mesh generation algorithms,
to the notable exception of \cite{Shew} and \cite{Bois}.\\

The approximation result \iref{eqCSX} is proved to be optimal if we restrict our attention to triangulations
which satisfy a technical condition defined as follows. 
We say that a sequence $(\cT_N)_{N \geq N_0}$ of triangulations of a polygonal domain $\Omega \subset \R^2$ is admissible if  $\#(\cT_N) \leq N$ and if 
\be
\label{defAdmiIntro}
\sup_{N \geq N_0} \left(\sqrt N \max_{T \in \cT_N} \diam(T)\right) <\infty.
\ee
For any admissible sequence $(\cT_N)_{N \geq N_0}$ of triangulations of $\Omega$, one can establish the lower bound 
\be
\label{eqCSXOpt}
\liminf_{N\to \infty} N \|f-\interp^1_{\cT_N} f\|_{L^p(\Omega)} \geq c \left\| \sqrt{|\det (d^2f)|}\right\|_{L^\tau(\Omega)},
\ee
where the constant $c>0$ is universal (see also Theorem \ref{thLowerLPIntro} in the outline
of the thesis below). 
Furthermore the admissibility condition is not too restrictive : for each $\ve>0$ there exists an admissible sequence of triangulations $(\cT_N)_{N \geq N_0}$ which satisfies the upper estimate \iref{eqCSX} up to the additive constant $\ve$ in the right hand side.\\

Similar results to \iref{eqCSX} and \iref{eqCSXOpt} can be developed for isotropic meshes, in which the triangles may vary in size but are constrained to be isotropic in the sense that the measure of degeneracy $\rho(T) := \diam(T)^2/|T|$ is uniformly bounded by a constant $\rho_0$, see for instance \cite{CMi2}. In that case the counterpart to  \iref{eqCSX} has the form
\be
\label{eqCSXiso}
\limsup_{N\to \infty} N \|f-\interp^1_{\cT_N} f\|_{L^p(\Omega)} \leq C \| d^2f\|_{L^\tau(\Omega)},
\ee
for the same value of $\tau$ and similarly for the counterpart to \iref{eqCSXOpt}, with
with constants $c$ and $C$ that depend on the bound $\rho_0$ on the measure of degeneracy, see for instance \cite{CMi2}. Therefore the 
nonlinear quantity
$\sqrt{|\det (d^2f)|}$ is replaced by
the linear $d^2f$ in the $L^\tau$ norm, and these results are now very similar
to those of best $N$-term wavelet approximation \cite{Co}.

In terms of the eigenvalues $\lambda_1(z), \lambda_2(z)$ of the symmetric matrix $d^2 f(z)$ we thus replace the geometric mean $\sqrt{|\lambda_1(z) \lambda_2(z)|}$ by $\max \{|\lambda_1(z)|, |\lambda_2(z)|\}$ which may be significantly larger when these eigenvalue have different order 
of magnitude. This is typically the case if the approximated function $f$ exhibits strongly anisotropic features, and therefore we may expect a significant improvement when using anisotropic meshes
for such functions.\\ 

The result \iref{eqCSX} gives a sharp account of the improvement that anisotropic triangulations can provide compared to isotropic triangulations, however in a rather restrictive setting which has
motivated our work:
\begin{enumerate}[I.]
\item The original result only applies to piecewise linear interpolation error measured in the $L^p$ norm, while higher order finite elements and Sobolev $W^{1,p}$ norms are equally relevant. In particular
the $W^{1,2}$ (or $H^1$) norm appears very naturally in the context of second order elliptic problems.
How should the elements be optimally adapted and 
how should the resulting error bound \iref{eqCSX} be modified in this more general context ?
\item The approximated function $f$ needs to be $C^2$, while the most interesting applications of adaptive approximation involve non smooth and even discontinuous functions, such as those appearing in Figures \ref{figEllTri0} and \ref{figInria}. Can we define in some sense a 
quantity such as $\|\sqrt{|\det (d^2f)|}\|_{L^\tau(\Omega)}$ when $f$ is a discontinuous function ?
\item In numerical applications, the riemannian metric $z\mapsto H(z)$ is used 
as an intermediate tool for mesh generation, but this approach lacks a precise equivalence result between these continuous objects and the different classes of anisotropic triangulations. Under which conditions on the metric $z\mapsto H(z)$ can we associate a triangulation $\cT$ which agrees with this
metric in the sense of \iref{triellipse} ?
\item The already mentionned anisotropic mesh generation algorithms 
are not hierarchical, in the sense that better accuracy is not achieved by
local refinement but by a global re-design of the mesh. Can we propose
a hierarchical mesh refinement algorithm that meets the optimal 
approximation bound  \iref{eqCSX} ?
\end{enumerate}

We also need to mention two fundamental issues which are discussed in this thesis,
but have remained open problems and
will be the subject of future research. First the error estimate \iref{eqCSX} only gives asymptotical information, as the number of triangles $N$ tends to $\infty$, while an error estimate valid for all values of $N$ is highly desirable. Second the extension of this result to functions defined on domains of dimension higher than two is hindered by a difficult problem of computational geometry : the production of {\it anisotropic and conforming meshes} in dimension $3$ or higher. This problem is addressed by some numerical methods such as \cite{Inria} but not solved from a theoretical point of view.

\DontWriteThisInToc
\section{Outline of the thesis}

This thesis is composed of four parts which attempt to solve four key issues, numbered I to IV, encountered in earlier results on adaptive and anisotropic finite element approximation.
These four parts are mostly self consistent and can be read independently.
The chapters that constitute each part should preferably be read sequentially,
to the exeption of chapter 1 which may be skipped by the reader who wishes to
go more directly to the heart of the matter. Chapters
1, 2, 3, 4, 7, 8, as well as the third part of Chapter 9, come from the papers
\cite{BLM}, \cite{Mi}, \cite{Mi2}, \cite{CM2}, \cite{CDHM}, \cite{CM} 
and \cite{CMi2}, respectively. Notations have been unified for the purpose of the thesis.
Chapters 5 and 6 are the object of future papers.

\subsection*{Part I. \ Finite elements of arbitrary degree, and Sobolev norms} 
\label{subsecOptAnisoIntro}

In this part, we discuss generalizations of the asymptotic estimate \iref{eqCSX} to finite elements of arbitrary degree, and to the Sobolev $W^{1,p}$ semi-norms as measure of error.
Through this analysis, we are led to introduce and study some key 
concepts for optimal mesh adaptation.\\

\begin{figure}
\centerline{
\includegraphics[width=4cm,height=4cm]{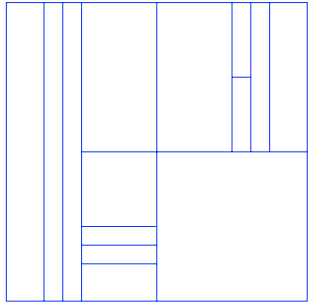}
}
\caption{Anisotropic partition of a domain into rectangles.\label{figRectPart}}
\end{figure}

As a starter we consider in Chapter \ref{chapBlocks} partitions of a rectangular domain into rectangles aligned with the coordinate axes,  as illustrated on Figure \ref{figRectPart}, instead of triangles of arbitrary direction. Such partitions are relevant when the coordinate axes play a privileged role, in such way that the anisotropic features of the approximated function $f$ are likely to be aligned with the coordinate axes.  In this setting, we obtain an optimal asymptotic error estimate.
Our results apply to piecewise polynomials of arbitrary degree, when the approximation error is measured in the $L^p$ norm, and in any dimension $d>1$. Here, we do not consider $W^{1,p}$ norms, since the approximants are generally discontinuous.

The biggest advantage of this setting is that the mathematical technicalities required for the construction of an anisotropic partition of a domain, as well as the error analysis, are simplified by the presence of preferred directions. We thus take advantage of this simple setting to introduce and study a 
key concept named the {\it shape function}, which governs the local approximation
error after optimal adaptation of the elements to the local properties of the approximated function.
The shape function is also defined and used for anisotropic triangulations
which are the object of Chapter 2 and 3. We give below its precise definition
in this setting. \\

Chapter \ref{chapOptAniso} is devoted triangular finite elements of arbitrary degree $m-1$
for $m\geq 2$, with $L^p$ as error norm.
In order to present our results we need to introduce some notations.
We denote by $\H_m$ the space of homogeneous polynomials of degree $m$:
$$
\H_m := \Span\{x^k y^l \sep k+l = m\}.
$$
A key ingredient of our approach is the shape function $K_{m,p} : \H_m \to \R_+$, where $1\leq p\leq \infty$ is an exponent. 
This function is defined by an optimization of the $L^p$ local interpolation error among the collection of triangles of unit area and of all possible shapes: for all $\pi\in \H_m$,
\be
\label{defKMPIntro}
K_{m,p}(\pi) := \inf_{|T| = 1} \|\pi-\interp_T^{m-1} \pi\|_{L^p(T)},
\ee
where $\interp_T^{m-1}$ is the local interpolation operator on $T$.
Our main result is the following asymptotic error estimate.
\begin{theoremIntro}
\label{thOptiLPIntro}
Let $\Omega\subset \R^2$ be a bounded polygonal domain, let $f\in C^m(\overline \Omega)$ and let  $1\leq p < \infty$.
There exists a sequence $(\cT_N)_{N \geq N_0}$, $\#(\cT_N) \leq N$, of conforming triangulations of $\Omega$ such that 
\be
\label{eqLPIntro}
\limsup_{N \to \infty} N^{\frac m 2}\|f-\interp_{\cT_N}^{m-1} f\|_{L^p(\Omega)} \leq \left\|K_{m,p}\left(\frac{d^m f}{m!}\right) \right\|_{L^\tau(\Omega)}
\ee
where the exponent $\tau \in (0, \infty)$ is defined by $\frac 1 \tau := \frac m 2 + \frac 1 p$.
\end{theoremIntro}

In the error estimate \iref{eqLPIntro}, we identify the collection $d^m f(z)$ of derivatives of order $m$ of $f$ at the point $z\in \Omega$ to the corresponding homogeneous polynomial in the Taylor development of $f$ at $z$. In mathematical form 
$$
\frac {d^m f(z)}{m!} \sim \sum_{k+l = m} \frac{\partial^m f}{\partial^k x\, \partial^l y}(z) \, \frac {x^k}{k!} \frac {y^l}{l!}.
$$

Theorem \ref{thOptiLPIntro} extends the known result \iref{eqCSX} to finite elements of arbitrary degree $m-1$. The quality of the adaptive and anisotropic approximation of $f$ is similarly determined by a non-linear quantity of its derivatives: the shape function $K_{m,p}(d^m f(z))$ is the ``generalization'' of the determinant $\sqrt{|\det (d^2f(z))|}$ appearing in \iref{eqCSX} to derivatives of arbitrary order.


This remark raises the need for an in depth study of the shape function $K_{m,p}$. For $m=2$ we have established that $K_{2,p}(\pi)$ is proportional to the square root of the determinant of the quadratic polynomial $\pi = a x^2 + 2 b xy + cy^2\in \H^2$:
$$
K_{2,p}(\pi) = c_{2,p} \sqrt{|\det \pi|} = c_{2,p}\sqrt {|ac-b^2|},
$$
where the constant $c_{2,p}>0$ only depends on the sign of $\det \pi$. We thus recover the earlier result \iref{eqCSX}.
For $m=3$, we establish that the shape function $K_{3,p}$ is the fourth root of the discriminant of the homogenous cubic polynomial $\pi = a x^3+3bx^2y+3cxy^2+dy^3\in \H_3$:
$$
K_{3,p}(\pi) = c_{3,p} \sqrt[4]{|\disc \pi|} = c_{3,p}\sqrt[4]{|4(ac-b^2)(bd-c^2) - (ad-bc)^2|},
$$
where the constant $c_{3,p}>0$ only depends on the sign of $\disc \pi$. For larger values $m\geq 4$ we have not obtained an explicit expression of the shape function $K_{m,p}$, but an explicit quantity uniformly equivalent to this function. This equivalent has the following form : there exists a polynomial $Q_m(a_0, \cdots, a_m)$ of the $m+1$ variables $a_0, \cdots, a_m$ and a constant $C_m\geq 1$ such that for all $\pi = a_0 x^m+ a_1 x^{m-1}y+ \cdots a_m y^m\in \H_m$ one has with $r_m := \deg Q_m$,
\be
\label{equivKQ}
C_m^{-1} \sqrt[r_m]{Q_m(a_0, \cdots , a_m)} \leq K_{m,p}(\pi) \leq  C_m\sqrt[r_m]{Q_m(a_0, \cdots , a_m)}.
\ee
The polynomial $Q_m$ is obtained using the theory of invariant polynomials developed by Hilbert in \cite{Hilbert}. We also characterise the zeros of the shape function, and thus 
the possible cases of ``super-convergence'' : $K_{m,p}(\pi)=0$ if and only if $\pi$ has a linear factor $ax+by$ of multiplicity $s>m/2$, in other words the homogeneous polynomial $\pi$ is sufficiently degenerated.\\

The proof of Theorem \ref{thOptiLPIntro} is based on a ``local patching strategy'', a two scale mesh generation procedure which proceeds as follows: we consider a first ``macro-triangulation'' $\cR$ of the polygonal domain $\Omega$ which is sufficiently fine in such way that the derivatives of order $m$ of $f$ vary little on each triangle $R\in \cR$ around an average value denoted by $\pi_R\in \H_m$. To each polynomial $\pi_R$, $R\in \cR$, we then associate a triangle $T_R$ which is a minimizer, or a near minimizer, of the optimization problem defining $K_{m,p}(\pi_R)$. We then tile each $R\in \cR$ using the triangle $T_R$ properly scaled and its symmetric with respect to the origin. Eventually, as illustrated on Figure \ref{figLocPatch} the triangulation $\cT_N$ is obtained by gluing together the periodic tilings on each $R\in \cR$, with a few additional triangles in order obtain a globally conforming mesh.\\


The next theorem establishes that the asymptotic error estimate \iref{eqLPIntro} is sharp, at least if we restrict our attention to admissible sequences of triangulations, a condition defined in \iref{defAdmiIntro}. The approximation result \iref{eqLPEpsIntro} establishes furthermore that the admissibility condition is not too restrictive.

\begin{theoremIntro}
\label{thLowerLPIntro}
Let  $\Omega\subset \R^2$ be a bounded polygonal domain, let $f\in C^m(\overline \Omega)$ and let $1\leq p\leq \infty$.
Let $(\cT_N)_{N \geq N_0}$, $\#(\cT_N) \leq N$, be an admissible sequence of triangulations of $\Omega$. Then
\be
\label{eqLPLowerIntro}
\liminf_{N \to \infty} N^{\frac m 2} \|f-\interp_{\cT_N}^{m-1} f\|_{L^p(\Omega)} \geq \left\|K_{m,p}\left(\frac{d^m f}{m!}\right) \right\|_{L^\tau(\Omega)},
\ee
where  $\frac 1 \tau := \frac m 2 + \frac 1 p$. Furthermore for any $\ve>0$ there exists an admissible sequence $(\cT_N^\ve)_{N \geq N_0}$ of triangulations of $\Omega$, satisfying $\#(\cT_N^\ve) \leq N$, and such that: 
\be
\label{eqLPEpsIntro}
\limsup_{N \to \infty} N^{\frac m 2} \|f-\interp_{\cT^\ve_N}^{m-1} f\|_{L^p(\Omega)} \leq \left\|K_{m,p}\left(\frac{d^m f}{m!}\right) \right\|_{L^\tau(\Omega)}+\ve.
\ee
\end{theoremIntro}

Chapter \ref{chapW1P} is devoted to the counterparts of Theorems \ref{thOptiLPIntro} and \ref{thLowerLPIntro} when the interpolation error is measured in the Sobolev $W^{1,p}$ semi-norms, $1\leq p < \infty$. These estimates involve the counterpart $L_{m,p} : \H_m \to \R_+$ of the shape function $K_{m,p}$ which is defined as follows: for all $\pi \in \H_m$ 
$$
L_{m,p}(\pi) := \inf_{|T|=1} \|\nabla (\pi-\interp_T^{m-1} \pi)\|_{L^p(T)}.
$$
We also provide some explicit equivalents of $L_{m,p}(\pi)$, defined similarly to \iref{equivKQ} by an algebraic expression in the coefficients of $\pi\in \H_m$.
Our results for the $W^{1,p}$ semi-norm are thus extremely similar to the results obtained for the $L^p$ norm.

\begin{figure}
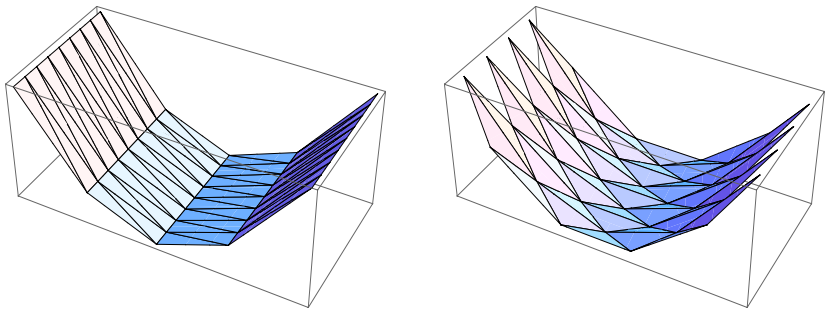

\centerline{
\includegraphics[width=5cm,height=4cm]{\pathPic/Triangles/ParaboleAigue.pdf}
\hspace{0.5cm}
\includegraphics[width=5cm,height=4cm]{\pathPic/Triangles/ParaboleObtuse.pdf}
}
\caption{Interpolation of a bi-dimensional parabolic function on a mesh built of acute (left) or ``flat'', strongly obtuse, triangles (right).\label{figParaAO}}
\end{figure}

The adaptation of the proofs however is not straightforward, due to the following new phenomenon: the presence of {\it strongly obtuse} triangles in a mesh (with one of their angles
close to $\cPi$) may cause the gradient of the interpolated function to oscillate, as illustrated in Figure \ref{figParaAO}. This phenomenon deteriorates the interpolation error in the $W^{1,p}$ semi-norm, but not in the $L^p$ norm. These ``flat'' triangles therefore need to be carefully avoided. In summary, long and thin triangles may be desirable but
they should not be strongly obtuse.

Before turning to description of the rest of this thesis we recall to the reader that the three chapters composing Part \ref{partAsymptApprox} are self consistent and can be read independently.

\subsection*{Part II. \ Anisotropic smoothness classes and image models}

In this part that consists of the sole Chapter 4, we discuss the extension of
our approximation results to non-smooth functions.\\

There exists various ways of measuring the smoothness 
of functions on a domain $\Omega\subset \R^2$, generally through the definition of an appropriate {\it smoothness space} with an associated norm. Classical instances are Sobolev and Besov spaces. Such spaces  are of common use when describing the regularity of solutions to partial differential equations. From a numerical perspective they are also useful in order to sharply characterize at which rate a function $f$ may be approximated by simpler functions such as Fourier series, finite elements (on isotropic triangulations), splines or wavelets.


The result of adaptive anisotropic approximation \iref{eqCSX} and its generalisation Theorem \ref{thOptiLPIntro} establish that the quality of the approximation of a function $f$ by finite elements on anisotropic triangulations is governed by a non-linear quantity of the derivatives of $f$, at least from an asymptotical point of view.
In the case of finite elements of degree $1$, and of the approximation in $L^2$ norm, the relevant quantity is the following
$$
A(f) := \left\|\sqrt {|\det(d^2f)|}\right\|_{L^{2/3}(\Omega)}.
$$
The functional $A$ strongly differs from the Sobolev, Holder or Besov norms, because of its nonlinear behavior: $A$ does not satisfy the triangle inequality, or any quasi-triangle inequality. In other words for any constant $C$ there exists $f,g\in C^2(\overline \Omega)$ such that 
\be
\label{eqNoTri}
A(f+g) > C(A(f)+A(g)).
\ee
The lack of a triangular inequality forbids to use the classical techniques of linear analysis to define a smoothness space attached to the functional $A$.
The extension of the approximation result  \iref{eqCSX} to functions which are not $C^2$
is therefore not straightforward.\\ 

Functions arising in concrete applications, for instance in image processing or in the solution of hyperbolic PDEs, often 
exhibit areas of smoothness separated by local discontinuities.
A mathematical model for this type of behavior is the collection of cartoon functions, which are smooth except for a jump discontinuity accross a collection of smooth curves.
A heuristical analysis presented in Chapter \ref{chapCartoon} suggests that for any cartoon function $f$ defined on a bounded polygonal domain $\Omega$, there exists a sequence $(\cT_N)_{N \geq N_0}$ of anisotropic triangulations of $\Omega$, satisfying $\#(\cT_N)\leq N$, and such that 
\be
\label{eqApproxCartoon}
\sup_{N\geq N_0} N \|f-\interp_{\cT_N}^1 f\|_{L^2(\Omega)} < \infty. 
\ee
As illustrated on Figure \ref{figEllTri0}, the triangulations $(\cT_N)_{N \geq N_0}$ combine highly anisotropic triangles aligned with the discontinuities of $f$, and large triangles in the regions where $f$ is smooth.
The approximation result \iref{eqApproxCartoon} raises the hope that a precise and quantitative asymptotic error estimate for the adaptive anisotropic approximation of cartoon functions extends the known result 
\be
\label{eqApproxA}
\limsup_{N \to \infty} N \|f-\interp_{\cT_N}^1 f\|_{L^2(\Omega)} \leq C A(f),
\ee
which holds if $f$ is $C^2$.

So far we have only fulfilled one part of this program: the extension of the functional $A$ to cartoon functions. More precisely, consider a radial and compactly supported function $\vp\in C^\infty(\R^2)$ of unit integral. 
For any $\delta>0$ we define the mollifier $\vp_\delta(z) := \frac 1 {\delta^2} \vp\left(\frac z \delta\right)$ and we denote by $f_\delta := f*\vp_\delta$ the convolution of $f$ with $\vp_\delta$. 
We establish in Chapter \ref{chapCartoon} that if $f$ is a cartoon function, then
$A(f_\delta)$ converges as $\delta \to 0$ to an explicit expression:
\begin{eqnarray*}
\lim_{\delta \to 0} A(f_\delta)^{2/3} &=& \left\|\sqrt{|\det(d^2 f)|}\right\|_{L^{2/3}(\Omega\sm\Gamma)}^{2/3} + C(\vp) \left\|\sqrt{|\kappa|} \gamma\right\|_{L^{2/3}(\Gamma)}^{2/3}\\ 
&=& \int_{\Omega\sm\Gamma} \sqrt[3]{|\det (d^2 f(z))|}dz+ C(\vp)\int_\Gamma |\kappa(s)|^{1/3} |\gamma(s)|^{2/3}ds.
\end{eqnarray*} 
Here, we have denoted by $\Gamma$ the collection of curves where the cartoon function $f$ is discontinuous, by $\gamma(s)$ the jump of $f$ at the point $s\in \Gamma$, and by $\kappa(s)$ the curvature of $\Gamma$ at $s$. The constant $C(\vp)$ is positive and explicit in terms of $\vp$.\\

The extension of the nonlinear functional $A$ to cartoon functions puts in light a link, explored in  
\S\ref{secCartoon4}, 
between anisotropic finite element approximation and several other mathematical fields. We think in particular of affine invariant image processing \cite{Ca} and of the definition in \cite{DPW} of smoothness spaces through the regularity of the level lines of a function.\\

\begin{figure}
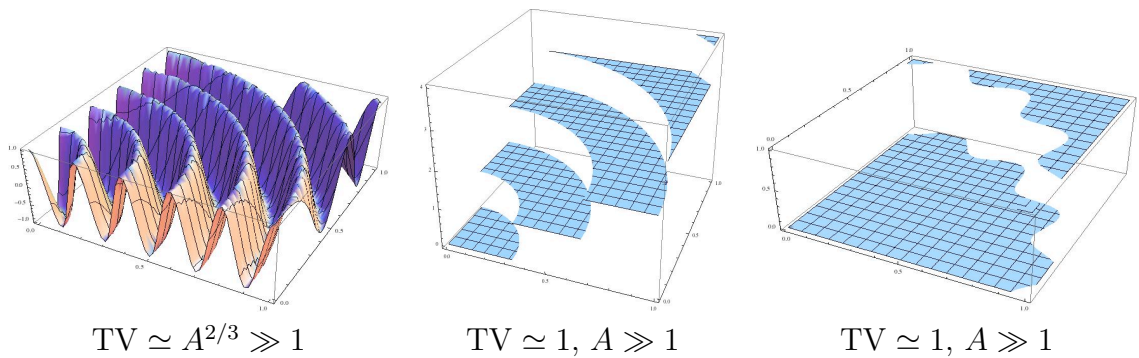

\centering
\begin{tabular}[c]{ccc}
\includegraphics[height=3.5cm,width=5cm]{\pathPic/Bayesien/SinaNoText.pdf}
&\includegraphics[height=4cm,width=4cm]{\pathPic/Bayesien/StairaNoText.pdf}
&\includegraphics[height=3.5cm,width=5cm]{\pathPic/Bayesien/FloweraNoText.pdf}\\
 $\TV \simeq A^{2/3} \gg 1$
& $\TV \simeq 1$, $A \gg 1$
& $\TV \simeq 1$, $A \gg 1$
\end{tabular}
\caption{Behavior of the functionals $\TV$ and $A$ on different types of images.\label{figTVA}}
\end{figure}

We also regard the quantity $A$ as a ``second order'' counterpart to the total variation semi-norm $\TV$, a measure of smoothness defined in terms of the first derivatives which is also finite for all cartoon functions. The total variation plays a central role in image processing and in the analysis of transport equations, two domains in which piecewise smooth functions naturally appear. For any cartoon function $f$, the total variation $\TV(f)$ is given by the following formula  
$$
\TV(f) = \int_{\Omega\sm\Gamma} |\nabla f(z)|dz+\int_\Gamma |\gamma(s)| ds.
$$
We compare the behavior of the two quantities $A(f)^{2/3}$ and $\TV(f)$ for different families of cartoon functions $f$. It appears that these quantities are equivalent when $f$ is a smoothly oscillating function $f(z) := \cos(\omega |z|)$ with large $\omega$, see figure \ref{figTVA} (i). In contrast the discontinuities are penalized in a different manner by these two functionals. For a staircase-like function, as illustrated on Figure \ref{figTVA} (ii), the total variation remains bounded while $A$ tends to infinity as the number of steps grows, due contribution $|\gamma(s)|^{2/3}$ of the jump term. Furthermore, due to the curvature term $|\kappa(s)|^\frac 1 3$, $A$ is much larger than $\TV$ for characteristic functions of sets with a complex or oscillating boundary as illustrated in Figure \ref{figTVA}. 

The functional $A$ can thus be thought of as a new quantitative image model : a monochromatic image, described by a brightness function $f: [0,1]^2 \to [0,1]$, is plausible if $A(f)$ is sufficiently small. Generalising the work \cite{LMo} we discuss an image recovery procedure using this model as a bayesian prior. At the present time, this
algorithm is not satisfactory due to very high convergence time, and for this reason we present some
numerical illustrations in a simplified one-dimensional setting. \\

Eventually we discuss the extension to cartoon functions of other non-linear quantities appearing in anisotropic finite element approximation such as the norm $\|K_{m,p}(d^m f)\|_{L^\tau}$ of the shape function for $m\geq 2$, or the counterpart to this quantitiy for functions $f$ of more than two variables.

\subsection*{Part III. \ Mesh generation and riemannian metrics} 

Triangulations are discrete objects of combinatorial nature : they can be described by the collection of their vertices and of the edges between them. This description is fruitful for the demonstration of algebraic results such as the Euler formula, or for computer processing.
In contrast, as explained earlier,
many approaches towards anisotropic mesh adaptation \cite{Shew,Bois,A} 
are based on a continuous equivalent object, namely 
a riemannian metric $z\mapsto H(z)$, in other words a continuous function $H$ from $\Omega$ to the set $S_2^+$ of symmetric positive definite matrices. Once this metric has been
properly designed, it is the task of the mesh generation algorithm to design
a triangulation that agrees with this metric according to \iref{triellipse}.
The purpose of Part \ref{partRiemann} is to formulate 
precise equivalence results between some classes of triangulations and of riemannian metrics. 
This equivalence translates some geometrical constraints satisfied by the triangulations into the form regularity properties of the equivalent riemannian metrics.\\

In order to state our results we need to introduce some notations. We associate to each triangle $T$ its barycenter $z_T\in \R^2$ and the symmetric positive definite matrix $\cH_T\in S_2^+$ such that the ellipse $\cE_T$ defined by 
$$
\cE_T := \{z\in \R^2 \sep (z-z_T)^\trans \cH_T (z-z_T)\},
$$
is the ellipse of minimal area containing $T$. The point $z_T$ thus encodes the position of $T$, while the matrix $\cH_T\in S_2^+$ encodes the area, the aspect ratio and the orientation of the triangle $T$. 

\begin{figure}
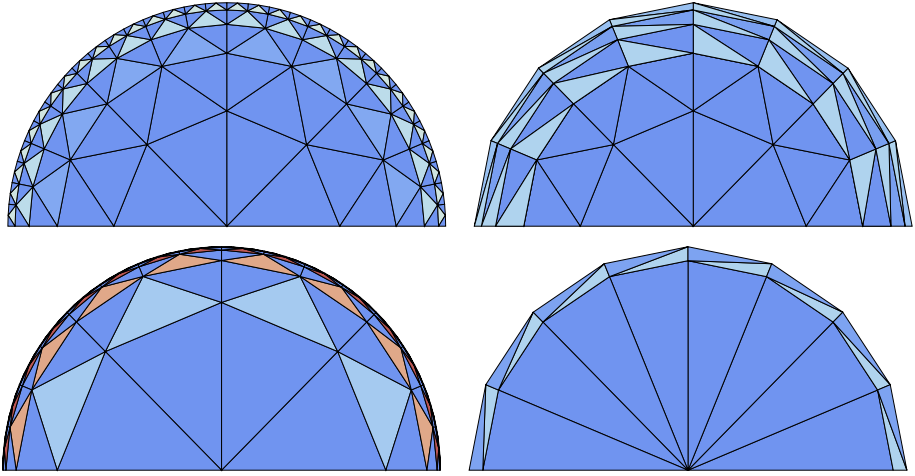

	\centering
		\includegraphics[width=6cm,height=3.2cm]{\pathPic/CoucheLimite/DiskIso.pdf}
		\includegraphics[width=6cm,height=3.2cm]{\pathPic/CoucheLimite/DiskAigu.pdf}
		\includegraphics[width=6cm,height=3.2cm]{\pathPic/CoucheLimite/DiskGen.pdf}
		\includegraphics[width=6cm,height=3.2cm]{\pathPic/CoucheLimite/DiskNoReg.pdf}
	\caption{\label{figDisk} Generation of a layer of simplices of width $\delta$ along the unit circle using $\cO(\delta^{-1})$ isotropic triangles (top left), $\cO(\delta^{-\frac 1 2} |\ln \delta|)$ ``quasi-acute'' triangles (top right), $\cO(\delta^{-\frac 1 2})$ triangles satisfying the grading condition
	or not (bottom left and bottom right).}
\end{figure}

We denote by $\bT$ the collection of conforming triangulations of the infinite domain $\R^2$.
The choice to consider infinite triangulations is guided by simplicity, 
and further work will be devoted to triangulations of bounded polygonal domains.
We denote by $\bH := C^0(\R^2, S_2^+)$ the collection of riemannian metrics on $\R^2$. A metric $H\in \bH$ continuously associates to each point $z\in \R^2$ a symmetric positive matrix $H(z)\in S_2^+$.
For $C\geq 1$, we say that a triangulation $\cT\in \bT$ is $C$-equivalent to a metric $H \in \bH$ if for all $T\in \cT$ and for all $z\in T$ we have 
$$
C^{-2} H(z) \leq \cH_T \leq C^2 H(z),
$$
where these inequalities are meant in the sense of symmetric positive matrices.
We say that a class of triangulations $\bT_*\subset \bT$ is equivalent to a class of metrics $\bH_*\subset \bH$ if there exists a uniform constant $C\geq 1$ such that 
\begin{itemize}
\item For any triangulation $\cT\in \bT_*$ there exists a metric $H\in \bH_*$ such that $\cT$ and $H$ are $C$-equivalent.
\item For any metric $H\in \bH_*$ there exists a triangulation $\cT\in \bT_*$ such that $\cT$ and $H$ are $C$-equivalent.\\
\end{itemize}
We consider three particularly relevant classes of triangulation of $\R^2$
$$
\bT_{i,C} \subset \bT_{a,C} \subset \bT_{g,C}
$$
which have a minor dependency on a parameter $C\geq 1$.
The class $\bT_{g,C}$ of \emph{graded} triangulations, is defined by the following condition which imposes a minimum of consistency among the shapes of neighbouring triangles. A triangulation $\cT$ belongs to $\bT_{g,C}$ if for all $T, T'\in \cT$ one has
$$
T \cap T' \neq \emptyset \quad \Ra \quad C^{-2} \cH_T \leq \cH_{T'} \leq C^2 \cH_T.
$$
A triangulation $\cT$ belongs to the class $\bT_{i,C}$ of \emph{isotropic} triangulations if $\cT$ is graded, $\cT\in \bT_{g,C}$,
and if the elements of $\cT$ are sufficiently close to the equilateral triangle, as expressed by the following condition: for all $T\in \cT$
$$
\|\cH_T\| \|\cH_T^{-1}\|\leq C^2.
$$
The class $\bT_{a,C}$ of \emph{quasi-acute} triangulations is defined by a slightly more technical condition on the maximal angles of the triangles, see Chapter \ref{chapMeshMet}. 
As explained in the earlier description of Chapter 3, avoiding flat angles is needed
to guarantee the stability of the gradient when applying the interpolation operator.
One of our key results is the reformulation of this condition 
under the form of a regularity assumption on the equivalent riemannian metric. 
From a practical point of view, this condition is unfortunately not guaranteed
by existing anisotropic mesh generation software.
The constraints satisfied by these classes of triangulations of $\R^2$ are illustrated by some bounded triangulations on Figure \ref{figDisk}.\\

The results of Chapter \ref{chapMeshMet} establish that when $C$ is sufficiently large, these 
three classes of triangulations are equivalent to three classes of metrics respectively
$$
\bH_i \subset \bH_a \subset \bH_g,
$$
that are defined by precise smoothness conditions on the function $z\mapsto H(z)$.
In particular, the class $\bH_i$ consists of those metrics $H$ which are proportional to the identity
$
H(z) = \Id /s(z)^2 ,
$
and such that the proportionality factor $s$ satisfies one of the following surprisingly
\emph{equivalent} properties
\begin{enumerate}[$\bullet$]
\item Euclidean Lipschitz Property: $|s(z) - s(z')| \leq |z-z'|$ for all $z,z'\in \R^2.$
\item Riemannian Lipschitz Property: $|\ln s(z) - \ln s(z')| \leq d_H(z,z')$ for all $z,z'\in \R^2.$
\end{enumerate}
We recall that the riemannian distance $d_H(z,z')$ measures the length, distorted by the metric $H$, of the smallest path joining $z$ and $z'$. 
$$
d_H(z,z') := \inf_{\substack{\gamma(0) = z\\ \gamma(1) = z'}} \int_0^1 \|\gamma'(t) \|_{H(\gamma(t))} dt
$$
where $\|u\|_M := \sqrt {u^\trans M u}$ and where the infimum is taken among all paths $\gamma\in C^1([0,1], \R^d)$ joining $z$ and $z'$.

The two lipschitz properties presented above have natural extensions 
to general anisotropic riemannian metrics, {\it which are not anymore
equivalent}. Considering $H\in \bH$ and defining $S(z) := H(z)^{-\frac 1 2}$ for all $z\in \R^2$, these Lipschitz properties are namely
\be
\label{eqDPlusIntro}
d_+(H(z), H(z') )  \leq  |z-z'| \stext{ for all } z,z'\in \R^2,
\ee
where $d_+(H(z), H(z') ) := \|S(z) - S(z')\|$, and
\be
\label{eqDTimesIntro}
d_\times(H(z), H(z') )  \leq d_H(z,z') \stext{ for all } z,z'\in \R^2,
\ee
where $d_\times(H(z), H(z') ) := \ln \max \{\|S(z) S(z')^{-1}\| , \|S(z') S(z)^{-1}\| \}$.
Note that $d_+$ and $d_\times$ are distances on $S_2^+$.

The main result of Chapter \ref{chapMeshMet} characterizes
the classes of metrics $\bH_g$ and $\bH_a$ (associated to
the classes of graded triangulations $\bT_{g,C}$ and of quasi-acute triangulations $\bT_{a,C}$)
in terms of the above Lipschitz properties.
\begin{theoremIntro}
\label{thMeshMetIntro}
If the constant $C$ is sufficiently large, then the class $\bT_{g,C}$ of graded triangulations is equivalent to the class $\bH_g$ of metrics satisfying \iref{eqDTimesIntro}, and the class $\bT_{a,C}$ of quasi-acute triangulations is equivalent to the class $\bH_a$ of metrics satisfying simultaneously \iref{eqDTimesIntro} and \iref{eqDPlusIntro}.
\end{theoremIntro}

We give in Chapter \ref{chapApproxMet} some applications of these results in the contexts of approximation theory and of constrained mesh generation. Let us describe the latter example. 
For any closed set $E\subset \R^2$ and any mesh $\cT\in \bT$ we denote by $V_\cT(E)$ the neighborhood of $E$ in the triangulation $\cT$, which is defined as follows
$$
V_\cT(E) := \bigcup_{\substack{T\in \cT\\ T \cap E \neq \emptyset}} T.
$$
We say that a mesh $\cT$ separates two disjoint closed sets $X$ and $Y$ if $V_\cT(X) \cap V_\cT(Y) = \emptyset$. The counterpart for metrics of this property is the following : a metric $H\in \bH$ separates $X$ and $Y$ if $d_H(x,y)\geq 1$ for all $x\in X$ and $y\in Y$. We show that these two properties are rigorously equivalent, and we use this reformulation to compute the smallest number of triangles (up to periodicity) required for the separation of some (periodic) sets using a (periodic) isotropic, quasi-acute or graded triangulation. As illustrated on Figure \ref{figDisk}, 
imposing more constraints on the triangulation typically raises the number
of triangles needed to achieve the same task.

The rest of this chapter is devoted to the control of the finite element approximation error of a function $f$ on a triangulation $\cT$, in $L^p$ norm or $W^{1,p}$ semi-norm, by a quantity $e_H(f)_p$, $e_H^a(\nabla f)_p$ or $e_H^g(\nabla f)_p$ attached to a metric $H$ equivalent to $\cT$. This analysis puts in light, in the case of Sobolev norms, the role played by the angle or regularity conditions which define the elements of $\bT_a$ or $\bH_a$ respectively. Finally, we present a counterpart for metrics of the asymptotic approximation results developed for triangulations Chapters 2 and 3.

\subsection*{Part IV. \ Hierarchical refinement  algorithms}

\begin{figure}
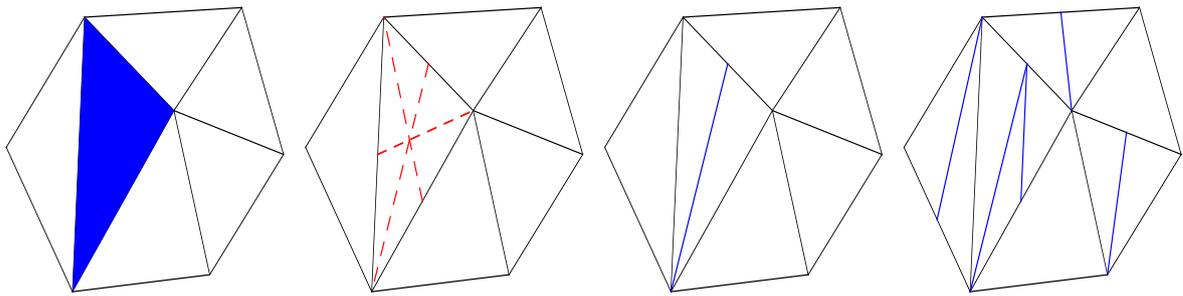

	\centering
		\includegraphics[width=3.8cm,height=4cm]{\pathPic/Subdivision/Algo1.pdf}
		\includegraphics[width=3.8cm,height=4cm]{\pathPic/Subdivision/Algo2.pdf}
		\includegraphics[width=3.8cm,height=4cm]{\pathPic/Subdivision/Algo2p.pdf}
		\includegraphics[width=3.8cm,height=4cm]{\pathPic/Subdivision/Algo3.pdf}
	\caption{\label{figRefinement}Refinement algorithm : selection of the triangle which maximizes the local approximation error (i, in dark), choice of a bisection among the three possibilities (ii,iii), iteration of these two steps (iv).}
\end{figure}

The last part of this thesis is devoted to the study an algorithm proposed by Cohen, Dyn and Hecht which produces hierarchical sequences $(\cT_N)_{N\geq N_0}$ of (non-conforming) anisotropic triangulations adapted to a given function $f$.
Given a triangulation $\cT$ of a domain $\Omega$ and a function $f\in L^p(\Omega)$, this algorithm creates in one step a triangulation $\cT'$ of $\Omega$, of cardinality $\#(\cT') = \#(\cT)+1$, proceeding as follows:
\begin{enumerate}
\item (Greedy triangle selection) A triangle $T\in \cT$ which has a maximal contribution to the error is selected 
$$
T := \underset {T'\in \cT} \argmax \, \|f-\cA_{T'} f\|_{L^p(T')},
$$
where $\cA_T : L^p(T) \to \P_{m-1}$ is a projection operator, for instance the $L^2(T)$ orthogonal projection onto $\P_{m-1}$.
\item (Decision of a bisection) An edge $e\in \{a,b,c\}$ of $T$ is selected by minimizing a 
a given decision function $e\mapsto d_T(e,f)$ among
the three edges. The triangle is refined by joining the mid-point of this edge
to the opposite vertex creating the sub-triangles $T_e^1$ and $T_e^2$. 
The new triangulation is thus
$$
\cT' := \cT-\{T\} + \{T_e^1, T_e^2\}.
$$
\end{enumerate}
Starting from a triangulation $\cT_0$ of cardinality $N_0$ of the domain $\Omega$, the algorithm produces step after step a sequence $(\cT_N)_{N\geq N_0}$, see Figure \ref{figRefinement}, of triangulations ``adapted'' to a given function $f\in L^p(\Omega)$. The properties of these triangulations strongly depend on the approximated function $f$ and on the choice of the decision function $d_T(e,f)$, which governs the creation of anisotropy. In contrast the projection operator $\cA_T$ plays a rather minor role. A typical choice for the decision function $e\mapsto d_T(e,f)$
is the local error after bisection of the edge $e$, namely the algorithm selects the 
bisection that most reduces the error. \\

We present this algorithm in more detail in Chapter \ref{chapCDHM}, and we establish its convergence 
in the sense that the (discontinuous) piecewise polynomial approximation $\cA_{\cT_N} f$ defined by
$$
\cA_{\cT_N} f(z)=\cA_T(z),\; \; z\in T,\;\; T\in \cT_N
$$
converges towards $f$ in $L^p$ as $N\to +\infty$ for any $f\in L^p$, 
under some assumptions on the decision function $e\mapsto d_T(e,f)$.
We also discuss the possibility of using the multiscale hierachical structure to
define multiresolution approximation, wavelets and CART algorithms.
We illustrate the anisotropic adaptation of the algorithm
on several types of functions and images presenting 
sharp transitions along curved edges.\\

\begin{figure}
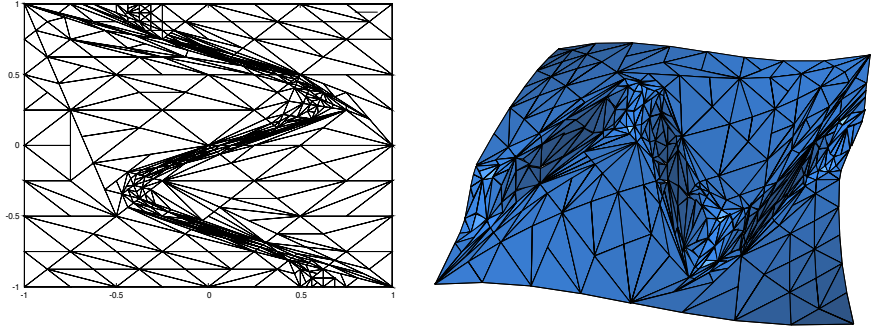

	\centering
		\includegraphics[width=6cm,height=5cm]{\pathPic/Triangles/f1-512.pdf}
		\includegraphics[width=6cm,height=4cm]{\pathPic/Triangles/f1-512_2.pdf}
	\caption{\label{figSinGreedy}
	Triangulation produced by the hierarchical refinement algorithm, and interpolation, for a function exhibiting a sharp transition close to a sinuosidal curve.}
\end{figure}

In Chapter \ref{chapBisecOpt} we make a deeper convergence analysis
of the algorithm in the case $m=2$ of piecewise linear elements.
Our main result states that when $f$ is a strictly convex $C^2$ function, then
for a particular choice of the decision function based on the
$L^1$ local interpolation error, 
the sequence $(\cT_N)_{N\geq N_0}$ of triangulations generated by this algorithm satisfies the optimal asymptotic convergence estimate
$$
\limsup_{N \to \infty} N \| f-\cA_{\cT_N} f\|_{L^p(\Omega)} \leq C \|\sqrt{|\det (d^2f)|}\|_{L^\tau(\Omega)},
$$
where $\frac 1 \tau := 1+ \frac 1 p$. The key observation leading to this
result is that when $f$ is a convex quadratic polynomial, the minimization of the
decision function selects the {\it longest edge} in the metric associated to 
the homogeneous quadratic part of this polynomial. This is used to prove that
the triangles generated by the algorithm tend in majority to adopt an optimal aspect ratio.
This good behaviour of the algorithm is also observed on more general non-convex
$C^2$ functions, as illustrated on Figure \ref{figSinGreedy}.
However, proving the above optimal convergence bound remains an open problem
in this general setting.\\

Finally we study in Chapter \ref{chapVariations} some variants of the hierarchical refinement algorithm presented above. We first consider a decision function based on the $L^2$ local 
projection error, which is particularly well suited to the numerical implementation since it can be evaluated in small computing time.

We then focus on the behaviour of the algorithm when applied to cartoon functions. 
The original algorithm does not satisfy the best possible convergence estimate
for such functions. We show that the optimal convergence rate may be
recovered when replacing the bisection from the mid-point of the
edge towards the opposite vertex by more general geometric
splitting procedure.

Eventually we consider another variant of the algorithm, which 
is based on rectangles aligned with the coordinate axes instead of triangles of arbitrary direction,
in the spirit of the rectangular partitions studied in Chapter 1.
This simplification leads to a result which guarantees the best possible convergence estimate for all $C^1$ functions, in the case of piecewise constant approximations.

%


\part{Optimal mesh adaptation for finite elements of arbitrary order} 
\label{partAsymptApprox}

\chapter[Sharp asymptotics of the $L^p$ approximation error on rectangles]{Sharp asymptotics of the $L^p$ approximation error on rectangular partitions} 
\minitoc
\label{chapBlocks}

\section {Introduction}

The purpose of this chapter is to study the adaptive anisotropic approximation of a function by interpolating splines defined over block partitions in $\RR^d$. We use the word block as a synonym for ``$d$-dimensional rectangle''.
Our analysis applies to an arbitrary projection operator in arbitrary dimension. We then apply the obtained estimates to several different interpolating schemes most commonly used in practice.

Our approach is to introduce the ``shape function'' which reflects the interaction of approximation procedure with polynomials. Throughout the chapter we shall study the asymptotic behavior of the approximation error and, whenever possible, the explicit form of the shape function which plays a major role in finding the constants in the formulae for exact asymptotics.

\subsection{The projection operator}

Let us first introduce the definitions that will be necessary to state the main problem and the results of this chapter.

We consider a fixed integer $d\geq 1$ and we denote by $x=(x_1, \cdots , x_d)$ the elements of $\R^d$. A block $R$ is a subset of $\R^d$ of the form
$$
R = \prod_{1\leq i\leq d} [a_i,b_i] 
$$
where $a_i< b_i$, for all $1\leq i\leq d$. 
For any block $R\subset \R^d$, by $L^p(R)$, $1\le p\le \infty$, we denote the space of measurable functions $f:R \to\RR$ for which the value
$$
  \|f\|_p = \|f\|_{L^p(R)} : = \left\{\begin{array}{ll}
                \left(\displaystyle\int\limits_{R} |f (x) | ^p dx\right)^{\frac 1p},   &{\rm if}\;\;\;    1\leq  p < \infty, \\ [10pt]
                {\rm ess sup} \{|f (x) | \sep x \in R \},                                    &{\rm if}\;\;\;       p =\infty.
\end{array}\right.
$$
is finite. 
We also consider the space $C^0(R)$ of continuous functions on $R$ equipped
with the uniform norm $\|\cdot\|_{L^\infty(R)}$.
We shall make a frequent use of the canonical block $\mI^d$, where $\mI$ is the interval
$$
\mI := \left[-\frac 1 2, \frac 1 2\right].
$$
Next we define the space $V := C^0(\mI^d)$ and the norm $\|\cdot\|_V := \|\cdot\|_{L^\infty(\mI^d)}$.
Throughout this chapter we consider a linear and bounded (hence, continuous) operator
$
\interp : V\to V. 
$
This implies that there exists a constant $C_{\interp}$ such that
\be
\label{contI}
\|\interp u\|_V \leq C_{\interp} \|u\|_V \text{ for all } u\in V.
\ee
We assume furthermore that $\interp$ is a projector, which means that it satisfies
\be
\label{projAxiom}
\interp\circ \interp = \interp.
\ee
Let $R$ be an arbitrary block. 
It is easy to show that there exists a unique $x_0\in \R^d$ and a unique diagonal matrix $D$ with positive diagonal coefficients such that the transformation
\be
\label{defPhi}
\phi(x) := x_0+ Dx \ \text{ satisfies } \ \phi(\mI^d) = R.
\ee
The volume of $R$, denoted by $|R|$, is equal to $\det(D)$. For any function $f\in C^0(R)$ we then define
\be
\label{interpR}
\interp_R f := \interp(f\circ\phi)\circ \phi^{-1}.
\ee
Note that
\be
\label{changeRect}
\|f-\interp_R f\|_{L^p(R)} = (\det D)^{\frac 1 p}\|f \circ \phi - \interp(f\circ\phi)\|_{L^p(\mI^d)}.
\ee
A block partition $\cR$ of a block $R_0$ is a finite collection of blocks such that their union covers $R_0$ and which pairwise intersections have zero Lebesgue measure.
If $\cR$ is a block partition of a block $R_0$ and if $f\in C^0(R_0)$, by
$
\interp_\cR f \in L^\infty(R_0)
$
we denote the (possibly discontinuous) function which coincides with $\interp_R f$ on the interior of each block $R\in \cR$.

{\bf Main Question.}
The purpose of this chapter is to understand the asymptotic behavior of the quantity
$$
\|f-\interp_{\cR_N} f \|_{L^p(R_0)}
$$
{\it {for each given function $f$}} on $R_0$ from some class of smoothness,
where $(\cR_N)_{N\geq 1}$ is a sequence of block partitions of $R_0$ that are optimally adapted to $f$. 

Note that the exact value of this error can be explicitly computed only in trivial cases. Therefore, the natural question is to  study the asymptotic behavior of the shape function, i.e. the behavior of the error as the number of elements of the partition ${\cR_N}$ tends to infinity.

Most of our results hold with only assumptions \iref{contI} of continuity of the operator $\interp$, the projection axiom \iref{projAxiom}, and the definition of $\interp_R$ given by \iref{interpR}. Our analysis therefore applies to various projection operators $\interp$, such as the $L^2$-orthogonal projection on a space of polynomials, or spline interpolating schemes described in \S\ref{exampleI}.

\subsection{History}

The main problem formulated above is interesting for functions
of arbitrary smoothness as well as for various classes of splines
(for instance, for splines of higher order, interpolating splines,
best approximating splines, etc.). In the univariate case general
questions of this type have been investigated by many authors. The
results are more or less complete and have numerous applications
(see, for example, ~\cite{LSh}).

Fewer results are known in the multivariate case. Most of them are
for the case of approximation by splines on triangulations (for
review of existing results see, for instance ~\cite{Gr3, KL, BBLS, CMi2} and Chapter \ref{chapOptAniso}).
However, in applications where preferred directions exist, box
partitions are sometimes more convenient and efficient.

The first result on the error of interpolation on
rectangular partitions by bivariate splines linear in each variable
(or bilinear) is due to D'Azevedo ~\cite{Daz1} who obtained
local (on a single rectangle) error estimates. In ~\cite{PhD} Babenko obtained the exact
asymptotics for the error (in $L^1$, $L^2$, and $L^{\infty}$ norms)
of interpolation of $C^2(\mI^d)$ functions by bilinear splines.

In ~\cite{JAT} Babenko generalized the result to interpolation and
quasiinterpolation of a function $f\in C^2(\mI^d)$ with arbitrary but fixed throughout the domain signature (number of positive and negative second-order partial derivatives). However, the
norm used to measure the error of approximation was uniform.

In this chapter we use a different, more abstract, approach which allows us to obtain the exact asymptotics of the error  in a more general framework which can be applied to many particular interpolation schemes by an appropriate choice of the interpolation operator. In general, the constant in the asymptotics is implicit. However, imposing additional assumptions on the interpolation operator allows us to compute the constant explicitly.

The chapter is organized as follows. Section \ref{subsecApprox}
contains the statements of main approximation results. The closer
study of the shape function, as well as its explicit formulas under
some restrictions, can be found in Section \ref{studyK}. The proofs
of the theorems about asymptotic behavior of the error are contained
in Section \ref{proofTh}.

\subsection{Polynomials and the shape function}

In order to obtain the asymptotic error estimates we need to study the interaction of the projection operator $\interp$ with polynomials.

The notation ${\alpha}$ always refers to a  $d$-vector
of non-negative integers
$$
\alpha = (\alpha_1,\cdots , \alpha_d) \in \ZZ_+^d.
$$
For each $\alpha$ we define
the following quantities
$$
|\alpha|:= \sum_{1\leq i\leq d} \alpha_i, \quad \alpha! := \prod_{1\leq i\leq d} \alpha_i!, \quad \max(\alpha) := \max_{1\leq i\leq d} \alpha_i.
$$
We also define the monomial
$$
X^\alpha := \prod_{1\leq i\leq d} X_i^{\alpha_i},
$$
where the variable is $X=(X_1,...,X_d)\in \RR^d$.
Finally, for each integer $k\geq 0$  we define the following three vector spaces of polynomials
\be
\label{defPk}
\begin{array}{lcl}
\P_k    & :=& \Vect\{X^\alpha \sep |\alpha| \leq k\}, \\
\P_k^* &:=& \Vect\{X^\alpha \sep \max(\alpha) \leq k \text{ and } |\alpha| \leq k+1\}, \\
\P_k^{**} &:=& \Vect\{X^\alpha \sep \max(\alpha) \leq k\}.
\end{array}
\ee 
Note that clearly $\dim(\P_k^{**}) = (k+1)^d$.
In addition, a classical combinatorial argument shows that 
$$
\dim \P_k = \binom {k+d}{d} \stext{ and }\dim\P_k^* = \dim\P_{k+1} - d = \binom {k+d+1}{d} - d.
$$

By $V_{\interp}$ we denote the image of $\interp$, which is a subspace of $V = C^0(\mI^d)$. Since
$\interp$ is a projector \iref{projAxiom}, we have \be \label{defVI}
V_{\interp} = \{\interp(f): \;\; f\in V\} = \{f \in V : \;\;
f=\interp(f)\}. \ee From this point on, the integer $k$ is
fixed and defined as follows \be \label{defk} k = k(\interp) := \max
\{k'\geq 0 \sep \P_{k'} \subset V_{\interp} \} \ee Hence, the
operator $\interp$ reproduces polynomials of total degree less or
equal than $k$. (If $k=\infty$ then we obtain, using the density of
polynomials in $V$ and the continuity of $\interp$, that $\interp(f)
= f$ for all $f\in V$. We exclude this case from now on.)

In what follows, by $m$ we denote the integer defined by
\be
\label{defm}
m = m(\interp) :=k+1,
\ee
where $k = k(\interp)$ is defined in \iref{defk}.
By $\H_m$ we denote the space of homogeneous polynomials of degree $m$
$$
\H_m := \Vect\{ X^\alpha\sep |\alpha| = m\}.
$$
We now introduce a function $K_I$ on $\H_m$, further referred to as the ``shape function''.
\begin{definition}[Shape Function]
For all $\pi \in \H_m$ 
\be 
\label{defK} K_I(\pi) := \inf_{|R|=1}
\|\pi - \interp_R \pi\|_{L^p(R)}, 
\ee 
where the infimum is taken
over all blocks $R$ of unit $d$-dimensional volume.
\end{definition}

The shape function $K$ plays a major role in our asymptotical error estimates developed in the next subsection.
Hence, we dedicate \S \ref{studyK} to its close study, and we provide its explicit form in various cases.
Note that for any $A>0$, since $\pi$ is homogeneous of degree $m$
\be
\label{infAreaA}
\inf_{|R|=A} \|\pi - \interp_R \pi\|_{L^p(R)} = A^{\frac m d+\frac 1 p} K_I(\pi).
\ee

The optimization \iref{defK} among blocks can be rephrased into an  optimization among diagonal matrices.
Indeed for any block $R$ there exists a unique $x_0\in \R^d$ and a unique diagonal matrix with positive coefficients such that $R = \phi(\mI^d)$ with $\phi(x) = x_0+ Dx$. Furthermore, the homogeneous component of degree $m$ is the same in both $\pi\circ \phi$ and $\pi\circ D$, hence $\pi\circ \phi - \pi\circ D \in \P_k$ (recall that $m = k+1$) and therefore this polynomial is reproduced by the projection operator $\interp$. Using the linearity of $\interp$, we obtain
$$
\pi \circ \phi - \interp (\pi \circ \phi) = \pi \circ D - \interp (\pi \circ D).
$$
Combining this with \iref{changeRect}, and observing that $\det D=|R|$, we obtain that
\be
\label{KD}
K_I(\pi) = \inf_{\substack{\det D = 1\\D\geq 0}} \|\pi\circ D - \interp(\pi\circ D)\|_{L^p(\mI^d)},
\ee
where the infimum is taken over the set of diagonal matrices with non-negative entries and unit determinant.

\subsection{Examples of projection operators}
\label{exampleI}

In this section we define several possible choices for the projection operator $\rm I$ which are consistent with \iref{defk} and, in our opinion, are most useful for practical purposes.
However, many other possibilities could be considered.
\begin{definition}[$L^2(\mI^d)$ orthogonal projection]
\label{L2Proj}
We may define $\interp(f)$ as the $L^2(\mI^d)$ orthogonal projection of $f$ onto one of the spaces of polynomials $\P_k$,
$\P_k^*$ or $\P_k^{**}$ defined in \iref{defPk}.
\end{definition}

If the projection operator $\interp$ is chosen as in Definition \ref{L2Proj}, then a simple change of variables shows that for any block $R$, the operator $\interp_R$ defined by \iref{interpR} is the $L^2(R)$ orthogonal projection onto the same space of polynomials.



To introduce several possible interpolation schemes for which we obtain the estimates using our approach,
we consider a set $U_k\subset \mI$ of cardinality $\#(U_k) = k+1$ (special cases are given below).
For any $\u = (u_1, \cdots u_d) \in U_k^d$ we define an element of $\P_k^{**}$ as follows
$$
\mu_\u (X):= \prod_{1 \leq i\leq d} \left( \prod_{\substack{v\in U_k\\ v\neq u_i}} \frac{X_i - v}{u_i-v}\right) \in \P_k^{**}.
$$

Clearly, 
$
\mu_\u (\u) = \mu_\u (u_1, \cdots , u_d) = 1
$
and 
$
\mu_\u(\v) = \mu_\u(v_1, \cdots , v_d) = 0
$
if $\v = (v_1, \cdots , v_d)\in U_k^d$ and $\v\neq \u$. 

It follows that the elements of  $B := (\mu_\u)_{\u\in U_k^d}$ are linearly independent. Since $\# (B) = \#(U_k^d) = (k+1)^d = \dim(\P_k^{**})$, $B$ is a basis of $\P_k^{**}$.

Therefore, any element of $\mu \in \P_k^{**}$ can be written in the form
$$
\mu(X) = \sum_{\u\in U_k^d} \lambda_\u \mu_\u(X).
$$
It follows that, for any given $f\in C^0(\mI^d)$, there exists a unique element of $\mu\in \P_k^{**}$ such that $\mu(\u) = f(\u)$ for all $\u\in U_k^d$. We define $\interp f := \mu$, namely 
$$
(\interp f)(X) := \sum_{\u\in U_k^d} f(\u) \mu_\u(X) \in \P_k^{**}.
$$
We may take $U_k$ to be the set of $k+1$ equi-spaced points on $\mI$
\be
\label{interpEqui}
U_k = \left\{-\frac 1 2 + \frac n k \sep 0 \leq n \leq k\right\}.
\ee
We obtain a different, but equally relevant, operator $\interp$ by choosing $U_k$ to be the set of Tchebychev points on $\mI$
\be
\label{interpTche}
U_k = \left\{\frac 1 2 \cos\left( \frac {n\pi} k\right) \sep 0 \leq n \leq k\right\}.
\ee
Different interpolation procedures can be used to construct $\interp$. Another convenient interpolation scheme is to take
$$
\interp(f) \in \P_k^*
$$
and $\interp(f) = f$ on a subset of $U_k^d$. This subset  contains $\dim \P_k^*$ points, which are convenient to choose first on the boundary of $\mI^d$ and then (if needed) at some interior lattice points.
Note that since $\dim \P_k^*< \#(U_k^d) =(k+1)^d$, it is always possible to construct such an operator.

If the projection operator $\interp$ is chosen as described above, then for any block $R$ and any $f\in C^0(R)$, $\interp_R(f)$ is the unique element of respective space of polynomials which coincides with $f$ at the image $\phi(p)$ of the points $p$ mentioned in the definition of $\interp$, by the transformation $\phi$ described in \iref{defPhi}.

\subsection{Main results}
\label{subsecApprox}
In order to obtain the approximation results we often impose a slight technical restriction (which can be removed, see for instance ~\cite{BBLS}) on sequences of block partitions, which is defined as follows.
\begin{definition}[Admissibility]
\label{defAdmiRect}
We say that a sequence $(\cR_N)_{N\geq 1}$ of block partitions of a block $R_0$ is \emph{admissible} if $\#(\cR_N) \leq N$ for all $N \geq 1$, and
\be
\label{defadmi}
\sup_{N\geq 1} \left( N^{\frac 1 d} \sup_{R \in \cR_N} \diam(R)\right) < \infty
\ee
\end{definition}

%

We recall that the approximation error is measured in $L^p$ norm, where the exponent $p$ is fixed and $1\leq p \leq \infty$.
We define $\tau\in (0, \infty)$ by
\be
\label{defTau}
\frac 1 \tau := \frac m d+ \frac 1 p.
\ee
In the following estimates we identified $d^m f(x)$ with an element of $\H_m$ according to
\be
\label{dmfHm}
\frac{d^m f(x)}{m!} \sim \sum_{|\alpha| = m} \frac{\partial^m f(x)}{\partial x^\alpha}  \frac{X^\alpha}{\alpha!}.
\ee
We now state the asymptotically sharp lower bound for the approximation error of a function $f$ on an admissible sequence of block partitions.
\begin{theorem}
\label{thLower}
Let $R_0$ be a block and let $f\in C^m(R_0)$. For any admissible sequence of block partitions $(\cR_N)_{N\geq 1}$ of $R_0$
$$
\liminf_{N\to \infty} N^{\frac m d}\|f-\interp_{\cR_N} f\|_{L^p(R_0)} \geq \left\|K_I\left(\frac{d^m f}{m!}\right)\right\|_{L^\tau(R_0)}.
$$
\end{theorem}

The next theorem provides an upper bound for the projection error of a function $f$ when an optimal sequence of block partitions is used. It confirms the sharpness of the previous theorem.

\begin{theorem}
\label{thUpper}
Let $R_0$ be a block and let $f\in C^m(R_0)$. Then there exists a (perhaps non-admissible) sequence $(\cR_N)_{N\geq 1}$, $\#(\cR_N) \leq N$, of block partitions of $R_0$ satisfying
\be
\label{upperEstim}
\limsup_{N\to \infty} N^{\frac m d}\|f-\interp_{\cR_N} f\|_{L^p(R_0)} \leq \left\|K_I\left(\frac{d^m f}{m!}\right)\right\|_{L^\tau(R_0)}.
\ee

Furthermore, for all $\ve>0$ there exists an admissible sequence $(\cR_N^\ve)_{N\geq 1}$ of block partitions of $R_0$ satisfying
\be
\label{upperEstimEps}
\limsup_{N\to \infty} N^{\frac m d}\|f-\interp_{\cR_N^\ve} f\|_{L^p(R_0)} \leq \left\|K_I\left(\frac{d^m f}{m!}\right)\right\|_{L^\tau(R_0)}+\ve.
\ee
\end{theorem}
An important feature of these estimates is the ``$\limsup$''.
Recall that the upper limit of a sequence $(u_N)_{N\geq N_0}$ is defined by
$$
\limsup_{N\to \infty} u_N := \lim_{N\to \infty} \sup_{n\geq N} u_n,
$$
and is in general strictly smaller than the supremum $\sup_{N\geq
N_0} u_N$. It is still an open question to find an appropriate upper
estimate of $\sup_{N\geq N_0} N^{\frac{m} d}  \|f-\interp_{\cR_N}
f\|_{L^p(R_0)}$ when optimally adapted block partitions are used.

In order to have more control of the quality of approximation on various parts of the domain we introduce a positive weight function $\Omega\in C^0(R_0)$.
For $1\leq p\leq \infty$ and for any $u\in L^p(R_0)$ as usual we define
\be 
\label{defWeight}
\|u\|_{L^p(R_0, \Omega)} := \|u\Omega\|_{L^p(R_0)}. 
\ee

\begin{remark}
\label{remarkWeight} Theorems \ref{thLower}, \ref{thUpper} and
\ref{thNoEps} below also hold when the norm $\| \cdot\|_{L^p(R_0)}$ (resp
$\| \cdot\|_{L^\tau(R_0)}$) is replaced with the weighted norm $\| \cdot\|_{L^p(R_0,
\Omega)}$ (resp $\| \cdot\|_{L^\tau(R_0, \Omega)}$) defined in \iref{defWeight}.
\end{remark}

In the following section we shall use some restrictive hypotheses on the interpolation operator in order to obtain an explicit formula for the shape function.
In particular, 
Propositions \ref{propOdd}, \ref{propEven}, and equation \iref{Kdmf} show that, under some assumptions, there exists a constant $C=C(\interp)>0$ such that
$$
\frac 1 {C} K_I\left(\frac{d^m f}{m!}\right) \leq \sqrt[d]{\left|\prod_{1\leq i\leq d} \frac{\partial^m f}{\partial x_i^m}\right|} \leq C K_I\left(\frac{d^m f}{m!}\right).
$$
These restrictive hypotheses also allow to improve slightly the estimate \iref{upperEstimEps} as follows.
\begin{theorem}
\label{thNoEps}
If the hypotheses of Proposition \ref{propOdd} or \ref{propEven} hold, and if
$
K_I\left(\frac{d^m f}{m!}\right) >0
$
everywhere on $R_0$, then there exists an \emph{admissible} sequence of block partitions $(\cR_N)_{N\geq 1}$ of $R_0$ which satisfies the optimal estimate \iref{upperEstim}.
\end{theorem}
The proofs of the Theorems \ref{thLower}, \ref{thUpper} and \ref{thNoEps} are given in \S\ref{proofTh}.
Each of these proofs can be adapted to weighted norms, hence establishing Remark \ref{remarkWeight}. Some details on how to adapt proofs for the case of weighted norms are provided at the end of each proof.

\section{Study of the shape function}
\label{studyK}

In this section we perform  a close study of the shape function $K_I$, since it plays a major role in our asymptotic error estimates. In the first subsection \S\ref{studyKGen} we investigate general properties which are valid for any continuous projection operator $\interp$.
However, we are not able to obtain an explicit form of $K_I$ under such general assumptions. 
Recall that in \S\ref{exampleI} we presented several possible choices of projection operators $\interp$ that seem more likely to be used in practice. In \S\ref{subsecproperties} we identify four important properties shared by these examples. 
These properties are used in \S \ref{subsecExact} to obtain an explicit form of $K_I$.

\subsection{General properties}
\label{studyKGen}
The shape function $K$ obeys the following important invariance property with respect to diagonal changes of coordinates.
\begin{prop}
\label{propInv}  For all $\pi\in \H_m$ and all diagonal matrices $D$
with non-negative coefficients
$$
K_I(\pi\circ D) = (\det D)^{\frac m d} K_I(\pi).
$$
\end{prop}

\proof
We first assume that the diagonal matrix $D$ has positive diagonal coefficients. Let $\overline D$ be a diagonal matrix with positive diagonal coefficient and which satisfies $\det \overline D = 1$.
Let also $\pi\in \H_m$. Then, since $\pi$ is homogeneous and of degree $m$
$$
\pi \circ (D \overline D) = \pi \circ ((\det D)^{\frac 1 d} \tilde D) = (\det D)^{\frac m d} \pi \circ \tilde D,
$$
where
$
\tilde D := (\det D)^{- \frac 1 d} D\overline D
$
satisfies $\det \tilde D = \det \overline D=1$ and is uniquely determined by $\overline D$.
According to \iref{KD} we therefore have
\begin{eqnarray*}
K_I(\pi \circ D) &=& \inf_{\substack{\det \overline D = 1\\ \overline D \geq 0}} \|\pi\circ (D \overline D) - \interp(\pi\circ (D \overline D))\|_{L^p(\mI^d)}\\
&=& (\det D)^{\frac m d} \inf_{\substack{\det \tilde D = 1\\ \tilde D \geq 0}} \|\pi\circ \tilde D - \interp(\pi\circ \tilde D)\|_{L^p(\mI^d)}\\
&=& (\det D)^{\frac m d} K_I(\pi),
\end{eqnarray*}
which concludes the proof in the case where $D$ has positive diagonal coefficients.\\
Let us now assume that $D$ is a diagonal matrix with non-negative diagonal coefficients and such that $\det(D) = 0$.
Clearly there exists a sequence $(D_n)_{n \geq 0}$ of diagonal matrices with positive coefficients, such that $\det D_n = 1$ for all $n \geq 0$ and that $DD_n \to 0$ as $n \to \infty$.
Therefore $\pi \circ (DD_n) \to 0$, which implies that $K(\pi \circ D) = 0$ and concludes the proof of this proposition.
\sq

The next proposition shows that the exponent $p$ used for measuring the approximation error plays a rather minor role.
By $K_p$ we denote the shape function associated with the exponent $p$.
\begin{prop}
\label{propEquiv} There exists a constant $c>0$ such that for all
$1\leq p_1\leq p_2\leq \infty$ we have on $\H_m$
$$
c K_\infty\leq K_{p_1} \leq K_{p_2} \leq K_{\infty}.
$$
\end{prop}

\proof 
For any function $f\in V = C^0(\mI^d)$ and for any $1\leq p_1\leq p_2\leq \infty$ by a standard convexity argument we obtain that
$$
\|f\|_{L^1(\mI^d)} \leq \|f\|_{L^{p_1}(\mI^d)} \leq \|f\|_{L^{p_2}(\mI^d)} \leq \|f\|_{L^{\infty}(\mI^d)}.
$$
Using \iref{KD}, it follows that
$$
K_1\leq K_{p_1}\leq K_{p_2}\leq K_\infty
$$
on $\H_m$.
Furthermore, the following semi norms on $\H_m$
$$
|\pi|_1 := \|\pi-\interp \pi\|_{L^1(\mI^d)} \text{ and } |\pi|_\infty := \|\pi-\interp \pi\|_{L^\infty(\mI^d)}
$$
vanish precisely on the same subspace of $\H_m$, namely $V_{\interp} \cap H_m = \{\pi \in
\H_m \sep \pi = \interp \pi\}$. Since $\H_m$ has finite dimension,
it follows that they are equivalent. Hence, there exists a constant
$c>0$ such that $c|\cdot |_\infty\leq |\cdot |_1$ on $\H_m$. Using
\iref{KD}, it follows that $cK_\infty\leq K_1$, which concludes the
proof. \sq

\subsection{Desirable properties of the projection operator}
\label{subsecproperties}
The examples of projection operators presented in \S\ref{exampleI} share some important properties which allow to obtain the explicit expression of the shape function $K_I$. 
These properties are defined below and called $H_\pm$, $H_\sigma$, $H_*$ or $H_{**}$.
They are satisfied when the operator $\interp$ is the interpolation at equispaced points (Definition \ref{interpEqui}), at Tchebychev points (Definition \ref{interpTche}), and usually on the most interesting sets of other points. They are also satisfied when $\interp$ is the $L^2(\mI^d)$ orthogonal projection onto $\P_k^*$ or $\P_k^{**}$ (Definition \ref{L2Proj}).

The first property reflects the fact that a coordinate $x_i$ on $\mI^d$ can be changed to $-x_i$, independently of the projection process.
\begin{definition}[$H_\pm$ hypothesis]
We say that the interpolation operator $\interp$ satisfies the $H_\pm$ hypothesis if for any diagonal matrix $D$ with entries in $\pm 1$ we have for all $f\in V$
$$
 \interp(f\circ D) = \interp(f) \circ D.
$$
\end{definition}

The next property implies that the different coordinates $x_1, \cdots, x_d$ on $\mI^d$ play symmetrical roles with respect to the projection operator.
\begin{definition}[$H_\sigma$ hypothesis]
If $M_\sigma$ is a permutation matrix, i.e. $(M_\sigma)_{ij} :=
\delta_{i \sigma(j)}$ for some permutation $\sigma$ of $\{1,\cdots,
d\}$, then for all $f\in V$
$$
\interp (f\circ M_\sigma) = \interp( f) \circ M_\sigma.
$$
\end{definition}

According to \iref{defk}, the projection operator $\interp$ reproduces the space of polynomials $\P_k$. However, in many situations the space $V_{\interp}$ of functions reproduced by $\interp$ is larger than $\P_k$. In particular $V_{\interp} = \P_k^{**}$ when $\interp$ is the interpolation on equispaced or Tchebychev points, and $V_{\interp} = \P_k$ (resp $\P_k^*$,  $\P_k^{**}$) when $\interp$ is the $L^2(\mI^d)$ orthogonal projection onto $\P_k$ (resp $\P_k^*$,  $\P_k^{**}$).

It is particularly useful to know whether the projection operator $\interp$ reproduces the elements of $\P_k^*$, and we therefore give a name to this property. Note that it clearly does not hold for the $L^2(\mI^d)$ orthogonal projection onto $\P_k$.
\begin{definition}[$H_*$ hypothesis]
The following inclusion holds :
$$
P_k^* \subset V_{\interp}.
$$
\end{definition}

On the contrary it is useful to know that some polynomials, and in particular pure powers $x_i^m$, are not reproduced by $\interp$.
\begin{definition}[$H_{**}$ hypothesis]
$$
\text{If } \ \sum_{1\leq i\leq d} \lambda_i x_i^m \in V_{\interp} \ \text{ then } \ (\lambda_1, \cdots, \lambda_d) = (0, \cdots , 0).
$$
\end{definition}
This condition obviously holds if $\interp(f)\in \P_k^{**}$ (polynomials of degree $\leq k$ in each variable) for all $f$. Hence, it holds for all the examples of projection operators given in the previous subsection \S\ref{exampleI}.

\subsection{Explicit formulas}
\label{subsecExact}

In this section we provide the explicit expression of the shape function $K$ when some of the hypotheses $H_\pm$, $H_\sigma$, $H_*$ or $H_{**}$ hold.
Let $\pi\in \H_m$ and let $\lambda_i$ be the coefficient of $X_i^m$ in $\pi$, for all $1\leq i\leq d$.
We define
$$
K_*(\pi) := \sqrt[d]{\prod_{1\leq i\leq d} |\lambda_i|}
$$
and
$$
s(\pi) := \#\{1\leq i\leq d\sep \lambda_i >0\}.
$$
If $\frac{d^m f(x)}{m!}$ is identified by  \iref{dmfHm}  to an element of $\H_m$, then one has 
\be
\label{Kdmf}
K_*\left(\frac{d^m f(x)}{m!}\right) = \frac 1 {m!} \sqrt[d]{\left|\prod_{1\leq i\leq d} \frac{\partial^m f}{\partial x_i^m}(x)\right|}.
\ee

\begin{prop}
\label{propOdd}
If $m$ is odd and if $H_\pm$, $H_\sigma$ and $H_*$ hold, then
$$
K_p(\pi) = C(p) K_*(\pi),
$$
where
$$
 C(p) := \left\|\sum_{1\leq i\leq d} X_i^m - \interp \left(\sum_{1\leq i\leq d} X_i^m\right)\right\|_{L^p(\mI^d)}>0.
$$
\end{prop}

\begin{prop}
\label{propEven}
If $m$ is even and if $H_\sigma$, $H_*$ and $H_{**}$ hold then
$$
K_p(\pi) = C(p,s(\pi)) \, K_*(\pi).
$$
Furthermore,
\be
\label{eqEven}
C(p,0) = C(p,d) =
\left\|\sum_{1\leq i\leq d} X_i^m - \interp
\left(\sum_{1\leq i\leq d} X_i^m\right)\right\|_{L^p(\mI^d)}>0.
\ee
Other constants $C(p,s)$ are positive and obey $C(p,s) =
C(p,d-s)$.\\
\end{prop}

\noindent
Next we turn to the proofs of Propositions \ref{propOdd} and \ref{propEven}.
\paragraph{Proof of Proposition \ref{propOdd}}
Let $\pi\in \H_m$ and let $\lambda_i$ be the coefficient of $X_i^m$ in $\pi$.
Denote by
$$
\pi_* := \sum_{ 1\leq i\leq d} \lambda_i X_i^m
$$
so that $\pi-\pi_* \in \P_k^*$ and, more generally, $\pi \circ D - \pi_* \circ D \in \P_k^*$ for any diagonal matrix $D$.
The hypothesis $H_*$ states that the projection operator $\interp$ reproduces the elements of $\P_k^*$, and therefore 
$$
\pi\circ D - \interp (\pi\circ D) = \pi_*\circ D - \interp (\pi_*\circ D).
$$
Hence,
$
K_I(\pi) = K_I(\pi_*)
$
according to \iref{KD}.
If there exists $i_0$, $1\leq i_0\leq d$, such that $\lambda_{i_0} = 0$, then we denote by $D$ the diagonal matrix of entries $D_{ii} = 1$ if $i\neq i_0$ and $0$ if $i=i_0$. Applying Proposition \ref{propInv} we find 
$$
K_I(\pi) = K_I(\pi_*) = K_I(\pi_*\circ D) = (\det D)^{\frac m d} K_I(\pi_*) = 0.
$$
which concludes the proof.
We now assume that all the coefficients $\lambda_i$, $1\leq i\leq d$, are different from $0$, and we denote by $\ve_i$ be the sign of $\lambda_i$. Applying Proposition
\ref{propInv} to the diagonal matrix $D$ of entries $|\lambda_i|^{\frac
1 m}$ we find that
$$
K_I(\pi) = K_I(\pi_*) = (\det D)^{\frac m d} K_I(\pi_*\circ D^{-1}) = K_*(\pi) \, K_I\left(\sum_{1\leq i\leq d} \ve_i X_i^m\right).
$$
Using the $H_\pm$ hypothesis with the diagonal matrix $D$ of entries
$D_{ii} = \ve_i$, and recalling that $m$ is odd, we find that
$$
K_I\left(\sum_{1\leq i\leq d} \ve_i X_i^m\right) = K_I\left(\sum_{1\leq i\leq d}  X_i^m\right).
$$
We now define the functions
$$
g_i := X_i^m - \interp(X_i^m) \text{ for } 1\leq i \leq d.
$$
It follows from \iref{KD} that
$$
K_I\left(\sum_{1\leq i\leq d}  X_i^m\right) = \inf_{\prod_{1\leq
i\leq d} a_i = 1} \left\|\sum_{1\leq i\leq d} a_i
g_i\right\|_{L^p(\mI^d)},
$$
where the infimum is taken over all $d$-vectors of positive reals of product $1$.
Let us consider such a $d$-vector $(a_1, \cdots , a_d)$, and a permutation $\sigma$ of the set $\{1, \cdots,
d\}$.
The $H_\sigma$ hypothesis implies that the quantity
$$
\left\|\sum_{1\leq i\leq d} a_{\sigma(i)} g_i\right\|_{L^p(\mI^d)}
$$
is independent of $\sigma$. Hence, summing over all permutations, we obtain using the triangle inequality
\be
\label{sumPerm}
\begin{array}{rcl}
\displaystyle
\left\|\sum_{1\leq i\leq d} a_i
g_i\right\|_{L^p(\mI^d)} 
&=&
\displaystyle
 \frac 1 {d!} \sum_\sigma
\left\|\sum_{1\leq i\leq d} a_{\sigma(i)} g_i\right\|_{L^p(\mI^d)}\\
&\geq &
\displaystyle
\frac 1 {d!}
\left\| \sum_{1\leq i\leq d} \left(\sum_\sigma a_{\sigma(i)}\right)g_i\right\|_{L^p(\mI^d)}\\
&=&
\displaystyle
\frac 1 d \left\|\sum_{1\leq i\leq d} g_i\right\|_{L^p(\mI^d)}
\sum_{1\leq i\leq d} a_i. 
\end{array}
\ee
The right-hand side is minimal when $ a_1 =
\cdots = a_d = 1 $, which shows that
$$
\left\|\sum_{1\leq i\leq d} a_i g_i\right\|_{L^p(\mI^d)} \geq  \left\|\sum_{1\leq i\leq d} g_i\right\|_{L^p(\mI^d)} = C(p)
$$
with equality when $a_i=1$ for all $i$. 
Note as a corollary that 
\be
\label{existsRPos}
K_I(\pi_\ve) = \|\pi_\ve - \interp(\pi_\ve)\|_{L^p(\mI^d)} = C(p) \ \text{ where } \ \pi_\ve = \sum_{1\leq i\leq d} \ve_i X_i^m.
\ee
It remains to prove that $C(p)>0$. Using the hypothesis $H_\pm$, we find that for all $\mu_i \in\{\pm 1\}$ we have
$$
 \left\|\sum_{1\leq i\leq d}\mu_i g_i\right\|_{L^p(\mI^d)} = C(p).
$$
In particular, for any $1\leq i_0\leq d$ one has 
$$
2\|g_{i_0}\|_{L^p(\mI^d)} \leq  \left\|\sum_{1\leq i\leq d} g_i\right\|_{L^p(\mI^d)} +  \left\|2g_{i_0} - \sum_{1\leq i\leq d} g_i\right\|_{L^p(\mI^d)} \leq 2 C(p).
$$
If $C(p) = 0$, it follows
that $g_{i_0} = 0$ and therefore that $X_{i_0}^m = \interp(X_{i_0}^m)$, for any $1\leq
i_0\leq d$. Using the assumption $H_*$, we find that the projection operator
$\interp$ reproduces all the polynomials of degree $m= k+1$, which
contradicts the definition \iref{defk} of the integer $k$.
\sq

\paragraph{Proof of proposition \ref{propEven}}
We define $\lambda_i$, $\pi_*$ and $\ve_i\in\{\pm 1\}$ as before and
we find, using similar reasoning, that
$$
K_I(\pi) = K_*(\pi) K_I\left(\sum_{1\leq i\leq d} \ve_i X_i^m\right).
$$
\noindent
For $1\leq s\leq d$ we define
$$
C(p, s) := K_I\left(\sum_{1\leq i\leq s}  X_i^m - \sum_{s+1\leq i\leq d}  X_i^m\right).
$$
From the hypothesis $H_\sigma$ it follows that $K_I(\pi) = K_*(\pi)
C(p,s(\pi))$.

Using again $H_\sigma$ and the fact that $K_I(\pi) = K_I(-\pi)$ for all $\pi \in \H_m$, we find that
\begin{eqnarray*}
C(p,s) &=& K_I\left(\sum_{1\leq i\leq s}  X_i^m - \sum_{s+1\leq i\leq d}  X_i^m\right)\\
&=& K_I\left(-\left(\sum_{1\leq i\leq d-s}  X_i^m - \sum_{d-s+1\leq i\leq d}  X_i^m\right)\right)\\
&=& C(p,d-s).
\end{eqnarray*}

We define $g_i := X_i^m - \interp(X_i^m)$, as in the proof of Proposition \ref{propOdd}.
We obtain the expression for $C(p,0)$ by summing over all permutations as in  \iref{sumPerm}
$$
C(p,0) = \left\|\sum_{1\leq i\leq d} g_i\right\|_{L^p(\mI^d)}.
$$
This concludes the proof of the first part of Proposition \ref{propEven}.
We now prove that $C(p,s)>0$ for all $1\leq p\leq \infty$ and all $s\in \{0, \cdots, d\}$.
To this end we define the following quantity on $\R^d$
$$
 \|a\|_K := \left\|\sum_{1\leq i\leq d}a_i  g_i\right\|_{L^p(\mI^d)}.
$$
Note that $\|a\|_K = 0$ if and only if
$$
\sum_{1\leq i\leq d} a_i X_i^m = \sum_{1\leq i\leq d} a_i \interp(X_i^m),
$$
and the hypothesis $H_{**}$ precisely states that this equality occurs if and only if $a_i = 0$, for all $1\leq i\leq d$. Hence, $\|\cdot \|_K$ is a norm on $\R^d$.
Furthermore, let
$$
E_s := \left\{a\in \R_+^s\times \R_-^{d-s}\sep \prod_{1\leq i \leq
d} |a_i| = 1\right\}
$$
Then
$$
C(p,s) = \inf_{a\in E_s} \|a\|_K.
$$
Since $E_s$ is a closed subset of $\R^d$, which does not contain the
origin, this infimum is attained. It follows that $C(p,s)>0$, and that there exists a rectangle $R_\ve$ of unit volume such that
\be
\label{existsREven}
K_I(\pi_\ve) = \|\pi_\ve - \interp \pi_\ve\|_{L^p(R_\ve)} = C(p, s(\pi_\ve)) \ \text{ where } \ \pi_\ve = \sum_{1\leq i\leq d} \ve_i X_i^m,
\ee
which concludes the proof of this proposition.
\sq

\section{Proof of the approximation results}
\label{proofTh}

In this section, let the
block $R_0$, the integer $m$, the function $f\in C^m(R_0)$ and the exponent $p$ be fixed. We conduct our proofs for $1\leq p<\infty$ and provide comments on how to adjust our arguments for the case $p=\infty$.

For each $x\in \R_0$ by $\mu_x\in \P_m$ we denote
 the  $m$-th degree Taylor
polynomial of $f$ at the point $x$
\be
\label{defmux}
\mu_x = \mu_x(X) := \sum_{|\alpha| \leq m}\frac{\partial^m
f}{\partial x^\alpha}(x) \, \frac{(X-x)^\alpha}{\alpha!},
\ee
and we define $\pi_x\in \H_m$ to be the homogeneous component of degree $m$ in
$\mu_x$,
\be
\label{defpix}
\pi_x = \pi_x(X) :=  \sum_{|\alpha| = m}\frac{\partial^m f}{\partial x^\alpha}(x) \, \frac{X^\alpha}{\alpha!}.
\ee
Since $\pi_x$ and $\mu_x$ are polynomials of degree $m$, their $m$-th derivative is constant, and clearly $d^m \pi_x=d^m \mu_x = d^m
f(x)$. 
In particular, for any $x\in R_0$ the polynomial $\mu_x - \pi_x$ belongs to $\P_k$ (recall that $k=m-1$) and is therefore reproduced by the projection operator $\interp$. It follows that for
any $x\in R_0$ and any block $R$
\be
\label{muzPiz} \pi_x- \interp_R \pi_x = \mu_x
-\interp_R \mu_x.
\ee
In addition, we introduce a measure $\rho$ of
the degeneracy of a block $R$
$$
\rho(R) := \frac{\diam(R)^d}{|R|}.
$$
Given any function $g\in C^m(R)$ and any $x\in R$ we can define, similarly to \iref{defpix}, a polynomial $\hat \pi_x\in \H_m$ associated to $g$ at $x$.
We then define 
\be
\label{normdmg}
\|d^m g\|_{L^\infty(R)}:=\sup_{x\in R} \left(\sup_{|u| = 1} |\hat \pi_x(u)|\right).
\ee
\begin{prop}
\label{propLocalIso}
There exists a constant $C = C(m,d)>0$ such that for any block $R$ and any function $g\in C^m(R)$
\be
\label{localIso}
\|g-\interp_R g\|_{L^p(R)} \leq C |R|^{\frac 1 \tau} \rho(R)^{\frac m d} \|d^m g\|_{L^\infty(R)}.
\ee
\end{prop}

\proof
Let $x_0\in R$ and let $g_0$ be the Taylor polynomial for $g$ of degree $m-1$ at point $x_0$ which is defined as follows
$$
g_0(X) := \sum_{|\alpha| \leq m-1}\frac{\partial^\alpha f(x_0)}{\partial x^\alpha}\frac{(X-x_0)^\alpha}{\alpha!}.
$$
Let $x\in R$ and let $x(t) =  x_0+ t (x-x_0)$.
We have
$$
g(x) = g_0(x)+ \int_{t=0}^1 d^m g_{x(t)}(x-x_0) \frac{(1-t)^m}{m!} dt.
$$
Hence,
\be
\label{ineqGG0}
\begin{array}{rcl}
 |g(x) - g_0(x)| &\leq&  \displaystyle\int_{t=0}^1 \|d^m g\|_{L^\infty(R)} |x-x_0|^m \frac{(1-t)^m}{m!}dt\\
&\leq&\displaystyle \frac 1 {(m+1)!} \|d^m g\|_{L^\infty(R)} \diam(R)^m.
\end{array}
\ee
Since $g_0$ is a polynomial of degree at most $m-1$, we have $g_0 = \interp g_0$. Hence,
\begin{eqnarray*}
\|g-\interp_R g\|_{L^p(R)} &\leq &  |R|^{\frac 1 p} \|g-\interp_R g\|_{L^\infty(R)}\\
& =& |R|^{\frac 1 p} \|(g-g_0)-\interp_R (g-g_0)\|_{L^\infty(R)}\\
&\leq & (1+C_{\interp}) |R|^{\frac 1 p} \|g-g_0\|_{L^\infty(R)},
\end{eqnarray*}
where $C_{\interp}$ is the operator norm of $\interp : V\to V$. Combining this estimate with \iref{ineqGG0}, we obtain \iref{localIso} which concludes the proof.
\sq

\subsection{Proof of Theorem \ref{thLower} (Lower bound)}

The following lemma allows us to bound the interpolation error of $f$ on the block $R$ from below.
\begin{lemma}
\label{lemmaLower} For any block $R\subset R_0$ and $x\in R$ we have
$$
\|f- \interp_R f\|_{L^p(R)} \geq |R|^{\frac 1 \tau} \left(K_I(\pi_x) - \omega(\diam R) \rho(R)^{\frac m d}\right),
$$
where the function $\omega$ is positive, depends only on $f$ and $m$, and satisfies $\omega(\delta) \to 0$ as $\delta\to 0$.
\end{lemma}

\proof
Let $h:= f-\mu_x$, where $\mu_x$ is defined in \iref{defmux}. Using \iref{muzPiz} and \iref{infAreaA}, we obtain
\begin{eqnarray*}
\|f- \interp_R f\|_{L^p(R)} & \geq & \|\pi_x - \interp_R \pi_x\|_{L^p(R)} - \|h - \interp_R h\|_{L^p(R)}\\
& \geq &  |R|^{\frac 1 \tau} K_I(\pi_x) - \|h - \interp_R
h\|_{L^p(R)},
\end{eqnarray*}
and according to Proposition \ref{propLocalIso} we have
$$
 \|h - \interp_R h\|_{L^p(R)} 
\leq C_0 |R|^{\frac 1 \tau}
\rho(R)^{\frac m d} \|d^m h\|_{L^\infty(R)}.
$$
Observe that
$$
\|d^m h\|_{L^\infty(R)} = \|d^m f - d^m\pi_x\|_{L^\infty(R)}  =
\|d^m f - d^m f(x)\|_{L^\infty(R)}. 
$$
We introduce the modulus of continuity $\omega_*$ of the $m$-th derivatives of $f$. 
\be
\label{defOmega}
\omega_*(r) := \sup_{\substack{x_1,x_2\in R_0 :\\ |x_1 - x_2|\leq r}} \|d^m f(x_1)-d^m f(x_2)\|= \sup_{\substack{x_1,x_2\in R_0 :\\ |x_1-x_2|\leq r}} \left( \sup_{|u|\leq 1} |\pi_{x_1}(u) - \pi_{x_2} (u)|\right)
\ee
By setting $\omega = C_0\, \omega_*$
we conclude the proof of
this lemma.
\sq

We now consider an admissible sequence of block partitions $(\cR_N)_{N\geq 0}$. For all $N\geq 0$,
$R\in \cR_N$ and $x\in \interior(R)$, we define
$$
\phi_N(x) := |R| \qquad \hbox{and} \qquad
\psi_N(x) := \left(K_I(\pi_x) - \omega(\diam(R)) \rho(R)^{\frac m d} \right)_+,
$$
where $\lambda_+ := \max\{\lambda,0\}$.
We now apply Holder's inequality
$$
\int_{R_0} f_1 f_2 \leq \|f_1\|_{L^{p_1}(R_0)}\|f_2\|_{L^{p_2}(R_0)}
$$
with the functions
$$
f_1 = \phi_N^{\frac {m \tau} d} \psi_N^\tau \ \text{ and } f_2 = \phi_N^{-\frac {m \tau} d}
$$
and the exponents
$
p_1 = \frac p \tau \ \text{ and } \ p_2 = \frac d {m\tau}.
$
Note that $\frac 1 {p_1}+ \frac 1 {p_2} =  \tau \left(\frac 1 p + \frac m d\right) = 1$. Hence,
\be
\label{holderPsi}
\int_{R_0} \psi_N^\tau \leq \left(\int_{R_0} \phi_N^{\frac{m p} d}
\psi_N^{p}\right)^{\frac \tau p} \left(\int_{R_0}
\phi_N^{-1}\right)^{\frac {m\tau} d}.
\ee

\noindent
Note that $\int_{R_0} \phi_N^{-1} = \# (\cR_N)\leq N$. Furthermore, if
$R\in \cR_N$ and $x\in \interior(R)$ then according to Lemma \ref{lemmaLower}
$$
 \phi_N(x)^{\frac m d} \psi_N(x)
=|R|^{\frac 1 \tau-\frac 1 p}  \psi_N(x) \leq|R|^{-\frac 1 p} \|f-\interp_R f\|_{L^p(R)}.
$$
Hence,
\be
\label{intphipsi}
\left[\int_{R_0} \phi_N^{\frac{m p} d} \psi_N^{p}\right]^{\frac 1 p} \leq \left[\sum_{R\in \cR_N} \frac 1 {|R|}  \int_R  \|f-\interp_R f\|_{L^p(R)}^p\right]^{\frac 1 p}=
\|f-\interp_R f\|_{L^p(R_0)}.
\ee
Inequality \iref{holderPsi} therefore leads to
\be
\label{upperPsiRect}
 \|\psi_N\|_{L^\tau(R_0)} \leq \|f- \interp_{\cR_N} f\|_{L^p(R_0)} N^{\frac m d}.
\ee

\noindent
Since the sequence $\seqR$ is admissible, there exists a constant $C_A>0$ such that for all $N$ and all $R\in \cR_N$ we have $\diam(R)\leq C_AN^{-\frac 1 d}$.
We introduce a subset of $\cR'_N\subset \cR_N$ which collects the most degenerate blocks
$$
\cR'_N = \{ R\in \cR_N \sep \rho(R)\geq \omega(C_AN^{-\frac 1
d})^{-\frac{1} m}\},
$$
where $\omega$ is the function defined in Lemma \ref{lemmaLower}. 
By $R'_N$ we denote the portion of $R_0$ covered by $\cR'_N$.
For all $x\in R_0\sm R'_N$
we obtain
$$
\psi_N(x)\geq K_I(\pi_x) -\omega(C_A N^{-\frac 1 d})^{1-\frac 1 d}.
$$
We define $\ve_N := \omega(C_A N^{-\frac 1 d})^{1-\frac 1 d}$ and we notice that $\ve_N \to 0$ as $N \to \infty$.
Hence,
$$
\begin{array}{ll}
 \|\psi_N\|_{L^\tau(R_0)}^\tau & \geq \left \|\left(K_I(\pi_x) -\ve_N\right)_+\right\|_{L^\tau(R_0\sm R'_N)}^\tau\\
 & \geq
 \left \|\left(K_I(\pi_x) -\ve_N\right)_+\right\|_{L^\tau(R_0)}^\tau
 -C^\tau |R'_N|,
 \end{array}
 $$
where $C:=\max_{x\in R_0}K_I(\pi_x)$. The last expression involves a slight abuse of notations, since $(K_I(\pi_x) -\ve_N)_+$ stands for the function $x\in R_0 \mapsto (K_I(\pi_x) -\ve_N)_+ \in \R$. Next we observe
that $|R'_N|\to 0$ as $N\to +\infty$: indeed
for all $R\in \cR'_N$ we have
$$
|R| = \diam(R)^d \rho(R)^{-1} \leq C_A^d N^{-1} \omega(C_A N^{-\frac 1 d})^{\frac 1 m}.
$$
Since $\#(\cR'_N)\leq N$, we obtain $|R'_N|\leq C_A^d \omega(C_A N^{-\frac 1 d})^{\frac 1 m}$, and the right-hand side tends to $0$ as $N\to \infty$. We thus obtain
$$
\liminf_{N\to \infty} \|\psi_N\|_{L^\tau(R_0)}
\geq \lim_{N\to \infty}   \left \|\left(K_I(\pi_x) -\ve_N\right)_+\right\|_{L^\tau(R_0)}
 = \|K_I(\pi_x)\|_{L^\tau(R_0)}.
$$
Combining this result with \iref{upperPsiRect}, we conclude the proof of the announced estimate.

Note that this proof also works with the exponent $p = \infty$ by changing
$$
\left(\int_{R_0} \phi_N^{\frac{m p} d}
\psi_N^{p}\right)^{\frac \tau p} \ \text{ into } \ \|\phi_N^{\frac m d}
\psi_N\|_{L^\infty(R_0)}^\tau
$$
in \iref{holderPsi} and performing the standard modification in \iref{intphipsi}.
\begin{remark}
As announced in Remark \ref{remarkWeight}, this proof can be adapted to the weighted norm $\|\cdot\|_{L^p(R_0, \Omega)}$ associated to a positive weight function $\Omega\in C^0(R_0)$ and defined in \iref{defWeight}. For that purpose let $r_N := \sup \{ \diam(R) \sep R \in \cR_N\}$ and let
$$
\Omega_N(x) := \inf_{\substack{x'\in R_0\\ |x-x'|\leq r_N}} \Omega(x').
$$
The sequence of functions $\Omega_N$ increases with $N$ and tends uniformly to $\Omega$ as $N\to \infty$.
If $R\in \cR_N$ and $x\in R$, then
$$
\|f-\interp_R f\|_{L^p(R,\Omega)}
\geq  \Omega_N(x) \|f-\interp_R f\|_{L^p(R)}.
$$
The main change in the proof is that the function $\psi_N$ should be replaced with $\psi'_N := \Omega_N \psi_N$. Other details are left to the reader.
\end{remark}
\sq

\subsection{Proof of the upper estimates}

The proof of Theorems \ref{thUpper} (and \ref{thNoEps}) is based on the actual construction of an asymptotically optimal sequence of block partitions. To that end we introduce the notion of a local block specification.

\begin{definition}{\bf{(local block specification)}}
\label{defBlockSpec}
A local block specification on a block $R_0$ is a (possibly discontinuous) map $x \mapsto R(x)$ which associates to each point $x\in R_0$ a block $R(x)$, 
and such that
\begin{itemize}
\item
The volume
$
|R(x)|
$
is a positive continuous function of the variable $x \in R_0$.
\item The diameter is bounded : $\sup \{\diam(R(x))\sep x \in R_0\}<\infty$. 
\end{itemize}
\end{definition}

The following lemma shows that it is possible to build sequences of block partitions of $R_0$ adapted in a certain sense to a given local block specification.

\begin{lemma}
\label{lemmaSeqBlock}
Let $R_0$ be a block in $\RR^d$ and let $x\mapsto R(x)$ be a local block specification on $R_0$. Then there exists a sequence $(\cP_n)_{n\geq 1}$ of block partitions of $R_0$,
$
\cP_n = \cP_n^1 \cup \cP_n^2,
$
satisfying the following properties.
\begin{itemize}
\item (The number of blocks in $\cP_n$ is asymptotically controlled)
\be \label{limCardRn} \lim_{n \to \infty} \frac{\#(\cP_n)}{n^{2d}} =
\int_{R_0} |R(x)|^{-1} dx. \ee
\item (The elements of $\cP_n^1$ follow the block specifications)
For each $R\in \cP_n^1$ there exists $y\in R_0$ such that
\be
\label{n2Ry}
R \text{ is a translate of } n^{-2} R(y), \text{ and } |x-y| \leq \frac{\diam(R_0)} n \text{ for all } x\in R.
\ee
\item (The elements of $\cP_n^2$ have a small diameter)
\be \label{smallDiam} \lim_{n \to \infty } \left( n^2 \sup_{R\in
\cP_n^2} \diam(R)\right) =0. \ee
\end{itemize}
\end{lemma}

\proof
See Appendix.
\sq

We recall that the block $R_0$, the exponent $p$ and the function $f\in C^m(R_0)$ are fixed, and that at each point $x\in R_0$ the polynomial $\pi_x\in \H_m$ is defined by \iref{defpix}.
The sequence of block partitions described in the previous lemma is now used to obtain an asymptotical error estimate. 
\begin{lemma}
\label{lemmaNn} Let $x \mapsto R(x)$ be a local block specification
such that for all $x\in R_0$ \be \label{unitError} \|\pi_x
-\interp_{R(x)}(\pi_x)\|_{L^p(R(x))} \leq 1. \ee Let $(\cP_n)_{n\geq
1}$ be a sequence of block partitions satisfying the properties of
Lemma \ref{lemmaSeqBlock}, and let for all $N\geq 0$
$$
n(N) := \max\{ n\geq 1  \sep \# (\cP_n) \leq N\}.
$$
Then $\cR_N := \cP_{n(N)}$ is an admissible sequence of block partitions and
\be
\label{limsupR}
\limsup_{N\to \infty} N^{\frac m d} \|f-\interp_{\cR_N} f \|_{L^p(R_0)} \leq \left(\int_{R_0} R(x)^{-1} dx\right)^{\frac 1 \tau}.
\ee
\end{lemma}

\proof
Let $n \geq 0$ and let $R\in \cP_n$.
If $R\in \cP_n^1$ then let $y\in R_0$ be as in \iref{n2Ry}. Using Proposition \ref{propLocalIso} and \iref{infAreaA} we find
\begin{eqnarray*}
\|f-\interp_R f\|_{L^p(R)} &\leq& \|\pi_y-\interp_R \pi_y\|_{L^p(R)}  + \|(f-\pi_y)-\interp_R (f-\pi_y)\|_{L^p(R)} \\
&\leq & n^{-\frac {2d} \tau}  \|\pi_y-\interp_{R(y)} \pi_y\|_{L^p(R(y))} + C |R|^{\frac 1 p} \diam(R)^m \|d^m f-d^m \pi_y\|_{L^{\infty(R)}}\\
&\leq &  n^{-\frac {2d} \tau}  + C n^{-\frac {2d} \tau} |R(y)|^{\frac 1 p} \diam(R(y))^m \|d^m f-d^m f(y)\|_{L^\infty(R)}\\
&\leq & n^{-\frac {2d} \tau}  (1+C' \omega_*(n^{-1}\diam(R_0))),
\end{eqnarray*}
where we defined $C' := C \sup_{y\in \R_0} |R(y)|^{\frac 1 p} \diam(R(y))^m$, which is finite by Definition \ref{defBlockSpec}. 
We denoted by $\omega_*$ the modulus of continuity of the $m$-th derivatives of $f$ which is defined at \iref{defOmega}.
We now define for all $n \geq 1$, 
$$
\delta_n := n^2 \sup_{R \in \cP_n^2} \diam(R).
$$
According to \iref{smallDiam} one has $\delta_n \to 0$ as $n \to \infty$.
If $R\in \cP_n^2$, then $\diam(R)\leq n^{-2} \delta_n$ and therefore $|R|\leq \diam(R)^d\leq n^{-2d} \delta_n^d$. Using again \iref{localIso}, and recalling that $\frac 1 \tau = \frac m d+ \frac 1 p$ we find
$$
\|f-\interp_R f\|_{L^p(R)} \leq C |R|^{\frac 1 p} \diam(R)^m \|d^m f\|_{L^{\infty(R_0)}} \leq C'' n^{-\frac {2d} \tau} \delta_n^{\frac d \tau}
$$
where $C'' = C \|d^m f\|_{L^{\infty(R_0)}}$.
From the previous observations it follows that
\begin{eqnarray*}
\|f-\interp_{\cP_n} f\|_{L^p(R_0)} &\leq& \#(\cP_n)^{\frac 1 p}
\max_{R\in \cP_n} \|f-\interp_R f\|_{L^p(R)}\\
&\leq&  \#(\cP_n)^{\frac
1 p}  n^{-\frac {2d} \tau}  \max\{1+
C'\omega_*(n^{-1}\diam(R_0)), \, C'' \delta_n^{\frac d
\tau}\}.
\end{eqnarray*}
Hence,
$$
\limsup_{n\to \infty} \#(\cP_n)^{-\frac 1 p}  n^{\frac {2d} \tau}\|f-\interp_{\cP_n} f\|_{L^p(R_0)}  \leq 1.
$$
Combining the last equation with \iref{limCardRn}, we obtain
$$
\limsup_{n\to \infty} \#(\cP_n)^{\frac m d} \|f-\interp_{\cP_n} f\|_{L^p(R_0)}  \leq \left(\int_{R_0} R(x)^{-1} dx\right)^{\frac 1 \tau}.
$$
The sequence of block partitions $\cR_N := \cP_{n(N)}$
clearly satisfies $\#(\cR_N)/N\to 1$ as $N
\to \infty$ and therefore leads to the announced equation \iref{limsupR}.
Furthermore, it follows from the boundedness of $\diam(R(x))$ on
$R_0$ and the properties of $\cP_n$ described in Lemma
\ref{lemmaSeqBlock} that
$$
\sup_{n\geq 1}\left( \#(\cP_n)^{\frac 1 d} \sup_{R \in \cP_n}
\diam(R) \right)< \infty
$$
which implies that $\cR_N$ is an admissible sequence of partitions.
\sq

We now choose adequate local block specifications in order to obtain
the estimates announced in Theorems \ref{thUpper} and \ref{thNoEps}.
For any $M\geq \diam(\mI^d) = \sqrt d$ 
we define the modified shape function 
\be 
\label{defKM}
K_M(\pi) := \inf_{\substack{|R| = 1,\\ \diam(R)\leq M}} \|\pi - \interp_R \pi\|_{L^p(R)}, 
\ee
where the infimum is taken on blocks of unit volume and diameter smaller that $M$.
It follows from a compactness argument that this infimum is attained and that $K_M$ is a continuous function on $\H_m$. Furthermore, for any fixed $\pi\in \H_m$, $M\mapsto K_M(\pi)$ is a decreasing function of $M$ which tends to $K_I(\pi)$ as $M \to \infty$.

For all $x\in R_0$ we denote by $R_M^*(x)$ a block which realises the infimum in $K_M(\pi_x)$. Hence,
$$
|R_M^*(x)| = 1, \ \diam(R_M^*(x))\leq M, \text{ and } K_M(\pi_x) = \|\pi_x - \interp_{R_M^*(x)} \pi_x\|_{L^p(R_M^*(x))}
$$
We define a local block specification on $R_0$ as follows
\be
\label{defRM}
R_M(x) := (K_M(\pi_x) + M^{-1})^{- \frac \tau d} R_M^*(x).
\ee
We now observe that using a change of variables and the homogeneity of $\pi_x$, as in \iref{infAreaA}, that 
$$
 \|\pi_x - \interp_{R_M(x)} \pi_x\|_{L^p(R_M(x))}  = K_M(\pi_x) (K_M(\pi_x)+ M^{-1})^{-1} \leq 1.
$$
Hence, according to Lemma \ref{lemmaNn}, there exists a sequence $(\cR_N^M)_{N\geq 1}$ of block partitions of $R_0$ such that
$$
\limsup_{N\to \infty} N^{\frac m d} \|f-\interp_{\cR^M_N} f \|_{L^p(R_0)} \leq \|K_M(\pi_x)+M^{-1}\|_{L^\tau (R_0)}.
$$
Using our previous observations on the function $K_M$, we see that
$$
\lim_{M \to \infty} \|K_M(\pi_x)+M^{-1}\|_{L^\tau (R_0)} = \|K_I(\pi_x)\|_{L^\tau (R_0)}.
$$
Hence, given $\ve >0$ we can choose $M(\ve)$ large enough in such a way that
$$
\|K_{M(\ve)}(\pi_x)+M(\ve)^{-1}\|_{L^\tau (R_0)} \leq \|K_I(\pi_x)\|_{L^\tau (R_0)}+ \ve,
$$
which concludes the proof of the estimate \iref{upperEstimEps} of Theorem \ref{thUpper}.

For each $N$ let $M(N)$ be such that
$$
N^{\frac m d} \|f-\interp_{\cR^{M(N)}_N} f \|_{L^p(R_0)} \leq \|K_{M(N)}(\pi_x)+M(N)^{-1}\|_{L^\tau (R_0)}+M(N)^{-1}
$$
and $M(N) \to \infty$ as $N \to \infty$. Then the (perhaps non admissible) sequence of block partitions $\cR_N := \cR_N^{M(N)}$ satisfies \iref{upperEstim} which concludes the proof of Theorem \ref{thUpper}.
\sq

We now turn to the proof of Theorem \ref{thNoEps}, which follows the same scheme for the most. 
There exists $d$ functions $\lambda_1(x)
, \cdots, \lambda_d(x) \in C^0(R_0)$, and a function $x \mapsto \pi_*(x) \in \P^*_k$ such that
for all $x\in R_0$ we have
$$
\pi_x = \sum_{1\leq i\leq d} \lambda_i(x) X_i^m + \pi_*(x).
$$
The hypotheses of Theorem \ref{thNoEps} state that $K_I\left(\frac {d^m f(x)}{m!}\right) = K_I(\pi_x)\neq 0$ for all $x\in \R_0$. It follows from Propositions \ref{propOdd} and \ref{propEven} that the product $\lambda_1(x)\cdots \lambda_d(x)$ is nonzero for all $x\in R_0$.
We denote by $\ve_i\in \{\pm 1\}$ the sign of $\lambda_i$, which is therefore constant over the block $R_0$, and we define
$$
\pi_\ve := \sum_{1\leq i\leq d} \ve_i X_i^m
$$
The proofs of Propositions \ref{propEven} and \ref{propOdd} show that there exists a block $R_\ve$, satisfying $|R_\ve| = 1$, and such that $K_I(\pi_\ve) = \|\pi-\interp_{R_\ve} \pi\|_{L^p(R_\ve)}$.
By $D(x)$ we denote the diagonal matrix of entries $|\lambda_1(x)|, \cdots , |\lambda_d(x)|$, and we define
$$
\phi_x :=  (\det D(x))^{\frac 1 {md} } D(x)^{- \frac 1 m}.
$$
Clearly $\det \phi_x = 1$ and $\pi_x \circ \phi_x = (\det D(x))^\frac 1 d \pi_\ve+ \pi_*(x) \circ \phi_x$, and $\pi_*(x) \circ \phi_x\in \P_k^*$.
Hence using \iref{changeRect} we obtain 
\begin{eqnarray*}
 \|\pi_x - \interp_{\phi_x(R_\ve)} \pi_x\|_{L^p(\phi_x(R_\ve))} &=&   \|\pi_x\circ \phi_x - \interp_{R_\ve} (\pi_x\circ \phi_x)\|_{L^p(R_\ve)} \\
 &=& (\det D(x))^\frac 1 d \|\pi_\ve - \interp_{R_\ve} \pi_\ve\|_{L^p(R_\ve)} \\
 &=&  (\det D(x))^{\frac 1 d} K_I(\pi_\ve)\\
  &=& K_I(\pi_x).
\end{eqnarray*}
We then define the local block specification
\be
\label{defR}
R(x) := K_I(\pi_x)^{-\frac \tau d}\phi_x(R_\ve),
\ee
in such way that $ \|\pi_x - \interp_{R(x)} \pi_x\|_{L^p(R(x))} = 1$ for all $x\in R_0$, using the homogeneity of $\pi_x$ and an isotropic change of variables.
The admissible sequence $(\cR_N)_{N \geq 1}$ of block partitions constructed in Lemma \ref{lemmaNn} then satisfies the optimal upper estimate \iref{upperEstim}, which concludes the proof of Theorem \ref{thNoEps}. \sq

\begin{remark}[Adaptation to weighted norms]
Lemma \ref{lemmaNn} also holds if \iref{unitError} is replaced with
$$\Omega(x) \|\pi_x -\interp_{R(x)}(\pi_x)\|_{L^p(R(x))} \leq 1$$
and if the $L^p(R_0)$
norm is replaced with the weighted $L^p(R_0, \Omega)$ norm in \iref{limsupR}.
Replacing the block $R_M(x)$ defined in \iref{defRM} with
$$
R'_M(x) := \Omega(x)^{- \frac \tau d}R_M(x),
$$
one can easily obtain the extension of Theorem \ref{thUpper} to
weighted norms. Similarly, replacing $R(x)$ defined in \iref{defR}
with $R'(x) := \Omega(x)^{- \frac \tau d}R(x)$, one obtains the
extension of Theorem \ref{thNoEps} to weighted norms.
\end{remark}

\section{Appendix : Proof of Lemma \ref{lemmaSeqBlock}}
\begin{figure}
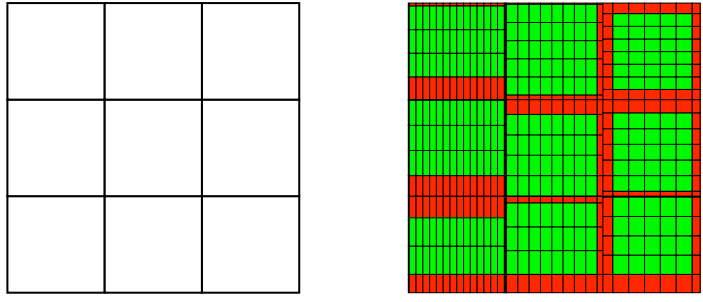

    \centering
    \includegraphics[width=4cm,height=4cm]{\pathPic/PaperRectangles/frame.pdf}
    \hspace{1cm}
    \includegraphics[width=4cm,height=4cm]{\pathPic/PaperRectangles/tile.pdf}
    \caption{\label{fig1Rectangles}(Left) the initial uniform (coarse) tiling $\cQ_3$ of $R_0$. (Right) the set of blocks $\cP_n^1$ in green and the set of blocks $\cP_n^{2*}$ in red.}
\end{figure}

By $\cQ_n$ we denote  the standard partition of $R_0\in \RR^d$ in $n^d$ identical blocks of diameter $n^{-1} \diam(R_0)$ illustrated on the left in Figure \ref{fig1Rectangles}.
For each $Q \in \cQ_n$ by $x_Q$ we denote the barycenter of $Q$ and we consider the tiling $\cT_Q$ of $\R^d$ formed with the block $n^{-2}R(x_Q)$ and its translates.
We define $\cP_n^1(Q)$ and $\cP_n^1$ as follows
$$
\cP_n^1(Q) := \{R\in \cT_Q \sep R \subset Q\} \ \stext{ and } \ \cP_n^1 := \bigcup_{Q\in \cQ_n} \cP_n^1(Q).
$$
Comparing the areas, we obtain
$$
\# (\cP_n^1) = \sum_{Q\in \cQ_n} \cP_n^1(Q) \leq \sum_{Q \in \cQ_n} \frac{|Q|}{|n^{-2} R(x_Q)|} = n^{2 d} \sum_{Q \in \cQ_n} |Q| |R(x_Q)|^{-1}.
$$
From this point, using the continuity of $x \mapsto |R(x)|$, one easily shows that
$
\frac {\# (\cP_n^1)}{n^{2d}} \to \int_{R_0} |R(x)|^{-1} dx
$
as $n \to \infty$. Furthermore, the property \iref{n2Ry} clearly holds. In order to construct $\cP_n^2$, we first define two sets of blocks $\cP_n^{2*}(Q)$ and $\cP_n^{2*}$ as follows
$$
\cP_n^{2*}(Q) := \{R\cap Q \sep R \in \cT_Q \text{ and }\interior(R)\cap \partial Q\neq \emptyset\} \ \text{ and } \ \cP_n^{2*} := \bigcup_{Q\in \cQ_n} \cP_n^{2*}(Q).
$$
Comparing the surface of $\partial Q$ with the dimensions of $R(x_Q)$, we find that
$$
\#(\cP_n^{2*}(Q)) \leq C n^{d-1}
$$
where $C$ is independent of $n$ and of $Q\in \cQ_n$.
Therefore, $\#(\cP_n^{2*})\leq C n^{2d-1}$. The set of blocks $\cP_n^2$ is then obtained by subdividing each block of $\cP_n^{2*}$ into $o(n)$ (for instance, $\lfloor \ln(n)\rfloor^d$) identical sub-blocks, in such a way that $\#(\cP_n^2)$ is $o(n^{2d})$ and that the requirement \iref{smallDiam} is met.

\chapter[Sharp asymptotics of the $L^p$ interpolation error]{Sharp asymptotics of the $L^p$ interpolation error on optimal triangulations} 
\minitoc
\label{secOptAniso}
\label{chapOptAniso}

%
%

\section{Introduction.}
\label{secOptAniso1}
\subsection{Optimal mesh adaptation}

In finite element approximation, a usual distinction is between {\it uniform} and
{\it adaptive} methods. In the latter, the elements defining the mesh 
may vary strongly in size and shape for a better adaptation
to the local features of the approximated function $f$. This naturally raises
the objective of characterizing and constructing an {\it optimal mesh}
for a given function $f$. 

Note that depending on the context, the function $f$ may be fully
known to us: either through an explicit formula or a discrete sampling; 
or observed through noisy measurements;
or implicitly defined as the solution of a given
partial differential equation.

In this chapter, we assume that $f$ is a function 
defined on a polygonal bounded domain $\Omega\subset \R^2$.
For a given conforming triangulation $\cT$ of $\Omega$
and an arbitrary but fixed 
integer $m>1$, we denote by $\interp^{m-1}_\cT$ 
the standard interpolation operator on the 
space of Lagrange finite elements of degree $m-1$ associated to $\cT$.
Given a norm $X$ of interest and a number
$N>0$, the objective of finding the optimal mesh for $f$ can be
formulated as solving the optimization problem
$$
\min_{\#(\cT)\leq N} \|f-\interp^{m-1}_\cT f\|_X,
$$
where the minimum is taken over all conforming triangulations of 
cardinality $N$. We denote by $\cT_N$ the minimizer of the
above problem.

Our first objective is to establish
sharp asymptotic error estimates that precisely describe the
behavior of $\|f-\interp^{m-1}_\cT f\|_X$ as $N\to +\infty$. Estimates of that type
were obtained in \cite{BBLS,CSX} in the particular case of linear
finite elements ($m-1=1$) and with the error measured in $X=L^p$. They have the form
\be
\limsup_{N\to +\infty} \(N\min_{\#(\cT)\leq N}\|f-\interp^1_\cT f\|_{L^p}\)  \leq C\left\| \sqrt{|\det(d^2f)|}\right\|_{L^\tau},\;\; \frac 1 \tau:=\frac 1 p+1,
\label{optiaffine}
\ee
which reveals that the convergence rate is governed by the quantity $\sqrt{|\det(d^2f)|}$, which depends
nonlinearly the on Hessian $d^2f$. This is heavily tied to the fact that we allow
triangles with possibly highly anisotropic shape.
In the present work, the polynomial degree $m-1$ is arbitrary,
and the quantities governing the convergence rate will therefore depend nonlinearly on the 
$m$-th order derivative $d^mf$.

Our second objective is to propose simple and practical ways
of designing meshes which behave similarly to the optimal one,
in the sense that they satisfy the sharp error estimate up
to a fixed multiplicative constant.

\subsection{Main results and layout}
\label{secOptAniso1p2}
We denote by $\H_m$ the space of homogeneous polynomials of degree $m$:
$$
\H_m:={\rm Span}\{ x^ky^l\sep k+l=m\}.
$$
For any triangle $T$, we denote by $\interp^{m-1}_T$ the local interpolation
operator acting from $C^0(T)$ onto the space $\P_{m-1}$ of 
polynomials of total degree $m-1$. The image of $v\in C^0(T)$ by this operator is defined by the
conditions
$$
\interp^{m-1}_T v(\gamma)=v(\gamma)
$$
for all points $\gamma\in T$ with barycentric coordinates in
the set $\{0, \frac 1 {m-1},\frac 2 {m-1},\cdots,1\}$. We denote by
$$
e_{m,T}(v)_p:=\|v-\interp^{m-1}_T v\|_{L^p(T)}
$$
the interpolation error measured in the norm $L^p(T)$. We also denote
by
$$
e_{m,\cT}(v)_p:=\|v-\interp^{m-1}_\cT v\|_{L^p}=\left(\sum_{T\in\cT}e_{m,T}(v)_p^p\right)^{\frac 1 p}
$$
the global interpolation error for a given triangulation $\cT$, with the standard modification if $p=\infty$.

A key ingredient in this chapter is a function 
defined by a {\it shape optimization problem}:
for any fixed $1\leq p\leq \infty$ and for any $\pi\in \H_m$, we define
\be
K_{m,p}(\pi):=\inf_{|T|=1}e_{m,T}(\pi)_p.
\label{shapefunction}
\ee
Here, the infimum is taken over all triangles of area $|T|=1$. 
Note that from the homogeneity of $\pi$, we find that
\be
\inf_{|T|=A}e_{m,T}(\pi)_p=K_{m,p}(\pi)A^{\frac m 2 + \frac 1 p}.
\label{shapeA}
\ee
This optimization problem thus gives the shape of the triangles 
of a given area which is best adapted
to the polynomial $\pi$ in the sense of minimizing the interpolation error
measured in $L^p$.  
We refer to $K_{m,p}$ as the {\it shape function}.
We discuss in \S \ref{secOptAniso2} the main properties of this function. 
\nl
\nl
Our asymptotic error estimate for the optimal triangulation is given by
the following theorem:
\begin{theorem}
\label{maintheorem}
For any bounded polygonal domain $\Omega\subset \R^2$ and any function 
$f\in C^m(\overline \Omega)$, there exists a sequence 
of triangulations $\seqT$, conforming if $p<\infty$, with $\#(\cT_N)\leq N$, such that
\be
\label{Cf}
\limsup_{N\ra \infty} N^{\frac m 2} e_{m,\cT_N}(f)_p\leq 
\left\|K_{m,p}\left(\frac{d^m f}{m!}\right)\right\|_{L^\tau(\Omega)},\ \frac 1 \tau := \frac m 2 +\frac 1 p .
\ee
\end{theorem}
An important feature of this estimate is the ``$\limsup$'' asymptotical operator. Recall that the upper limit of a sequence $(u_N)_{N\geq N_0}$ is defined by 
$$
\limsup_{N\to \infty} u_N := \lim_{N\to \infty} \sup_{n\geq N} u_n
$$
and is in general strictly smaller than the supremum $\sup_{N\geq N_0} u_N$. It is still an open question to find an appropriate upper estimate for $\sup_{N\geq N_0} N^{\frac m 2} e_{m,\cT_N}(f)_p$ when optimally adapted anisotropic triangulations are used.

In the estimate \iref{Cf}, the $m$-th derivative $d^mf(z)$ is identified to 
a homogeneous polynomial in $\H_m$:
\be
\label{defDMFPol}
\frac{d^m f(z)}{m!}\sim \sum_{k+l=m}\frac{\partial^m f}{\partial^kx\partial^l y}(z)\frac{x^k}{k!}\frac{y^l}{l!}.
\ee
In order to illustrate the sharpness of \iref{Cf}, we introduce
a slight restriction on sequences of triangulations, following 
an idea in \cite{BBLS}: a sequence $\seqT$ of triangulations, such that $\#(\cT_N) \leq N$, is said to be \emph{admissible} if
$$
\sup_{N \geq N_0} \left (N^\frac 1 2 \sup_{T\in \cT_N} \diam(T)\right) < \infty.
$$
In other words
\be
\label{admissibilitycond}
\sup_{T\in \cT_N} \diam(T) \leq C_AN^{-1/2}
\ee
for some $C_A>0$ independent of $N$. The following theorem shows that the estimate
\iref{Cf} cannot be improved when we restrict our attention to admissible sequences.
It also shows that this class is reasonably large in the sense that 
\iref{Cf} is ensured to hold up to small perturbation.

\begin{theorem}
\label{optitheorem}
Let $\Omega\subset \R^2$ be a bounded polygonal domain, and $f\in C^m(\overline \Omega)$. 
Set $ \frac 1 \tau := \frac m 2 +\frac 1 p$.
For all \emph{admissible} sequences of triangulations $\seqT$, conforming or not, one has
$$
\liminf_{N\ra \infty} N^{\frac m 2} e_{m,\cT_N}(f)_p \geq \left\|K_{m,p}\left(\frac{d^m f}{m!}\right)\right\|_{L^\tau(\Omega)}.
$$
For all $\ve>0$, there exists an \emph{admissible} sequence of triangulations $(\cT_N^\ve)_{N\geq N_0}$, conforming if $p<\infty$, such that
$$
\limsup_{N\ra \infty} N^{\frac m 2} e_{m,\cT_N^\ve}(f)_p \leq \left\|K_{m,p}\left(\frac{d^m f}{m!}\right)\right\|_{L^\tau(\Omega)}+\ve.
$$
\end{theorem}
Note that the sequences $(\cT_N^\ve)_{N\geq N_0}$ satisfy the admissibility condition \iref{admissibilitycond} with a constant $C_A(\ve)$ which may explode as $\ve\to 0$.
The proofs of both theorems are given in \S \ref{secOptAniso3}. These proofs reveal that the construction of the 
optimal triangulation obeys two principles:
(i) the triangulation should {\it equidistribute} the local approximation error $e_{m,T}(f)_p$ between 
each triangle, and (ii) the aspect ratio of a triangle $T$ should be 
{\it isotropic} with respect to a distorted
metric induced by the local value of $d^mf$ on $T$
(and therefore anisotropic in the sense of the Euclidean metric). 
Roughly speaking, the quantity $\|K_{m,p}\left(\frac{d^m f}{m!}\right)\|_{L^\tau(T)}$
controls the local interpolation $L^p$-error estimate on a triangle $T$ once this triangle  
is optimized with respect to the local properties of $f$.
This type of estimate differs from those obtained in \cite{Apel}, which hold for any $T$, optimized
or not, and involve the partial derivatives of $f$ in a local coordinate system which is 
adapted to the shape of $T$.

The proof of the upper estimates in Theorem \ref{optitheorem} involves the construction
of an optimal mesh based on a patching strategy similar to \cite{BBLS}. 
However, inspection of the proof reveals that
this construction becomes effective only when 
the number of triangles $N$ becomes very large. Therefore it 
may not be useful in practical applications.

A more practical approach consists in deriving the above mentioned distorted metric from 
the exact or approximate data of $d^mf$ using the following procedure. 
To any $\pi\in \H_m$, we associate a 
symmetric positive definite matrix $h_\pi\in S_2^+$. If $z\in \Omega$, and $d^m f(z)$ is close to $\pi$, then the triangle $T$ containing $z$ should be isotropic in the metric $h_\pi$.
The global metric is given
at each point $z$ by
$$
h(z)=s(\pi_z) h_{\pi_z},\;\; \pi_z=d^mf(z),
$$
where $s(\pi_z)$ is a scalar factor which depends on the desired accuracy
of the finite element approximation. 
Once this metric has been properly identified, fast algorithms such as in \cite{Inria, Bamg, Peyre}
can be used to design a near-optimal mesh based on it. Recently in \cite{Shew,Bois}, several algorithms have been rigorously proven to terminate and produce good quality meshes, see also Chapter 5. 
Computing the map
\be
\pi\in \H_m \mapsto h_\pi\in S_2^+
\label{mappi}
\ee
is therefore of key use in applications. This problem is well understood
in the case of linear elements ($m=2$): the matrix $h_\pi$ is then defined as the
absolute value (in the sense of symmetric matrices) of the matrix associated to
the quadratic form $\pi$. In contrast, the exact form of this map in the case $m\geq 3$
is not well understood. 

In this chapter, we propose
algebraic strategies for computing the map \iref{mappi} for $m=3$ which corresponds to
quadratic elements. These strategies
have been implemented in an open-source 
Mathematica code \cite{sitejm}. 
In a similar manner, we address the algebraic computation 
of the shape function $K_{m,p}(\pi)$
from the coefficients of $\pi\in \H_m$, when $m\geq 3$. 
All these questions are addressed in \S \ref{secOptAniso4}, \ref{secOptAniso5} and \ref{secOptAniso6}.

In \S \ref{secOptAniso4}, we discuss the particular case of linear ($m=2$)
and quadratic ($m=3$) elements. In this case, it is possible
to obtain explicit formulas for $K_{m,p}(\pi)$ from the coefficients
of $\pi$. In the case $m=2$, this formula is of the form
$$
K_{2,p}(ax^2+2bxy+cy^2)=\sigma\sqrt{|b^2-ac|},
$$
where the constant $\sigma$ only depends on $p$ and the sign of $b^2-ac$,
and we therefore recover the known estimate \iref{optiaffine} from Theorem \ref{maintheorem}. The formula for $m=3$ involves the
discriminant of the third degree polynomial $d^3f$. Our analysis
also leads to an algebraic computation of the map \iref{mappi}.
We want to mention that a different strategy for 
the construction of the distorted metric and
the derivation of the error estimate for a finite element
of arbitrary order was proposed in
\cite{C3}. In this approach, the distorted
metric is obtained at a point $z\in\Omega$ by finding the largest
ellipse contained in a level set of the polynomial associated to $d^mf(z)$ by \iref{defDMFPol}.
This optimization problem has connections with the one
that defines the shape function in \iref{shapefunction},
as we shall explain in \S \ref{secOptAniso2}. The approach 
proposed in the present work in the case $m=3$ has the advantage of
avoiding the use of numerical optimization,
the metric being directly derived from the
coefficients of $d^mf$. 

In \S \ref{secOptAniso5}, we address the case $m>3$. In this case, explicit formulas for $K_{m,p}(\pi)$
seem out of reach. However, we can introduce explicit functions 
$\Kpol_m(\pi)$ which are polynomials in the coefficients of $\pi$, 
and are equivalent to $K_{m,p}(\pi)$, leading therefore to similar
asymptotic error estimates up to multiplicative constants. At the current stage,
we did not obtain a simple solution to the
algebraic computation of the map \iref{mappi} in the case $m>3$.
The derivation of $\Kpol_m$ is based on the 
theory of invariant polynomials due to Hilbert. Let us mention
that this theory was also recently applied in \cite{OST} to image processing tasks
such as affine invariant edge detection and denoising.

We finally discuss in \S \ref{secOptAniso6} the possible extension of our analysis to simplicial
elements in higher dimension. This extension is not straightforward except in the case of
linear elements $m=2$.

\section{The shape function}
\label{secOptAniso2}
In this section, we establish several properties of the 
function $K_{m,p}$ which will be of key use subsequently.
We assume that $m\geq 2$ is an integer, and $p\in[1,\infty]$.
We equip the finite dimensional vector space $\H_m$ with a norm $\|\cdot\|$
defined as follows
\be
\label{normdef}
\|\pi\| := \sup_{x^2+y^2 \leq 1} |\pi(x,y)|.
\ee
Our first result shows that the function $K_{m,p}$
vanishes on a set of polynomials which has a simple algebraic characterization.
\begin{prop}
\label{vanishprop}
We denote by $\mhalf:= \lfloor \frac m 2 \rfloor +1$ the smallest integer strictly larger than $m/2$.
The vanishing set of $K_{m,p}$ is the set of polynomials which have a generalized 
root of multiplicity at least $\mhalf$:
$$
K_{m,p}(\pi)=0 \Leftrightarrow  \pi(x,y) = (\alpha x +\beta y)^\mhalf \tilde \pi, \mbox{ for some }
\alpha,\beta \in \R\mbox{ and } \ti\pi \in \H_{m-\mhalf}.
$$
\end{prop}

\proof
We denote by $\TEq$ a fixed equilateral triangle of unit area, centered at $0$.

We first assume that $\pi(x,y)= (\alpha x+\beta y)^\mhalf \tilde \pi$. Then there exists a rotation  
$R\in \cO_2$ and $\hat \pi\in H_{m-\mhalf}$ such that 
$$
\pi\circ R(x,y)= x^\mhalf \hat \pi(x,y)=x^\mhalf \left( \sum_{i=0}^{m-s_m}a_ix^iy^{m-s_m-i} \right).
$$
Therefore, denoting by $\phi_\ve$ the linear transform $\phi_\ve(x,y) = R\left(\ve x ,\frac y \ve\right)$, 
we obtain 
$$
\pi \circ \phi_\ve = (\ve x)^\mhalf \left( \sum_{i=0}^{m-s_m}a_i (\ve x)^i (y/\ve)^{m-s_m-i} \right) = \ve^{2\mhalf-m} x^\mhalf \left( \sum_{i=0}^{m-s_m}\ve^{2i}a_ix^iy^{m-s_m-i} \right),
$$
hence $\pi \circ \phi_\ve \to 0$ as $\ve \to 0$.
Since $|\det \phi_\ve| = 1$, the triangles $\phi_\ve(\TEq)$ have unit area. Consequently, 
$$
e_{m,\phi_\ve(\TEq)} (\pi)_p = e_{m,\TEq} (\pi\circ \phi_\ve)_p\to 0\; \; {\rm as}\;\; \ve\to 0,
$$
and therefore $K_{m,p}(\pi) = 0$.

Conversely, let $\pi\in \H_m\sm\{0\}$ be such that $K_{m,p}(\pi)=0$. Then there exists a sequence $(T_n)_{n\geq 0}$ of triangles with unit area such that $e_{m,T_n}(\pi)_p\to 0$.
We remark that the interpolation error $e_T(\pi)_p$ of $\pi\in \H_m$ is 
invariant by a translation $\tau_h: z\mapsto z+h$ of the triangle $T$. Indeed, 
$\pi-\pi\circ\tau_h\in \P_{m-1}$, so that
\be
\|\pi- \interp^{m-1}_T \pi\|_{L^p(\tau_h(T))}=
\|\pi\circ \tau_h - \interp^{m-1}_T (\pi\circ\tau_h)\|_{L^p(T)}=
\|\pi- \interp^{m-1}_T \pi\|_{L^p(T)}.
\label{transinv}
\ee
Hence we may assume that the barycenter of $T_n$ is $0$
and write $T_n = \phi_n(\TEq)$ for some linear transform $\phi_n$ with $\det \phi_n = 1$.
Since $e_{m,\TEq}(\cdot)_p$ is a norm on $\H_m$, it follows 
that $\pi\circ \phi_n\ra 0$.

The linear transform $\phi_n$ has a singular value decomposition
$$
\phi_n = U_n \circ D_n \circ V_n, \text{ where } U_n,V_n\in \cO_2,\text{ and } D_n = \left(
\begin{array}{cc}
\ve_n & 0\\
0 & 1/\ve_n
\end{array}
\right)
,\ 0<\ve_n\leq 1.
$$
Note that $\|\pi \circ V\| = \|\pi\|$ for any $\pi \in \H_m$ and $V\in \cO_2$.
Therefore, 
$$
\|\pi\circ U_n \circ D_n\| = \|\pi\circ U_n \circ D_n \circ V_n \circ V^{-1}_n\|
=\|\pi\circ\phi_n\| \to 0.
$$
Denoting by $a_{i,n}$ the coefficient of $x^i y^{m-i}$ in
$\pi\circ U_n$, we find that $a_{i,n}\ve_n^{2i-m}$ tends to $0$ as $n\to +\infty$.
In the case where $i<s_m$, this implies that $a_{i,n}$ tends to $0$ as $n\to +\infty$.

By compactness of $\cO_2$
we may assume, up to a subsequence extraction, that $U_n$ converges to some $U\in \cO_2$. 
Denoting by $a_i$ the coefficient of  $x^i y^{m-i}$ in
$\pi\circ U$, we thus find that $a_i=0$ if $i<s_m$. 
This implies that $\pi\circ U (x,y) = x^\mhalf \hat \pi(x,y)$ which concludes the proof.
\sq

\begin{remark}
In the simple case $m=2$, we infer from Proposition \ref{vanishprop}
that $K_{2,p}(\pi)=0$ if and only if $\pi$ is of the form $\pi(x,y)=x^2$ up to a rotation, 
and therefore a one-dimensional function. For such a function, the optimal
 triangle $T$ degenerates to a segment in the $y$ direction, i.e.,
optimal triangles of a fixed area tend to be infinitely long in one direction. 
This situation also holds when $m>2$. Indeed, we see in the second part
in the proof of Proposition \ref{vanishprop} that if $\pi$ is a nontrivial polynomial such that
$K_{m,p}(\pi)=0$, then $\e_n$ must tend to $0$ as $n\to +\infty$. This shows that 
$T_n=\phi_n(T)$ tends to be infinitely flat in the direction $Ue_y$ with $e_y=(0,1)$.
However, $K_{m,p}(\pi)=0$ does not any longer mean that
$\pi$ is a polynomial of one variable.
\end{remark}

\noindent
Our next result shows that the function $K_{m,p}$ is homogeneous
and obeys an invariance property with respect to linear change of variables.
\begin{prop}
\label{propinvarK}
For all $\pi \in \H_m$, $\lambda\in \R$, and $\phi \in \cL(\R^2)$,
\begin{eqnarray}
\label{homogK}
K_{m,p}(\lambda \pi) &=& |\lambda| K_{m,p}(\pi),\\
\label{invK}
K_{m,p}(\pi\circ \phi) &=& |\det\phi|^{m/2} K_{m,p}(\pi).
\end{eqnarray}
\end{prop}

\proof
The homogeneity property \iref{homogK} is a direct consequence of the definitions of $K_{m,p}$. In order to prove the invariance property \iref{invK}, we assume in a first part that $\det \phi\neq 0$, and
we define $\ti T:=\frac{\phi(T)}{\sqrt{|\det \phi|}}$
and $\ti \pi (z):=\pi(\sqrt{|\det \phi|}z)= |\det \phi|^{m/2} \pi(z)$. 

We now remark that the local interpolant $\interp^{m-1}_T $ commutes 
with linear change of variables in the sense that, when $\phi$ is an invertible linear transform,
\be
\interp^{m-1}_T (v\circ \phi)=(\interp^{m-1}_{\phi(T)} v) \circ \phi
\label{intercommut}
\ee
for all continuous functions $v$ and triangles $T$. Using this commutation formula we obtain 
\begin{eqnarray*}
e_{m,T}(\pi\circ \phi)_p &=& |\det \phi|^{-1/p} e_{m,\phi(T)}(\pi)_p\\
&=& e_{m,\ti T}(\ti \pi)_p\\
&=& |\det \phi|^{m/2} e_{m,\ti T} (\pi)_p.
\end{eqnarray*} 

Since the map $T\mapsto \ti T$ is a bijection of the set of triangles onto itself, leaving the area invariant, we obtain the relation \iref{invK} when $\phi$ is invertible.
When $\det \phi = 0$, the polynomial $\pi \circ \phi$ can be written $(\alpha x+\beta y)^m$
so that $K_{m,p}(\pi\circ \phi) = 0$ by Proposition \ref{vanishprop}.
\sq

The functions $K_{m,p}$ are not necessarily continuous, but the 
following properties will be sufficient for our purposes.

\begin{prop}
\label{propsemicont}
The function $K_{m,p}$ is upper semi-continuous in general and continuous if $m=2$ or $m$ is odd.
Moreover, the following property holds:
\be
\text{If } \pi_n\ra \pi \text{ and } K_{m,p}(\pi_n) \ra 0, \text{ then } K_{m,p}(\pi) = 0.
\label{contzero}
\ee
\end{prop}

\proof
The upper semi-continuity property comes from the fact that the infimum of a family of upper
semi-continuous functions is an upper semi-continuous function. 
We apply this fact to the functions $\pi \mapsto e_{m,T}(\pi)_p$ indexed by triangles which are obviously continuous.

For any polynomial $\pi \in \H_2$, $\pi = a x^2+2 b xy+ c y^2$, we define $\det \pi = ac-b^2$. It will be shown in \S \ref{secOptAniso4} that $K_{2,p}(\pi)= \sigma_p\sqrt{|\det \pi|}$, where $\sigma_p$ only depends on the sign of $\det \pi$. This clearly implies the continuity of $K_{2,p}$.
We next turn to the proof of the continuity of $K_{m,p}$ for odd $m$. 
Consider a polynomial $\pi\in \H_m$. If $K_{m,p}(\pi)=0$, then the upper semi-continuity of $K_{m,p}$, combined with its nonnegativity, implies that it is continuous at $\pi$. 
Otherwise, assume that $K_{m,p}(\pi)>0$. Consider a sequence $\pi_n\in \H_m$ converging to $\pi$ and a sequence $\phi_n$ of linear transformations satisfying $\det \phi_n = 1$, and such that
$$
\lim_{n\to +\infty} e_{\phi_n(\TEq)}(\pi_n) =  \liminf_{\pi^*\to \pi} K_{m,p}(\pi^*)
:= \lim_{r\to 0} \inf_{\|\pi^*- \pi\|\leq r} K_{m,p}(\pi^*).
$$
If the sequence $\phi_n$ admits a converging subsequence $\phi_{n_k}\to \phi$, it follows that 
$$
K_{m,p}(\pi) \leq e_{\phi(\TEq)}(\pi) = \lim_{k\to+\infty} e_{\phi_{n_k}(\TEq)}(\pi_{n_k}) = \liminf_{\pi^*\to \pi} K_{m,p}(\pi^*).
$$
This asserts that $K_{m,p}$ is lower semi-continuous at $\pi$, and therefore continuous at $\pi$ since we already know that $K_{m,p}$ is upper semi-continuous.

If $\phi_n$ does not admit any converging subsequence, then we invoke the Singular Value Decomposition (SVD)   
$\phi_n = U_n\circ D_n \circ V_n$, where $U_n,V_n\in \cO_2$ and $D_n = \diag(\ve_n,\frac 1 {\ve_n})$, where $0<\e_n\leq 1$. (Here and below, we use the shorthand $\diag(a,b)$ to denote the diagonal matrix with entries $a$ and $b$.) 
The compactness of $\cO_2$ implies that $U_n$ admits a converging subsequence $U_{n_k}\to U$. 
In particular, $\pi_{n_k}\circ U_{n_k}$ converges to $\pi\circ U$. Therefore,
denoting by $a_{i,n}$ the coefficient of $x^i y^{m-i}$ in $\pi_n \circ U_n$, the subsequence $a_{i,n_k}$ converges to the coefficient $a_i$ of $x^i y^{m-i}$ in $\pi\circ U$.
Observe also that $\e_n\to 0$; otherwise, some converging 
subsequence could be extracted from $\phi_n$.
Since $e_{\phi_n(\TEq)}(\pi_n) = e_{\TEq}(\pi_n\circ\phi_n)$, the sequence of polynomials $\pi_n\circ\phi_n$ is uniformly bounded and so is the sequence $\pi_n\circ U_n\circ D_n$. Therefore,
the sequences $(a_{i,n} \ve_n^{2i-m})_{n\geq 0}$ are uniformly bounded.
It follows that $a_i = 0$ when $i< \frac m 2$. Since $m$ is odd, this implies that $\pi\circ U(x,y) = x^{s_m} \tilde\pi(x,y)$, and Proposition \ref{vanishprop} implies that $K_{m,p}(\pi) = 0$, which contradicts the hypothesis $K_{m,p}(\pi)>0$.

Finally, we prove property \iref{contzero}.
The assumption $K_{m,p}(\pi_n)\to 0$ is equivalent to the existence of 
a sequence $T_n=\phi_n(\TEq)$ with $\det \phi_n = 1$ such that
$e_{m,T_n}(\pi_n)_p \to 0$. Reasoning in a similar way as in
the proof of Proposition \ref{vanishprop}, we first obtain that $\pi_n\circ \phi_n\to 0$,
and we then invoke the SVD decomposition of $\phi_n$
to build a converging sequence of orthogonal matrices 
$U_n\ra U$ and a sequence $0<\ve_n \leq 1$ such that if $a_{i,n}$ is the coefficient of 
$x^i y^{m-i}$ in $\pi_n\circ U_n$, we have $a_{i,n} \ve_n^{2i-m} \to 0$. 
When $i<s_m$, it follows that $a_{i,n}\to 0$, and therefore $\pi\circ U (x,y)= x^{s_m} \hat \pi(x,y)$. The result follows from Proposition \ref{vanishprop}.
\sq

We finally make a connection between the shape function
and the approach developed in \cite{C3}.
For all $\pi\in \H_m$, we denote by $\Lambda_\pi$ the level set of $|\pi|$ for the value $1$:
\be
\Lambda_\pi = \{(x,y)\in \R^2,\ |\pi (x,y)|\leq 1\}.
\label{lambdapi}
\ee
We now define 
\be
\label{defKE}
K^\cE_m(\pi) = \left(\sup_{E\in \cE,\ E\subset \Lambda_\pi} |E|/ \cPi \right)^{-m/2},
\ee
where the supremum is taken over the set $\cE$ of all ellipses centered at $0$. We use a bold font to denote the numerical constant $\cPi = 3.14159...$
The optimization among ellipses defining $K^\cE_m$ can be rephrased as an optimization on the cone 
 $S_2^+$ of $2\times 2$ positive symmetric matrices.
\be
\label{eqOptimEll}
K_m^\cE(\pi) = \inf \{(\det M)^{\frac m 4}\sep  M\in S_2^+ \text{ and } \forall z \in \R^2, \<Mz,z\> \geq |\pi(z)|^{2/m}\}.
\ee
The minimizing ellipse $E^*$ is then given by $\{ \<Mz,z\>\leq 1\}$.
The optimization problem described in \iref{eqOptimEll} 
is quadratic in dimension $2$ and subject to (infinitely many) linear constraints.
This apparent simplicity is counterbalanced by the fact that it is nonconvex. In
particular, it does not have unique solutions and may also have no solution.  We construct in Chapter 6, \S \ref{subsecWellPosed}, a variant $K^{(\alpha)}$ of $K^\cE$ which is defined by by a well posed optimization problem, for which the optimal matrix is unique and continuously depends on the parameter $\pi$.

\begin{prop}
\label{propequivEllTri}
On $\H_m$, one has the equivalence
$$
cK^\cE_m \leq K_{m,p} \leq C K^\cE_m,
$$
with constant $0<c\leq C$ independent of $p$.
\end{prop}

\proof 
We consider a fixed triangle $T_*$ of unit area, for instance an equilateral triangle, for each exponent $1\leq p \leq \infty$ on $\H_m$ we define a norm $\|\cdot \|_p$ on $\H_m$ as follows 
$$
\|\pi\|_p := \|\pi - \interp_{T_*}^{m-1} \pi\|_{L^p(T_*)},
$$
in such way that 
\be
\label{eqKMPPhi}
K_{m,p}(\pi) = \inf_{|\det \phi|=1} \|\pi\circ \phi\|_p,
\ee
where the infimum is taken among the collection of linear changes of variables $\phi$ satisfying $|\det \phi|=1$.
Using \iref{eqOptimEll} and the homogeneity of $\pi$ we obtain a similar expression for the shape function $K_m^\cE$ based on ellipses
\be
\label{eqKEPhi}
K_m^\cE(\pi) = \inf_{|\det \phi|=1} \|\pi\circ \phi\|.
\ee
Indeed if $\phi\in \GL_d$ satisfies $|\det \phi| = 1$, then $M := \|\pi\circ \phi\|^\frac 2 m (\phi^{-1})^\trans \phi^{-1}$ satisfies $(\det M)^\frac m 4 = \| \pi \circ \phi\|$ and $\<Mz, z\> \geq |\pi(z)|^\frac 2 m$. Conversely to any $M\in S_2^+$ we associate $\phi = M^{-\frac 1 2} (\det M)^\frac 1 4$ which satisfies $\det \phi = 1$.

Since the vector space $\H_m$ has finite dimension there exists $C\geq c>0$ such that for any $\pi \in \H_m$ and any exponent $1\leq p\leq \infty$ one has 
$$
c \|\pi\|\leq \|\pi\|_1 \leq \|\pi\|_p \leq \|\pi\|_\infty \leq C \|\pi\|.
$$
Combining this with \iref{eqKMPPhi} and \iref{eqKEPhi} we obtain the announced result.
\sq

\begin{remark}
Since $K_{m,p}$ and $K_m^\cE$ are equivalent, they vanish on the
same set, and therefore Proposition \ref{vanishprop} is also valid for $K_m^\cE$.
It also easy to see that $K_m^\cE$ satisfies the homogeneity and invariance
properties stated for $K_{m,p}$ in \iref{homogK} and \iref{invK}, as well as
the continuity properties stated in Proposition \ref{propsemicont}.
\end{remark}

\begin{remark}
The continuity of the functions $K_{m,p}$ and $K_m^\cE$ can be established
when $m$ is odd or equal to $2$, as shown by Proposition \ref{propsemicont}, but seems to fail otherwise. In particular, direct computation shows that $K_4^\cE(x^2y^2-\ve y^4)$ is independent of 
$\ve>0$ and strictly smaller than $K_4^\cE(x^2 y^2)$. 
Therefore, $K_4^\cE$ is upper semi-continuous but discontinuous at the point $x^2 y^2\in \H_4$.
\end{remark}

\section{Optimal estimates}
\label{secOptAniso3}
This section is devoted to the proofs of our main theorems, starting with the lower estimate of Theorem \ref{optitheorem}, and continuing with the upper estimates involved in both Theorem \ref{maintheorem} and \ref{optitheorem}.

Throughout this section, for the sake of notational simplicity, 
we fix the parameters $m$ and $p$ and use the shorthand
$$
K=K_{m,p}\;\; {\rm and}\;\; e_T(\pi) = e_{m,T}(\pi)_p.
$$
For each point $z\in \Omega$ we define
$$
\pi_z := \frac{d^mf_z}{m!} \in \H_m,
$$
where $f\in C^m(\overline \Omega)$ is the function in the statement of the theorems. We 
denote by
$$
\omega(r):=\sup_{\|z-z'\|\leq r} \|\pi_z-\pi_{z'}\|
$$
the modulus of continuity of $z\mapsto \pi_z$
with the norm $\|\cdot\|$ defined by \iref{normdef}. 
Note that $\omega(r)\to 0$ as $r\to 0$.

\subsection{Lower estimate}

In this proof we will use an estimate from below of the local interpolation error.

\begin{prop}
\label{errorfPz}
Assume that $1\leq p<\infty$. There exists a constant $C>0$, depending on $f$ and $\Omega$,
such that for all triangles $T\subset \Omega$ and $z\in T$,
\be
e_T(f)^p \geq K^p(\pi_z) |T|^{\frac{mp} 2+1} - C(\diam T)^{mp} |T| \omega(\diam T).
\label{localbelow}
\ee
\end{prop}

\proof
Denoting by $\mu_z\in \P_m$ the Taylor development of $f$ at the point $z$ up to degree $m$, we obtain
$$
f(z+u) - \mu_z(z+u) = m \int_{t=0}^1 (\pi_{z+tu}(u)-\pi_z(u))(1-t)^{m-1} dt,
$$
and therefore
$$
\|f-\mu_z\|_{L^\infty(T)} \leq C_0\diam(T)^m \omega(\diam(T)),
$$
where $C_0$ is a fixed constant.
By construction $\pi_z$ is the homogenous part of $\mu_z$ of degree $m$, and therefore $\mu_z-\pi_z\in \P_{m-1}$. It follows that for any triangle $T$, we have 
\be
\label{eqMuPi}
\mu_z - \interp^{m-1}_T \mu_z = \pi_z -\interp^{m-1}_T  \pi_z.
\ee
We therefore obtain 
\begin{eqnarray*}
|e_T(f)- e_T(\pi_z)| &\leq & \|(f-\interp^{m-1}_T f) - (\pi_z-\interp^{m-1}_T \pi_z)\|_{L^p(T)} \\
&\leq & |T|^{1/p} \|(f-\interp^{m-1}_T f) - (\mu_z-\interp^{m-1}_T \mu_z)\|_{L^\infty(T)} \\
&= & |T|^{1/p} \|(I-\interp^{m-1}_T )(f-\mu_z)\|_{L^\infty(T)} \\
&\leq & C_1|T|^{1/p} \|f-\mu_z\|_{L^\infty(T)} \\
&\leq & C_0C_1 |T|^{1/p} \diam(T)^m \omega(\diam(T)),
\end{eqnarray*}
where $C_1$ is the norm of the operator $I-\interp^{m-1}_T : C^0(T) \to C^0(T)$ in $L^\infty(T)$ norm
which is independent of $T$. 

From \iref{shapeA} we know that $e_T(\pi_z) \geq |T|^{\frac m 2+\frac 1 p} K(\pi_z)$, and therefore
$$
e_T(f) \geq K(\pi_z) |T|^{\frac m 2+\frac 1 p} - C_0C_1 |T|^{1/p} \diam(T)^m \omega(\diam(T)).
$$
We now remark that for all $p\in[1,\infty)$ the function $r\mapsto r^p$ is convex, and therefore if $a,b,c$ are positive numbers, and $a\geq b-c$, then $a^p \geq \max\{0,b-c\}^p \geq b^p -p c b^{p-1}$. Applying this to our last inequality we obtain
$$
e_T(f)^p \geq K^p(\pi_z) |T|^{\frac {mp} 2+1} -
pC_0C_1 (K(\pi_z))^{p-1} |T|^{(p-1)(\frac m 2+\frac 1 p)+\frac 1 p} \diam(T)^m \omega(\diam T).
$$
Since $|T|^{(p-1)(\frac m 2+\frac 1 p)+\frac 1 p}=|T|^{(p-1)\frac m 2} |T| \leq (\diam T)^{m(p-1)}|T|$, this leads to
$$
e_T(f)^p \geq K^p(\pi_z) |T|^{\frac{mp} 2+1} - C(\diam T)^{mp} |T| \omega(\diam T),
$$
where $C :=pC_0C_1 (\sup_{z\in\Omega}K(\pi_z))^{p-1}$. 
\sq

We now turn to the proof of the lower estimate in Theorem \ref{optitheorem} in the case
where $p<\infty$.
Consider a sequence $\seqT$ of triangulations which is admissible in the sense of equation 
\iref{admissibilitycond}. Therefore, there exists a constant $C_A$ such that 
$$
\diam T\leq C_A N^{-1/2},\; N\geq N_0,\; T\in \cT_N.
$$
For $T\in \cT_N$, we combine this estimate with \iref{localbelow}, which gives
$$
e_T(f)^p \geq K^p(\pi_z) |T|^{\frac{mp} 2+1} - (C_AN^{-1/2})^{mp} |T| C\omega(C_AN^{-1/2}).
$$
Averaging over $T$, we obtain
$$
e_T(f)^p \geq \int_T K^p(\pi_z) |T|^{\frac{mp} 2} dz - |T|  N^{-\frac{mp} 2}C_A^{mp}
C\omega(C_A N^{-1/2}).
$$
Summing on all $T\in \cT_N$, and denoting by $T_z^N$ the triangle in $\cT_N$ containing the point $z\in\Omega$, we obtain the estimate
\be
\label{lowerint}
e_{\cT_N}(f)^p \geq \int_\Omega K(\pi_z) |T_z^N|^{\frac{mp} 2} dz - N^{-\frac{mp} 2}\e(N),
\ee
where $\e(N):=|\Omega|C_A^{mp}C\omega(C_AN^{-1/2})\to 0$ as $N\to +\infty$.
The function $z\mapsto |T_z^N|$ is linked with the number of triangles in the following way:
$$
\int_\Omega \frac {dz} {|T_z^N|} = \sum_{T\in \cT_N} \int_T \frac 1{|T|} \leq N.
$$
On the other hand, with $\frac 1 \tau = \frac m 2 +\frac 1 p$, we have by H\"older's inequality,
\be
\int_\Omega K^\tau(\pi_z) dz\leq \left(\int_\Omega K^p(\pi_z) |T_z^N|^{\frac{mp} 2} dz\right)^{\tau/p}
\left(\int_\Omega \frac 1 {|T_z^N|}dz\right)^{1-\tau/p}.
\label{holder}
\ee
Combining the above, we obtain a lower bound for the integral term in \iref{lowerint}
which is independent of $\cT_N$:
$$
\int_\Omega K^p(\pi_z) |T_z^N|^{\frac{mp} 2} dz\geq \left(\int_\Omega K^\tau(\pi_z)dz\right)^{p/\tau} N^{-m p/2}.
$$
Inserting this lower bound into \iref{lowerint} we obtain
$$
e_{\cT_N}(f)^p \geq \left[\left(\int_\Omega K^\tau(\pi_z)dz\right)^{p/\tau} -\e(N) \right] N^{-m p/2}.
$$
This allows us to conclude
\be
\liminf_{N\to +\infty} N^{\frac m 2} e_{\cT_N}(f) \geq \left(\int_\Omega K^\tau(\pi_z)dz\right)^{\frac 1 \tau},
\label{lowerest}
\ee
which is the desired estimate.
\nl
\nl
The case $p=\infty$ follows the same ideas. 
Adapting Proposition \ref{errorfPz}, one proves that
$$
e_T(f) \geq K(\pi_z) |T|^{\frac m 2} - C(\diam T)^{m} \omega(\diam T),
$$
and therefore
\be
e_{\cT_N}(f) \geq \left\|K(\pi_z) |T_z^N|^{\frac m 2} \right\|_{L^\infty(\Omega)} - N^{-\frac m 2} \e(N), 
\label{lowerint2}
\ee
where $\e(N):= C_A^m C\omega(C_A N^{-\frac 1 2})\to 0$ as $N\to +\infty$.
The Hölder inequality now reads:
$$
\int_\Omega K(\pi_z)^{\frac 2 m} dz \leq \left\|K(\pi_z)^{\frac 2 m} |T_z^N|\right\|_{L^\infty(\Omega)} \left\|\frac 1 {|T_z^N|}\right\|_{L^1(\Omega)} ,
$$
equivalently,
$$
\left\|K(\pi_z) |T_z^N|^{\frac m 2} \right\|_{L^\infty(\Omega)} \geq \left(\int_\Omega K(\pi_z)^{\frac 2 m} dz\right)^{\frac m 2} N^{- \frac m 2}.
$$
Combining this with \iref{lowerint2}, we obtain to the desired estimate \iref{lowerest}
with $p=\infty$ and $\tau=\frac 2 m$.

\begin{remark}
This proof reveals the two principles which characterize the optimal triangulations.
Indeed, the lower estimate \iref{lowerest}
becomes an equality only when both inequalities in \iref{localbelow} and \iref{holder}
are equality. The first condition - equality in \iref{localbelow} - is met when each 
triangle $T$ has an optimal shape, in the sense that $e_T(\pi_z)=K(\pi_z)|T|^{\frac m 2+\frac 1 p}$
for some $z\in T$. The second condition - equality in \iref{holder} - is met when the ratio
between 
$K^p(\pi_z)|T_z^N|^{\frac{mp} 2}$ and $|T_z^N|^{-1}$ is constant, or equivalently
$K(\pi_z)|T|^{\frac m 2+\frac 1 p}$ is independent of the triangle $T$. 
Combined with the first condition, this means that the error $e_T(f)^p$ is equidistributed
over the triangles, up to the perturbation by  $(\diam T)^{mp} |T| \omega(\diam T)$ which
becomes negligible as $N$ grows.
\end{remark}

\subsection{Upper estimate}
\label{secOptAniso3p2}
We first remark that the upper estimate in Theorem \ref{optitheorem}. 
implies the upper estimate in Theorem \ref{maintheorem}
by a sub-sequence extraction argument: if the upper estimate in Theorem \ref{optitheorem}
holds, then for all $\ve>0$ there exists a sequence
$(\cT_N^\ve)_{N>N_0}$, $\#(\cT_N) \leq N$, such that 
$$
\limsup_{N\to +\infty}\( N^{\frac m 2} e_{\cT_N^\ve}(f)\) \leq \left\| K\left(\frac {d^mf}{m!}\right)\right\|_{L^\tau}+ \ve,
$$
with $\frac 1 \tau=\frac 1 p+\frac m 2$. We then choose a sequence $(\ve_N)_{N\geq N_0}$ such that 
$$
N^{\frac m 2} e_{\cT_N^{\ve_N}}(f) \leq \left\| K\left(\frac {d^mf}{m!}\right)\right\|_{L^\tau}+ 2\ve_N
$$
for all $N \geq N_0$, and $\ve_N \to 0$ as $N \to \infty$.
Defining $\cT_N:=\cT_N^{\ve_N}$ we thus obtain 
$$
\limsup_{N\to +\infty}\( N^{\frac m 2} e_{\cT_N}(f)\) \leq \left\| K\left(\frac {d^mf}{m!}\right)\right\|_{L^\tau}
$$
which concludes the proof of Theorem \ref{maintheorem}.
We are thus left with proving the upper estimate in Theorem \ref{optitheorem}. 
We begin by fixing a (large) number $M>0$. We shall take the limit $M\to \infty$ in the very last step of our proof. We define 
$$
\mT_M = \{T \text{ triangle}\sep |T|=1,\ \bary(T)=0 \text{ and } \diam(T)\leq M\},
$$
the set of triangles centered at the origin, of unit area and diameter smaller than $M$.
This set is compact with respect to the Hausdorff distance.
This allows us to define a ``tempered'' version of $K = K_{m,p}$ that we denote by $K_M$: for all $\pi \in \H_m$
$$
K_M(\pi) := \inf_{T\in \mT_M} e_T(\pi).
$$
Since $\mT_M$ is compact, the above infimum
is attained on at least a triangle, that we denote by $T_M(\pi)$.
Note that the map $\pi\mapsto T_M(\pi)$ need not be continuous.
It is clear that $K_M(\pi)$ decreases as $M$ grows.
Note also that the restriction to triangles $T$ centered at $0$ is artificial, since
the error is invariant by translation as noticed in 
\iref{transinv}. Therefore, $K_M(\pi)$ converges to $K(\pi)$ as $M\to+\infty$.
Since $\mT_M$ is compact, the map $\pi \mapsto \max_{T\in\mT_M}e_T(\pi)$ 
defines a norm on $\H_m$, and is therefore bounded by $C_M\|\pi\|$ for some $C_M>0$.
One easily sees that the functions $\pi \mapsto e_T(\pi)$ are uniformly $C_M$-Lipschitz
for all $T\in \mT_M$, and so is $K_M$.

We now use this new function $K_M$ to obtain a local upper error estimate that is closely related to the local lower estimate in Proposition \ref{errorfPz}.
\begin{prop}
For $z_1\in\Omega$, let $T$ be a triangle which is obtained from
$T_M(\pi_{z_1})$ by rescaling and translation ( $T=tT_M(\pi_{z_1})+z_0$).
Then for any $z_2\in T$, 
\be
\label{localabove}
e_T(f) \leq \(K_M(\pi_{z_2}) + B_M\omega(\max\{|z_1-z_2|,\diam(T)\})\) |T|^{\frac m 2+\frac 1 p},
\ee
where $B_M>0$ is a constant which depends on $M$.
\end{prop}

\proof
For all $z_1,z_2\in \Omega$, we have
\begin{eqnarray*}
\label{eTMPxPy}
e_{T_M(\pi_{z_1})}(\pi_{z_2}) &\leq & e_{T_M(\pi_{z_1})} (\pi_{z_1}) + C_M \|\pi_{z_1}-\pi_{z_2}\| \\
& = & K_M(\pi_{z_1}) + C_M \|\pi_{z_1}-\pi_{z_2}\|,\\
&\leq & K_M(\pi_{z_2}) + 2 C_M\|\pi_{z_1}-\pi_{z_2}\|,\\
&\leq & K_M(\pi_{z_2}) + 2C_M\omega (|z_1-z_2|).
\end{eqnarray*}
Therefore, if $T$ is of the form $T=tT_M(\pi_{z_1})+z_0$, we obtain by a change of variable that
$$
e_T(\pi_{z_2})\leq \(K_M(\pi_{z_2}) + 2C_M\omega(|z_1-z_2|)\) |T|^{\frac m 2+\frac 1 p}.
$$
Let $\mu_z\in \P_m$ be the Taylor polynomial of $f$ at the point $z$ up to degree $m$.
Using \iref{eqMuPi} we obtain 
\begin{eqnarray*}
e_T(f) &\leq &  e_T(\mu_{z_2})+ e_T(f-\mu_{z_2})\\
&=& e_T(\pi_{z_2})+ e_T(f-\mu_{z_2})\\
&\leq &  \(K_M(\pi_{z_2}) + 2C_M\omega(|z_1-z_2|)\) |T|^{\frac m 2+\frac 1 p} + e_T(f-\mu_{z_2}).\\
\end{eqnarray*}
By the same argument as in the proof of Proposition \ref{errorfPz}, we derive that
$$
e_T(f-\mu_{z_2})\leq C|T|^{\frac 1 p} \diam(T)^m \omega(\diam T),
$$
and thus 
$$
e_T(f) \leq \(K_M(\pi_{z_2}) + 2C_M\omega(|z_1-z_2|)\) |T|^{\frac m 2+\frac 1 p} +
C|T|^{\frac 1 p} \diam(T)^m \omega(\diam T).
$$
Since $T$ is the scaled version of a triangle in $\mT_M$, it obeys $\diam(T)^2 \leq M^2 |T|$. Therefore, 
$$
e_T(f) \leq \left(K_M(\pi_{z_2}) + (2C_M+CM^m)\omega(\max\{|z_1-z_2|,\diam (T)\})\right) |T|^{\frac m 2 +\frac 1 p},
$$
which is the desired inequality with $B_M:=2C_M+CM^m$.
\sq

For some $r>0$ to be specified later, we now choose an arbitrary triangular mesh $\cR$ of $\Omega$ satisfying 
$$
r \geq \sup_{R\in \cR} \diam(R).
$$
Our strategy to build a triangulation that satisfies the optimal upper estimate
is to use the triangles $R$ as {\it macro-elements} in the sense that each of them
will be tiled by a locally optimal uniform triangulation.  This strategy was already used in \cite{BBLS}.

For all $R\in \cR$ we consider the triangle
$$
T_R:=(K_M(\pi_{b_R})+2B_M\omega(r))^{-\frac \tau 2} T_M(\pi_{b_R}),
$$
which is a scaled version of $T_M(\pi_{b_R})$
where $b_R$ is the barycenter of $R$. We use this triangle
to build a periodic tiling $\cP_R$
of the plane: there exists $c\in \R^2$ such that 
$T_R\cup T'_R$ forms a parallelogram
of side vectors $a$ and $b$, with $T'_R=c-T_R$. We
then define 
\be
\cP_R:=\{T_R+ma+nb\sep m,n\in\Z^2\} \cup  \{T'_R +ma+nb\sep m,n\in\Z^2\}.
\label{defTiling}
\ee
Observe that for all $\pi\in \H_m$ and all triangles $T,T'$ such that $T'=-T$, one has $e_T(\pi) = e_{T'}(\pi)$, since $\pi$ is either an even polynomial when $m$ is an even integer, or an odd polynomial when $m$ is odd. Since we already know that $e_T(\pi)$ is invariant by translation of $T$,
we find that the local error $e_T(\pi)$ is constant on all $T\in \cP_R$.

\begin{figure}
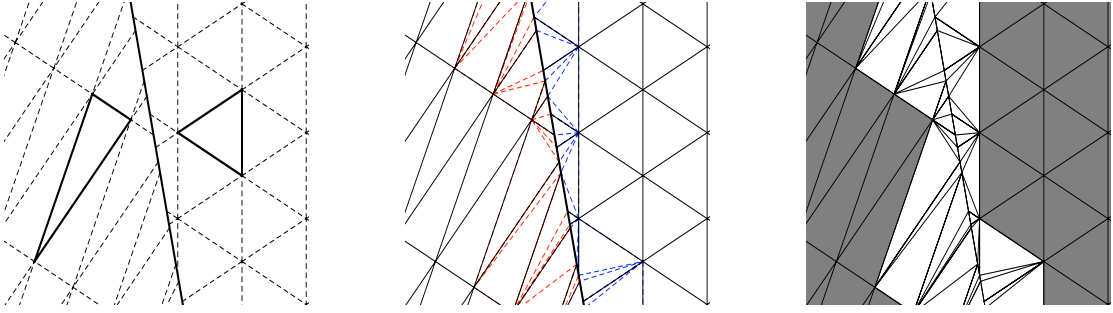

	\centering
		\includegraphics[width=4cm,height=4cm]{\pathPic/PaperOptAniso/NonConformTiling.pdf}
		\hspace{1cm}
		\includegraphics[width=4cm,height=4cm]{\pathPic/PaperOptAniso/ConformTiling.pdf}
		\hspace{1cm}
		\includegraphics[width=4cm,height=4cm]{\pathPic/PaperOptAniso/TriangleClasses.pdf}
	\caption{\label{fig1OptAniso}a. An edge (thick) of the macro-triangulation $\cR$ separating to uniformly paved regions ($T_R$ is thick, $\cP_R$ is dashed). b. Additional edges (dashed) are added near the interface in order to preserve conformity. c. The sets of triangles $\cTr$ (gray) and $\cTb$ (white)}
\end{figure}

We now define as follows a family of triangulations $\cT_s$ of the domain 
$\Omega$, for $s>0$. For every $R\in \cR$, we consider the elements $T\cap R$
for $T\in s\cP_R$, where $s\cP_R$ denotes the triangulation $\cP_R$
scaled by the factor $s$. Clearly, $\{T\cap R,\; T\in s\cP_R,\; R\in\cR\text{ and } \interior (T \cap R) \neq \emptyset\}$ constitute
a partition of $\Omega$, up to a set of zero Lebesgue measure. In this partition, we distinguish the interior 
elements 
$$
\cTr:=\{T\in s\cP_R \sep T\subset \interior(R) \; , R\in \cR\},
$$
which define pieces of a conforming triangulation, and the 
boundary elements $T\cap R$ for $T\in s\cP_R$ such that $\interior(T)\cap \partial R\neq \emptyset$.
These last elements might not be triangular, nor conformal with the 
elements on the other side. Note that for $s>0$ small enough, each $R\in\cR$
contains at least one triangle in $\cTr$, and therefore the boundary elements
constitute a layer around the edges of $\cR$. In order to obtain a conforming
triangulation, we proceed as follows: for each boundary element $T\cap R$,
we consider the points on its boundary which are either its vertices or
those of a neighboring element. We then build the Delaunay triangulation
of these points, which is a triangulation of $T\cap R$ since it is a convex set.
We denote by $\cTb$ the set of all triangles obtained by this procedure,
which is illustrated in Figure \ref{fig1OptAniso}.

Our conforming triangulation is given by
$$
\cT_s=\cTr \cup \cTb.
$$
As $s\ra 0$, clearly,
$$
\#(\cTb)\leq  C_{\rm bd}s^{-1} \text{ and } \sum_{T\in \cTb} |T| \leq C_{\rm bd} s
$$
for some constant $C_{\rm bd}$ which depends on the macro-triangulation $\cR$. 
We do not need to estimate $C_{\rm bd}$ since $\cR$ is fixed and the contribution due to $C_{\rm bd}$ in the following estimates is negligible as $s\to 0$.
We therefore obtain that the number of triangles
in $\cTb$ is dominated by the number of triangles in $\cTr$. More precisely, we
have the equivalence
\be
\label{cardcTs}
\#( \cT_s )\sim \#(\cTr) \sim  \sum_{R\in \cR} \frac{ |R|}{s^2|T_R|}
= s^{-2} \sum_{R\in \cR} |R| (K_M(\pi_{b_R})+2B_M \omega(r))^\tau,
\ee
in the sense that the ratio between the above quantities tends to $1$ as $s\to 0$.
The right-hand side in \iref{cardcTs} can be estimated through an integral:
\begin{eqnarray*}
s^2 \#(\cTr) &\leq & \sum_{R\in \cR} |R| (K_M(\pi_{b_R})+2B_M \omega(r))^\tau\\
& = & \sum_{R\in \cR} \int_R (K_M(\pi_{b_R})+2B_M \omega(r))^\tau dz\\
& \leq & \sum_{R\in \cR} \int_R (K_M(\pi_z)+C_M\|\pi_z-\pi_{b_R}\|+2B_M \omega(r))^\tau dz\\
&\leq & \int_\Omega (K_M(\pi_z)+(2B_M+C_M) \omega(r))^\tau dz.
\end{eqnarray*}
Therefore, since $C_M\leq B_M$,
\be
\label{cardcTsUpper}
\#(\cT_s)\leq s^{-2}\left(\int_\Omega (K_M(\pi_z)+3B_M \omega(r))^\tau dz +C_{\rm bd} s \right).
\ee
Observe that the construction of $\cT_s$ gives a bound on the diameter of its elements
$$
\sup_{T\in \cT_s}\diam(T)\leq s C_a, \;\; C_a:=\max_{R\in\cR}{\rm diam}(T_R).
$$
Combining this with \iref{cardcTs}, we obtain that   
$$
\sup_{T\in \cT_s}\diam(T) \leq C_A \#(\cT_s)^{-1/2} \text{ for all } s>0,
$$
which is analogous to the admissibility condition 
\iref{admissibilitycond}. 
\nl
\nl
We now estimate the global interpolation error 
$\|f-\interp^{m-1}_{\cT_s}f\|_{L^p}:=(\sum_{T\in\cT_s}e_T(f)^p)^{\frac 1 p}$,
assuming first that $1\leq p<\infty$. We first estimate
the contribution of $\cTb$, which will eventually be negligible. 
Denoting by $\nu_z\in \P_{m-1}$ the Taylor polynomial of $f$ up to degree $m-1$ at $z$, we remark that
$$
\|f-\interp^{m-1}_T  f\|_{L^\infty(T)} = \|(I-\interp^{m-1}_T  )(f-\nu_{b_T})\|_{L^\infty(T)}
 \leq C_1\|f-\nu_{b_T}\|_{L^\infty(T)} \leq C_0C_1 \diam(T)^m,
$$
where $C_1$ is the norm of $I-\interp^{m-1}_T $
in $L^\infty(T)$ which is independent of $T$,
and $C_0$ only depends on the $L^\infty$ norm
of $d^mf$. Remarking that $e_T(f) = \|f-\interp^{m-1}_T  f\|_{L^p(T)} \leq |T|^{\frac 1 p} \|f-\interp^{m-1}_T  f\|_{L^\infty(T)}$, we obtain an upper bound for the contribution of $\cTb$ to the error:
\begin{eqnarray*}
\sum_{T\in \cTb}  e_{T}(f)^p &\leq& C_0^pC_1^p\sum_{T\in \cTb} |T|  \diam(T)^{mp}\\
&\leq & C_0^pC_1^p \left(\sum_{T\in \cTb} |T|\right) \sup_{T\in \cTb} \diam(T)^{mp} \\
&\leq & C_0^pC_1^p C_{\rm bd} s \sup_{T\in \cTb} \diam(T)^{mp}\\
&\leq & C^*_{\rm bd}s^{mp+1},
\end{eqnarray*}
with $C^*_{\rm bd}=C_0^pC_1^p C_a^{mp}C_{\rm bd}$.
We next turn to the the contribution of $\cTr$ to the error.
If $T\in \cTr$, $T\subset R\in \cR$, we consider any point $z_1 = z\in T$ and define $z_2 = b_R$
the barycenter of $R$. With such choices, the estimate \iref{localabove} reads 
$$
e_{T}(f) \leq \left(K_M(\pi_z) + B_M \omega(\max\{r,C_A s\})\right) |T|^{\frac m 2+\frac 1 p}.
$$
We now assume that $s$ is chosen small enough such that 
$C_A s\leq r$. Geometrically, this condition ensures that the ``micro-triangles'' constituting $\cT_s$ actually have a smaller diameter than the ``macro-triangles'' constituting $\cR$. This implies
\be
\label{estimabove}
e_{T}(f)^p \leq \left(K_M(\pi_z) + B_M \omega(r)\right)^p |T|^{\frac {mp} 2+1}.
\ee
Given a triangle $T\in \cTr$, $T\subset R \in \cR$, and a point $z\in T$, one has 
\begin{eqnarray*}
|T|  &= &s^2 \left(K_M(\pi_{b_R})+2B_M \omega(r)\right)^{-\tau}\\
& \leq & s^2 \left(K_M(\pi_{z})-C_M\|\pi_z-\pi_{b_R}\|+2B_M \omega(r)\right)^{-\tau}\\
&\leq& s^2\left(K_M(\pi_z)+(2B_M-C_M) \omega(r)\right)^{-\tau}.
\end{eqnarray*}
Observing that $B_M\geq C_M$, and that $p-\tau\frac{mp} 2 = \tau$, we insert the above inequality into the estimate \iref{estimabove}, which yields
$$
e_{T}(f)^p \leq s^{mp} \left(K_M(\pi_z) + B_M \omega(r)\right)^\tau |T|.
$$
Averaging on $z\in T$, we obtain
$$
e_{T}(f)^p \leq s^{mp} \int_T  \left(K_M(\pi_z) + B_M \omega(r)\right)^\tau dz.
$$
Adding up contributions from all triangles in $\cT_s$, we find
$$
 e_{\cT_s}(f)^p =  \sum_{T\in \cTr} e_{T}(f)^p + \sum_{T\in \cTb} e_{T}(f)^p \leq s^{mp} \int_\Omega \left(K_M(\pi_z) + B_M \omega(r)\right)^\tau dz + C^*_{\rm bd} s^{mp+1}.
$$
Combining this with the estimate \iref{cardcTsUpper} we obtain
$$
e_{\cT_s} \#(\cT_s)^{\frac m 2} \leq \left(\int_\Omega \left(K_M(\pi_z) + B_M \omega(r)\right)^\tau dz+C^*_{bd} s\right)^{\frac 1 p} 
\left(\int_\Omega \left(K_M(\pi_z)+3B_M \omega(r)\right)^\tau dz +C_{\rm bd} s \right)^{\frac m 2},
$$
and therefore, since $\frac 1 \tau = \frac m 2+\frac 1 p$,
$$
\limsup_{s\to 0}\(\#(\cT_s)^{\frac m 2} e_{\cT_s}\) \leq \left(\int_\Omega \left(K_M(\pi_z)+3B_M \omega(r)\right)^\tau dz \right)^{\frac 1 \tau}.
$$
It is now time to observe that for any fixed $M$,
$$
\lim_{r\ra 0} \int_\Omega \left(K_M(\pi_z) + 3B_M \omega(r)\right)^\tau dz = \int_\Omega K_M^\tau(\pi_z) dz,
$$
and that, since $K_M(\pi)$ converges decreasingly to $K(\pi) := K_{m,p}(\pi)$ for any $\pi \in \H_m$
$$
\lim_{M\ra +\infty} \int_\Omega K_M^\tau(\pi_z)dz = \int_\Omega K^\tau(\pi_z) dz.
$$
Therefore, for all $\e>0$, we can choose $M$ sufficiently large and $r$
sufficiently small, such that
$$
\limsup_{s\to 0}\(\#(\cT_s)^{\frac m 2} e_{\cT_s}\)\leq \(\int_\Omega K^\tau(\pi_z) dz\)^{\frac 1 \tau}+\e.
$$
This gives us the previously mentionned statement of Theorem \ref{optitheorem}, by defining
$$
s_N:=\min\{s>0 \sep \#(T_s)\leq N\},
$$
and by setting $\cT_N=\cT_{s_N}$.
\nl
\nl
The adaptation of the above proof in the case $p=\infty$ is not straightforward
due to the fact that the contribution to the error of $\cTb$ is no longer negligible
with respect to the contribution of $\cTr$. For this reason, one needs to modify
the construction of $\cTb$. Here, we provide a simple construction
but for which the resulting triangulation $\cT_s$ is nonconforming, as we do not know how to produce a satisfying conforming triangulation.

More precisely, we define $\cTr$ in a similar way as for $p<\infty$,
and add to the construction of $\cTb$ a post-processing step
in which each triangle is split into $4^j$ similar triangles according
to the midpoint rule. Here we take for $j$ the smallest integer
which is larger than $-\frac {\log s}{4\log 2}$. With such an additional splitting,
we thus have
$$
\max_{T\in \cTb}{\rm diam}(T)\leq s^{\frac 1 4}\max_{R\in\cR}{\diam}(sT_R)
=C_as^{1+\frac 1 4}.
$$
The contribution of $\cTb$ to the $L^\infty$ interpolation error is bounded by 
$$
e_{\cTb}(f) \leq C_0C_1\max_{T\in \cTb}{\rm diam}(T)^m \leq C^*_{\rm bd}s^{\frac {5m} 4},
$$
with $C^*_{\rm bd}:=C_0C_1C_a^m$. We also have
$$
\#(\cTb)\leq C_{\rm bd} s^{-3/2},
$$
which remains negligible compared to $s^{-2}$. We therefore obtain
\be
\label{cardcTsInf}
\#(\cT_s)\leq s^{-2}\left(\int_\Omega (K_M(\pi_z)+3B_M \omega(r))^{\frac 2 m} dz +C_{\rm bd} s^{1/2} \right).
\ee
Moreover, if $T\in \cTr$ and $T\subset R\in \cR$, we have according to the estimate \iref{localabove},
$$
e_{T}(f) \leq \left(K_M(\pi_{b_R}) + B_M \omega(\max\{r,C_A s\})\right) |T|^{\frac{m} 2}.
$$
By construction $|T| = s^2 (K_M(\pi_{b_R})+2B_M \omega(r))^{-2/m}$. This implies $e_{T}(f) \leq s^m$ when $C_A s\leq r$. Therefore,
$$
e_{\cT_s}(f) = \max\{e_\cTr,e_\cTb\} \leq s^m\max\{1, C^*_{\rm bd} s^{\frac m 4}\}.
$$
Combining this estimate with \iref{cardcTsInf} yields
$$
\limsup_{s\to 0}\(\#(\cT_s)^{\frac m 2} e_{\cT_s}\) \leq \left(\int_\Omega \left(K_M(\pi_z)+3B_M \omega(r)\right)^{\frac  2 m} dz \right)^{\frac  m 2},
$$
and we conclude the proof in a similar way as for $p<\infty$.

\section{The shape function and the optimal metric for linear and quadratic elements}
\label{secOptAniso4}
This section is devoted to linear ($m=2$) and quadratic ($m=3$) elements,
which are the most commonly used in practice.
In these two cases, we are able to derive an exact expression for $K_{m,p}(\pi)$ in terms
of the coefficients of $\pi$. Our analysis also gives us access to the distorted metric
which characterizes the optimal mesh.
While the results concerning linear elements
have strong similarities with those of \cite{BBLS}, those
concerning quadratic elements are to our knowledge
the first of this kind, although \cite{C1} analyzes a similar setting. 

\subsection{Exact expression of the shape function}

In order to give the exact expression of $K_{m,p}$, we define the determinant
of a homogeneous quadratic polynomial by
$$
\det (ax^2+2bxy+cy^2)=ac-b^2,
$$
and the discriminant of a homogeneous cubic polynomial by
$$
\disc(a x^3+ b x^2 y+ c x y^2+ d y^3) = b^2 c^2 - 4 a c^3 - 4 b^3 d + 18 a b c d - 27 a^2 d^2.
$$
The functions $\det$ on $\H_2$ and $\disc$ on $\H_3$ are homogeneous in the sense that
\be
\label{homogDisc}
\det(\lambda \pi) = \lambda^2\det \pi,\quad \disc(\lambda \pi) = \lambda^4  \disc \pi.
\ee
Moreover, it is well known that they obey an invariance property with respect to linear changes of coordinates $\phi$:
\be
\label{invDisc}
\det(\pi\circ \phi) = (\det \phi)^2 \det \pi,\quad \disc( \pi\circ \phi) =  (\det \phi)^6 \disc \pi.
\ee
Our main result relates $K_{m,p}$ to these quantities.
\begin{theorem}
\label{equal23}
We have for all $\pi\in\H_2$,
$$
K_{2,p}(\pi) = \sigma_p(\det\pi) \sqrt{|\det \pi|},
$$
and for all $\pi\in \H_3$,
$$
K_{3,p}(\pi) = \sigma^*_p(\disc \pi) \sqrt[4]{|\disc \pi|},
$$
where $\sigma_p(t)$ and $\sigma_p^*(t)$ are constants that only depend on the sign of $t$.
\end{theorem}

The proof of Theorem \ref{equal23} relies on the possibility of
mapping an arbitrary polynomial $\pi\in \H_2$ such that 
${\rm det}(\pi)\neq 0$ or $\pi\in\H_3$ such that $\disc(\pi)\neq 0$ onto 
two fixed polynomials $\pi_-$ or $\pi_+$ by a linear change of variable
and a sign change.

In the case of $\H_2$, it is well known
that we can choose $\pi_-=x^2-y^2$ and $\pi_+=x^2+y^2$.
More precisely, to  all $\pi\in H_2$, we associate a
symmetric matrix $Q_\pi$ such that $\pi(z)=\<Q_\pi z,z\>$. This
matrix can be diagonalized according to 
$$
Q_\pi=U^\trans 
\left(
\begin{array}{cc}
\lambda_1 & 0\\
0 & \lambda_2
\end{array}
\right)
U,\quad U\in \cO_2, \; \lambda_1,\lambda_2\in\R.
$$
Then, defining the linear transform
$$
\phi_\pi:=U^\trans 
\left(
\begin{array}{cc}
|\lambda_1|^{-\frac 1 2} & 0\\
0 & |\lambda_2|^{-\frac 1 2}
\end{array}
\right)
$$
and $\lambda_\pi={\rm sign}(\lambda_1)\in \{-1,1\}$, 
it is readily seen that 
$$
\lambda_\pi \pi\circ \phi_\pi =
\left\{
\begin{array}{cl}
x^2+y^2 &\text{ if } \det \pi>0,\\
x^2-y^2 &\text{ if } \det \pi<0.
\end{array}
\right.
$$

In the case of $\H_3$, a similar result holds, as shown by the following lemma: 

\begin{lemma}
\label{lemmaChgVar}
Let $\pi \in \H_3$ be such that $\disc \pi \neq 0$. There exists a linear transform
$\phi_\pi$ such that
\be
\label{eqChgVar}
\pi\circ \phi_\pi =
\left\{
\begin{array}{cl}
x(x^2-3y^2) &\text{ if } \disc \pi>0,\\
x(x^2+3y^2) &\text{ if } \disc \pi<0.
\end{array}
\right.
\ee
\end{lemma}

\proof
Let us first assume that $\pi$ is not divisible by $y$ so that it can be factorized
as 
$$
\pi = \lambda (x-r_1y) (x-r_2 y) (x-r_3 y),
$$
with $\lambda\in\R$ and $r_i\in\C$. If $\disc \pi>0$, then the $r_i$ are real
and we may assume $r_1<r_2<r_3$. Then defining
$$
\phi_\pi = \lambda (2\disc \pi)^{-1/3} 
\left(\begin{array}{cc} 
r_1 (r_2 + r_3) -2 r_2 r_3 & (r_2-r_3) r_1 \sqrt 3\\
2 r_1-(r_2+r_3) & (r_2-r_3) \sqrt 3,
\end{array}\right).
$$
an elementary computation shows that  $\pi\circ \phi_\pi = x(x^2 - 3 y^2)$.
If $\disc \pi < 0$, then we may assume that $r_1$ is real and $r_2$ and $r_3$
are complex conjugates with ${\rm Im}(r_2)>0$. Then defining
$$
\phi_\pi = \lambda (2\disc \pi)^{-1/3} 
\left(\begin{array}{cc} 
r_1 (r_2 + r_3) -2 r_2 r_3 & \mi(r_2-r_3) r_1 \sqrt 3\\
2 r_1-(r_2+r_3) & \mi(r_2-r_3) \sqrt 3 
\end{array}\right),
$$
an elementary computation shows that  $\pi\circ \phi_\pi = x(x^2 +3 y^2)$.
Moreover, it is easily checked that $\phi_\pi$ has real entries
and is therefore a change of variable in $\R^2$.

In the case where $\pi$ is divisible by $y$, there exists a rotation $U\in \cO_2$
such that $\ti \pi:=\pi\circ U$ is not divisible by $y$. By the invariance property
\iref{invDisc} we know that $\disc \pi=\disc \ti \pi$. Thus, we reach the same conclusion
with the choice $\phi_\pi :=U\circ \phi_{\ti \pi}$.
\sq

\noindent
{\bf Proof of Theorem \ref{equal23}} \, For all $\pi\in \H_2$ such that $\det\pi\neq 0$
and for all change of variable $\phi$ and $\lambda\neq 0$,
we may combine the properties of the determinant in \iref{homogDisc} and \iref{invDisc} 
with those of the shape function established in
Proposition \ref{propinvarK}. This gives us
$$
\frac {K_{2,p}(\pi)} {\sqrt{|\det\pi|}}
=\frac {K_{2,p}(\lambda \pi\circ \phi)} {\sqrt{|\det(\lambda \pi\circ\phi)|}}.
$$
Applying this with $\phi=\phi_\pi$ and $\lambda=\lambda_\pi$, we therefore obtain
$$
K_{2,p}(\pi) =
 \sqrt{|\det \pi|} 
\left\{\begin{array}{cc}
K_{2,p}(x^2+y^2) & \text{ if } \det \pi>0,\\
K_{2,p}(x^2-y^2) & \text{ if } \det \pi<0.
\end{array}\right.
$$
This gives the desired result with $\sigma_p(t)=K_{2,p}(x^2+y^2)$ for $t>0$
and $\sigma_p(t)=K_{2,p}(x^2-y^2)$ for $t<0$. In the case
where $\det\pi=0$, $\pi$ is of the
form $\pi(x,y) = \lambda (\alpha x +\beta y)^2$, and we conclude by
Proposition \ref{vanishprop} that $K_{2,p}(\pi)=0$.

For all $\pi\in \H_3$ such that $\disc\pi\neq 0$, a similar reasoning yields
$$
K_{3,p}(\pi) = \sqrt[4]{|\disc \pi|} 108^{-\frac{1} 4}
\left\{\begin{array}{cc}
K_{3,p}(x(x^2-3y^2)) & \text{ if } \disc \pi>0,\\
K_{3,p}(x(x^2+3y^2)) & \text{ if } \disc \pi<0,
\end{array}\right.
$$
where the constant $108$ comes from the fact that
$\disc(x(x^2-3y^2)) = - \disc(x(x^2-3y^2)) = 108$. This gives the
desired result with $\sigma_p^*(t)=108^{-\frac{1} 4}K_{3,p}(x(x^2-3y^2))$ for $t>0$
and $\sigma_p^*(t)=108^{-\frac{1} 4}K_{3,p}(x(x^2+3y^2))$ for $t<0$.
In the case where $\disc\pi=0$, $\pi$ is of the
form $\pi(x,y) =  (\alpha x +\beta y)^2(\gamma x+\delta y)$, and we conclude by
Proposition \ref{vanishprop} that $K_{3,p}(\pi)=0$.\sq 

\begin{remark}
We do not know any simple analytical expression for the constants involved in 
$\sigma_p$ and $\sigma^*_p$, but these can be found by numerical optimization. These constants are known for some special values of $p$ in the case $m=2$, see for example \cite{BBLS}.
\end{remark}

\subsection{Optimal metrics}
\label{secOptAniso4p2}
Practical mesh generation techniques such as in \cite{Shew,Bois,Peyre,Bamg,Inria} 
are based on the data of a Riemannian metric, by which we mean
a field $h$ of symmetric definite positive matrices
$$
x \in \Omega \mapsto h(x) \in S_2^+.
$$
Typically, the mesh generator takes the metric $h$ as an input and hopefully
returns a triangulation $\cT_h$ adapted to it in the sense that all triangles are close
to equilateral of unit side length with respect to this metric. 
Recently, it has been rigorously proved in \cite{Shew2, Bois}, see also Chapter 5, that some algorithms produce bidimensional meshes obeying these constraints under certain conditions. This must be contrasted with algorithms based on heuristics, such as \cite{Bamg} in two dimensions and \cite{Inria} in three dimensions, which have been available for some time and offer good performance \cite{A} but no theoretical guaranties. See \cite{FG} for a review of these mesh generation techniques. 

For a given function $f$ to be approximated, the field of metrics given 
as input should be such that the 
local errors are equidistributed and the aspect ratios are optimal for the 
generated triangulation. Assuming that the error is measured in $X=L^p$
and that we are using finite elements of degree $m-1$, we can construct
this metric as follows, provided that some estimate of $\pi_z = \frac{d^mf(z)}{m!}$ is available
for all points $z\in \Omega$.
An ellipse $E_z$ such that $|E_z|$ is equal
or close to
\be
\sup_{E\in \cE, E\subset \Lambda_{\pi_z}} |E|
\label{maxellips}
\ee
is computed, where $\Lambda_{\pi_z}$ is defined as in \iref{lambdapi}.
We denote by $h_{\pi_z}\in S_2^+$ the associated symmetric definite
positive matrix such that
$$
E_z=\left\{(x,y) \sep (x,y)^T h_{\pi_z} (x,y) \leq 1\right\}.
$$
Let us notice that the supremum in \iref{maxellips} might not
always be attained or even be finite. This particular case is discussed
in the end of this section. Denoting by $\nu>0$ the desired
order of the $L^p$ error on each triangle, we then define the
metric by rescaling $h_{\pi_z}$ according to
$$
h(z)= \frac 1 {\alpha_z^2} h_{\pi_z}\; \; {\rm where}\;\; 
\alpha_z:=\nu^{\frac {p}{mp+2}} |E_z|^{-\frac {1}{mp+2}}.
$$
With such a rescaling, any triangle $T$ designed by the
mesh generator should be comparable to the ellipse
$z+\alpha_z E_z$ centered around $z$ the barycenter of $T$,
in the sense that
\be
z+c_1\alpha_z E_z \subset T\subset z+c_2\alpha_z E_z
\ee
for two fixed constants $0<2c_1\leq c_2$ independent of $T$ (recall that for any ellipse $E$
there always exists a triangle $T$ such that $E\subset T \subset 2E$).

Such a triangulation heuristically 
fulfills the desired properties of optimal aspect ratio and 
error equidistribution when the level of refinement
is sufficiently small. Indeed, we then have
\begin{eqnarray*}
e_{m,T}(f)_p & \approx & e_{m,T}(\pi_z)_p \\
&=& \|\pi_z - \interp^{m-1}_T \pi_z\|_{L^p(T)},\\
&\sim & |T|^{\frac 1 p}\|\pi_z - \interp^{m-1}_T \pi_z\|_{L^\infty(T)},\\
&\sim & |T|^{\frac 1 p} \|\pi_z\|_{L^\infty(T)},\\
&\sim &  |\alpha_z E_z|^{\frac 1 p} \|\pi_z\|_{L^\infty(\alpha_z E_z)},\\
&=&	 \alpha_z^{m+\frac 2 p} |E_z|^{\frac 1 p} \|\pi_z\|_{L^\infty(E_z)},\\
&= & \nu,
\end{eqnarray*}
where we have used the fact that $\pi_z\in\H_m$ is homogeneous of degree $m$. 

Leaving aside these heuristics on error estimation and mesh generation, 
we focus on the main computational issue in the design of the metric $h(z)$,
namely the solution to the problem \iref{maxellips}: to any given $\pi\in \H_m$,
we want to associate $h_\pi\in S_2^+$ such that
the ellipse $E_\pi$ defined by $h_\pi$ has area equal
or close to $\sup_{E\in \cE, E\subset \Lambda_{\pi}} |E|$.

When $m=2$ the computation of the optimal matrix $h_\pi$ 
can be done by elementary algebraic means. In fact, as will
be recalled below, $h_\pi$ is simply the absolute value
(in the sense of symmetric matrices) of the symmetric
matrix $[\pi]$ associated to the quadratic form $\pi$. These facts
are well known and used in mesh generation algorithms
for $\P_1$ elements.

When $m\geq 3$ no such algebraic derivation of $h_\pi$ from
$\pi$ has been proposed up to now, and current 
approaches instead consist in numerically solving 
the optimization problem \iref{eqOptimEll}, see \cite{C3}. 
Since these computations have to be done extremely frequently 
in the mesh adaptation process, a simpler algebraic procedure 
is highly valuable. In this section, we propose a simple and 
algebraic method in the case $m=3$, corresponding to quadratic elements. 
For purposes of comparison, the results already known in the case $m=2$ are recalled.

\begin{prop}
\label{propEllipseMax}
\begin{enumerate}
\item
Let $\pi\in \H_2$ be such that $\det(\pi)\neq 0$, and consider its associated $2\times 2$ matrix which can be 
written as
$$
[\pi] = U^\trans 
\left(
\begin{array}{cc}
\lambda_1 & 0\\
0 & \lambda_2
\end{array}
\right)
U,\quad U\in \cO_2.
$$
Then an ellipse of maximal volume inscribed in $\Lambda_\pi$
is defined by the matrix
$$
h_\pi = U^\trans 
\left(
\begin{array}{cc}
|\lambda_1| & 0\\
0 & |\lambda_2|
\end{array}
\right)
U.
$$

\item
Let $\pi \in \H_3$ be such that $\disc \pi > 0$, and let $\phi_\pi$ be a matrix satisfying \iref{eqChgVar}. Define
\be
\label{hPiPos}
h_\pi =(\phi_\pi^{-1})^\trans \phi_\pi^{-1}.
\ee

Then $h_\pi$ defines an ellipse of maximal volume inscribed in $\Lambda_\pi$. Moreover, $\det h_\pi =\frac {2^{-2/3}} 3 (\disc \pi)^{\frac 1 3}$. 

\item
Let $\pi\in \H_3$ be such that $\disc \pi < 0$, and $\phi_\pi$ a matrix satisfying \iref{eqChgVar}. Define  
$$
h_\pi = 2^{\frac 1 3} (\phi_\pi^{-1})^\trans \phi_\pi^{-1}.
$$
Then $h_\pi$ defines an ellipse of maximal volume inscribed in $\Lambda_\pi$. Moreover, $\det h_\pi = \frac 1 3 |\disc \pi|^{\frac 1 3}$ .
\end{enumerate}
\end{prop}

\begin{figure}
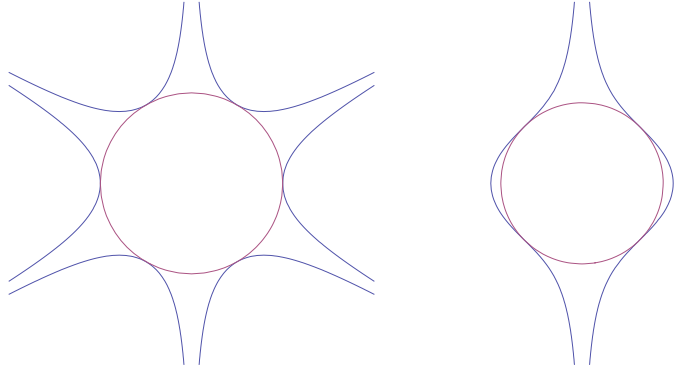

	\centering
		\includegraphics[width=5cm,height=5cm]{\pathPic/PaperOptAniso/EllMaxDiscPos.pdf}
		\includegraphics[width=5cm,height=5cm]{\pathPic/PaperOptAniso/EllMaxDiscNeg.pdf}
	\caption{\label{fig2OptAniso}Maximal ellipses inscribed in $\Lambda_\pi$, $\pi = x(x^2-3y^2)$ or $\pi = x(x^2+3y^2)$.}
\end{figure}

\proof
Clearly, if the matrix $h_\pi$ defines an ellipse of maximal volume in the set $\Lambda_\pi$, then for any linear change of coordinates $\phi$, the matrix $(\phi^{-1})^\trans h_\pi \phi^{-1}$ defines an ellipse of maximal volume in the set $\Lambda_{\pi\circ \phi}$.
When $\pi\in \H_2$, we know that $\lambda_\pi \pi \circ \phi_\pi = x^2+y^2$ when $\det \pi >0$, and $x^2-y^2$ when $\det \pi <0$, where $|\lambda_\pi|=1$. When $\pi \in \H_3$, 
we know from Lemma \ref{lemmaChgVar} that $\pi \circ \phi_\pi = x(x^2-3y^2)$ 
when $\disc \pi>0$, and $x (x^2+3y^2)$ when $\disc \pi <0$. 
Hence it only remains to prove that when $\pi \in \{x^2+y^2,\ x^2-y^2,\ x(x^2-3y^2)\}$, then $h_\pi=\Id$, which means that the disc of radius $1$ is an ellipse of maximal volume inscribed in $\Lambda_\pi$, while when $\pi = x(x^2+3y^2)$ we have $h_\pi=2^{1/3}\Id$.

The case $\pi = x^2+y^2$ is trivial. We next concentrate on the case $\pi = x(x^2+3y^2)$,
the treatment of the two other cases being very similar.
Let $E$ be an ellipse included in $\Lambda_\pi$, $\pi = x(x^2+3y^2)$. 
Analyzing the variations of the function $\pi(\cos\theta,\sin\theta)$, it is not hard to see that we can rotate $E$ into another ellipse $E'$, also satisfying the inclusion $E'\subset \Lambda_\pi$, and whose principal axes are $\{x=0\}$ and $\{y=0\}$.
We therefore only need to consider ellipses of the form $k x^2 + h y^2 \leq 1$. 
For a given value of $h$, we denote by $k(h)$ the minimal value of $k$ for which this ellipse is included in $\Lambda_\pi$. Clearly, the boundary of the ellipse, defined by $k(h) x^2 + h y^2 = 1$, must be tangent to the curve defined by $\pi(x,y) = 1$ at some point $(x,y)$. This translates into the following system of equations:
\be 
\left\{
\begin{array}{ccc}
\pi(x,y) &=& 1,\\
h x^2 + k y^2 &=& 1,\\
k y \partial_x \pi(x,y) - hx \partial_x \pi(x,y) &=& 0.
\end{array}
\right.
\label{tangencyEqn}
\ee
Eliminating the variables $x$ and $y$ from this system, as well as negative
or complex-valued solutions, we find that $k(h) = \frac{4+h^3}{3h^2}$ when $h\in (0,2]$, and $k(h) = k(2) = 1$ when $h\geq 2$.
The minimum of the determinant $h k(h) = \frac 1 3 \left(\frac 4 h + h^2\right)$ is attained for $h=2^{\frac 1 3}$. Observing that $k(2^{\frac 1 3}) = 2^{\frac 1 3}$, we obtain, as previously stated, $h_\pi = 2^{1/3}\Id$ and that the ellipse of largest area included in $\Lambda_\pi$ is the disc
of equation $2^{1/3}(x^2+y^2)\leq 1$, as illustrated on Figure  \ref{fig2OptAniso}.b.

The same reasoning applies to the other cases. For $\pi = x^2-y^2$ we obtain $k(h) = \frac 1 h$, $h\in (0,\infty)$. In this case the determinant $h k(h)$ is independent of $h$, and we simply choose $h=1 = k(1)$.
For $\pi = x(x^2-3y^2)$ we obtain $k(h) = \frac{4-h^3}{3h^2}$ when $h\in (0,1]$ and $k(h) = k(1) = 1$ when $h>1$. The maximal volume is attained when $h=1$, corresponding to the
unit disc, as illustrated on Figure  \ref{fig2OptAniso}.a.
\sq

\begin{remark}
When $\pi\in \H_3$ and $\disc \pi>0$ a surprising simplification happens: the matrix 
{\rm \iref{hPiPos}} has entries which are symmetric functions of the roots $r_1,r_2,r_3$.
Using the relation between the roots and the coefficients of a polynomial, we find the following expression:
if $\pi = a x^3+ 3 b x^2 y+ 3 c x y^2+ d y^3$, then
$$
h_\pi = 2^{-\frac 1 3} 3 (\disc \pi)^{\frac{-1} 3} 
\left(
\begin{array}{cc}
2 (b^2-ac) & bc - ad\\
bc - ad & 2 (c^2-bd)
\end{array}
\right).
$$
This yields a direct expression of the matrix as a function of the coefficients. Unfortunately, there is no such expression when $\disc \pi<0$.
\end{remark}

At first sight, Proposition \ref{propEllipseMax} might seem to be a complete solution to the problem of building an appropriate metric for mesh generation. However, some difficulties arise at points $z\in \Omega$ where $\det \pi_z=0$ or $\disc \pi_z=0$.
If $\pi\in \H_2\sm \{0\}$ and $\det \pi = 0$, then up to a linear change of coordinates, and a change of sign, we can assume that $\pi = x^2$. The minimization problem clearly yields the degenerate matrix
$h_\pi = \diag(1,0)$, the $2\times 2$ diagonal matrix with entries $1$ and $0$.
If $\pi\in \H_3\sm \{0\}$ and $\disc \pi = 0$, then up to a linear change of coordinates either $\pi = x^3$ or $\pi = x^2 y$. In the first case the minimization problem gives again $h_\pi = \diag(1,0)$. In the second case a wilder behavior appears, in the sense that minimizing sequences
for the problem \iref{maxellips} are of the type
$h_\pi = \diag(\ve^{-1},\ve^2)$ with $\ve\to 0$. The minimization process 
therefore gives a matrix which is not only degenerate, but also unbounded.

These degenerate cases appear generically and constitute a problem for
mesh generation since they mean that the adapted triangles are not well defined. 
Current anisotropic mesh generation algorithms for linear elements 
often solve this problem by fixing a small parameter $\delta>0$ and working with the modified matrix
$\tilde h_\pi := h_\pi+\delta\Id$ which cannot degenerate. However, this procedure cannot be extended to quadratic elements, since $h_{x^2 y}$ is both degenerate and unbounded.

\begin{figure}
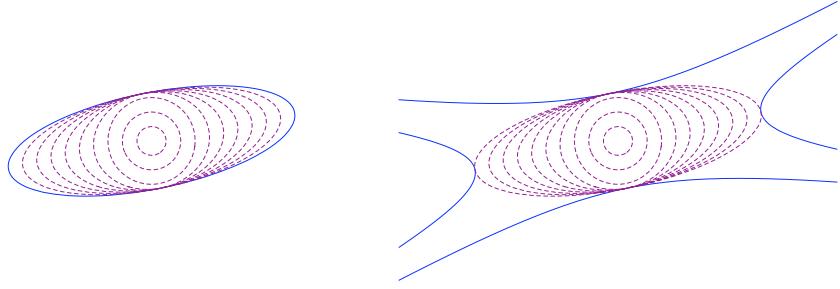

	\centering
		\includegraphics[width=6cm,height=4cm]{\pathPic/PaperOptAniso/AllEllipsesDetPos.pdf}
		\includegraphics[width=6cm,height=4cm]{\pathPic/PaperOptAniso/AllEllipsesDetNeg.pdf}
	\caption{\label{fig3OptAniso}The set $\Lambda_\pi$ (solid) and the ellipses $E_{\pi,\alpha}$ (dashed)
    for various values of $\alpha>0$ when $\pi\in \H_2$.}
\end{figure}

In the theoretical construction of an optimal mesh which
was discussed in \S \ref{secOptAniso3p2}, we  tackled this problem by imposing a bound $M>0$ on the diameter of the triangles. This was the purpose of the modified shape function $K_M(\pi)$ and of the triangle $T_M(\pi)$ of minimal interpolation error among the triangles of diameter smaller than $M$.
We follow a similar idea here, looking for the ellipse of largest area included in $\Lambda_\pi$ with constrained diameter. This provides matrices which are both positive definite and bounded
and which vary continuously with respect to the data $\pi\in \H_3$. A similar construction is presented in \S \ref{subsecWellPosed} of Chapter 6 in arbitrary degree $m$ and dimension $d$, however the expression of the matrix in terms of the polynomial $\pi$ is less explicit in that general context.
The constrained problem, depending on $\alpha>0$, is the following:
\be
\sup \{|E|\sep E\in \cE,\; E\subset \Lambda_\pi \text{ and }\diam E\leq 2\alpha^{-1/2}\},
\label{ellipsconstrained}
\ee
or equivalently,
\be
\label{hConstrained}
\inf \{\det H\sep  H\in S_2^+ \;\;{\rm s.t.}\;\; \<Hz,z\> \geq |\pi(z)|^{2/m}, z\in\R^2,\; \text{ and } H\geq \alpha \Id\}.
\ee
We denote by $E_{\pi,\alpha}$ and $h_{\pi,\alpha}$ the solutions
to \iref{ellipsconstrained} and \iref{hConstrained} respectively.
In the remainder of this section, we show that this solution can also be computed 
by a simple algebraic procedure, avoiding any kind of numerical optimization. In the case
where $\pi\in \H_2$, it can easily be checked that
\be
[h_{\pi,\alpha}] = U^\trans 
\left(
\begin{array}{cc}
\max\{|\lambda_1|,\alpha\} & 0\\
0 & \max\{|\lambda_2|,\alpha\},
\end{array}
\right)
U,
\label{family2}
\ee
as illustrated in Figure  \ref{fig3OptAniso}. 

When $\pi \in \H_3$, the problem is
more technical, and the matrix $h_{\pi,\alpha}$ takes different forms depending on the value of $\alpha$ and the sign of $\disc \pi$.
In order to describe these different regimes, we introduce three real numbers 
$0\leq \beta_\pi\leq \alpha_\pi\leq \mu_\pi$ and a matrix 
$U_\pi\in \cO_2$ which are defined as follows. We first define $\mu_\pi$ by
$$
\mu_\pi^{-1/2}:=\min\{\|z\|\sep |\pi(z)|=1\},
$$ 
the radius of the largest
disc $D_\pi$ inscribed in $\Lambda_\pi$. For $z_\pi$ such that $|\pi(z_\pi)|=1$
and $\|z_\pi\|=\mu_\pi^{-1/2}$, we define $U_\pi$
as the rotation which maps $z_\pi$ to the vector $(\|z_\pi\|,0)$.
We then define $\alpha_\pi$ by
$$
2\alpha_\pi^{-1/2}:=\max \{{\rm diam}(E)\sep E\in \cE\; ;\; D_\pi\subset E\subset\Lambda_\pi\},
$$
the diameter of the largest ellipse inscribed in $\Lambda_\pi$
and containing the disc $D_\pi$. In the case where $\pi$ is
of the form $(ax+by)^3$, this ellipse is infinitely long,
and we set $\alpha_\pi=0$. We finally define $\beta_\pi$ by
$$
2\beta_\pi^{-1/2}:={\rm diam}(E_\pi),
$$ 
where $E_\pi$ is the optimal
ellipse described in Proposition \ref{propEllipseMax}. In the case
where $\disc\pi=0$, the ``optimal ellipse'' is infinitely long,
and we set $\beta_\pi=0$. It is readily seen that 
$0\leq \beta_\pi\leq \alpha_\pi\leq \mu_\pi$.

All these quantities can be algebraically computed
from the coefficients of $\pi$ by solving equations of degree
at most $4$, as well as the other quantities involved in the
description of the optimal $h_{\pi,\alpha}$ and $E_{\pi,\alpha}$
in the following result.

\begin{prop}
\label{family3}
For $\pi\in \H_3$ and $\alpha>0$, the matrix $h_{\pi,\alpha}$ and ellipse $E_{\pi,\alpha}$ are described as follows:
\begin{enumerate}
\item
If $\alpha\geq \mu_\pi$, then $h_{\pi,\alpha}=\alpha\Id$
and $E_{\alpha,\pi}$ is the disc of radius $\alpha^{-1/2}$.
\item
If $\alpha_\pi\leq \alpha\leq \mu_\pi$, then 
\be
\label{hBigAlpha}
h_{\pi,\alpha} = U_\pi^\trans 
\left(
\begin{array}{cc}
\mu_\pi & 0\\
0 & \alpha
\end{array}
\right)
U_\pi,
\ee
and $E_\alpha$ is the ellipse of diameter $2\alpha^{-1/2}$
which is inscribed in $\Lambda_\pi$ and contains $D_\pi$.
It is tangent to $\partial\Lambda_\pi$ at the two points $z_\pi$ and $-z_\pi$.
\item
If $\beta_\pi\leq \alpha\leq \alpha_\pi$ then $E_{\pi,\alpha}$ is tangent
to $\partial\Lambda_\pi$ at four points
and has diameter $2\alpha^{-1/2}$. There are at most 
three such ellipses, and $E_{\pi,\alpha}$ is the one
of largest area. The matrix $h_{\pi,\alpha}$ has a form 
which depends on the sign of $\disc\pi$.
\nl
(i) If $\disc \pi<0$,
then
$$
h_{\pi,\alpha} = (\phi_\pi^{-1})^\trans 
\left(
\begin{array}{cc}
\lambda_\alpha & 0\\
0 & \frac{4+ \lambda_\alpha^3}{3\lambda_\alpha^2}
\end{array}
\right)
\phi_\pi^{-1},
$$
where $\phi_\pi$ is the matrix defined in Proposition {\rm \ref{propEllipseMax}} 
and $\lambda_\alpha$ is determined by $\det(h_{\pi,\alpha}-\alpha\Id)=0$.
\nl
(ii)
If $\disc \pi >0$, then
$$
h_{\pi,\alpha} = (\phi_\pi^{-1})^\trans V^\trans
\left(
\begin{array}{cc}
\lambda_\alpha & 0\\
0 & \frac{4-\lambda_\alpha^3}{3\lambda_\alpha^2}
\end{array}
\right)
V
\phi_\pi^{-1},
$$
where $\phi_\pi$ and $\lambda_\alpha$ are given as in the case $\disc\pi<0$ and
where $V$ is chosen between the three rotations by $0$, $60$ or $120$ degrees
so as to maximize $|E_{\alpha,\pi}|$.\nl
(iii)
If $\disc \pi = 0$ and $\alpha_\pi>0$, then there exists a linear change 
of coordinates $\phi$ such that $\pi\circ\phi =x^2 y$ and we have 
$$
h_{\pi,\alpha} = (\phi^{-1})^\trans
\left(
\begin{array}{cc}
\lambda_\alpha & 0\\
0 & \frac 4 {27\lambda_\alpha^2}
\end{array}
\right)
\phi^{-1},
$$
where $\lambda_\alpha$ is determined by $\det(h_{\pi,\alpha}-\alpha\Id)=0$.
\item
If $\alpha\leq \beta_\pi$, then $h_{\pi,\alpha}=h_\pi$, and
$E_{\pi,\alpha}=E_\pi$ is the solution of the unconstrained problem.
\end{enumerate}
\end{prop}

\proof
See Appendix.
\sq

\begin{figure}
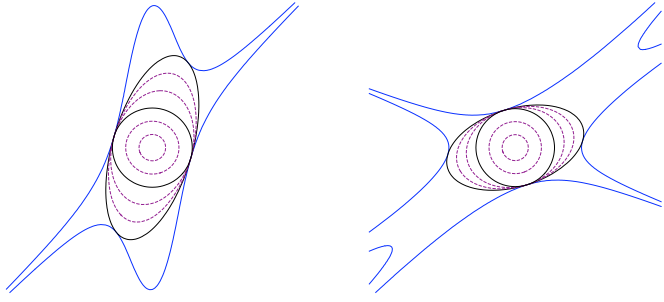

	\centering
		\includegraphics[width=4cm,height=4cm]{\pathPic/PaperOptAniso/EllNeg1.pdf} \hspace{5mm}
		\includegraphics[width=4cm,height=4cm]{\pathPic/PaperOptAniso/EllPos1.pdf}
	\caption{\label{fig4OptAniso}The set $\Lambda_\pi$ (solid), the disc $E_{\pi,\mu_\pi}=D_\pi$ (solid), 
    the ellipse $E_{\pi,\alpha_\pi}$ (solid), and the ellipses $E_{\pi,\alpha}$ (dashed)
    for various values of $\alpha>0$ when $\pi\in \H_3$ and $\alpha \in (\alpha_\pi,\infty)$.
    Left: $\disc\pi <0$. Right: $\disc\pi >0$}
\end{figure}

\begin{figure}
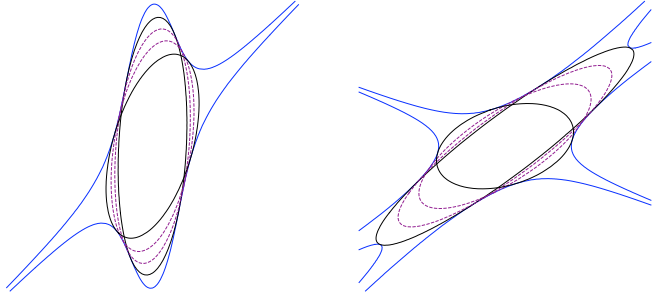

	\centering
		\includegraphics[width=4cm,height=4cm]{\pathPic/PaperOptAniso/EllNeg2.pdf}\hspace{5mm}
		\includegraphics[width=4cm,height=4cm]{\pathPic/PaperOptAniso/EllPos2.pdf}
	\caption{\label{fig5OptAniso}The set $\Lambda_\pi$ (solid), the ellipse $E_{\pi,\alpha_\pi}$ (solid), 
    the ellipse $E_{\pi,\beta_\pi}=E_\pi$ (solid), and the ellipses $E_{\pi,\alpha}$ (dashed)
    for various values of $\alpha>0$ when $\pi\in \H_3$ and $\alpha \in (\beta_\pi,\alpha_\pi)$.
    Left: $\disc\pi <0$. Right: $\disc\pi >0$}
\end{figure}

\noindent
Figure  \ref{fig4OptAniso} illustrates the ellipses $E_{\pi,\alpha}$, $\alpha \in (\alpha_\pi,\infty)$ when $\disc \pi>0$ ( \ref{fig4OptAniso}.a) or $\disc\pi<0$ ( \ref{fig4OptAniso}.b). Figure  \ref{fig5OptAniso} illustrates the ellipses $E_{\pi,\alpha}$, $\alpha \in (\beta_\pi,\alpha_\pi)$ when $\disc \pi>0$ ( \ref{fig5OptAniso}.a) or $\disc\pi<0$ ( \ref{fig5OptAniso}.b). Note that when $\alpha\geq \alpha_\pi$,
the principal axes of $E_{\pi,\alpha}$ are independent of $\alpha$, since $U_\pi$
is a rotation that only depends on $\pi$, while these axes generally
vary when $\beta_\pi\leq \alpha \leq \alpha_\pi$, since the matrix $\phi_\pi$ is not
a rotation.

\begin{remark}
\label{overfitting}
For interpolation by cubic or higher-degree polynomials ($m\geq 4$), an additional difficulty arises that can be summarized as follows: one should be careful not to ``overfit'' the polynomial $\pi$ with the matrix $h_\pi$. An approach based on exactly solving the optimization problem \iref{maxellips} might indeed
lead to a metric $h(z)$ with unjustified strong variations with respect to $z$ 
and/or bad conditioning, and jeopardize the mesh generation process. 
As an example, consider the one-parameter family of polynomials 
$$
\pi_t=x^2y^2+ty^4\in\H_4,\;\; t\in [-1,1].
$$ 
It can be checked that when $t>0$, the supremum
$S_+=\sup_{E\in\cE,E\subset\Lambda_{\pi_t}} |E|$ is finite and independent of $t$, but not attained, and that
any sequence $E_n\subset\Lambda_{\pi_t}$ of ellipses such that $\lim_{n\to\infty} |E_n| = S_+$ becomes infinitely elongated in the $x$ direction as $n\to\infty$. For $t<0$, the supremum
$S_-=\sup_{E\in\cE,E\subset\Lambda_{\pi_t}} |E|$ is independent of $t$ 
and attained for the optimal ellipse of equation 
$|t|^{-1/2}\frac {\sqrt 2-1} 2 x^2+|t|^{1/2}y^2\leq 1$. This ellipse becomes infinitely elongated
in the $y$ direction as $t\to 0$. This example shows the instability of the optimal matrix $h_\pi$ with respect
to small perturbations of $\pi$. However, for all values of $t\in [-1,1]$, 
these extremely elongated ellipses could be discarded in favor, for example,
of the unit disc $D=\{x^2+y^2\leq 1\}$, which obviously satisfies $D\subset \Lambda_{\pi_t}$ and is a near-optimal
choice in the sense that $2|D|= S_+\leq S_-= |D|\sqrt{2(\sqrt 2+1)}$. 
\end{remark}

\section{Polynomial equivalents of the shape function in higher degree}
\label{secOptAniso5}
In degrees $m\geq 4$, we could not find analytical expressions of $K_{m,p}$ or $K_m^\cE$ and do not expect them to exist. However, equivalent quantities with analytical expressions
are available, under the same general form as in Theorem \ref{equal23}: the root of a polynomial in the coefficients of the polynomial $\pi\in\H_m$. This result improves on the analysis of \cite{C2}, where a similar setting is studied.

In the following, we say that a function $\rm \b R$ is a polynomial on $\H_m$ if there exists a polynomial $P$ of $m+1$ variables such that for all $(a_0,\cdots, a_m)\in \R^{m+1}$,
$$
{\rm \b R}\left(\sum_{i=0}^m a_i x^i y^{m-i}\right) := P(a_0,\cdots, a_m) ,
$$
and we define $\deg {\rm \b R} := \deg P$.

The object of this section is to prove the following theorem:
\begin{theorem}
\label{thequiv}
For all degree $m\geq 2$, there exists a polynomial $\Kpol_m$ on $\H_m$, 
and a constant $C_m>0$ such that for all $ \pi \in \H_m$ and all $1\leq p \leq \infty$,
$$
\frac 1 {C_m} \sqrt[r_m]{\Kpol_m(\pi)} \leq K_{m,p}(\pi )\leq C_m \sqrt[r_m]{\Kpol_m(\pi)},
$$
where $r_m = \deg \Kpol_m$.
\end{theorem}

Since for fixed $m$ all functions $K_{m,p}$, $1\leq p\leq \infty$, are equivalent on $\H_m$, there is no need to keep track of the exponent $p$ in this section, and we use below the notation $K_m = K_{m,\infty}$. In this section, the reader should not confuse the functions $K_m$ and $\Kpol_m$, nor the polynomials $Q_d$ and $\b Q_d$ below, which notations are only distinguished by their case.

Theorem \ref{thequiv} is a generalization of Theorem \ref{equal23}, and the polynomial $\Kpol_m$ involved should be seen as a generalization of the determinant on $\H_2$ and of the discriminant on $\H_3$.
Let us immediately stress that the polynomial $\Kpol_m$ is not unique. In particular, we shall propose
two constructions that lead to different $\Kpol_m$ with different degree $r_m$.
Our first construction is simple and intuitive but leads to a polynomial 
of degree $r_m$ that grows quickly with $m$. Our second construction
uses the tools of invariant theory to provide a polynomial
of much smaller degree, which might be more useful in practice.

We first recall that there is a strong connection between the roots of a polynomial in $\H_2$ or $\H_3$ and its determinant or discriminant:
\begin{eqnarray*}
\det\left(\lambda \prod_{1\leq i\leq 2} (x-r_i y)\right) & = & \frac {-1} 4 \lambda^2(r_1-r_2)^2,\\ 
\disc\left(\lambda\prod_{1\leq i\leq 3} (x-r_i y)\right)  & = & \lambda^4 (r_1-r_2)^2(r_2-r_3)^2(r_3-r_1)^2.
\end{eqnarray*}

We now fix an integer $m>3$. Observing that these expressions are a ``cyclic'' product of the squares of differences of roots, we define
$$
\cyc(\lambda,r_1,\cdots,r_m) := \lambda^4 (r_1-r_2)^2\cdots (r_{m-1}-r_m)^2 (r_m-r_1)^2.
$$
Since $m>3$, this quantity is not invariant anymore under reordering of the $r_i$. For any
positive integer $d$, we
introduce the symmetrized version of the $d$-powers of the cyclic product
$$
Q_d(\lambda,r_1,\cdots,r_m) := \sum_{\sigma\in \Sigma_m} \cyc(\lambda,r_{\sigma_1},\cdots,
r_{\sigma_m})^d,
$$
where $\Sigma_m$ is the set of all permutations of $\{1,\cdots, m\}$. 

\begin{prop}
For all $d>0$ there exists a homogeneous polynomial $\b Q_d$ of degree $4 d$ on $\H_m$, with integer coefficients, and such that:
$$
\text{If } \pi = \lambda \prod_{i=1}^m (x-r_i y),
\text{ then } \b Q_d \left(\pi \right) = Q_d(\lambda,r_1,\cdots,r_m).
$$ 
In addition, $\b Q_d$ obeys the invariance property
\be
\label{invprop}
\b Q_d(\pi\circ\phi) = (\det \phi)^{2md} \b Q_d(\pi).
\ee
\end{prop}

\proof
We denote by $\sigma_i$ the elementary symmetric functions in the $r_i$, in such way that
$$
\prod_{i=1}^m (x-r_i y) = x^m - \sigma_1 x^{m-1} y+ \sigma_2 x^{m-2} y^2 -\cdots +(-1)^m \sigma_m y^m.
$$
A well-known theorem of algebra (see, e.g., Chapter IV.6 in \cite{Lang}) asserts that any symmetrical polynomial in the $r_i$ can be reformulated as a polynomial in the $\sigma_i$. Hence for any $d$ there exists a polynomial $\tilde Q_d$ such that
$$
Q_d(1,r_1,\cdots,r_m) = \tilde Q_d(\sigma_1,\cdots,\sigma_m).
$$
In addition, it is known that the total degree of $\tilde Q_d$ is the partial degree of $Q_d$ in the variable $r_1$, in our case $4d$, and that $\tilde Q_d$ has integer coefficients since $Q_d$ has.

Given a polynomial $\pi\in H_m$ not divisible by $y$, we write it under the two equivalent forms:
$$
\pi = a_0 x^m + a_1 x^{m-1}y+\cdots + a_m y^m = \lambda \prod_{i=1}^m(x-r_i y).
$$
Clearly $a_0 = \lambda$ and $\sigma_i = (-1)^i \frac{a_i}{a_0}$. It follows that 
$$
Q_d(\lambda,r_1,\cdots,r_m)= \lambda^{4d}\tilde Q_d(\sigma_1, \cdots,\sigma_m) = a_0^{4d} \tilde Q_d\left(\frac{-a_1}{a_0},\cdots,\frac{(-1)^m a_m}{a_0}\right).
$$
Since $\deg \tilde Q_d = 4d$, the negative powers of $a_0$ due to the denominators are cleared 
by the factor $a_0^{4d}$, and the right-hand side is thus 
a polynomial in the coefficients $a_0,\cdots,a_m$ that we denote by $\b Q_d(\pi)$.

We now prove the invariance of $\b Q_d$ with respect to linear changes of coordinates; this proof is adapted from \cite{Hilbert}. 
By continuity of $\b Q_d$, it suffices to prove this invariance property for pairs $(\pi,\phi)$ such that $\phi$ is an invertible linear change of coordinates and neither $\pi$ or $\pi\circ \phi^{-1}$ is divisible by $y$.

Under this assumption, we observe that if $\pi = \lambda \prod_{i=1}^m(x-r_i y)$ and $\phi = \left(\begin{array}{cc} \alpha &\beta\\ \gamma &\delta \end{array}\right)$, then $\pi \circ \phi^{-1} = \tilde \lambda \prod_{i=1}^m (x-\tilde r_i y)$, where
$$
\tilde \lambda = \lambda (\det\phi)^{-m} \prod_{i=1}^m(\gamma+\delta r_i) \stext{ and } \tilde r_i =\frac {\alpha r_i+\beta }{\gamma r_i+\delta}.
$$
Observing that 
$$
\tilde r_i - \tilde r_j = \frac{\det\phi}{(\gamma r_i+\delta)(\gamma r_j+\delta)} (r_i-r_j),
$$
it follows that 
$$
\cyc(\tilde \lambda,\tilde r_1,\cdots,\tilde r_m) = (\det \phi)^{-2m} \cyc(\lambda,r_1,\cdots,r_m).
$$
The invariance property \iref{invprop} follows readily.
\sq
\nl
We now define $r_m = 2 {\rm lcm}\{\deg \b Q_d\sep  1\leq d \leq m!\}$,
where ${\rm lcm}\{a_1,\cdots,a_k\}$ stands for the lowest common multiple of 
$\{a_1,\cdots,a_k\}$, and we consider the following polynomial on $\H_m$: 

$$
\Kpol_m := \sum_{d=1}^{m!} \b Q_d^{\frac{r_m}{\deg \b Q_d} }.
$$
Clearly, $\Kpol_m$ has degree $r_m$ and obeys the invariance property $\Kpol_m(\pi\circ\phi) = (\det\phi)^{\frac{r_m m} 2} \Kpol_m(\pi)$. 
\begin{lemma}
Let $\pi\in\H_m$. If $\Kpol_m(\pi) = 0$, then $K_m(\pi)=0$.
\end{lemma}

\proof
We assume that $K_m(\pi)\neq 0$ and intend to prove that $\Kpol_m(\pi) \neq 0$. Without loss of generality, we may assume that $y$ does not divide $\pi$, since $K_m(\pi\circ U)= K_m(\pi)$ and $\Kpol_m(\pi\circ U)= \Kpol_m(\pi)$ for any rotation $U$. We thus write $\pi = \lambda\prod_{i=1}^m(x-r_i y)$, where $r_i\in\C$. Since $K_m(\pi)\neq 0$, we know
from Proposition \ref{vanishprop} that there is no group of $\mhalf:= \lfloor \frac m 2 \rfloor +1$ equal roots $r_i$.
 
We now define a permutation $\sigma^* \in \Sigma_m$ such that $r_{\sigma^*(i)}\neq r_{\sigma^*(i+1)}$ for $1\leq i\leq m-1$ and $r_{\sigma^*(m)} \neq r_{\sigma^*(1)}$.
In the case where $m=2m'$ is even and $m'$ of the $r_i$ are equal, any permutation $\sigma^*$ such that $r_{\sigma^*(1)} = r_{\sigma^*(3)} = \cdots = r_{\sigma^*(2m'-1)}$ satisfies this condition. 
In all other cases let us assume that the $r_i$ are sorted by equality:
if $i<j<k$ and $r_i=r_k$, then $r_i=r_j=r_k$. If $m=2m'$ is even, we set $\sigma^*(2i-1) = i$ and $\sigma^*(2i) = m'+i$, $1\leq i\leq m'$. If $m=2m'+1$ is odd we set $\sigma^*(2i) = i$, $1\leq i\leq m'$ and $\sigma^*(2i-1) = m'+i$, $1\leq i\leq m'+1$. For example, $\sigma^* = (4\ 1\ 5\ 2\ 6\ 3\ 7)$
when $m=7$ and  $\sigma^* =(1\ 5\ 2\ 6\ 3\ 7\ 4\ 8)$ when  $m=8$.
With such a construction, we find that $|\sigma^*(i) -\sigma^*(i+1)|\geq m'$
if $m$ is odd and $|\sigma^*(i) -\sigma^*(i+1)|\geq m'-1$ if $m$ is even,
for all $1\leq i\leq m$, where we have set $\sigma^*(m+1):=\sigma^*(1)$.
Hence $\sigma\*$ satisfies the required condition, and therefore $\cyc(\lambda,r_{\sigma^*(1)},\cdots,r_{\sigma^*(m)})\neq 0$.

It is well known that if $k$ complex numbers $\alpha_1,\cdots,\alpha_k\in\C$ are such that $\alpha_1^d+\cdots+\alpha_k^d=0$, for all $1\leq d \leq k$, then $\alpha_1=\cdots =\alpha_k=0$.
Applying this property to the $m!$ complex numbers $\cyc(\lambda,r_{\sigma(1)},\cdots,r_{\sigma(m)})$, $\sigma\in\Sigma_m$ and noticing that the term corresponding to $\sigma^*$ is non zero, we see that there exists $1\leq d\leq m!$ such that $\b Q_d(\pi) = Q_d(\lambda,r_1,\cdots,r_m)\neq 0$. Since 
$\b Q_d$ has real coefficients, the numbers $\b Q_d(\pi)$ are real. Since the exponent $r_m/\deg \b Q_d$ is even, it follows that $\Kpol_m(\pi)>0$, which concludes the proof of this lemma.
\sq
\nl
The following proposition, when applied to the function $\Keq=\sqrt[r_m]{\Kpol_m}$, concludes the proof of Theorem \ref{thequiv}:
\begin{prop}
\label{propequiv}
Let $m\geq 2$, and let $\Keq : \H_m\to \R_+$ be a continuous function obeying the following properties:
\begin{enumerate}
\item \emph{Invariance property} : $\Keq(\pi\circ\phi) = |\det\phi|^{\frac m 2} \Keq(\pi)$. 
\item \emph{Vanishing property} : for all $\pi\in\H_m$, if $\Keq(\pi)=0$, then $K_m(\pi) = 0$.
\end{enumerate}
Then there exists a constant $C>0$ such that $\frac 1 C \Keq\leq K_m\leq C \Keq$ on $\H_m$.
\end{prop}

\proof
We first remark that $\Keq$ is homogeneous in a similar way as 
$K_m$: if $\lambda\geq 0$, then applying the invariance property to 
$\phi = \lambda^{\frac 1 m}\Id$ yields $\Keq(\pi\circ(\lambda^{\frac 1 m}\Id)) = \Keq(\lambda\pi)$ and $|\det\phi|^{\frac m 2} = \lambda$. Hence $\Keq(\lambda\pi) = \lambda \Keq(\pi)$.

Our next remark is that a converse of the vanishing property
holds: if $K_m(\pi) = 0$, then there exists a sequence $\phi_n$ of linear changes of coordinates, $\det\phi_n=1$, such that $\pi\circ\phi_n\to 0$ as $n\to\infty$. Hence $\Keq(\pi) = \Keq(\pi\circ\phi_n) \to \Keq(0)$. Furthermore, $\Keq(0) =  0$ by homogeneity. Hence $\Keq(\pi) = 0$.

We define the set $\NF_m:=\{\pi\in \H_m\sep K_m(\pi) = 0\}$. 
We also define a set $A_m\subset \H_m$ by a property ``opposite'' to the property defining $\NF_m$. A polynomial $\pi\in \H_m$ belongs to $A_m$ if and only if
$$
\|\pi\| \leq \|\pi\circ \phi\|\text{ for all } \phi \text{ such that } \det \phi = 1.
$$
The sets $\NF_m$ and $A_m$ are closed by construction, and clearly $\NF_m\cap A_m = \{0\}$. 
We now denote by $\underline {K_m}$ the lower semi continuous envelope of $K_m$, which is defined by  
$$
\underline{K_m}(\pi) = \lim_{r\to 0} \inf_{\|\pi'-\pi\|\leq r} K_m(\pi')
$$
the lower semi-continuous envelope of $K_m$. If $\underline{K_m}(\pi) = 0$, then there exists a converging sequence $\pi_n\to \pi$ such that $K_m(\pi_n)\to 0$. According to Proposition \ref{propsemicont}, it follows that $K_m(\pi)=0$ and hence $\pi\in\NF_m$. Therefore, the lower semi-continuous function $\underline{K_m}$ and the continuous function $\Keq$ are bounded below by a positive constant on the compact set $\{\pi\in A_m, \|\pi\|=1\}$. Since
in addition $\Keq$ is continuous and $K_m$ is upper semi-continuous, we find that
the constant
$$
C = \sup_{\pi\in A_m, \|\pi\|=1} \max\left\{ \frac{\Keq(\pi)}{\underline{K_m}(\pi)},\; \frac{K_m(\pi)}{\Keq(\pi)} \right\}
$$
is finite. By homogeneity of $K_m$ and $\Keq$, we infer that on $A_m$,
\be
\label{equivAm}
\frac 1 C \Keq \leq \underline{K_m} \leq K_m \leq C \Keq.
\ee
Now, for any $\pi \in \H_m$, we consider $\hat \pi$ of minimal norm in the
closure of the set 
$\{\pi \circ \phi \sep \det \phi = 1\}$. By construction, we have $\hat \pi \in A_m$, and there exists a sequence $\phi_n$, $\det \phi_n=1$ such that $\pi \circ \phi_n \to \hat \pi$ as $n\to \infty$.
If $\hat \pi = 0$, then $K_m(\pi) = \Keq(\pi) = 0$.
Otherwise, we observe that 
$$
\underline{K_m}(\hat \pi)\leq K_m(\pi) \leq K_m(\hat \pi) \text{ and } \Keq(\hat \pi) = \Keq(\pi),
$$
where we used the fact that $\underline{K_m}$, $K_m$, and $\Keq$ are respectively lower semi-continuous, upper semi-continuous, and continuous on $\H_m$.
Combining this with inequality \iref{equivAm} concludes the proof.
\sq

A natural question is to find the polynomial of smallest degree satisfying Theorem \ref{thequiv}.
This leads us to the theory of {\it invariant polynomials} 
introduced by Hilbert \cite{Hilbert} (we also refer to \cite{dix} for a survey on this subject). 
A polynomial $R$ on $\H_m$ is said to be invariant if $\mu = \frac{m\deg R} 2$ is a positive integer and for all $\pi\in \H_m$ and any linear change of coordinates $\phi$, one has
\be
R(\pi\circ\phi) = (\det\phi)^\mu R(\pi).
\label{invarpolR}
\ee
We have seen for instance that $\Kpol_m$ and $\b Q_d$ are ``invariant polynomials'' on $\H_m$.

Nearly all the literature on invariant polynomials is concerned with the case of complex coefficients, both for the polynomials and the changes of variables. It is known in particular \cite{dix} that for all $m\geq 3$, there exists $m-2$ invariant polynomials $R_1,\cdots R_{m-2}$ on $\H_m$ such that for any $\pi$ (complex coefficients are allowed) and any other invariant polynomial $R$ on $\H_m$:
\be
\text{If } R_1(\pi) = \cdots = R_{m-2}(\pi) = 0,\text{ then } R(\pi) = 0.
\label{generators2}
\ee
A list of such polynomials with minimal degree is known explicitly at least when $m\leq 8$, and they have real coefficients. 
Defining $r=2{\rm lcm} (\deg R_i)$ and $\Keq := \sqrt[r]{\sum_{i=1}^{m-2} \b R_i^{\frac{r}{\deg R_i} }}$, we see that $\Keq(\pi)=0$ implies $\Kpol_m(\pi) = 0$, and hence $K_m(\pi)=0$. According to proposition \ref{propequiv}, we have constructed a new, possibly simpler, equivalent of $K_m$.

For example, when $m=2$ the list $(R_i)$ is reduced to the polynomial $\det$, and for $m=3$ to the polynomial $\disc$. 
For $m=4$, given $\pi = ax^4+4 b x^3 y+6c x^2 y^2+4 d x y^3+ey^4$, the list consists of the two polynomials 
$$
I = ae-4bd+3c^2,\quad J=\left|\begin{array}{ccc} a&b&c\\ b&c&d\\ c&d&e\end{array}\right|;
$$
therefore $K_4(\pi)$ is equivalent to the quantity $\sqrt[6]{|I(\pi)|^3+J(\pi)^2}$.
As $m$ increases, these polynomials unfortunately become more and more complicated, and their number $m-2$ obviously increases. According to \cite{dix}, for $m=5$ the list consists of three polynomials of 
degrees $4,8,12$, while for $m=6$ it consists of $4$ polynomials of degrees $2,4,6,10$.

\section{Extension to higher dimension}
\label{secOptAniso6}
The function $K_{m,p}$  can be generalized to higher dimension $d>2$ in the following way.
We denote by $\H_{m,d}$ the set of homogeneous polynomials of degree $m$ in $d$ variables. For all $d$-dimensional simplex $T$, we define the interpolation operator $\interp^{m-1}_T $ acting from $C^0(T)$ onto the space $\P_{m-1,d}$ of polynomials of total degree $m-1$ in $d$ variables. This operator is defined by the conditions $\interp^{m-1}_T v(\gamma) = v(\gamma)$ for all points $\gamma \in T$ with barycentric coordinates in the set $\{0,\frac 1 {m-1},\frac 2 {m-1},\cdots,1\}$. Following Section \S \ref{secOptAniso1p2}, and generalizing definition \iref{shapefunction}, we define the local interpolation error on a simplex, the global interpolation error on a mesh,
and the shape function.

For all $\pi \in \H_{m,d}$,
$$
K_{m,p,d}(\pi) := \inf_{|T|=1} \|\pi -\interp^{m-1}_T \pi\|_p,
$$
where the infimum is taken on all $d$-dimensional simplices $T$ of volume $1$.
The variant $K_m^\cE$ introduced in \iref{defKE} also generalizes in higher dimension and was introduced by Weiming Cao in \cite{C3}. Denoting by $\cE_d$ the set of $d$-dimensional ellipsoids, 
we define 
\begin{eqnarray*}
K_{m,d}^\cE(\pi) &=& \left( \sup_{E\in \cE_d, E\subset \Lambda_\pi} |E|/|B| \right)^{-\frac {m} d} \\
&=& \inf\{(\det M)^{\frac m {2d}} \sep M\in S_d^+ \text{ and } \forall z\in \R^d, \<Mz,z\> \geq |\pi(z)|^{\frac m 2}\},
\end{eqnarray*}
where $\Lambda_\pi := \{z\in \R^d \sep |\pi(z)|\leq 1\}$, $B := \{z\in \R^d \sep |z|\leq 1\}$ is the unit euclidean ball, and where $S_d^+$ denotes the set of symmetric positive definite $d\times d$ matrices. 
Similarly to Proposition \ref{propequivEllTri}, it is not hard to show that the functions $K_{m,p,d}(\pi)$ and $K_{m,d}^\cE(\pi)$ are equivalent: there exist constants $0<c\leq C$ depending only on $m,d$, such that
$$
c K_{m,d}^\cE \leq K_{m,p,d} \leq C K_{m,d}^\cE.
$$
Let $(\cT_N)_{N\geq 0}$ be a sequence of simplicial meshes (triangles if $d=2$, tetrahedrons
if $d=3$, \ldots) of a bounded $d$-dimensional, polygonal open set $\Omega$. Generalizing 
\iref{admissibilitycond}, we say that $(\cT_N)_{N\geq N_0}$ is admissible if $\#(\cT_N) \leq N$ and if there exists a constant $C_A$ satisfying
$$
\sup_{T\in \cT_N} \diam(T) \leq C_A N^{-1/d}.
$$
The lower estimate in Theorem \ref{optitheorem} can be generalized, with straightforward adaptations in the proof. If $f\in C^m(\overline \Omega)$ and $\seqT$ is an admissible sequence of simplicial meshes, then
\be
\label{eqLowerLPD}
\liminf_{N\to\infty} N^{\frac m d} e_{m,\cT_N}(f)_p \geq \left\|K_{m,d,p}\left(\frac{d^m f}{m!}\right)\right\|_{L^\tau(\Omega)},
\ee
where $\frac 1 \tau:= \frac m d+\frac 1 p$.

The upper estimate in Theorem \ref{optitheorem} 
however does not generalize. The reason is that we used 
in its proof a tiling of the plane consisting of translates of a single triangle 
and of its symmetric with respect to the origin. This construction is not possible anymore in higher dimension, for example it is well known that one cannot tile the space $\R^3$, with equilateral tetrahedra. 

The generalization of the second part of Theorem \ref{optitheorem} is therefore the following.
For all $m$ and $d$, there exists a constant $C=C(m,d)>0$ such that for any bounded polygonal open set
$\Omega \subset \R^d$ and $f\in C^m(\overline \Omega)$ the following holds: for all $\e>0$, there exists an admissible sequence $(\cT_N)_{N \geq N_0}$ of simplicial meshes $\Omega$, possibly non conforming, such that  
$$
\limsup_{N\to\infty} N^{\frac m d} e_{m,\cT_N}(f)_p \leq C\left\|K_{m,d,p}\left(\frac{d^m f}{m!}\right)\right\|_{L^\tau(\Omega)}+\e.
$$ 
The ``tightness'' of Theorem \ref{optitheorem} is partially lost due to the constant $C$. This upper bound is not new and can be found in \cite{C3}. 
In the proof of the bidimensional theorem, we define by \iref{defTiling} a tiling $\cP_R$ of the plane made of a triangle $T_R$ and some of its translates and of their symmetric with respect to the origin. In dimension $d$, the tiling $\cP_R$ cannot be constructed by the same procedure. 
The idea of the proof is to first consider a fixed tiling $\cP_0$ of the space constituted of simplices
bounded diameter, and of volume bounded below by a positive constant, as well as a reference equilateral simplex $\TEq$ of volume $1$. We then set $\cP_R = \phi(\cP_0)$, where $\phi$ is a linear change of coordinates such that $T_R=\phi(\TEq)$. This procedure can be applied in any dimension and yields all subsequent estimates ``up to a multiplicative constant,'' which concludes the proof.

Since this upper bound is not tight anymore, and since the functions $K_{m,p,d}$ are all equivalent to $K_{m,d}^\cE$ as $p$ varies (with equivalence constants independent of $p$), there is no real need to keep track of the exponent $p$. We therefore denote by $K_{m,d}$ the function $K_{m,\infty,d}$.

For practical as well as theoretical purposes, it is desirable to have an efficient way to compute the shape function $K_{m,d}$, and an efficient algorithm to produce adapted triangulations.
The case $m=2$, which corresponds to piecewise linear elements, has been extensively studied, see for instance \cite{BBLS,CSX}. In that case there exists constants $0<c<C$, depending only on $d$, such that for all $\pi\in \H_{2,d}$,
$$
c\sqrt[d]{|\det \pi|} \leq K_{2,d}(\pi) \leq C\sqrt[d]{|\det\pi|},
$$
where $\det \pi$ denotes the determinant of the symmetric matrix associated to $\pi$.
Furthermore, similarly to Proposition \ref{propEllipseMax}, the optimal metric for mesh refinement is given by the absolute value of the matrix of second derivatives, see \cite{BBLS,CSX}, which is constructed in a similar way as in dimension $d=2$: with $U$ and $D=\diag(\lambda_1,\cdots,\lambda_d)$
the orthogonal and diagonal matrices such that $[\pi] = U^\trans D U$ and 
with $|D|:=\diag(|\lambda_1|,\cdots,|\lambda_d|)$, we set $h_\pi = U^\trans |D| U$. It can be shown that the matrix $h_\pi$ defines an ellipsoid of maximal volume included
in the set $\Lambda_\pi$. The case $m=2$ can therefore be regarded as solved.

For values $(m,d)$ both larger than $2$, the question of computing the shape function as well as the optimal metric is much more difficult, but we have partial answers, in particular for quadratic elements in dimension $3$. Following \S\ref{secOptAniso5}, we need fundamental results from the theory of invariant polynomials, developed in particular by Hilbert \cite{Hilbert}. In order to apply these
results to our particular setting, we need to introduce a compatibility condition
between the degree $m$ and the dimension $d$.

\begin{definition}
We call the pair of integers $m\geq 2$ and $d\geq 2$ ``compatible'' if and only if the following holds:
For all $\pi \in \H_{m,d}$ such that there exists a sequence $(\phi_n)_{n\geq 0}$ of $d\times d$ matrices with \emph{complex} coefficients, satisfying $\det \phi_n=1$ and $\lim_{n\to\infty} \pi\circ \phi_n = 0$,
there also exists a sequence $\psi_n$ of $d\times d$ matrices with \emph{real} coefficients, satisfying $\det \psi_n=1$ and $\lim_{n\to\infty} \pi\circ \psi_n = 0$.
\end{definition}

Following Hilbert \cite{Hilbert}, we say that
a polynomial $Q$ of degree $r$ defined on $\H_{m,d}$ is 
\emph{invariant} if  $\mu=\frac{mr} d$ is a positive integer and 
if for all $\pi\in\H_{m,d}$ and all linear changes of coordinates $\phi$,
\be
Q(\pi\circ\phi) = (\det \phi)^\mu Q(\pi).
\label{invpolQ}
\ee
This is a generalization of \iref{invarpolR}. We denote by $\I_{m,d}$ the set of invariant 
polynomials on $\H_{m,d}$. It is easy to see that if $\pi\in \H_{m,d}$
is such that $K_{m,d}(\pi)=0$, then $Q(\pi)=0$ for all
$Q\in \I_{m,d}$. Indeed, as seen in the proof of Proposition
\ref{vanishprop}, if $K_{m,d}(\pi)=0$, then there exists a sequence $\phi_n$ such that
$\det \phi_n=1$ and $\pi\circ\phi_n\to 0$. Therefore, \iref{invpolQ}
implies that $Q(\pi)=0$. 
The following lemma shows that
the compatibility condition for the pair $(m,d)$ 
is equivalent to a converse of this property.

\begin{lemma}
\label{lemmacomp}
The pair $(m,d)$ is compatible if and only if for all $\pi \in \H_{m,d}$,
$$
K_{m,d}(\pi) = 0 \text{ if and only if } Q(\pi) = 0 \text{ for all } Q\in\I_{m,d}. 
$$
\end{lemma}

\proof
We first assume that the pair $(m,d)$ is not compatible. Then there exists a polynomial $\pi_0\in\H_{m,d}$ such that there exists a sequence $\phi_n$, $\det \phi_n=1$ of matrices with \emph{complex} coefficients such that $\pi\circ\phi_n \to 0$, but there exists no such sequence with \emph{real} coefficients. This last property indicates that $K_{m,d}(\pi)>0$. On the contrary, let $Q\in \I_{m,d}$ be an invariant polynomial, and set $\mu = \frac{m\deg Q} d$. The identity
$$
Q(\pi_0\circ\phi) =(\det\phi)^\mu Q(\pi_0)
$$
is valid for all $\phi$ with real coefficients and is a \emph{polynomial identity} in the coefficients of $\phi$. Therefore, it remains valid if $\phi$ has complex coefficients. If follows that $Q(\pi_0) = Q(\pi_0\circ\phi_n)$ for all $n$, and therefore $Q(\pi_0)=0$, which concludes the proof in the case where the pair $(m,d)$ is not compatible.

We now consider a compatible pair $(m,d)$.
Following Hilbert \cite{Hilbert}, we say that a polynomial $\pi\in H_{m,d}$ is a {\it null form} if and only if there exists a sequence of matrices $\phi_n$ with \emph{complex} coefficients such that $\det \phi_n=1$
and $\pi\circ \phi_n \to 0$. We denote by $\NF_{m,d}$ the set of such polynomials.
Since the pair $(m,d)$ is compatible, note that $\pi\in\NF_{m,d}$ if and only if there exists 
a sequence $\phi_n$ of matrices with \emph{real} coefficients such that $\det \phi_n=1$ and $\pi\circ\phi_n\to 0$. Hence, we find that
$$
\NF_{m,d}=\{\pi\in\H_{m,d}\sep K_{m,d}(\pi) = 0\}.
$$
Denoting by $\I_{m,d}^\sC$ the set of invariant 
polynomials on $\H_{m,d}$ with {\it complex} coefficients, 
a difficult theorem of \cite{Hilbert} states that 
$$
\NF_{m,d} = \{\pi\in\H_{m,d}\sep Q(\pi) =0 \text{ for all } Q \in \I_{m,d}^\sC\}.
$$
It is not difficult to check that if $Q = Q_1+ i  Q_2$ where $Q_1$ and $Q_2$ have real coefficients,
then \iref{invpolQ} holds for $Q$ if and only if it holds for both $Q_1$ and $Q_2$,
i.e., $Q_1$ and $Q_2$ are also invariant polynomials. 
Hence denoting by $\I_{m,d}$ the set of invariant polynomials on $\H_{m,d}$ with real coefficients, we have obtained that 
$$
\NF_{m,d} = \{\pi\in\H_{m,d}\sep Q(\pi) =0  \text{ for all }  Q \in \I_{m,d}\},
$$
which concludes the proof.
\sq

\begin{theorem}
\label{ThEquivMd}
If the pair $(m,d)$ is compatible, then there exists a polynomial $\Kpol$ on $\H_{m,d}$ (we set $r=\deg \Kpol$) and a constant $C>0$ such that for all $\pi \in \H_{m,d}$,
\be 
\frac 1 C \sqrt[r]{\Kpol(\pi)}\leq K_{m,d}(\pi) \leq C\sqrt[r]{\Kpol(\pi)}.
\label{equivKmd}
\ee
Furthermore any polynomial $\Kpol$ satisfying the above equivalence needs to be an invariant polynomial on $\H_m$.
If the pair $(m,d)$ is not compatible, then there does not exist such a polynomial $\Kpol$.
\end{theorem}

\proof
The proof of the invariance of any polynomial $\Kpol$ satisfying \iref{equivKmd}, and of the nonexistence property when the pair $(m,d)$ is not compatible, is reported in the appendix.
Assume that the pair $(m,d)$ is compatible. We follow a reasoning very similar to \S\ref{secOptAniso5} to prove the equivalence \iref{equivKmd}.

We use the notations of Lemma \ref{lemmacomp} and we consider the set 
$$
\NF_{m,d} = \{\pi\in\H_{m,d}\sep K_{m,d}(\pi) = 0\} = \{\pi\in\H_{m,d}\sep Q(\pi) =0, \; Q \in \I_{m,d}\}. 
$$
The ring of polynomials on a field is known to be Noetherian. 
This implies that there exists a finite family $Q_1,\cdots,Q_s\in \I_{m,d}$ of invariant polynomials on $\H_{m,d}$ such that any invariant polynomial is of the form $\sum P_i Q_i$ where 
$P_i$ are polynomials on $\H_{m,d}$. We therefore obtain
$$
\NF_{m,d} = \{\pi\in\H_{m,d}\sep Q_1(\pi) =\cdots =Q_s(\pi) = 0\},
$$
which is a generalization of \iref{generators2}, however with no clear bound on $s$.

We now fix such a set of polynomials, set $r:= 2 {\rm lcm}_{1\leq i\leq s} \deg Q_i$, and define 
$$
\Kpol = \sum_{i=1}^s \b Q_i^{\frac r {\deg \b Q_i}}\;\;{\rm and}\;\;\Keq := \sqrt[r]{\Kpol}.
$$
Clearly, $\Kpol$ is an invariant polynomial on $\H_{m,d}$, and $\NF_{m,d} = \{\pi\in\H_{m,d}\sep \Kpol(\pi)=0\}$. 
Hence the function $\Keq$ is \emph{continuous} on $\H_{m,d}$, obeys the invariance property $\Keq(\pi\circ\phi) = |\det\phi| \Keq(\pi)$, and for all $\pi\in\H_m$, $\Keq(\pi)=0$ implies $\Kpol(\pi)=0$ and therefore $K_{m,d}(\pi)=0$.
We recognize here the hypotheses of Proposition \ref{propequiv}, except that the dimension $d$ has changed.
Inspection of the proof of Proposition \ref{propequiv} shows that we use only once the fact that $d=2$, when we refer to Proposition \ref{propsemicont} and state that if $(\pi_n) \in\H_m$, $\pi_n\to \pi$ and $K_m(\pi_n)\to 0$, then $K_m(\pi)=0$. 
This property also applies to $K_{m,d}$, when the pair $(m,d)$ is compatible. Assume that $(\pi_n)\in\H_{m,d}$, $\pi_n\to \pi$ and that $K_{m,d}(\pi_n)\to 0$. Then there exists a sequence of linear changes of coordinates $\phi_n$, $\det\phi_n=1$, such that $\pi_n\circ\phi_n\to 0$. Therefore,
$$
\Kpol(\pi) = \lim_{n\to \infty} \Kpol(\pi_n) = \lim_{n\to\infty} \Kpol(\pi_n\circ \phi_n) = 0.
$$
It follows that $\pi\in\NF_{m,d}$, and therefore $K_{m,d}(\pi)=0$. 
Since the rest of the proof of Proposition \ref{propequiv} never uses that $d=2$, this concludes the proof of equivalence \iref{equivKmd}.
\sq

Hence there exists a ``simple'' equivalent of $K_{m,d}$ for all compatible pairs $(m,d)$, 
while equivalents of $K_{m,d}$ for incompatible pairs need to be more sophisticated, 
or at least different from the root of a polynomial. This theorem leaves open several questions. 
The first one is to identify the list of compatible pairs $(m,d)$. It is easily shown that the pairs  $(m,2)$, $m\geq 2$, and $(2,d)$, $d\geq 2$ are compatible, but this does not provide any new results since we already derived equivalents of the shape function in these cases. More interestingly, we show in the next corollary that the pair $(3,3)$ is compatible, which corresponds to approximation by quadratic elements in dimension $3$.
There exist two generators $S$ and $T$ 
of $\I_{3,3}$, whose expressions are given in \cite{Salmon}
and which have respectively degree $4$ and $6$.

\begin{corollary}
$\sqrt[6]{|S|^3+T^2}$ is equivalent to $K_{3,3}$ on $\H_{3,3}$.
\end{corollary}

\proof
The invariants $S$ and $T$ obey the invariance properties $S(\pi\circ\phi) = (\det\phi)^4 S(\pi)$ and $T(\pi\circ\phi) = (\det\phi)^6 T(\pi)$. We intend to show that if $\pi\in\H_{3,3}$ and $S(\pi) = T(\pi) = 0$, then $K_{3,3}(\pi) = 0$. Let us first admit this property and see how to conclude the proof of this corollary. 
According to Lemma \ref{lemmacomp} the pair $(3,3)$ is compatible.
The function $\Keq := \sqrt[6]{|S|^3+T^2}$ is continuous on $\H_{3,3}$, obeys the invariance property 
$\Keq(\pi\circ\phi) = |\det\phi|\Keq(\pi)$, 
and is such that $\Keq(\pi) = 0$ implies $K_{3,3}(\pi) = 0$.
We have seen in the proof of Theorem \ref{ThEquivMd} that these properties imply the desired equivalence of $\Keq$ and $K_{3,3}$.

We now show that $S(\pi) = T(\pi) = 0$ implies $K_{3,3}(\pi) = 0$.
A polynomial $\pi\in\H_{3,3}$ can be of two types. Either it is \emph{reducible}, meaning that there exist $\pi_1\in \H_{1,3}$ (linear) and $\pi_2\in \H_{2,3}$ (quadratic) such that $\pi = \pi_1 \pi_2$, or it is \emph{irreducible}.
In the latter case, according to \cite{Hartshorne}, there exists a linear change of coordinates $\phi$ and two reals $a,b$ such that 
$$
\pi\circ\phi = y^2 z - (x^3+ 3 a x z^2 + b z^3).
$$
A direct computation from the expressions given in \cite{Salmon} shows that $S(\pi\circ\phi) = a$ and $T(\pi\circ\phi) = -4b$.  If $S(\pi)=T(\pi) = 0$, then $S(\pi\circ\phi)=T(\pi\circ\phi) = 0$ and $\pi\circ \phi = y^2 z -x^3$. Therefore for all $\lambda\neq 0$, $\pi\circ \phi(\lambda x,\lambda^2 y, \lambda^{-3} z) = \lambda y^2 z -\lambda ^3 x^3$, which tends to $0$ as $\lambda\to 0$. We easily construct from this point a sequence $\phi_n$, $\det \phi_n = 1$, such that $\pi\circ\phi_n\to 0$. Therefore, $K_{3,3}(\pi) = 0$.

If $\pi$ is reducible, then $\pi = \pi_1 \pi_2$ where $\pi_1$ is linear and $\pi_2$ is quadratic. Choosing a linear change of coordinates $\phi$ such that $\pi_1\circ\phi = z$, we obtain
$$
\pi\circ\phi = 3 z (a x^2+ 2 b xy + c y^2)+ z^2 ( u x + v y + w z)
$$ 
for some constants $a,b,c,u,v,w$.
Again, a direct computation from the expressions given in \cite{Salmon} shows that $S(\pi\circ\phi) = -(ac -b^2)^2$ (and $T(\pi\circ\phi) = 8 (ac -b^2)^3$). Therefore, if $S(\pi) = T(\pi)= 0$ then the quadratic function $a x^2+ 2 b xy + c y^2$ of the pair of variables $(x,y)$ is degenerate. Hence there exists a linear change of coordinates $\psi$, altering only the variables $x,y$, and reals $\mu,u',v'$ such that
$$
\pi\circ\phi\circ\psi = \mu z x^2 + z^2 ( u' x + v' y + w z).
$$
It follows that $\pi\circ\phi\circ\psi(x,\lambda^{-1} y,\lambda z)$ tends to $0$ as $\lambda\to 0$. Again, this implies that $K_{3,3}(\pi)=0$, and concludes the proof of this proposition.
\sq

We could not find any example of an incompatible pair $(m,d)$, which leads us to formulate the conjecture that all pairs $(m,d)$ are compatible (hence providing ``simple'' equivalents of $K_{m,d}$ in full generality). Another even more difficult problem is to derive a polynomial $\Kpol$ of minimal degree for all pairs $(m,d)$ which are compatible and of interest. 

Last but not least, efficient algorithms are needed to compute metrics, from which effective triangulations are built that yield the optimal estimates. A possibility is to follow the approach proposed in \cite{C3}, i.e., solve numerically the optimization problem
$$ 
\inf \{\det M \sep  M\in S_d^+ \text{ and } \forall z \in \R^d, \<Mz,z\> \geq |\pi(z)|^{2/m}\},
$$
which amounts to minimizing a degree $d$ polynomial under an infinite set of linear constraints. When $d>2$, this minimization problem is not quadratic, which makes it rather delicate. Furthermore, numerical instabilities similar to those described in Remark \ref{overfitting} can be expected to appear.

\section{Conclusion and Perspectives}
\label{secOptAniso7}
In this chapter, we have introduced asymptotic estimates for the
finite element interpolation error measured in the $L^p$ norm when the mesh is 
optimally adapted to the interpolated function.
These estimates are asymptotically sharp for functions of two variables, see Theorem \ref{optitheorem}, and precise up to a fixed multiplicative constant in higher dimension, as described in \S \ref{secOptAniso6}. They involve a shape function $K_{m,p}$ (or $K_{m,d,p}$ if $d>2$) which generalizes the determinant which appears in estimates for piecewise linear interpolation \cite{CSX,BBLS}. This function can be explicitly computed in several cases, as Theorem \ref{equal23} shows, and has equivalents of a simple form in a number of other cases, see Theorems \ref{thequiv} and \ref{ThEquivMd}.

All our results are stated and proved for sufficiently smooth functions.
One of our future objectives is to extend these results to larger classes of functions, 
and in particular to functions exhibiting discontinuities along curves. This means that
we need to give a proper meaning to the nonlinear 
quantity $K_{m,p}\left(\frac{d^m f}{m!}\right)$
for nonsmooth functions. This question is addressed in Chapter 4.

This chapter also features a constructive algorithm (similar to \cite{BBLS}), that produces 
triangulations obeying our sharp estimates, and is described in \S\ref{secOptAniso3p2}. However, this algorithm
becomes asymptotically effective only for a highly refined triangulation. A more practical
way to produce quasi-optimal triangulations is to adapt them to a metric, see \cite{Bois,Shew,Peyre} and Chapter 5. 
This approach is discussed in \S \ref{secOptAniso4p2}. This raises the question of generating 
the appropriate metric from the (approximate) knowledge of the derivatives of the function to be interpolated. We addressed this question in the particular case of piecewise quadratic approximation in two dimensions in Theorems \ref{propEllipseMax} and \ref{family3}.

We plan to integrate this result in the PDE solver FreeFem++ in the near future, see Chapter 3 for a synthetic example. 
Note that a Mathematica source code is already available on the web \cite{sitejm}.
We also would like to derive appropriate metrics for other settings of 
degree $m$ and dimension $d$, although, as we pointed out in 
Proposition \ref{overfitting}, this might be a rather delicate matter.

We finally remark that in many applications, one seeks for error estimates in the Sobolev norms $W^{1,p}$ (or $W^{m,p}$) 
rather than in the $L^p$ norms. Finding the optimal triangulation for such norms requires a new error analysis, see Chapter 3. 
For instance, in the survey \cite{Shew2} on piecewise linear approximation,
it is observed that the metric $h_\pi = |d^2f|$ (evoked in Equation \iref{family2}) should be replaced with $h_\pi = (d^2f)^2$ for best adaptation in $H^1$ norm. In other words, the principal axes of the positive definite matrix $h_\pi$ remain the same, but its conditioning is squared.  

\section{Appendix}
\subsection{Proof of Proposition \ref{family3}}

We consider a fixed polynomial $\pi\in\H_3$, set a parameter $\alpha>0$, and look for an ellipse  $E_{\pi,\alpha}$ of maximal volume included in the set $\alpha^{-1/2} D \cap \Lambda_\pi$. Since this set is compact, a
standard argument shows that there exists at least one such ellipse. 

If $\alpha\geq \mu_\pi$, then $\alpha^{-1/2} D \subset \Lambda_\pi$, and therefore $\alpha^{-1/2} D\cap \Lambda_\pi = \alpha^{-1/2} D$. It follows that $E_{\pi,\alpha}=\alpha^{-1/2} D$, which proves part 1. 

In the following, we denote by $E'_\alpha$ the ellipse defined by the matrix \iref{hBigAlpha}. Note that any ellipse containing  $D_\pi$ and included in $\Lambda_\pi$ must be tangent to $\partial\Lambda_\pi$ at the point $z_\pi$, and hence of the form $E'_\delta$ for some $\delta>0$. Clearly, $E'_\delta\subset E'_\mu$ if and only if $\delta \geq \mu$. Therefore, $E'_\alpha \subset \Lambda_\pi$ if and only if $\alpha\geq \alpha_\pi$. 
Let $E$ be an arbitrary ellipse, $D_1$ the largest disc contained in $E$, and $D_2$ the smallest disc containing $E$. Then it is not hard to check that $|E| = \sqrt{|D_1| |D_2|}$.
For any $\alpha$ satisfying $\alpha_\pi\leq \alpha\leq \mu_\pi$, the ellipse $E'_\alpha$ is such that $D_1 = D_\pi$, which is the largest centered disc contained in $\Lambda_\pi$, and $D_2 = \alpha^{-1/2} D$, which corresponds to the bound $2\alpha^{-1/2}$ on the diameter of $E_{\pi,\alpha}$. 
It follows that $E'_\alpha$ is an ellipse of maximal volume 
included in $\alpha^{-1/2} D\cap \Lambda_\pi$, and this concludes the proof of part 2.

Part 4 is trivial; hence we concentrate on part $3$ and assume that  $\beta_\pi\leq \alpha\leq \alpha_\pi$.

An elementary observation is that $E_{\pi,\alpha}$ must be ``blocked with respect to rotations.'' Indeed, assume for contradiction that $R_\theta(E_{\pi,\alpha}) \subset \Lambda_\pi$ for $\theta \in [0,\ve]$ or $[-\ve,0]$, where we denote by $R_\theta$ the rotation of angle $\theta$. Observing that the set $\cup_{\theta \in [0,\ve]} R_\theta(E_{\pi,\alpha})$ contains an ellipse of larger area than $E_{\pi,\alpha}$ and of the same diameter, we obtain a contradiction.

In the following, we say that an ellipse $E$ is quadri-tangent to $\Lambda_\pi$ when there are at least four points of tangency between $\partial E$ and $\partial\Lambda_\pi$ (a tangency point being counted twice 
if the radii of curvature of $\partial E$ and $\partial\Lambda_\pi$ coincide at this point).

The fact that $E_{\pi,\alpha}$ is ``blocked with respect to rotations'' implies that it is either quadri-tangent to $\Lambda_\pi$ or tangent to $\partial\Lambda_\pi$ at the extremities of its small axis.
In the latter case, the extremities of the small axis must clearly be the points $z_\pi$ and $-z_\pi$, the closest points of $\partial \Lambda_\pi$ to the origin. It follows that $E_{\pi,\alpha}$ belongs to the family $E'_\delta$, $\delta \geq \alpha_\pi$ described above, and therefore is equal to $E'_{\alpha_\pi}$ since $\alpha\leq\alpha_\pi$. But $E'_{\alpha_\pi}$ is quadri-tangent to $\Lambda_\pi$, since otherwise we would have $E'_{\alpha_\pi-\ve}\subset \Lambda_\pi$ for some $\ve>0$.

We have now established that $E_{\pi,\alpha}$ is quadri-tangent to $\Lambda_\pi$ when $\beta_\pi\leq \alpha\leq \alpha_\pi$. This property is invariant by under linear change of coordinates: if an ellipse $E$ is quadri-tangent to $\Lambda_{\pi\circ\phi}$, then $\phi(E)$ is quadri-tangent to $\Lambda_\pi$. Furthermore, if $E$ is defined by a 
symmetric positive definite matrix $H$, then $\phi(E)$ is defined by $(\phi^{-1})^\trans H\phi^{-1}$.
This remark leads us to the problem of identifying the family of ellipses quadri-tangent to $\partial \Lambda_\pi$ when $\pi$ is among the four reference polynomials $x(x^2 - 3 y^2),\ x(x^2+3 y^2),\ x^2 y$, and $x^3$. In the case of $x^3$ there is no quadri-tangent ellipse and we have $\alpha_\pi=0$; therefore part 3 of the theorem is irrelevant.
In the three other cases, which respectively correspond to 
part 3 (i), (ii), and (iii), the quadri-tangent ellipses
are easily identified using the symmetries of these polynomials and the system of equations \iref{tangencyEqn}. 

The ellipses quadri-tangent to $x(x^2 + 3 y^2)$ are defined by matrices of the form $H_\lambda = \diag(\lambda,\frac{4+ \lambda^3}{3\lambda^2})$, where $0<\lambda\leq 2$. Note that $\det H_\lambda$ is decreasing on $(0,2^{\frac 1 3}]$ and increasing on $[2^{\frac 1 3},2]$. 
Given $\pi$ with $\disc\pi<0$, the optimization problem \iref{hConstrained} therefore becomes
$$
\min_\lambda\{\det H_\lambda \sep (\phi_\pi^{-1})^\trans H_\lambda \phi_\pi^{-1} \geq \alpha \Id\}.
$$
If the constraint is met for $\lambda = 2^{1/3}$, we obtain $E_{\pi,\alpha}=E_\pi$ and therefore $\alpha\leq \beta_\pi$. Otherwise, using the monotonicity of $\lambda\mapsto \det H_\lambda$ on each side of its minimum $2^{\frac 1 3}$ we see that the matrix $H_\lambda - \alpha \phi_\pi^\trans \phi_\pi$ must be singular. Taking the determinant, we obtain an equation of degree $4$ from which $\lambda$ can be computed, and this concludes the proof of part 3 (i).

The ellipses quadri-tangent to $x(x^2 - 3 y^2)$ are defined by  $H_{\lambda, V} = V^\trans \diag(\lambda,\frac{4- \lambda^3}{3\lambda^2}) V$, where $0<\lambda\leq 1$ and $V$ is a rotation by $0$, $60$, or $120$ degrees. Since $\det H_{\lambda,V}$ is a decreasing function of $\lambda$ on $(0,1]$, we can apply the same reasoning as above to polynomials $\pi$ such that $\disc \pi>0$. This concludes the proof of part 3 (ii). 

Finally, the ellipses quadri-tangent to $x y ^2$ are defined by $H_\lambda = \diag(\lambda,\frac 4 {27\lambda^2})$, $\lambda>0$. The determinant $\lambda \mapsto \det H_\lambda$ is a decreasing function of $\lambda$ with lower bound $0$ as $\lambda\to\infty$, and the same reasoning applies again, hence concluding the proof of part 3 (iii).

\subsection{Proof of the nonexistence property in Theorem \ref{ThEquivMd}}
\label{appendixBOA}
%

Let $\Kpol$ be a polynomial on $\H_{m,d}$ which satisfies the inequalities \iref{equivKmd}. We observe that $\Kpol$ takes nonnegative values on $\H_{m,d}$, and we assume in a first time that $\mu := \frac {mr} d$ is an integer.
We derive from \iref{equivKmd} and from the invariance of the shape function $K_{m,d}$ with respect to changes of variables the inequalities
\be
C^{-2 r} (\det \phi)^\mu\Kpol(\pi)\leq \Kpol(\pi\circ\phi)\leq C^{2 r} (\det \phi)^\mu\Kpol(\pi),
\label{equivInvarPol}
\ee
for all $\pi \in \H_{m,d}$ and all $\phi\in M_d$ (the space of $d\times d$ real matrices),
where $C$ is the constant appearing in inequalities \iref{equivKmd}.
We regard the function $\b Q(\pi,\phi) = \Kpol(\pi\circ\phi)$ as a polynomial on the vector space $V= \H_{m,d}\times M_d $ and we observe that it vanishes on the \emph{hypersurface} $V_{\det} := \{(\pi,\phi)\in V \sep\det \phi = 0\}$.
Since $\phi\mapsto\det(\phi)$ is an irreducible polynomial, as shown in \cite{Bochner}, it follows that $\b Q(\pi,\phi) = (\det \phi) \b Q_1(\pi,\phi)$ for some polynomial $\b Q_1$ on $V$. Injecting this expression in \iref{equivInvarPol} we obtain that $\b Q_1(\pi,\phi)$ also vanishes on the hypersurface $V_{\det}$ if $\mu>1$, and the argument can be repeated. By induction we eventually obtain a polynomial $\widehat \Kpol$ on $V$ such that $\Kpol(\pi\circ\phi) = (\det \phi)^\mu \widehat \Kpol(\pi,\phi)$. It follows from inequality \iref{equivInvarPol} that for all $(\pi,\phi)\in V$,
$$
C^{-2 r} \Kpol(\pi) \leq \widehat \Kpol(\pi,\phi)\leq C^{2 r} \Kpol(\pi).
$$
This implies that $\widehat \Kpol(\pi,\phi)$ does not depend on $\phi$. Otherwise, since it is a polynomial, we could find $\pi_1\in H_{m,d}$ and a sequence $\phi_n\in M_d$ such that $|\widehat \Kpol(\pi_1,\phi_n)|\to \infty$.
Therefore,
$$
\Kpol(\pi\circ\phi) = (\det\phi)^\mu \widehat \Kpol(\pi,\phi)=(\det\phi)^\mu \widehat \Kpol(\pi,\Id)= (\det\phi)^\mu \Kpol(\pi).
$$
This establishes the invariance property of $\Kpol$ under the hypothesis that $\mu$ is an integer.
Let us assume for contradiction that $\mu$ is not an integer. Then applying the previous reasoning to $\Kpol^d$ we obtain that $\Kpol(\pi\circ \phi)/ \Kpol(\pi) = \alpha(\pi, \phi) |\det \phi|^\mu$ for all $\pi \in \H_{m,d}$ and $\phi\in M_d$, where $\alpha : \H_{m,d} \times M_d \to \{-1,1\}$. Since $\phi \mapsto \det \phi$ is an irreducible polynomial we obtain  that $\mu$ is an integer and as before that $\Kpol(\pi\circ\phi) = (\det\phi)^\mu \Kpol(\pi)$, which concludes the proof of the invariance property.

Let $(m,d)$ be an incompatible pair. We know from Lemma \ref{lemmacomp} that there exists $\pi_0\in \H_{m,d}$ such that $K_{m,d}(\pi_0) >0$ and $Q(\pi_0)=0$ for any invariant polynomial $Q\in \I_{m,d}$. Therefore there exists no polynomial $\Kpol$ satisfying \iref{equivKmd}.

\chapter[Sharp asymptotics of the $W^{1,p}$ interpolation error]{Sharp asymptotics of the $W^{1,p}$ interpolation error on optimal triangulations} 
\minitoc
\label{chapW1P}
\section{Introduction}

We consider in this section the problem of finding the optimal mesh for the interpolation of a function by finite elements of a given degree, when the error is measured in the Sobolev $W^{1,p}$ norm.
Our purpose is to establish sharp asymptotic estimates, of the same type as the following in the case of the $L^p$ norm,
\be
\limsup_{N\to +\infty} \(N\min_{\#(\cT)\leq N}\|f-\interp_\cT^1 f\|_{L^p}\)  \leq C\| \sqrt{|\det(d^2f)|}\|_{L^\tau},\;\; \frac 1 \tau=\frac 1 p+1,
\label{optiaffineW1P}
\ee
which holds for a $C^2$ smooth function $f$ defined on a bounded polygonal domain $\Omega$, see  \cite{CSX, BBLS} and Chapter 2.
We refer to the introduction of this thesis and to Chapter \ref{secOptAniso} for a more detailed overview of existing results and of our motivations.
The convergence estimate \iref{optiaffineW1P} is extended 
to arbitrary approximation order in Chapter \ref{chapOptAniso}, where the quantity governing 
the convergence rate for finite elements of arbitrary degree $m-1$ is identified. 
This quantity, referred to as the shape function, depends nonlinearly on the $m$-th order derivative $d^mf$.

The purpose of the present chapter is to investigate this problem when the $L^p$-norm is
replaced by the $W^{1,p}$ semi-norm which plays a critical role in PDE analysis, and which is defined as follows
$$
|f|_{W^{1,p}(\Omega)} := \|\nabla f\|_{L^p(\Omega)} = \left(\int_\Omega |\nabla f|^p\right)^{1/p}.
$$
Our second objective is to propose simple and practical ways
of designing meshes which behave similar to the optimal one,
in the sense that they satisfy the sharp error estimate up
to a fixed multiplicative constant. 

\subsection{Main results and layout}

We denote by $\P_m$ the space of 
polynomials of total degree less or equal to $m$
and by $\H_m$ the space of homogeneous polynomials of total degree $m$,
$$
\P_m:={\rm Span}\{ x^ky^l\; ; \; k+l\leq m\} \;\;{\rm and}\;\; \H_m:={\rm Span}\{ x^ky^l\; ; \; k+l=m\}.
$$
For any triangle $T$, we denote by $\interp_T^m$ the local interpolation
operator acting from $C^0(T)$ onto $\P_m$. 
For any continuous fonction $\nu \in C^0(T)$, the interpolating polynomial $\interp_T^m \nu\in \P_m$ is defined by the conditions
$$
\interp_T^m \nu(\gamma)=\nu(\gamma),
$$
for all points $\gamma\in T$ with barycentric coordinates in
the set $\{0, \frac 1 m,\frac 2 m,\cdots,1\}$. 
This interpolation operator is invariant by translation, hence for any polynomial $\pi\in \H_m$, triangle $T$ and offset $z$ we have
\be
\label{transInv}
 |\pi - \interp_T^{m-1} \pi|_{W^{1,p}(T)} =  |\pi - \interp_T^{m-1} \pi|_{W^{1,p}(z+T)}.
\ee
If $\cT$ is a triangulation of a domain $\Omega$, then $\interp_\cT^m$ 
refers to the interpolation operator which coincides with $\interp_T^m$ on each triangle $T\in \cT$. 

A key ingredient of this chapter is the {\it shape function} $L_{m,p}$, which is  
defined by a {\it shape optimization problem}:
for any fixed $1\leq p \leq \infty$ and for any $\pi\in \H_m$, we define
\be
L_{m,p}(\pi):=\inf_{|T|=1} |\pi - \interp_T^{m-1} \pi|_{W^{1,p}(T)}.
\label{shapeFunctionL}
\ee
Here, the infimum is taken over all triangles of area $|T|=1$. 
From the homogeneity of $\pi$, it is easily checked that
$$
\inf_{|T|=A} |\pi - \interp_T^{m-1} \pi|_{W^{1,p}(T)} = L_{m,p}(\pi)A^{\frac {m-1} 2 + \frac 1 p}.
$$
The solution to this optimization problem thus describes the shape of the triangles 
of a given area which are best adapted
to the polynomial $\pi$ in the sense of minimizing the interpolation error
measured in $W^{1,p}$.  

The function $L_{m,p}$ is the natural generalisation of the function $K_{m,p}$ introduced in Chapter \ref{chapOptAniso} for the study of optimal anisotropic triangulations in the sense of the $L^p$ interpolation error
$$
K_{m,p}(\pi):=\inf_{|T|=1} \|\pi - \interp_T^{m-1} \pi\|_{L^p(\Omega)}.
$$
Our asymptotic error estimate for the optimal triangulation is given by
the following theorem.
\begin{theorem}
\label{mainTheorem}
For any bounded polygonal domain $\Omega\subset \R^2$, any function 
$f\in C^m(\overline \Omega)$ and any $1\leq p <\infty$, there exists a sequence 
$\seqT$, $\#(\cT_N)\leq N$, of triangulations of $\Omega$ such that
\be
\label{upperEstimW1P}
\limsup_{N\ra \infty} N^{\frac {m-1} 2} |f-\interp_{\cT_N}^{m-1} f|_{W^{1,p}(\Omega)} \leq 
\left\|L_{m,p}\left(\frac{d^m f}{m!}\right)\right\|_{L^\tau(\Omega)},\text{ where } \frac 1 \tau := \frac {m-1} 2 +\frac 1 p.
\ee
\end{theorem}
In the above convergence estimate, the number $N_0$ is independent of $f$
and refers to the 
minimal cardinality of a conforming triangulation of $\Omega$.  
The $m$-th derivative $d^mf(z)$ at each point $z$ is identified to 
a homogeneous polynomial $\pi_z\in \H_m$:
\be
\label{dmfHmW1P}
\frac{d^mf(z)}{m!}\sim \pi_z = \sum_{k+l=m}\frac{\partial^m f}{\partial x^k\partial y^l}(z)\frac{x^k}{k!}\frac{y^l}{l!}. 
\ee
An important feature of this estimate is the ``$\limsup$''. 
Recall that the upper limit of a sequence $(u_N)_{N\geq N_0}$ is defined by 
$$
\limsup_{N\to \infty} u_N := \lim_{N\to \infty} \sup_{n\geq N} u_n,
$$
and is in general stricly smaller than the supremum $\sup_{N\geq N_0} u_N$. It is still an open question to find an appropriate upper estimate of $\sup_{N\geq N_0} N^{\frac{m-1} 2}  |f-\interp_{\cT_N}^{m-1} f|_{W^{1,p}(\Omega)}$ when optimally adapted anisotropic triangulations are used.

In order to illustrate the sharpness of \iref{upperEstimW1P}, we introduce
a slight restriction on sequences of triangulations, following 
an idea in \cite{BBLS}: a sequence $\seqT$ of triangulations is said to be \emph{admissible} if
$\#(\cT_N) \leq  N$ and $\sup_{N \geq N_0} (N^\frac 1 2 \sup_{T\in \cT_N} \diam(T)) < \infty$. In other words if
\be
\label{admissibilityCond}
\sup_{T\in \cT_N} \diam(T) \leq C_AN^{-\frac 1 2}
\ee
for some constant $C_A>0$ independent of $N$. The following theorem shows that the estimate
\iref{upperEstimW1P} cannot be improved when we restrict our attention to admissible sequences.
It also shows that this class is reasonably large in the sense that 
\iref{upperEstimW1P} is ensured to hold up to small perturbation.

\begin{theorem}
\label{optiTheorem}
Let $\Omega\subset \R^2$ be a bounded polygonal domain, let $f\in C^m(\overline\Omega)$ and $1\leq p <\infty$. We define $ \frac 1 \tau := \frac {m-1} 2 +\frac 1 p$.
For any \emph{admissible} sequence $\seqT$ of triangulations of $\Omega$, one has
\be
\label{lowerEstim}
\liminf_{N\ra \infty} N^{\frac {m-1} 2} |f-\interp_{\cT_N}^{m-1} f|_{W^{1,p}(\Omega)} \geq \left\|L_{m,p}\left(\frac{d^m f}{m!}\right)\right\|_{L^\tau(\Omega)}.
\ee
Furthermore, for all $\ve>0$ there exists an \emph{admissible} sequence $(\cT_N^\ve)_{N\geq N_0}$ of triangulations of $\Omega$ such that
\be
\label{upperEstimEpsW1P}
\limsup_{N\ra \infty} N^{\frac {m-1} 2}  |f-\interp_{\cT_N^\ve}^{m-1} f|_{W^{1,p}(\Omega)} \leq \left\|L_{m,p}\left(\frac{d^m f}{m!}\right)\right\|_{L^\tau(\Omega)}+\ve.
\ee
\end{theorem}
Note that the sequences $(\cT_N^\ve)_{N\geq N_0}$ satisfy the admissibility condition \iref{admissibilityCond} with a constant $C_A(\ve)$ which may grow to $+\infty$ as $\ve\to 0$.
The proofs of these two theorems are given in \S \ref{secW1P3}. 
Theorem \ref{mainTheorem} is a direct consequence of Theorem \ref{optiTheorem}, by
considering a sequence of triangulations of 
the type $\cT_N^{\e_N}$ with $\e_N\to 0$ as $N\to +\infty$.
The proof of the upper estimate in Theorem \ref{optiTheorem} involves the construction
of an optimal mesh based on a patching strategy adapted from the one encountered in \cite{BBLS}. 
However, inspection of the proof reveals that
this construction only becomes effective as
the number of triangles $N$ becomes very large. Therefore it 
may not be useful in practical applications.

\begin{remark}
It can easily be shown that if $(\cT_N)_{N\geq N_0}$ is an admissible sequence of triangulations and
$f\in\cC^m(\Omega)$, then $\|f-\interp_{\cT_N}^m f\|_{L^p(\Omega)}$
decays with the rate $N^{-m/2}$ which is faster than the 
decay rate obtained for the $W^{1,p}$ error. Therefore, our convergence estimates
are also valid in the $W^{1,p}$ norm 
$$
\|f\|_{W^{1,p}(\Omega)} := \(\|f\|_{L^p(\Omega)}^p+|f|_{W^{1,p}(\Omega)}^p\)^{1/p}.
$$ 
\end{remark}

We show in \S \ref{secW1P2} that in order to satisfy the optimal estimate \iref{upperEstimW1P} up to a fixed multiplicative constant, it suffices to build a triangulation which obeys four general principles:

(i) The interpolation error should be evenly distributed on all triangles,

(ii) The triangles should adopt locally a specific aspect ratio, dictated by the local value of $d^m f$, 

(iii) the largest angle of the triangles should be bounded away from $\cPi = 3.14159\ldots$

(iv) the triangulation $\cT$ should be sufficiently refined in order to adapt to the local features of $f$.
\nl

The third point (iii) is the main new ingredient of this chapter compared to Chapter \ref{chapOptAniso}, and is necessary for obtaining $W^{1,p}$ error estimates (but not for $L^p$ error estimates). Roughly speaking, two triangles having the same optimized aspect ratio imposed
by (ii) may greatly differ in term of their largest angle, and the most acute triangle should
be preferred when the interpolation error is measured in $W^{1,p}$ rather than $L^p$.
The influence of large angles in mesh adaptation has already been studied in \cite{Ba,Ja,Shew2}, see also Chapter 6. 
The heuristic guideline is that large angles should be avoided in general, since they lead to oscillations of the gradient of the interpolant. On the contrary, extremely thin triangles and very small angles can be necessary for optimal mesh adaptation.  

A practical approach for mesh generation is discussed in \S \ref{secW1P4}, and consists in deriving a distorted metric from 
the exact or approximate data of $d^m f$ at each point $x\in\Omega$. We restrict in this section to the case of linear and quadratic finite elements, and we provide simple mesh generation procedures and numerical results. 
To any $\pi\in \H_m$, $m\in \{2,3\}$, we associate a 
symmetric positive definite matrix $\cM_m(\pi)\in S_2^+$ (strictly speaking, this matrix is degenerated if $\pi$ is univariate, a detail that needs to be taken care of in the numerical implementation). If $z\in \Omega$ and $d^m f(z)$ is close to $\pi$, then the triangle $T$ containing $z$ should be isotropic in the metric $\cM_m(\pi)$.
The requirements (i) and (ii) above, which are respectively linked to the size and shape of the triangles, can then be summarized through a global metric on $\Omega$ given by
\be
\label{defPiZ}
h(z)=s(\pi_z) \cM_m(\pi_z),\;\; \pi_z=\frac{d^mf(z)}{m!},
\ee
where $s(\pi_z)$ is a scalar factor which depends on the desired accuracy
of the finite element approximation.
Once this metric has been properly identified, fast algorithms such as in \cite{Inria, FreeFem, Peyre}
can be used in order to design a near-optimal mesh based on it. Recently it has been rigorously proved in \cite{Shew,Bois}, that several algorithms terminate and produce good quality meshes, under certain conditions. Although we are not aware that the angle constraint (iii) is guaranteed
in such algorithms, it seems to hold in practice.
Computing the map
$$
\pi\in \H_m \mapsto \cM_m(\pi)\in S_2^+,
$$
is therefore of key use in applications ($S_2^+$ refers to the set of $2\times 2$ symmetric and non-negative matrices). This problem is solved in \cite{Shew},
in the case of linear elements ($m=2$): the matrix $\cM_m(\pi)$ is then defined as the
square of the matrix associated to
the quadratic form $\pi$. 
We give a simple expression of $\cM_m(\pi)$ for piecewise quadratic finite elements ($m=3$). The optimality of this construction is proved theoretically, and numerical experiments confirm its adequacy.
An open source implementation for the mesh generator FreeFem++ \cite{FreeFem} is provided
at \cite{sitejm}.

The shape function $L_{m,p}$ does not always have a simple analytic expression from the coefficients of $\pi$. For this reason we introduce in \S \ref{secW1P5} explicit functions $\pi\mapsto \Lpol_m(\pi)$ which are defined as the root of a polynomial in the coefficients of $\pi$, and are equivalent to $L_{m,p}$, leading therefore to similar
asymptotic error estimates up to multiplicative constants.
We finally discuss in \S \ref{secW1P6} the possible extension of our analysis to simplicial
elements in dimension $d>2$. 

\subsection*{Notations}
Throughout this chapter, we define $L_m := L_{m,\infty}$, where $L_{m,p}$ is defined at Equation \iref{shapeFunctionL}. 
We prove further in Lemma \ref{lemmaEquivLm} that for all $1\leq p\leq \infty$
\be
\label{equivLm}
c L_m\leq L_{m,p}\leq L_m \text{ on } \H_m,
\ee
where the constant $c>0$ depends only on $m$.
For any compact set $E\subset \R^d$, of non-zero Lebesgue measure, we denote by $\bary(E)$ its barycenter. 
For any pair of vectors $u,v\in\R^d$,
we denote by $\<u,v\>$ their inner product, and by
$$
|u|:=\sqrt{\<u,u\>},
$$
the euclidean norm of $u$. When $g\in L^p(E, \R^d)$ is a vector valued function,
we denote by $\|g\|_{L^p(E)}$ the $L^p$ norm of $x\mapsto |g(x)|$ on $E$.

We denote by $\M_d(\R)$ the set of all $d\times d$ real matrices, equiped with the norm
$$
\|A\|:=\max_{|u|\leq 1}|A u|.
$$
We denote by $\GL_d\subset \M_d(\R)$ the linear group of invertible matrices and by $\SL_d\subset \GL_d$ the special linear group of matrices of determinant $1$.
$$
\GL_d := \{A\in \M_d(\R)\sep \det A \neq 0\} \text{ and } \SL_d:=\{A\in \M_d(\R) \sep \det A = 1\}.
$$
For $A\in \GL_d$, we denote by
\be
\label{condnumber}
\kappa(A):=\|A\|\, \|A^{-1}\|,
\ee
its condition number.
We denote by $S_d\subset \M_d(\R)$ the subset of symmetric matrices, by $S_d^\oplus \subset S_d$ the subset of non-negative symmetric matrices and by $S_d^+$ the subset of positive definite symmetric matrices.
For any two symmetric matrices $S,S'\in S_d$, we write $S\leq S'$ if and only if $S'-S\in S_d^\oplus$.

We equip the spaces $\P_m$ and $\H_m$ with the norm
\be
\label{normhm}
\|\pi\|:=\max_{|u|\leq 1}|\pi(u)|.
\ee
Note that the greek letter $\pi$ always refers to an homogeneous polynomial $\pi\in \H_m$, 
while the large and bold notation $\cPi$ refers to the mathematical constant $\cPi= 3.14159\ldots$.

Recall that if $g$ is a $C^m$ function, we identify $d^mg(x)$ to a polynomial in $\H_m$, by \iref{dmfHmW1P}. We
then denote
\be
\|d^m g\|_{L^\infty(E)}:=\max_{z\in E} \|d^m g(z)\|
\label{normdmgW1P}
\ee
with $\|\cdot\|$ the previously defined norm on $\H_m$.

\section{The shape function $L_{m,p}$ and local error estimates}
\label{secW1P2}
In this section, we study the function $L_{m,p}$ and obtain local $W^{1,p}$ error estimates for functions of two variables. These estimates naturally give rise to a heuristic method for the design of ``near optimal'' triangulations adapted to a function $f$, in other words triangulations satisfying the estimate \iref{upperEstimW1P} up to a fixed multiplicative constant, and it is put into practice in \S \ref{secW1P4} in the case of piecewise linear and piecewise quadratic finite elements.
The results of this section are also useful to the proof, in \S \ref{secW1P3}, of the optimal error 
estimates presented in Theorems \ref{mainTheorem} and \ref{optiTheorem}.

We first introduce the measure of sliverness $S(T)$ of a triangle $T$. Given two triangles $T,T'$, there are precisely $6$ affine transformations $\Psi$ such that $\Psi(T) = T'$. Each of these affine transformations $\Psi$ defines a linear transformation $\psi$ and we set 
\be
\label{defDist}
d(T,T') := \ln \left( \inf\{\kappa(\psi)\sep \Psi(T) = T'\}\right),
\ee
where $\kappa(\psi)$ is the condition number defined in \iref{condnumber}.
Clearly $d(T,T')\geq 0$, $d(T,T')=d(T',T)$ and $d(T,T'')\leq d(T,T')+ d(T',T'')$. 
Furthermore $d(T,T') = 0$ if and only if $T$ can be transformed into $T'$ through a translation, a rotation and a dilatation. Therefore $d(\cdot,\cdot)$ defines a distance between shapes of triangles. The heuristic guideline of the papers \cite{Ja,Ba} and Chapter 6 is that obtuse shapes should be avoided when possible in the design of Finite Element meshes for the approximation of a function in the $W^{1,p}$ norm. We therefore introduce the set $\mathbb A$ of all acute triangles
and we define the {\it measure of sliverness} $S(T)$ of a triangle $T$ as follows
\be
\label{defS}
S(T) := \exp d(T, \mathbb A) = \inf \{ \kappa(\psi) \sep \Psi(T)\in \mathbb A\}.
\ee
This quantity reflects the distance from $T$ to the set of acute triangles:
in particular $S(T)=1$ if and only if $T\in\mathbb A$, and $S(T)>1$ otherwise.
It has an analytic  expression, which is given in the following proposition.

\begin{prop} For any triangle $T$ with largest angle $\theta$, one has 
\label{propSTan}
$S(T) = \max\{1, \tan \frac \theta 2\}.$
\end{prop}

\proof
The result of this proposition is trivial if the triangle $T$ is acute, and we therefore assume that $T$ is obtuse. We can assume without loss of generality that the vertices of $T$ are $0$, $\alpha u$ and $\beta v$, where $\alpha, \beta>0$, $u,v\in \R^2$, $|u|=|v|=1$ and $\<u,v\> = \cos \theta$. Note that $|u-v|= 2 \sin (\theta/2)$ and $|u+v|= 2 \cos (\theta/2)$. Let $\Psi$ be such that $\Psi(T)\in \mathbb A$, and let $\psi$ be the associated linear transform. Since $\Psi(T)$ is acute we have $\<\psi(u),\psi(v)\>\geq 0$ and therefore $|\psi(u)-\psi(v)| \leq |\psi(u)+\psi(v)|$. It follows that 
$$
\kappa(\psi) = \|\psi\|\; \|\psi^{-1}\|\geq\frac {|u-v|} {|u+v|} \times \frac {|\psi(u)+\psi(v)|} {|\psi(u)-\psi(v)|} \geq \frac {2 \sin (\theta/2)}{2\cos(\theta/2)}  = \tan \frac \theta 2.
$$ 
Therefore $S(T)\geq \tan \frac \theta 2$.
Furthermore, let $\psi$ be defined by $\psi(u)=(0,1)$ and $\psi(v) = (1,0)$. Obviously 
$\psi(T)$ has one of its angles equal to $\cPi/2$ and is therefore acute. On
the other hand, one easily checks that $\kappa(\psi) = \tan (\theta/2)$ and therefore
$S(T)\leq \tan \frac \theta 2$. This concludes the proof of this proposition.
\sq

The previous proposition implies in particular that $S(T)$ is equivalent to the quantity $\frac 1 {\sin \theta}$, where $\theta$ denotes the largest angle of $T$, which is used in \cite{Ja,Ba}.  
The following lemma shows that the interpolation process
is stable with respect to the $L^\infty$ norm of the gradient
if the measure of sliverness $S(T)$ is controlled. Let us mention
that a slightly different formulation of this result was already proved in \cite{Ja},
yet not exactly adapted to our purposes.

\begin{lemma}
\label{lemmaS}
There exists a constant $C=C(m)$ such that 
for any triangle $T$ and any $f\in W^{1, \infty}(T)$ one has
\be
\label{ineqS}
\|\nabla \interp^{m-1}_T f \|_{L^\infty(T)} \leq C S(T) \|\nabla f \|_{L^\infty(T)}, 
\ee
\end{lemma}

\begin{figure}
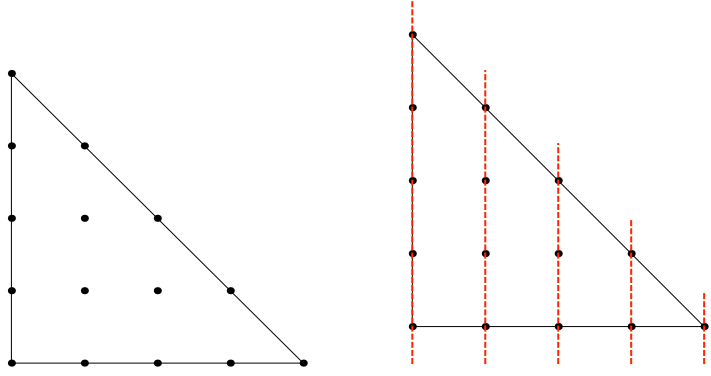

	\centering
		\includegraphics[width=4cm,height=4cm]{\pathPic/PaperW1P/Lagrange4.pdf}
		\hspace{1cm}
		\includegraphics[width=4cm,height=5cm]{\pathPic/PaperW1P/Lagrange4Lines.pdf}
	\caption{\label{fig1W1P}The interpolation points on the triangle $T_0$ are aligned vertically.}
\end{figure}

\proof
Let $T_0$ be the triangle of vertices $(0,0)$, $(1,0)$ and $(0,1)$, and let $g\in W^{1,\infty}(T_0)$. We define $\ti g(x,y):= g(x,0)$ and $h(x,y):=g(x,y)-g(x,0)$. 
Since $\ti g$ does not depend on $y$ and since the Lagrange interpolation points on $T_0$ are aligned vertically, as illustrated on Figure \ref{fig1W1P}, the Lagrange interpolant $\interp_{T_0}^{m-1} \ti g$ does not depend on $y$ either. Futhermore, for all $(x,y)\in T_0$, we have 
$|h(x,y)| = |\int_{s = 0}^y \frac {\partial g}{\partial y}(x,s) ds | \leq \|\frac {\partial g}{\partial y}\|_{L^\infty(T_0)}$.
Hence
\begin{eqnarray*}
\left\|\frac {\partial \interp_\TRect^{m-1} g} {\partial y}   \right\|_{L^\infty(\TRect)}
&=& \left\| \frac {\partial \interp_\TRect^{m-1} h}{\partial y}  \right\|_{L^\infty(\TRect)}\\
&\leq & C_0 \|\interp_\TRect^{m-1} h\|_{L^\infty(\TRect)}\\
&\leq & C_0 C_1 \|h\|_{L^\infty(\TRect)} \\
&\leq & C_0 C_1  \left\|\frac {\partial g}{\partial y}\right\|_{L^\infty(\TRect)},
\end{eqnarray*}
where the constants $C_0$ and $C_1$ are the $L^\infty(T_0)$ norms of the operators $g\mapsto \frac{\partial g}{\partial y}$
restricted to $\P_{m-1}$ and $g\mapsto \interp_{T_0}^{m-1} g$ respectively.

Let $e$ be an edge vector of $T$. There exists an affine change of coordinates $\Psi_e$, with linear part $\psi_e$, such that $T_0 = \Psi_e(T)$ and $\psi_e(e) = e_0$ where $e_0=(0,1)$ is the vertical edge vector of $T_0$.
Noticing that 
$$
\<e,\nabla \interp_T^{m-1} (g\circ \Psi_e)\> = \<e,\nabla ((\interp_\TRect^{m-1} g)\circ \Psi_e)\> = \<e_0, (\nabla \interp_\TRect^{m-1} g) \circ \Psi_e\> = \frac{\partial  \interp_\TRect^{m-1} g}{\partial y} \circ \Psi_e, 
$$
we obtain
\be
\label{edgeGrad}
\begin{array}{rcl}
\|\<e, \nabla \interp_T^{m-1} (g\circ \Psi_e)\>\|_{L^\infty(T)} &=& \left\|\frac{\partial \interp_{T_0}^{m-1} g}{\partial y}  \right\|_{L^\infty(\TRect)} \\
&\leq& 
C_0C_1 \left\|\frac{\partial g}{\partial y} \right\|_{L^\infty(\TRect)}\\
&=&C_0C_1 \|\<e, \nabla (g\circ \Psi_e)\>\|_{L^\infty(T)}.
\end{array}
\ee
Applying this inequality to $g=f\circ \Psi_e^{-1}$ we obtain that 
\be
\label{ineqSfe}
\|\<e, \nabla \interp_T^{m-1} f\>\|_{L^\infty(T)} \leq C_0C_1 \|\<e, \nabla f\>\|_{L^\infty(T)},
\ee
for all edge vectors $e\in \{a,b,c\}$ of $T$. We next define a norm on $\R^2$ as follows
$$
|v|_T := |a|^{-1} |\<a,v\>|+|b|^{-1} |\<b,v\>|+|c|^{-1} |\<c,v\>|.
$$
It follows from inequality \iref{ineqSfe}, that 
\be
\label{gradNormT}
\| \ |\nabla \interp_T^{m-1} f|_T\ \|_{L^\infty(T)} \leq 3C_0C_1 \| |\nabla f|_T \|_{L^\infty(T)}.
\ee
We next observe that if $\theta$ denotes the maximal angle of $T$, 
$$
\cos(\theta/2) |v| \leq |v|_T\leq 3 |v|,
$$
where $|\cdot|$ is the euclidean norm: the upper inequality is trivial
and the lower one is implied by the fact that at least one of the edge vectors makes
an angle less than $\theta/2$ with $v$. Combining this with \iref{gradNormT}, we obtain
$$
\| \nabla \interp_T^{m-1} f\ \|_{L^\infty(T)} \leq \frac {9 C_0C_1}{\cos(\theta/2)}
 \| \nabla f \|_{L^\infty(T)}.
$$
Since $\theta> \pi/3$ we have $\frac 1 {\cos(\theta/2)}\leq 2 \tan(\theta/2) \leq 2S(T)$
according to Proposition \ref{propSTan}, which concludes the proof with $C=18C_0C_1$.
\sq

\begin{remark}The following example in the simple case of piecewise linear approximation
illustrates the sharpness of inequality \iref{ineqS}. Let $T$ be a triangle having an obtuse angle 
$\theta$ at a vertex $v$, and edges neighbouring $v$ of length $l$ and $l'$. Let $f(z) := |z-v|^2$. A simple computation shows that 
$$
\|\nabla \interp_T^1 f\|_{L^\infty(T)} = \frac {\diam T}{\sin \theta} \ \text{ and } \ \|\nabla f\|_{L^\infty(T)} = 2\max(l,l').
$$ 
It follows that
$$
\|\nabla \interp_T^1 f\|_{L^\infty(T)} =\lambda(T) S(T) \|\nabla f\|_{L^\infty(T)},
$$
with
$$
\lambda(T) := \frac{\diam(T)}{ 2\sin(\theta)S(T) \max\{l,l'\}}=\frac{\diam(T)}{4 \sin(\theta/2)^2 \max\{l,l'\}}\in [1/4,\, 1],
$$
which shows the sharpness of Lemma \ref{lemmaS} in this context.
\end{remark}

We now introduce for each polynomial $\pi\in \H_m$, a special set $\cA_\pi \subset \M_2(\R)$ of linear maps.  
\be
\label{defA}
\cA_\pi := \{A\in \M_2(\R)\sep |\nabla \pi(z)| \leq |Az|^{m-1} \text{ for all } z\in \R^2\}. 
\ee
This set has a geometrical interpretation : since $\nabla\pi$ is homogeneous
of degree $m-1$, we find that $A\in \cA_\pi$ if and only if the ellipse $\{z\in \R^2 \sep |A z|\leq 1 \}$ is included in the set $\{z\in \R^2 \sep |\nabla \pi(z)|\leq 1\}$ which is limited by
the algebraic curve $\{|\nabla \pi(z)|= 1\}$. 
If $T$ is a triangle that contains the origin and if $A\in \cA_\pi$, we observe that 
\be
\label{ineqPi}
\|\nabla \pi \|_{L^\infty(T)}   \leq \diam(A(T))^{m-1}.
\ee
We define
$$
\gamma_m(\pi):= \inf \{|\det A|\sep A\in \cA_\pi\},
$$
so that $\frac {\cPi}{\gamma_m(\pi)}$ is the maximal area of an ellipse 
contained in $\{z\in \R^2 \sep |\nabla \pi(z)|\leq 1\}$.

\begin{remark} 
Similar concepts have been
introduced in \cite{C3}
for the purpose of studying the $L^p$ interpolation error
of anisotropic finite elements, therefore with $|\pi(z)|$ in place of $|\nabla\pi(z)|$. 
\end{remark}
\noindent
The following result shows that a certain power of $\gamma_m$ is equivalent
to the shape function $L_m$. For any domain $\Omega\subset \R^d$ and any $f,g\in L^p(\Omega)$ we use the shorthand
\be
\label{defNormFG}
\|(f,g)\|_{L^p(\Omega)} := \left(\int_{\Omega} |(f(z),g(z))|^p dz \right)^\frac 1 p =  \left(\int_{\Omega} (f(z)^2+g(z)^2)^{p/2} dz \right)^\frac 1 p.
\ee

\begin{lemma}
\label{lemmaA}
There exists a constant $C=C(m)$ such that for all $\pi\in \H_m$
\be
\label{eqLA}
C^{-1} L_m(\pi)\leq 
\gamma_m(\pi)^{\frac{m-1} 2}
\leq C L_m(\pi).
\ee
\end{lemma}

\proof
We first prove the left part of \iref{eqLA}.
Let $\pi \in \H_m$ and $A\in \cA_\pi$ such that $A$ is invertible. 
The matrix $A$ admits a singular value decomposition
$$
A = UDV,
$$
where $U,V$ are unitary and $D={\rm diag}(\lambda_1,\lambda_2)$ with $\lambda_i>0$
such that $\lambda_i^2$ are the eigenvalues of $A^{\trans}A$.
Let $\TRect$ be the triangle of vertices $(0,0)$, $(0,\sqrt 2)$ and $(\sqrt 2,0)$. 
We define the triangle
$$
T :=  \sqrt{|\det A|} V^\trans D^{-1} (\TRect),
$$ 
which satisfies $|T| = |\TRect| = 1$ and has an angle of $\cPi/2$ at the origin so that $S(T) = 1$. 
Denoting by $C$ the constant in Lemma \ref{lemmaS} and using \iref{ineqPi}, we obtain
$$
(1+C)^{-1} \|\nabla \pi-\nabla \interp_T^{m-1} \pi \|_{L^\infty(T)} \leq \|\nabla \pi\|_{L^\infty(T)}
\leq \diam(A(T))^{m-1} 
=  |\det A|^{\frac {m-1} 2} \diam(\TRect)^{m-1}.
$$
Taking the infimum over all invertible $A\in \cA_\pi$
and remarking that this set is
dense in $\cA_\pi$, we conclude the proof of the left part of \iref{eqLA}.
For the right part, we define
for all $q_1,q_2\in \H_{m-1}$, and any triangle $T$,
\be
\label{defNormInf}
\|(q_1,q_2)\|_T := \inf_{r_1,r_2\in \sP_{m-2}} \|(q_1,q_2) - (r_1,r_2)\|_{L^\infty(T)}.
\ee
We denote by $\TEq$ an equilateral triangle centered at the origin and of area $1$. 
Since the functions $ \|\cdot \|_{\TEq}$ and $\|\cdot\|_{L^\infty(\TEq)}$ are norms on $\H_{m-1}\times \H_{m-1}$ there exists a constant $C_*$ such that $\|\cdot\|_{L^\infty(\TEq)}\leq C_* \|\cdot \|_{\TEq}$.
Let $T$ be a triangle satisfying $|T|=1$ and $\bary(T) = 0$. Then there exists a linear change of coordinates $\phi\in \SL_2$ such that $T=\phi(\TEq)$. We then obtain
$$
 \|(q_1,q_2)\|_{L^\infty(T)} = \|(q_1\circ\phi, q_2 \circ \phi)\|_{L^\infty(\TEq)} \leq C_* \|(q_1\circ\phi, q_2 \circ \phi)\|_{\TEq} =C_* \|(q_1,q_2)\|_T 
$$
We now choose a polynomial $\pi\in \H_m$ and we set $(q_1,q_2) := \nabla \pi$. It follows from the previous equation and \iref{defNormInf} that
$$
\|\nabla \pi\|_{L^\infty(T)} \leq C_* \|\nabla \pi\|_T \leq C_* \|\nabla \pi-\nabla \interp_T^{m-1} \pi\|_{L^\infty(T)}
$$
We define a linear map $A\in \GL_2$ associated to $\pi$ and $T=\phi(\TEq)$ as follows
$$
A := \|\nabla \pi\|_{L^\infty(T)}^{\frac 1 {m-1} } \ \lambda^{-1} \phi^{-1},
$$ 
where $\lambda = 3^{-3/4}$ is the minimal distance from $0$ to $\partial \TEq$.
Then for all $z\in \partial T$ we have
$\phi^{-1}(z)\in \partial \TEq$ and hence $|\phi^{-1}(z)|\geq \lambda$. Therefore, 
$$
|A(z)|^{m-1} =\|\nabla \pi\|_{L^\infty(T)}\ (\lambda^{-1} |\phi^{-1}(z)|)^{m-1} \geq \|\nabla \pi\|_{L^\infty(T)} \geq |\nabla\pi(z)|.
$$
By homogeneity, we thus find that for all $z\in\R^2$
$$
|\nabla\pi(z)|\leq |A(z)|^{m-1}, 
$$
which means that $A\in \cA_\pi$.  
Furthermore, since $\det \phi=1$, we have
$$
|\det A|^{\frac {m-1} 2 } = \lambda^{-(m-1)} \|\nabla \pi\|_{L^\infty(T)} \leq C_* \lambda^{-(m-1)} \|\nabla \pi-\nabla \interp_T^{m-1} \pi\|_{L^\infty(T)}.
$$
Hence taking the infimum over all triangles $T$ satisfying $|T|=1$ and $\bary(T) = 0$ we obtain
\be
\label{infbary}
\gamma_m(\pi)^{\frac {m-1} 2} \leq C \inf_{\substack{|T| = 1 \\ \bary(T) = 0}}  \|\nabla \pi-\nabla \interp_T^{m-1} \pi\|_{L^\infty(T)} 
\ee
with $C = C_* \lambda^{-(m-1)}$. Using the invariance of the interpolation error under translation, as expressed by  \iref{transInv}, we find that the right hand side of \iref{infbary} is $C L_m(\pi)$, which concludes the proof. 
\sq


We next introduce a measure of the isotropy of a triangle $T$
 with respect to the euclidean metric:
\be
\label{defRho}
\rho(T) := \frac{\diam(T)^2}{|T|}.
\ee
If $T$ is an obtuse triangle, an elementary computation shows that $4 S(T)\leq \rho(T)$. 
Indeed, if the largest angle of $T$ is $\theta\geq \cPi/2$, and if the edges neighbouring the angle $\theta$ have length $l_1, l_2$, we obtain using $l_1^2+l_2^2\geq 2 l_1 l_2$ that 
$$
\rho(T) = \frac {l_1^2+l_2^2 - 2 l_1 l_2\cos\theta}{ \frac 1 2 l _1 l_2 \sin\theta} \geq 4\frac {1-\cos\theta}{\sin \theta} = 4 \tan\frac \theta 2 = 4 S(T)
$$
Since the minimal value of $\rho$ is $4/\sqrt 3$ (for the equilateral triangle), and since $S(T)=1$ for acute triangles, we obtain that for any triangle $T$ 
\be
\label{ineqRhoS}
\rho(T) \geq \frac {4 }{\sqrt 3} S(T).
\ee
The functions $S$ and $\rho$ have a different behavior : $\rho(T)$ increases as $T$ becomes thinner, while $S(T)$ increases only if an angle of $T$ approaches $\cPi$.

In the follow up of this chapter, we frequently distort the measure of isotropy $\rho$ by a linear transform. If $A\in \GL_2$, then $\rho(A(T))$ reflects the isotropy of $T$ \emph{measured in the metric} $A^\trans A$. In particular $\rho(A(T))$ is minimal, i.e. equal to $4/\sqrt 3$, if and only if the  ellipse $\cE$ containing $T$ and of minimal area is of the form 
$$
\cE=\{z\in \R^2\sep |A(z-\bary(T))|\leq r\}
$$ 
for some $r>0$.


We may now state the main theorem of this section.
\begin{theorem}
\label{thLocal}
There exists a constant $C=C(m)$ such that for all $\pi\in \H_m$, all $A \in \cA_\pi$ 
and any triangle $T$, we have
\be
\label{localPi}
|\pi - \interp^{m-1}_T \pi|_{W^{1,p}(T)} \leq C |T|^{\frac 1 \tau} \ S(T) \ \rho(A (T))^{\frac{m-1} 2}  |\det A |^{\frac {m-1} 2}
\ee
where $\frac 1 \tau := \frac {m-1} 2 +\frac 1 p$.
Furthermore
 for any triangle $T$ and any $g\in C^m(T)$, we have
\be
\label{localIsoW1P}
|g- \interp^{m-1}_T g|_{W^{1,p}(T)} \leq C |T|^{\frac 1 \tau} \ S(T) \ \rho(T)^{\frac{m-1} 2} \|d^m g\|_{L^\infty(T)},
\ee
where $\|d^m g\|_{L^\infty(T)}$ is defined by \iref{normdmgW1P}.
\end{theorem}

Before proving this result, we make some observations on its consequences.
Combining the two estimates contained in this theorem, 
we obtain a mixed anisotropic-isotropic estimate, that can 
be used as a guideline for producing triangulations adapted to a function $f\in C^m(\Omega)$.
Let $T$ be a triangle, let $f\in C^m(T)$, $\pi \in \H_m$ and $A\in \cA_\pi$. Then
\be
\label{localMixed}
\begin{array}{cl}
\vspace{1mm}
&|f - \interp^{m-1}_T f|_{W^{1,p}(T)} \\
\leq&   |\pi - \interp^{m-1}_T \pi|_{W^{1,p}(T)} + |(f-\pi) - \interp^{m-1}_T (f-\pi)|_{W^{1,p}(T)} ,\\
\leq& C |T|^{\frac 1 \tau} S(T) \left( \rho(A(T))^{\frac{m-1} 2} |\det A|^{\frac {m-1} 2}  +  \rho(T)^{\frac{m-1} 2} \|d^m f-d^m\pi\|_{L^\infty(T)} \right),
\end{array}
\ee
where $C = C(m)$. Note that the left term in the parenthesis is an ``anisotropic'' contribution to the error, while the right term is an ``isotropic'' contribution.

Let $\ve>0$ and $1\leq p<\infty$. We now explain how the requirements (i), (ii), (iii) and  (iv) heuristically exposed in the introduction can be mathematically stated, 
and show that 
the estimate
$$
\# (\cT)^{\frac{m-1} 2}|f - \interp^{m-1}_\cT f|_{W^{1,p}(\Omega)} \leq C\left\|L_m(\pi_z)+\ve\right\|_{L^\tau(\Omega)},
$$
is met 
when the triangulation satisfies these requirements.
Consider a polygonal and bounded domain $\Omega$, 
a function $f\in C^m(\overline \Omega)$
and a triangulation $\cT$. 
For each $z\in \Omega$, we denote by $T_z\in \cT$ the triangle containing $z$ and define $\pi_z = \frac{d^mf(z)}{m!}\in \H_m$.  The adaptation of $\cT$ with respect to $f$ for the $W^{1,p}$ semi-norm,
can be measured by the smallest constant $C_\cT\geq 1$ such that the following criterions are met:

\begin{enumerate}[(i)]
\item (Equilibrated errors) There exists a constant $\delta>0$ such that for all $z\in \Omega$,
\be
\label{eqError}
C_\cT^{-1}\delta\leq |T_z|^{\frac 1 \tau} (L_m(\pi_z)+\ve)\leq C_\cT \delta.
\ee
\item (Optimized shapes) 
For all $z\in \Omega$, there exists $A_z\in \cA_{\pi_z}$, such that 
\be
\label{optShape}
\rho(A_z(T_z)) \leq C_\cT \text{ and } |\det A_z|^{\frac {m-1} 2} \leq C_L (L_m(\pi_z)+\ve), 
\ee
where $C_L$ is the constant that appears in Lemma \iref{lemmaA}.
According to this lemma, such an $A_z$ always exists for any $\ve>0$.
\item (Bounded sliverness in average) The averaged $l^p(\cT)$ norm of $S$ is bounded as follows
\be
\label{boundedS}
\left(\frac 1 {\#(\cT)} \sum_{T\in \cT} S(T)^p\right)^{\frac 1 p}\leq C_\cT .
\ee
This condition is less stringent than asking that $S(T)\leq C_\cT$ for {\it all} $T\in\cT$,
and turns out to be sufficient for proving the optimal error estimate.

\item (Sufficient refinement)  The mesh $\cT$ is sufficiently fine in such way that the local interpolation error estimate\iref{localMixed} is controlled by the ``anisotropic'' component. More precisely, for all $z\in \Omega$, 
\be
\label{sufficientRefinement}
\rho(T_z)^{\frac{m-1} 2} \|d^m f-d^m f(z)\|_{L^\infty(T_z)} \leq C_\cT (L_m(\pi_z)+\ve).
\ee
This condition is ensured by sufficient refinement of the triangulation due to the following 
observation: If $T'_z$ is the image of $T_z$ by a homothetic size reduction around $z$, then $\rho(T'_z) = \rho(T_z)$ while $\|d^m f-d^m f(z)\|_{L^\infty(T'_z)}$ tends to zero due to the continuity of $d^m f$. 
\end{enumerate}

We now produce a global error estimate from these four assumptions.
For a given $z\in \Omega$ we inject successively $\pi = \pi_z$, \iref{optShape}, \iref{sufficientRefinement} and \iref{eqError} into the estimate \iref{localMixed} and obtain 
\begin{eqnarray*}
& &|f - \interp^{m-1}_{T_z} f|_{W^{1,p}(T_z)} \\
&\leq & C |T_z|^{\frac 1 \tau} S(T_z) \left( \rho(A_z(T_z))^{\frac{m-1} 2} |\det A_z|^{\frac {m-1} 2}  +  \rho(T_z)^{\frac{m-1} 2} \|d^m f-d^m\pi_z\|_{L^\infty(T)} \right)\\
&\leq & C|T_z|^{\frac 1 \tau} S(T_z) (C_\cT^{\frac {m-1} 2} C_L (L_m(\pi_z)+\ve)+ \rho(T_z)^{\frac{m-1} 2} \|d^m f - d^m f(z)\|_{L^\infty(T)})\\
&\leq &  C|T_z|^{\frac 1 \tau} S(T_z) (C_\cT^{\frac {m-1} 2} C_L +  C_\cT)(L_m(\pi_z)+\ve)\\
&\leq & C_0 \delta S(T_z)
\end{eqnarray*}
where $C_0 = C_0(m,C_\cT,C_L)$. 
Using \iref{boundedS} we obtain
\be
\label{lowerDelta}
|f - \interp^{m-1}_\cT f|_{W^{1,p}(\Omega)}^p  = \sum_{T\in \cT} |f - \interp^{m-1}_T f|_{W^{1,p}(T)}^p \leq C_0^p \delta^p \sum_{T\in \cT} S(T)^p \leq (C_0 C_\cT)^p \delta^p \# (\cT) .
\ee
On the other hand, the left side of inequality \iref{eqError} provides an upper estimate of $\delta$ as follows.
\be
\label{upperDelta}
C_\cT^{-\tau} \delta^\tau \# (\cT)  = C_\cT^{-\tau} \delta^\tau \int_\Omega \frac {dz}{|T_z|} \leq \int_\Omega (L_m(\pi_z)+\ve)^\tau dz =\left\|L_m(\pi_z)+\ve\right\|_{L^\tau(\Omega)}^\tau.
\ee
Combining \iref{lowerDelta} with \iref{upperDelta} we eliminate the variable $\delta$ and obtain 
\be
\label{nonAsympt}
\# (\cT)^{\frac{m-1} 2}|f - \interp^{m-1}_\cT f|_{W^{1,p}(\Omega)} \leq C\left\|L_m(\pi_z)+\ve\right\|_{L^\tau(\Omega)},
\ee
where $C = C(m,C_\cT)$. Hence the optimal asymptotic estimate \iref{upperEstimW1P} is satisfied up to a multiplicative constant depending only on the degree $m-1$ of interpolation and the quality of the mesh reflected by $C_\cT$. 
Note however that the properties \iref{eqError}, \iref{optShape}, \iref{boundedS} and \iref{sufficientRefinement} required for $C_\cT$ may lead to a very pessimistic constant $C=C(m,C_\cT,C_L)$ in inequality \iref{nonAsympt}. 
Finer estimates and weaker conditions on the mesh $\cT$ can be obtained from \iref{localMixed}. 


In the context of the $H^1 = W^{1,2}$ semi-norm and of piecewise linear and quadratic elements we present numerical results in \S \ref{secW1P4p2} and discuss the quality of a numerical mesh $\cT$ using three quantities $\sigma(\cT)$, $\rho(\cT)$ and $S(\cT)$ that are related to the conditions (i), (ii) and (iii) respectively. 
We also discuss in \S \ref{secW1P4} a reformulation of the requirements of size \iref{eqError} and shape \iref{optShape} for the triangles $T$ of the mesh $\cT$ in terms of Riemannian metrics, a more convenient form for mesh generation. 

The construction of a mesh which satisfies both the requirements \iref{optShape} of optimized shapes and \iref{boundedS} of bounded measure of sliverness is a difficult problem. The construction presented in this chapter, for the proof of Theorems \ref{mainTheorem} and \ref{optiTheorem}, is based on a local patching strategy. A small portion of the triangulations, which can be neglected as the cardinality tends to infinity, does not satisfy conditions.
Another approach to this mesh generation problem is presented in Chapter 5, where the requirements  \iref{boundedS} is replaced as follows : the measure of sliverness needs to be uniformly bounded on a refinement of  $\cT$. The role of the measure of sliverness in the $W^{1,p}$ approximation error is again discussed in Chapter 6.

\paragraph{Proof of Theorem \ref{thLocal} :} 
Let $T$ be a triangle and let $h\in C^1(T)$. Using lemma \ref{lemmaS}, we obtain
\be
\label{ineqH}
\begin{array}{ccc}
|h -\interp^{m-1}_T h|_{W^{1,p}(T)} & = & \|\nabla h - \nabla \interp^{m-1}_T h\|_{L^p(T)}\\
&\leq & |T|^{\frac 1 p} \|\nabla h - \nabla \interp^{m-1}_T h\|_{L^\infty(T)}\\ 
&\leq & |T|^{\frac 1 p} (1+CS(T)) \|\nabla h \|_{L^\infty(T)}.\\ 
\end{array}
\ee
Replacing $h$ with $\pi$ in inequality \iref{ineqH} and combining it with \iref{ineqPi}, we 
obtain that if $T$ contains the origin, then for all $A\in\cA_\pi$
$$
|\pi -\interp^{m-1}_T \pi|_{W^{1,p}(T)} \leq |T|^{1/p} (1+CS(T))\diam(A(T))^{m-1}.
$$
The left and right quantities in the above inequality are invariant by 
translation of $T$ and therefore this inequality remains valid for any $T$.
Combining it with the identity
$$
\diam(A(T))^2=|T| \, |\det A| \,\rho(A(T)),
$$
this leads to the first inequality \iref{localPi} of Theorem \ref{thLocal}. For the second
inequality, we take $g\in C^m(T)$ and $z_0 = (x_0,y_0)\in T$. We now take for $h$ 
the remainder of the Taylor development of $g$ at $z_0$,
$$
h(x,y) := g(x,y) - \sum_{ k+l\leq m-1} 
\frac{\partial^{k+l} g}{\partial x^k \partial y^l}(z_0) \frac{(x-x_0)^k}{k!} \frac{(y-y_0)^l}{l!}.
$$ 
Therefore $h\in C^m(T)$ and  
\be
\label{dh0}
h(z_0) = dh(z_0) = \cdots = d^{m-1} h(z_0) = 0.
\ee
It follows that
\be
\label{ineqH0}
\|\nabla h\|_{L^\infty(T)} \leq C_1 (\diam T)^{m-1} \|d^m h\|_{L^\infty(T)} =  C_1 |T|^{\frac{m-1} 2} \rho(T)^{\frac{m-1} 2}  \|d^m h\|_{L^\infty(T)},
\ee
where $C_1=C_1(m)$. Combining \iref{ineqH} and \iref{ineqH0}, we obtain 
\be
\label{localIsoH}
|h -\interp^{m-1}_T h|_{W^{1,p}(T)} \leq C |T|^{\frac 1 \tau} S(T) \rho(T)^{\frac {m-1} 2}\|d^m h\|_{L^\infty(T)},
\ee
where $C=C(m)$. We now observe that $g-h\in \P_{m-1}$, hence $d^mg = d^m h$ and $h -\interp^{m-1}_T h = g -\interp^{m-1}_T g$. Injecting this into the last equation we conclude the proof of inequality \iref{localIsoW1P} and of Theorem \ref{thLocal}
\sq
\nl
\noindent
As a conclusion to this section we prove inequality  \iref{equivLm}, which links the functions $L_{m,p}$
and $L_m=L_{m,\infty}$.

\begin{lemma}
\label{lemmaEquivLm}
There exists a constant $c = c(m)>0$ such that for all $1\leq p_1 \leq p_2 \leq \infty$, 
\be
\label{compareL}
c L_m\leq L_{m,p_1}\leq L_{m,p_2} \leq L_m \text{ on } \H_m,
\ee
\end{lemma}
\proof
Let $\TEq$ be an equilateral triangle of area $1$. Since all norms are equivalent on the finite dimensional space $\P_{m-1}^2$, there exists a constant $c=c(m)>0$ such that for all $(q_1,q_2)\in \P_{m-1}\times \P_{m-1}$, 
\be
\label{ineqInf1}
c \|(q_1,q_2)\|_{L^\infty(\TEq)} \leq \|(q_1,q_2)\|_{L^1(\TEq)},
\ee
Furthermore, since $\TEq$ has area $1$, we have 
\be
\label{ineqp1p2}
 \|(q_1,q_2)\|_{L^{p_1}(\TEq)} \leq \|(q_1,q_2)\|_{L^{p_2}(\TEq)},
\ee
for all $1\leq p_1\leq p_2 \leq \infty$. If $T$ is a triangle satisfying $|T|=1$, there exists an affine change of coordinates $\Psi$ such that $T = \Psi(\TEq)$ and we have $\|(q_1,q_2)\|_{L^p(T)} = \|(q_1\circ \Psi, q_2\circ \Psi)\|_{L^p(\TEq)}$ for all $(q_1,q_2)\in \P_{m-1}\times \P_{m-1}$ and $1\leq p\leq \infty$. 
Combining this invariance property with inequalities \iref{ineqInf1} and \iref{ineqp1p2} we obtain
\be
\label{ineqQ1Q2}
c\|(q_1,q_2)\|_{L^\infty(T)} \leq \|(q_1,q_2)\|_{L^{p_1}(T)} \leq \|(q_1,q_2)\|_{L^{p_2}(T)}\leq \|(q_1,q_2)\|_{L^\infty(T)}.
\ee
We now choose a polynomial $\pi\in \H_m$, we set $(q_1,q_2) := \nabla \pi - \nabla\interp_T^{m-1} \pi$, and we take the infimum of \iref{ineqQ1Q2} among all triangles $T$ of area $1$.
This leads to the announced inequality \iref{compareL} which concludes the proof.
\sq

\section{Proof of Theorems \ref{mainTheorem} and \ref{optiTheorem}}
\label{secW1P3}
The polygonal domain $\Omega$, the integer $m$, the function $f\in C^m(\overline \Omega)$ and the exponent $1\leq p<\infty$ are fixed in this section which is devoted to the proof of the lower estimate \iref{lowerEstim} and the upper estimates \iref{upperEstimW1P} and \iref{upperEstimEpsW1P} which are stated in Theorems  \ref{mainTheorem} and \ref{optiTheorem}.

We denote by $\mu_{z_0}$ the Taylor polynomial of $f$ and of degree $m$ at the point $z_0 = (x_0,y_0)\in \Omega$ 
$$
\mu_{z_0}(x,y) := \sum_{k+l \leq m}\frac{\partial^{k+l} f(z_0)}{\partial x^k\, \partial y^l}\frac{(x-x_0)^k}{k!} \frac{(y-y_0)^l}{l!}.
$$
Note that $\pi_z$ is the homogeneous component of degree $m$
in $\mu_z$. Therefore $d^m \pi_z=d^m \mu_z = d^m f(z)$ for any $z\in \overline \Omega$, and for any triangle $T$  
\be
\label{muzPizW1P}
\pi_z- \interp_T^{m-1} \pi_z = \mu_z -\interp_T^{m-1} \mu_z.
\ee

\subsection{Proof of the lower estimate \iref{lowerEstim}}
\label{secW1P3p1}
The following lemma allows to bound by below the interpolation error of $f$ on a triangle $T$. 
\begin{lemma}
\label{lemmaLowerW1P}
Let $\frac 1 \tau := \frac {m-1} 2 +\frac 1 p$. For any triangle $T\subset \Omega$ and $z\in T$ we have
$$
|f- \interp^{m-1}_T f|_{W^{1,p}(T)} \geq |T|^{\frac 1 \tau} \left(L_{m,p}(\pi_z) - \omega(\diam T) \rho(T)^{\frac{m-1} 2} S(T)\right),
$$
where the function $\omega$ is positive, depends only on $f$ and $m$, and satisfies $\omega(\delta) \to 0$ as $\delta\to 0$.
\end{lemma}

\proof
Let $h:= f-\mu_z$. Using Equation \iref{muzPizW1P} we obtain
\begin{eqnarray*}
|f- \interp^{m-1}_T f|_{W^{1,p}(T)} & \geq & |\pi_z - \interp^{m-1}_T \pi_z|_{W^{1,p}(T)} - |h - \interp^{m-1}_T h|_{W^{1,p}(T)}\\
& \geq &  |T|^{\frac 1 \tau} L_{m,p}(\pi_z) - |h - \interp^{m-1}_T h|_{W^{1,p}(T)}.
\end{eqnarray*}
and we have seen in Theorem \ref{thLocal} that
$$
 |h - \interp^{m-1}_T h|_{W^{1,p}(T)} \leq C_0 |T|^{\frac 1 \tau}  S(T) \rho(T)^{\frac {m-1} 2} \|d^m h\|_{L^\infty(T)}.
$$
for some constant $C_0>0$ depending only on $m$.
We then remark that
$$
\|d^m h\|_{L^\infty(T)} = \|d^m f - d^m\pi_z\|_{L^\infty(T)}  =  \|d^m f - d^m f(z)\|_{L^\infty(T)}.
$$
Therefore, defining 
$$
\omega(\delta) := C_0\sup_{z,z'\in \Omega\sep |z-z'|\leq \delta} \|d^mf(z) - d^mf(z')\|,
$$
we conclude the proof of this lemma.
\sq

We now consider an admissible sequence of triangulations $(\cT_N)_{N\geq N_0}$. For all $N\geq N_0$, 
$T\in \cT_N$ and $z\in T$, we define $\phi_N(z) := |T|$ and 
$$
\psi_N(z) := \left(L_{m,p}(\pi_z) - \omega(\diam(T)) \rho(T)^{\frac {m-1} 2} S(T)\right)_+,
$$
where $\lambda_+ := \max\{\lambda,0\}$.
Holder's inequality gives, with $\frac 1 \tau := \frac{m-1} 2 + \frac 1 p$,
\be
\label{holderPsiW1P}
\int_\Omega \psi_N^\tau \leq \left(\int_\Omega \phi_N^{\frac{(m-1)p} 2} \psi_N^{p}\right)^{\frac \tau p} \left(\int_\Omega \phi_N^{-1}\right)^{\frac {(m-1)\tau} 2}.
\ee

\noindent
Note that $\int_\Omega \phi_N^{-1} = \# \cT_N\leq N$. Furthermore if $T\in \cT_N$ and $z\in T$ then according to Lemma \ref{lemmaLowerW1P}
$$
 \phi_N(z)^{\frac{(m-1)p} 2} \psi_N(z)^p
=|T|^{\frac p \tau-1}  \psi_N(z)^p \leq\frac 1 {|T|} |f-\interp_T^{m-1} f|_{W^{1,p}(T)}^p,
$$
hence 
$$
\int_\Omega \phi_N^{\frac{(m-1)p} 2} \psi_N^{p}\leq \sum_{T\in \cT_N} \frac 1 {|T|}  \int_T  |f-\interp_T^{m-1} f|_{W^{1,p}(T)}^p= |f-\interp_T^{m-1} f|_{W^{1,p}(\Omega)}^p
$$
Inequality \iref{holderPsiW1P} therefore leads to
\be
\label{upperPsiW1P}
 \|\psi_N\|_{L^\tau(\Omega)} \leq |f- \interp^{m-1}_{\cT_N} f|_{W^{1,p}(\Omega)} N^{\frac {m-1} 2}.
\ee

\noindent
Since the sequence $\seqT$ is admissible, there exists a constant $C_A>0$ such that for all $N$ and all $T\in \cT_N$ we have $\diam(T)\leq C_AN^{-\frac 1 2}$. 
We introduce a subset of $\cT'_N\subset \cT_N$ which gathers the most degenerate triangles
$$
\cT'_N = \{ T\in \cT_N \sep \rho(T)\geq \omega(C_AN^{-\frac 1 2})^{\frac{-1}{m+1}}\},
$$
where $\omega$ is the function from Lemma \ref{lemmaLowerW1P}. 
We denote by $\Omega'_N$ the portion of $\Omega$ covered by $\cT'_N$.
For all $z\in \Omega\sm\Omega'_N$, recalling from \iref{ineqRhoS} that $\rho \geq S$, 
we obtain 
$$
\psi_N(z)\geq L_{m,p}(\pi_z) -\sqrt{\omega(C_A N^{-\frac 1 2})}.
$$
Hence
$$
\begin{array}{ll}
 \|\psi_N\|_{L^\tau(\Omega)}^\tau & \geq \left \|\left(L_{m,p}(\pi_z) -\sqrt{\omega(C_A N^{-\frac 1 2})}\right)_+\right\|_{L^\tau(\Omega\sm\Omega'_N)}^\tau\\
 & \geq 
 \left \|\left(L_{m,p}(\pi_z) -\sqrt{\omega(C_A N^{-\frac 1 2})}\right)_+\right\|_{L^\tau(\Omega)}^\tau
 -C^\tau |\Omega'_N|,
 \end{array}
 $$
where $C:=\max_{z\in\Omega}L_{m,p}(\pi_z)$. We next observe
that $|\Omega'_N|\to 0$ as $N\to +\infty$: indeed
for all $T\in \cT'_N$ we have 
$$
|T| = \diam(T)^2 \rho(T)^{-1} \leq C_A ^2 N^{-1} \omega(C_A N^{-\frac 1 2})^{\frac 1 {m+1}}.
$$
Since $\#\cT'_N\leq N$, we obtain $|\Omega'_N|\leq C_A^2 \omega(C_A N^{-\frac 1 2})^{\frac 1 {m+1}}$, which tends to $0$ as $N\to \infty$. We thus obtain
$$
\liminf_{N\to \infty} \|\psi_N\|_{L^\tau(\Omega)}
\geq \lim_{N\to \infty}   \left \|\left(L_{m,p}(\pi_z) -\sqrt{\omega(C_A N^{-\frac 1 2})}\right)_+\right\|_{L^\tau(\Omega)}
 = \|L_{m,p}(\pi_z)\|_{L^\tau(\Omega)}.
$$
Combining this result with \iref{upperPsiW1P} we conclude the proof of the announced estimate \iref{lowerEstim}.

\subsection{Proof of the upper estimates \iref{upperEstimW1P} and \iref{upperEstimEpsW1P}}
\label{secW1P3p2}
The proof of these upper estimates is based on an explicit construction of the triangulations $\cT_N$, which is adapted from the construction in \cite{BBLS}. 
Roughly speaking, the idea of this construction is to produce a first mesh $\cR$ of the domain $\Omega$, composed of elements sufficiently small so that $f$ can be regarded as a polynomial $\pi_R$ on each triangle $R\in \cR$. Each element $R\in \cR$ is then tiled with small triangles optimally adapted to $\pi_R$, and some technical manipulations are done in order to preserve the conformity at the interfaces of the elements of $\cR$.
The main difference with the construction first proposed in \cite{BBLS}, and used later in Chapter \ref{chapOptAniso}, is that the measure of sliverness $S$ of the generated triangles should be kept under control. 

Let $T$ be a triangle with vertices $(z_0,z_1,z_2)$.
We define the symmetrized
triangle $\ti T$ of vertices $(z_1,z_2,z_1+z_2-z_0)$ so that
$T\cup \ti T$ is a parallelogram.
We define a tiling $\cP_T$ of the plane $\R^2$ as follows
\be
\label{defPT}
\cP_T := \{\alpha (z_1-z_0)+\beta (z_2-z_0)+ T'\sep \alpha,\beta\in \ZZ,  T'\in \{T,\ti T\}\}.
\ee
A homogeneous polynomial $\pi\in \H_m$ is either even or odd (depending on the parity of $m$). Combining this observation with the translation invariance \iref{transInv} we obtain that $|\pi-\interp_{T'}^{m-1} \pi|_{W^{1,p}(T')}$ is constant among all triangles $T'\in \cP_T$. 
We also define
\be
\label{defPTN}
\cP_{T,n}:=\frac 1 n \cP_T
\ee
the tiling obtained by rescaling $\cP_T$
by a factor $\frac 1 n$. We use this rescaled tiling in order to subdivide
an arbitrary triangle $R$, up to a few additional triangles located near the
boundary of $R$, as expressed by the following lemma.
\begin{lemma}
\label{lemmaTile}
Let $R$ and $T$ be two triangles. There exists a family $(\cP_{T,n}(R))_{n\geq 0}$, 
of conforming triangulations of $R$ such that the following holds
\begin{enumerate}
\item Nearly all the elements of $\cP_{T,n}(R)$ belong to $\cP_{T,n}$, which is defined by \iref{defPTN}, in the sense that
\be
\label{limCardRN}
\lim_{n\to \infty} \frac {\# (\cP_{T,n}^1(R))} {n^2} = \frac {|R|}{|T|} \ \text{ and } \  \lim_{n\to \infty}  \frac {\# (\cP_{T,n}^2(R))} {n^2}   = 0.
\ee
where
\be
\label{defP12}
\cP_{T,n}^1(R) := \cP_{T,n}(R)\cap \cP_{T,n} \ \text{ and } \quad \cP_{T,n}^2(R) := \cP_{T,n}(R)\sm \cP_{T,n}
\ee

\item The vertices of $\cP_{T,n}(R)$ on the boundary of $R$ are exactly those of the form $\frac k n a+ (1-\frac k n) b$, where $0\leq k \leq n$ and $a,b$ are vertices of $R$.
\item There exists constants $C_1 = C_1(R,T)$ and $C_2=C_2(R,T)$ such that 
\be
\label{admiRN}
\sup_{n\geq 0}  \; \( n \max_{T'\in \cP_{T,n}(R)} \diam(T') \)\leq C_1
\ \text{ and } \ \sup_{n\geq 0} \;  \max_{T'\in \cP_{T,n}(R)} S(T')\leq C_2.
\ee
\end{enumerate}
\end{lemma}

\proof
See appendix.
\sq

For any $M>0$, we define the compact set of triangles
$$
\mT_M := \{T  \sep |T| = 1, \ \diam(T)\leq M \text{ and } \bary(T) = 0\}. 
$$
Note that for all $T\in \mT_M$,
$$
\rho(T)\leq M^2.
$$
We also define the function
\be
\label{defLM}
L_M(\pi) := \min_{T\in \mT_M} |\pi - \interp_T^{m-1} \pi|_{W^{1,p}(T)}.
\ee
Since $\mT_M$ is compact for the Hausdorff distance
between sets and since $T\mapsto |\pi - \interp_T^{m-1} \pi|_{W^{1,p}(T)}$
is continuous with respect to this distance on the set of all triangles, we find
that this minimum is indeed attained and that $L_M$ is continuous. 
We also observe that $M \mapsto L_M(\pi)$ is a decreasing function of $M$ 
and that
$$
\lim_{M\to \infty}L_M(\pi)=
 \inf_{|T|=1, \bary(T) = 0} |\pi - \interp_T^{m-1} \pi|_{W^{1,p}(T)}
 = \inf_{|T|=1} |\pi - \interp_T^{m-1} \pi|_{W^{1,p}(T)}
 =L_{m,p}(\pi),
 $$
 where we have used the invariance under translation of the interpolation error \iref{transInv}
 for the second equality.

The constant $M>0$ is now fixed until the last step of this proof.
Let $\pi \in \H_m$ and let $T\subset \Omega$ be homothetic to a triangle achieving the minimum in
the definition of $L_M(\pi)$. Then, 
\be
\label{eqHomoth}
{\everymath{\displaystyle\everymath{}}
\begin{array}{ccl}
\vspace{0.5ex}
|f-\interp_T^{m-1} f|_{W^{1,p}(T)} &\leq &|\pi - \interp_T^{m-1} \pi|_{W^{1,p}(T)} + |(f-\pi) -\interp_T^{m-1} (f-\pi) |_{W^{1,p}(T)}\\
\vspace{0.5ex}
&\leq & |T|^{\frac 1 \tau} L_M(\pi) + C |T|^{\frac 1 \tau} \rho(T)^{\frac {m-1} 2}  S(T) \| d^m f -d^m \pi\|_{L^\infty(T)}\\
&\leq &  |T|^{\frac 1 \tau} ( L_M(\pi) +CM^{ m+1} \| d^m f -d^m \pi\|_{L^\infty(T)}),
\end{array}
}
\ee
where we have used inequality \iref{localIsoW1P} in the second line,
and in third line the fact that 
$$
\rho(T)^{\frac {m-1} 2}  S(T)\leq \rho(T)^{\frac {m+1} 2} \leq M^{m+1},
$$
since $S\leq \rho$ and $T$ is homothetic to an element of $\mT_M$. 

Let $\delta >0$ which value will be specified later. Since $d^m f$ is continuous, we can choose a sufficiently fine mesh $\cR = \cR(M,\delta)$ of $\Omega$
 in such way that, 
\be
\label{epsChoice}
C M^{m+1} \| d^m f(x) -d^m f(y)\|_{L^\infty(T)}\leq \delta, \text{ for all } R\in\cR \text{ and } x,y\in R.
\ee
For any triangle $R\in \cR$ we define
\be
\label{defZR}
z_R := \underset{z\in R}{\argmin} \, L_M(\pi_z) \ \text{ and } \ \pi_R := \pi_{z_R}. 
\ee
We also define 
\be
\label{defTR}
T_R := (L_M(\pi_R)+\delta)^{-\frac \tau 2} T_*,
\ee
where $T_*\in \mT_M$ achieves the minimum in the definition of $L_M(\pi_R)$.
We denote by $\cP_n(R)=\cP_{T_R,n}(R)$ the triangulation 
of Lemma \ref{lemmaTile} built from the two triangles $R$ and $T_R$, and similarly $\cP_n^1(R)=\cP_{T_R,n}^1(R)$ and $\cP_n^2(R)=\cP_{T_R,n}^2(R)$.
We define for all $n$
the global mesh of $\Omega$
$$
\cT_n^{M,\delta}=\bigcup_{R\in\cR} \cP_n(R),
$$
which coincides with $\cP_n(R)$ on each $R\in \cR$. Since all the meshes $\cP_n(R)$ are conforming, and since $\cP_n(R)$ has by construction $n+1$ equispaced vertices on each edge of $R$, the mesh $\cT_n^{M,\delta}$ is also conforming. 
According to Equations \iref{limCardRN} and \iref{defZR}, we have
\be
\label{limCardTN}
\begin{array}{rcl}
\displaystyle\lim_{n\to \infty} \frac{\# \left(\cT_n^{M,\delta}\right)}{n^2} &=& 
\displaystyle\sum_{R\in \cR}\left( \lim_{n\to \infty}  \frac{\# (\cP_n(R))}{n^2} \right)\\
&=&\displaystyle \sum_{R\in \cR}  |R| (L_M(\pi_R)+\delta)^\tau\\
&\leq& \displaystyle \int_\Omega  (L_M(\pi_z)+\delta)^\tau dz.
\end{array}
\ee
For $T\in \cP_n^1(R)$, we combine \iref{eqHomoth}, \iref{epsChoice} and \iref{defTR} to obtain 
\be
\label{erTiled}
|f-\interp_T^{m-1} f|_{W^{1,p}(T)}\leq n^{-\frac 2 \tau}
\text{ for all } T\in \cP_n^1(R).
\ee
For $T\in  \cP_n^2(R)$, we invoke the isotropic estimate \iref{localIsoW1P} 
to obtain
\be
\label{erGarbage1}
\begin{array}{rcl}
|f-\interp_T^{m-1} f|_{W^{1,p}(T)} 
&\leq & C |T|^{\frac 1 p} S(T) \diam(T)^{m-1} \|d^m f\|_{L^\infty(\Omega)}\\
&\leq & CS(T) \diam(T)^{\frac 2 \tau} \|d^m f\|_{L^\infty(\Omega)}
\end{array}
\ee
where $C$ is the constant from \iref{localIsoW1P}. Using 
the third item in Lemma \ref{lemmaTile}, we find that
that there exists constants $C_1=C_1(M,\delta)$
and $C_2=C_2(M,\delta)$ such that
\be
\sup_{n\geq 0}\left(n \max_{T\in \cT_n^{M,\delta}} {\rm diam}(T)\right)\leq C_1\;\; {\rm and}
\;\; \sup_{n\geq 0} \left(\max_{T\in \cT_n^{M,\delta}} S(T)\right)\leq C_2,
\label{admiT}
\ee
so that, combining with \iref{erGarbage1}, we have
for all $T\in \cP_n^2(R)$
\be
 \label{erGarbage}
 |f-\interp_T^{m-1} f|_{W^{1,p}(T)}^p  \leq C_0 n^{-\frac{2} \tau},
 \ee
with $C_0=C_0(M,\delta)$.
Combining \iref{erTiled} and \iref{erGarbage}, and using the
first item in Lemma \ref{lemmaTile}, we obtain
\begin{eqnarray*}
 |f-\interp_{\cT_n^{M,\delta}}^{m-1} f|_{W^{1,p}(\Omega)}^p &=& \sum_{T \in \cT_n^{M,\delta}} |f-\interp_T^{m-1} f|_{W^{1,p}(T)}^p\\
 &\leq&  \sum_{R\in \cR} \left(\#(\cP_n^1(R)) n^{-\frac {2p} \tau} + \#(\cP_n^2(R)) C_1^p n^{-\frac {2p} \tau}\right)
\end{eqnarray*}
therefore
\begin{eqnarray*}
\limsup_{n\to \infty} n^{\frac {2p} \tau -2}  |f-\interp_{\cT_n^{M,\delta}}^{m-1} f|_{W^{1,p}(\Omega)}^p 
& \leq &  \sum_{R\in \cR} \lim_{n\to \infty} \frac{\#(\cP_n^1(R))+ \#(\cP_n^2(R)) C_1^p}{n^2}\\
& = & \sum_{R\in \cR} |R|(L_M(\pi_R)+\delta)^\tau\\
&\leq &\int_\Omega (L_M(\pi_z)+\delta)^\tau dz.
\end{eqnarray*}
Combining this with \iref{limCardTN} we obtain 
\be
\label{limMD}
\limsup_{n\to \infty} \# (\cT_n^{M,\delta})^{\frac{m-1} 2} \left|f-\interp_{\cT_n^{M,\delta}}^{m-1} f\right|_{W^{1,p}(\Omega)} \leq \left\| L_M(\pi_z)+\delta\right\|_{L^\tau(\Omega)}.
\ee
Let $\ve >0$. Since
$$
\lim_{M\to \infty} \lim_{\delta\to 0} \| L_M(\pi_z)+\delta\|_{L^\tau(\Omega)} = \lim_{M\to \infty} \| L_M(\pi_z)\|_{L^\tau(\Omega)} = \| L_{m,p}(\pi_z)\|_{L^\tau(\Omega)}
$$
we can choose adequately $M$ and $\delta$ in such way that $\| L_M(\pi_z)+\delta\|_{L^\tau(\Omega)}\leq \| L_{m,p}(\pi_z)\|_{L^\tau(\Omega)}+ \ve$. 

Let $n=n(N,M,\delta)$ be the largest integer 
such that $\# (\cT_{n}^{M,\delta})\leq N$, we define
$$
\cT_N^\ve := \cT_{n}^{M,\delta}.
$$
so that $N^{-1}\#( \cT_N^\ve)\to 1$ as $N\to \infty$. 
If follows from \iref{limCardTN} and \iref{admiT} that the sequence of triangulations $(\cT_N^\ve)$ is admissible, and inequality \iref{limMD} gives 
$$
\limsup_{N\to \infty} N^{\frac {m-1} 2}  |f-\interp_{\cT_N^\ve}^{m-1} f|_{W^{1,p}(\Omega)} \leq  \| L_{m,p}(\pi_z)\|_{L^\tau(\Omega)}+ \ve.
$$
which is the upper estimate \iref{upperEstimEpsW1P} announced.
Last we choose for all $N$ large enough $\ve(N)>0$ such that 
$$
N^{\frac{m-1} 2} |f-\interp_{\cT_N^{\ve(N)}}^{m-1} f|_{W^{1,p}(\Omega)} \leq  \| L_{m,p}(\pi_z)\|_{L^\tau(\Omega)}+ 2\ve(N).
$$
and such that $\ve(N)\to 0$ as $N\to \infty$. The sequence of triangulations $\cT_N := \cT_N^{\ve(N)}$ fulfills the estimate \iref{upperEstimW1P} which concludes the proof.

\section{Optimal metrics for linear and quadratic elements}
\label{secW1P4}
The proof of the upper estimate \iref{upperEstimEpsW1P} exposed in the previous section involves
the construction of meshes $\cT_N^\ve$ by tiling each element $R$ of the ``coarse'' triangulation
$\cR$ using the finer mesh $\cP_n(R)$. In practice, such a construction may require a very large
number of triangles in order to match the optimal error estimate. More commonly used strategies
for mesh generation are based on the prescription of a non-euclidean metric depending on $f$ 
for which each triangle should be
isotropic. In this section, we explain how to design such metric in order to
derive near-optimal error estimates
and we give analytic expressions 
in the particular case of $\P_1$ and $\P_2$ finite elements.

\subsection{Optimal metrics}
\label{secW1P4p1}
As a first step, we express the requirements (i), (ii) and (iv) of mesh adaptation
in terms of metrics. 
We therefore use the following notations :
we consider a polygonal domain $\Omega$, an integer $m\geq 2$, an exponent $1\leq p<\infty$, and a function $f\in C^m(\Omega)$ to be approximated  in the $W^{1,p}$ semi-norm by $\P_{m-1}$ finite element interpolation on a triangulation of $\Omega$. We also consider two real numbers $\ve>0$ and $\delta>0$.


We define for all $\pi \in \H_m$
\be
\label{defA'}
\cA'_\pi := \{M\in S_2^+\sep |\nabla \pi(z)|^2\leq (z^\trans M z)^{m-1} \text{ for all } z\in \R^2\} = \{A^\trans A\sep A \in \cA_\pi\},
\ee
and we consider a continuous field $\cM : \Omega \to S_2^+$ of symmetric positive definite matrices satisfying $\cM(z)\in \cA'_{\pi_z}$ for all $z\in \Omega$ and 
\be
\label{defMz}
 C_2^{-1} (L_m(\pi_z)+\ve) \leq  (\det \cM(z))^{\frac {m-1} 4} \leq C_2 (L_m(\pi_z)+\ve),
\ee
where $C_2$ is an absolute constant. The existence of such a field is established in full generality in Chapter 6.
We explain in the sequel of this section a practical construction in the case of piecewise linear and bilinear elements. 
We then define a field of symmetric positive definite matrices $h$ on $\Omega$ by
\be
\label{optMetric}
h(z) := \delta^{-\tau}(\det \cM(z))^{\frac{-\tau}{2p}} \cM(z)
\ee
where $\frac 1 \tau := \frac {m-1} 2 + \frac 1 p$.
Such a field $h$ is called a Riemannian metric. 
Under some assumptions on the metric $h$ and on the domain $\Omega$, which are discussed in \cite{Shew, Bois} and Chapter 5 for the infinite domain $\R^d$, it is possible to produce a triangulation $\cT$ of $\Omega$ satisfying for all $T\in \cT$ and $z\in T$
\be
\label{adaptTh}
C_1^{-1} \leq |T|\sqrt{\det h(z)}\leq C_1 \ \text{ and } \ \rho\left(\sqrt{h(z)}(T)\right)\leq C_1
\ee
where the constant $C_1\geq 1$ reflects the quality of the adaptation of the mesh to the metric $h$.  
(In the second inequality the square root is meant in the sense of symmetric positive matrices).
Examples of such mesh generators are \cite{FreeFem,Inria,Peyre}. 
We shall not discuss in this chapter the conditions under which such a mesh can be generated. Let us only mention that, if one ignores a few outliers at the corners of $\Omega$, these conditions hold if $\delta$ is small enough.

Note that $(\det h(z))^{\frac 1 {2\tau}} = \delta^{-1} (\det \cM(z))^{\frac {m-1} 4}$, therefore if \iref{adaptTh} holds we find that for all $T\in \cT$ and $z\in T$,
$$
(C_1C_2)^{-1} \delta \leq |T|^{\frac 1 \tau}( L_m(\pi_z)+\ve)\leq C_1C_2 \delta,
$$
hence condition (i) of equilibrated errors, as stated in \iref{eqError}, holds provided $C_\cT\geq C_1C_2$. Furthermore for all $z\in \Omega$ let us define $A_z := \sqrt{\cM(z)}$ and note that $A_z\in \cA_{\pi_z}$ and $\det A_z = \sqrt{\det \cM(z)}$. Using \iref{defMz} and \iref{adaptTh} we find that condition (ii) of optimal shapes, as stated in \iref{optShape}, holds provided $C_\cT\geq C_2$.
Condition (iv) holds when the mesh $\cT$ is sufficiently refined, which is the case if $\delta$ is small enough.

In summary, given a map $\cM : \Omega\to S_2^+$ satisfying $\cM(z)\in \cA'_{\pi_z}$, \iref{defMz} and such that $\cM(z)$ is positive definite, state of the art mesh generators allow us to build triangulations
$\cT$ that match the conditions (i), (ii) and (iv). In order to prove the
near-optimal estimate 
$$
\#(\cT)^{\frac {m-1} 2} |f-\interp^{m-1}_\cT|_{W^{1,p}(\Omega)} \leq C \|L_m(\pi_z)+\ve\|_{L^\tau(\Omega)},
$$
it is also necessary that the generated meshes satisfy condition (iii) of bounded measure of sliverness, as stated in \iref{boundedS}. Unfortunately, the author has not heard of theoretical results 
that would guarantee this condition when a mesh is built by such algorithms, apart from Theorems \ref{thEquiv} and \ref{thPer} which only apply to the infinite domain $\R^2$ or the periodic domain $(\R/\Z)^2$ respectively.
We discuss in  \S \ref{secW1P4p2} the observed behaviour of $S(T)$ when
using the mesh generation software \cite{FreeFem}.

For $m \in\{ 2,3\}$, which correspond to $\P_1$ and $\P_2$ elements, we 
give in the sequel a simple expression of a continuous 
map $\cM_m : \H_m\to S_2^\oplus$ satisfying $\cM_m(\pi)\in \cA'_\pi$ for all $\pi \in \H_m$ and 
\be
\label{defMm}
K^{-1} L_m(\pi) \leq (\det \cM_m(\pi))^{\frac {m-1} 4} \leq K L_m(\pi)
\ee
for some absolute constant $K\geq 1$. 
It is not hard to build from $\cM_m(\pi_z)$ a matrix $\cM(z)$ satisfying \iref{defMz}.
For practical uses one usually takes
$$
\cM(z) := \cM_m(\pi_z) +\mu \Id.
$$
where the constant $\mu\geq 0$ is here to avoid degeneracy problems. 


Let us mention that there exists radically different approaches to anisotropic mesh generation, which are not based on Riemannian metrics. For example the hierarchical refinement procedure exposed in  Chapter \ref{chapCDHM}, which is proved in Chapter 8 to yield the best possible estimate \iref{optiaffineW1P} in the case of piecewise linear interpolation of bidimensional convex functions with the error measured in $L^p$ norm. This approach does not seem to adapt well to the $W^{1,p}$ norm: the main problem arises again from condition (iii) of bounded measure of sliverness, and from the lack of conformity of the triangulations generated by this procedure.

\subsection{The case of linear and quadratic elements}
\label{subSecOptMet}
\label{secW1P4p2}

We now give analytic expression of matrix fields $\cM_2$ and $\cM_3$
satisfying \iref{defMm}, which correspond to linear and quadratic elements.
In the simplest and already well established case of $\P_1$ elements, 
a more detailed analysis can be found in \cite{Shew2}. 
In contrast, the results for quadratic elements are new.
\nl
\nl
For any homogeneous quadratic polynomial $\pi\in \H_2$, $\pi = a x^2 + 2 b x y + c y^2$, we define the symmetric matrix
$$
 [\pi] = 
\left(\begin{array}{cc}
a & b\\
b & c
\end{array}\right).
$$
We define 
\be
\label{defMP1}
\cM_2(\pi) := 4 [\pi]^2 = 
4 \left(\begin{array}{cc}
a & b\\
b & c
\end{array}\right)^2
= 
4 \left(\begin{array}{cc}
a^2+b^2 & ab+bc\\
ab+bc & b^2+c^2
\end{array}\right)
\ee
For all $z\in \R^2$ one has $\nabla \pi(z) = 2 [\pi] z$, and therefore $|\nabla \pi(z)| ^2 = z^\trans \cM_2(\pi) z$. It follows that $\cM_2(\pi)\in \cA'_\pi$ and 
$$
\det \cM_2(\pi) =  \inf\{\det M\sep M\in \cA'_\pi\}
$$
which implies \iref{defMm} according to Lemma \ref{lemmaA}. 
\nl
\nl
For any $\pi\in \H_3$ we define
\be
\label{defMP2}
\cM_3(\pi) := \sqrt{[\partial_x \pi]^2+[\partial_y \pi]^2} 
\ee
If $\pi = a x^3+ 3 b x^2 y+ 3 c x y^2 + d y^3$ 
then in terms of the coefficients of $\pi$
\begin{eqnarray*}
\cM_3(\pi)  &=& 
3 \sqrt{
\left(\begin{array}{cc}
a & b\\
b & c
\end{array}\right)^2
+
\left(\begin{array}{cc}
b & c\\
c & d
\end{array}\right)^2
}\\
&=& 3
\sqrt{
\left(\begin{array}{cc}
a^2+2b^2+ c^2 & ab+2 bc + cd\\
ab+2 bc + cd & b^2+2c^2+d^2
\end{array}\right)
}
\end{eqnarray*}

In the sense of symmetric matrices, we have
$$
\cM_3(\pi) = \sqrt{[\partial_x \pi]^2+[\partial_y \pi]^2} \geq \sqrt{[\partial_x \pi]^2} = |[\partial_x \pi]|,
$$
where we used the fact that the square root $\sqrt : S_2^\oplus\to S_2^\oplus$ is increasing.
It follows that 
$$
|\nabla \pi(z)| ^2 =  |\partial_x\pi(z)|^2 + |\partial_y\pi(z)|^2 \leq 2 (z^\trans \cM_3(\pi) z)^2,
$$
hence $\sqrt 2 \cM(\pi) \in \cA'(\pi)$.
Note that
\be
\label{detMPi}
\det \cM_3(\pi) = 9 \sqrt{(a^2+2b^2+c^2)(b^2+2c^2+d^2) - (ab+2bc+cd)^2}.
\ee
It remains to establish \iref{defMm}.
This point is postponed to \S \ref{secW1P5}, right after \iref{L3eq}, as we develop a general method for obtaning simple equivalents of the functions $L_m$. 
Let us finally mention the work \cite{Kuate} in which approximate solutions to the optimization problem $\inf \{\det M\sep M\in \cA'_\pi\}$ are obtained through numerical optimization. This approach works for general $m$ but is harder to use than the algebraic expressions of $\cM_2(\pi)$ and $\cM_3(\pi)$ given here. 

\subsection{Limiting the anisotropy in mesh adaptation}

The measure of non degeneracy of a triangle and of its image by a linear transform
can be linked by the following result.
\begin{prop}
There exists an absolute constant $c>0$ such that for any triangle $T$ and any $A\in \GL_2$,
\be
\label{rhoA}
\frac c {\rho(A(T))} \leq \frac{\rho(T)}{\|A\| \|A^{-1}\|} \leq  \rho(A(T)).
\ee
\end{prop}

\proof
We use in this proof the identity $|\det B| = \|B\| \|B^{-1}\|^{-1}$ which holds for all $B\in \GL_2$.
Let $T'$ be a triangle and let $A'\in \GL_2$, then
$$
\rho(A'(T')) = \frac{\diam(A'(T'))^2}{|A'(T')|}\leq \frac{\|A'\|^2}{|\det A'|} \frac{\diam(T')^2}{|T'|} = \|A'\|\|A'^{-1}\| \rho(T').
$$
with the particular choice $A'=A^{-1}$ and $T'=A(T)$ we obtain the right side of \iref{rhoA}.
Let $\TEq$ be an equilateral triangle of area $1$, and let $\mu$ be the diameter of the largest ball included in $\TEq$. Up to a translation on $\TEq$ we can assume that there exists $B\in \GL_2$ such that $T = B(\TEq)$.
We then have
\begin{eqnarray*}
\diam(T) \diam(A(T)) &=& 
\diam(B(\TEq))\diam(AB(\TEq))\\ 
&\geq& \mu^2 \|B\| \|AB\| \\
&\geq& \mu^2 \|B\| \|A\| \|B^{-1}\|^{-1} \\
&=& \mu^2 \|A\| |\det B|.
\end{eqnarray*}
Hence, since $|T|=|B(\TEq)| = |\det B|$, 
$$
\rho(T) \rho(A(T)) = \frac{\diam(T)^2 \diam(A(T))^2}{|T| |A(T)|} \geq \frac{(\mu^2 \|A\| |\det B|)^2}{|\det B|^2 |\det A|} = \mu^4 \|A\| \|A^{-1}\|
$$
which establishes the left part of \iref{rhoA} with $c=\mu^4 = \frac{2^4}{3^3}$.
\sq

A consequence of the above lemma
is that if $\cT$ is a mesh adapted to a metric $h$ in the
sense of \iref{adaptTh}, then for all $T\in \cT$ and $z\in T$ we have  
$$
c C_1^{-1}\sqrt{\|h(z)\|\|h(z)^{-1}\|} \leq \rho(T) \leq C_1 \sqrt{\|h(z)\|\|h^{-1}(z)\|}.
$$
The measure of non-degeneracy $\rho(T)$ is thus large when $h(z)$ is ill conditioned. Although this property is desirable in order to adapt to highly anisotropic features of the function $f$ to be approximated, excessive degeneracy can cause mesh generation problems, which are discussed in \S \ref{secW1P4p2}. 
In the following, we explain how to slightly modify the construction of $\cM_2$ and $\cM_3$
in order to control the value of $\rho(T)$.

According to \iref{optMetric} we have $\|h(z)\|\|h^{-1}(z)\| = \|\cM(z)\| \|\cM(z)^{-1}\|$, 
and thus 
$$
\rho(T)\leq C_1\sqrt{\|\cM(z)\| \|\cM(z)^{-1}\|}.
$$
This leads us to define for all $\alpha \geq 1$,
$$
\cA'_{\pi,\alpha} := \{M \in \cA'_\pi\sep \|M\| \|M^{-1}\|\leq \alpha ^2\}.
$$

Let $M\in S_2^+$, let $R$ be a rotation and let $\lambda\geq \mu \geq 0$ be the eigenvalues of $M$ in such way that 
$
M = R^\trans \diag(\lambda, \mu) R.
$
We define for any $\alpha\geq 1$
\be
\label{defMAlpha}
M^{(\alpha)} := 
 R^\trans \left(
 \begin{array}{ccc}
 \lambda & 0\\
 0 & \max(\lambda \alpha^{-2}, \mu)
 \end{array}
\right) R.
\ee
Clearly $M^{(\alpha)}\geq M$, and if $M\neq 0$ then $\|M^{(\alpha)}\| \|\left(M^{(\alpha)}\right)^{-1}\|\leq \alpha^2$. Hence for all $M\in \cA'_\pi$ we have $M^{(\alpha)}\in \cA'_{\pi,\alpha}$.

In the case of piecewise linear elements we therefore have $\cM_2^{(\alpha)}(\pi)\in \cA'_{\pi,\alpha}$ for all $\pi \in \H_2$, and one easily shows that $\det \cM_2^{(\alpha)}(\pi) = \inf \{\det M\sep M\in \cA'_{\pi,\alpha}\}$. This suggests that constructing $\cM(z)$ from $\cM_2^{(\alpha)}(\pi_z)$ instead of $\cM_2(\pi_z)$ leads to a near-optimal mesh adaptation to the function $f$, under the constraint $\rho(T)\leq C_1 \alpha$ for all triangles $T$ in the triangulation.
The following proposition implies the same in the case of piecewise quadratic finite elements.

\begin{prop}
Let $\pi\in \H_3$ and $\alpha \geq 1$. Then $\sqrt 2 \cM_3^{(\alpha)}(\pi)\in \cA'_{\pi, \alpha}$ and 
\be
\label{optRAM3}
\det \cM_3^{(\alpha)}(\pi) \leq K \inf \{\det M\sep M\in \cA'_{\pi,\alpha}\}
\ee
where the constant $K$ is independent of $\pi$ and $\alpha$.
\end{prop}

\proof
We already know that $\sqrt 2 \cM_3^{(\alpha)}(\pi)\in \cA'_{\pi,\alpha}$. If $\cM_3^{(\alpha)}(\pi) = \cM_3(\pi)$ then \iref{optRAM3} holds as a consequence of \iref{defMm} and Lemma \ref{lemmaA}. We therefore assume in the following that $\cM_3^{(\alpha)} (\pi)\neq \cM_3(\pi)$.
Let 
$$
\lambda_*(\pi) := \|\nabla\pi\|_{L^\infty(D)} = \sup_{|z|\leq 1} |\nabla\pi(z)|,
$$
where $D=\{z\in \R^2\sep |z|\leq 1\}$ is the unit disc of $\R^2$.
The largest ball inscribed in $\{ z\in\R^2 \sep |\nabla \pi(z)|\leq 1\}$ is $\lambda_*(\pi)^{-\frac 1 2} D$. 
Let $M\in \cA'_{\pi,\alpha}$ and let $\lambda_1\geq \lambda_2>0$ be its eigenvalues. 
The ellipse $\{z\in \R^2 \sep z^\trans M z\leq 1\}$ contains the ball $\lambda_1^{-\frac 1 2} D$, hence $\lambda_1\geq \lambda_*(\pi)$. Furthermore $\lambda_2\geq \alpha^{-2} \lambda_1$, hence
\be
\label{infDetAlpha}
\det M = \lambda_1 \lambda_2\geq \alpha^{-2}\lambda_1^2 \geq \alpha^{-2} \lambda_*(\pi)^2. 
\ee
Let $\lambda(\pi)$ be the largest eigenvalue of $\cM_3(\pi)$, and assume that $\pi =  a x^2+3 b x^2 y + 3 c x y^2+ d y^3$. We obtain from \iref{defMP2} that
$$
\lambda(\pi)\leq \sqrt{\Tr\cM_3(\pi)^2} = \sqrt{a^2+3 b^2+ 3c^2+ d^2}. 
$$
Since the norms $\|\nabla\pi\|_{L^\infty(D)}$ and $\sqrt{a^2+3 b^2+ 3c^2+ d^2}$ are equivalent on the vector space $\H_3$, there exists a constant $C_0>0$ independent of $\pi\in \H_3$ such that $\lambda(\pi) \leq C_0  \lambda_*(\pi)$. Since $\cM_3^{(\alpha)}(\pi) \neq \cM_3(\pi)$, the eigeinvalues of $\cM_3^{(\alpha)}(\pi)$ are $\lambda(\pi)$ and $\alpha^{-2} \lambda(\pi)$. Hence
$$
\det \cM_3^{(\alpha)}(\pi) = \alpha^{-2} \lambda(\pi)^2 \leq C_0^2 \alpha^{-2} \lambda_*(\pi)^2.
$$
Combining this with \iref{infDetAlpha} we conclude the proof, with $K = C_0^2$.
\sq

Let us finally mention that, although they are derived from the
coefficients of $\pi$, the maps $\pi \mapsto \cM_m(\pi)$ and $\pi \mapsto \cM_m^{(\alpha)}(\pi)$ 
for $m\in\{2,3\}$ are invariant under rotation,
and therefore not tied to the chosen system of
coordinate $(x,y)$, as expressed by the following result.

\begin{prop}
\label{propInvMRot}
For any $m\in\{2,3\}$, any $\pi\in \H_m$ and any unitary matrix $U\in \cO_2$, one has
$$
\cM_m(\pi\circ U) = U^\trans \cM_m(\pi) U.
$$
Furthermore, for any $\alpha \geq 1$ one has $\cM_m^{(\alpha)} (\pi\circ U) = U^\trans \cM_m^{(\alpha)} (\pi) U$.
\end{prop}

\proof
We only prove the invariance under unitary transformation of $\cM_3$, since the proof
for $\cM_2$ is elementary, as well at the result for $\cM_m^{(\alpha)}$.
Let $\pi\in \H_3$, let $D_x=[\partial_x \pi]$ and $D_y = [\partial_y \pi]$. Let 
$
U = 
\left(
\begin{array}{cc}
u_{11} & u_{12}\\
u_{21} & u_{22}
\end{array}
\right)
$
be unitary,
then 
$$
[\partial_x(\pi\circ U)] = u_{11} U^\trans D_x U + u_{12} U^\trans D_y U \quad \text{ and } \quad [\partial_y(\pi\circ U)] = u_{21} U^\trans D_x U + u_{22} U^\trans D_y U
$$
Hence
\begin{eqnarray*}
[\partial_x(\pi\circ U)]^2 +  [\partial_y(\pi\circ U)]^2 &=& (u_{11}^2+u_{21}^2) U^\trans D_x^2 U \\
&&+ (u_{11} u_{12}+u_{21} u_{22}) U^\trans (D_x D_y+ D_y D_x) U \\
&&+ (u_{12}^2+u_{22}^2) U^\trans D_y^2 U
\end{eqnarray*}
which equals $U^\trans D_x^2 U+U^\trans D_y^2 U$ since $U$ is unitary. Eventually
$$
\cM_3(\pi\circ U) = \sqrt{U^\trans D_x^2 U+U^\trans D_y^2 U} = U^\trans \sqrt{D_x^2+D_y^2} \,U = U^\trans \cM_3(\pi) U
$$
which concludes the proof.
\sq

\subsection{Numerical results}

\begin{figure}
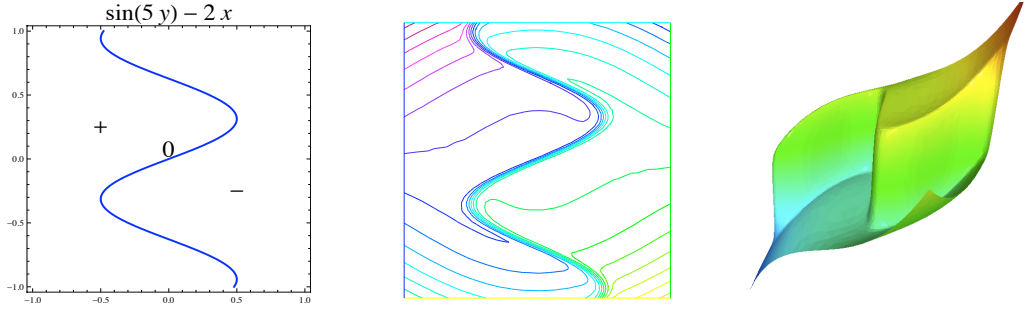

	\centering
		\includegraphics[width=4cm,height=4cm]{\pathPic/PaperW1P/SignSin.pdf}
		\hspace{-0.5cm}	
		\includegraphics[width=6.4cm,height=4cm]{\pathPic/PaperW1P/PicFFine.pdf}
		\hspace{-1cm}	
		\includegraphics[width=4cm,height=4cm]{\pathPic/PaperW1P/PicF3DFine.pdf}
	\caption{\label{fig2W1P}Description of the function $f_\delta$, $\delta = 0.1$.}
\end{figure}

The envisionned applications for the theory developped in this chapter are mainly in the field of partial differential equations that exhibit ``shocks'', and strongly anisotropic features, in particular conservation laws and fluid dynamics. 
We therefore test the quality of our meshes 
on a synthetic function that mimics the typical behavior 
of functions encountered in these contexts. For all $\delta>0$, 
our test function $f_\delta : [-1,1]^2 \to \R$ is defined as follows 
$$
f_\delta(x,y) := \tanh\left(\frac{2x - \sin(5y)} \delta\right) + x^3+x y^2.
$$
In all numerical results, we choose $\delta := 0.1$.
This function $f_\delta$, although smooth, exhibits a ``smoothed jump'' of height $2$ along to the curve defined by the equation $2x = \sin(5y)$,
on a layer of width $\delta$.
On the rest of the domain, $f_\delta$ is dominated by the polynomial part $x^3+x y^2$. The level lines and a 3D plot of $f_\delta$ are presented on the two rightmost pictures of Figure \ref{fig2W1P}.

Our purpose is to produce four triangulations $\cT_{H^1,\sP_1}$, $\cT_{H^1,\sP_2}$, $\cT_{L^2,\sP_1}$ and $\cT_{L^2,\sP_2}$ containing $2000$ triangles each and which, for this cardinality, produce respectilvely the smallest possible interpolation errors 
 $\|\nabla f-\nabla \interp_\cT^1 f\|_2$, $\|\nabla f-\nabla \interp_\cT^2 f\|_2$, $\|f-\interp_\cT^1 f\|_2$ and $\|f-\interp_\cT^2 f\|_2$.
It is clearly out of reach to find the triangulations leading exactly to the smallest error. Following the analysis developed in the beginning of this section we have generated $\cT_{H^1,\sP_1}$ and $\cT_{H^1,\sP_2}$ based on the metrics
\be
\label{defHNum}
\begin{array}{cclcc}
h_{H^1,\sP_1}(z) & = & \lambda_1 (\det \cM_2^{(100)}(\pi_z))^{-\frac{1} 4} \cM_2^{(100)}(\pi_z) & \text{ where } & \pi_z := \frac{d^2f_\delta(z)} 2,\\
h_{H^1,\sP_2}(z) &=&  \lambda_2 (\det \cM_3(\pi_z))^{-\frac{1} 6} \cM_3(\pi_z) & \text{ where } & \pi_z := \frac{d^3f_\delta(z)} 6,
\end{array}
\ee
where the positive constants $\lambda_1, \lambda_2$ are adjusted in such way that the meshes generated have $2000$ elements. Mesh generation was performed by the open source program FreeFEM++ \cite{FreeFem} and results are illustrated on Figure \ref{fig3W1P}. Note that we have used
$\cM_2^{(100)}$ (defined as in \iref{defMAlpha}) instead of $\cM_2$
which would lead to a different triangulation $\cT^*_{H^1,\sP_1}$,
also displayed on Figure \ref{fig3W1P}, and
associated to the metric
$$
h^*_{H^1,\sP_1} (z) := \lambda^*_1 (\det \cM_2(\pi_z))^{-\frac{1} 4} \cM_2(\pi_z)  \text{ where }  \pi_z := \frac{d^2f_\delta(z)} 2,
$$
with again $\lambda_1^*$ adjusted to obtain $2000$ elements. 
The use of $\cM_2^{(100)}$ in place of $\cM_2$ is justified by 
mesh generation issues which are discussed in the next subsection.

Similarly, and following the study developed in Chapter \ref{chapOptAniso}, we have generated $\cT_{L^2,\sP_1}$ and $\cT_{L^2,\sP_2}$ from the metrics 
 $$
\begin{array}{cclcc}
h_{L^2,\sP_1}(z) & = & \mu_1 (\det \cN_2(\pi_z))^{-\frac{1} 6} \cN_2(\pi_z) & \text{ where } & \pi_z := \frac{d^2f_\delta(z)} 2,\\
h_{L^2,\sP_2}(z) &=&  \mu_2 (\det \cN_3(\pi_z))^{-\frac{1} 8} \cN_3(\pi_z) & \text{ where } & \pi_z := \frac{d^3f_\delta(z)} 6,
\end{array}
$$
where again $\mu_1,\mu_2$ are positive constants adjusted in order to generate a mesh with $2000$ elements. Here $\cN_2(\pi)  :=\sqrt {\cM_2(\pi)}$ and 
$$
\cN_3(\pi) := \argmin\{\det M\sep M\in S_2^+ \text{ and } \ |\pi(z)| \leq (z^\trans M z)^{\frac 3 2} \text{ for all } z\in \R^2\}.
$$
We have obtained the following results, which confirm that the use of the metric adapted to a given norm and interpolation degree produces the triangulation that yields the smallest interpolation error in this case (at least among these four triangulations). 
\be
\label{comparativeRes}
\begin{array}{r|cccc}
\# \cT = 2000  &  \cT_{H^1,\sP_1} & \cT_{H^1,\sP_2} &  \cT_{L^2,\sP_1} &\cT_{L^2,\sP_2}\\
\hline
\tb{| f_\delta-\interp_\cT^1 f_\delta|_{H^1_0}} & \tr{1.35} & 1.47 &  1.43 &  1.63\\
10 \tb{| f_\delta-\interp_\cT^{\tr 2} f_\delta|_{H^1_0}} & 1.66 &  \tr{1.17} &  1.89 &   1.47\\
10^2 \tb{\|f_\delta-\interp_\cT^1 f_\delta\|_{L^2}} & 1.54 & 2.73 & \tr{0.759} & 1.18 \\
10^4 \tb{\|f_\delta-\interp_\cT^{\tr2} f_\delta\|_{L^2}} & 6.64 & 6.61 & 4.73 & \tr{3.17}
\end{array}
\ee 


\subsection{Quality of a triangulation generated from a metric}
\begin{figure}
\centering
\hspace{-9mm}
\hspace{-9mm}
\includegraphics[width=6.4cm,height=4cm]{\pathPic/PaperW1P/PicP1H1Mesh.pdf}
\hspace{-9mm}
\hspace{-7mm}
\hspace{-9mm}
\includegraphics[width=6.4cm,height=4cm]{\pathPic/PaperW1P/PicP1H1CondMesh.pdf}
\hspace{-9mm}
\hspace{-7mm}
\hspace{-9mm}
\includegraphics[width=6.4cm,height=4cm]{\pathPic/PaperW1P/PicP2H1Mesh.pdf}
\hspace{-9mm}
\hspace{-9mm}
\caption{\label{fig3W1P}The meshes $\cT^*_{H^1,\sP_1}$,  $\cT_{H^1,\sP_1}$ and $\cT_{H^1,\sP_2}$  adapted to  $f_\delta$ ($500$ triangles only).}
\end{figure}
Given a metric $h : \Omega\to S_2^+$, there does not always exists a triangulation $\cT$ adapted to $h$, i.e. satisfying \iref{adaptTh} for some constant $C_1 \geq 1$ not too large. Such a triangulation exists only if $h$ satisfies some constraints which are analyzed in \cite{Shew}, see also Chapter 5.
Instead of analysing the metric $h$ prior to the process of mesh generation, we choose here the simpler option of evaluating a posteriori the quality of a triangulation $\cT$. 

Since we are interested in the $H^1 = W^{1,2}$ semi norm we define following \iref{boundedS} 
$$
S(\cT) := \left(\frac 1 {\# \cT} \sum_{T\in \cT} S(T)^2\right)^{\frac 1 2}.
$$
For all $T\in \cT$ we define $h_T := h(\bary(T))\in S_2^+$. We also define the sets
$$
E := \left\{\ln\left(|T| \sqrt {\det h_T}\right)\sep T\in \cT\right\} \text{ and } F := \left\{\rho\left(\sqrt{h_T} (T)\right)\sep T \in \cT\right\}.
$$
According to \iref{adaptTh}, the quality of $\cT$ is reflected by the quantities
$$
\exp(\max E - \min E) \ \text{ and } \ \max F.
$$
However these quantities give a rather pessimistic account of the adaptation of $\cT$ to $h$, and \emph{heuristically} we find it more fruitful to consider averages. We therefore define
$$
\rho(\cT,h) :=\frac 1 {\# (\cT)} \sum_{T \in \cT} \rho\left(\sqrt{h_T} (T)\right).
$$
and 
$$
\sigma(\cT,h) :=  \exp\left(\frac 1 {\#(E)} \sum_{e\in E} \left|e- \frac 1 {\# (E)} \sum_{e\in E}e\right|\right).
$$ 
The following table shows that the quantities $S(\cT)$, $\rho(\cT,h)$ and $\sigma(\cT,h)$ are abnormally large for the triangulation $ \cT^*_{H^1,\sP_1}$ generated from the metric 
$h^*_{H^1,\sP_1}$
but reasonable for the triangulations $\cT_{H^1,\sP_1}$ and $\cT_{H^1,\sP_2}$ generated from the metrics \iref{defHNum}. 
$$
\begin{array}{c|ccc}
\# \cT = 2000  & \cT^*_{H^1,\sP_1} &  \cT_{H^1,\sP_1} & \cT_{H^1,\sP_2} \\ 
\hline
S(\cT) & 14.2 & 3.14 & 4.04\\
\rho(\cT,h) & 10.6 & 6.02 & 4.18\\
\sigma(\cT,h) & 2.39 &   2.25 & 1.70
\end{array}
$$

In practice $\cT^*_{H^1,\sP_1}$ led to a poor interpolation error, contrary to $\cT_{H^1,\sP_1}$.
 We believe that the poor quality of  $\cT^*_{H^1,\sP_1}$ is due to the excessively wild behavior of the metric $h^*_{H^1,\sP_1}$ and not to a deficiency of the excellent mesh generator BAMG \cite{FreeFem}.

\section{Polynomial equivalents of the shape function}
\label{secW1P5}
The optimal error estimates established in Theorem \ref{mainTheorem}
involve the quantity $L_{m,p}(\frac {d^m f}{m!})$. The function 
$\pi \mapsto L_{m,p}(\pi)$ is obtained by solving an optimization problem,
and it does not have an explicit analytic expression 
in terms of the coefficients of $\pi\in \H_m$. In this section, we introduce quantities which
are equivalent to $L_{m}(\pi)$, and therefore to $L_{m,p}(\pi)$ for
all $1 \leq p\leq \infty$, and  which can be written in analytic form in terms of the coefficients of $\pi\in \H_m$.

Given a pair of non negative functions $Q$ and $R$ on $\H_m$ we write $Q\sim R$ if and only if there exists a constant $C>0$ such that $C^{-1} Q\leq R \leq C Q$ uniformly on $\H_m$. 
We sometimes slightly abuse notations and write $Q(\pi)\sim R(\pi)$.
We say that a function $Q$ is a polynomial on $\H_m$ if there exists a polynomial $P$ of $m+1$ real variables such that for all $a_0, \cdots,  a_m\in \R$,
$$
Q\left(\sum_{i=0}^m a_i x^i y^{m-i}\right) = P(a_0,\cdots, a_m).
$$
We define $\deg Q := \deg P$, and we say that $Q$ is homogeneous
if $P$ is homogeneous. For all $m\geq 2$, we shall build an
homogeneous polynomial $Q$ on $\H_m$ such that 
\be
\label{rootPol}
L_m \sim \sqrt [r]{|Q|}\;\; {\rm with}\;\; r:=\deg Q,
\ee
where the constants in the equivalence only depend on $m$. 

We first introduce for all  $\pi \in \H_m$ the set 
$$
\cB_\pi := \{B\in \M_2(\R)\sep |\pi(z)|\leq |B z|^m \text{ for all } z\in \R^2\}, 
$$
and the function
$$
K_m^\cE(\pi) := \inf\{ |\det B|^{\frac m 2}\sep B\in \cB_\pi\}.
$$
According to Lemma \ref{lemmaA} we have for any $m\geq 2$ 
\be
\label{equivLK}
L_m(\pi) \sim \sqrt{K_{2m-2}^\cE( |\nabla \pi|^2)} 
\ee
where $ |\nabla \pi|^2 = (\partial_x \pi)^2+(\partial_y \pi)^2\in \H_{2m-2}$.
The function $K_m^\cE$ is 
extensively studied in Chapter 2. 
In particular we know that 
\be
\label{K2}
K_2^\cE(\pi) \sim \sqrt{|\det [\pi]|},
\ee
and
\be
\label{K3}
K_3^\cE(\pi)\sim \sqrt[4]{|\disc(\pi)|}
\ee
where $\disc(\pi)$ denotes the discriminant of a polynomial $\pi\in \H_3$, namely
$$
\disc(a x^3 + b x^2  y+ c x y^2 + dy^3) = b^2c^2 - 4ac^3 - 4b^3d + 18abcd - 27a^2 d^2.
$$
More generally, it is proved in Chapter \ref{chapOptAniso} that for 
all $m\geq 2$, the function $K_m^\cE$ has an equivalent 
of the form  $\sqrt [r]{|Q|}$, where $Q$ is an homogeneous polynomial
of degree $r$ on $\H_m$. Combining this result with 
\iref{equivLK} we obtain the main result of this section.

\begin{prop}
\label{thPolEq}
Let $m\geq 2$ and let $Q$ be an homogeneous polynomial on $\H_{2m-2}$ such that 
$K_{2m-2}^\cE \sim \sqrt[r]{|Q|}$, where $r = \deg Q$. 
Let $Q_*$ be the polynomial on $\H_m$ defined by 
$$
Q_*(\pi) := Q(|\nabla \pi|^2).
$$
Then $L_m \sim \sqrt[2r]{Q_*}$ on $\H_m$.
\end{prop}

\noindent
Let $\pi\in \H_2$ and let us observe that $|\nabla \pi(z)|^2 = |2[\pi]z|^2 = 4z^\trans [\pi]^2 z$. 
Using \iref{K2} we therefore obtain
\be
\label{L2eq}
L_2(\pi) \sim \sqrt{K_2^\cE(|\nabla\pi|^2)} \sim \sqrt{\sqrt{\det(4 [\pi]^2)}} =2 \sqrt{|\det [\pi]|}. 
\ee

The construction suggested by Theorem \ref{thPolEq} uses an equivalent of $K_{2m-2}^\cE$ to produce an equivalent to $L_m$.
Unfortunately, as $m$ increases, the practical construction of $Q$ such that $\sqrt [r]{|Q|}$
is equivalent to $K_{m}^\cE$ becomes more involved and the degree
$r$ quickly raises. In the following theorem, we build 
an equivalent to $L_{m}$ from an equivalent of $K_{m-1}^\cE$ instead of $K_{2m-2}^\cE$, which is therefore
simpler.

\begin{theorem}
\label{thPolEq2}
Let $m\geq 3$ and let $Q$ be an homogeneous polynomial on $\H_{m-1}$ such that $K_{m-1}^\cE \sim \sqrt[r]{|Q|}$, 
where $r=\deg Q$. Let $(Q_k)_{0\leq k\leq r}$ be the 
homogeneous polynomials of degree $r$ on $\H_{m-1}\times \H_{m-1}$ such that 
for all $u,v\in \R$ and all $\pi_1,\pi_2\in \H_m$ we have
\be
\label{defQk}
Q(u \pi_1+v \pi_2) = \sum_{0\leq k\leq r} \binom r k u^k v^{r-k} Q_k(\pi_1,\pi_2),
\ee
where $\binom r k := \frac {r!}{k!(r-k)!}$.
Let $Q_*$ be the polynomial defined for all $\pi\in\H_{m}$ by 
$$
Q_*(\pi) := \sum_{0\leq k\leq r} \binom r k Q_k\left(\partial_x \pi,\partial_y \pi\right)^2
$$ 
then $L_{m} \sim \sqrt[2r]{Q_*}$ on $\H_{m}$.
\end{theorem}

\proof
See Appendix.
\sq

Using this construction and \iref{K2} we obtain an equivalent of $L_3$ as follows.
Let $\pi_1 = a x^2+ 2b x y + c y^2$ and $\pi_2 = a' x^2+ 2b' x y + c' y^2$ be two elements
of $\H_2$. We obtain
\begin{eqnarray*}
\det( [u \pi_1 + v \pi_2]) &=& (ua+va') ( uc+vc') - (ub+vb')^2\\
& =& u^2 (ac-b^2) + u v (ac'+a'c - 2 b b') + v^2 (a'c' - b'^2).
\end{eqnarray*}
Applying the construction of Theorem \ref{thPolEq2} to $\pi = a x^3 + 3 b x^2  y+ 3 c x y^2 + dy^3\in \H_3$ we obtain
\be
\label{L3eq}
L_3(\pi) \sim 3 \sqrt[4]{(ac-b^2)^2+(ad-bc)^2/2+(bd-c^2)^2}.
\ee
Remarking that 
$$
2 [(ac-b^2)^2+(ad-bc)^2/2+(bd-c^2)^2 ]= (a^2+2b^2+c^2)(b^2+2c^2+d^2) - (ab+2bc+cd)^2,
$$
 and using equation \iref{detMPi} we obtain that $L_3(\pi) \sim \sqrt{\det \cM_3(\pi)}$. 
This point was announced in \S \ref{secW1P4p1} and establishes that the map $\cM_3$ defined in \iref{defMP2} can be used for optimal mesh adaptation for quadratic finite elements.
\newline
\newline
\noindent
Using \iref{K3} and the construction of Theorem \ref{thPolEq2}, we 
also obtain an equivalent of $L_4(\pi)$ 
\begin{eqnarray*}
L_4(\pi)^8 & \sim &(3 b^2 c^2 - 4 a c^3 - 4 b^3 d + 6 a b c d - a^2 d^2)^2 \\
&+& (2 b c^3 - 6 a c^2 d + 4 a b d^2 - 4 b^3 e + 6 a b c e - 2 a^2 d e)^2/4 \\
&+& (3 c^4 - 6 b c^2 d + 8 b^2 d^2 - 6 a c d^2 - 6 b^2 c e + 6 a c^2 e + 2 a b d e - a^2 e^2)^2/6\\
&+& (2 c^3 d - 4 a d^3 - 6 b c^2 e + 4 b^2 d e + 6 a c d e - 2 a b e^2)^2/4 \\
&+& (3 c^2 d^2 - 4 b d^3 - 4 c^3 e + 6 b c d e - b^2 e^2)^2.
\end{eqnarray*}

The following proposition identifies the polynomials $\pi\in \H_m$ for which $L_m(\pi) = 0$, and therefore the 
values of $d^{m}f$ for which anisotropic mesh adaptation may lead to \emph{super-convergence}.

\begin{prop}
Let $m\geq 2$ and let $t_m := \left\lfloor \frac{m+3} 2 \right\rfloor$. Then for all $\pi \in \H_m$,
\be
\label{vanishL}
L_m(\pi) = 0 \text{ if and only if } \pi = (\alpha x+ \beta y)^{t_m} \tilde \pi \text{ for some } \alpha, \beta\in \R \text{ and } \tilde \pi \in \H_{m-t_m}.
\ee
\end{prop}
\proof
According to \iref{equivLK},  $L_m(\pi) = 0$ if and only if $K_{2m-2}^\cE(|\nabla \pi|^2)=0$. 
On the other hand, it is proved in Chapter \ref{chapOptAniso} that for any $\pi_*\in \H_{2m-2}$ one has $K_{2m-2}^\cE(\pi_*)=0$ if and only if $\pi_*$ has a linear factor of multiplicity $m$. Therefore $L_m(\pi) = 0$ if and only
$|\nabla \pi|^2$ is a multiple of $l^m$, where $l$ is of the form $l=\alpha x +\beta y$.

Let us first assume that $|\nabla \pi|^2 = (\partial_x \pi)^2+ (\partial_y \pi)^2$ has such a form. Since they are non-negative, $(\partial_x \pi)^2$ and $(\partial_y \pi)^2$ are both multiples of $l^m$. Therefore $\partial_x \pi$ and $\partial_y \pi$ are multiples of $l^s$ where $s$ is an integer such that $2s\geq m$, hence $s\geq t_m-1$. We therefore have 
$$
\partial_x \pi = l^s \pi_1 \ \text{ and } \ \partial_y \pi = l^s \pi_2 \text{ where } \pi_1,\pi_2 \in \H_{m-1-s}
$$
Recalling that $l = \alpha x+\beta y$ we obtain
$$
0 = \partial_{yx}^2 \pi -\partial_{xy}^2 \pi = l^s(\partial_y \pi_1 - \partial_x\pi_2) + sl^{s-1} (\beta\pi_1-\alpha\pi_2),
$$
hence $\beta\pi_1-\alpha\pi_2$ is a multiple of $l$. Since $\pi$ is homogenous of degree $m$ it obeys
the Euler identity $m \pi(z) = \<z, \nabla \pi(z)\>$ for all $z\in \R^2$. 
Assuming without loss of generality that $\alpha \neq 0$,
we therefore obtain
$$
m \pi(x,y) =  l^s (x\pi_1 +y \pi_2) = l^s\left( (\alpha x+\beta y)\frac {\pi_1} \alpha + \frac y \alpha (\alpha \pi_2 - \beta \pi_1)\right)
$$
which shows that $\pi$ is a multiple of $l^{s+1}$, hence of $l^{t_m}$.  

Conversely if $\pi$ is a multiple of $l^{t_m}$ then $\partial_x \pi$ and $\partial_y \pi$ are both multiples of $l^{t_m-1}$. Since $2(t_m-1)\geq m$ the polynomial $|\nabla \pi|^2$ is a multiple of $l^m$ which concludes the proof.
\sq

\section{Extension to higher dimension}
\label{secW1P6}
This section partially extends the results exposed in the previous sections to functions of $d$ variables. 
We give in \S \ref{secW1P6p1} the generalisations of the shape function $L_{m,p}$ and of the measure of sliverness $S$. 

Subsection \S \ref{secW1P6p2} is devoted to interpolation error estimates. We prove a local $d$-dimensional error estimate in Theorem \ref{thLocalD} which generalises Theorem \ref{thLocal}. We then establish an asymptotic lower error estimate in Theorem \ref{optiTheoremD} which generalises Theorem \ref{optiTheorem}. We give sufficient conditions under which the interpolation on a $d$-dimensional mesh $\cT$ achieves this optimal lower bound up to a multiplicative constant. 
However due to technical issues linked to the measure of sliverness $S$ we were not able to construct such meshes, and we therefore state the upper bound as a conjecture.

We discuss in subsection \S \ref{secW1P6p3} the construction of optimal metrics for practical mesh generation. We partially extend the results of \S \ref{secW1P4} and raise open questions.

\subsection{Generalisation of the shape function and of the measure of sliverness.}
\label{secW1P6p1}
We extend in this section the tools used in our analysis of optimally adapted triangulations to arbitrary dimension $d$. We begin with the spaces of polynomials. Let 
$$
\H_{m,d} := {\rm Span}\{x_1^{\alpha_1}\cdots x_d ^{\alpha_d}\sep |\alpha | = m\} \text{ and } \P_{m,d} := {\rm Span}\{x_1^{\alpha_1}\cdots x_d ^{\alpha_d}\sep |\alpha | \leq m\},
$$
where $\alpha = (\alpha_1, \cdots, \alpha_d)$ denotes a d-plet of non-negative integers, $|\alpha|:= \alpha_1+\cdots+ \alpha_d$.
For any simplex $T$ the Lagrange interpolation operator $\interp_T^m : C^0(T) \to \P_{m,d}$ is defined 
by imposing $f(\gamma) = \interp_T^m f(\gamma)$ for all points $\gamma$ with barycentric coordinates in the set $\{0,\frac 1 m, \frac 2 m,\cdots , 1\}$ with respect to the vertices of $T$.
For all $\pi \in \H_{m,d}$ we define 
$$
L_{m,d,p}(\pi) := \inf_{|T| = 1} \|\nabla \pi - \nabla \interp_T^{m-1} \pi\|_{L^p(T)},
$$
where the infimum is taken on the set of $d$-dimensional simplices of unit volume. Similarly to \iref{equivLm} the functions $L_{m,d,p}$, $1\leq p \leq \infty$, are uniformly equivalent on $\H_{m,d}$. We define $L_{m,d} := L_{m,d,\infty}$.

The distance defined at \iref{defDist} between triangles extends easily to simplices. Given two $d$-dimensional simplices $T$, $T'$ there are precisely $(d+1)!$ affine transformations $\Psi$ such that $\Psi(T) = T'$. For
each such $\Psi$, we denote by $\psi$ its linear part and we define
$$
d(T,T') := \ln \(\inf\{\kappa(\psi) \sep \Psi(T) = T'\}\).
$$
We say that a $d$-dimensional simplex $T$ is acute if the exterior normals $n,n'$ to any two distinct faces $F,F'$ of $T$ have a negative scalar product $\<n, n'\>$. In other words if all pairs of faces of $T$ form acute dihedral angles. We denote the set of acute simplices by $\mathbb A$ and we generalise the measure of sliverness to arbitrary dimension $d$ as follows
\be
\label{defSD}
S(T) := \exp d (T,\mathbb A) = \inf\{\kappa(\psi) \sep \Psi(T) \in \mathbb A\}.
\ee
Similarly to \iref{defS}, the quantity $S(T)$ reflects the distance from a simplex $T$ to the set of acute simplexes $\mathbb A$. 
The definition \iref{defSD} of $S(T)$ raises a legitimate question : how to produce an affine transformation $\Psi$ such that $\Psi(T)$ has acute angles, and $\kappa(\psi)$ is comparable to $S(T)$?
This question is answered by the following proposition. 

For any $d$-dimensional simplex $T$ with vertices $(v_i)_{0\leq i\leq d}$, we define the symmetric matrix
\be
\label{defM}
M_T := \sum_{0\leq i<j \leq d} e_{ij} e_{ij}^\trans, \text{ where } e_{ij} := \frac{v_i - v_j}{|v_i - v_j|}.
\ee
Observe that 
\be
\label{alphaD}
1\leq \|\sqrt{M_T}\| = \sqrt{\|M_T\|} \leq \alpha_d\stext{ where} \alpha_d := \sqrt{\frac {d(d+1)} 2},
\ee
since $\alpha_d^2$ is the number of distinct pairs $(i,j)$ satisfying $0\leq i<j \leq d$.

\begin{prop}
\label{propMS}
For any simplex $T$, the simplex $M_T^{-\frac 1 2}(T)$ is acute and  
\be
\label{eqMS}
S(T)\leq  \kappa(\sqrt{M_T}) \leq \alpha_d S(T).
\ee 
\end{prop}

\proof
See appendix.
\sq

\begin{remark}
In the paper \cite{Ja} an alternative measure of sliverness $S'(T)$ of a simplex $T$ is introduced, and defined as
$$
S'(T) := \left(\inf_{|u| = 1} \max_{ i<j } |\<u, e_{i,j}\>|\right)^{-1}.
$$
This quantity is equivalent to $S(T)$. Indeed, for any $u\in \R^d$ we have
$$
\max_{ i<j } |\<u, e_{i,j}\>| \leq \sqrt{\sum_{i<j} \<e_{ij}, u\>^2} =  \sqrt{u^\trans M_T u}  = \left|M_T^{\frac 1 2} u\right| \leq \alpha_d \max_{ i<j } |\<u, e_{i,j}\>|,
$$
which implies that $\alpha_d^{-1} S'(T)\leq \|M_T^{-\frac 1 2}\| \leq S'(T)$, hence $\alpha_d^{-1} S(T)\leq S'(T) \leq \alpha_d^2 S(T)$ using \iref{alphaD}. 
Our approach therefore introduces a new geometrical interpretation to the quantity $S'$ introduced in \cite{Ja}, as the distance from a given simplex to the set of acute simplices.
\end{remark}
The following lemma generalises Lemma \ref{lemmaS} and shows that the interpolation process is stable in the $L^\infty$ norm of the gradient if the measure of sliverness is controlled. Let us mention that a slightly different version of this lemma can be found in \cite{Ja}, yet not exactly adapted to our purposes. 
\begin{lemma}
\label{lemmaSD}
For all $m\geq 2$ and all $d\geq 2$ there exists a constant $C=C(m,d)$ such that 
for any $d$-dimensional simplex $T$ and any $f\in W^{1, \infty}(T)$, one has
\be
\label{ineqSD}
\|\nabla \interp^m_T f \|_{L^\infty(T)} \leq C S(T) \|\nabla f \|_{L^\infty(T)}.
\ee
\end{lemma}

\proof
The proof this lemma is extremely similar to the proof of Lemma \iref{lemmaS}.
Let $\TRect$ be the simplex which vertices are the origin and the canonical basis of $\R^d$. For the same reason as in Lemma \ref{lemmaS}, if a function $\ti g(x_1,x_2,\cdots ,x_d)\in C^0(\TRect)$ does not depend on the coordinate $x_d$, then $\interp_\TRect^{m-1} \ti g$ does not depend on $x_d$ either.
Using the same reasonning as in Lemma \ref{lemmaS} we obtain that there exists a constant $C_0 = C_0(m,d)$ such that for all $g\in W^{1, \infty}(\TRect)$
$$
\left\|\frac{\partial\interp_\TRect^m g}{\partial x_d}\right\|_{L^\infty(\TRect)} \leq C_0 \left\|\frac{\partial g}{\partial x_d}\right\|_{L^\infty(\TRect)}.
$$
Again similarly to the proof of Lemma \ref{lemmaS} we obtain using a change of variables that for any simplex $T$, any $f\in W^{1,\infty}(T)$ and any edge vector $u$ of $T$
$$
\|\<u,\nabla \interp_T^m f\>\|_{L^\infty(T)}\leq C_0 \|\<u, \nabla f\>\|_{L^\infty(T)}.
$$
We use the notations of Proposition \ref{propMS} and we define a norm $|v|_T$ on $\R^d$ by 
$$
|v|_T^2 := v^\trans M_T v = \sum_{0\leq i<j\leq d} \<v,e_{ij}\>^2.
$$ 
Observe that 
\be
\label{euclMT}
\|M_T^{-\frac 1 2}\|^{-1} |v|\leq |v|_T \leq \|M_T^{\frac 1 2}\| |v|.
\ee
Then, since $e_{ij}$ is proportional to an edge vector of $T$, 
\begin{eqnarray*}
\| \ |\nabla \interp_T^m f|_T \|_{L^\infty(T)}^2 &\leq& \sum_{0\leq i<j\leq d} \|\<e_{ij}, \nabla \interp_T^m f\>\|_{L^\infty(T)}^2\\
&\leq&  C_0 \sum_{0\leq i<j\leq d} \|\<e_{ij}, \nabla f\>\|_{L^\infty(T)}^2 \leq C_0 \alpha_d^2 \| \ |\nabla f|_T \|_{L^\infty(T)}^2.
\end{eqnarray*}
Combining this result with \iref{euclMT} we obtain 
$$
\|M_T^{-\frac 1 2}\|^{-1}  \| \nabla \interp_T^m f \|_{L^\infty(T)} \leq C_0 \alpha_d^2 \|M_T^{\frac 1 2}\| \|  \nabla f \|_{L^\infty(T)}
$$
and we conclude the proof using \iref{eqMS}.
\sq

The oscillation of the gradient of the interpolated function is an important problem encountered by numerical methods that try to take advantage of highly anisotropic meshes, see the discussion in \cite{Shew2}. As the previous lemma shows, such oscillations are kept under control if $S(T)$ is bounded on the mesh of interest. For checking this property in pratical situations one needs an equivalent of the sliverness $S$ that can be computed at low numerical cost. The formula \iref{defSD} is clearly not adapted, since it involves a complicated optimisation procedure. Instead we propose to use
\be
\label{equivS}
\hat S(T) := \sqrt{\Tr(M_T^{-1})}.
\ee
Observing that $\|M_T^{-\frac 1 2} \| \leq \hat S(T) \leq \sqrt d \|M_T^{-\frac 1 2} \|$, and recalling that $1\leq \|M_T\|\leq \alpha_d$ we obtain 
$$
\alpha_d^{-1} \hat S(T) \leq S(T) \leq \sqrt d \alpha_d \hat S(T).
$$
Note that $\hat S(T)$ has an analytic expression in terms of the coordinates of $T$ : the square root of the ratio of two polynomials in the positions of the vertices of $T$.\\

\begin{remark} We illustrate the sharpness of inequality \iref{ineqSD} in a simple example.
Let $x,y,z$ be the coordinates on $\R^3$ and let $\pi_0 := x^2 \in \H_{2,3}$.
Let $T_\lambda$ be the tetrahedron of vertices $(-\lambda,0,0)$, $(\lambda,0,0)$, $(\lambda,1,0)$ and $(0,0,1)$. 
Simple computations show that
$$
\|\nabla \interp_{T_\lambda}^1 \pi_0\|_{L^\infty(T_\lambda)} = \lambda^2, \ 
\|\nabla \pi_0\|_{L^\infty(T_\lambda)} = 2 \lambda  
\text{ and } \lim_{\lambda\to \infty}  \frac{\hat S(T_\lambda)} \lambda = \sqrt{\frac 5 7}.
$$
Let $T'_\lambda$ be defined by replacing the vertex $(-\lambda,0,0)$ of $T_\lambda$ with $(0,0,0)$. Then 
$$
\|\nabla \interp_{T_\lambda}^1 \pi_0\|_{L^\infty(T'_\lambda)} = \lambda, \ 
\|\nabla \pi_0\|_{L^\infty(T'_\lambda)} = 2 \lambda  
\text{ and } \lim_{\lambda\to \infty}  \hat S(T'_\lambda)  = \frac 3 2.
$$
Hence the simplices $T_\lambda$ and $T'_\lambda$ have very different interpolation properties for large $\lambda$, although they have a similar aspect ratio. They are representatives of ``bad'' and ``good'' anisotropy respectively. The tetrahedrons $T_{\frac 3 2}$ and $T'_{\frac 3 2}$ are illustrated on the left of Figure 4, bottom and top respectively.
 \end{remark}
\begin{figure}
	\centering
		\vspace{-3cm}
		\includegraphics[width=15cm,height=9cm]{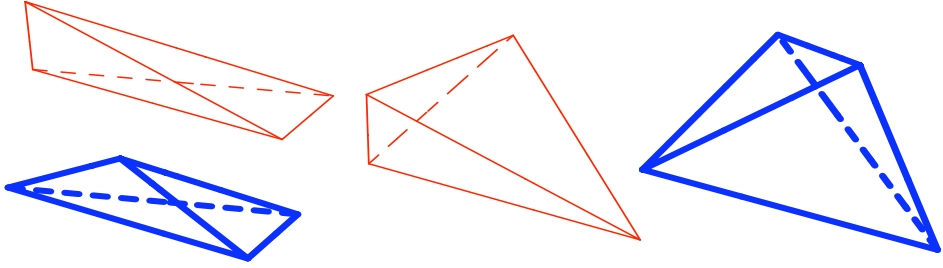}
		\vspace{-3cm}
	\caption{\label{fig4W1P}Examples of Good anisotropy (Thin lines, $S(T) \sim 1$), and Bad anisotropy (Thick lines, $S(T)\gg1$) .}
\end{figure}

For any $d$-dimensional simplex $T$, we define its measure of non degeneracy by  
$$
\rho(T) := \frac{\diam(T) ^d}{|T|}.
$$
Let $T_*$ be a fixed $d$-dimensional acute simplex, for instance the reference equilateral simplex. For any $d$-dimensional simplex $T$ let $\psi\in \GL_d$ and $z\in \R^d$ be such that $T = z+\psi(T_*)$. 
Since $T_*$ is acute, we obtain a generalization of \iref{ineqRhoS}
$$
S(T) \leq \kappa(\psi)\leq \|\psi\|^d |\det \psi|^{-1} \leq \frac{\diam(T)^d}{\mu(T_*)^d}\frac{ |T_*|}{ |T|} = C(d) \rho(T).
$$
where $\mu(T_*)$ is the diameter of the largest ball incribed in $T_*$,
and where we have used the inequality $|\det(\psi^{-1})| \geq \|\psi^{-1}\| \|\psi\|^{-(d-1)}$. This last inequality can be
derived by using the singular value decomposition $\psi=U{\rm diag}(\lambda_1,\cdots,\lambda_d)V$
with $0<\lambda_1\leq \cdots \leq \lambda_d$
and noting that $\|\psi\|=\lambda_d$ and $\|\psi^{-1}\|=\lambda_1^{-1}$.
\nl
\nl
Generalizing \iref{defA}, we define for all $\pi\in \H_{m,d}$,
$$
\cA_\pi := \{A\in \M_d(\R)\sep |\nabla \pi(z)| \leq |Az|^{m-1} \text{ for all } z\in \R^d\}. 
$$
Geometrically, one has $A \in \cA_\pi$ if and only if  the ellipsoid $\{z\in \R^d \sep |A z|\leq 1 \}$ is included in the algebraic set $\{z\in \R^d\sep |\nabla \pi(z)|\leq 1\}$.
This leads us to the generalisation of Lemma \ref{lemmaA}.
\begin{lemma}
\label{lemmaAD} 
For all $m\geq 2$ and all $d\geq 2$ 
there exists a constant $C=C(m,d)$ such that 
for all $\pi \in \H_{m,d}$, we have
$$
C^{-1} L_{m,d}(\pi)\leq \inf \{|\det A|^{\frac{m-1} d} \sep A \in \cA_\pi\} \leq C L_{m,d}(\pi).
$$
\end{lemma}

\proof
The proof of this lemma is completely similar to the proof of its bidimensional version Lemma \ref{lemmaA}.
The only point that needs to be properly generalized is the following : given a matrix $A\in \GL_d$, how to construct an acute simplex $T = T(A)$ such that $\rho(A(T))$ is bounded independently of $A$?

The following construction is not the simplest but will be useful in our subsequent analysis. Let $A=U D V$, be the singular value decomposition of $A$, where $U,V$ are orthogonal matrices and $D$ is a diagonal matrix with positive diagonal entries $(\lambda_i)_{1\leq i\leq d}$. We define the Kuhn simplex $T_0$
$$
\TRect := \{x\in [0,1]^d\sep x_1\geq x_2\geq \cdots \geq x_d\},
$$
and $T :=  V^\trans D^{-1} \TRect$.
Then $\rho(A(T)) = \rho(U(\TRect)) = \rho(\TRect) = d! d^{d/2}$ which is independent of $A$. 
We now show that $T$ is an acute simplex.

Let  $(e_1,\cdots, e_d)$ be the canonical basis of $\R^d$, and let by convention $e_0 = e_{d+1} = 0$.
For $0\leq i\leq d$, an easy computation shows that the the exterior normal to the face $F_i$ of $\TRect$, opposite to the vertex $v_i =  \sum_{0\leq k \leq i} e_k$, is 
$
n_i = \frac{ e_i - e_{i+1}}{\| e_i - e_{i+1}\|}.
$
It follows that the exterior normal $n'_i$ to the face $D^{-1}(F_i)$ of the simplex $D^{-1}(T_0)$ is 
 $$
 n'_i = \frac{D(n_i)}{|D (n_i)|} = \frac{\lambda_i e_i - \lambda_{i+1} e_{i+1}}{|\lambda_i e_i - \lambda_{i+1} e_{i+1}|}.
 $$
Hence $\<n'_i, n'_j\> = 0$ if $|i -j|>1$, and $\<n'_i, n'_{i+1}\> < 0$ for all $0\leq i\leq d-1$. It follows that the simplex $D^{-1}(T_0)$ is acute, and therefore $T = V^\trans D^{-1} (T_0)$ is also acute since $V$ is a rotation.
 \sq

\subsection{Generalisation of the error estimates}
\label{secW1P6p2}
We present in this section the generalisation to higher dimension of our anisotropic error estimates. We prove a local error estimate in theorem \ref{thLocalD} and an asymptotic lower estimate in \ref{optiTheoremD}. We also point out in conjecture \ref{conjTileD} a technical point which, if proved, would lead to the optimal asymptotic upper estimates \iref{upperEstimD} and \iref{upperEstimDEps}.


\begin{theorem}
\label{thLocalD}
For all $m\geq 2$ and $d\geq 2$ there exists a constant $C=C(m,d)$ such that for all $\pi\in \H_{m,d}$, all $A \in \cA_\pi$ and any simplex $T$ we have
\be
\label{localPiD}
|\pi - \interp^{m-1}_T \pi|_{W^{1,p}(T)} \leq C |T|^{\frac 1 \tau} \, S(T) \, \rho(A (T))^{\frac{m-1} d} |\det A|^{\frac{m-1} d}, 
\ee
where $\frac 1 \tau := \frac {m-1} d +\frac 1 p$. 
Furthermore for any $g\in C^m(T)$ we have 
$$
|g - \interp^{m-1}_T g|_{W^{1,p}(T)} \leq C |T|^{\frac 1 \tau} \, S(T) \, \rho(T)^{\frac{m-1} d} \|d^m g\|_{L^\infty(T)}.
$$
\end{theorem}

\proof
It is a straightforward generalization 
of the proof of Theorem \ref{thLocal}.
\sq

Combining these two estimates, we can obtain a mixed estimate similar to \iref{localMixed}, with the new value of $\tau$ and the generalised $S$ and $\rho$. 
For all $m\geq 2$ and $d\geq 2$ there exists a constant $C=C(m,d)$ such that for any simplex $T$, any $f\in C^m(T)$, any $\pi\in \H_m$ and any $A\in \cA_\pi$
\be
\label{localMixedD}
\begin{array}{l}
\displaystyle |f - \interp^{m-1}_T f|_{W^{1,p}(T)} \\
\displaystyle\quad \leq  C |T|^{\frac 1 \tau} S(T) \left( \rho(A(T))^{\frac{m-1} d} |\det A|^{\frac {m-1} d} +  \rho(T)^{\frac{m-1} d} \|d^m f-d^m\pi\|_{L^\infty(T)} \right).
\end{array}
\ee
This leads us to a straightforward generalisation of the points (i) to (iv) exposed in \iref{eqError}. Similarly to the bidimensional case \iref{nonAsympt} if a triangulation $\cT$ meets these requirements, then it satisfies the error estimate 
\be
\label{nonAsymptD}
\# (\cT)^{\frac{m-1} d}|f - \interp^{m-1}_\cT f|_{W^{1,p}(\Omega)} \leq C\left\|L_m(\pi_z)+\ve\right\|_{L^\tau(\Omega)}.
\ee


Generalizing \iref{admissibilityCond}, we say that a sequence $\seqT$ of simplicial meshes of a $d$-dimensional polygonal domain is admissible if $\#(\cT_N)\leq N$ and if there exists a constant $C_A>0$ such that 
$$
\sup_{T\in \cT_N} \diam(T) \leq C_A N^{-\frac 1 d}.
$$
Similarly to \iref{lowerEstim}, it can be shown
that \iref{nonAsymptD} cannot be improved for an admissible sequence of triangulations,
in the following asymptotical sense.

\begin{theorem}
\label{optiTheoremD}
Let $\seqT$ be an admissible sequence of triangulations of a domain $\Omega$, let $f\in C^m(\Omega)$ and $1\leq p<\infty$. Then
\be
\label{lowerEstimD}
\liminf_{N\to \infty} N^{\frac{m-1} d}|f-\interp^{m-1}_{\cT_N}  f |_{W^{1,p}(\Omega)} \geq \left\|L_{m,d,p}\left(\frac {d^m f}{m!}\right)\right\|_{L^\tau(\Omega)}
\ee
where $\frac 1 \tau := \frac {m-1} d+\frac 1 p$.
\end{theorem}

\proof
It is identical to the proof of the bidimensional estimate \iref{lowerEstim}, which is exposed in \S \ref{secW1P3p1}.
\sq

In contrast, the upper estimates \iref{upperEstimW1P} and \iref{upperEstimEpsW1P} do not generalize easily to
higher dimension.
A first problem is that the bidimensional mesh $\cP_T$ defined in Equation \iref{defPT} has no equivalent in higher dimension,
in the sense that we cannot exactly tile the space by simplices of optimal shape. 
We may however build a tiling made of near optimal simplices, based on the following procedure: 
for any permutation $\sigma\in \Sigma_d$ of $\{1,\cdots,d\}$ we define
$$
T_\sigma := \{x\in [0,1]^d\sep x_{\sigma(1)} \geq \cdots \geq x_{\sigma(d)}\}.
$$
Let $A\in \GL_d(\R)$, and let $A=UDV$ be the singular value decomposition of $A$, where $U$ and $V$ are unitary and $D$ is diagonal. We define 
$$
\cP_A := \{V^\trans D^{-1}(T_\sigma+z)\sep \sigma \in \Sigma_d, \ z\in \ZZ^d\},
$$
which is a tiling of $\R^d$ built of \emph{acute} simplices $T$ satisfying $\rho(A(T)) = d! d^{d/2}$ (these properties are established in the proof of Lemma \ref{lemmaSD}). Using such a tiling, we would like to build partitions $\cP_{A,n}(R)$
of any $d$-dimensional simplex $R$, with properties similar to those expressed in Lemma \ref{lemmaTile}
for the triangulations $\cP_{T,n}(R)$. At the present stage we do not know how to properly
adapt the construction of $\cP_{A,n}(R)$ near the boundary of $R$ in order to respect
the condition on the measure of sliverness.
The following conjecture, if established, would serve as a generalisation of Lemma \ref{lemmaTile}.
\begin{conjecture}
\label{conjTileD}
Let $R$ be a d-dimensional simplex, and let $A\in \GL_d(\R)$. There exists a sequence $(\cP_{A,n}(R))_{N\geq 0}$, of conformal triangulations of $R$ such that
\begin{itemize}

\item Nearly all the elements of $\cR_N$ belong to $\cP_{A,n} := \frac 1 n \cP_A$, in the sense that 
$$
\lim_{n\to \infty} \frac {\# (\cP_{A,n}^1(R))} {n^d} = \frac {d! |R|}{|\det A|} \ \text{ and } \  \lim_{n\to \infty}  \frac {\# (\cP_{A,n}^2(R))} {n^d}   = 0.
$$
where
$$
\cP_{A,n}^1(R) := \cP_{A,n}(R)\cap \cP_{A,n} \ \text{ and } \quad \cP_{A,n}^2(R) := \cP_{A,n}(R)\sm \cP_{A,n}
$$

\item The restriction of $\cP_{A,n}(R)$ to a face $F$ of $R$ is its standard periodic tiling with $n^{d-1}$ elements.
\item 
The sequence $( \cP_{A,n}(R))_{n\geq 0}$ satisfies 
$$
\sup_{n\geq 0} \left(n \max_{T\in \cP_{A,n}(R)} \diam(T)\right) <\infty \ \text{ and } \ \sup_{n\geq 0} \, \max_{T\in \cP_{A,n}(R)} S(T)< \infty. 
$$
\end{itemize}
\end{conjecture}
The validity of this conjecture would imply the following result using the same proof as
for the estimates \iref{upperEstimW1P} and \iref{upperEstimEpsW1P} established in \S \ref{secW1P3p1}.
\begin{conjecture}
\label{upperPropD}
For all $m\geq 2$ there exists a constant $C=C(m,d)$ such that the following holds.
Let $\Omega\subset \R^d$ be polygonal domain, let $f\in C^m(\Omega)$ and $1\leq p<\infty$. Then there exists a sequence $\seqT$ of simplicial meshes of $\Omega$ such that $\# (\cT_N)\leq N$ and  
\be
\label{upperEstimD}
\limsup_{N\to \infty} N^{\frac{m-1} d}|f-\interp^{m-1}_{\cT_N}f |_{W^{1,p}(\Omega)} \leq C \left\|L_{m,d}\left(\frac {d^m f}{m!}\right)\right\|_{L^\tau(\Omega)}
\ee
where $\frac 1 \tau := \frac {m-1} d +\frac 1 p$. 
Furthermore, for all $\ve>0$, there exists an admissible sequence of simplicial meshes $(\cT^\ve_N)_{N\geq N_0}$ of $\Omega$ such that $\# (\cT^\ve_N)\leq N$ and
\be
\label{upperEstimDEps}
\limsup_{N\to \infty} N^{\frac{m-1} d}|f-\interp^{m-1}_{\cT^\ve_N} f |_{W^{1,p}(\Omega)} \leq  C\left\|L_{m,d}\left(\frac {d^m f}{m!}\right)\right\|_{L^\tau(\Omega)}+\ve.
\ee
\end{conjecture}

\subsection{Optimal metrics and algebraic expressions of the shape function.}
\label{secW1P6p3}
The theory of anisotropic mesh generation in dimension three or higher is only at its infancy. 
However efficient software already exists such as \cite{Inria} for tetrahedral 
mesh generation in domains of $\R^3$. A description of such an algorithm can be found in \cite{BFGLS} as well as some applications to computational mechanics. These software take as input a field $h(z)$
of symmetric positive definite matrices and attempt to create a mesh satisfying \iref{adaptTh}. 
Defining the set of symmetric matrices $\cA'_\pi$ in a similar way as in the
two dimensional case \iref{defA'}, let us
 consider a continuous function $\cM_{m,d}^{(\ve)} : \H_{m,d}\to S_d^+$ such that 
\be
\label{defMD}
\cM_{m,d}^{(\ve)}(\pi) \in \cA'_\pi \text{ and } \det \cM_{m,d}^{(\ve)}(\pi) \leq K(\inf \{\det M\sep M\in \cA'_\pi\}+\ve),
\ee
where $\ve>0$ can be chosen arbitrarily small and where $K$ is an absolute constant independent of $\ve$. The existence of such a function is established in Chapter 6.
Let $\cM(z) := \cM_{m,d}^{(\ve)}(d^m f(z))$,  
let $\delta>0$ and let 
\be
\label{defHd}
h(z) := \delta^{-\tau}(\det \cM(z))^{\frac{-\tau}{dp}} \cM(z).
\ee
A heuristic analysis similar to the one developed in \S \ref{secW1P4} suggests that mesh generation based on this metric leads to a mesh $\cT$ of $\Omega$ optimally adapted for approximating $f$ with $\P_{m-1}$ elements in the $W^{1,p}$ semi norm. This justifies the search for functions $\cM_{m,d}$ satisfying \iref{defMD}.

The form of $\cM_{2,d}$, which corresponds to piecewise linear finite elements,
is already established, see for instance \cite{Shew2}, but we recall it for completeness. 
The same analysis as in \S \ref{secW1P4p2} shows that 
$$
\cM_{2,d}(\pi) := 4 [\pi]^2
$$
satisfies $\cM_{2,d}(\pi)\in \cA'_\pi$ and $\det \cM_{2,d}(\pi) = \inf \{\det M \sep M\in \cA'_\pi\}$. As a byproduct we obtain from Lemma \ref{lemmaAD} that there exists a constant $C=C(d)$ such that for all $\pi\in \H_{2,d}$
$$
C^{-1} \sqrt[d]{|\det [\pi]|}\leq L_{2,d}(\pi) \leq C \sqrt[d]{|\det [\pi]|}.
$$
For piecewise quadratic elements, we generalise \iref{defMP2} and define 
$$
\cM^*_{3,d}(\pi) := \sqrt{[\partial_{x_1} \pi]^2+\cdots+[\partial_{x_d} \pi]^2}.
$$
Then $\sqrt d \cM_{3,d}^*(\pi) \in \cA'_\pi$, but we have found that for each $K>0$ there exists $\pi \in \H_{3,d}$ such that 
$$
\det \cM^*_{3,d}(\pi)> K \inf \{\det M\sep M\in \cA'_\pi\}.
$$
The map $\cM^*_{3,d}$ may still be used for mesh adaptation 
through the formula \iref{defHd} but this metric may not be 
optimal in the area where $\pi_z=\frac {d^3f(z)} 6$ is such that
$\det \cM^*_{3,d}(\pi_z)$ is not well controlled by $\inf \{\det M\sep M\in \cA'_{\pi_z}\}$.

\section{Final remarks and conclusion}
\label{secW1P7}

In this chapter, we have introduced asymptotic estimates for the
finite element interpolation error measured in the $W^{1,p}$ semi-norm, when the mesh is 
optimally adapted to a function of two variables and the degree of interpolation $m-1$ is arbitrary.
The approach used is an adaptation of the ideas developped in Chapter \ref{chapOptAniso} for the $L^p$ interpolation error, and leads to asymptotically sharp error estimates, exposed in Theorems \ref{mainTheorem} and \ref{optiTheorem}.
These estimates involve a shape function $L_{m,p}$ which generalises the determinant which appears in estimates for piecewise linear interpolation. The shape function has equivalents of polynomial 
form for all values of $m$, as established in theorems \ref{thPolEq} and \ref{thPolEq2}. 
Up to a fixed multiplicative constant, our estimates can therefore be written 
under analytic form in terms of the derivatives of the function to be approximated.
 
In the case of piecewise linear and piecewise quadratic finite elements, we have presented in \S \ref{secW1P4} metrics which allow to produce near optimal meshes. This metric is new in the case of quadratic elements.
Some numerical experiments presented in \S \ref{secW1P4p2} illustrate the efficiency of this procedure, and the C++ source code is freely available on the internet \cite{sitejm}.

We have partially extended these results to higher dimension, in particular we provide a local error estimate \iref{localMixedD} which leads to sufficient conditions for building meshes that satisfy the best possible estimate up to a multiplicative constant.
A multidimensional asymptotical lower error estimate is proved in Theorem \ref{optiTheoremD} and generalises the bidimensional study. The corresponding asymptotical upper estimate is presented in \ref{upperPropD} but not proved.

One of the main tools used throughout this chapter 
for the construction of an optimal partition 
is the measure of sliverness $S(T)$ of a simplex, defined in \iref{defSD}, 
which has a geometrical interpretation as the distance from $T$ to the set of acute simplices. 
This measure accurately distinguishes between good anisotropy, 
that leads to optimal error estimates, and bad anisotropy 
that leads to oscillation of the gradient of the interpolated function. 
Equivalent quantities can be found in \cite{Ba,Ja}, but had not been used in the context of optimal mesh adaptation.

\section{Appendix}
\subsection{Proof of Lemma \ref{lemmaTile}}

\begin{figure}[h]
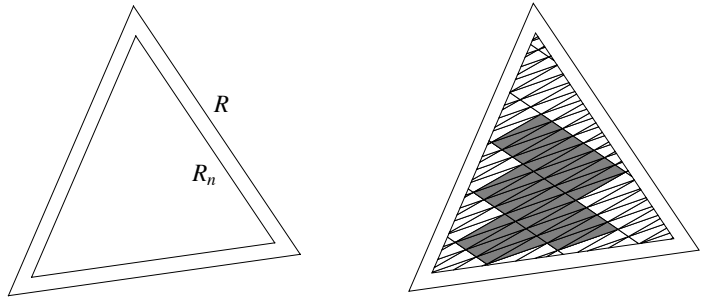

	\centering
		\includegraphics[width=4cm,height=4cm]{\pathPic/PaperW1P/RandRn.pdf}
		\hspace{1cm}
		\includegraphics[width=4cm,height=4cm]{\pathPic/PaperW1P/TilingPure.pdf}
	\caption{\label{fig5W1P}Left : The triangles $R$ and $R_n$. Right : The partition $\cP'_n$ of $R_n$.} 
\end{figure}

\begin{figure}[b]
	\centering
		\includegraphics[width=4cm,height=4cm]{\pathPic/PaperW1P/TilingBlue.pdf}
		\hspace{1cm}
		\includegraphics[width=4cm,height=4cm]{\pathPic/PaperW1P/TilingRed.pdf}
	\caption{\label{fig6W1P}Left : detail of the partition $\cP''_n$ of $R_n$. \ Right :  the partition $\cP_n = \cP_{T,n}(R) = \cP''_n \cup \widetilde \cP_n$ of $R$.} 
\end{figure}

Let $R_n$ be the homothetic contraction of $R$ by the factor $1-n^{-1}$ and with the same barycenter. We define a partition $\cP'_n$ of $R_n$ into convex polygons as follows 
$
\cP'_n := \{ R_n \cap T' \sep T'\in \cP_{T,n}\}.
$
The triangles $R$, $R_n$, and the partition $\cP'_n$ are illustrated on Figure \ref{fig5W1P}.
Note that the normals to the faces of the polygons building the partition $\cP'_n$ belong to a family of only $6$ elements $(\n_i)_{1\leq i\leq 6}$ : the normals to the faces of of $R$, and the normals to the faces of $T$. Hence only  $6\times 5$ different angles can appear in $\cP'_n$, and we denote the largest of these by $\alpha < \cPi$.

We now partition into triangles each convex polygon $C\in \cP'_n$ using the Delaunay triangulation of its vertices. Note that the angles of the triangles partitionning a convex polygon $C$ are smaller than the maximal angle of $C$, hence than $\alpha$.
We denote by $\cP''_n$ the resulting triangulation of $R_n$, as illustrated on the left of Figure \ref{fig6W1P}.

We denote by $E_n$ the collection of $n$ equidistributed points on each edge of $R$, described in item 2 of Lemma \ref{lemmaTile}.
We denote by $E'_n$ the set of vertices of the triangles in $\cP''_n$ that fall on $\partial R'_n$. 
For each point $p\in E_n$, we draw an edge between $p$ and the point of $E'_n$ which is the closest to $p$. This produces a partition of $R\sm R_n$ into triangles and convex quadrilaterals. 
Eventually we partition each of these polygons $C$ into triangles using the Delaunay triangulation of the point set $\overline C\cap (E_n \cup E'_n)$, which produces a triangulation $\widetilde \cP_n$ of $R\sm R_n$ illustrated on the right of Figure \ref{fig6W1P}. 
The triangles $T'\in \widetilde \cP_n$ obey 
$$
\diam(T')\leq (2\diam R +\diam T) n^{-1} = C n^{-1}.
$$ 
Furthermore let $L$ be the length of the edge of $T'$ included in $\partial R\cup \partial R_n$, and let $H$ be the height of the triangle $T'$  such that $L H = 2 |T'|$. 
Then 
$$
H\geq \min \{ |z-z'|\sep z\in \partial R, \ z'\in \partial R_n\} = cn^{-1}.
$$
where $c>0$ is independent of $n$.
Let $L'$ be another edge of $T'$, and let $\theta$ be the angle of $T'$ between the edges $L$ and $L'$. Then
$$
2 |T'| = L L' \sin \theta = L H.
$$
hence $\sin \theta \geq \frac H {\diam(T')} \geq \frac c C$, and therefore $\arcsin(\frac c C) \leq \theta \leq \pi - \arcsin(\frac c C)$. It follows that all the angles of $T'$ are smaller than $\pi - \arcsin(\frac c C)$.
We eventually define $\cP_{T,n}(R) := \cP''_n\cup \widetilde \cP_n$ and we observe that the largest angle of a triangle in $\cP_{T,n}(R)$ is bounded by the constant $\beta(R,T) = \max\{\alpha, \pi-\arcsin(\frac c C)\}<\cPi$ which is independent of $n$. 
Hence 
$$
\sup_{n\geq 1} \sup_{T\in \cP_{T,n}(R)} S(T) \leq \tan \left(\frac {\beta(R,T)} 2 \right) <\infty
$$
The other properties of $\cP_{T,n}(R)$ mentionned in \ref{lemmaTile} are easily checked.

\subsection{Proof of Theorem \ref{thPolEq2}}

Let $m\geq 2$ be arbitrary 
and let $s_m := \lfloor\frac m 2\rfloor +1$. We have proved in chapter \S\ref{chapOptAniso}, Proposition \ref{vanishprop}, that for all $\pi \in \H_m$ the three following properties are equivalent  
\be
\label{vanishKE}
\left [ 
\begin{array}{l}
K_m^\cE(\pi) = 0, \\
\text{There exists } \alpha, \beta\in \R \text{ and } \tilde \pi \in \H_{m-s_m} \text{ such that } \pi = (\alpha x+\beta y)^{s_m} \tilde \pi,\\
\text{There exists a sequence } (\phi_n)_{n\geq 0}, \phi_n \in \SL_2 \text{ such that } \pi\circ \phi_n\to 0.
\end{array}
\right.
\ee
We also proved in chapter \S\ref{chapOptAniso}, Theorem \ref{ThEquivMd}, the following invariance property : Let $Q$ be a polynomial on $\H_m$ such that $K_m^\cE \sim \sqrt[r] {|Q|}$ where $r=\deg Q$. Then
\be
\label{invQ}
Q(\pi\circ \phi) = (\det \phi)^{\frac{rm} 2}  Q(\pi) \text{ for all }  \pi \in \H_m \text{ and } \phi\in \M_2(\R).
\ee
It immediately follows that the polynomials $(Q_k)_{0\leq k \leq r}$, defined in \iref{defQk}, satisfy for all $\pi_1,\pi_2\in \H_m$ and $\phi\in \M_2(\R)$
\be
\label{invQk}
Q_k(\pi_1\circ \phi, \pi_2\circ\phi) = (\det\phi)^{\frac{rm} 2} Q_k(\pi_1, \pi_2) \text{ for all } \pi_1, \pi_2 \in \H_m \text{ and all } \phi\in \M_2(\R).
\ee
We define two functions on $\H_m \times \H_m$
$$
\KEq(\pi_1,\pi_2) := \sqrt[2r] {\sum_{0\leq k\leq r} Q_k(\pi_1,\pi_2)^2} \text{ and } K(\pi_1,\pi_2) := \sqrt[2 \tilde r]{\tilde Q(\pi_1^2+\pi_2^2)},
$$
where $\tilde Q$ is such that $K_{2m}^\cE\sim \sqrt[\tilde r]{\tilde Q}$, $\tilde r := \deg \tilde Q$. 
We show below that $K\sim K_*$ on $\H_m \times \H_m$. This result combined with  \iref{equivLK} concludes the proof of Theorem \ref{thPolEq2}. (Note that $m$ is replaced with $m+1$ in the statement of this theorem.)
Using \iref{invQk} and remarking the invariance property $\ti Q (\pi \circ \phi) = (\det \phi)^{\ti r m} Q(\pi)$, for the same reasons as \iref{invQ}, we obtain
\be
\label{invPair}
\text{for all } \pi_1, \pi_2 \in \H_m \text{ and all } \phi\in \M_2(\R) , 
\left\{\begin{array}{lcl}
K(\pi_1\circ \phi, \pi_2\circ\phi) &=& |\det \phi|^{\frac m 2} K(\pi_1, \pi_2),\\
\KEq(\pi_1\circ \phi, \pi_2\circ\phi) &=& |\det \phi|^{\frac m 2} \KEq(\pi_1, \pi_2).
\end{array}\right.
\ee

If $K(\pi_1, \pi_2) = 0$, then $\pi_1^2+ \pi_2^2\in \H_{2m}$ has a linear factor of multiplicity $s_{2m} = m+1$  according to \iref{vanishKE}, and therefore $\pi_1$ and $\pi_2$ have a common linear factor of multiplicity $s_m$. 

If $\KEq(\pi_1,\pi_2) = 0$, then for all $k$, $0\leq k\leq r$, we have $Q_k(\pi_1,\pi_2) = 0$. Using \iref{defQk} we obtain that for all $u,v\in \R$ we have $K_m(u \pi_1+ v\pi_2) = 0$. It follows from \iref{vanishKE}  that for all $u,v\in \R$ the polynomial $u\pi_1+ v\pi_2\in \H_m$ has a linear factor of multiplicity $s_m$, hence that $\pi_1$ and $\pi_2$ have a common linear factor of multiplicity $s_m$. 

Hence the following properties are equivalent 
\be
\label{commonRoot}
\left[\begin{array}{l}
K(\pi_1,\pi_2) = 0,\\
\KEq(\pi_1,\pi_2) = 0,\\
\text{There exists } \alpha, \beta\in \R \text{ and } \tilde \pi_1, \tilde \pi_2 \in \H_{m-s_m} \\
\qquad \text{ such that } \pi_1 = (\alpha x+\beta y)^{s_m} \tilde \pi_1,
\text{ and } \pi_2 = (\alpha x+\beta y)^{s_m} \tilde \pi_2.\\
\end{array}\right.
\ee
Using \iref{vanishKE} we find that these properties are also equivalent to 
\be
\label{commonRoot2}
\left[\begin{array}{l}
K_{2m}^\cE(\pi_1^2+\pi_2^2) = 0,\\
\text{There exists a sequence } (\phi_n)_{n\geq 0}, \phi_n \in \SL_2, \text{ such that } (\pi_1\circ \phi_n)^2 + (\pi_2\circ \phi_n)^2 \to 0,\\
\text{There exists a sequence } (\phi_n)_{n\geq 0}, \phi_n \in \SL_2, \text{ such that } \pi_1\circ \phi_n \to 0 \text{ and }\pi_2\circ \phi_n \to 0.\\
\end{array}\right.
\ee
We now define the norm $\|(\pi_1,\pi_2)\| := \sup_{|u|\leq 1} |(\pi_1(u),\pi_2(u))|$ on $\H_m\times \H_m$ and 
$$
\cF := \{(\pi_1, \pi_2)\in \H_m\times \H_m\sep \|(\pi_1, \pi_2)\| =1 \text{ and } \|(\pi_1\circ\phi, \pi_2\circ\phi)\| \geq 1 \text{ for all } \phi\in \SL_2\}.
$$ 
$\cF$ is compact subset of $\H_m \times \H_m$ and $K$ as well as $K_*$ do not vanish on $\cF$ according to \iref{commonRoot} and \iref{commonRoot2}. Since these functions are continuous, there exists a constant $C_0>0$ such that 
\be
\label{equivKA}
C_0^{-1} K\leq K_*\leq C_0 K \text{ on } \cF.
\ee
Let $\pi_1, \pi_2\in \H_m$. If there exists a sequence $(\phi_n)_{n\geq 0}$, $\phi_n \in \SL_2$, such that $\pi_1\circ \phi_n \to 0$ and $\pi_2 \circ \phi_n \to 0$, then $K(\pi_1, \pi_2) = \KEq(\pi_1, \pi_2) = 0$. Otherwise, consider a sequence $(\phi_n)_{n\geq 0}$, $\phi_n\in\SL_2$ such that 
$$
\lim_{n\to \infty}\|(\pi_1\circ \phi_n, \pi_2\circ \phi_n)\| = \inf_{\phi\in \SL_2} \|(\pi_1\circ\phi, \pi_2\circ\phi)\|.
$$
By compactness there exists a pair $(\tilde \pi_1, \tilde \pi_2)\in \H_m \times \H_m$ and a subsequence $(\phi_{n_k})_{k\geq 0}$ such that  
$$
(\pi_1\circ \phi_{n_k}, \pi_2\circ \phi_{n_k})\to (\tilde \pi_1, \tilde \pi_2).
$$ 
One easily checks that $\frac{(\ti \pi_2, \ti \pi_2)}{\|(\ti \pi_2,\ti \pi_2)\|}\in \cF$.
Using \iref{invPair} we obtain
$$
\frac{K(\pi_1,\pi_2)}{K_*(\pi_1,\pi_2)} =
\lim_{n\to \infty}  \frac{K(\pi_1 \circ\phi_n,\pi_2\circ\phi_n)}{K_*(\pi_1\circ\phi_n,\pi_2\circ\phi_n)} = \frac{K(\tilde\pi_1,\tilde\pi_2)}{K_*(\tilde \pi_1,\tilde \pi_2)} 
$$
Using  \iref{equivKA} and the homogeneity of $K$ and $K_*$, we obtain that $C_0^{-1} K\leq K_*\leq C_0 K$ on $\H_m \times \H_m$ which concludes the proof.

\subsection{Proof of Proposition \ref{propMS}}

We denote by $\TEq$ a $d$-dimensional equilateral simplex such that $\bary(\TEq) = 0$, where $\bary$ denotes the barycenter, and such that its vertices  $q_i$, $0\leq i\leq d$, belong to the unit sphere, i.e. $|q_i| = 1$. Since the vertices of $\TEq$ play symmetrical roles there exists a constant $\xi \in \R$ such that 
\be
\label{eqTEq}
\text{ For all } \ 0\leq i\leq d , \ \ 0 \leq j\leq d, \  \ 0 \leq k \leq d, \text{ one has } \<q_i-q_j,q_k\> = \xi(\delta_{ik} - \delta_{jk}),
\ee
where $\delta$ is the Kronecker symbol : $\delta_{ij} = 1$ if $i=j$, and $0$ otherwise. 
Using the relation $q_0+ \cdots + q_d = 0$ we obtain $\xi d = \sum_{j=0}^d \<q_0-q_j,q_0\> = d+1$ hence $\xi = 1+ \frac 1 d$.
Note also that the unit exterior normal to the face of $\TEq$ opposite to the vertex $q_i$ is $-q_i$. 

We recall the following property : if $A\in \GL_d$ and if $\n$ is the exterior normal to a face $F$ of a simplex $T$, then the exterior normal to the face $A(F)$ of $A(T)$ is 
\be
\label{extNorm}
\n ' = 
\frac{(A^{-1})^\trans \n}{\left|(A^{-1})^\trans \n\right|}.
\ee 


We first establish that for any simplex $T$, the simplex $M_T^{-\frac 1 2} (T)$ is acute.
Without loss of generality we can assume that $\bary(T) = 0$, hence there exists $A\in \GL_d$ such that $T = A(\TEq)$.
Since the vertices of $T$ are $v_i = A q_i$ for $0\leq i\leq d$, we obtain from definition \iref{defM} that  
$$
M_T = A \left(\sum_{0\leq i< j \leq d}  \frac{(q_i - q_j) (q_i -q_j)^\trans}{|A(q_i-q_j)|^2}\right) A^\trans.
$$
According to \iref{extNorm}, the exterior normal to the face of the simplex 
$\TAc := M_T^{-\frac 1 2}(T) = M_T^{-\frac 1 2} A(\TEq)$
opposite to the vertex $v_i$ is 
$$
\n_i = -\nu_i M_T^{\frac 1 2} (A^{-1})^\trans q_i
$$ 
where $\nu_i >0$.
For all $0\leq a< b\leq d$, we therefore obtain using \iref{eqTEq}
\begin{eqnarray*}
\nu_a \nu_b \<\n_a,\n_b\> &=& \left\<M_T^{\frac 1 2} (A^{-1})^\trans q_a, \ M_T^{\frac 1 2} (A^{-1})^\trans q_b\right\> \\
&=& q_a^\trans A^{-1} M_T (A^{-1})^\trans q_b \\
&=& q_a^\trans  \left(\sum_{0\leq i< j \leq d}  \frac{(q_i - q_j) (q_i -q_j)^\trans}{|A(q_i-q_j)|^2}\right) q_b\\
&=& \xi^2 \sum_{0\leq i< j \leq d}  \frac{(\delta_{a i} - \delta_{a j}) (\delta_{b i} - \delta_{b j})}{|A(q_i-q_j)|^2}\\
&=& \frac{-\xi^2}{|A(q_a-q_b)|^2} <0.
\end{eqnarray*}
This establishes that the simplex  $\TAc := M_T^{-\frac 1 2}(T)$ is acute, and therefore that $S(T) \leq \kappa(\sqrt{M_T})$ since $1 \leq \|\sqrt{M_T}\| \leq \alpha_d$.

The rest of this appendix is devoted to the proof that $\|M_T^{-\frac 1 2}\|\leq S(T)$, which implies that $\kappa(\sqrt{M_T})\leq \alpha_d S(T)$ and thus concludes the proof of Theorem \ref{propMS}. For this we need the following lemma.

\begin{lemma}
For any acute simplex $\TAc$ one has $M_{\TAc} \geq \Id$.
\end{lemma}
\proof
Without loss of generality we can assume that $\bary(\TAc) = 0$, hence there exists $A\in \GL_d$ such that $\TAc = A(\TEq)$. The vertices of $\TAc$ are $c_i = A q_i$ for $0\leq i\leq d$, and the exterior normal to the face of $\TAc$ opposite $c_i$ is 
$$
\m_i = - \mu_i (A^{-1})^\trans q_i.
$$
where $\mu_i>0$.
We define for all $0\leq i<j\leq d$
$$
\lambda_{i j} := \frac{-|c_i - c_j|^2 \<\m_i, \m_j\>}{\xi^2\mu_i \mu_j}.
$$ 
Since $\TAc$ is acute we have $\<\m_i,\m_j\>\leq 0$ and therefore $\lambda_{i j}\geq 0$. We now introduce the symmetric matrix
\be
\label{deff}
M := \sum_{0\leq i<j\leq d} \lambda_{i j} f_{i j} f_{i j}^\trans \text{ where } f_{i j} := \frac{c_i-c_j}{|c_i-c_j|}.
\ee
For all $0\leq a<b\leq d$ we obtain using the relation $\<\m_i, c_j\> = - \mu_i\<q_i,q_j\>$ that 
\begin{eqnarray*}
\m_a^\trans M \m_b&=&\mu_a \mu_b  \sum_{0\leq i<j\leq d}  \lambda_{ij} \frac{\<q_a, q_i - q_j\> \<q_b,q_i - q_j\>}{|c_i-c_j|^2} \\
&=& \mu_a \mu_b \xi^2 \sum_{0\leq i<j\leq d}  \lambda_{i j} \frac{(\delta_{a i}-\delta_{a j})(\delta_{b i} - \delta_{b j})}{|c_i-c_j|^2} \\
&=& \frac{-\mu_a \mu_b \xi^2\lambda_{a b}}{|c_a-c_b|^2} =  \<\m_a, \m_b\>.
\end{eqnarray*}
Therefore $\m_a^\trans M \m_b = \m_a^\trans \m_b$ for all $0\leq a<b\leq d$, which implies that $M = \Id$. 
Furthermore for all $0\leq a<b\leq d$, we have 
$$
1 = |f_{a b}|^2 = f_{a b}^\trans M f_{a b}  = \sum_{0\leq i<j\leq d} \lambda_{i j} \< f_{a b}, \, f_{i j}\>^2 \geq \lambda_{a b}.
$$
It follows that in the sense of symmetric matrices, 
$$
M_\TAc := \sum_{0\leq i<j\leq d} f_{i j} f_{i j}^\trans \geq \sum_{0\leq i<j\leq d} \lambda_{i j} f_{i j} f_{i j}^\trans =M = \Id,
$$
which concludes the proof of this lemma.
\sq

We now conclude the proof of inequality \iref{eqMS}.
Let $T$ be an arbitrary simplex, and 
 let $\psi\in \GL_d$ be such that the simplex $\TAc := \psi(T)$ is acute. Let $v_i$, $0\leq i\leq d$ be the vertices of $T$ and $c_i = \psi(v_i)$ the vertices of $\TAc$. We define the vectors $e_{i j}$ and $f_{i j}$ similarly to \iref{defM} and \iref{deff}
$$
e_{i j} := \frac{v_i - v_j} {|v_i-v_j|}, \quad f_{i j} := \frac{c_i-c_j}{|c_i-c_j|} = \frac{\psi(e_{i j})}{|\psi(e_{i j})|}.
$$
for all $0\leq i<j\leq d$. 
For any $v\in \R^d$ we therefore have 
\begin{eqnarray*}
 v^\trans M_{\TAc} v &=&  \sum_{0\leq i<j\leq d} \frac{\<\psi(e_{i j}),v\>^2}{|\psi(e_{i j})|^2} \\
&\leq &  \|\psi^{-1}\|^2 \sum_{0\leq i<j\leq d} \<\psi(e_{i j}),v\>^2 \\
&=& \|\psi^{-1}\|^2 \, (\psi^\trans v)^\trans M_T (\psi^\trans v).
\end{eqnarray*}
Using the previous lemma and defining $u:=\psi^\trans v$ we obtain 
$$
|u|^2 \leq \|\psi\|^2 |v|^2 \leq \|\psi\|^2 v^\trans M_\TAc v \leq \|\psi\|^2 \|\psi^{-1}\|^2 u^\trans M_T u = \left(\|\psi\| \|\psi^{-1}\| | M_T^{\frac 1 2} u|\right)^2,
$$
hence $\|M_T^{\frac{-1} 2}\| \leq \kappa(\psi)$.
Recalling that 
$\|\sqrt{M_T}\|\leq \alpha_d$ we obtain
$$
\kappa(\sqrt{M_T}) = \|M_T^{\frac 1 2}\| \|M_T^{\frac{-1} 2}\| \leq \alpha_d \kappa(\psi).
$$
We conclude the proof by taking the infimum among all $\psi$ such that $\psi(T)$ is acute.


\part{Anisotropic smoothness classes}
\label{partCartoon}

\chapter{From finite element approximation to image models}
\label{chapCartoon}
\minitoc

\section{Introduction}
There exists various ways of measuring the smoothness of functions
on a domain $\Omega\subset \RR^d$, 
generally through the definition of an appropriate {\it smoothness space}. 
Classical instances are Sobolev, H\"older and Besov spaces. Such spaces are of common use
when describing the regularity of solutions to partial differential equations.
From a numerical perspective, they are also useful in order to sharply characterize at which
rate a function $f$ may be approximated by simpler functions such as Fourier series,
finite elements, splines or wavelets (see \cite{Co,DL,De} for surveys on
such results).

Functions arising in concrete applications may have
inhomogeneous smoothness properties, in the sense that they exhibit
area of smoothness separated by localized discontinuities. Two typical
instances are (i) edge in functions representing real images and 
(ii) shock profiles in solutions to non-linear hyperbolic PDE's.
The smoothness space that is best taylored to take such features
into account is the space $BV(\Omega)$ of bounded variation functions. This space
consists of those $f$ in $L^1(\Omega)$ such that $\nabla f$ is a bounded measure,
i.e. such that their total variation
$$
\TV(f)=|f|_{BV}:=\max \left\{\int_\Omega f {\rm div}(\vp)\; ; \; \vp\in\cD(\Omega)^d,\; \|\vp\|_{L^\infty}\leq 1\right\}
$$
is finite. Functions of bounded variation are allowed to have jump discontinuities
along hypersurfaces of finite measure. In particular, the characteristic function of a smooth
subdomain $D\subset \Omega$ has finite total variation equal to the $d-1$-dimensional
Hausdorff measure of its boundary:
\be
|\Chi_D|_{BV}=\cH_{d-1}(\partial D).
\label{haus}
\ee
It is well known that $BV$ is a regularity space for certain 
hyperbolic conservation laws \cite{GR,Le}, in the sense that the total variation
of their solutions remains finite for all time $t>0$. This space also plays an
important role in image processing since the seminal paper \cite{FOR}.
Here, a small total variation is used as a prior to describe the mathematical 
properties of ``plausible images'', when trying to restore an unknown
image $f$ from an observation $h=Tf+e$ where $T$ is a known operator
and $e$ a measurement noise of norm $\|e\|_{L^2}\leq \e$. The restored image
is then defined as the solution to the minimization problem
\be
\min_{g\in BV} \{|g|_{BV}\; ; \; \|Tg-h\|_{L^2}\leq \e\}.
\label{minBV}
\ee
From the point of view of approximation theory, it was shown in \cite{CDPX,CDDD2} that
the space $BV$ is almost characterized by expansions in wavelet bases. For example,
in dimension $d=2$, if $f=\sum d_\lambda\psi_\lambda$ 
is an expansion in a tensor-product $L^2$-orthonormal wavelet basis, one has
$$
(d_\lambda)\in \ell^1 \Rightarrow f\in BV \Rightarrow (d_\lambda)\in w\ell^1,
$$
where $w\ell^1$ is the space of weakly summable sequences. The fact that 
the wavelet coefficients of a $BV$ function are weakly summable implies
the convergence estimate
\be
\|f-f_N\|_{L^2}\leq CN^{-1/2}|f|_{BV},
\label{wavN}
\ee
where $f_N$ is the {\it nonlinear} approximation of $f$ obtained by retaining the
$N$ largest coefficients in its wavelet expansion. Such approximation results 
have been further used in order to justify the performance of compression or
denoising algorithms based on wavelet thresholding \cite{CDDD1,DJKP,Do}.

In recent years, it has been observed that the space $BV$ (and more generally
classical smoothness spaces) do not provide a fully satisfactory description of
piecewise smooth functions arising in the above mentioned applications.
Indeed, formula \iref{haus} reveals that the total variation only takes into account
the {\it size} of the sets of discontinuities and not their geometric {\it smoothness}.
In image processing, this means that the set of bounded variation images does
not make the distinction between smooth and non-smooth edges as long as they have finite length.

The fact that edges have some geometric smoothness can be exploited in order to study
approximation procedures which outperform wavelet thresholding in terms of 
convergence rates. For instance, it is easy to prove that if $f=\Chi_D$ 
where $D$ is a bidimensional domain with smooth boundary, one can find a sequence
of triangulations $\cT_N$ with $N$ triangles such that the convergence estimate
\be
\|f-\interp_{\cT_N}f\|_{L^2} \leq CN^{-1},
\label{triN}
\ee
holds, where $\interp_\cT$ denotes the 
piecewise linear interpolation operator on a triangulation $\cT$.
Other methods are based on thresholding a decomposition of the 
function in bases or frames which differ from classical wavelets, see e.g. \cite{CD,LM,ACDDM}.
These methods also yield improvements over \iref{wavN} similar to \iref{triN}. The common feature in
all these approaches is that they achieve 
{\it anisotropic refinement} near the edges. For example,
in order to obtain the estimate \iref{triN}, 
the triangulation $\cT_N$ should include a thin layer of triangles
which approximates the boundary $\partial D$. These triangles typically have
size $N^{-2}$ in the normal direction to $\partial D$ and $N^{-1}$ in the tangential direction,
and are therefore highly anisotropic.  

Intuitively, these methods are well adapted to
functions which have anisotropic smoothness properties
in the sense that their local variation is significantly stronger
in one direction. Such properties are not well described by
classical smoothness spaces such as $BV$, and a natural
question to ask is therefore:
\nl
\nl
{\it What type of smoothness properties 
govern the convergence rate of anisotropic refinement
methods and how can one quantify these properties ?}
\nl
\nl
The goal of this chapter is to answer this question,
by proposing and studying measures of smoothness which are suggested
by recent results on anisotropic finite element approximation
\cite{BBLS, CSX} and Chapter 2. Before going further, 
let us mention several existing approaches which have been developed
for describing and quantifying anisotropic smoothness,
and explain their limitations.
\begin{enumerate}
\item
The so-called {\it mixed smoothness} classes have been introduced 
and studied in order to
describe functions which have a different order of
smoothness in each coordinate, see e.g. \cite{Ni,Te}. These spaces
are therefore not adapted to our present goal since
the anisotropic smoothness that we want to describe 
may have preferred directions that are not 
aligned with the coordinate axes and that may vary from
one point to another (for example an image with a curved edge).
\item
Anisotropic smoothness spaces
with more general and locally varying directions have
been investigated in \cite{KP}. Yet, in such spaces
the amount of smoothness in different directions
at each point is still fixed in advance and therefore
again not adapted to our goal, since this amount
may differ from one function to another (for example
two images with edges located at different positions).
\item
A  class of functions which is often used to 
study the convergence properties of anisotropic approximation
methods is the family of $C^m-C^n$ {\it cartoon images}, i.e.
functions which are $C^m$ smooth on a finite number
of subdomains $(\Omega_i)_{i=1,\cdots,k}$ separated by a union of
discontinuity curves 
$(\Gamma_{j})_{j=1,\cdots,l}$
that are $C^n$ smooth. The defects of this class are revealed
when searching for simple expression that quantifies the amount of
smoothness in this sense. A natural choice is
to take the supremum of all $C^m(\Omega_i)$ norms of $f$ and
$C^n$ norms of the normal parametrization of $\Gamma_{j}$.
We then observe that this quantity is unstable 
in the sense that it becomes extremely large
for blurry images obtained by convolving
a cartoon image by a mollifier $\vp_\delta=\frac 1{\delta^{2}}\vp(\frac \cdot h)$
as $\delta \to 0$.  In addition, this quantity does not
control the number of subdomains in the partition.
\item
A recent approach proposed in \cite{DPW} defines
anisotropic smoothness through the geometric 
smoothness properties of the {\it level sets} of the function $f$.
In this approach the measure of smoothness 
is not simple to compute directly from $f$ since it involves each of its level sets
and a smoothness measure of their local parametrization.
\end{enumerate}

The results of \cite{BBLS,CSX} and Chapter \ref{chapOptAniso} describe 
the $L^p$-error of piecewise linear interpolation by 
an optimally adapted triangulation of at most $N$ elements,
when $f$ is a $C^2$ function of two variables. This error is defined
as
$$
\sigma_N(f)_p:=\inf_{\#(\cT)\leq N}\|f-\interp_\cT f\|_{L^p}.
$$
It is shown in \cite{BBLS} for $p=\infty$
and in Chapter \ref{chapOptAniso} for all $1\leq p\leq \infty$ that
that 
\be
\limsup_{N\to +\infty} N\sigma_N(f)_p\leq CA_p(f),
\label{limsupest}
\ee
where $C$ is an absolute constant and
\be
A_p(f):=\|\sqrt {|{\rm det}(d^2f)|}\|_{L^\tau},\;\; \frac 1\tau:= 1+\frac 1 p.
\label{nonlinnorm}
\ee
Moreover, this estimate is known to be optimal in the sense that
$\liminf_{N\to +\infty} N\sigma_N(f)_p\geq cA_p(f)$ also holds, under some mild restriction
on the class of triangulations in which one selects the optimal one.
These results are extended in Chapter \ref{chapOptAniso} to the case of
higher order finite elements and space dimension $d>2$,
for which one can identify similar measures $f\mapsto A(f)$ governing the
convergence estimate. Such quantities thus constitute natural 
candidates to measure anisotropic smoothness properties. Note
that $A_p(f)$ is not a semi-norm due to the presence of the determinant
in \iref{nonlinnorm}, and in particular the
quasi-triangle inequality $A_p(f+g)\leq C(A_p(f)+A_p(g))$ 
does not hold even with $C>1$.

This chapter is organized as follows.
We begin in \S \ref{secCartoon2} by a brief account of the
available estimates on anisotropic finite element 
interpolation, and we recall in particular the argument that
leads to \iref{limsupest} with the quantity $A_p(f)$ defined by \iref{nonlinnorm}.

Since $A_p(f)$ is not a norm, we cannot associate
a linear smoothness space to it by a standard completion
process. We are thus facing a difficulty in extending
the definition of $A_p(f)$ to functions which are not
$C^2$-smooth and in particular to cartoon images
such as in item 3 above. Since we know from \iref{triN}
that for such cartoon images the $L^2$ error of adaptive
piecewise linear interpolation decays like $N^{-1}$, we would expect that the
quantity 
$$
A_2(f)=\|\sqrt {|{\rm det}(d^2f)|}\|_{L^{2/3}},
$$
corresponding to the case $p=2$ can be properly defined for
piecewise smooth functions. We address this difficulty in 
\S \ref{secCartoon3} by a regularization process:
if $f$ is a cartoon image we introduce its regularized version
\be
f_\delta:=f*\vp_{\delta},
\ee
where $\vp_\delta=\frac 1{\delta^{2}}\vp(\frac \cdot h)$ 
is a standard mollifier. Our
main result is the following: for any cartoon image $f$ of $C^2-C^2$ type, the quantity $A_2(f_\delta)$
remains uniformly bounded as $\delta\to 0$ and one has
\be
\lim_{\delta \to 0} A_2(f_\delta)^{2/3}= \sum_{i=1}^k \int_{\Omega_i}\left |\sqrt {|{\rm det}(d^2f)|}\right |^{2/3}
+C(\vp)\sum_{j=1}^l \int_{\Gamma_j} |[f](s)|^{2/3} |\kappa(s)|^{1/3} ds,
\label{lima2}
\ee
where $[f](s)$ and $\kappa(s)$ respectively denote
the jump of $f$ and the curvature of $\Gamma_j$ at the point $s$,
and where $C(\vp)$ is a constant that only depends on the choice of the mollifier.
This constant can be shown to be uniformly bounded by below for the class of
radially decreasing mollifiers. 
This result reveals that $A_2(f)$ is stable under regularization of cartoon images
(in contrast to the measure of smoothness described in item 3 above). 
We also discuss the behaviour of $A_p$ when $p\neq 2$.

These results lead us in \S \ref{secCartoon4} to a comparison
between the quantity $A_2(f)$ and the total variation $\TV(f)$. We also
make some remarks on the existing links between the limit expression
in \iref{lima2} and classical results on adaptive approximation
of curves, as well as
with operators of affine-invariant image processing which also
involve the power $1/3$ of the curvature. 

We devote \S \ref{secCartoon5} to numerical tests performed on cartoon images that illustrate the validity 
of our results. In addition, we compare the quantity $A_2(f)$
with the total variation $\TV(f)$ as a model for plausible images.

One could be tempted to use the quantity $A_2(f)$
in place of the total variation $\TV(f)$ as a prior 
in image processing, and in particular to replace $|g|_{BV}$
by $A_2(g)$ in a restoration procedure such as \iref{minBV}. 
In \S \ref{secCartoon5.5}, we show that this optimization program is
ill-posed in the sense that the restored image would be the noisy image itself.
An alternate strategy is to use $A_2(f)$ in the framework of a
bayesian least-square estimator, as proposed in \cite{LMo} in the
case of the total variation. At the present stage, the algorithm 
implementing this approach did not give satisfactory results due
to its very slow convergence. For this reason, we only present
some numerical illustration of this approach in a simplified one-dimensional setting.

Eventually we describe in \S \ref{secCartoon6} the extension of our results
to finite elements of higher degree and higher space dimensions.
Concluding remarks and perspectives of our work are given in \S \ref{secCartoon7}.

\section{Anisotropic finite element approximation}
\label{secCartoon2}
A standard estimate in finite element approximation states that 
if $f\in W^{2,p}(\Omega)$ then
$$
\|f-\interp_{\cT_h} f\|_{L^p} \leq C  h^{2}\|d^2f\|_{L^{p}},
$$
where $\cT_h$ is a triangulation of mesh size
$h:=\max_{T\in\cT_h} {\rm diam}(T)$. If we restrict our attention
to a family {\it quasi-uniform} triangulations, $h$ is linked with 
the complexity $N:=\#(\cT_h)$ according to
$$
C_1 h^{-2} \leq N \leq  C_2 h^{-2}
$$
Therefore, denoting by $\sigma^{\rm unif}_N(f)_{L^p}$ the $L^p$ approximation
by quasi-uniform triangulations of cardinality $N$, we can re-express the above
estimate as
\be
\sigma^{\rm unif}_N(f)_{L^p} \leq CN^{-1}\|d^2f\|_{L^p}.
\ee
In order to explain how this estimate can be improved
when using adaptive partitions, we first give some
heuristic arguments which are based on the assumption that on each triangle $T$
the relative variation of $d^2f$ is small so that it can be considered
as constant over $T$, which means that 
$f$ coincides with a quadratic function $q_T$ on each $T$.
Denoting by $\interp_T$ the local interpolation operator
on a triangle $T$ and by $e_{T}(f)_p:=\|f-\interp_T f\|_{L^p(T)}$
the local $L^p$ error, we thus have according to
this heuristics
$$
\|f-\interp_\cT f\|_{L^p} =\(\sum_{T\in \cT} e_T(f)_p^p\)^{\frac 1 p}=\(\sum_{T\in \cT} e_T(q_T)_p^p\)^{\frac 1 p}
$$
We are thus led to study the local interpolation error $e_T(q)_p$ when $q\in \P_2$
is a a quadratic polynomial. Denoting by $\bq$ the homogeneous part of $q$, we
remark that
$$
e_T(q)_p=e_T(\bq)_p.
$$
We optimize the shape of $T$ with respect to the quadratic form $\bq$
by introducing a function $K_p$ defined on the space of quadratic forms by
$$
K_p(\bq):=\inf_{|T|=1}e_T(\bq)_p,
$$
where the infimum is taken among all triangles of area $1$. 
It is easily seen that $e_T(\bq)_p$ is invariant
by translation of $T$ and so is therefore the minimizing triangle
if it exists. By homogeneity, it is also easily seen
that 
$$
\inf_{|T|=a} e_T(\bq)_p=a^{\frac 1 \tau}K_p(\bq),\;\; \frac 1 \tau=\frac 1 p+1,
$$
and that the minimizing triangle of area $a$ is obtained by rescaling the minimizing
triangle of area $1$ if it exists. Finally, it is easily seen that if $\vp$ is an invertible
linear transform 
$$
K_p(\bq \circ \vp)=|{\rm det} (\vp)|K_p(\bq),
$$
and that the minimizing triangle of area $|{\rm det} (\vp)|^{-1}$ for $\bq \circ \vp$
is obtained by application of $\vp^{-1}$ to the minimizing triangle of 
area $1$ for $\bq$ if it exists. If ${\rm det} (\bq)\neq 0$, there exists 
a $\vp$ such that $\bq \circ \vp$ is either $x^2+y^2$ or $x^2-y^2$ up to 
a sign change, and we have $|{\rm det} (\bq)|=|{\rm det} (\vp)|^{-2}$.
It follows that $K_p(\bq)$ has the simple form
\be
K_p(\bq)=\sigma  |{\rm det} (\bq)|^{1/2},
\label{Kpexp}
\ee
where $\sigma$ is a constant 
equal to $K_p(x^2+y^2)$ 
if ${\rm det} (\bq)>0$ and to  $K_p(x^2-y^2)$ if ${\rm det} (\bq)<0$. 
One easily checks that this equality
also holds when ${\rm det} (\bq)= 0$ in which case $K_p(\bq)=0$.

Assuming that the triangulation $\cT$ is
such that all its triangles $T$ have optimized shape in the above sense
with respect to the quadratic form $\bq_T$ associated with $q_T$, we thus have
for any triangle $T\in\cT$
$$
e_T(f)_p=e_T(\bq_T)_p=|T|^{\frac 1 \tau}K_p(\bq_T)=\left\|K_p\(\frac {d^2f}2\)\right\|_{L^\tau(T)}.
$$
since we have assumed $\frac {d^2f} 2=\bq_T$ on $T$. In order to optimize the trade-off
between the global error and the complexity $N=\#(\cT)$,
we apply the principle of {\it error equidistribution}: the triangles $T$ have
area such that all errors $e_T(\bq_T)_p$ are equal
i.e. $e_T(\bq_T)_p= \eta$ for some $\eta>0$ independent of $T$.
It follows that
$$
N\eta^\tau\leq \left\|K_p\(\frac {d^2f}2\)\right\|_{L^\tau(\Omega)}^\tau,
$$
and therefore
$$
\sigma_N(f)_p \leq \|f-\interp_\cT f\|_{L^p} \leq N^{1/p} \eta \leq\left\|K_p\(\frac {d^2f}2\)\right\|_{L^\tau(\Omega)}N^{-1},
$$
which according to \iref{Kpexp} implies 
\be
\sigma_N(f)_p \leq CN^{-1}A_p(f),
\label{false}
\ee
with $A_p$ defined as in \iref{nonlinnorm}.

The estimate \iref{false} is too optimistic to be correct:
if $f$ is a univariate function then $A_p(f)=0$ while $\sigma_N(f)_p$
may not vanish. In a rigorous 
derivation such as in \cite{BBLS} and Chapter 2, one observes that
if $f\in C^2$, the replacement of $d^2f$ by a constant over $T$
induces an error which becomes negligible only when the
triangles are sufficiently small, and therefore a
correct statement is that for any $\e>0$ there exists 
$N_0=N_0(f,\e)$ such that 
\be
\sigma_N(f)_p \leq N^{-1}\left(\left\|K_p\(\frac {d^2f}2\)\right\|_{L^\tau(\Omega)}+\e\right),
\label{true}
\ee
for all $N\geq N_0$, i.e.
\be
\limsup_{N\to +\infty} N\sigma_N(f)_p\leq \left\|K_p\(\frac {d^2f}2\)\right\|_{L^\tau(\Omega)}
\label{limsupest1}
\ee
which according to \iref{Kpexp} implies \iref{limsupest}.

\section{Piecewise smooth functions and images}
\label{secCartoon3}
As already observed, the quantities $A_p(f)$ are well defined for
functions $f\in C^2$, but 
we expect that they should in some sense also be well defined for
functions representing $C^2-C^2$ ``cartoon images'' when $p\leq 2$. We first give 
a precise definition of such functions.

\begin{definition} 
\label{defcartoon}
A cartoon function on an open set $\Omega$ is a function almost everywhere of the form 
$$
f=\sum_{1\leq i\leq k} f_i \Chi_{\Omega_i},
$$
where the $\Omega_i$ are disjoint open sets with piecewise $C^2$ boundary, no cusps (i.e. satisfying an interior and exterior cone condition), and such that $\overline \Omega = \cup_{i=1}^k \overline \Omega_i$, and where the function $f_i$ is $C^2$ on $\overline \Omega_i$ for each $1\leq i\leq k$. 
\end{definition}

Let us consider a fixed cartoon function $f$ on an open polygonal domain $\Omega$ (i.e.
$\Omega$ is such that $\overline\Omega$ is a closed polygon)
associated with a decomposition $(\Omega_i)_{1\leq i\leq k}$. We define 
\be
\Gamma :=\bigcup_{1\leq i\leq k} \partial \Omega_i,
\ee
the union of the boundaries of the $\Omega_i$. The above definition
implies that $\Gamma$ is the disjoint union of a finite set of points $\cP$ and a finite number of open curves $(\Gamma_i)_{1\leq i\leq l}$.
$$
\Gamma = \(\bigcup_{1\leq i\leq l} \Gamma_i\) \cup \cP.
$$
Furthermore for all $1\leq i < j \leq l$, we may impose that $\overline \Gamma_i\cap \overline \Gamma_j \subset \cP$ (this may be ensured by a splitting of some of the $\Gamma_i$ if necessary).

We now consider the piecewise linear interpolation
$\interp_{\cT_N}f$ of $f$ on a triangulation $\cT_N$ of cardinality $N$. 
We distinguish two types of elements of $\cT_N$. A 
triangle $T\in \cT_N$ is called ``regular'' if $T\cap \Gamma=\emptyset$, 
and we denote the set of such triangles by $\cT_N^r$. 
Other triangles are called ``edgy'' and their set is denoted by $\cT_N^e$.
We can thus split $\Omega$ according to
$$
\Omega:=\left(\cup_{T\in \cT_N^r}T\right) \cup \left(\cup_{T\in \cT_N^e}T\right)=\Omega_N^r \cup \Omega_N^e.
$$
We split accordingly the $L^p$ interpolation error into
$$
\|f-\interp_{\cT_N}f\|_{L^p(\Omega)}^p = \int_{\Omega_N^r} |f-\interp_{\cT_N}f|^p +\int_{\Omega_N^e} |f-\interp_{\cT_N}f|^p.
$$
We may use $\cO(N)$ triangles in $\cT_N^e$ and $\cT_N^r$ (for example $N/2$
in each set). Since $f$ has discontinuities along $\Gamma$, the $L^\infty$ interpolation error on 
$\Omega_N^e$ does not tend to zero and $\cT_N^e$ should be chosen so 
that $\Omega_N^e$ has the aspect of a thin layer around $\Gamma$. 
Since $\Gamma$  is a finite union of $C^2$ curves, we can build this layer
of width $\cO(N^{-2})$ and therefore of global area $|\Omega_N^e|\leq C N^{-2}$,
by choosing long and thin triangles in $\cT_N^e$. 
On the other hand, since $f$ is uniformly $C^2$ on $\Omega_N^r$, we may
choose all triangles in $\cT_N^r$ of regular shape 
and diameter $h_T\leq C N^{-1/2}$.
Hence we obtain the following heuristic error estimate, 
for a well designed anisotropic triangulation:
\begin{eqnarray*}
\|f-\interp_{\cT_N}f\|_{L^p(\Omega)} &=& \(\|f-\interp_{\cT_N}f\|_{L^p(\Omega_N^r)}^p +\|f-\interp_{\cT_N}f\|_{L^p(\Omega_N^e)}^p\)^{1/p}\\
&\leq&\( \|f-\interp_{\cT_N}f\|_{L^\infty(\Omega_N^r)}^p |\Omega_N^r|+  
\|f-\interp_{\cT_N}f\|_{L^\infty(\Omega_N^e)}^p |\Omega_N^e|\)^{1/p}\\
& \leq & C(N^{-p}+ N^{-2})^{1/p},
\end{eqnarray*}
and therefore
\be
 \|f-\interp_{\cT_N}f\|_{L^p(\Omega)} \leq C N^{-\min\{1,2/p\}},
 \label{heurist}
 \ee
where the constant $C$ depends on $\|d^2f\|_{L^\infty(\Omega\sm \Gamma)}$, $\|f\|_{L^\infty(\Omega)}$
and on the number, length and maximal curvature of the $C^2$ curves which constitute $\Gamma$.

Observe in particular that the error
is dominated by the edge term $\|f-\interp_{\cT_N}f\|_{L^p(\Omega_N^e)}$
when $p>2$ and by the smooth term $\|f-\interp_{\cT_N}f\|_{L^p(\Omega_N^r)}$
when $p<2$. For the critical value $p=2$ the two terms have the same order. 

For $p\leq 2$, we obtain the approximation rate $N^{-1}$ 
which suggests that approximation results such as \iref{limsupest}
should also apply to cartoon functions and
that the quantity $A_p(f)$ should be
finite. We would therefore like to bridge the gap 
between anisotropic approximation of cartoon functions and smooth functions. 
For this purpose, we first need to give a proper meaning
to $A_p(f)$ when $f$ is a cartoon function. This is not 
straightforward, due to the fact that
the product of two distributions has no meaning in general.
Therefore, we cannot define $\det (d^2 f)$ in a distributional sense,
when the coefficients of $d^2 f$ are distributions without sufficient smoothness.
Our approach will rather be based on regularisation. This is additionally justified by the fact that sharp curves of discontinuity are a mathematical idealisation. In real world applications, such as photography, several physical limitations (depth of field, optical blurring) impose a certain level of blur on the edges.

In the following, we consider a fixed
radial nonnegative function $\vp$ of unit integral and supported in the unit ball,
and we define for all $\delta>0$ and $f$ defined on $\Omega$,
\be
\vp_\delta(z) := \frac 1 {\delta^2} \vp\left(\frac z \delta\right) \text{ and } f_\delta = f * \vp_\delta.
\label{defphih}
\ee
Our main result gives a meaning to $A_p(f)$ based on this regularization.
If $f$ is a cartoon function on a set $\Omega$, 
and if $x\in \Gamma\sm\cP$, we denote by $[f](x)$ the jump of $f$
at this point.
We also denote $\bt(x)$ and $\bn(x)$ the unit tangent and normal
vectors to $\Gamma$ at $x$ oriented in such way that $\det(\bt,\bn)=+1$,
and by $\kappa(x)$ the curvature
at $x$ which is defined by the relation
$$
\partial_{\bt(x)}\bt(x)=\kappa(x)\bn(x).
$$
For $p\in [1,\infty]$ and $\tau$ defined by $\frac 1 \tau := 1+\frac 1 p$,
we introduce the two quantities
\begin{eqnarray*}
S_p(f) &:= & \|\sqrt{|\det (d^2 f) |}\|_{L^\tau(\Omega\sm \Gamma)}=A_p(f_{|\Omega\sm \Gamma}),\\
E_p(f) &:= & \|\sqrt{|\kappa|} [f]\|_{L^\tau(\Gamma)},\\
\end{eqnarray*}
which respectively measure the ``smooth part'' and the ``edge part'' of $f$. We
also introduce the constant
\be
C_{p,\vp}:=\|\sqrt {|\Phi\Phi'|}\|_{L^\tau(\R)},\;\; \text{ where } \ \Phi(x) := \int_{y\in\R} \vp(x,y)dy.
\label{defchi}
\ee
Note that $f_\delta$ is only properly defined on the set
$$
\Omega^\delta:=\{z\in \Omega\; ; \; B(z,\delta)\subset \Omega\},
$$
and therefore, we define
$A_p(f_\delta)$ as the $L^\tau$ norm of $\sqrt{|\det(d^2f_\delta)|}$ on this set.

\begin{theorem}
\label{threg}
For all cartoon functions $f$, the quantity $A_p(f_\delta)$ behaves as follows:
\begin{itemize}
\item If $p<2$, then
$$
\lim_{\delta\to 0} A_p(f_\delta) = S_p(f).
$$
\item If $p=2$, then $\tau=\frac 2 3$ and 
\be
\lim_{\delta\to 0} A_2(f_\delta) = \(S_2(f)^\tau+ E_2(f)^\tau C_{2,\vp}^\tau\)^{\frac 1 \tau}.
\label{Acartoon}
\ee
\item If $p>2$, then $A_p(f_\delta) \to \infty$ according to
\be
\lim_{\delta\to 0}  \delta^{\frac 1 2 - \frac 1 p} A_p(f_\delta)= E_p(f) C_{p,\vp}.
\label{equivAphih}
\ee
\end{itemize}
\end{theorem}

\begin{remark}
This theorem reveals that as $\delta\to 0$, the 
contribution of the neighbourhood of $\Gamma$ to $A_p(f_\delta)$
is negligible when $p<2$ and 
dominant when $p>2$, which was already remarked in
the heuristic computation leading to {\rm \iref{heurist}}.
\end{remark}

\begin{remark}
It seems to be possible to eliminate the ``no cusps'' condition
in the definition of cartoon functions, while still retaining the validity of
this theorem. It also seems possible to take the more natural choice 
$\vp(z) = \frac 1 \pi e^{-\|z\|^2}$, which is not compactly supported.
However, both require higher technicality in the proof
which we avoid here.
\end{remark}

Before attacking the proof of Theorem \ref{threg},
we show below that the constant $C_{p,\vp}$ involved in 
the result for $p\geq 2$ is uniformly bounded by below
for a mild class of mollifiers.

\begin{proposition}
\label{propLowerCPhi}
Let $\vp$ be a radial and positive function supported on the unit ball
such that $\int \vp=1$ and that $\vp(x)$ decreases as $|x|$ increases. For any $p\geq 2$ we have 
$$
C_{p, \vp} \geq \frac 2 \pi \left(\frac 4 {\tau+2}\right)^{\frac 1 \tau}.
$$
and this lower bound is optimal. 
There is no such bound if $p<2$, but note that Theorem \ref{threg} does not involve $C_{p, \vp}$ for $p<2$.
\end{proposition}

\proof
Let $D$ be the unit disc of $\R^2$. We define a non smooth mollifier $\psi$ and a function $\Psi$ as follows
$$
\psi := \frac {\Chi_D} \pi \ \text{ and } \ \Psi(x) :=  \int_\R\psi(x,y) dy.
$$
One easily obtains that
$
\Psi(x) = \frac 2 \pi \sqrt{1-x^2} \Chi_{[-1,1]}(x)\text{ and } \Psi'(x) = \frac {-2x}{\pi \sqrt{1-x^2}} \Chi_{[-1,1]}(x), 
$
hence 
$$
\Psi(x) \Psi' (x)=  \frac {-4x} {\pi^2} \Chi_{[-1,1]}(x).
$$
For all $\delta>0$ we define $\psi_\delta := \delta^{-2} \psi(\delta^{-1} \cdot)$, and  $\Psi_\delta(x) :=  \int_\R\psi_\delta(x,y) dy$. Similarly we obtain 
$$
\Psi_\delta(x) \Psi'_\delta (x) = \frac {-4x} {\pi^2}  \Chi_{[-\delta,\delta]} \delta^{-4}.
$$
Hence 
\begin{eqnarray*}
C_{p,\psi_\delta} &=& \left\|\sqrt{\Psi_\delta(x) |\Psi'_\delta (x)|}\right\|_{L^\tau(\R)} \\
&=& \frac 2 {\pi \delta^2} \left(\int_{-\delta}^\delta |x|^{\frac\tau 2} dx\right)^{\frac 1 \tau} \\
&=& \frac 2 {\pi \delta^2} \left(\frac{2 \delta^{\frac \tau 2 +1}}{\frac \tau 2 +1}\right)^{\frac 1 \tau}\\
&=& \frac 2 \pi \left(\frac 4 {\tau+2}\right)^{\frac 1 \tau} \delta^{\frac 1 p -\frac 1 2}.
\end{eqnarray*}
Note that
\be
\label{ineqPsiDelta}
\text{ If } p\geq 2 \text{ and }\delta \in (0,1] \text{ then }  C_{p,\psi_\delta} \geq C_{p,\psi} =   \frac 2 \pi \left(\frac 4 {\tau+2}\right)^{\frac 1 \tau}.
\ee
The mollifier $\vp$ of interest is radially decreasing, has unit integral and is supported on the unit ball. It follows that there exists a Lebesgue measure $\mu$ on $(0,1]$, of mass $1$, such that 
$$
\vp = \int_0^1\psi_\delta \ d\mu(\delta).
$$
Hence
$
\Phi(x) := \int_\R \vp(x, y) dy = \int_0^1\Psi_\delta(x) \ d\mu(\delta)
$,
 for any $x\in \R$.
Since $s\mapsto s^\tau$ is concave on $\R_+$ when $0<\tau \leq 1$, we obtain
$$
\Phi(x)^\tau = \left(\int_0^1 \Psi_\delta(x) \ d\mu(\delta)\right)^\tau \geq \int_0^1 \Psi_\delta(x)^\tau \ d\mu(\delta)
$$
Similarly, since the sign of $\Psi'_\delta(x)$ is independent of $\delta$, 
$
|\Phi'(x)|^\tau = \left(\int_0^1 |\Psi_\delta'(x)| \ d\mu(\delta)\right)^\tau \geq \int_0^1 |\Psi_\delta'(x)|^\tau \ d\mu(\delta).
$
Applying the Cauchy-Schwartz inequality we obtain
\begin{eqnarray*}
\sqrt{\Phi(x) |\Phi'(x)|}^{\, \tau} &\geq& 
\sqrt{\left(\int_0^1 \Psi_\delta(x)^\tau \ d\mu(\delta)\right)\left(\int_0^1 |\Psi_\delta'(x)|^\tau \ d\mu(\delta)\right)} \\
&\geq& \int_0^1 \sqrt{\Psi_\delta(x) |\Psi_\delta'(x)|}^{\, \tau} \ d\mu(\delta)
\end{eqnarray*}
Eventually we obtain using the previous equation and \iref{ineqPsiDelta} that
$$
C_{p,\vp}^\tau = \int_\R \sqrt{\Phi |\Phi'|}^{\,\tau}
 \geq \int_0^1 \left(\int_\R \sqrt{\Psi_\delta |\Psi_\delta'|}^{\,\tau}\right) d\mu(\delta)\\
=  \int_0^1  C_{p,\psi_\delta}^\tau d\mu(\delta) \geq  C_{p,\psi}^\tau,
$$
which concludes the proof of this lemma.
\sq

The rest of this section is devoted to the proof of Theorem \ref{threg}. Since
it is rather involved, we split its presentation into several main steps.
\nl
\nl
{\bf Step 1: decomposition of $A_p(f_\delta)$.}
Using the notation $K(M):= \sqrt{|\det M|}$, we can write
\be
\label{tointegrate2d}
A_p(f_\delta)^\tau = \int_{\Omega^\delta} K(d^2 f_\delta)^\tau.
\ee
We decompose this quantity based on a partition of
$\Omega^\delta$ into three subsets 
$$
\Omega^\delta=\Omega_\delta \cup \Gamma_\delta \cup \cP_\delta.
$$
The first set $\Omega_\delta$ corresponds to the
{\it smooth part}:
$$
\Omega_\delta :=  \bigcup_{1\leq i\leq k} \Omega_{i,\ssdelta}, \;\;{\rm where}\;\; 
\Omega_{i,\ssdelta} := \{z\in \Omega_i\sep d(z,\Omega\sm\Omega_i)>\delta\}.
$$
Note that $\Omega_\delta$ is strictly contained in $\Omega^\delta$.
The second set corresponds to the {\it edge part}:  we first define
$$
 \Gamma_\delta^0 := \bigcup_{1\leq j\leq l} \Gamma_{j,\ssdelta}^0, \;\;{\rm where}\;\; 
\Gamma_{j,\ssdelta}^0 := \{z\in \Gamma_j\sep d(z,\Gamma\sm\Gamma_j)>2\delta\},
$$
and then set
$$
\Gamma_\delta := \bigcup_{1\leq j\leq l} \Gamma_{j,\ssdelta}\;\;{\rm where}\;\; 
\Gamma_{j,\ssdelta} :=  \{z\in \Omega \sep d(z,\Gamma)< \delta \text{ and } \pi_\Gamma(z)\in \Gamma_{j,\ssdelta}^0\} 
$$
where $\pi_\Gamma(z)$ denotes the point of $\Gamma$ which is the closest to $z$.
The third set corresponds the {\it corner part}:
$$
\cP_\delta :=  \Omega^\delta \sm (\Omega_\delta\cup\Gamma_\delta).
$$
The measures of the sets $\Gamma_\delta$ and $\cP_\delta$ tends to
$0$ as $\delta\to 0$, while $|\Omega_\delta|$ tends to $|\Omega|$. More
precisely, we have
$$
 |\Gamma_\delta|\leq C\delta\;\;{\rm and}\;\; |\cP_\delta|\leq C  \delta^2
$$
where the last estimate exploits the ``no cusps'' property of the cartoon function.
We analyze separately the contributions of these three sets to 
\iref{tointegrate2d}.
\nl
\nl
{\bf Step 2: Contribution of the smooth part $\Omega_\delta$.}
The contribution of $\Omega_\delta$ to the integral \iref{tointegrate2d} is easily measured. Indeed, let us define
$$
Q_\delta(z):= \left\{
\begin{array}{cc}
K(d^2 f_\delta (z))^\tau & \text{ if } z\in \Omega_\delta,\\
0 & \text{ otherwise.}
\end{array}
\right.
$$
Then we have pointwise convergence $Q_\delta(z)\to K(d^2 f(z))^\tau$ on $\Omega\sm \Gamma$.
Since the  $\delta$-neighbourhood of $\Omega_\delta$ is included in $\Omega\sm\Gamma$, we have 
$$
\|d^2 (f*\vp_\delta) \|_{L^\infty(\Omega_\delta)} = \|(d^2 f)*\vp_\delta \|_{L^\infty(\Omega_\delta)}\leq \|d^2 f\|_{L^\infty(\Omega\sm \Gamma)} \|\vp_\delta\|_{L^1} = \|d^2 f\|_{L^\infty(\Omega\sm \Gamma)} 
$$
Since 
$$
K(M) = \sqrt{|\det M|} \leq \|M\|
$$ 
we have
$
K(d^2 f_\delta)\leq \|d^2 f\|_{L^\infty(\Omega\sm \Gamma)} 
$ on $\Omega_\delta$,
and we conclude by dominated convergence that 
$$
\lim_{\delta \to 0}\int_{\Omega_\delta} K(d^2 f_\delta)^\tau = \lim_{\delta\to 0} \int_{\Omega \sm \Gamma} Q_\delta =  \int_{\Omega \sm \Gamma} K(d^2 f)^\tau.
$$
\nl
{\bf Step 3: Contribution of the corner part $\cP_\delta$.}
We only need a rough upper estimate of the contribution of $\cP_\delta$ 
to the integral \iref{tointegrate2d}.
We observe that 
$$
\|d^2 (f*\vp_\delta) \|_{L^\infty(\Omega)} = \|f*(d^2  \vp_\delta)\|_{L^\infty(\Omega)}\leq \|f\|_{L^\infty(\Omega)} \|d^2 \vp_\delta\|_{L^1(\R^2)} =   \frac M {\delta^2},
$$
where $M:=\|f\|_{L^\infty(\Omega)} \|d^2 \vp\|_{L^1(\R^2)}$.
It follows that 
$$
\int_{\cP_\delta} K(d^2 f_\delta)^\tau \leq |\cP_\delta| \left(\frac M {\delta^2}\right)^\tau \leq C  \delta^{2-2 \tau}.
$$
If $\tau<1$, this quantity tends to $0$ and is therefore negligible compared to the 
contribution of the smooth part. If $\tau=1$, which corresponds to $p=\infty$, our
further analysis shows that the contribution of the edge part tends
to $+\infty$, and therefore the contribution of the corner part is always negligible.
\nl
\nl
{\bf Step 4: Contribution of the edge part $\Gamma_\delta$.} {
This step is the main difficulty of the proof. We make a key use
of an asymptotic analysis of $f_\delta$ on $\Gamma_\delta$, which relates
its second derivatives to the jump $[f]$ and the curvature $\kappa$ as $\delta\to 0$. 
We first define for all $\delta>0$ the map
$$
\begin{array}{ccc}
U_\delta : \Gamma\sm\cP \times [-1,1] &\to& \Omega\\
(x,u) &\mapsto & x+\delta u \bn (x).
\end{array}
$$
We notice that according to our definitions,
for $\delta$ small enough, the map $U_\delta$ induces a diffeomorphism 
between $\Gamma_\delta^0\times [-1,1]$ and $\Gamma_\delta$, such
that $\pi_\Gamma(U_\delta(x,u))=x$ and $d(U_\delta(x,u),\Gamma)=|U_\delta(x,u)-x|=\delta |u|$.
We establish asymptotic estimates on the
second derivatives of $f_\delta$ which have the following form:
\begin{eqnarray}
\label{estimnn2d}
\left |\partial_{\bn,\bn} f_\delta(z) -  \frac 1 { \delta^2} [f](x) \Phi'(u)\right| &\leq & \frac C { \delta}\\ 
\label{estimnt2d}
|\partial_{\bn,\bt} f_\delta(z)| &\leq &  \frac C \delta\\
\label{estimtt2d}
\left |\partial_{\bt,\bt} f_\delta(z) +  \frac 1 \delta [f](x) \kappa(x)\Phi(u)\right| &\leq &\frac {\omega(\delta)} \delta
\end{eqnarray}
where $\lim_{\delta\to 0} \omega(\delta) = 0$ and with the notation
$z=U_\delta(x,u)$. The constant $C$ and the function $\omega$ depend only on $f$.
The proof of these estimates is given in the appendix. As an immediate consequence,
we obtain an asymptotic estimate of $K(d^2 f_\delta)=\sqrt{|\det(d^2 f_\delta)|}$ of the form
\be
\label{estimK2d}
\left|  \delta^{\frac 3 2} K(d^2 f_\delta(z)) - \sqrt{|\kappa(x)|} \, |[f](x)| \sqrt{|\Phi(u) \Phi'(u)|} \right|\leq \ti \omega(\delta),
\ee
where $\lim_{\delta\to 0} \ti \omega(\delta) = 0$, and the function $\ti \omega$ depends only on $f$. 
Using the notations
$$
g_\delta(z):=\delta^{\frac 3 2} K(d^2 f_\delta(z)),\;\;  \lambda(x):=\sqrt{|\kappa(x)|}\, |[f](x)|,\;\; \mu(u):=\sqrt{|\Phi(u) \Phi'(u)|},
$$
we thus have
\be
|g_\delta(z) - \lambda(x) \mu(u)|\leq \ti \omega(\delta), 
\label{glambmu}
\ee
for all $x\in \Gamma^0_\delta$, $u\in [-1,1]$ and $\delta>0$
sufficiently small, with $z=U_\delta(x,u)$. We claim that
for any continuous functions $(g_\delta,\lambda,\mu)$ satisfying \iref{glambmu},
we have for any $\tau>0$, 
\be
\lim_{\delta\to 0}  \delta^{-1} \int_{\Gamma_\delta} g_\delta^\tau =\int_{\Gamma} \lambda(x)^\tau dx  \int_{-1}^1 \mu(u)^\tau du,
\label{claim}
\ee
which is in our case equivalent to the estimate
\be
\lim_{\delta\to 0}  \delta^{ \frac 3 2 \tau-1} \int_{\Gamma_\delta} K(d^2 f_\delta)^\tau = 
 \int_\Gamma |[f]|^\tau |\kappa|^{\tau/2} \int_\R |\Phi \Phi'|^{\tau/2}.
\label{estimgammadelta}
\ee
In order to prove \iref{claim}, we may assume without loss of generality that $\tau =1$
up to replacing $(g_\delta,\lambda,\mu)$  by $(g_\delta^\tau,\lambda^\tau,\mu^\tau)$.
We first express the jacobian matrix of $U_\delta$
using the bases $B_1 = ((\bt(x),0),(0,1))$ and $B_2 =(\bt(x),\bn(x))$ 
for the tangent spaces of $\Gamma\times [-1,1]$ and $\Omega$.
This gives the expression
$$
[d U_\delta(x,u)]_{B_1,B_2} = \left(
\begin{array}{c|c}
1- \delta u\kappa(x)& 0\\
\hline
0 & \delta
\end{array}
\right)
$$
and therefore $|\det([dU_\delta(x,u)]_{B_1,B_2})| = \delta -\delta^2 u\kappa(x)$. Since
$B_1$ and $B_2$ are orthonormal bases, this quantity is the jabobian of $U_\delta$ at $(x,u)$,
and therefore
$$
\int_{\Gamma_\delta} g_\delta =\delta\int_{\Gamma_\delta^0 \times [-1,1]} g_\delta(x+\delta u \bn (x)) 
(1 -\delta u\kappa(x)) dx\, du
$$
Combining with \iref{glambmu}, and using dominated
convergence we obtain \iref{claim}.
\nl
\nl
{\bf Step 5: summation of the different contributions.}
Summing up the contributions of $\Omega_\delta$, $\cP_\delta$ and $\Gamma_\delta$,
we reach the estimate
\begin{eqnarray*}
\int_\Omega K(d^2 f_\delta)^\tau &=&  \int_{\Omega_\delta} K(d^2 f_\delta)^\tau+ \int_{\cP_\delta} K(d^2 f_\delta)^\tau+ \int_{\Gamma_\delta} K(d^2 f_\delta)^\tau\\
&=& \left(\int_{\Omega\sm\Gamma} K(d^2 f)^\tau +\e_1(\delta) 
\right) + B(\delta) \delta^{2-2\tau} \\
& &+  \delta^{1-\frac 3 2 \tau} \left( \int_\Gamma |[f]|^\tau |\kappa|^{\tau/2}\int_\R |\Phi \Phi'|^{\tau/2} +\e_2(\delta)\right),\\
&=& (S_p(f)^\tau +\e_1(\delta)) + B(\delta) \delta^{2-2\tau}+  \delta^{1-\frac 3 2 \tau} (E_p(f)^\tau C_{p,\phi}^\tau+\e_2(\delta)),\\
\end{eqnarray*}
where $\lim_{\delta\to 0}\e_1(\delta)=\lim_{\delta\to 0}\e_2(\delta)=0$ and $B(\delta)$ is uniformly bounded.
This concludes the proof of Theorem \ref{threg}.

\section{Relation with other works}
\label{secCartoon4}
Theorem \ref{threg} allows us to extend
the definition of $A_2(f)$ when $f$ is a cartoon
function, according to
\be
A_2(f):= \(S_2(f)^{2/3}+ E_2(f)^{2/3} C_{2,\vp}^{2/3} \)^{3/2}.
\label{addita2}
\ee
We first compare this additive form
with the total variation $\TV(f)$. If $f$ is a
cartoon function,
its total variation has the additive form
\be
\TV(f):=\|\nabla f\|_{L^1(\Omega\sm\Gamma)}+ \|[f]\|_{L^1(\Gamma)},
\label{additv}
\ee
Both \iref{addita2} and \iref{additv}
include a ``smooth term'' and an ``edge term''.
It is interesting to compare the edge term
of $A_2(f)$, which is given by
$$
E_2(f)=\|\sqrt {|\kappa|}[f]\|_{L^{2/3}(\Gamma)},
$$
up to the multiplicative constant $C_{2,\vp}$, with the
one of $\TV(f)$ which is simply the integral of the jump
$$
J(f):=\|[f]\|_{L^{1}(\Gamma)},
$$
Both terms are $1$-homogeneous with the value of the jump of the
function $f$. In particular, if the value of this jump is $1$
(for example when $f$ is the characteristic function of a set of boundary $\Gamma$), we have
\be
E_2(f)=\(\int_\Gamma |\kappa|^{1/3}\)^{3/2},
\label{powercurv}
\ee
while $J(f)$ coincides with the length of $\Gamma$. 
In summary, $A_2(f)$ takes into account the {\it smoothness} of
edges, through their curvature $\kappa$, while 
$\TV(f)$ only takes into account their {\it length}.

Let us now investigate more closely the measure of
smoothness of edges which is incorporated in $A_2(f)$. According to 
\iref{powercurv}, this smoothness is meant in the sense
that the arc length parametrizations of the curves that
constitute $\Gamma$ admit  
second order derivatives in $L^{\frac 1 3}$. In the following,
we show that this particular measure of smoothness is 
naturally related to some known results in two different areas:
adaptive approximation of curves and affine-invariant image
processing.

We first revisit the derivation of the heuristic estimate
\iref{heurist} for the error between a cartoon function
and its linear interpolation on an optimally adapted 
triangulation. In this computation, the contribution of the
``edgy triangles'' was estimated by the area of the layer $\Omega_N^e$ according to
$$
\|f-\interp_{\cT_N} f\|_{L^p(\Omega^{e}_N)} \leq \|f\|_{L^\infty} |\Omega_N^e|^{1/p}.
$$
Then we invoke the fact that $\Gamma$ is a finite union of $C^2$ curves $\Gamma_j$
in order to build a layer of global area $|\Omega_N^e| \leq CN^{-2}$, 
which results in the case $p=2$ into a contribution to the $L^p$ error
of the order $\cO(N^{-1})$. The area of the layer $\Omega_N^e$ 
is indeed of the same order as the area between the edge $\Gamma$
and its approximation by a polygonal line with $\cO(N)$ segments.

Each of the curves $\Gamma_j$
can be identified to the graph of a $C^2$ function
in a suitable orthogonal coordinate system.
If $\gamma$ is one of these functions, the area
between $\Gamma$
and its polygonal approximation can thus be 
locally measured by the $L^1$ error
between the one-dimensional function
$\gamma$ and a piecewise linear approximation
of this function.
Since $\gamma$ is $C^2$, it is obvious
that it can be approximated by a piecewise linear function on $\cO(N)$ intervals
with accuracy $\cO(N^{-2})$ in the $L^\infty$ norm and therefore
in the $L^1$ norm. However, we may ask whether such a rate
could be achieved under weaker conditions on the smoothness of $\gamma$.
The answer to this question is a chapter of {\it nonlinear} approximation theory
which identifies the exact conditions for a function $\gamma$ to be approximated
at a certain rate by piecewise polynomial functions on adaptive one-dimensional
partitions. We refer to \cite{De} for a detailed
treatment and only state the result which is of interest to us. We say that a function
$\gamma$ defined on a bounded interval $I$ belongs to the approximation
space ${\cal A}^s(L^p)$ if and only if there exists a sequence $(p_N)_{N>0}$ of functions where
each $p_N$ is piecewise affine on a partition of $I$ by $N$ intervals
such that 
$$
\|\gamma-p_N\|_{L^p} \leq CN^{-s}.
$$
For $0<s\leq 2$, it is known that 
$\gamma\in {\cal A}^s(L^p)$ provided that $\gamma\in B^{s}_{\tau,\tau}(I)$ with $\frac 1 \tau:=\frac 1 p+s$,
where $B^s_{\tau,\tau}(I)$ is the standard Besov space
that roughly describes those functions having $s$ derivatives in $L^{\tau}$.
In the case $s=2$ and $p=1$ which is of interest to us, we find $\tau=\frac 1 3$ and therefore $\gamma$ should
belong to the Besov space $B^2_{\frac 1 3,\frac 1 3}(I)$. Note that in our definition 
of cartoon functions, we assume much more than $B^2_{\frac 1 3,\frac 1 3}$ smoothness on $\gamma$,
and it is not clear to us if Theorem \iref{threg} can be derived under this minimal smoothness 
assumption. However it is striking to see that the quantity $E_2(f)$
that is revealed by Theorem \iref{threg} precisely measures the second derivative of the
arc-length parametrization of $\Gamma$ in the $L^{\frac 1 3}$ norm, up to
the multiplicative weight $|[f]|^{2/3}$. Let us also mention
that Besov spaces have been used in \cite{DPW} in order to describe
the smoothness of functions through the regularity of their level sets.
Note that edges and level sets are two distinct concepts, which
coincide in the case of piecewise constant cartoon functions.

The quantity $|\kappa|^{1/3}$ is also encountered in
mathematical image processing, for the design of
simple smoothing semi-groups that respect {\it affine invariance}
with respect to the image. Since these semi-groups should also
have the property of {\it contrast invariance}, they can be defined
through curve evolution operators acting on the level sets
of the image. The simplest curve evolution operator
that respects affine invariance is given by the equation
$$
\frac {d \Gamma}{d t}=-|\kappa|^{1/3}{\rm n},
$$
where ${\rm n}$ is the outer normal, see e.g. \cite{Ca}.
Here the value $1/3$ of the exponent plays a critical
role. The fact that we also find it in $E_2(f)$ suggests
that some affine invariance property also holds for this quantity
as well as for $A_2(f)$. 
We first notice that if $f$ is a compactly supported
$C^2$ function of two variables and $T$ is a 
bijective affine transformation, then with
$\ti f$ such that
$$
f=\ti f\circ T,
$$
we have
the property
$$
d^2 f(z)=L^\trans d^2 \ti f(Tz) L,
$$
where $L$ is the linear part of $T$ and $L^\trans$ its transpose, so that
$$
\sqrt {|\det (d^2f(z))|}=|\det L|\sqrt{|\det (d^2 \ti f(Tz))|}.
$$
By change of variable, we thus find that
\be
A_p(\ti f)=|\det L|^{1/\tau-1}A_p(f)=|\det L|^{1/p}A_p(f).
\label{invap}
\ee
A similar invariance property can be  derived on the interpolation error
$\sigma_N(f)_p=\|f-\interp_{\cT_N}f\|_{L^p}$ where $\cT_N$ is a triangulation which is optimally adapted
to $f$  in the
sense of minimizing the linear interpolation error in the $L^p$ norm
among all triangulations of cardinality $N$. We indeed remark that
an optimal triangulation for $\ti f$ is then given by applying $T$ to
all elements of $\cT_N$. For such a triangulation 
$\ti \cT_N:=T(\cT_N)$, one has the commutation formula
$$
\interp_{\cT_N}f=(\interp_{\ti \cT_N} \ti f)\circ T,
$$
and therefore we obtain by a change of variable that
\be
\sigma_N(\ti f)_p=\|\ti f-\interp_{\ti\cT_N}\ti f\|_{L^p}=|\det L|^{1/p}\|f-\interp_{\cT_N}f\|_{L^p}=|\det L|^{1/p}\sigma_N(f)_p.
\label{invsigmap}
\ee
Let us finally show that if $f$ is a cartoon function, then $E_2(f)$ satisfies a similar invariance property
corresponding to $p=2$, namely
\be
E_2(\ti f)=|\det L|^{1/2} E_2(f).
\label{inve2}
\ee
Note that this cannot be derived by arguing that $A_2(f)$ satisfies this
invariance property when $f$ and $\ti f$ are smooth, since we lose
the affine invariance property as we introduce the convolution by $\vp_\delta$: 
we do not have
$$
(f\circ T) * \vp_\delta = (f* \vp_\delta)\circ T.
$$
unless $T$ is a rotation or a translation. 

Let $\Gamma_j$ be one of the $C^2$ pieces of $\Gamma$
and  $\gamma_j: [0,B_i]\to \Omega$ a regular parametrisation of 
$\Gamma_j$. The curvature 
of $\Gamma$ on $\Gamma_j$ at the point $\gamma_j(t)$ is therefore given by
\be
\label{exprKappa}
\kappa(\gamma_j(t))= \frac{\det(\gamma_j'(t),\gamma_j'' (t))}{\|\gamma'_j(t)\|^{3}}
\ee
Since $f=\ti f\circ T$, the discontinuity curves of $\ti f$ are the images of those of $f$ by $T$:
$$
\ti \Gamma_j=T(\Gamma_j).
$$
The curvature of $\ti \Gamma_j$ at the point 
$T(\gamma_j(t))$ is therefore given by 
$$
\ti \kappa(T(\gamma_j(t))) = \frac{\det(L\gamma'_j(t), \ L\gamma''_j(t))}{\|L\gamma'_j(t)\|^3} = \det(L) \frac{\det (\gamma'_j(t),\gamma''_j(t))}{\|L\gamma'_j(t)\|^3} .
$$
This leads us to the relation : 
\be
\label{relKappa}
|\det(L)|^{1/3} |\kappa(\gamma_j(t))|^{1/3} \ \|\gamma'_j(t)\| =  |\ti\kappa(T(\gamma_j(t)))|^{1/3}
\|L\gamma'_j(t)\| ,
\ee
and therefore
\begin{eqnarray*}
\int_{\ti \Gamma_j} |[\ti f]|^{2/3} |\ti \kappa|^{1/3} &=& \int_0^{B_j}  |[\ti f](T(\gamma_j(t)))|^{2/3}
|\ti\kappa(T(\gamma_j(t)))|^{1/3} \ \|L\gamma'_j(t)\| dt,\\
&=& |\det L|^{1/3} \int_0^{B_j}  |[f](\gamma_j(t))|^{2/3} |\kappa(\gamma_j(t))|^{1/3} \  \|\gamma'_j(t)\|dt,\\
&=& |\det L|^{1/3} \int_{ \Gamma_i} |[ f]|^{2/3} |\kappa|^{1/3}.
\end{eqnarray*}
Summing over all $j=1,\cdots,l$ and elevating to the $3/2$ power
we obtain \iref{inve2}.

\section{Numerical tests}
\label{secCartoon5}

We first validate our previous results
by  numerical tests applied to a simple cartoon image:
the Logan-Shepp phantom. We use a $256\times 256$ pixel
version of this image, with a slight modification which is motivated
further. This image is iteratively smoothed
by the numerical scheme 
\be
\label{numSmoothing}
u_{i,j}^{n+1} = \frac {u_{i,j}^n} 2 + \frac {u_{i+1,j}^n+u_{i-1,j}^n+u_{i,j+1}^n+u_{i,j-1}^n} 8.
\ee
This scheme is an explicit discretization of the heat equation.
Formally, as $n$ grows,
$u^n$ is a discretization of 
\be
\label{Gaussian}
u * \vp_{\lambda \sqrt n} \text{ with } \vp_\delta (z) := \frac 1 { \pi\delta^2 } e^{-\frac {\|z\|^2}{\delta^2}},
\ee 
where $u$ stands for the continuous image. The determinant of the hessian 
is discretised by the following $9$-points formula
\be
\begin{array}{lll}
d_{i,j}^n &:=& (u_{i,j+1}^n-2u_{i,j}^n+u_{i,j-1}^n) (u_{i+1,j}^n-2u_{i,j}^n+u_{i-1,j}^n) \vspace{2mm}\\
& &\displaystyle- \frac {(u_{i+1,j+1}^n+u_{i-1,j-1}^n-u_{i+1,j-1}^n-u_{i-1,j+1}^n)^2}{16}
\end{array}
\label{ninepoint}
\ee
For each value of $n$, we then compute the $\ell^\tau$ norm of the array
$\left(\sqrt{|d^n_{i,j}|}\right)$ for $\tau \in [\frac 1 2,1]$, which corresponds to
$p\in [1,\infty]$ with $\frac 1 \tau := 1 + \frac 1 p$. This norm is thus a discretization of the quantity
$$
\left \|\sqrt{|\det(d^2(u*\vp_{\lambda \sqrt n}))|}\right \|_{L^\tau}.
$$
For each value of $n$ we obtain a function $\tau\in  [\frac 1 2,1]\to D_n(\tau) \in\R_+$.

As $n$ grows, three consecutive but potentially overlapping phases appear 
in the behaviour of the functions $D_n$, which are illustrated on Figure \ref{fig2Cartoon}.

\begin{enumerate}

\item For small $n$, the $9$-points discretisation is not a good approximation of the determinant of the hessian
due to the fact that the pixel discretization is too coarse compared to the 
smoothing width. During this phase, the functions $D_n$ decay rapidly for all values of $\tau$. 

\item For some range of $n$, the edges have been smoothed by the action of \iref{numSmoothing}, 
but the parameter $\lambda \sqrt n$ in \iref{Gaussian} remains rather small. Our previous analysis applies and we observe that $D_n(2/3)$ is (approximately) constant while $D_n(\tau)$ increases for 
$\tau < 2/3$ and decreases for $\tau > 2/3$. 

\item For large $n$, the details of the picture fade and begin to disappear. The picture begins to resemble a constant picture. Therefore the functions $D_n$ decay for all values of $\tau$, and eventually tend to $0$. 

\end{enumerate}

\begin{figure}
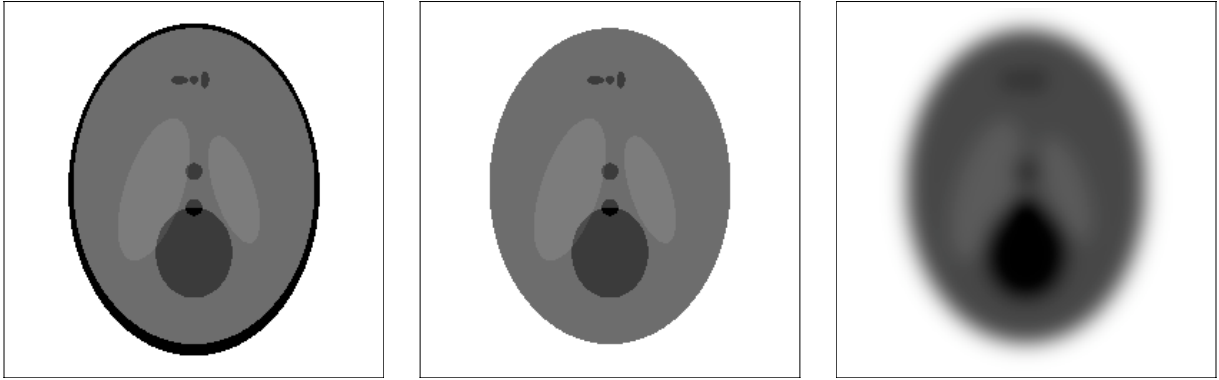

	\centering
		\includegraphics[width=5cm,height=5cm]{\pathPic/PaperCartoon/PicShepp.pdf}
		\hspace{0.2cm}
		\includegraphics[width=5cm,height=5cm]{\pathPic/PaperCartoon/PicSheppSimpler.pdf}
		\hspace{0.2cm}
		\includegraphics[width=5cm,height=5cm]{\pathPic/PaperCartoon/PicSmoothedShepp.pdf}
	\caption{\label{fig1Cartoon}The Logan-Shepp image (left), modified image (center), smoothed image (right).}
\end{figure}

\begin{figure}
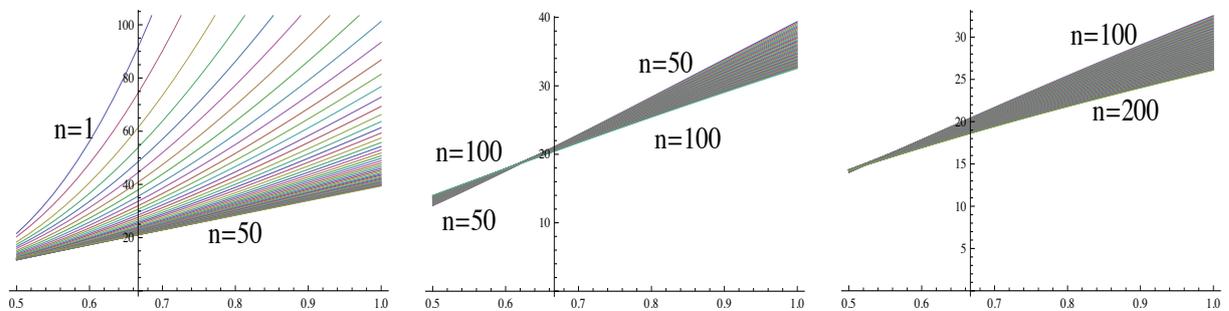

	\centering
		\includegraphics[width=5cm,height=4cm]{\pathPic/PaperCartoon/PicBeginShepp.pdf}
		\hspace{0.2cm}
		\includegraphics[width=5cm,height=4cm]{\pathPic/PaperCartoon/PicMidShepp.pdf}
		\hspace{0.2cm}
		\includegraphics[width=5cm,height=4cm]{\pathPic/PaperCartoon/PicEndShepp.pdf}
	\caption{\label{fig2Cartoon}The curves $D_n(\tau)$ for $n\leq 50$ (left),
	$50\leq n\leq 100$ (center) and 
	$100\leq n\leq 200$ (right).}
\end{figure}

In our numerical experiments, we used the well known Shepp-Logan Phantom  with a slight modification
as shown on Figure \ref{fig1Cartoon}: we have removed the thin layer around the head, which represents the skull, 
because it disappears too quickly by the smoothing procedure
and causes phases 1 and 3 to overlap, masking phase 2.
For more complicated functions $f$, such as most photographic pictures, phases $1$ and $3$ 
also tend to overlap for similar reasons. 
Indeed, these pictures often have details at the pixel scale, including electronic noise due to the captor. 
Since these details disappear early phase 3 begins immediately, therefore phase 2 cannot be observed.\\

Our next discussion aims at comparing the quantities $A_2(f)$ with the
total variation $\TV(f)$ on numerical images. In particular, we want to compare how these quantities 
measure the geometric complexity of images. Here, we consider non-discretized images
$$
z=(x,y)\in [0,1]^2 \mapsto f(z),
$$
and we compare the numerical behaviour of  $A_2(f)$ and $\TV(f)$ for some relevant cases.
Let us recall that for any $g\in C^2(\R_+)$, if $f$ is the
radial function
$$
f(z)=g(|z|),
$$
one has 
$$
\det \left(d^2f (z)\right) = \frac 1 {|z|} g'(|z|) g''(|z|).
$$
As a result we find that the two functionals $\TV$ and $A_2^{2/3}$ behave similar 
on the oscillatory images $f_\omega(z):= \cos(\omega |z|)$ illustrated on Figure \ref{figTVDa} (left),
as $\omega \to +\infty$:
\be
\label{eqComplexaCos}
\TV( f_\omega ) \simeq \omega \stext{ and } A_2(f_\omega) \simeq \omega^{3/2}.
\ee
We next consider the cartoon images defined by 
$g_\omega(z) = \frac {\lfloor \omega |z|\rfloor} \omega$, which
have the form of a circular staircase with approximately $\omega$ steps of height $1/\omega$,
as illustrated on Figure \ref{figTVDa} (center). 
In this case, we may give a meaning to $A_2(g_\omega)$ based on \iref{addita2}, 
which since $g_\omega$ is piecewise constant gives
$$
A_2(g_\omega)^{2/3}=C\int_\Gamma  |[g_\omega]|^{\frac 2 3} |\kappa|^{\frac 1 3},
$$
where $\Gamma$ are the circular curves of discontinuities. From this we find that
as $\omega\to +\infty$,
\be
\label{eqComplexaStair}
\TV(g_\omega) \simeq 1 \stext{ and } A_2(g_\omega) \simeq \sqrt \omega.
\ee
This shows that in contrast to $\TV$, the functional $A_2$ penalizes images which 
have the appearance of a staircase.

Finally we consider the cartoons functions $h_\omega(z)= \chi_{S_\omega}$ 
where the set $S_\omega$ is defined by
$$
S_\omega := \left\{z=(r\cos \theta, r \sin \theta) \sep 0 \leq \theta \leq \frac \cPi 2, \text{ and } 0 \leq r \leq 1 + \frac 1 \omega \cos(\omega \theta)\right\},
$$
see Figure \ref{figTVDa} (right). These functions have therefore only one step
but its geometry becomes more and more oscillatory as $\omega$ grows. More precisely 
its length remains finite, but its curvature behaves like $\omega$. In turn, we obtain
that as $\omega\to +\infty$,
\be
\label{eqComplexaFlower}
\TV(h_\omega) \simeq 1 \stext{ and } A_2(h_\omega) \simeq \sqrt \omega.
\ee
This shows that in contrast to $\TV$, the functional $A_2$ penalizes 
the fact that the curve of discontinuity $\partial S_\omega$ strongly oscillates.

These examples suggest that the functional $A_2$ 
gives a more relevant account of the complexity of cartoon images 
than the total variation $\TV$. In particular staircasing effects and 
oscillatory curves of discontinuity are penalized.\\


\begin{figure}
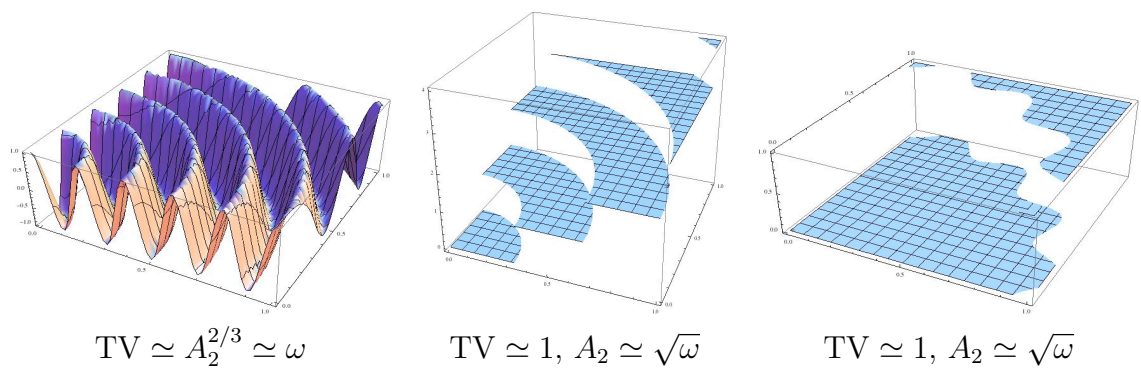

\centering
\begin{tabular}[c]{ccc}
\includegraphics[height=3.5cm,width=5cm]{\pathPic/Bayesien/SinaNoText.pdf}
&\includegraphics[height=4cm,width=4cm]{\pathPic/Bayesien/StairaNoText.pdf}
&\includegraphics[height=3.5cm,width=5cm]{\pathPic/Bayesien/FloweraNoText.pdf}\\
 $\TV \simeq A_2^{2/3} \simeq \omega$
& $\TV \simeq 1$, $A_2 \simeq \sqrt \omega$
& $\TV \simeq 1$, $A_2 \simeq \sqrt \omega$
\end{tabular}
\caption{Some representatives of families $f_\omega$, $g_\omega$ and $h_\omega$ of images, and behavior of the functionals $\TV$ and $A_2$.\label{figTVDa}}
\end{figure}


\section{Applications to image restoration} 
\label{secCartoon5.5}

The previous discussion suggests that we could use $A_2$ as an alternative
to $\TV$ as a prior in image restoration. In particular, we may try to replace $|g|_{BV}$
by $A_2(g)$ in \iref{minBV}, and therefore consider the minimization problem
$$
\min_{g} \{A_2(g)\; ; \; \|Tg-h\|_{L^2}\leq \e\},
$$
or its formulation using a Lagrange multiplier 
\be
\label{eqMinLagrange}
\min_g \|Tg-h\|_{L^2}^2 + t A_2(g)^{2/3}.
\ee
Our first observation is that, in contrast to \iref{minBV}, these problems are non-convex.

Moreover, it is easily seen that even in the very simple case of image denoising corresponding to 
$T=Id$, the above problems are ill-posed, in the following sense.

\begin{prop}
\label{propNotWellPosed}
For any $f\in L^2([0,1]^2)$, there exists a sequence $(f_n)_{n\geq 0}$ of $C^\infty$ functions
such that 
$$
\|f-f_n\|_{L^2}\to 0\;\;{\rm and}\;\; A_2(f_n) \to 0.
$$
Consequently the infimums of $A_2(g)$ over those $g\in C^\infty$ such that $\|g-h\|_{L^2}\leq \e$
and of $\|g-h\|_{L^2}^2 + t A_2(g)^{2/3}$ over all $g\in C^\infty$ are both equal to $0$, and are
not attained in general.
\end{prop}

\proof
If $\cT$ is a triangulation of the domain $[0,1]^2$, and if $g$ is piecewise affine on $\cT$,
we may consider its regularized version $g_\delta:=g*\vp_\delta$. From Theorem \ref{threg},
we find that $A_2(g_\delta)\to 0$ as $\delta \to 0$. The proof follows by observing that 
piecewise affine functions on triangulations are dense in $L^2$.
\sq

Note that the above result is formulated in the setting of non-discretized
images, and does not exactly hold for numerical images if we define $A_2$ 
using the $9$-points formula \iref{ninepoint} for the discretization of the 
determinant of the hessian. However, we expect that the ill-posedness
in the continuous setting is reflected by a bad behaviour of the solution
of the discrete optimization problem.

In summary, a straightforward generalization of the minimization
approach \iref{minBV} with $A_2$ in place of $\TV$ is doomed to fail.
In order to circumvent this difficulty, we propose a different strategy 
based on a restoration algorithm introduced in \cite{LMo}
and \cite{LMo2}. This algorithm is based on a bayesian framework that we briefely recall.

\subsection{A bayesian approach to image denoising} 

For simplicity, we focus from this point on the discrete setting, in which images are discretized on a $N \times N$ grid where $N\geq 1$.
We denote by 
$$
\cF := \R^{N \times N}.
$$
the collection of all discrete images. Bayesian restauration 
algorithms are based on a prior probability distribution on $\cF$ describing
how certain images are more ``plausible'' than others. One way to
build such priors is through a functional 
$$
J : \cF \to \R_+
$$
which reflects the plausibility of an image: typically
the value of $J$ is large for complex images and zero only for few very simple images. 
The functional $J$ is typically obtained as the discretization of a measure of smoothness, such as the total variation $\TV$ or the quantity $A_2^{2/3}$ (discretized with the nine points formula \iref{ninepoint}) which is shown in \S\ref{secCartoon5} to give a good account of the complexity of cartoon images. 

Assuming that the quantity 
$$
Z_J := \int_\cF \exp(-J(F)) dF
$$
is finite, our bayesian prior for \emph{plausible images} is the probability distribution
\be
\label{defProbaF}
\frac{e^{-J(F)}}{Z_J} dF
\ee
where $dF$ stands for the Lebesgue measure on $\cF$.

We use the simplistic, although popular, model of additive gaussian noise in which random drafts from $\cF$ with respect to the probability 
\be
\label{defProbaOmega}
\frac{e^{-\alpha\|\Omega\|^2}}{Z_\alpha} d\Omega,
\ee
where $\alpha>0$ is a parameter and $Z_\alpha$ is a normalizing factor, represent \emph{typical noise}.

We regard a \emph{corrupted image} $G$ as a random variable defined as the sum 
$$
G = F+\Omega
$$
of a random variable $F$ of distribution \iref{defProbaF}, the \emph{original image}, and a random variable $\Omega$ of distribution \iref{defProbaOmega}, the \emph{corruption} by additive gaussian noise.

By construction the conditional probability density $p(G|F)$ of the corrupted image $G$ with respect to the original image $F$ is 
$$
p(G|F) = \frac {e^{-\alpha\|\Omega\|^2}}{Z_\alpha} = \frac{e^{-\alpha\|G-F\|^2}}{Z_\alpha}.
$$
The conditional probability density $p(F|G)$ of the original image $F$ with respect to the corrupted one $G$ is obtained by the bayesian rule of conditional probabilities
$$
p(F|G)p(G) = p(G|F) p(F).
$$
Replacing $p(G|F)$ and $p(F)$ with their explicit expressions
we obtain 
\begin{eqnarray*}
p(F|G) &=&  \frac 1 {p(G)} p(G|F) \, p(F)\\
&=& \frac 1 {p(G)} \left(\frac 1 {Z_\alpha} e^{-\alpha \|G-F\|^2_2}\right) \left(\frac 1 {Z_J} e^{-J(F)} \right)\\
&=& \frac 1 {p(G) Z_\alpha Z_J} e^{-\alpha \|G-F\|^2_2 - J(F)}.
\end{eqnarray*}
We define
$$
\sigma(F,G) = \exp( -\alpha \|G-F\|^2_2 - J(F)).
$$ 
For any fixed $G$, 
the explicit function $F\in \cF \mapsto \sigma(F,G)\in \R_+^*$ 
is therefore proportional to the conditional probability density $F\mapsto p(F|G)$.

In order to recover the original image $F$ from the corrupted one $G$ a first approach, called the maximum a posteriori (MAP), consists in maximizing the conditional probability density 
\begin{eqnarray*}
F_J^* &:=& \underset {F\in \cF} \argmax \ p(F|G)\\
&=& \underset {F\in \cF}\argmax \ \sigma(F,G)\\
&=& \underset{F\in \cF}\argmin\ \alpha \|F-G\|^2 + J(F)
\end{eqnarray*}
We thus recover the optimization procedure \iref{eqMinLagrange} in the case $T=\Id$ of image denoising, and of the functionnal $J = A_2^{2/3}$. As observed in Proposition \ref{propNotWellPosed} this approach is doomed to fail for the functional $A_2^{2/3}$. For the total variation functional $J = \TV$, the MAP approach gives 
good results, but is also known to produce visual artifacts:  the restored image is exactly constant on large regions, delimited by sharp discontinuities which do not correspond to a feature of the original image.
The heuristical reason of this problem is that the maximum of a probability density is generally not a good representative of a random draft for this probability, as discussed in \cite{LMo}.

A second approach, called the minimum mean square error (MMSE), consists in finding $F$ which minimizes the empirical quadratic risk $\|F-F'\|^2$ with respect to a random image $F'$ distributed according to the conditional probability density $p(F'|G)$:
\begin{eqnarray*}
F_J &:=& \underset {F\in \cF} \argmin \int_\cF \|F-F'\|^2 \, p(F'|G) \, dF',\\
&=&\underset {F\in \cF} \argmin \int_\cF \|F-F'\|^2 \, \sigma (F',G) \, dF',\\
&=& \frac 1 {Z} \int_\cF F\, \sigma(F,G) \, dF. 
\end{eqnarray*}
where 
$
Z = \int_\cF \sigma(F, G) \, dF
$ 
 is a normalizing factor. This method is studied in depth in \cite{LMo} in the case $J = \TV$ where it
 is shown that it has the advantage of suppressing the visual artifacts observed in the MAP approach, and we describe in the next section its adaptation to the case $J = A_2^{2/3}$.

\subsection{The restauration algorithm}
\label{secAlgoBayes}
%
%

The computation of the estimator $F_J$ is a numerical challenge, since it involves an integration on the space $\cF = \R^{N \times N}$ which has dimension at least $512 \times 512 \simeq 250000$ in realistic cases.
We use an algorithm which was first proposed in the thesis \cite{LMo}, and which is a variation on the Metropolis Hastings algorithm. We shall only sketch its main principles and we refer to the original thesis for more details. 

The central point is to define a Markov chain $(F_k)_{k\geq 0}$ of images in $\cF$ which is recurrent with respect to the probability measure $\sigma(F,G) dF/Z$, where $G$ is the denoised image and $Z$ is a normalizing factor.
Then, almost surely,
\be
\label{eqConAvg}
\lim_{K\to \infty} \frac 1 K \sum_{0\leq k \leq K-1} F_k = \frac 1 Z \int_\cF F \, \sigma(F,G) \, dF.
\ee
Two successive images $F_k$ and $F_{k+1}$ generated by the algorithm proposed in \cite{LMo} only differ by the value of a single pixel. The cost of a step of the algorithm, generating a new image $F_k$, comes mainly from the computation of the ratio
\be
\label{eqRatioSigma}
\frac {\sigma(U,G)}{\sigma(V,G)} = \exp\( \|V-G\|^2 - \|U-G\|^2 +J(V) - J(U)\),
\ee
where $U$ and $V$ are two elements of $\cF$ which differ at a single position $(i,j)$. 

This algorithm applies to any continuous functional $J \in C^0(\cF, \R_+)$, hence in particular to {\it non-convex functionals} such as $J=A_2^{2/3}$.
For numerical applications, in order to have a good speed of convergence, the ratio \iref{eqRatioSigma} needs to be extremely cheap to compute in terms of computer time. 
The discretisation of the total variation $\TV$ or of the functional $A_2^{2/3}$, are given for an image $F$ by respectively
$$
\frac 1 N \sum_{i,j} |F_{i,j}-F_{i+1,j}|+|F_{i,j} - F_{i,j+1}|\stext{ and } N^{-2/3} \sum_{i,j} |d_{i,j}|^\frac 1 3
$$
where $d_{i,j}$ is given by the formula \iref{ninepoint}.
The ratio \iref{eqRatioSigma} is cheap to compute numerically, as required by the algorithm, since it has an explicit algebraic expression which only involves the 
the values of the images $U,V\in \cF$ at the positions $(i+k, j+l)$ where $\max\{|k|, |l|\} \leq 2$, where $(i,j)$ denotes the single position at which these images differ.

From a theoretical point of view, the averaging procedure \iref{eqConAvg} converges at the speed $\cO(K^{-\frac 1 2})$. 
Unfortunately the \emph{constant} in front of the convergence rate plays a key role in practice, and we did not manage to ``reach the convergence'' for realistic $512 \times 512$ images with the functional $J = A_2^{2/3}$. This issue might be solved by a parallelization of the algorithm, or a modification of the algorithm in which a group of neighboring pixels is modified at each step instead of a single pixel. For the time being, we thus only present one-dimensional results.

\subsection{Numerical illustration in the 1D case}

%
%

We present numerical results in one dimension only, using some counterparts of the total variation $\TV$ and of the functional $A_2^{2/ 3}$ in that context.
We define for each function $f\in C^2([0,1])$ the quantity
\be
\label{defDf}
D(f) := \int_{[0,1]} \sqrt{|f''|}.
\ee
The next proposition shows that $D(f)$ can be given a meaning when $f$ has localized discontinuites, similarly to the functional $A_2^{2/3}$ for cartoon images.
\begin{prop}
Let $f : [0,1] \to \R$ be piecewise $C^2$, with a finite set $E$ of discontinuity points in $]0,1[$.
Let $\vp$ be a $C^2$ mollifier supported in $[-1,1]$, satisfying $\int_\R \vp = 1$ and $\vp(x) = \vp(-x)$ for all $x\in \R$. For all $\delta>0$ let $\vp_\delta := \frac 1 \delta \vp(\frac \cdot \delta)$ 
 and let $f_\delta := f* \vp_\delta$.
Then $f_\delta \in C^2([\delta, 1-\delta])$ and satisfies
\be
\label{defDDisc}
\lim_{\delta\to 0} \int_{[\delta, 1-\delta]} \sqrt{|f_\delta ''|} = \int_{]0,1[\sm E} \sqrt{|f''|} + C(\vp) \sum_{e\in E} \sqrt{|[f](e)|}
\ee
where $C(\vp) := \int_\R \sqrt{|\vp'|}$. 
\end{prop}

\proof
We denote by $0 <e_1< e_2 < \cdots< e_{n-1}< 1$ the points of $E$
and we define $e_0 = 0$ and $e_n = 1$.  
The function $f$ can written as the sum 
$$
f=J+A+S,
$$
where $J$ (the Jump part) is piecewise constant, and $A$ is continuous and piecewise Affine, 
with respect to the partition $\left(\, ]e_i,e_{i+1}[\,\right)_{0\leq i\leq n-1}$ of the interval $[0,1]$. 
The function $S$ is $C^1$ on $]0,1[$, and $S''$ is uniformly bounded. 
Then, on the interval $[\delta, 1-\delta]$, 
$$
f_\delta = f*\vp_\delta = (J+A+S)*\vp_\delta 
= J_\delta+ A_\delta+ S_\delta.
$$
We now assume that the parameter $\delta$ satisfies
$$
0< \delta <  \min_{0 \leq i \leq n-1} \frac{e_{i+1}-e_i} 2.
$$
For any $0 \leq i \leq n-1$ one has on the interval $[e_i+\delta,e_{i+1}-\delta]$
$$
J_\delta = J, 
\quad 
A_\delta = A
\stext{ and }
S_\delta = f_\delta - A- J.
$$
Hence $J_\delta'' = A_\delta'' = 0$, and $S_\delta'' = f_\delta''$ converges uniformly to $f''$ as $\delta \to 0$, on this interval. Therefore 
\be
\label{eqIntSmooth}
\lim_{\delta \to 0} \sum_{0 \leq i \leq n-1} \int_{e_i+\delta}^{e_{i+1}-\delta} \sqrt{|f_\delta''|} = \int_{]0,1[\sm E} \sqrt{|f''|}.
\ee
Furthermore for any $1\leq i\leq n-1$ and any $|x|\leq 1$ one easily checks that 
$$
J_\delta''(e_i+\delta x) = \frac 1 {\delta^2}[f](e_i) \, \vp'(x)
\stext{ and } 
A_\delta''(e_i+\delta x) = \frac 1 \delta [f'](e_i) \, \vp(x)
$$
where
$$
[f](e_i) := \lim_{\ve\to 0^+} f(e_i+\ve) - f(e_i - \ve) \stext{ and } [f'](e_i) := \lim_{\ve\to 0^+} f'(e_i+\ve) - f'(e_i-\ve).
$$
On the other hand $S''_\delta(e_i+ \delta x)$ is uniformly bounded independently of $x$ and $\delta$. Hence the contribution of $J_\delta''$ is dominant on the intervals $[e_i-\delta, e_i+\delta]$, and we obtain 
\be
\label{eqIntDisc}
\lim_{\delta\to 0} \sum_{1\leq i\leq n-1} \int_{e_i-\delta}^{e_i+\delta} \sqrt{|f_\delta''|} = \sum_{1\leq i \leq n-1} |[f](e_i)| \int_{-1}^1 \sqrt{|\vp'|}.
\ee
Combining \iref{eqIntSmooth} and \iref{eqIntDisc} we conclude the proof of this proposition.
\sq

For the simple oscillating function $f_\omega(x) := \cos(\omega x)$, illustrated on Figure \ref{figTVDf} (left), the two functionals $\TV$ and $D$ behave similarly as $\omega \to \infty$:
\be
\label{eqComplexfCos}
\TV( f_\omega) \simeq \omega \stext{ and } D(f_\omega) \simeq \omega.
\ee
On the contrary for the function $g_\omega(x) := \frac {\lfloor \omega x\rfloor} \omega$, illustrated on Figure \ref{figTVDf} (right) 
we obtain as $\omega\to \infty$
\be
\label{eqComplexfStair}
\TV( g_\omega) \simeq 1 \stext{ and } D( g_\omega)\simeq \sqrt \omega,
\ee
where $D(g_\omega)$ is understood in the sense of \iref{defDDisc}.
Hence the functionals $\TV$ and $D$ treat discontinuities very differently, and the latter penalizes functions which have the appearance of a staircase. 

\begin{figure}
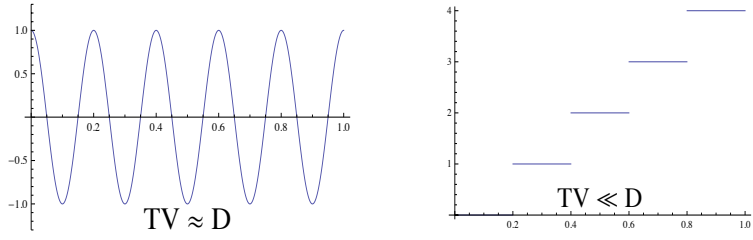

\centering
\includegraphics[height=3cm,width=4.5cm]{\pathPic/Bayesien/Sinf.pdf}
\hspace{1cm}
\includegraphics[height=3cm,width=4cm]{\pathPic/Bayesien/Stairf.pdf}
\caption{ Behavior of the functionals $\TV$ and $D$ on different types of functions.\label{figTVDf}}
\end{figure}

We denote by 
$$
\cF := \R^N.
$$
the collection of discretized one-dimensional functions. For any $F\in \cF$, we denote by $\TV(F)$ and $D(F)$ the discretization of the functionals $\TV$ and $D$ using finite differences, as follows 
\be
\label{eqTVF}
\TV(F) := \sum_{1\leq i\leq N-1} |F_i - F_{i+1}|.
\ee
and 
\be
\label{eqDF}
D(F) := \sum_{2\leq i\leq N-1} \sqrt{|F_{i+1} - 2 F_i+ F_{i-1}|},
\ee

We now turn to numerical results. Given a one dimensional discretized function $F\in \cF := \R^N$ we compute using the algorithm described in \S \ref{secAlgoBayes} three different regularizations of $F$.
\begin{eqnarray*}
F_{\TV}^* &:= & \underset {F'\in \R^N} \argmin \ \alpha_1 \|F-F'\|^2+ \TV(F')\\
F_{\TV} &:= & \frac 1 {Z_2} \int_{\R^N} F' \exp(-\alpha_2 \|F-F'\|^2 - \beta_2 \TV(F') )dF'\\
F_{D} &:= & \frac 1 {Z_3} \int_{\R^N} F' \exp(-\alpha_3 \|F-F'\|^2 - \beta_3 D(F') )dF',\\
\end{eqnarray*}
where $Z_2$ and $Z_3$ are normalizing coefficients. The parameters $\alpha_1, \alpha_2, \alpha_3$ and $\beta_2, \beta_3$ can be freely chosen.

We have not attempted to denoise some data and to compare the $\PSNR$ of the restoration by the different methods.  Indeed such a comparison should, in order to be fair, involve an in depth analysis of the role of the parameters $\alpha_1, \alpha_2, \alpha_3$ and $\beta_2, \beta_3$, that we did not have the time to do. Note also that $F_{\TV}^*$ involves a single parameter, while two parameters are needed for $F_{\TV}$ and $F_D$. We refer to the thesis \cite{LMo} for the comparison of $F_{\TV}^*$ and $F_{\TV}$ from the point of view of the $\PSNR$.

Instead we regard the three procedures as regularization methods, that we apply to different test functions.
We focus on the qualitative differences of the regularizations $F_{\TV}^*$, $F_{\TV}$ and $F_D$ of different functions $F$, and we discuss on how these differences traduce the different bayesian priors implicit in the three methods. \\

\begin{figure}
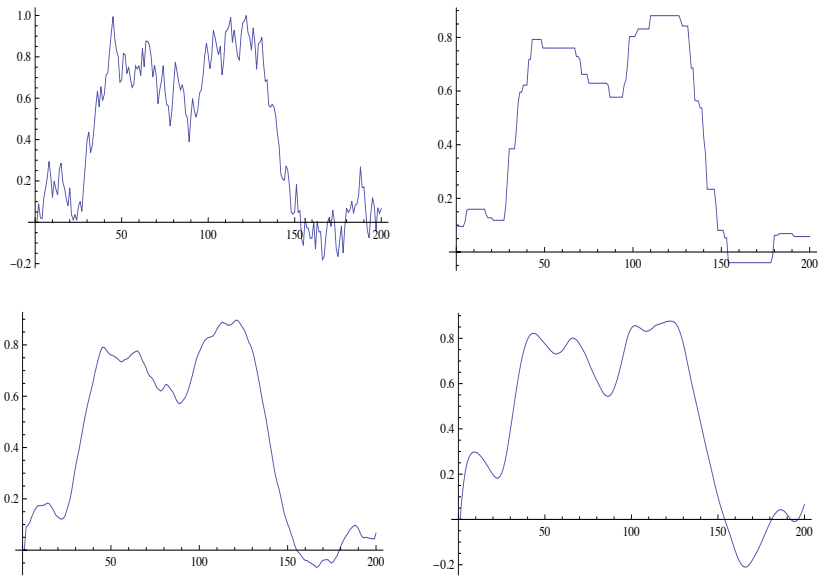

	\centering
	\includegraphics[height = 3.5cm, width = 5cm]{\pathPic/Bayesien/Noise.pdf}\hspace{0.5cm}
	\includegraphics[height = 3.5cm, width = 5cm]{\pathPic/Bayesien/NoiseTVMap.pdf}\vspace{0.5cm}\\
	
	\includegraphics[height = 3.5cm, width = 5cm]{\pathPic/Bayesien/NoiseTV.pdf}\hspace{0.5cm}
	\includegraphics[height = 3.5cm, width = 5cm]{\pathPic/Bayesien/NoiseD.pdf}
	\caption{\label{figRandomWalk}Regularization of a random walk. (Top, Left) Original $F$, (Top, Right) $F^*_{\TV}$, (Bottom,Left) $F_{\TV}$, (Bottom, Right) $F_D$.}
\end{figure}

Our first experiment, see Figure \ref{figRandomWalk}, is the regularization of a random walk. As remarked in \cite{LMo}, $F_{\TV}^*$ is constant over large intervals and also has several sharp discontinuities. These two types of features were not present in the original function $F$ and are undesirable. They are avoided in the two other regularizations $F_{\TV}$ and $F_D$. \\

\begin{figure}
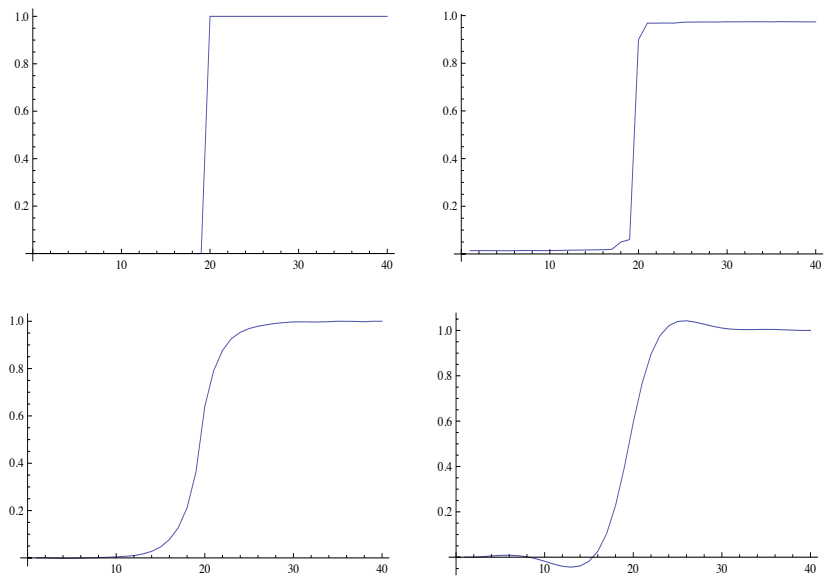

	\centering
	\includegraphics[height = 3.5cm, width = 5cm]{\pathPic/Bayesien/Step.pdf}\hspace{0.5cm}
	\includegraphics[height = 3.5cm, width = 5cm]{\pathPic/Bayesien/StepTVMap.pdf}\vspace{0.5cm}\\
	
	\includegraphics[height = 3.5cm, width = 5cm]{\pathPic/Bayesien/StepTV.pdf}\hspace{0.5cm}
	\includegraphics[height = 3.5cm, width = 5cm]{\pathPic/Bayesien/StepD.pdf}
	\caption{\label{figHeaviside}Heaviside function. (Top, Left) Original $F$, (Top, Right) $F^*_{\TV}$, (Bottom,Left) $F_{\TV}$, (Bottom, Right) $F_D$.}
\end{figure}

Our second experiment, see Figure \ref{figHeaviside} is the regularization of a Heaviside function. Up to numerical artifacts, the functions $F$ and $F_{\TV}^*$ are identical, which is precisely what is wanted in this situation.
The discontinuity disappears in $F_{\TV}$ and $F_D$, and is replaced with a smooth but sharp and well localized transition. Numerical experiments on real (bidimensional) images in \cite{LMo} show that generally, and with appropriate parameters, the regularization $F_{\TV}$ of an image $F$ does not cause a perceptible blurring of the edges which appear in the original image. 

This example puts in light a strong flaw of the regularization method $F_D$ : it does not obey the maximum principle. Indeed, it is clear on Figure \ref{figHeaviside} that $F_{\TV}^*$ and $F_{\TV}$ take their values in the interval $[\min F, \max F]$, but $F_D$ does not.\\

\begin{figure}
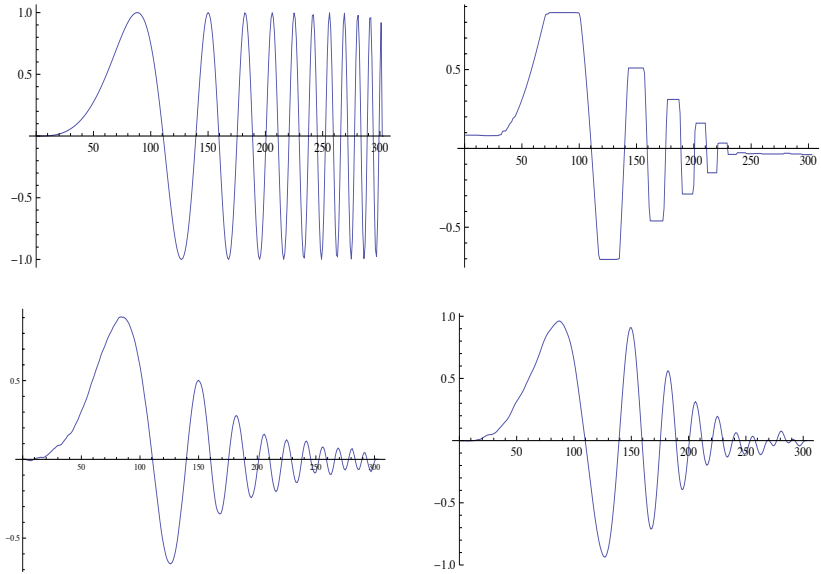

	\centering
	\includegraphics[height = 3.5cm, width = 5cm]{\pathPic/Bayesien/Churp.pdf}\hspace{0.5cm}
	\includegraphics[height = 3.5cm, width = 5cm]{\pathPic/Bayesien/ChurpTVMap.pdf}\vspace{0.5cm}\\
	
	\includegraphics[height = 3.5cm, width = 5cm]{\pathPic/Bayesien/ChurpTV.pdf}\hspace{0.5cm}
	\includegraphics[height = 3.5cm, width = 5cm]{\pathPic/Bayesien/ChurpD.pdf}
	\caption{\label{figChurp}Churp : $\sin(x^3)$ where $x\in [0,4]$. (Top, Left) Original $F$, (Top, Right) $F^*_{\TV}$, (Bottom,Left) $F_{\TV}$, (Bottom, Right) $F_D$.}
\end{figure}

Our third and last experiment, see Figure \ref{figChurp} is the regularization of a ``Churp'' : an oscillating function of increasing frequency, here $\sin (x^3)$ where $x\in [0,4]$.
Qualitatively, most regularization algorithms regard slow oscillations as valuable data which should be preserved, and fast oscillations as noise which should be eliminated. The regularization of a churp allows to analyse the behavior of a regularization algorithm in the intermediate regimes.

From a qualitative point of view, the MAP regularization $F_{\TV}^*$ again introduces undesirable visual artefacts : $F_{\TV}^*$ is exactly constant on some regions, and $\nabla F_{\TV}^*$ has sharp discontinuities which do not correspond to a feature of the original function. In contrast the MMSE regularizations $F_{\TV}$ and $F_D$ do not exhibit these artefacts. As anticipated, the slow oscillations of the original function $F$ are preserved, while the fast oscillations are attenuated. It seems qualitatively that fast oscillations are more strongly attenuated in $F_D$ than in $F_{\TV}$. The author could not prove this last property, but a heuristical analysis suggests that sinuosidal oscillations of large frequency $\omega\gg 1$ are attenuated by a factor $\omega^{-2}$ in 
$F_{\TV}$ and $\omega^{-4}$ in $F_D$.

\section{Extension to higher dimensions and higher order elements}
\label{secCartoon6}
The results on approximation by anisotropic bidimensional piecewise linear
finite elements that we have exposed in \S \ref{secCartoon2} have been generalized in Chapter \ref{chapOptAniso}
to the case of elements of arbitrary order $m-1$
defined on partitions of $\Omega\subset \RR^d$ by simplices.
Here, the local error is defined as
$$
e_{m,T}(f)_p:=\|f-\interp^{m-1}_T f\|_{L^p(T)},
$$
where $\interp^{m-1}_T$ denotes the local interpolation operator
on $\P_{m-1}$ for a $d$-dimensional simplex $T$. This operator
is defined by the condition
$$
\interp^{m-1}_T v(\gamma)=v(\gamma),
$$
for all points $\gamma\in T$ with barycentric coordinates in
the set $\{0, \frac 1 {m-1},\frac 2 {m-1},\cdots,1\}$. We denote by $\H_m\subset \P_m$ the collection of homogeneous polynomials of degree $m$, and we define for any $\bq\in \H_m$ the quantity
$$
K_{m,p}(\bq):=\inf_{|T|=1}e_{m,T}(\bq)_p.
$$
We refer to $K_{m,p}$ as the {\it shape function}. 
For piecewise linear elements in dimension two, i.e. $m=d=2$, 
we have observed that $K_p=K_{2,p}$ has the
special form given by \iref{Kpexp} which justifies the
introduction of the quantity $A_p(f)$. In a similar way, it can 
easily be proved that for piecewise linear elements in higher dimension,
i.e. $m=2$ and $d>2$, one has
$$
c_1 |{\rm det} (\bq)|^{1/d}\leq K_{2,p}(\bq) \leq  c_2 |{\rm det} (\bq)|^{1/d}.
$$
For piecewise quadratic elements in dimension two, i.e. $m=3$ and $d=2$, 
it is proved in Chapter \ref{chapOptAniso} that
$$
c_1 |{\rm disc} (\bq)|^{1/4}\leq K_{3,p}(\bq) \leq  c_2 |{\rm disc} (\bq)|^{1/4}.
$$
for any {\it homogeneous} polynomial $\bq\in \H_3$, where
$$
{\rm disc}(a x^3+ b x^2 y+ c x y^2+ d y^3):= b^2 c^2 - 4 a c^3 - 4 b^3 d + 18 a b c d - 27 a^2 d^2.
$$
For other values of $m$ and $d$, equivalent expressions
of $K_{m,p}(\bq)$ in terms of polynomials in the coefficients of $\bq$
are available but of less simple form, see Chapter 2.

Defining the finite element interpolation error by an optimally adapted partition
$$
\sigma_N(f)_p:=\inf_{\#(\cT)\leq N}\|f-\interp^{m-1}_\cT f\|_{L^p},
$$
where $\interp^{m-1}_\cT$ is the global interpolation operator for the 
simplicial partition (possibly non conforming) $\cT$, the following generalization of 
\iref{limsupest1} is proved in Chapter \ref{chapOptAniso}:
\be
\limsup_{N\to +\infty} N^{\frac m d}\sigma_N(f)_p\leq C_d\left\|K_{m,p}\(\frac {d^mf}{m!}\)\right\|_{L^\tau(\Omega)},
\;\; \frac 1 \tau=\frac 1 p+\frac m d.
\label{limsupest2}
\ee
The constant $C_d$ is equal to $1$ when $d=2$ but larger
than $1$ when $d>2$ due to the impossibility of exactly tiling the
space with locally optimized simplices.
If $f$ is a $C^m$ function of $d$ variables, it is therefore natural to 
consider the quantity 
\be
A_{m,p}(f):=\|K_{m,p}(d^m f)\|_{L^\tau(\Omega)},
\;\; \frac 1 \tau=\frac 1 p+\frac m d,
\ee
as a possible way to measuring anisotropic smoothness. For $d=2$ and piecewise
linear elements, we have seen in \S \ref{secCartoon2} that $A_{2,p}(f)$ is equivalent to the quantity $A_p(f)$.

Similarly to $A_p$ we are interested in the possible extension of $A_{m,p}$
to cartoon functions. We first introduce a generalisation 
of the notion of cartoon functions to higher piecewise smoothness $m$ and dimension $d$.
\begin{definition} 
\label{defcartoongen}
Let $m\geq 2$ and $d\geq 2$ be two integers. Let $\Omega\subset \R^d$ be an open set.
We say that a function $f$ defined
on $\Omega$ is a $C^m$ cartoon function if it is almost everywhere of the form 
$$
f=\sum_{1\leq i\leq k} f_i \Chi_{\Omega_i},
$$
where the $\Omega_i$ are disjoint open sets with piecewise $C^2$ boundary, no cusps (i.e. satisfying an interior and exterior cone condition), and such that $\overline \Omega = \cup_{i=1}^k \overline \Omega_i$.
Additionally, for each $1\leq i\leq k$, the function $f_i$ is assumed to be $C^m$ on $\overline \Omega_i$. 
\end{definition}

Let us consider a fixed cartoon function $f$ on a polyhedral domain $\Omega\subset\R^d$ (i.e. $\Omega$ is
such that $\overline \Omega$ is a closed polyhedron), and a decomposition 
$(\Omega_i)_{1\leq i\leq k}$ of $\Omega$ as in definition \ref{defcartoongen}. As before we define 
$ 
\Gamma := \bigcup_{1\leq i\leq k} \partial \Omega_i,
$
the union of the boundaries of the $\Omega_i$. Our assumptions on the sets $(\Omega_i)_{1\leq i\leq k}$ imply that $\Gamma$ 
is the union of a finite number of open hypersurfaces $(\Gamma_j)_{1\leq j\leq l}$, and of a set $\cP$ of dimension $d-2$.

As in \S \ref{secCartoon3}, we now consider a sequence $f_N$ of piecewise linear approximations of $f$ on 
simplicial partitions $\cT_N$ of cardinality $N$. We distinguish two types of elements of $\cT_N$. A simplex $T\in \cT_N$ is called ``regular'' if $T\cap \Gamma=\emptyset$, and we denote the set of these simplices by $\cT_N^r$. Other simplices are called ``edgy'' and their set is denoted by $\cT_N^e$.
We can again split $\Omega$ according to
$$
\Omega:=(\cup_{T\in \cT_N^r}T) \cup (\cup_{T\in \cT_N^e}T)=\Omega_N^r \cup \Omega_N^e.
$$
Heuristically, if the partitions $\cT_N$ are built with approximation error minimisation in mind, the number of elements should be balanced between $\cT_N^r$ and $\cT_N^e$. The partition $\cT_N^r$ tends to cover most of the surface of $\Omega$, with simplices of diameter $\leq C  N^{-\frac 1 d}$, and $L^\infty$ approximation error $|f-f_N|\leq C N^{-\frac m d}$ (since we use $\P_{m-1}$ elements).
On the other hand, since $f$ has discontinuities along $\Gamma$, the $L^\infty$ approximation error on $\cT_N^e$ does not tend to zero, and $\cT_N^e$ should thus be chosen so as to produce a thin layer around $\Gamma$. 
Let $h$ be the typical diameter of an element of $\cT_N^e$. Since the $\Gamma_j$ has
bounded curvature, this layer can be made of width $\cO(h^2)$
and therefore the layer around $\Gamma$ has volume bounded by
$h^2 \cH_{d-1}(\Gamma)$  up to a fixed multiplicative constant, 
where $\cH_{d-1}(\Gamma)$ is the $d-1$ dimensional
Hausdorff measure of $\Gamma$. On the other hand the minimal number of 
such elements of diameters $h$ needed to cover $\Gamma$ is 
bounded by $h^{1-d}\cH_{d-1}(\Gamma)$ up to a fixed multiplicative constant.
Eventually, we find that the layer around $\Gamma$ has volume bounded by $C N^{-\frac 2 {d-1}}$. 

Hence we have the following heuristic error estimate, for a well designed anisotropic partition:
\begin{eqnarray*}
\|f-f_N\|_{L^p(\Omega)} &\leq& \|f-f_N\|_{L^p(\Omega_N^r)} +\|f-f_N\|_{L^p(\Omega_N^e)}\\
&\leq& \|f-f_N\|_{L^\infty(\Omega_N^r)} |\Omega_N^r|^{\frac 1 p}+  \|f-f_N\|_{L^\infty(\Omega_N^e)} |\Omega_N^e|^{\frac 1 p}\\
& \leq & C (N^{-\frac m d}  + N^{\frac {-2} {p (d-1)}})
\end{eqnarray*}
This leads us to define a critical exponent
$$
p_c = p_c(m,d):=\frac{2d}{m (d-1)}.
$$
If one measures the error in $L^p$ norm with $p>p_c(m,d)$, then the contribution of the edge neighbourhood
$\Omega_N^e$ dominates, while if $p<p_c(m,d)$ it is negligible compared to the contribution
of the smooth region $\Omega_N^r$. 
For the critical exponent $p=p_c(m,d)$ the two terms have the same order, which makes the situation more interesting. Note in particular that $p_c(2,2) =2$, which is consistent with our previous analysis.

For $p\leq p_c(m,d)$, we obtain the approximation rate $N^{-m/d}$ which suggests that approximation
results such as \iref{limsupest2} should also apply to cartoon functions and that the 
quantity $A_{m,p}(f)$ should be finite for such functions. We again need to
use a regularization approach, for the same reasons as in \S \ref{secCartoon3}.
For a given dimension $d$, we consider a radial nonnegative function $\vp$ 
of unit integral and supported in the unit ball of $\R^d$, and we define for $\delta >0$
\be
\vp_\delta(z) := \frac 1 { \delta^d} \vp\left(\frac z \delta\right) \stext{ and } f_\delta = f * \vp_\delta.
\ee

In order to define the quantities of involved in our conjecture, we need to introduce the second fundamental form of an hypersurface.
At any point $x\in \Gamma\sm \cP$ we denote by $\bn(x)$ the unit normal to $\Gamma$. Note that since 
$\Gamma$ is piecewise $C^2$, the map $x\mapsto \bn(x)$ is 
$C^1$ on $\Gamma\sm \cP$.
We define $T_x\Gamma := \bn(x)^\perp$, the tangent space to $\Gamma$ at $x$.
In a neighbourhood of $x\in \Gamma\sm \cP$, the hypersurface $\Gamma$ admits a parametrization of the form
$$
u\in T_x\Gamma \mapsto x + u+ \lambda(u) \bn(x) \in \Gamma \sm \cP,
$$
where $\lambda$ is a scalar valued $C^2$ function.
By definition, the second fundamental form of $\Gamma$ at the point $x$ is the quadratic form $\II_x$ associated
to $d^2\lambda(0)$ which is defined on $T_x \Gamma\times T_x \Gamma$.
Alternatively, for all $u,v\in T_x\Gamma$ we have $\II_x(u,v) := -\<\partial_u \bn, v\>$.
The Gauss curvature $\kappa(x)$ is the determinant of $\II_x$, in any orthonormal basis of $T_x \Gamma$,
$$
\kappa(x):= \det \II_x.
$$
For example, in two space dimensions the tangent space $T_x\Gamma$ is one dimensional, and we simply have 
$\II_x(u,v) = \kappa(x) \<u,v\>$.
We also denote by $\sigma(x)\in \{0, \cdots, d-1\}$ the signature of the quadratic form $\II_x$, which is defined as the number of its positive eigenvalues. 

With $\tau$ such that $\frac 1 \tau :=\frac m d+\frac 1 p$, we 
define
$$
S_p(f) :=  \|K(d^m f )\|_{L^\tau(\Omega\sm \Gamma)}=A_p(f_{|\Omega\sm \Gamma}).
$$
We conjecture the following generalization to Theorem \ref{threg}.

\begin{conjecture}
There exists $d$ positive constants $C(k)$, $k\in \{0,\cdots, d-1\}$, that depend on  $\vp,p,m,d$, such that,
with 
$$
E_p(f):= \| C(\sigma) |\kappa|^{\frac m {2d}} [f] \|_{L^\tau(\Gamma)}=
\(\int_\Gamma \left |C(\sigma(x)) |\kappa(x)|^{\frac m {2d}} [f(x)] \right |^\tau dx\)^{\frac 1 \tau},
$$
we have
\label{thAmdp}
\begin{itemize}
\item If $p<p_c$ then 
$$
\lim_{\delta\to 0} A_{m,p}(f_\delta) = S_p(f).
$$ 
\item If $p=p_c$ then 
$$
\lim_{\delta \to 0} \left(A_{m,p}(f_\delta)\right)  = \left(S_p(f)^\tau +  E_p(f)^\tau\right)^{1/\tau}.
$$ 
\item If $p>p_c$ then
$$
\lim_{\delta \to 0}  \delta^{\frac 1 {p_c} - \frac 1 p} A_{m,p}(f_\delta) =  E_p(f).
$$ 
\end{itemize}
\end{conjecture}

In the remainder of this section, we give some 
arguments that justify this conjecture. Given a cartoon function $f$, 
we define the sets $\Omega_\delta$, $\Gamma_\delta^0$, 
$\Gamma_\delta$ and $\cP_\delta$ similarly to \S \ref{secCartoon3}.
We need to perform an asymptotic analysis of the integral
\be
\label{tointegrate}
\int_\Omega K(d^m f_\delta)^\tau = \int_{\Omega_\delta} K(d^m f_\delta)^\tau+\int_{\cP_\delta} K(d^m f_\delta)^\tau
+\int_{\Gamma_\delta} K(d^m f_\delta)^\tau.
\ee
As in the proof of Theorem \ref{threg} the contribution of $\cP_\delta$ can be proved to be negligible 
compared with those of $\Omega_\delta$ and $\Gamma_\delta$ as $\delta \to 0$. The contribution of
$\Omega_\delta$ satisfies
$$
\lim_{\delta \to 0} \int_{\Omega_\delta} K(d^m f_\delta)^\tau =\int_{\Omega\sm \Gamma} K(d^m f)^\tau.
$$ 
The main difficulty lies again in the contribution of $\Gamma_\delta$. Let us define $\tau_c$ by 
$$
\frac 1 {\tau_c} := \frac m d+ \frac 1 {p_c}.
$$
The contribution of $\Gamma_\delta$ can be 
computed if one can establish an estimate generalizing \iref{estimK2d} according to 
\be
\label{estimK}
\left|  \delta^{\frac 1 {\tau_c}} K_{m,p}(d^m f_\delta(z)) - |[f](x)| |\kappa(x)|^{\frac m {2 d}} \Phi_{m,d,\sigma(x)}(u) \right|\leq \omega(\delta)
\ee
where $\omega(\delta)\to 0$ as $\delta\to 0$, $x\in \Gamma_\delta^0$, $u\in [-1,1]$, $z = x+\delta u \bn(x)$,
and where the function $\Phi_{m,d,k} : [-1,1]\to \R$ only depends on $m,d,k$ and $\vp$. If
\iref{estimK} holds, we then easily derive that
$$
\lim_{\delta \to 0} \delta^{\frac {\tau}{\tau_c}-1}\int_{\Gamma_\delta} K(d^m f_\delta)^\tau=\int_\Gamma C(\sigma) |\kappa|^{\frac {\tau m} {2d}} |[f]|^\tau,
$$
with $C(k):=\int_{-1}^1|\Phi_{m,d,k}(u)|^\tau du$, which leads to the proof of the conjecture.

We do not have a general proof of \iref{estimK} for any $m$, $p$ and $d$. 
In the following, we justify its validity in
two particular cases for which the explicit expression of $K_{m,p}$ is known to us:
piecewise quadratic in two space dimensions ($d=2$ and $m=3$)
and piecewise linear  in any dimension ($m=2$).

\paragraph{Piecewise quadratic elements in two dimensions.}
For all $\delta>0$, $x\in \Gamma_\delta^0$ and $u\in [-1,1]$, let $\pi_{x,\delta,u}\in \H_3$
 be the homogeneous cubic polynomial on $\R^2$ corresponding to $d^3 f_\delta (x+\delta u \bn(x))$.
Let also $\pi_{x,u}\in \H_3$ be the homogeneous cubic polynomial on $\R^2$ defined by 
\be
\label{eqThirdCartoon}
\pi_{x,u}(\lambda \bn(x)+\mu \bt(x)) = -\lambda(\Phi''(u) \lambda^2- 3\Phi'(u) \kappa(x) \mu^2)
\ee
for all $(\lambda,\mu)\in \R^2$, where $\Phi$ is defined
by \iref{defchi}. For all $x\in \Gamma$, we denote by $M_{x,\delta}$ the (symmetric) linear map defined by 
$$
M_{x,\delta}\bn(x) = \delta \bn(x) \text{ and } M_{x,\delta} \bt(x)= \sqrt \delta \bt(x).
$$
Then, using a reasoning similar to the one used in the appendix of this chapter, it can be proved that
\be
\label{limdiscD}
\|\pi_{x,\delta,u} \circ M_{x,\delta} -[f](x)\pi_{x,u}\|\leq \omega(\delta).
\ee
where $\lim_{\delta\to 0} \omega(\delta) = 0$ and the function $\omega$ depends only on $f$.
Furthermore, it is proved in Chapter \ref{chapOptAniso} that for all $\bq\in \H_3$
$$
K_{3,p}(\bq) = C \sqrt[4]{|\disc \bq|}
$$
where the positive constant $C$ depends on $p$ and the sign of $\disc \bq$. Combining this expression with \iref{limdiscD} proves \iref{estimK} and thus the conjecture in the case $m=3$ and $d=2$.

\paragraph{Piecewise linear elements in any dimension.} We use the second 
fundamental form of the discontinuity set $\Gamma$ in order to evaluate
$d^m f_\delta$ on $\Gamma_\delta$.
Characteristic functions are one of the simplest types of cartoon functions. In that case, it is possible
to establish a simple relation between the second fundamental form of
the edge set and the second derivatives of $f$ in a distributional sense:
if $\Omega\subset \R^d$ is a bounded domain with smooth boundary $\Gamma$ and inward normal $\bn$,
we then have for all $C^2$ test function $\psi$ and $u,v\in \R^d$ 
\be
-\int_\Omega \partial^2_{u,v}\psi = \int_\Gamma \<u,\bn\>\<v,\bn\> (\partial_\bn \psi - \Tr(\II) \psi) +\II'(u,v)\psi,
\label{d2charR2eqd}
\ee
where $\II'_x(u,v)$ is the second fundamental form $\II_x$ applied to the orthogonal projection of $u$ and $v$ on $T_x\Gamma$.
The proof of this formula (that generalizes the simpler bidimensional case \iref{d2charR2eq} which is proved in the appendix) is given further below. For all $x\in \Gamma$, 
we denote by $M_{x,\delta}$ the (symmetric) linear map defined by 
$$
M_{x,\delta}\bn(x) = \delta \bn(x) \text{ and } M_{x,\delta} \bt = \sqrt \delta \bt  
$$
for all $\bt\in T_x\Gamma$.
For all $\delta>0$, $x\in \Gamma_\delta^0$ and $u\in [-1,1]$, let $\pi_{x,\delta,u}\in \H_2$ be the homogeneous quadratic polynomial on $\R^d$ corresponding to $d^2 f_\delta (x+\delta u \bn(x))$.
Let also $\pi_{x,u}\in \H_2$ be the homogeneous quadratic polynomial on $\R^d$ defined by 
$$
\pi_{x,u}(\lambda \bn(x)+\bt) = \Phi'(u) \lambda^2- \Phi(u) \II_x(\bt,\bt)
$$
for all $\lambda\in \R$ and $\bt\in T_x \Gamma$, where $\Phi(x):=\int_{\R^{d-1}}\vp(x,y)dy$.
Then, using \iref{d2charR2eqd} and a reasoning analogous to the one presented in the appendix, 
it can be proved that 
\be
\label{limQuadD}
\|\pi_{x,\delta,u} \circ M_{x,\delta} -[f](x)\pi_{x,u}\|\leq \omega(\delta).
\ee
where $\lim_{\delta\to 0} \omega(\delta) = 0$ and the function $\omega$ depends only on $f$.
Furthermore, it is proved in Chapter \ref{chapOptAniso} that 
$$
K_{2,p}(\bq) = C \sqrt[d]{|\det \bq|}
$$
where the positive constant $C$ depends on $d,p$ and the signature of $\bq\in \H_2$. Combining this expression with \iref{limQuadD} proves the estimate \iref{estimK} and thus the conjecture in the case $m=2$ in any dimension $d>1$.
\nl
\nl
{\bf Proof of \iref{d2charR2eqd}:}
Let $\proj_\Gamma$ be the orthogonal projection onto $\Gamma$, and for all $x\in \Gamma$ let $\proj_x$ be the orthogonal projection onto $T_x\Gamma$. 
We consider a vector $u\in \R^d$ and we define $\move : \Gamma\to\Gamma$ by $\move (x) := \proj_\Gamma(x+u)$. If $\|u\| \|\II\|_{L^\infty(\Gamma)}<1$, then $\move$ is smooth and its differential $d_x \move : T_x\Gamma\to T_{x'}\Gamma$, where 
$$
x' = \move(x),
$$ 
is given by the following formula
$$
d_x \move = (\Id - \<u, \bn(x')\>\II_{x'})^{-1} \proj_{x'} 
$$
The determinant of $d_x U$ (more precisely the determinant of the matrix of $d_x U$ in direct orthogonal bases of $T_x\Gamma$ and $T_{x'} \Gamma$) is 
$$
\det (d_x \move) = \det(\Id - \<u, \bn(x')\>\II_{x'})^{-1}\<\bn(x),\bn(x')\> = 1+\<u, \bn(x')\> \Tr(\II_{x'})+ \|u\|\omega_1(u,x).
$$
where 
$\omega_1(u,x)$ tends uniformly to $0$ as $u\to 0$.
Furthermore, it is easy to show that 
$$
|\psi(x+u) - \psi(x') -  \<u, \bn(x)\> \partial_{\bn(x')} \psi(x')|\leq C\|u\|^2,
$$
and $\|n(x') - n(x) - \II_{x'}(\proj_{x'}(u))\|\leq \|u\|\omega_2(u)$, where $C$ and $\omega_2$ are independent of $x\in \Gamma$ and $\omega_2(u) \to 0$ as $u \to 0$ (the quadratic form $\II_x(r,s) := -\<\partial_r \n(x),s\>$ for all $r,s\in T_x\Gamma$ is identified here to the differential of $\n$). Combining these results, we obtain
\begin{eqnarray*}
& &\int_\Gamma \psi(x+u)\<\bn(x),v\> dx\\
&=& \int_\Gamma \psi(x+u)\<\bn(x),v\> \det (d_{x'} \move_u)^{-1} dx'\\
&=&\int_\Gamma \<\bn(x'), v\> \psi(x')dx'\\
& &+ \int_\Gamma \<\bn(x'), v\> \<\bn(x'),u\>(\partial_{\bn(x')} \psi- \Tr(\II_x')\psi(x'))
 +\<v,\II_x'(\proj_{x'}(u))\>\psi(x') dx' \\
& &+ \|u\|\omega_3(u).
\end{eqnarray*}
where $\omega_3(u)\to 0$ as $u \to 0$.
We conclude the proof of \iref{d2charR2eqd} using the formula
$$
- \int_\Omega \partial^2_{u,v}\psi = \lim_{h\to 0} h^{-1}\int_\Gamma (\psi(x+hu)-\psi(x))\<\bn(x),v\> dx.
$$ 
\sq

\begin{remark}
Similarly to the results presented in \S \ref{secCartoon4}, there is an affine invariance property associated to $\kappa$: if $T$ is an affine transformation of $\R^d$ with linear part $L$, 
and if $f = \ti f\circ T$, $\ti \Gamma = T(\Gamma)$ and $\ti \kappa$ is the Gauss curvature of $\ti \Gamma$, then
one has for any $s \geq 0$, 
$$
(\det L)^{\frac{d-1}{d+1}}\int_\Gamma |C(\sigma) [f]|^s |\kappa|^{\frac 1 {d+1}} = \int_{\ti\Gamma} |C(\ti \sigma)[\ti f]|^s |\ti \kappa|^{\frac 1 {d+1}}.
$$
It follows from this observation that when $p=p_c$, the contribution of the edges is affine invariant in the sense
that 
$$
E_{p_c}(\ti f) = (\det L)^{\frac {d-1}{d+1}} E_{p_c}(f).
$$
Since one also has $A_{m,p}(\ti f) = (\det L)^{\frac {d-1}{d+1}} A_{m,p} (f)$ this comforts the conjecture.
Let us mention that the quantity $|\kappa|^{\frac 1 {d+1}}$ has been used in \cite{Ol} 
in order to define surface smoothing operators that are invariant under affine change of coordinates.
\end{remark}

\section{Conclusion}
\label{secCartoon7}
In this chapter we have investigated the quantity $A_p(f)$ which
governs the rate of approximation by anisotropic $\P_1$ finite elements as a way to
describe anisotropic smoothness of functions. This quantity is not a semi-norm
due to the presence of the non-linear quantity $\det(d^2f)$ 
and cannot be defined in a straightforward manner for general distributions.
We nevertheless have shown that this quantity can be defined
for cartoon images with geometrically smooth edges when $p\leq 2$.
A theoretical issue remains to give a satisfactory meaning to 
the full class of functions for which this quantity is finite.

From a more applied perspective, it could be interesting to
investigate the role of $A_p(f)$ in problems where anisotropic
features naturally arise: 

\begin{enumerate}
\item
Approximation of PDE's: in the case of one dimensional hyperbolic conservation laws, it was proved
in \cite{DLu} that despite the appearance of discontinuities
the solution has high order smoothness in Besov spaces
that govern the rate of adaptive approximation by piecewise polynomials.
A natural question is to ask wether similar results hold
in higher dimension, which corresponds to understanding
if $A_p(f)$ remains bounded despite the appearance of shocks.
\item
Image processing: as illustrated in \S \ref{secCartoon5}, the quantity $A_p(f)$ can easily be discretized
and defined for pixelized images. It is therefore tempting to use
$A_2(f)$ in a similar way
as the total variation in \iref{minBV}, by solving a problem of the form
\be
\min_{g\in BV} \{A_2(g)\; ; \; \|Tg-h\|_{L^2}\leq \e\},
\label{minA2}
\ee
with the objective of promoting images with piecewise smooth edges. The main
difficulty is that $A_2$ is not a convex functional. One way to solve this
difficulty could be to reformulate \iref{minA2} in a Bayesian framework as the search of a maximum
of an a-posteriori probability distribution (MAP) as an estimator of $f$. In this framework, we
may instead search for a minimal mean-square error estimator (MMSE),
and this search can be implemented by stochastic algorithms which
do not require the convexity of $A_2$, see \cite{LMo} and \S\ref{secCartoon5.5}.
\end{enumerate}

\section{Appendix: proof of the estimates \iref{estimnn2d}-\iref{estimnt2d}-\iref{estimtt2d}.}

It is known since
the work of Whitney on extension theorems 
(see in particular \cite{Whitney}) 
that for any open set $U\subset \R^d$, and any $g\in C^2(\overline U)$ there exists $\ti g \in C^2(\R^d)$ such that $\ti g_{|U} = g$.
It follows that for each $1\leq i \leq k$, there exists $\ti f_i\in C^2(\R^2)$, compactly supported, and such that $\ti f_{i|\Omega_i} = f_i$. 

Let $\Gamma_j$ be one of the pieces of $\Gamma$, between the domains $\Omega_k$ and $\Omega_l$, and let  $s=\ti f_k$ and $t = \ti f_l- \ti f_k$. Although the domains $\Omega_k$ and $\Omega_l$ are only piecewise smooth, there exists an open set $\Omega'$ with $C^2$ boundary such that for $\delta_0>0$ small enough
$$
f = s \Chi_{\Omega'} + t \ \text{ on } \ \bigcup_{0<\delta\leq \delta_0} (\Gamma_{j,\delta}+B_\delta),
$$
where $B_\delta$ is the ball of radius $\delta$ centered at $0$.
Note that $\Gamma_j\subset \Gamma' := \partial \Omega'$ and that $s = [f]$ on $\Gamma_j$.
In the following, the variables $x,z$ are always subject to the restriction
\be
\label{restrictXZ}
x\in \Gamma_{j,\delta}^0 \text{ and } z = U_\delta(x,u) = x+\delta u \bn(x) \text{ where } 0<\delta\leq \delta_0 \text{ and } |u|\leq 1,
\ee
note that $z\in \Gamma_{j,\delta}$ and $\|x-z\|\leq \delta$. We therefore have 
$$
f_\delta(z) = \int_{\Omega'} s(\ti x) \vp_\delta(z-\ti x) d\ti x + t_\delta(z),
$$
where $t_\delta:=t * \vp_\delta$. The second derivatives of $t_\delta$ are uniformly bounded, and are therefore negligible in regard of all three estimates \iref{estimnn2d}, \iref{estimnt2d} and \iref{estimtt2d}, indeed
$$
\|d^2 t_\delta\|_{L^\infty} = \|(d^2 t)* \vp_\delta\|_{L^\infty} \leq \|d^2 t\|_{L^\infty} \|\vp_\delta\|_{L^1} = \|d^2 t\|_{L^\infty} \|\vp\|_{L^1} <\infty.
$$
We now define the $2\times 2$ symmetric matrices
$$
I(z,x) := \int_{\Omega'} (s(\ti x)-s(x)) \ d^2 \vp_\delta(z-\ti x) \ d\ti x \text{ and } J(z) := \int_{\Omega'} d^2 \vp_\delta(z-\ti x) d\ti x 
$$
so that 
\be
d^2 f_\delta(z) =  d^2 t_\delta + I(z,x) + [f](x)  J(z).
\label{d2fdelta}
\ee
We already know that the contribution of $d^2t_\delta$ is negligible. We now prove that
the same holds for the contribution of $I(z,x)$. Since $\vp_\delta(z-\ti x)$ is non-zero
only if $\|\ti x-z\| \leq \delta$ and therefore $\|\ti x-x\| \leq 2\delta$, we
can bound the norm of the matrix $I(z,x)$ by
\be
\label{d2f1}
\|I(z,x)\| \leq 2\delta \|ds\|_{L^\infty} \|d^2 \vp_\delta\|_{L^1} \leq 2\delta \|ds\|_{L^\infty} \|d^2\vp\|_{L^1} \delta^{-2} = C\delta^{-1}.
\ee
This proves that the contribution of $I(z,x)$ is negligible for the two estimates
\iref{estimnn2d} and \iref{estimnt2d}. In order to prove that it is
also negligible in the estimate \iref{estimtt2d},
we need a finer analysis of $\bt(x)^\trans I(z,x) \bt(x)$.
For this purpose we fix a unit vector $u$ and the pair $(x,z)$. We introduce
$$
\Lambda(\ti x) := (s(\ti x)-s(x)) \ \partial_u \vp_\delta(z-\ti x) + \partial_u s(\ti x) \ \vp_\delta(z-\ti x),
$$
so that by Leibniz rule
$$
(s(\ti x)-s(x)) \ \partial^2_{u,u} \vp_\delta(z-\ti x)= \partial^2_{u,u} s(\ti x) \ \vp_\delta(z-\ti x) - \partial_u \Lambda(\ti x).
$$
Therefore
\begin{eqnarray*}
u^\trans I(z,x)u &=& \int_{\Omega'} \left( \partial^2_{u,u} s(\ti x) \ \vp_\delta(z-\ti x) - \partial_u \Lambda(\ti x)\right) d\ti x, \\
&=&  \int_{\Omega'} \partial^2_{u,u} s(\ti x) \ \vp_\delta(z-\ti x) d\ti x - \int_{\Gamma'} \Lambda(\ti x) \<\bn(\ti x), u\> d \ti x.
\end{eqnarray*}
The first integral clearly satisfies 
$$
\left |\int_{\Omega'} \partial^2_{u,u} s(\ti x) \ \vp_\delta(z-\ti x) d\ti x \right |\leq \|d^2 s\|_{L^\infty} \|\vp_\delta\|_{L^1},
$$
and is therefore bounded independently of $\delta$. We estimate the second integral 
for the special case $u = \bt(x)$, remarking that  $|\<\bn(\ti x), \bt(x)\>|\leq C_1 \delta$ on the domain of integration.
Therefore
$$
\left| \int_{\Gamma'} \Lambda(\ti x) \<\bn(\ti x), \bt(x)\> d \ti x\right| \leq C_1 \delta |\Gamma'\cap B(z,\delta)| \|\Lambda\|_{L^\infty},
$$
where, slightly abusing notations, we denote by $|\Gamma'\cap B(z,\delta)|$ the length ($1$-dimensional Hausdorff measure) of the curve $\Gamma'\cap B(z,\delta)$.
Clearly $\Lambda(\ti x) = 0$ if $\|z-\ti x\|\geq \delta$. If $\|z-\ti x\|\leq \delta$ we have 
\be
\label{upperPsi}
|\Lambda(\ti x)|\leq (\|x-z\|+\|z-\ti x\|) \|ds\|_{L^\infty} \|d\vp\|_{L^\infty} \delta^{-3} + \|ds\|_{L^\infty} \|\vp\|_{L^\infty} \delta^{-2} \leq C_0\delta^{-2}.
\ee
Since in addition $|\Gamma'\cap B(z,\delta)|\leq C_2 \delta$, we finally find that
$$
\left| \int_{\Gamma'} \Lambda(\ti x) \<\bn(\ti x), \bt(x)\> d \ti x\right| \leq C_0C_1C_2.
$$
We have therefore proved that 
 $$
|\bt(x)^\trans I(z,x) \bt(x)|\leq C,
 $$
 where the constant $C$ is independent of $\delta$, which shows that 
 the contribution of $I(z,x)$ is negligible in \iref{estimtt2d}.

We now analyze the contribution the quantity $[f](x)J(z)$ in \iref{d2fdelta}.
For this purpose, we use an expression of the second derivative of 
the characteristic function $\Chi_{\Omega'}$ of a smooth set $\Omega'$
in the distribution sense. 
We assume without loss of generality that $\Gamma'$ is parametrized 
in the trigonometric sense, and therefore that $\bn$ is the inward normal to $\Omega$.
For all test function $\psi$, we have
$$
-\int_{\Omega'} \partial^2_{u,v} \psi = \int_{\Gamma'} \partial_u \psi \<v,\bn\> =  \int_{\Gamma'} (\partial_\bn \psi \<u,\bn\> + \partial_\bt\psi \<u,\bt\>) \<v,\bn\>
$$
and, by integration by parts, 
$$
\int_{\Gamma'} \partial_\bt \psi \<u,\bt\>\<v,\bn\> = - \int_{\Gamma'} \psi (\<u,\kappa\bn\>\<v,\bn\>- \<u,\bt\>\<v,\kappa\bt\>). 
$$
Therefore, we have
\be
\label{d2charR2eq}
-\int_{\Omega'} \partial^2_{u,v} \psi = \int_{\Gamma'} \<u,\bn\>\<v,\bn\> (\partial_\bn \psi-\kappa \psi) +\kappa \<u,\bt\>\<v,\bt\> \psi
\ee
Applying this formula to $\psi(\ti x):=\vp_{\delta}(z-\ti x)$ we obtain
\be
\begin{array}{rcl}
\displaystyle-u^\trans J(z) v &=&  \displaystyle\int_{\Gamma'} \<u,\bn(\ti x)\>\<v,\bn(\ti x)\> (\partial_\bn \vp_\delta(z-\ti x) - \kappa(\ti x)\vp_\delta(z-\ti x))d \ti x\\
& & \displaystyle+\int_{\Gamma'} \kappa(\ti x) \<u,\bt(\ti x)\>\<v,\bt(\ti x)\> \vp_\delta(z-\ti x) d\ti x.
\end{array}
\label{uJv}
\ee
Since $\Gamma_j$ is $C^2$, there exists a constant $C_0$ such that for all $x_1,x_2\in \Gamma_j$, we have 
$$
|\<\bt(x_1),\bn(x_2)\>|\leq C_0 \|x_1- x_2\|,
$$
and 
$$
|1-\<\bn(x_1),\bn( x_2)\>| = |1-\<\bt(x_1),\bt( x_2)\>|\leq C_0 \|x_1- x_2\|^2.
$$
We finally remark that $ |\Gamma'\cap B(z,\delta)|\leq C_1 \delta$, and that 
$\|\vp_\delta\|_{L^\infty}\leq \|\vp\|_{L^\infty}\delta^{-2}$ and  $\|\partial_\bn \vp_\delta\|_{L^\infty} \leq \|d\vp\|_{L^\infty} \delta^{-3}$. 

Taking the vectors $\bt(x)$ or $\bn(x)$ as possible values of $u$ and $v$ in \iref{uJv} and using the
above remarks, we obtain the estimates
\begin{eqnarray}
\label{threeineq1}
\left|\bn(x)^\trans J(z)\bn(x) + \int_{\Gamma'} \partial_\bn \vp_\delta(z-\ti x) d \ti x\right| &\leq& C \delta^{-1},\\
\label{threeineq2}
|\bt(x)^\trans J(z)\bn(x)| &\leq& C \delta^{-1},\\
\label{threeineq3}
\left|\bt(x)^\trans J(z)\bt(x) + \int_{\Gamma'} \kappa(\ti x) \vp_\delta(z-\ti x)d \ti x\right| &\leq& C,
\end{eqnarray}
where the constant $C$ depends only on $f$. In view of \iref{d2fdelta}
we can immediately derive estimate
\iref{estimnt2d} from \iref{threeineq2}.

In order to derive the estimate \iref{estimtt2d} from \iref{threeineq3}, we first
introduce the modulus of continuity  $\omega$ of $\kappa$ on $\Gamma_j$, 
$$
\omega(\delta) := \sup_{x_1, x_2\in \Gamma_j\sep \|x_1- x_2\|\leq \delta} |\kappa(x_1)-\kappa( x_2)|.
$$
Therefore
\begin{eqnarray*}
\left|\bt(x)^\trans J(z)\bt(x) + \int_{\Gamma'} \kappa(\ti x) \vp_\delta(z-\ti x)d \ti x\right| 
&\leq& \left|\bt(x)^\trans J(z)\bt(x) + \kappa(x)\int_{\Gamma'}  \vp_\delta(z-\ti x)d \ti x\right| \\
& &+ C\omega(\delta)\delta^{-1}.
\end{eqnarray*}
We now claim that 
\be
\left| \int_\Gamma \vp_\delta( z-\ti x) d\ti x - \delta^{-1} \Phi(u) \right| \leq C,
\label{phiPhi}
\ee
holds with $C$ independent of $\delta$ which implies the validity of  \iref{estimtt2d}. 
In order to prove \iref{phiPhi}, we use a local parametrization of $\Gamma'$: let 
$\lambda:\R\to \R$ be such that for $h$ small enough we have, $x+h\bt(x)+\lambda(h)\bn(x)\in \Gamma$. 
Note that we have $|\lambda(h)|\leq C_0 h^2$ and $|\lambda'(h)|\leq C_0h$ for $h$ small enough. 
Then for $\delta$ small enough,
\begin{eqnarray*}
\left| \int_{\Gamma'} \vp_\delta( z-\ti x) d\ti x - \delta^{-1} \Phi(u) \right|  &\leq&  \left|\int_{\R} \vp_\delta(h\bt(x)+ (\delta u-\lambda(h))\bn(x))\sqrt{1+\lambda'(h)^2} dh \right.\\
& &\left.- \int_\R  \vp_\delta(h\bt(x)+ \delta u\bn(x)) dh \right|\\
&\leq & C\delta(\|\vp_\delta\|_{L^\infty} (\sqrt{1+(C_0\delta)^2}-1)+  \|d\vp_\delta\|_{L^\infty} C_0 \delta^2 ) \leq C
\end{eqnarray*}
Finally, we can derive the estimate \iref{estimnn2d} from \iref{threeineq1} using
the inequality
\be
\left| \int_\Gamma \partial_\bn \vp_\delta(z - \ti x) d\ti x + \delta^{-2} \Phi'(u)\right|  \leq C \delta^{-1}
\label{phiPhin}
\ee
which proof is very similar to the one of \iref{phiPhi}. 



\part{Mesh adaptation and riemannian metrics} 
\label{partRiemann}

\chapter{Are riemannian metrics equivalent to simplicial meshes ?} 
\minitoc
\label{chapMeshMet}

\section{Introduction}

Triangulations and meshes are finite objects of combinatorial nature: they can be described by a collection of vertices and of connections between these. This description well adapted to the demonstration of algebraic results, such as the Euler formula, or for computer processing.
In contrast, many approaches towards anisotropic mesh generation are based on a continuous object, namely a riemannian metric $z\mapsto H(z)$, in other words a continuous function $H$ from the domain $\Omega\subset \R^d$ to the set $S_d^+$ of symmetric positive definite matrices. Once this metric has been properly designed, it is the task of a mesh generation algorithm such as \cite{Shew,Bois,FreeFem,Inria} to generate a triangulation that agrees with this metric, in a sense specified below \iref{defEquiv}.
The purpose of this chapter is to formulate precise equivalence results between some classes of triangulations and of riemannian metrics. This equivalence translates some geometrical constraints satisfied by the triangulations into the form of regularity properties of the equivalent riemannian metrics.\\

Our results are so far limited to meshes and metrics defined on the entire infinite domain $\R^d$, and the dimension $d\geq 2$ is fixed throughout this chapter. 
The choice of the unbounded domain $\R^d$ is guided by simplicity, since the curvature and the singularities of
the boundary of a bounded domain induce additional difficulties 
from the point of view of computational mesh generation.
Bounded domains will be the object of future work.
We denote by $\bT$ the collection of conforming simplicial meshes of $\R^d$, and we introduce in \S\ref{subsecDefTStar} three embedded subsets of $\bT$ of particular interest 
$$
\bT_{i,C} \subset \bT_{a,C} \subset \bT_{g,C},
$$
where $C\geq 1$ is a parameter which plays a minor role.
The set $\bT_{i,C}$ collects all \emph{isotropic} meshes $\cT$ which are heuristically defined as follows: the simplices $T\in \cT$ may have strongly varying volumes, but their aspect ratio is uniformly bounded. 
The largest set $\bT_{g,C}$ collects \emph{graded} meshes, which only satisfy a condition of local consistency: all the  geometrical features of the simplices $T\in \cT$ may vary strongly, volume, aspect ratio and orientation, but two \emph{neighboring} simplices should not be excessively different.
The intermediate set $\bT_{a,C}$ of \emph{quasi-acute} meshes, is defined by a condition which involves the measure of sliverness $S(T)$ of a simplex $T$ introduced in Chapter \ref{chapW1P} and also discussed in Chapter \ref{chapApproxMet}. 
The measure of sliverness plays an important role in the finite element approximation of a function on a mesh $\cT$ when the error is measured in the Sobolev $W^{1,p}$ norm.\\

We associate to each simplex $T$ a symmetric positive definite matrix $\cH_T$ such that the ellipsoid
\be
\label{caractHT}
\cE_T := \{z\in \R^d \sep (z-z_T) \cH_T (z-z_T)^\trans \leq 1\},
\ee
which is centered at the barycenter $z_T$ of $T$,
is the ellipsoid of minimal volume containing $T$. Two triangles $T$ and their associated ellipses $\cE_T$ are illustrated on Figure \ref{figET}. We give in \S \ref{subsecIntroSimplex} the explicit expression of the matrix $\cH_T$ in terms of the coordinates of vertices of the simplex $T$, and we discuss its main properties. The matrix $\cH_T\in S_d^+$ encodes the volume, the aspect ratio and the orientation (but not the angles) of the simplex $T$. 

\begin{figure}
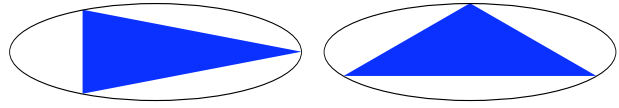

\centering
	\includegraphics[width = 4cm, height=1.5cm]{\pathPic/Triangles/TriangleAigu.pdf}
	\includegraphics[width = 4cm, height=1.5cm]{\pathPic/Triangles/TriangleObtus.pdf}
	\caption{Two triangles $T$ and their associated ellipse $\cE_T$.\label{figET}}
\end{figure}

A riemannian metric on $\R^d$ is a continuous map $H$ which associates to any $z\in \R^d$ a symmetric positive definite matrix $H(z)\in S_d^+$. We denote by $\bH = C^0(\R^d,S_d^+)$ the collection of riemannian metrics.
The following definition introduces a notion of equivalence between meshes and metrics
\begin{definition}
\label{defEqHT}
We say that a mesh $\cT\in \bT$ is $C_0$-equivalent to a given metric $H\in \bH$, where $C_0\geq 1$ is a fixed constant, if for all $T \in \cT$ and all $z\in T$ one has 
\be
\label{defEquiv}
C_0^{-2} H(z) \leq \cH_T \leq C_0^2 H(z)
\ee
We say that a collection of simplicial meshes $\bT_*\subset \bT$ is equivalent to a collection $H_*$ of metrics if there exists a uniform constant $C_0\geq 1$ such that the following holds:
\begin{itemize}
\item For any mesh $\cT\in \bT_*$ there exists a metric $H\in \bH_*$ such that $\cT$ and $H$ are $C_0$-equivalent.
\item For any metric $H\in \bH_*$ there exists a mesh $\cT\in \bT_*$ such that $\cT$ and $H$ are $C_0$-equivalent.
\end{itemize}
\end{definition}
Our main result (Theorem \ref{thEquiv}) states that when the constant $C$ is sufficiently large, and if dimension is $d=2$, the three classes of meshes $\bT_{i,C}\subset \bT_{a,C} \subset \bT_{g,C}$ are respectively equivalent to three classes of metrics
$$
\bH_i \subset \bH_a \subset \bH_g.
$$
that are defined by precise smoothness conditions on the function $z\mapsto H(z)$ in \S \ref{subsecDefHStar}.
The most striking point of this result is that, although the matrix $\cH_T$ does not encode the angles of a simplex $T$, the fact that a mesh $\cT\in \bT_{a,C}$ is \emph{quasi-acute} can be translated in a precise regularity condition on the metric $H$ equivalent to $\cT$.

In order to state our results more precisely we need to introduce some notations. We first recall some notations of linear algebra. We denote by $\M_d$ the collection of $d\times d$ matrices with real entries. We denote by $\GL_d$ the standard linear group, by $\SL_d$ the special linear group and by $\cO_d$ the orthogonal group
\begin{eqnarray*}
\GL_d &:=& \{A\in \M_d\sep \det A\neq 0\}\\
\SL_d &:=& \{A\in \M_d\sep \det A=1\}\\
\cO_d &:=& \{U\in \M_d\sep U^\trans U = \Id\}.
\end{eqnarray*}
We denote by $S_d$ the collection of symmetric matrices, by $S_d^\oplus$ the subset of non-negative symmetric matrices, and by $S_d^+$ the collection of symmetric positive definite matrices.
\begin{eqnarray*}
S_d &:=& \{M\in \M_d\sep M = M^\trans\}\\
S_d^\oplus &:=& \{M\in S_d\sep z^\trans M z \geq 0 \text{ for all } z\in \R^d\}\\
S_d^+ &:=& S_d^\oplus \cap \GL_d.
\end{eqnarray*}
For any $M,M'\in S_d$, we write $M\leq M'$ if $M'-M\in S_d^\oplus$, and $M< M'$ if $M'-M\in S_d^+$. We recall that for any $M,M'\in S_d$ such that $M\leq M'$ and any $\phi\in \M_d$ we have $\phi^\trans M \phi \leq \phi^\trans M' \phi$. Furthermore if $M,M' \in S_d^+$ and $M\leq M'$ then $M'^{-1} \leq M^{-1}$. For each symmetric positive definite matrix $M\in S_d^+$ we define a norm $\|\cdot\|_M$ on $\R^d$ as follows: for all $u\in \R^d$
$$
\|u\|^2_M := u^\trans M u.
$$

\subsection{Geometry of the simplex}
\label{subsecIntroSimplex}
We briefly recall the definition of a $d$-dimensional simplex, and of the collection of its faces.
\begin{definition}
A $d$-dimensional simplex $T$ is the convex envelope of a set $V\subset \R^d$ of $d+1$ 
vertices, not contained in a $d-1$ dimensional affine subspace of $\R^d$ : 
$$
T = \Cvx(V).
$$
We denote by $\cF(T)$ the collection of faces of $T$ of any dimension,
\be
\label{defFT}
\cF(T) := \{\Cvx(V') \sep V'\subset V\}.
\ee
\end{definition}
Observe that $\emptyset$ and $T$ are among the faces of a simplex $T$, since they respectively correspond to $V'=\emptyset$ and $V'=V$.
We denote by $z_T\in \R^d$ the barycenter of a simplex $T$ of vertices $V$
$$
z_T := \frac 1 {d+1} \sum_{v\in V} v,
$$
and we define a symmetric positive definite matrix $\cH_T\in S_d^+$ as follows
\be
\label{defcH}
\cH_T^{-1} = \frac d {d+1} \sum_{v\in V} (v-z_T) (v-z_T)^\trans. 
\ee
The next proposition shows that this definition is consistent with the characterization of $\cH_T$ given in \iref{caractHT}, and gives some of the key properties of the matrix $\cH_T$.
Throughout this chapter, we denote by $\TEq$ a $d$-dimensional equilateral simplex, centered at the origin
i.e. such that $z_\TEq=0$, and having its vertices on the unit sphere. One easily checks that 
$$
\cH_\TEq = \Id.
$$
%
\begin{prop}
\label{propHT}
The following holds for any $d$-dimensional simplex $T$. 
\begin{itemize}
\item
For any an affine change of coordinates $\Phi$, $\Phi(z) := \phi z+z_0$ where $\phi\in \GL_d$ and $z_0$ in $\R^d$, we have 
\be
\label{invcH}
\cH_{\Phi^{-1}(T)} = \phi^\trans \cH_T \phi.
\ee
\item
There exists a rotation $U\in \cO_d$, depending on $T$, such that 
\be
\label{eqcHTEq}
\cH_T^\frac 1 2 (T-z_T) = U(\TEq),
\ee
hence 
\be
\label{eqVolTHT}
|T| \sqrt{\det \cH_T} = |\TEq|.
\ee
\item
We have the inclusions 
\be
\label{eqTET}
\{z\in \R^d \sep \|z-z_T\|_{\cH_T}\leq 1/d\} \subset T \subset \cE_T := \{z\in \R^d \sep \|z-z_T\|_{\cH_T}\leq 1\},
\ee
and the these two ellipsoids are respectively the one of largest volume included in $T$, and of smallest volume containing $T$.
\end{itemize}
\end{prop}
\proof
For any vertex $\Phi^{-1}(v)$ of the simplex $\Phi^{-1}(T)$, we have 
$$
\Phi^{-1}(v)-z_{\Phi^{-1}(T)} = \Phi^{-1}(v)-\Phi^{-1}(z_T) = \phi^{-1}(v-z_T).
$$
Injecting this relation in \iref{defcH} we obtain \iref{invcH}. 
We now turn to the proof of \iref{eqcHTEq} and for that purpose we remark that the simplices $\TEq$ and $T' := \cH_T^\frac 1 2 (T-z_T)$ both have their barycenter at zero, hence there exists a linear map $\phi\in \GL_d$ such that $T' = \phi^{-1}(\TEq)$. We thus obtain, using the formula \iref{invcH} for the transformation of $\cH_T$ under affine change of coordinates 
$$
\begin{array}{ccccc}
\cH_{T'} &=& \cH_T^{-\frac 1 2} \cH_T \cH_T^{-\frac 1 2} &=& \Id \\ 
\cH_{\phi^{-1}(\TEq)} &=& \phi^\trans \cH_\TEq &=& \phi^\trans \phi.
\end{array}
$$
Hence $\phi^\trans \phi= \Id$ which establishes that $\phi\in \cO_d$ and concludes the proof of \iref{eqcHTEq}.
It is well known that the smallest ellipsoid containing $\TEq$ is the unit ball, and that the largest ellipsoid included in $\TEq$ is the ball of radius $1/d$. Combining this fact with the change of coordinates \iref{eqcHTEq} we obtain \iref{eqTET} which concludes the proof of this proposition.
\sq

Note that for any simplex $T$ and any $z,z'\in T$ we obtain using \iref{eqTET}
\be
\label{eqDiamTHT}
\|z-z'\|_{\cH_T} \leq \|z-z_T\|_{\cH_T} + \|z'-z_T\|_{\cH_T} \leq 1+1= 2,
\ee
and 
\be
\label{eqDiamTHT2}
 \frac 2 d \|\cH_T^{-\frac 1 2}\| \leq \diam(T) \leq 2 \|\cH_T^{-\frac 12}\|
\ee
We introduce a measure of degeneracy $\rho(T)\in [1,\infty)$ of a $d$-dimensional simplex $T$
\be
\label{defRhoMM}
\rho(T) := \sqrt{\|\cH_T\| \|\cH_T^{-1}\|}.
\ee
Note that $\rho(T)=1$ if and only if $\cH_T$ is proportional to the identity, which implies in view of \iref{eqcHTEq} that 
$
T = z_T+r U
$
for some $r>0$ and $U\in \cO_d$.

The measure of degeneracy $\rho$ is slightly different from the measure of degeneracy $\frac {\diam(T)^d} {|T|}$ used in Chapters \ref{chapOptAniso} and \ref{chapW1P} but both have the same role: they are minimal for equilateral simplices, and increase as the simplex becomes thinner.\\

The measure of sliverness $S(T)$ of a simplex $T$ is a quantity which plays an important role in the context of optimal mesh adaptation for the finite element approximation of a function in the Sobolev $W^{1,p}$ norm, see   \cite{Ba, Ja} and Chapters \ref{chapW1P} and \ref{chapApproxMet}. 
The measure of sliverness is defined in Chapter \ref{chapW1P} and illustrated by several examples and equivalent quantities. The measure of sliverness is defined by 
\be
\label{defSMM}
S(T) := \inf \{ \|\phi\| \|\phi^{-1}\| \sep \phi \in \GL_d \text{ and } \phi(T) \text{ is acute } \}.
\ee
We recall that a simplex is acute if and only if the exterior normals $\n, \n'$ to any two distinct $d-1$-dimensional faces have a negative scalar product $\<\n,\n'\>\leq 0$.
As observed in Chapter \ref{chapW1P}, the measure of sliverness $S(T)$ can be interpreted as the distance from $T$ to the collection of acute simplices.

The exterior normals $\n,\n'$ to any two distinct $d-1$-dimensional faces of $\TEq$ satisfy $\<\n,\n'\> =-1/d$, hence this simplex is acute. Recalling \iref{eqcHTEq}
 we thus obtain for any simplex $T$
$$
S(T) \leq \rho(T).
$$
Apart from this upper bound, the matrix $\cH_T$ does not contain any direct information on the measure of sliverness $S(T)$.
For any bidimensional triangle $T$ of largest angle $\theta$, Proposition \ref{propSTan} of Chapter \ref{chapW1P} states that :
$$
S(T) = \max\left\{1, \tan \frac \theta 2 \right\}.
$$

\subsection{Simplicial meshes}
\label{subsecDefTStar}

As a starter, we recall the definition of a conforming simplicial mesh of $\R^d$.

\begin{definition}
\label{defSimplMesh}
A conforming simplicial mesh of $\R^d$ is a collection $\cT$ of simplices which satisfy the conformity axiom :  for all $T,T'\in \cT$
\be
\label{confAxiom}
T\cap T' \in \cF(T) \cap \cF(T'),
\ee
as well as the following properties
\begin{itemize}
\item (Covering) The simplices $T\in \cT$ cover the whole infinite domain $\R^d$:
$$
\bigcup_{T\in \cT} T = \R^d.
$$
\item (Partition) The interiors of the simplices $T\in \cT$ are pairwise disjoint: for any $T,T'\in \cT$,
$$
\interior (T) \cap \interior(T') \neq \emptyset\stext{ implies }  T=T'. %
$$
\item (Local finiteness) For any compact set $K \subset \R^d$, only a finite collection of simplices $T\in \cT$ intersect $K$:
$$
\{T \in \cT \sep K \cap T \neq \emptyset \} \ \text{ is finite.}
$$
\end{itemize}
\end{definition}
The conformity axiom \iref{confAxiom} states that the intersection of any two simplices in $\cT$ needs to be a full common face. We denote by $\bT$ the collection of conforming simplicial meshes of $\R^d$.

For an optimal efficiency, numerous numerical simulation software use anisotropic meshes, in which the simplices may have an arbitrary shape. Anisotropic meshes are used for instance, to create a thin layer of simplices close to a geometric feature which is relevant in a numerical simulation.
In order to avoid excessively wild meshes, one often requires some consistency between the shapes of neighboring simplices. We therefore introduce for any constant $C \geq 1$ the class $\bT_{g,C} \subset \bT$ of \emph{graded} meshes, which is defined as follows.

\begin{definition}
A mesh $\cT$ belongs to $\bT_{g,C}$ if for any two simplices $T,T'\in \cT$ one has 
\be
\label{defTG}
T \cap T'\neq \emptyset \;\;  \Rightarrow \;\; C^{-2}\cH_T \leq \cH_{T'} \leq C^2 \cH_T.
\ee
\end{definition}

For any $C\geq 1$ we also introduce the class $\bT_{i,C}$ of \emph{isotropic} meshes, by requiring in addition to \iref{defTG} that the measure of degeneracy $\rho$ is uniformly bounded by $C$, which forbids any anisotropy.

\begin{definition}
A mesh $\cT$ belongs to $\bT_{i,C}$ if and only if $\cT\in \bT_{g,C}$ and for all $T \in \cT$
$$
\rho(T) \leq C. 
$$
\end{definition}

%

In order to introduce the intermediate collection $\bT_{a,C}$ of $C$-quasi-acute meshes, we first define the notion of refinement of a simplicial mesh.
\begin{definition}
Consider two meshes $\cT,\cT'\in \bT$ and a constant $C\geq 1$. We say that $\cT'$ is a $C$-refinement of $\cT$ if it satisfies the following properties.
\begin{itemize}
\item (Inclusion) Any simplex $T'\in \cT'$ is contained in a simplex $T\in \cT$. 
\item (Bounded refinement) Any simplex  $T\in \cT$, contains at most $C$ simplices $T'\in \cT$.
\end{itemize}
\end{definition}

\begin{definition}
A graded mesh $\cT \in \bT_{g,C}$ belongs to the collection $\bT_{a,C}$ of quasi-acute meshes if and only there exists a $C$-refinement $\cT'$ of $\cT$ which satisfies for all $T' \in \cT'$
$$
S(T') \leq C.
$$
\end{definition}

In particular if $\cT\in \bT_{i,C}$, then choosing $\cT' = \cT$, which is a $C$-refinement since $C\geq 1$, and recalling that $S(T) \leq \rho(T)$ for any simplex $T$, we obtain that $\cT\in \bT_{a,C}$. Thus $\bT_{i,C} \subset \bT_{a,C} \subset \bT_{g,C}$ as announced.

The constraints defining graded, quasi-acute and isotropic meshes are illustrated on Figure \ref{figDisk}
in the main introduction of the thesis. 
Figure \ref{figAcuteGen} gives an example of a quasi-acute triangulation which 
by the above described refinement process leads to a slightly finer triangulation on which the measure of 
sliverness $S(T)$ is uniformly bounded.

\begin{figure}
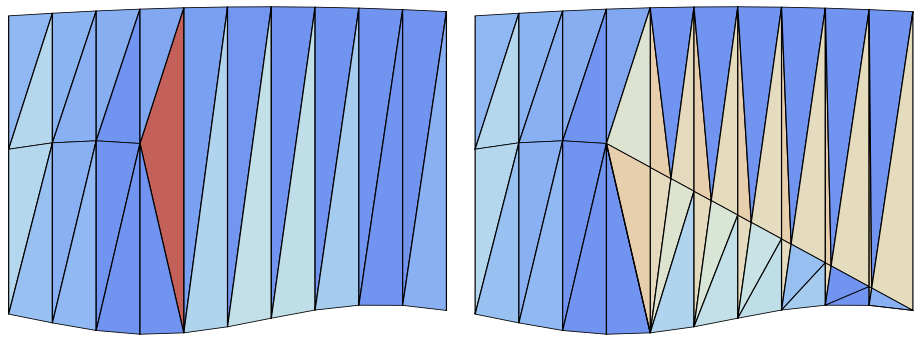

\centering
	\includegraphics[height = 4.5cm, width = 6cm]{\pathPic/Maillages/VarPteVap1.pdf}
	\includegraphics[height = 4.5cm, width = 6cm]{\pathPic/Maillages/VarPteVap2.pdf}
	\caption{\label{figAcuteGen}
Quasi-acute triangulation (left) and refined triangulation with uniformly bounded sliverness (right).
	}
\end{figure}

\subsection{Riemannian metrics}
\label{subsecDefHStar}

For any metric $H \in \bH := C^0(\R^d, S_d^+)$ and any path $\gamma \in C^1([a,b],\R^d)$, we define the riemmannian length $l_H(\gamma)$ by the integral
$$
l_H(\gamma) := \int_a^b \|\gamma'(t)\|_{H(\gamma(t))} dt.
$$
The riemmannian distance $d_H$ between two points  $z,z'\in \R^d$, is the infimum of the length of all paths joining $z$ and $z'$, 
\be
\label{defdH}
d_H(z,z') := \inf\left\{l_H(\gamma)\sep \gamma\in C^1\left([0,1],\R^d\right), \ \gamma(0) = z \text{ and } \gamma(1) = z'\right\}.
\ee
If $H, H'\in \bH$ are such that $H(z) \leq H'(z)$ for all $z\in \R^d$, then $d_H(z,z') \leq d_{H'}(z,z')$ for all $z,z'\in \R^d$.
By construction, the riemannian distance $d_H$ associated to a metric $H\in \bH$ is locally equivalent close to a point $z$ to the distance defined by the norm $\|\cdot \|_{H(z)}$: for all $z\in \R^d$
\be
\label{eqDistLoc}
\lim_{\ve \to 0} \sup_{p,q\in B(z, \ve)} |\ln d_H(p,q) - \ln \|p-q\|_{H(z)}| = 0,
\ee
where $B(z,\ve)$ denotes the euclidean ball of radius $\ve$ around $z$.\\

\begin{remark}[Riemannian geodesics]
When a metric $H$ is sufficiently smooth, a special family of curves named Riemannian geodesics can be defined and studied. A riemannian geodesic is a curve $\gamma\in C^0(\R, \R^d)$ which satisfies a second order differential equation 
called the equation of geodesics. This equation implies that for any sufficiently close $t,t'\in \R$ the curve $\gamma_{|[t,t']}$ is the path of smallest possible length joining the points $\gamma(t)$ and $\gamma(t')$ : 
$d_H(\gamma(t), \gamma(t')) = l_H(\gamma_{|[t,t']})$.

Geodesics are not defined for general continuous or Lipschitz riemannian metrics, such as those constituting the sets $\bH$ and $\bH_i \subset \bH_a\subset \bH_g$, since the coefficients of the equation of geodesics involve the second derivatives of the metric. Notions related to geodesics, such as the injectivity radius of a metric, therefore have no meaning in our context.
This is not an issue for our purposes since, following the point of view of \cite{Gromov},
we study Riemannian distances and not Riemannian geodesics. 


\end{remark}

\begin{definition}
We denote by $\bH_i$ the collection of metrics $H \in \bH$ which are proportionnal to the identity: for all $z\in \R^d$
$$
H(z) = \frac \Id {s(z)^2},
$$
and such that the proportionality factor $s$ satisfies one of the two following conditions conditions, which are surprisingly equivalent as shown in Proposition \ref{propEquivS}:
\begin{itemize}
\item (``Additive'' Lipschitz condition)  
\be
\label{eqLipSAdd}
|s(z) - s(z')| \leq |z-z'| \stext{ for all }z,z'\in\R^d.
\ee
\item (``Multiplicative'' Lipschitz condition)
\be
\label{eqLipSMult}
|\ln s(z) - \ln s(z')| \leq d_H(z,z') \stext{ for all } z,z'\in\R^d.
\ee
\end{itemize}
\end{definition}

The two Lipschitz properties \iref{eqLipSAdd} and \iref{eqLipSMult} have a natural generalisation to general anisotropic metrics but are not equivalent in that context.
In order to introduce the regularity conditions defining the classes $\bH_a$ and $\bH_g$ of metrics, we need to introduce two distances $d_\times$ and $d_+$ on the set $S_d^+$ of symmetric positive definite matrices. 
For any $M,M'\in S_d^+$ we define
\be
\label{defDTimes}
d_\times (M,M') := \min \{\delta \geq 0 \sep e^{-2\delta} M \leq M' \leq e^{2 \delta} M\}.
\ee
For instance, a mesh $\cT$ is $C$-equivalent to a metric $H$, as defined in \iref{defEquiv}, if and only if for all $T \in \cT$ and all $z\in \cT$ one has 
$$
d_\times (H(z), \cH_T) \leq \ln C.
$$
There exists various expressions of the distance $d_\times$ : for all $M,M'\in S_d^+$ one easily checks that  
\be
\label{eqDTimesLog}
\begin{array}{rcl}
d_\times (M,M') &=& \displaystyle\sup_{u\in \sR^d\sm \{0\}} \left|\ln \|u\|_M - \ln \|u\|_{M'}\right|\vspace{2mm}\\
&=&\displaystyle \sup_{|u|=1} \left|\ln \|u\|_M - \ln \|u\|_{M'}\right|,
\end{array}
\ee
and 
\be
\label{eqDTimesProd}
d_\times (M,M') = \ln \left(\max \{\|M^{-\frac 1 2} M'^{\frac 1 2}\|, \, \|M'^{-\frac 1 2} M^\frac 1 2\|\}\right).
\ee
This last expression implies that 
\be
\label{eqDTimesDiffLog}
d_\times (M,M') \geq \frac 1 2 \max\{ |\ln \|M\|-\ln\|M'\| |, \, |\ln \|M^{-1}\| - \ln \|M'^{-1}\| | \}
\ee
The exponential of the distance $d_\times$ was used in the earlier paper \cite{Shew} 
under the name ``relative deformation''. 
We define a second distance $d_+$ on $S_d^+$ as follows 
\be
\label{defDPlus}
d_+(M,M') := \|M^{-\frac 1 2} - M'^{-\frac 1 2}\|.
\ee
The expressions \iref{eqDTimesProd} and \iref{defDPlus} of the distances $d_\times$ and $d_+$ show that they are respectively tied to multiplicative or additive properties of matrices, which justifies their notations.
Note also that for any $s,s'>0$ one has
\be
\label{eqDIso}
d_\times (s^{-2} \Id, s'^{-2} \Id) = |\ln s - \ln s'| \stext{ and } d_+(s^{-2} \Id, s'^{-2} \Id) = |s-s'|.\\
\ee

\begin{definition}
\label{defHg}
We denote by $\bH_g$ the collection of metrics $H\in \bH$ which satisfy the natural extension of \iref{eqLipSMult} to general anisotropic metrics:
\be
\label{eqLipMult}
d_\times (H(z), H(z')) \leq d_H(z,z') \stext{ for all }z,z'\in\R^d.
\ee
\end{definition}

This condition can be heuristically described as follows: a metric $H\in \bH$ defines close to each point $z\in\R^d$ a norm $\|\cdot\|_{H(z)}$ which encodes an anisotropic notion of scale. The equation \iref{eqLipMult} states that the metric $H$ itself is consistent a the scale that it encodes.\\

\begin{definition}
\label{defHa}
We denote by $\bH_a$ the collection of metrics $H\in \bH$ which satisfy both \iref{eqLipMult} and the natural extension of \iref{eqLipSAdd} to general anisotropic metrics:
\be
\label{eqLipAdd}
d_+ (H(z), H(z')) \leq |z-z'| \stext{ for all }z,z'\in\R^d.
\ee
\end{definition}

In addition to the results of mesh generation presented in this chapter, we show in the next chapter, Lemma \ref{lemmaBoundedConv}, that this condition is critical in order to define a local averaging operator compatible with Sobolev norms. 
\begin{remark}
\label{remHomog}
Let $H\in \bH$, and let $\lambda>0$. Using the homogeneity properties of the distances $d_\times$, $d_+$ on $S_d^+$ and $d_H$ on $\R^d$ with respect to $H$, we obtain that $\lambda^2 H\in \bH_g$ if and only if 
\be
\label{eqLipTimesLambda}
d_\times (H(z), H(z')) \leq \lambda d_H(z,z') \stext{ for all } z,z'\in \R^d, 
\ee
and $\lambda^2 H\in \bH_a$ if and only if we have in addition to the previous property 
\be
\label{eqLipPlusLambda}
d_+ (H(z), H(z')) \leq \lambda |z-z'| \stext{ for all } z,z'\in \R^d.
\ee
In particular the collections $\bH_i \subset \bH_a \subset \bH_g$ of isotropic, quasi-acute and graded metrics are stable under multiplication by a constant larger than one : for any $\star \in \{a,b,c\}$, any $H\in \bH_\star$ and any $\lambda \geq 1$ we have $\lambda^2 H\in \bH_\star$.
\end{remark}

\subsection{Main results}

The main result of this chapter is the equivalence of the classes of meshes and metrics defined above, which is  announced in the first part of the introduction and proved in sections \S\ref{secMeshToMet} and  \S \ref{secMetToMesh}. In two cases we were only able to establish this equivalence when the dimension is $d=2$. 
\begin{theorem}
\label{thEquiv}
There exists $C_0 = C_0(d)$ such that for all $C\geq C_0$ the following holds.
\begin{enumerate}[i)]
\item The collections $\bT_{i,C}$ of triangulations and $\bH_i$ of metrics are equivalent.
\item If $d=2$, then the collections $\bT_{a,C}$ of triangulations and $\bH_a$ of metrics are equivalent.
\item If $d=2$, then the collections $\bT_{g,C}$ of triangulations and $\bH_g$ of metrics are equivalent.
\end{enumerate}
\end{theorem}

Let us comment on this result.
The theory of isotropic meshes is already well developed, and our result i) in this direction should be regarded as a reformulation of previous work, intended to put into perspective the results on anisotropic meshes. The main ingredient of the proof of i) is the mesh refinement procedure exposed in \cite{Nochetto}. 

In the case iii) of graded metrics the key ingredients used for the construction of a bidimensional triangulation $\cT\in \bT_{g,C}$ equivalent to a metric $H \in \bH_g$ come from the paper of computational geometry \cite{Shew}. In contrast the arguments used to produce a metric $H \in \bH_g$ from a mesh $\cT \in \bH_{g,C}$ hold in any dimension.

Last the author has not heard of any anterior results on the collection $\bT_{a,C}$ of quasi-acute meshes, and most of the techniques used in that context are new. 

\begin{conjecture}
The author conjectures that the points ii) and iii) of Theorem \ref{thEquiv} hold without restriction on the dimension $d\geq 2$.
\end{conjecture}

This chapter is organised as follows. We establish in section \S \ref{secMetProp} some general properties of the metrics in $\bH_a$ or $\bH_g$. 
We first prove in \S\ref{subsecMetInv} a property of invariance of the collection $\bH_g$ of graded metrics with respect to affine changes of coordinates.
The regularity assumptions \iref{eqLipMult} and \iref{eqLipAdd} defining graded and quasi-acute metrics simply mean that the map $H : \R^d \to S_d^+$ is Lipschitz with respect to some distances on $\R^d$ and $S_d^+$, and 
we therefore study metrics from this point of view in \S\ref{subsecMetLip}.
We focus in \S\ref{subsecMetGeom} on the ``geometrical properties'' of the space $\R^d$ equipped with the distance $d_H$ and the measure $\sqrt{\det H(z)} dz$ associated to a graded metric $H\in \bH_g$. In particular we compare the $d_H$ with the euclidean distance and we estimate the volume of balls.

Section \S\ref{secEigen} is devoted to the study of metrics $H\in \bH$ such that the symmetric matrix $H(z)$ has  an eigenspace of dimension at least $d-1$ at each point. This condition clearly holds if the dimension is $d=2$ but is also relevant in some applications to higher dimension such as the creation of a thin layer of simplices close to a $d-1$ dimensional surface, as discussed in \S\ref{secOptGeom} in the next chapter. We give an approximate translation of the regularity properties \iref{eqLipAdd} and \iref{eqLipMult} defining the collections $\bH_a$ and $\bH_g$ of metrics in terms of the two eigenvalues of $H$ and of the direction of the eigenspaces of such metrics. 

The last two sections \S\ref{secMeshToMet} and \S\ref{secMetToMesh} are devoted to the proof of Theorem \ref{thEquiv}. More precisely we depict in \S\ref{secMeshToMet} to construction of a metric which is equivalent to a given mesh. 
Conversely we address in \S\ref{secMetToMesh} the production of an isotropic, quasi-acute or graded mesh equivalent to a given metric $H\in \bH_i$, $\bH_a$ or $\bH_g$. 

\section{General properties of graded and quasi-acute metrics}
\label{secMetProp}

We study in this section the general properties of the  distinguished classes of metrics $\bH_i \subset \bH_a \subset \bH_g$. 
We first establish some invariance properties of these collections of metrics with respect to affine changes of coordinates.
We then regard metrics as Lipschitz functions $H : \R^d\to S_d^+$. 
The last subsection is devoted to the comparison of the distance $d_H$ and the measure $\sqrt{\det H(z)} dz$ associated to a graded metric $H\in \bH_g$, with respect to the standard euclidean distance and Lebesgue measure.

\subsection{Invariance properties and restriction of metrics}
\label{subsecMetInv}

We focus in this subsection on the properties of invariance of the collections of meshes $\bT_{i,C} \subset \bT_{a,C}\subset \bT_{g,C}$ and of metrics $\bH_i \subset \bH_a\subset \bH_g$ introduced in this chapter.

Let $\Phi : \R^d \to \R^d$ be an affine change of coordinates, $\Phi(z) := \phi z+z_0$ where $\phi\in \GL_d$ and $z_0\in \R^d$.
We have for any simplex $T$, as previously observed in \iref{invcH},
\be
\label{eqInvH}
\cH_{\Phi^{-1}(T)} = \phi^\trans \cH_T \phi.
\ee
Mimicking this property we define for any $H \in \bH$ a transported version $H_\Phi\in \bH$ as follows: for all $z\in \R^d$
\be
\label{defHPhi}
H_\Phi(z) = \phi^\trans H(\Phi(z)) \phi.
\ee
In the language of differential geometry,
the metric $H_\Phi$ is generally referred to as the ``pull back'' of $H$ by $\Phi$. If a mesh $\cT\in \bT$ is $C$-equivalent to a metric $H\in \bH$, then \iref{eqInvH} and \iref{defHPhi} clearly imply that the mesh $\Phi^{-1}(\cT)$ is $C$-equivalent to the metric $H_\Phi$.

The next proposition studies the invariance of the collections $\bT_{\star,C}$ of meshes and $\bT_\star$ of metrics, $\star \in \{i,a,g\}$, under the transformations $\cT \mapsto \Phi^{-1}(\cT)$ or $H \mapsto H_\Phi$ induced by an affine change of coordinates $\Phi$.

\begin{prop}
\label{propInvT}
Let $\Phi : \R^d \to \R^d$ be an affine change of coordinates, $\Phi(z) := \phi z+z_0$ where $\phi\in \GL_d$ and $z_0\in \R^d$. 
The following holds for any constant $C\geq 1$.
\begin{enumerate}[i)]
\item If $\cT \in \bT_{g,C}$ then $\Phi^{-1}(\cT) \in \bT_{g,C}$. If $H \in \bH_g$ then $H_\Phi\in \bH_g$.
\item Assume that $\Phi$ is an isometry, in other words $\phi\in \cO_d$.\\
If $\cT\in \bT_{a,C}$ then $\Phi^{-1}(\cT)\in \bT_{a,C}$, and similarly if $\cT\in \bT_{i,C}$ then $\Phi^{-1}(\cT)\in \bT_{i,C}$.\\
If $H\in \bH_a$ then $H_\Phi\in \bH_a$, and similarly if $H\in \bH_i$ then $H_\Phi\in \bH_i$.
\end{enumerate}
\end{prop}

\proof
We first establish the properties announced for meshes, and we then focus on metrics.
Let $M,M'\in S_d^+$, let $C\geq 1$ and let $\phi\in \GL_d$. Clearly
$$
\text{ if } C^{-2} M \leq M' \leq C^2 M \stext{ then } C^{-2} \phi^\trans M\phi \leq \phi^\trans M'\phi  \leq C^2 \phi^\trans M\phi.
$$
%
In view of the definition \iref{defTG} of $\bT_{g,C}$ we obtain that $\Phi^{-1}(\cT)\in \bT_{g,C}$ for any mesh $\cT\in\bT_{g,C}$ as announced in i).

If $\Phi$ is an isometry then for any simplex $T$ we obtain, using the definitions \iref{defRhoMM} and \iref{defSMM} of the measure of degeneracy $\rho$ and the measure of sliverness $S$, that
$$
\rho(\Phi^{-1}(T)) = \rho(T) \stext{ and } S(\Phi^{-1}(T)) = S(T).
$$
This implies the properties of $\bT_{a,C}$ and $\bT_{i,C}$ announced in b), which concludes the proof of the properties announced for meshes.

We therefore turn our attention to metrics. The fact that $H_\Phi\in \bH_g$ for any $H\in \bH_g$ and any affine change of coordinates $\Phi$ is a special case of Proposition \ref{propHG} below. We thus admit this property and we focus on isotropic and quasi-acute metrics and \emph{isometric} changes of coordinates.

If $H\in \bH_i$, then $H = \Id/s(z)^2$ where $s: \R^d \to \R_+^*$ is a Lipschitz function. Since $\phi\in \cO^d$ we have $\phi^\trans \phi = \Id$ and therefore 
$$
H_\Phi(z) = \phi^\trans H(\Phi(z)) \phi = \Id/s(\Phi(z))^2.
$$
The map $s \circ \Phi : \R_d\to \R_+^*$ is Lipschitz since it is the composition of a Lipschitz function with an isometry, which establishes that $H_\Phi\in \bH_i$ as announced.

Since $\phi\in \cO_d$ we have for any symmetric matrix $M\in S_d^+$ and any exponent $\alpha\in \R$
$$
(\phi^\trans M \phi)^\alpha = \phi M^\alpha \phi.
$$
Consider $H\in \bH_a$, then for any $z,z'\in \R^d$
\begin{eqnarray*}
\|H_\Phi(z)^{-\frac 1 2}-H_\Phi(z')^{-\frac 1 2}\| &=& \|\phi^\trans \left(H(\Phi(z))^{-\frac 1 2}-H(\Phi(z'))^{-\frac 1 2}\right) \phi\| \\
&=& \|H(\Phi(z))^{-\frac 1 2}-H(\Phi(z'))^{-\frac 1 2}\|\\
& \leq& \|\Phi(z)-\Phi(z')\|\\
& =& \|\phi(z-z')\|= |z-z'|,
\end{eqnarray*}
which establishes that $H_\Phi$ satisfies the Lipschitz regularity condition \iref{eqLipAdd}. We already know from the case of graded metrics that the other regularity condition \iref{eqLipMult} is satisfied by $H_\Phi$, which establishes that $H_\Phi\in \bH_a$ and concludes the proof of this proposition.
\sq

In practical applications, one often needs to mesh simultaneously a domain $\Omega$ and a surface $\Gamma$ embedded in $\Omega$. In terms of metric this raises the following question: given a metric $H$ on $\Omega$ which satisfies certain properties of regularity 
what is the regularity of the restriction of this metric to $\Gamma$?

The next proposition answers this question is the simplified setting where $\Omega$ is the infinite domain $\R^d$, $\Gamma$ is an affine subspace of $\R^d$, and properties of regularity are those defining the collection $\bH_g$ of graded metrics. For any $1\leq k \leq d$ we denote by $\bH_g(\R^k)$ the collection of graded metrics on $\R^k$ (with the convention that $\bH_g(\R)$ is the collection $\Lip(\R, \R_+^*)$ of Lipschitz functions from $\R$ to $\R_+^*$).

\begin{prop}
\label{propHG}
Let $1\leq k\leq d$, let $\Phi : \R^k\to \R^d$ be an affine injective map, $\Phi(z) = \phi z+z_0$ where $\phi\in M_{d,k}$ has rank $k$ and $z_0\in \R^d$. Let $H\in \bH$ and let for all $z\in \R^k$
$$
H_\Phi(z) := \phi^\trans H(\Phi(z)) \phi.
$$
If $H\in \bH_g$ then $H_\Phi\in \bH_g(\R^k)$. 
\end{prop}

\proof
For any $u\in \R^d$, one has 
\be
\label{chgNorm}
\|u\|_{H_\Phi(z)} = \|\phi(u)\|_{H(\Phi(z))}.
\ee
Hence for any $z,z'\in \R^k$ we obtain using \iref{eqDTimesLog},
\be
\label{dNHPhi}
\begin{array}{rcl}
d_\times(H_\Phi(z), \ H_\Phi(z'))  &=&\displaystyle \sup_{u\in \sR^k\sm\{0\}} \left|\ln \left\|\phi(u)\right\|_{H(\Phi(z))} -\ln \left\|\phi(u)\right\|_{H(\Phi(z'))}\right| \\
& \leq& d_\times(H(\Phi(z)), \  H(\Phi(z')))
\end{array}
\ee
For any path $\gamma\in C^1([0,1], \R^k)$, it follows from \iref{chgNorm} that 
$$
l_H(\Phi\circ \gamma) = \int_0^1 \|\phi(\gamma'(t))\|_{H(\Phi(\gamma(t)))}dt =  l_{H_\Phi}(\gamma).
$$
Taking the infimum among all paths $\gamma\in C^1([0,1], \R^k)$ joining two points $z,z'\in \R^k$, we obtain since $\Phi\circ \gamma\in C^1([0,1], \R^d)$ is a path joining $\Phi(z)$ and $\Phi(z')$
\be
\label{dHPhi}
d_H(\Phi(z), \Phi(z')) \leq d_{H_\Phi}(z,z').
\ee
Combining \iref{dNHPhi} and \iref{dHPhi}, we conclude the proof of this proposition.
\sq

\subsection{Metrics as Lipschitz functions}
\label{subsecMetLip}

\begin{definition}
Let $(X,d_X)$ and $(Y,d_Y)$ be metric spaces and
let $f \in C^0(X , Y)$ be a continuous function. We define 
the dilatation $\dil(f) \in [0, \infty]$ of $f$ as follows
$$
\dil(f) := \sup_{\substack{x,x'\in X \\ x \neq x'}} \frac{d_Y(f(x),f(x'))}{d_X(x,x')}
$$
\end{definition}

We say that $f$ is {\it  Lipschitz} if and only if $\dil(f) \leq 1$, and more generally that $f$ is 
$\lambda$-Lipschitz if and only if $\dil(f) \leq \lambda$.

We define the {\it local dilatation} $\dil_x(f)$ of a function $f\in C^0(X,Y)$ at a point $x\in X$ as follows
\be
\label{defDilX}
\dil_x(f) := \lim_{\ve\to 0} \dil \left(f_{|B(x,\ve)}\right) = 
\lim_{\ve\to 0} \left(\sup_{p,q\in B(x,\ve)} \frac{d_Y(f(q),f(q))}{d_X(p,q)}\right).
\ee

We specify the distance function on $X$ or $Y$ if it is not the expected ``canonical one'', for instance if $X$ or $Y$ is a subset of $\R^k$ and if the associated distance is not the standard euclidean distance but the riemannian distance $d_H$ associated to a metric, or if $X$ or $Y$ is the collection $S_d^+$ of symmetric positive definite matrices since none of the two distances $d_+$ or $d_\times$ is canonical. The local dilatation is in that case denoted as follows
$$
\dil_x(f\ssep d_X, d_Y), \quad \dil_x(f\ssep d_X)  \stext{ or } \dil_x(f\ssep d_Y).
$$

The local dilatation $\dil_x(f)$ only depends on the local properties of the metric spaces close to $x$ and $f(x)$. Consider the metric space $(\R^d, d_H)$, where $H \in \bH$ is a fixed riemannian metric, and an arbitrary metric space $(Y,d_Y)$. For any $f\in C^0(\R^d, Y)$ and any $z\in \R^d$ we obtain using \iref{eqDistLoc}
\be
\label{eqDilNorm}
\|H(z)^{-1}\|^{-\frac 1 2} \dil_z(f) \leq \dil_z(f\sep d_H) = \dil_z(f \ssep \|\cdot \|_{H(z)}) \leq \sqrt{\|H(z)\|} \dil_z(f)
\ee
We also introduce the lower local dilatation $\dil_z^*(f)\leq \dil_z(f)$ which is defined as follows
$$
\dil_z^*(f) := \lim_{\ve \to 0} \left(\inf_{p,q\in B(x,\ve)} \frac{d_Y(f(q),f(q))}{d_X(p,q)}\right).
$$
The next proposition establishes that local dilatations are sub-multiplicative.
\begin{lemma}
\label{lemmaSubMult}
Let $(X,d_X)$, $(Y,d_Y)$ and $(Z,d_Z)$ be three metric spaces, and let $f\in C^0(X,Y)$ and $g\in C^0(Y,Z)$. 
For any $x\in X$ one has 
\be
\label{eqDilUpper1}
\dil_x(g \circ f) \leq \dil_{f(x)}(g) \, \dil_x (f),
\ee
provided no indeterminate product $0\times \infty$ or $\infty \times 0$ appears in the right hand side. Similarly we have 
\be
\label{eqDilLower1}
\dil_x(g \circ f) \geq \dil^*_{f(x)}(g) \, \dil_x (f),
\ee
provided no indeterminate product $0\times \infty$ or $\infty \times 0$ appears in the right hand side.
\end{lemma}

\proof
We first establish \iref{eqDilUpper1}.
If $\dil_x(f) = \infty$ or $\dil_{f(x)}(g) = \infty$, then there is nothing to prove. We may therefore assume that $\dil_x(f) < \infty$ and $\dil_{f(x)}(g) < \infty$.
For any $\ve>0$ we define 
$$
F(\ve) = \sup_{p,q\in B(x, \ve)} \frac{d_Y(f(p),f(q))}{d_X(p,q)} \stext{ and } G(\ve) := \sup_{p,q\in B(f(x), \ve)} \frac{d_Z(g(p),g(q))}{d_Y(p,q)},
$$
we thus have
$$
\sup_{p,q\in B(x, \ve)} \frac{d_Z(g(f(p)),g(f(q)))}{d_X(p,q)} \leq G(F(\ve) \ve) \, F(\ve).
$$
Taking the limit as $\ve \to 0$ we obtain conclude the proof of \iref{eqDilUpper1}.
For the proof of \iref{eqDilLower1} we consider two sequences $(p_n)_{n \geq 0}$, $(q_n)_{n \geq 0}$ which both tend to $x$ and such that $d_Y(f(p_n),f(q_n))/d_X(p_n, q_n) \to \dil_x(f)$ as $n \to \infty$. We thus have 
$$
\frac{d_Z(g(f(p_n)),g(f(q_n)))}{d_X(p_n,q_n)} = \frac{d_Z(g(f(p_n)),g(f(q_n)))}{d_Y(f(p_n),f(q_n))} \times \frac{d_Y(f(p_n),g(f(q_n)))}{d_X(p_n,q_n)},
$$
and taking the limit as $n \to \infty$ we obtain \iref{eqDilLower1} which concludes the proof of this lemma.
\sq

Let $f: \Omega \to V$ be a $C^1$ function, where $\Omega\subset \R^d$ is an open set and $V$ is a Banach space. Then for any $z\in \Omega$ the local dilatation $\dil_z(f)$, and the lower dilatation $\dil_z^*(f)$ have an explicit expression in terms of the differential $d_z f : \R^d \to V$ of $f$ at $z$
\be
\label{eqDilDiff}
\dil_z(f) = \sup_{|u|=1} \|d_z f(u)\| \stext{ and } \dil_z^*(f) = \inf_{|u|=1} \|d_z f(u)\|
\ee
Applying Lemma \ref{lemmaSubMult} and \iref{eqDilDiff} to the function $g=\ln$ we obtain for any metric space $(X,d_X)$, any $f\in C^0(X,\R_+^*)$ and any $x\in X$
\be
\label{eqDilLog}
\dil_x(\ln f) = \frac {\dil_x(f)} {f(x)} .
\ee

The following proposition shows that, under certain circumstances, the Lipschitz property is a local property.
\begin{prop}
\label{propDilLoc}
Consider the metric space $(\R^d, d_H)$, where $H \in \bH$ is a fixed riemannian metric, and an arbitrary metric space $(Y,d_Y)$.
For any $f\in C^0(\R^d, Y)$ and any path $\gamma\in C^1([0,1], \R^d)$, $\gamma(0) = z$, $\gamma(1) = z'$, one has
\be
\label{eqLipLocPath}
d_Y(f(z), f(z')) \leq l_H(\gamma) \sup_{0\leq t \leq 1} \dil_{\gamma(t)}(f\ssep d_H).
\ee
It follows that 
\be
\label{eqDilLoc}
\dil(f\ssep d_H) = \sup_{z\in \sR^d} \dil_z(f\ssep d_H).
\ee
\end{prop}

\proof
We equip the segment $[0,1]$ with the distance $d_\gamma$ defined for $0 \leq t\leq t'\leq 1$ by 
\be
\label{defDGamma}
d_\gamma(t,t') = d_\gamma(t',t) = l_H(\gamma_{|[t,t']}).
\ee
Note that for any $t\leq t_* \leq t'$ one has 
\be
\label{exactTri}
d_\gamma(t,t') = d_\gamma(t,t_*)+d_\gamma(t_*,t')
\ee
The map $\gamma : ([0,1], d_\gamma) \to (\R^d, d_H)$ is clearly Lipschitz.
It thus follows from Lemma \ref{lemmaSubMult} that $F := f\circ \gamma :[0,1] \to Y$ satisfies for all $t\in [0,1]$
$$
\dil_t(F) \leq \dil_{\gamma(t)}(f) \leq \lambda := \sup_{0\leq t'\leq 1} \dil_{\gamma(t')}(f).
$$
We assume that $\lambda<\infty$, otherwise there is nothing to prove, and we consider a fixed $\delta>0$. For each $t\in [0,1]$ there exists therefore an interval $V_t\subset [0,1]$ containing $t$, relatively open in $[0,1]$, and such that $F_{|V_t}$ is $(\lambda+\delta)$-Lipschitz. Since the segment $[0,1]$ is compact, there exists a finite set $I_0\subset [0,1]$ such that 
\be
\label{finiteCover01}
[0,1] = \bigcup_{t\in I_0} V_t.
\ee
Let $V,V'\subset [0,1]$ be two intersecting intervals, and assume that $F$ is $\lambda+\delta$ Lipschitz on $V$ and on $V'$. For any $t\in V$ and $t'\in V'$, there exists $t_*\in V\cap V'$ which satisfies $t\leq t_*\leq t'$ or $t\geq t_* \geq t'$. Recalling \iref{exactTri} we thus obtain  
\begin{eqnarray*}
d_Y(F(t),F(t')) &\leq& d_Y(F(t), F(t_*))+d_Y(F(t_*), F(t')) \\
&\leq& (\lambda+\delta) d_\gamma(t,t_*)+( \lambda+\delta)d_\gamma(t_*,t') \\
&=& (\lambda+\delta)d_\gamma(t,t').
\end{eqnarray*}
Thus $F$ is $\lambda+\delta$ Lipschitz on the interval $V\cup V'$. Since $[0,1]$ is covered by the finite collection of open intervals \iref{finiteCover01} on which $F$ is $\lambda+\delta$ Lipschitz, we obtain proceeding by induction that $F$ is $(\lambda+\delta)$-Lipshitz on $[0,1]$ and therefore 
$$
d_Y(f(z),f(z')) = d_Y(F(0),F(1))\leq (\lambda+\delta) d_\gamma(0,1) = (\lambda+\delta)l_H(\gamma). 
$$
Letting $\delta\to 0$ we obtain \iref{eqLipLocPath}. Taking the infimum of \iref{eqLipLocPath} among all paths $\gamma\in C^1([0,1], \R^d)$ joining the points $z$ and $z'$ we obtain $d_Y(f(z),f(z')) \leq  \lambda d_H(z,z')$. Therefore $\dil(f\ssep d_H) \leq \lambda$. Conversely we have $\dil_z(f\ssep d_H) \leq \dil(f\ssep d_H)$ for any $z\in \R^d$, which establishes \iref{eqDilLoc} and concludes the proof of this proposition.\sq

The next corollary shows that the global dilatation of a function is controlled by the local dilatation on $\R^d\sm \Gamma$, if $\Gamma$ is a \emph{sufficiently thin} set. For instance the skeleton $\cup_{T\in \cT} \partial T$ of a mesh $\cT\in \bT$.

\begin{corollary}
\label{corolLipNoGamma}
Let $\Gamma\subset \R^d$ be a set which has the following property : for any path $\gamma\in C^1([0,1], \R^d)$ there exists a sequence of paths $(\gamma_n)_{n \geq 0}$ which converges to $\gamma$ in $C^1$ norm and such that  $\{t\in [0,1]\sep \gamma_n(t) \in \Gamma\}$ is finite for each $n$.\\
Let $H\in \bH$, let $(Y,d_Y)$ be an arbitrary metric space, and met let $f\in C^0(\R^d, Y)$.
Then 
$$
\dil(f\ssep d_H) = \sup_{z\in \R^d \sm \Gamma} \dil_z(f\ssep d_H).
$$
\end{corollary}

\proof
We denote $\lambda := \sup_{z\in \R^d \sm \Gamma} \dil_z(f)$.
We consider two points $z,z'\in \R^d$ and a path $\gamma\in C^1([0,1], \R^d)$ such that $\gamma(0) = z$ and $\gamma(1) = z'$. 

We first assume that the set $\{t\in [0,1]\sep \gamma(t) \in \Gamma\}$ is finite, and we denote its elements by $t_1< \cdots < t_k$. We set by convention $t_0 = 0$ and $t_{k+1} = 1$.
For each $0\leq i \leq k$ and each $t,t'$ such that $t_i <t<t'<t_{i+1}$ we have $d_Y(f(\gamma(t)),f(\gamma(t'))) \leq \lambda l_H(\gamma_{|[t,t']})$ according to Proposition \ref{propDilLoc}. Letting $t\to t_i$ and $t'\to t_{i+1}$ we obtain $d_Y(f(\gamma(t_i)),f(\gamma(t_{i+1}))) \leq \lambda l_H(\gamma_{|[t_i,t_{i+1}]})$. Adding up we obtain 
$$
d_Y(f(z), f(z')) \leq \sum_{0\leq i \leq  k} d_Y(f(\gamma(t_i)), f(\gamma(t_{i+1}))) \leq \lambda
\sum_{0\leq i \leq k} l_H(\gamma_{[t_i, t_{i+1}]}) = \lambda l_H(\gamma).
$$

We now consider the case where $\gamma$ meets the set $\Gamma$ an infinite number of times, and we choose a sequence $(\gamma_n)_{n \geq 0}$ as described is the statement of the corollary. We obtain $d_Y(f(\gamma_n(0)), f(\gamma_n(1))) \leq \lambda l_H(\gamma_n)$ for each $n\geq 0$. Taking the limit as $n \to \infty$ we obtain $d_Y(f(z), f(z')) \leq \lambda l_H(\gamma)$ since $\gamma_n$ converges to $\gamma$ is $C^1$ norm as $n \to \infty$.

Taking the infimum among all paths $\gamma\in C^1([0,1], \R^d)$ joining the points $z$ and $z'$, we obtain $d_Y(f(z), f(z')) \leq \lambda d_H(z,z')$ which concludes the proof of this proposition.
\sq

The next proposition establishes that the regularity conditions \iref{eqLipMult} and \iref{eqLipAdd} defining the collections of quasi-acute or graded metrics are equivalent in the case of a metric proportionnal to the identity. 

\begin{prop}
\label{propEquivS}
Let $s\in C^0(\R^d, \R_+^*)$ and let $H := \Id/s^2$. We have for all $z\in \R^d$
\be
\label{eqDilLnH}
\dil_z(\ln s\ssep d_H) = \dil_z(s).
\ee
Furthermore the following properties are equivalent.
\begin{enumerate}
\item The map $s : \R^d \to \R_+^*$ is Lipschitz. (i.e. $H\in \bH_i$) 
\item (``Additive'' Lipschitz property) For all $z,z'\in \R^d$ one has $d_+(H(z),H(z')) \leq |z-z'|$.
\item (``Multiplicative'' Lipschitz property) For all $z,z'\in \R^d$ one has $d_\times(H(z),H(z'))\leq d_H(z,z')$.
\end{enumerate}
\end{prop}

\proof 
Combining \iref{eqDilNorm} and \iref{eqDilLog} we obtain for any $z\in \R^d$
$$
\dil_z(\ln s\ssep d_H) = s(z) \dil_z(\ln s) = \dil_z(s),
$$
which establishes \iref{eqDilLnH}.
According to \iref{eqDIso} we have for any $z,z'\in \R^d$
$$
d_+(H(z),H(z')) = |s(z) - s(z')| \stext{ and } d_\times (H(z),H(z')) = |\ln s(z) - \ln s(z')|.
$$
The properties 1. and 2. are thus equivalent to $\dil(s)\leq 1$, and therefore also to :
\be
\label{eqDilsLoc}
\dil_z(s) \leq 1 \text{ for all } z\in \R^d,
\ee
according Proposition \ref{propDilLoc} applied to the constant metric $H_0 = \Id$ and to the function $f=s$.
On the other hand property 3. is equivalent to $\dil(\ln s \ssep d_H) \leq 1$, and thus to :
$$
\dil_z(\ln s\ssep d_H) \leq 1 \text{ for all } z\in \R^d,
$$
according to Proposition \ref{propDilLoc} applied to the metric $H$ and to the function $f=\ln s$.
This property equivalent to \iref{eqDilsLoc} according to \iref{eqDilLnH}, which concludes the proof of this Proposition.
\sq

The next corollary uses Proposition \ref{propEquivS} to analyse the regularity of the first and last eigenvalue of a metric.
\begin{corollary}
\label{corolLipNorm}
\begin{enumerate}
\item
For any $H\in \bH_g$, the map $z\mapsto \|H(z)\|^{-\frac 1 2}$ is Lipschitz.
\item
For any $H\in \bH_a$, the map $z\mapsto \|H(z)^{-\frac 1 2}\|$ is Lipschitz.
\end{enumerate}
\end{corollary}
\proof
We first establish Point 1.
We consider $H\in \bH_g$, and we define $s(z) := \|H(z)\|^{-\frac 1 2}$ and 
$$
H'(z) := \frac \Id {s(z)^2} = \|H(z)\|\Id,
$$
for each $z\in \R^d$.
Since $H(z) \leq H'(z)$ for all $z\in \R^d$, we have $d_H(z,z') \leq d_{H'} (z,z')$ for all $z, z'\in \R^d$.
It follows that, using \iref{eqDTimesDiffLog},
\begin{eqnarray*}
|\ln s(z)-\ln s(z')| &=& \frac 1 2 \left|\ln \|H(z)\| - \ln \|H(z')\|\right| \\ 
&\leq & d_\times(H(z),H(z')) \\
&\leq& d_H(z,z')\\
&\leq& d_{H'}(z,z').
\end{eqnarray*}
Applying Proposition \ref{propEquivS} to the metric $H'$ we obtain that $H'\in \bH_i$, and therefore that $s$ is Lipschitz as announced.\\

The proof of Point 2. is more straightforward, since for any metric $H\in \bH_a$ we have 
$$
\left|\|H(z)^{-\frac 1 2}\| - \|H(z')^{-\frac 1 2}\| \right| \leq \|H(z)^{-\frac 1 2}-H(z')^{-\frac 1 2}\| \leq |z-z'|
$$
for all $z,z'\in \R^d$,
which concludes the proof of this corollary.\sq

Our last proposition studies the local dilatation of the maximum and the minimum of two functions.
\begin{prop}
\label{propDilMax}
Let $(X,d_X)$ be a metric space and let $\lambda, \mu\in C^0(X, \R)$. Let $\gamma := \min \{\lambda, \mu\}$ and $\Gamma := \max\{\lambda, \mu\}$. Then for all $x\in X$
\be
\label{eqDilMax1}
\max\{\dil_x(\gamma), \dil_x(\Gamma)\} \leq \max\{\dil_x(\lambda), \dil_x(\mu)\}.
\ee
If $(X,d_X) = (\R^d, d_H)$, where $H \in \bH$ is a riemannian metric, then the above inequality is an equality 
\be
\label{eqDilMax2}
\max\{\dil_x(\gamma\ssep d_H),\, \dil_x(\Gamma\ssep d_H)\} = \max\{\dil_x(\lambda\ssep d_H),\, \dil_x(\mu\ssep d_H)\}.
\ee
\end{prop}

\proof
We have for any $p,q\in X$,
\begin{eqnarray*}
|\gamma(p)-\gamma(q)| &=& |\min\{\lambda(p), \mu(p)\} - \min \{\lambda(q), \mu(q)\}|\\
& \leq & \max\{|\lambda(p)-\lambda(q)|,\, |\mu(p)-\mu(q)|\}.
\end{eqnarray*}
For any $x\in X$ and any $\ve>0$ we thus obtain
\begin{eqnarray*}
\dil(\gamma_{|B(x,\ve)}) &=& \sup_{p,q\in B(x,\ve)} \frac{|\gamma(p)-\gamma(q)|}{d_X(p,q)} \\
&\leq& \sup_{p,q\in B(x,\ve)} \frac{\max\{|\lambda(p)-\lambda(q)|,\, |\mu(p)-\mu(q)|\}}{d_X(p,q)}\\
& =& \max\left\{\sup_{p,q\in B(x,\ve)}  \frac{|\lambda(p)-\lambda(q)|}{d_X(p,q)}, \sup_{p,q\in B(x,\ve)}  \frac{|\mu(p)-\mu(q)|}{d_X(p,q)}\right\}\\
&=& \max\left\{\dil(\lambda_{|B(x, \ve)}), \dil(\mu_{|B(x,\ve)})\right\}.
\end{eqnarray*}
Letting $\ve \to 0$ we obtain $\dil_x(\gamma) \leq \max \{\dil_x(\lambda), \dil_x(\mu)\}$. Proceeding likewise for $\Gamma$ we conclude the proof of \iref{eqDilMax1}.\\

We now turn to the proof of \iref{eqDilMax2}. For that purpose we consider a fixed $\ve>0$ and we define
$$
K := \max \{\dil(\gamma_{|B(x,\ve)}\ssep \|\cdot\|_{H(x)}), \ \dil(\Gamma_{|B(x,\ve)} \ssep \|\cdot\|_{H(x)})\},
$$
where $B(x,\ve)$ stands for the euclidean unit ball of radius $\ve$ centered at $x$.
Consider two points $p,q\in B(x,\ve)$. If $\lambda(p) \geq \mu(p)$ and $\lambda(q)\geq \mu(q)$, or if $\lambda(p) \leq \mu(p)$ and $\lambda(q)\leq \mu(q)$, then 
\begin{eqnarray*}
|\lambda(p) -\lambda(q)|&=& \max \{|\gamma(p)-\gamma(q)|, |\Gamma(p)-\Gamma(q)|\}\\
&\leq & K\|p-q\|_{H(x)}
\end{eqnarray*}
Otherwise if  $\lambda(p) \geq \mu(p)$ and $\lambda(q)\leq \mu(q)$, or $\lambda(p) \leq \mu(p)$ and $\lambda(q)\geq \mu(q)$, then there exists a point $r$ on the segment $[p,q]$ such that $\lambda(r)=\mu(r)$. We thus have 
\begin{eqnarray*}
|\lambda(p) -\lambda(q)| &\leq& \max \{|\gamma(p)-\gamma(r)|, |\Gamma(p)-\Gamma(r)|\}\\
& &+ \max \{|\gamma(r)-\gamma(q)|, |\Gamma(r)-\Gamma(q)|\}\\
& \leq & K\|p-r\|_{H(x)} + K\|r-q\|_{H(x)} = K\|p-q\|_{H(x)}
\end{eqnarray*}
If follows that $\dil(\lambda_{|B(x,\ve)}\ssep \|\cdot\|_{H(x)})\leq K$, and letting $\ve\to 0$ we obtain 
$$
\dil_x(\lambda\ssep \|\cdot\|_{H(x)}) \leq \max\{\dil_x(\gamma\ssep \|\cdot\|_{H(x)}), \ \dil_x(\Gamma \ssep \|\cdot \|_{H(x)})\}.
$$ 
Recalling that $\dil_x(f\ssep d_H) = \dil_x(f\ssep \|\cdot\|_{H(x)})$ for any $f\in C^0(\R^d, \R)$ and any $x\in \R^d$, see \iref{eqDilNorm}, we obtain $\dil_x(\lambda\ssep d_H) \leq \max\{\dil_x(\gamma\ssep d_H), \ \dil_x(\Gamma\ssep d_H)\}$. Proceeding likewise for $\mu$ we obtain \iref{eqDilMax2} which concludes the proof of this proposition.
\sq

\subsection{Geometric properties of the space $(\R^d, d_H)$} 
\label{subsecMetGeom}

We focus in this subsection on the ``geometrical properties'' of the space $\R^d$ equipped with the distance $d_H$ and the measure $\sqrt{\det H(z)} dz$ associated to a graded metric $H\in \bH_g$. The next proposition compares the riemannian distance $d_H$ between two points $z$ and $z+u$ with the norm $\|u\|_{H(z)}$ of their difference.


\begin{prop}
\label{propEuclidRiemann}
Let $H \in \bH_g$ and let $z\in \R^d$. For all $u\in \R^d$ one has
\be
\label{eqEuclidRiemann}
\ln(1+\|u\|_{H(z)}) \leq d_H(z,z+u) \leq -\ln(1-\|u\|_{H(z)}),
\ee
where the right hand side equals $\infty$ by convention when $\|u\|_{H(z)}\geq 1$.
\end{prop}
\proof
We consider a path $\gamma\in C^1([0,1], \R^d)$, and we define for each $t\in [0,1]$
\be
\label{defLGamma}
l(t) := \int_0^t \|\gamma'(t)\|_{H(\gamma(t))}dt.
\ee
Since $H\in \bH_g$ we have for any $t\in [0,1]$
$$
\exp(-2\, d_H(z, \gamma(t))) \, H(z) \leq H(\gamma(t)) \leq  \exp(2 \, d_H(z, \gamma(t))) \, H(z).
$$
Therefore, since $d_H(z, \gamma(t)) \leq l(t)$,
$$
\exp(-2l(t)) H(z)  \leq H(\gamma(t)) \leq  \exp(2l(t)) H(z).
$$
It follows that 
$$
\|\gamma'(t)\|_{H(z)} \exp(-l(t)) \leq l'(t) = \|\gamma'(t)\|_{H(\gamma(t))} \leq  \|\gamma'(t)\|_{H(z)} \exp(l(t)), 
$$
hence 
\be
\label{eqGammapLp}
-\frac d {dt} \exp(-l(t)) \leq \|\gamma'(t)\|_{H(z)} \leq \frac d {dt} \exp(l(t)).
\ee
We have 
\be
\label{eqGammaGammap}
\|u\|_{H(z)} =  \|\gamma(1)-\gamma(0)\|_{H(z)} \leq \int_0^1 \|\gamma'(t)\|_{H(z)} dt
\ee
with equality if $\gamma(t) = z+tu$ for all $t\in [0,1]$, i.e. $\gamma$ is a straight line.
For that specific path $\gamma$ we obtain 
$$
\|u\|_{H(z)} = \int_0^1 \|\gamma'(t)\|_{H(z)} dt \geq -\int_0^1 \left(\frac d {dt} \exp(-l(t))\right)dt = 1-\exp(-l_H(\gamma)).
$$
Rearranging the terms we obtain the right part of \iref{eqEuclidRiemann}:
$$
d_H(z,z+u) \leq l_H(\gamma) \leq -\ln(1-\|u\|_{H(z)}).
$$
We now consider an arbitrary path $\gamma$ joining $z$ and $z+u$ and we obtain integrating the right part of \iref{eqGammapLp} and recalling \iref{eqGammaGammap}
$$
\|u\|_{H(z)} \leq \int_0^1 \|\gamma'(t)\|_{H(z)} dt \leq \int_0^1 \left(\frac d {dt} \exp(l(t))\right) dt = \exp(l_H(\gamma)) -1,
$$
which is equivalent to
$
l_H(\gamma) \geq \ln(1+\|u\|_{H(z)}).
$
Taking the infimum among all paths $\gamma\in C^1([0,1], \R^d)$ satisfying $\gamma(0) = z$ and $\gamma(1) = z+u$, we obtain the left part of \iref{eqEuclidRiemann}, which concludes the proof of this proposition.
\sq

We discuss in the rest of this section the consequences of the comparison \iref{eqEuclidRiemann} of the riemannian distance $d_H$ with the distance associated to the norm $\|\cdot\|_{H(z)}$ when $H\in \bH_g$. 
Combining the right part of \iref{eqEuclidRiemann} with the expression \iref{eqDTimesLog} of the distance $d_\times$ we obtain an estimate of the local variations of the norm $\|\cdot \|_{H(z)}$ associated to the metric $H$ at a point $z$: let $H\in \bH_g$, and let $z,u,v\in \R^d$ be such that $\|u\|_{H(z)} < 1$, then
\be
\label{eqLocalNorm}
(1-\|u\|_{H(z)}) \|v\|_{H(z)} \leq \|v\|_{H(z+u)} \leq (1-\|u\|_{H(z)})^{-1} \|v\|_{H(z)}. 
\ee
Note also that 
\be
\label{eqLocalDet}
(1-\|u\|_{H(z)})^d \sqrt{\det H(z)} \leq \sqrt{\det H(z+u)} \leq (1-\|u\|_{H(z)})^{-d} \sqrt{\det H(z)}.
\ee
Another important consequence of \iref{eqEuclidRiemann} is that for any $z\in \R^d$ and any $r\geq 0$ the closed ball $\{z'\in \R^d \sep d_H(z,z')\leq r\}$ for the distance $d_H$, is a closed and bounded subset $\R^d$, hence is compact. 

The following example establishes that \iref{eqEuclidRiemann} is a sharp inequality. The metric $H$ defined by 
$$
H(z) := \frac \Id {(1+|z|)^2}
$$
belongs to $\bH_i$ since $z\mapsto 1+|z|$ is Lipschitz, hence $H$ also belongs to $\bH_a$ and $\bH_g$.
For any $v\in \R^d$, one easily checks that the path of minimal length joining $0$ to $v$ is the straight line, and that
$$
d_H(0, v) = \ln(1+|v|).
$$
Choosing $z=0$ and $u=v$ we obtain that the left part of \iref{eqEuclidRiemann} is sharp
$$
\ln (1+ \|u\|_{H(0)}) = \ln (1+ |u|) = d_H(0,0+u).
$$
Choosing $z=-v$ and $u = v$ we obtain that the right part of \iref{eqEuclidRiemann} is sharp
$$
-\ln\left(1-\|u\|_{H(-u)} \right) = -\ln \left(1-\frac {|u|}{1+|u|}\right) = \ln(1+|u|) = d_H(-u, -u+ u).
$$

The next corollary uses Proposition \ref{propEuclidRiemann} to obtain a lower bound for the mass of balls in the measure $z\mapsto \sqrt{\det H(z)} dz$ associated to a graded metric $H\in \bH_g$. 
For any metric $H\in \bH$, any $z\in \R^d$ and any $r>0$ we define the ellipse
$$
B_H(z,r) := \{ z+u\sep \|u\|_{H(z)} < r\},
$$
and we observe that 
$$
|B_H(z,r)| = \omega r^d (\det H(z))^{-\frac 1 2}
$$
where $\omega$ denotes the volume of the standard euclidean ball of radius one.
\begin{corollary}
\label{corolVol}
There exists $c=c(d)>0$ such that the following holds. For any $H \in \bH_g$, any $r\geq 1$ and any $z\in \R^d$ we have 
\be
\label{eqVolLog}
\int_{B_H(z,r)} \sqrt{\det H} \geq c \ln r. 
\ee
\end{corollary}

\proof
It follows from Proposition \ref{propEuclidRiemann} that $d_H(z_0, z) \geq \ln (r+1)$ for all $z\in \partial B_H(z,r)$. We define the integer 
$$
k = \left\lfloor \ln(r+1) -\frac 1 2\right\rfloor,
$$
and we consider $k$ points $z_1, \cdots , z_k \in E$ such that 
$
d_H(z_0, z_i)= i
$ 
for all $0 \leq i \leq k$.
We define 
$
r_0 := 1-e^{-\frac 1 2} = 0.39\cdots
$
in such way that 
$
-\ln (1-r_0) = 1/2.
$
We have according to \iref{eqLocalDet} for all $z\in B_H(z_i,r_0)$
\be
\label{eqzzi}
\sqrt{\det H(z)} \geq (1-\|z-z_i\|_{H(z_i)})^d  \sqrt{\det H(z_i)} \geq (1-r_0)^d \sqrt{\det H(z_i)}.
\ee
For any $z\in B_H(z_i,r_0)$ we have 
$$
d_H(z_0, z) \leq d_H(z_0, z_i) + d_H(z_i, z) \leq i- \ln (1-r_0) \leq k+\frac 1 2 \leq \ln (r+1),
$$
which implies that $B_H(z_i,r_0) \subset B_H(z,r)$ for all $0 \leq i \leq k$.
Furthermore for any $0\leq i<j\leq k$, $z\in B_H(z_i,r_0)$ and $z'\in B_H(z_j,r_0)$ we have
$$
d_H(z,z') \geq d_H(z_i, z_j) -d_H(z_i, z) - d_H(z_j,z') > |d_H(z_0, z_i) - d(z_0, z_j)| + 2\ln(1 -r_0) =0, 
$$
hence $B_H(z_i,r_0) \cap B_H(z_j,r_0) = \emptyset$.
It follows that 
\begin{eqnarray*}
\int_{B_H(z,r)} \sqrt{\det H} &\geq& \sum_{0 \leq i \leq k} \int_{B_H(z_i,r_0)} \sqrt{\det H}\\
&\geq& \sum_{0 \leq i \leq k} (1-r_0)^d \sqrt{\det H(z_i)} |B_H(z_i,r_0)|\\
&=& \omega r_0^d (1-r_0)^d (k+1)\\
& \geq & c_0 (\ln (r+1) - 1/2),
\end{eqnarray*}
where $c_0 :=  \omega r_0^d (1-r_0)^d$.
Since $r\geq 1$ we have 
$$
\ln(r+1) - \frac 1 2 \geq \left(\frac{\ln(r+1)}{2 \ln 2} -\frac 1 2\right)+  \left(1-\frac 1 {2 \ln 2}\right)\ln(r+1) \geq  \left(1-\frac 1 {2 \ln 2}\right)\ln r
$$
which concludes the proof with $c=(1- 1/(2 \ln 2))c_0$.
\sq

\section{Metrics having an eigenspace of dimension $d-1$ at each point}
\label{secEigen}
We focus in this section on metrics $H\in \bH$ such that the symmetric matrix $H(z)$ has an eigenspace of dimension at least $d-1$ at each point $z\in \R^d$. Note that this condition clearly holds if the dimension is $d=2$, but is also relevant in some applications to higher dimension as illustrated in \S \ref{secOptGeom} of the next chapter.
The main result of this section is Theorem \ref{th2Spaces} which characterises those of these metrics which belong to $\bH_a$ or $\bH_g$ in terms of the regularity of their eigenvalues and of their eigenvectors.

We denote by $\bbS := \{\theta\in \R^d \sep |\theta| = 1\}$ the euclidean unit sphere of $\R^d$, equipped with the distance 
$$
d_\bbS(\theta, \theta') := \arccos(\<\theta, \theta'\>).
$$
We denote by $\cA$ the space of parameters
$$
\cA := \R_+^*\times \R_+^*\times \bbS,
$$
and we define a map $\cS :\cA\to S_d^+$ as follows
\be
\label{defcS}
\cS(\lambda, \mu, \theta) := \lambda \theta \theta^\trans+ \mu(\Id -\theta \theta^\trans).
\ee
The matrix of $\cS$, in any orthonormal basis $\cB$ of $\R^d$ which begins with the vector $\theta$, has the following form
$$
[\cS(\lambda, \mu, \theta)]_\cB = 
\left(
\begin{array}{cccc}
\Lambda &0 & 0 & \cdots\\
0 & \mu & 0 & \cdots\\
0 & 0 & \mu & 0\\
\vdots & \vdots & 0 &\ddots
\end{array}
\right).
$$
We also define for any $a = (\lambda, \mu, \theta)\in \cA$
$$
\cH(a) := \cS(a)^{-2} = \lambda^{-2} \theta \theta^\trans+ \mu^{-2}(\Id-\theta \theta^\trans).
$$

For any metric $H \in \bH$ and any $z\in \R^d$ we denote by $\dil_z(H)_\times$ and $\dil_z(H)_+$ the local dilatations at $z$ associated to the Lipschitz conditions \iref{eqLipMult} and \iref{eqLipAdd} defining $\bH_g$ and $\bH_a$:
\be
\label{defDilH}
\dil_z(H)_\times := \dil_z(H \ssep d_\times, d_H) \stext{ and } \dil_z(H)_+ := \dil_z(H \ssep d_+).
\ee
Let $\Omega\subset \R^d$ be an open set, let $a\in C^0(\Omega, \cA)$, $a(z) = (\lambda(z), \mu(z), \theta(z))$, and let $z\in \Omega$ be such that $\lambda(z) \neq \mu(z)$ or $\dil_z (\theta) < \infty$. We define
\be
\label{eqDZA}
D_z(a)_\times := \max\left\{\dil_z(\ln \lambda\ssep d_H), \ \dil_z(\ln \mu\ssep d_H), \ \frac 1 2 \left| \frac {\lambda(z)} {\mu(z)} - \frac {\mu(z)} {\lambda(z)} \right|\dil_z(\theta\ssep  d_H)\right\}, 
\ee
and 
\be
\label{eqDZAP}
D_z(a)_+ := \max\left\{\dil_z(\lambda), \ \dil_z(\mu),\ |\lambda(z)-\mu(z)| \dil_z(\theta)\right\}. 
\ee
Note that these quantities are not defined if $\lambda(z) = \mu(z)$ and $\dil_z(\theta\ssep  d_H)=\infty$ simultaneously, since an indeterminate product $0\times \infty$ appears in \iref{eqDZA} and \iref{eqDZAP}.
\begin{theorem}
\label{th2Spaces}
Let $H \in \bH$. Assume that there exists an open set $\Omega\subset \R^d$ and a continuous function $a\in C^0(\Omega, \cA)$, $a(z) = (\lambda(z), \mu(z), \theta(z))$, such that $H = \cH\circ a$ on $\Omega$.
For each $z\in \Omega$ such that $\lambda(z) \neq \mu(z)$ or $\dil_z (\theta) < \infty$ one has  
\begin{eqnarray}
\label{eqDilTimes}
D_z(a)_\times &\leq \dil_z(H)_\times \leq& 2D_z(a)_\times\\
\label{eqDilPlus}
D_z(a)_+ &\leq \dil_z(H)_+ \leq& 2D_z(a)_+.
\end{eqnarray}
\end{theorem}

The rest of this section is devoted to the proof of this theorem and to Corollaries \ref{corolDilThetaG} and  \ref{corolHgToHa}. Our first intermediate result in the proof of Theorem \ref{th2Spaces} defines and estimates a quantity $\Delta(a,b,c)$ which appears repeatedly in the rest of the proof.

\begin{lemma}
\label{lemmaDelta}
For each $a,b,c\in \R$ we define 
\be
\label{defDelta}
\Delta(a,b,c) := 
\left\| 
\left(
\begin{array}{cccc}
a & b & 0 & \cdots\\
b & c& 0 & \cdots\\
0 & 0 & c & 0\\
\vdots & \vdots & 0 &\ddots
\end{array}
\right)
\right\|.
\ee
Then 
\be
\label{eqEstimM}
\max\{|a|, |b|, |c|\} \leq \Delta(a,b,c) \leq 2 \max\{|a|, |b|, |c|\}.
\ee
\end{lemma}

\proof
We consider fixed values of $a,b,c\in \R$ and we denote by $M$ the matrix appearing in \iref{defDelta}, in such way that $\Delta(a,b,c) = \|M\|$. We also define $\lambda := \max\{|a|,|b|,|c|\}$.
We have 
$$
\Delta(a,b,c) = \|M\|\geq \max_{1\leq i,j\leq d} |M_{ij}| = \lambda,
$$
which establishes the left part of \iref{eqEstimM}. For all $p\in [1,\infty]$ we denote by $l^p$ the usual norm on $\R^d$ of exponent $p$. We have 
$$
\|M\|_{l^\infty\to l^\infty} = \max_{1 \leq i \leq d} \sum_{1\leq j \leq d} |M_{ij}| \leq 2 \lambda \stext{ and } \|M\|_{l^1\to l^1} = \max_{1 \leq j \leq d} \sum_{1\leq 1 \leq d} |M_{ij}| \leq 2 \lambda.
$$
By interpolation we thus obtain  
$$
\|M\| = \|M\|_{l^2\to l^2} \leq \sqrt{\|M\|_{l^1\to l^1} \|M\|_{l^\infty\to l^\infty} } \leq 2 \lambda,
$$
which concludes the proof.
\sq

The parameter space $\cA = \bbS \times \R_+^* \times \R_+^*$ is a differential variety, and the regularity of functions defined on $\cA$ should be seen through local charts.
For each $a = (\lambda, \mu, \theta)\in \cA$ we introduce a local chart $\psi_a$ of a neigborhood of $a$ in $\cA$.
\be
\label{defPsia}
\psi_a : (-\lambda, \infty)\times (-\mu , \infty) \times (B(0,\cPi) \cap \theta^\perp) \to \cA,
\ee
where $\theta^\perp := \{\Theta\in \R^d \sep \<\theta, \Theta\> = 0 \}$ denotes the space orthogonal to $\theta$, and $B(0,\cPi)$ the open euclidean ball of radius $\cPi$.
We define 
$$
\psi_a(\Lambda, M, r\Theta) := (\lambda+ \Lambda,\, \mu+ M,\, \cos(r) \theta + \sin(r) \Theta),
$$
where $|\Theta| = 1$.
Note that $\psi_a$ is one to one and that $d_\bbS(\theta, \, \cos(r) \theta + \sin(r) \Theta) = r$.

Let $V$ be a banach space. We say that a function $\vp : \cA \to V$ is $C^1$ if for any $a\in \cA$ the function $\vp\circ \psi_a$ is $C^1$.
In that case for any $a\in \cA$ we define the differential 
$$
d_a\vp(A) : \R\times \R\times \theta^\perp \to V
$$
by the formula
$$
d_a\vp(A) := d_0(\vp\circ \psi_a)(A) = \lim_{t\to 0} \frac{\vp\circ \psi_a( tA) - \vp(a)} t. 
$$

\begin{prop}
\label{propPhiPlus}
Define $V_+ := S_d$, equipped with the standard norm $\|\cdot\|$, and $\vp_+ := \cS$ on $\cA$.
Then for all $a,b\in \cA$ we have 
\be
\label{eqDPPhiP}
d_+(\cH(a), \cH(b)) = \|\vp_+(a)-\vp_+(b)\|.
\ee
Furthermore $\vp_+ \in C^1(\cA, V_+)$ and for all $a=(\lambda, \mu, \theta)\in \cA$ and all $A = (\Lambda, M, r\Theta)\in T_a\cA$, with $|\Theta| = 1$, we have 
$$
\|d_a \vp_+(A)\|_\infty= 
\Delta(\Lambda, \ (\lambda-\mu) r ,\ M).
$$
\end{prop}

\proof
The identity \iref{eqDPPhiP} directly follows from the fact that $\cH(a)^{-\frac 1 2} = \cS(a) = \vp_+(a)$, and from the definition \iref{defDPlus} of $d_+$.
The function $\vp_+ = \cS$ is $C^1$ (in fact  $C^\infty$), since it has a polynomial expression 
$$
\cS(\lambda, \mu, \theta) := \lambda \theta \theta^\trans+ \mu(1-\theta \theta^\trans).
$$
We consider a fixed point $a = (\lambda,\mu,\theta)\in \cA$, and $A=(\Lambda, M, r\Theta) \in T_a\cA$, where $|\Theta|=1$. 
We thus have for $t\in \R$
\begin{eqnarray*}
\cS(\psi_a(tA))  &=& 
 (\lambda+ t\Lambda) (\theta+ tr\Theta) (\theta+tr\Theta)^\trans \\
 & &+ (\mu+tM)(\Id - (\theta+ tr\Theta) (\theta+tr\Theta)^\trans) +\cO(t^2)\\
 &=& S(\lambda, \mu, \theta) + t(\Lambda \theta \theta^\trans + M(\Id - \theta \theta^\trans) + r(\lambda-\mu) (\theta \Theta^\trans + \Theta \theta^\trans))+\cO(t^2)
\end{eqnarray*}
Therefore 
\be
\label{eqDiffS}
d_a\cS(A) =  \Lambda \theta \theta^\trans + M(\Id - \theta \theta^\trans) + r(\lambda-\mu) (\theta \Theta^\trans + \Theta \theta^\trans).
\ee
Choosing an orthonormal basis $\cB$ of $\R^d$ which begins with the unit orthogonal vectors $\theta$ and $\Theta$, we obtain
$$
[d_a\cS(A)]_\cB = 
\left(
\begin{array}{cccc}
\Lambda & r(\lambda-\mu) & 0 & \cdots\\
r(\lambda-\mu) & M & 0 & \cdots\\
0 & 0 & M & 0\\
\vdots & \vdots & 0 &\ddots
\end{array}
\right),
$$
which concludes the proof of this proposition.
\sq

\begin{prop}
\label{propPhiTimes}
Define $V_\times := C^0(\bbS, \R)$, equipped with the $\|\cdot\|_\infty$ norm, and define $\vp_\times : \cA \to V_\times$ as follows: for any $a\in \cA$ 
$$
\vp_\times (a) := (u \mapsto \ln \|u\|_{\cH(a)}).
$$
Then for all $a,b\in \cA$ we have 
\be
\label{eqDTPhiT}
\|\vp_\times(a)-\vp_\times(b)\|_\infty = d_\times(\cH(a), \cH(b)).
\ee
Furthermore $\vp_\times \in C^1(\cA, V_\times)$ and for all $a=(\lambda, \mu, \theta)\in \cA$ and all $A = (\Lambda, M, r\Theta)\in T_a\cA$, with $|\Theta| = 1$, we have 
\be
\label{eqNormDPhiPlus}
\|d_a \vp_\times(A)\|_\infty= 
\Delta\left(\frac \Lambda \lambda ,\ \frac r 2\left(\frac \lambda \mu - \frac \mu \lambda\right) ,\ \frac M \mu\right).
\ee
\end{prop}

\proof
The expression \iref{eqDTPhiT} immediately follows from the expression \iref{eqDTimesLog} of the distance $d_\times$.
It is well known that the inverse map ${\rm Inv} : \GL_d \to \GL_d$ is $C^1$ and has the the following differential : for any $\phi\in \GL_d$ and any $\Phi\in M_d$
$$
d_\phi \Inv (\Phi) = \phi^{-1} \Phi \phi^{-1}.
$$
Since $\cS : \cA\to S_d^+$ is $C^1$, as observed in the proof of Proposition \ref{propPhiPlus}, the composition $\Inv \circ\, \cS : \cA \to S_d^+$ is also $C^1$, and for any $a=(\lambda, \mu , \theta)\in \cA$ and any $A\in \R\times \R \times \theta^\perp$
$$
d_a(\Inv \circ \, \cS)(A) = S(a)^{-1} (d_a S(A) ) S(a)^{-1}.
$$
For any $u\in \bbS$ we have 
$$
\ln \|u\|_{\cH(a)} = \ln \|\cS(a)^{-1}u\| = \frac 1 2 \ln \left(\<\cS(a)^{-1}u,\cS(a)^{-1}u\>\right).
$$
Therefore, again by composition, 
\begin{eqnarray*}
d_a (\ln \|u\|_{\cH})(A) &=&  \frac 1 {2\|\cS(a)^{-1}u\|^2} (\<\cS(a)^{-1} (d_a\cS(A))\cS^{-1}(a)u,\ \cS^{-1}(a)u\>\\
& &\quad +\<\cS^{-1}(a)u,\ \cS(a)^{-1} (d_a\cS(A))\cS^{-1}(a)u\>)\\
&=& \frac{\<G(a,A) v,v\>}{\|v\|^2}
\end{eqnarray*}
where $v = \cS^{-1}(a)u$. The function $G$ has the following expression : for any $a = (\lambda, \mu, \theta)\in \R_+^* \times \R_+^* \times \bbS$ and any $A = (\Lambda, M, r\Theta)\in \R\times \R \times \theta^\perp$, where $|\Theta| = 1$,
\begin{eqnarray*}
G(a,A) &=& \frac{\cS(a)^{-1} (d_a\cS(A)) + (d_a\cS(A)) \cS(a)^{-1}} 2\\
&=& \frac \Lambda \lambda \theta \theta^\trans +  \frac M \mu (\Id - \theta \theta^\trans) + \frac r 2\left(\frac \lambda \mu - \frac \mu \lambda\right)(\theta \Theta^\trans + \Theta \theta^\trans),
\end{eqnarray*}
where we used the explicit expression \iref{eqDiffS} of $d_a\cS(A)$, and the fact that $\theta^\trans \Theta = 0$.
Choosing an orthonormal basis $\cB$ of $\R^d$ which begins with the unit orthogonal vectors $\theta$ and $\Theta$, we obtain
$$
[G(a,A)]_\cB = 
\left(
\begin{array}{cccc}
\frac \Lambda \lambda & \frac r 2\left(\frac \lambda \mu - \frac \mu \lambda\right) & 0 & \cdots\\
\frac r 2\left(\frac \lambda \mu - \frac \mu \lambda\right) & \frac M \mu & 0 & \cdots\\
0 & 0 & \frac M \mu & 0\\
\vdots & \vdots & 0 &\ddots
\end{array}
\right).
$$
The map $\vp_\times$ is differentiable, as the composition of differentiable maps, and $d_a(\vp_\times(A))$ is the element of $V_\times  := C^0(\bbS, \R)$ defined by 
$
u \mapsto d_a (\ln \|u\|_\cH)(A).
$
Therefore 
\begin{eqnarray*}
\|d_a \vp_\times(A)\|_\infty &=& \sup \left\{ \frac{|\<G(a,A) v,v\>|}{\|v\|^2} \sep u\in \bbS, \, v=\cS(a)^{-1} u\right\}\\
&=& \|G(a,A)\|,
\end{eqnarray*}
which establishes \iref{eqNormDPhiPlus} and concludes the proof of this proposition.
\sq

We now consider a fixed $z\in \Omega$ and we remark that 
$$
\dil_z(H)_\times = \dil_z(\cH\circ a)_\times  = \dil_z(\vp_\times \circ a\ssep d_H). 
$$
Therefore,  according to Lemma \ref{lemmaSubMult} for any norm $\|\cdot\|_{a(z)}$ defined on the tangent space $\R\times \R \times \theta(z)^\perp$ to $a(z)\in \cA$ we have 
\be
\label{eqDilPhiA}
\begin{array}{rcl}
& &\dil^*_0 (\vp_+\circ \psi_{a(z)} \ssep \|\cdot\|_{a(z)}) \dil_z(\psi_{a(z)}^{-1} a\ssep d_H, \|\cdot \|_{a(z)})\\
& \leq& \dil_z(H)_\times\\
& \leq& \dil_0 (\vp_+\circ \psi_{a(z)} \ssep \|\cdot\|_{a(z)}) \dil_z(\psi_{a(z)}^{-1} \circ a\ssep d_H, \|\cdot\|_{a(z)})
\end{array}
\ee
provided no indeterminate $0\times \infty$ or $\infty \times 0$ appears on the first and last term of this inequality.

We first consider the case where $\lambda(z) \neq \mu(z)$.
For each $a'=(\lambda', \mu', \theta')\in \cA$ we define a norm $\|\cdot \|_{a'}$ on the tangent space $\R\times \R \times \theta'^\perp$ to $\cA$ at the point $a'$ : for all $A = (\Lambda, M , \Theta)\in \R\times \R \times \theta'^\perp$
$$
\|A\|_{a'} = \max \left\{\frac {|\Lambda|} {\lambda'} , \, \frac {|M|} {\mu'}, \, \frac 1 2 \left| \frac {\lambda'} {\mu'} - \frac {\mu'} {\lambda'}\right| \|\Theta\|\right\}.
$$
We have by construction 
\be
\label{eqDilADZDH}
\dil_z(\psi_{a(z)}^{-1}\circ a\ssep d_H, \|\cdot \|_{a(z)}) = D_z(a)_\times
\ee
According to Lemma \ref{lemmaDelta} and Proposition \ref{propPhiTimes} we have for any $a'=(\lambda', \mu', \theta')\in \cA$ and any $A\in \R\times\R \times \theta'^\perp$
$$
\|A\|_{a'} \leq \|d_{a'} \vp_\times(A)\|_\infty = \|d_0 (\vp_\times \circ \psi_{a'})\|_\infty \leq 2 \|A\|_{a'}.
$$
Since $\vp_\times$ is $C^1$ it follows from \iref{eqDilDiff} that 
$$
1 \leq \dil^*_0 (\vp_\times\circ \psi_{a'}\ssep \|\cdot \|_{a'}) \leq \dil_0 (\vp_\times\circ \psi_{a'}\ssep \|\cdot \|_{a'}) \leq 2.
$$
Combining \iref{eqDilPhiA}, \iref{eqDilADZDH} and the last inequality we obtain
$$
D_z(a)_\times \leq \dil_z(H)_\times \leq 2 D_z(a)_\times,
$$
which establishes the announced result \iref{eqDilTimes} in the case $\lambda(z)\neq \mu(z)$. A similar reasoning establishes the counterpart \iref{eqDilPlus} of this inequality for the distance $d_+$.\\

We now consider the case $\lambda(z) = \mu(z)$.
For any $\ve>0$ and any $a'=(\lambda', \mu', \theta')\in \cA$ such that $\lambda' = \mu'$ we define a norm $\|\cdot\|_{a',\ve}$ on the tangent space $\R\times \R\times \theta'^\perp$ to $a'$ in $\cA$ as follows : 
$$
\|A\|_{a', \ve} = \max \left\{\frac {|\Lambda|} {\lambda'} , \, \frac {|M|} {\mu'}, \, \ve \|\Theta\|\right\}.
$$
This modification is required because the original norm $\|\cdot\|_{a'}$ used in the case $\lambda'\neq \mu'$ is only a semi-norm when $\lambda' = \mu'$.
Reasoning similarly to the case $\lambda(z) \neq \mu(z)$ we obtain the upper bound
$$
\dil_z(H)_\times \leq 2 \max\left\{\dil_z(\ln \lambda\ssep d_H), \ \dil_z(\ln \mu\ssep d_H), \ \ve\dil_z(\theta\ssep  d_H)\right\}.
$$
Since $\lambda(z) = \mu(z)$, the assumptions of the theorem state that $\dil_z(\theta \ssep d_H) < \infty$. Since $\ve>0$ is arbitrary we thus obtain 
$$
\dil_z(H)_\times \leq 2 \max\left\{\dil_z(\ln \lambda\ssep d_H), \ \dil_z(\ln \mu\ssep d_H)\right\} = 2D_z(a)_\times.
$$
We now turn to the lower bound on $\dil_z(H)_\times$. 
For that purpose we define the functions $\gamma := \min \{\lambda, \mu\}$ and $\Gamma := \max\{\lambda,\mu\}$. For all $z,z'\in \Omega$ we have according to \iref{eqDTimesDiffLog} 
$$
d_\times (H(z),H(z')) \geq \max \{ |\ln\gamma(z) - \ln \gamma(z')|,\, |\ln \Gamma(z) - \ln \Gamma(z')|\}.
$$
Therefore 
$$
\dil_z(H)_\times = \dil_z(H\ssep d_\times, d_H) \geq \max \{ \dil_z(\ln \gamma\ssep d_H), \, \dil_z(\ln \Gamma \ssep d_H)\}.
$$
Hence according to Proposition \ref{propDilMax}
$$
\dil_z(H)_\times \geq \max \{\dil_z(\ln \lambda\ssep d_H), \dil_z(\ln \mu\ssep d_H)\} = D_z(a)_\times.
$$
We have thus obtained $D_z(a)_\times \leq \dil_z(H)_\times \leq 2 D_z(a)_\times$. Proceeding likewise we obtain $D_z(a)_+ \leq \dil_z(H)_+ \leq 2 D_z(a)_+$, which concludes the proof of Theorem \ref{th2Spaces}. \\

Our first corollary compares $D_z(a)_\times$ with the dilatation $\dil_z(\theta)$. 

\begin{corollary}
\label{corolDilThetaG}
Under the hypotheses of Theorem \ref{th2Spaces} we have
$$ 
|\lambda(z) - \mu(z)| \dil_z(\theta) \leq 2 D_z(a)_\times 
$$
\end{corollary}
\proof
We have according to \iref{eqDilNorm} 
\begin{eqnarray*}
 |\lambda(z) - \mu(z)| \dil_z(\theta)
&\leq& \frac {|\lambda(z) - \mu(z)|}{\min\{\lambda(z), \mu(z)\}} \dil_z(\theta\ssep d_H)\\
&=& \left(\frac {\max \{\lambda(z), \mu(z)\}}{\min\{\lambda(z), \mu(z)\}}-1\right)\dil_z(\theta\ssep d_H)\\
& \leq & \left|\frac {\lambda(z)} {\mu(z)}-\frac {\mu(z)} {\lambda(z)}\right|\dil_z(\theta\ssep d_H) \\
&\leq & 2 D_z(a)_\times
\end{eqnarray*}
which concludes the proof.
\sq

The next corollary shows how to construct a quasi-acute metric from a graded one.
Let $M\in S_d$ and let $U\in \cO_d$ and $\lambda_1, \cdots , \lambda_d\in\R$ be such that 
\be
\label{eqDiagSd}
M = U^\trans \diag(\lambda_1, \cdots , \lambda_d) U,
\ee
where $\diag(\lambda_1, \cdots, \lambda_d)$ denotes the diagonal matrix of entries $\lambda_1, \cdots, \lambda_d$. 
For any $\alpha\in \R$ we define the matrix $\max \{\alpha, M\}\in S_d$ as follows
\be
\label{defMaxAlphaM}
\max\{\alpha, M\} := U^\trans \diag(\max\{\alpha, \lambda_1\}, \cdots, \max\{\alpha, \lambda_d\})U,
\ee
and we observe that $\max\{\alpha, M\}$ does not depend on the choice of the matrix $U$ and on the order of the eigenvalues $\lambda_1, \cdots, \lambda_d$ in the decomposition \iref{eqDiagSd} of $M$.
For any $(\lambda, \mu, \theta)\in \cA$ and any $s>0$ one has
\be
\label{maxsHa}
\max \{s^{-2}, \cH(a)\} = \min\{s, \lambda\}^{-2} \theta \theta^\trans+ \min\{s, \mu\}^{-2} (\Id-\theta \theta^\trans).
\ee

\begin{corollary}
\label{corolHgToHa}
Let $H\in \bH$ be such that $H(z)$ has an eigenvalue of multiplicity at least $d-1$ for all $z\in \R^d$. We define $s:\R^d\to \R_+^*$ and $H'\in \bH$ as follows : for all $z\in \R^d$
\begin{eqnarray*}
s(z) &:=& \inf_{z'\in \sR^d} |z-z'|+ \|H(z')^{-\frac 1 2}\|,\\
H'(z) &:=& \max\{s(z)^{-2}, H(z)\}.
\end{eqnarray*}
If $H\in \bH_g$, then $4^2H'\in \bH_a$.
\end{corollary}

\proof
We first observe that $s : \R^d \to \R_+^*$ is Lipschitz and that $s(z)\leq \|H(z)^{-\frac 1 2}\|$ for all $z\in \R^d$.
For each $\ve>0$ and each $z\in \R^d$ we define 
$$
H_\ve(z) := \max\{e^{2 \ve} s(z)^{-2}, H(z)\}.
$$

Consider a fixed point $z\in\R^d$. If $H(z)$ is proportional to $\Id$, which means that $\|H(z)\|^{-\frac 1 2} = \|H(z)^{-\frac 1 2}\| \geq s(z)$, then $H_\ve = e^{2 \ve} \Id/ s^2$ on a neighborhood of $z$. This implies according to \iref{eqDilLnH} 
$$
\dil_z(H_\ve)_\times =\dil_z(H_\ve)_+ = \dil_z(e^{-\ve} s) \leq e^{-\ve}\leq 1.
$$
On the other hand if $H(z)$ is not proportional to $\Id$ then there exists $a\in C^0(\R^d, \cA)$, $a(z) = (\lambda(z), \mu(z), \theta(z))$ such that $H = \cH\circ a$ on a neighborhood of $z$, and $D_z(a)_\times \leq 1$ according to Theorem \ref{th2Spaces}. 
We have according to \iref{maxsHa} on the same neighborhood of $z$, 
$$
H_\ve = \cH\circ a_\ve = \min\{e^{-\ve} s,\lambda\}^{-2} \theta \theta^\trans+ \min\{e^{-\ve} s,\mu\}^{-2} (\Id - \theta \theta^\trans),
$$
where $a_\ve=(\min\{e^{-\ve}s, \lambda\}, \min\{e^{-\ve} s, \mu\}, \theta) \in C^0(\R^d,\cA)$.

We define $H_\ve^i := e^{2\ve} \Id/ s^2$ and we observe that $H_\ve \geq H_\ve^i$ and $H_\ve \geq H$, hence for any $p,q\in \R^d$
$$
d_{H_\ve}(p,q) \geq \max \{d_H(p,q), \, d_{H_\ve^i}(p,q)\}.
$$
It follows from Proposition \ref{propDilLoc} that   
\begin{eqnarray*}
\dil_z(\ln \min\{e^{-\ve} s, \lambda\} \ssep d_{H_\ve}) &\leq& \max \{\dil_z(\ln (e^{-\ve} s)\ssep d_{H_\ve}), \,\dil_z(\ln \lambda\ssep d_{H_\ve})\}\\
&\leq & \max \{\dil_z(\ln(e^{-\ve} s) \ssep d_{H_\ve^i}), \,\dil_z(\ln \lambda\ssep d_H)\}\\
& \leq & \max \{\dil_z(e^{-\ve} s), \, \dil_z(H)_\times\}\\
& \leq & \max \{e^{-\ve}, 1\} =1,
\end{eqnarray*}
where we used Proposition \ref{propEquivS} and Theorem \ref{th2Spaces} in the third line.
Proceeding likewise for the other eigenvalue of $H_\ve$ we conclude that 
\be
\label{eqDilAeVal}
\dil_z(\ln \min\{e^{-\ve} s, \lambda\}\ssep d_{H_\ve}) \leq 1 \stext{ and }
\dil_z(\ln \min\{e^{-\ve} s, \mu\}\ssep d_{H_\ve}) \leq 1.
\ee
One easily checks that
$$
\left|\frac {\min\{e^{-\ve} s, \lambda\}}{\min\{e^{-\ve} s, \mu\}} - \frac {\min\{e^{-\ve} s, \mu\}}{\min\{e^{-\ve} s, \lambda\}}\right| \leq \left|\frac \lambda \mu - \frac \mu \lambda\right|,
$$
using that 
$
t \mapsto t/\mu - \mu/t 
$
is increasing on $\R_+^*$.
Therefore   
$$
\frac 1 2 \left|\frac {\min\{e^{-\ve} s, \lambda\}}{\min\{e^{-\ve} s, \mu\}} - \frac {\min\{e^{-\ve} s, \mu\}}{\min\{e^{-\ve} s, \lambda\}}\right|  \dil_z(\theta \ssep d_{H_\ve}) \leq \frac 1 2\left|\frac \lambda \mu - \frac \mu \lambda\right| \dil_z(\theta\ssep d_H) \leq 1.
$$
Combining the last estimate with \iref{eqDilAeVal} we conclude that $D_z(a_\ve)_\times \leq 1$ and therefore $\dil_z(H_\ve)_\times \leq 2$ according to Theorem \ref{th2Spaces}.

Our next objective is to estimate $\dil_z(H_\ve)_+$.
Let us assume that $\lambda(z)<\mu(z)$, then
\begin{eqnarray*}
\dil_z(\min\{e^{-\ve} s , \lambda\}) &\leq& \max \left\{\dil_z(e^{-\ve} s), \, \dil_z(\lambda)\right\}\\
& \leq & \max \left\{\dil_z(e^{-\ve} s), \, \frac{\dil_z(\lambda \ssep d_H)}{\min\{\lambda(z), \mu(z)\}}\right\}\\
& = & \max \left\{\dil_z(e^{-\ve} s), \, \frac{\dil_z(\lambda \ssep d_H)}{\lambda(z)}\right\}\\
&\leq& \max \{e^{-\ve}, \, \dil_z(\ln \lambda\ssep  d_H)\}\\
&\leq & \max \{e^{-\ve} , \, 1\} = 1,
\end{eqnarray*}
where we used \iref{eqDilNorm} in the second line.
Still assuming that $\lambda(z)<\mu(z)$, we have $s(z)\leq \|H(z)^{-\frac 1 2}\| = \mu(z)$ by construction of $s$, which implies that $\min\{e^{-\ve} s, \mu\} = e^{-\ve} s$ on a neighborhood of $z$, and therefore
$
\dil_z(\min\{e^{-\ve} s, \mu\}) = \dil_z(e^{-\ve} s)\leq e^{-\ve} \leq 1.
$
Proceeding similarly if $\lambda(z)>\mu(z)$ we also obtain 
\be
\label{eqDilPVal}
\dil_z(\min\{e^{-\ve} s, \lambda\})\leq 1 \stext{ and } \dil_z(\min\{e^{-\ve} s, \mu\}) \leq 1.
\ee
We now focus on the local dilatation of the orientation $\theta$, and for that purpose we observe that 
$$
|\min\{e^{-\ve}s(z),\lambda(z)\} - \min \{e^{-\ve}s(z), \mu(z)\}|  \leq  |\lambda(z) - \mu(z)|.
$$
Therefore according to Corollary \ref{corolDilThetaG} 
$$
|\min\{e^{-\ve}s(z),\lambda(z)\} - \min \{e^{-\ve}s(z), \mu(z)\}| \dil_z(\theta) \leq 2D_z(a)_\times \leq 2
$$
Combining this with \iref{eqDilPVal} we obtain $D_z(a_\ve)_+ \leq 2$ and therefore $\dil_z(H_\ve)_+ \leq 4$.\\

The above arguments show that 
$
\dil_z(H_\ve)_\times \leq 2
$
and 
$
\dil_z(H_\ve)_+ \leq 4,
$
which implies that $2^2H_\ve \in \bH_g$ and $4^2 H_\ve \in \bH_a$ according to Proposition \ref{propDilLoc} and Remark \ref{remHomog}.
Hence for any $z,z'\in \R^d$, since $H_\ve \leq e^{2\ve} H'$,
$$
d_\times(H_\ve(z), H_\ve(z')) \leq 2d_{H_\ve}(z,z') \leq 2e^\ve d_{H'}(z,z') \stext{ and } d_+(H_\ve(z), H_\ve(z')) \leq 4|z-z'|.
$$
Letting $\ve \to 0$ we obtain 
$$
d_\times(H'(z), H'(z')) \leq 2d_{H'}(z,z') \stext{ and } d_+(H'(z), H'(z')) \leq 4|z-z'|.
$$
Therefore $4^2 H\in \bH_a$ according to Remark \ref{remHomog}, which concludes the proof of this corollary.
\sq

\section{From mesh to metric}
\label{secMeshToMet}
The techniques presented in this section show how to produce an metric $H\in \bH_i$, $\bH_a$ or $\bH_g$ from an isotropic, quasi-acute (in two dimensions), or graded mesh $\cT$. These results are summarized in Proposition \ref{propMetGen} below.

Combining this proposition with its counterpart Proposition \ref{propMeshGen} in the next section, on the generation of a mesh from a metric, yields the main result of this chapter, Theorem \ref{thEquiv}, which states the equivalence of the classes $\bT_{i,C} \subset \bT_{a,C}\subset \bT_{g,C}$ of meshes and $\bH_i \subset \bH_a \subset \bH_g$ of metrics.
\begin{prop}
\label{propMetGen}
For any $C_0>0$ there exits $C = C(C_0,d)$ such that the following holds: 
\begin{enumerate}[i)]
\item For any $\cT \in \bT_{i, C_0}$ there exists $H \in \bH_i$ which is $C$-adapted to $\cT$.
\item If $d=2$, then for any $\cT \in \bT_{a, C_0}$ there exists $H \in \bH_a$ which is $C$-adapted to $\cT$.
\item For any $\cT \in \bT_{g,C_0}$ there exists $H \in \bH_g$ which is $C$-adapted to $\cT$.
\end{enumerate}
\end{prop}

The rest of this section is devoted to the proof of this proposition.
For any mesh $\cT\in \bT$, and for any vertex $v$ of $\cT$, we denote by $n_\cT(v)$ the number of simplices in the mesh $\cT$ containing $v$:
$$
n_\cT(v) := \#\{T\in \cT \sep v\in T\}. 
$$
We denote by $H_\cT^g\in \bH$ the piecewise linear metric on the mesh $\cT$ which satisfies for each vertex $v$ of $\cT$
\be
\label{eqHcTg}
H_\cT^g(v) := \frac 1 {n_\cT(v)} \sum_{\substack{T\in \cT\\\text{s.t. } v\in T}} \cH_T.
\ee

The next lemma establishes a property of the barycentric coordinates on a simplex, which is useful for the analysis of  the regularity of a piecewise linear function on an anisotropic mesh, such as the metric $H_\cT^g$.
\begin{lemma}
\label{lemmaNormLambda}
Let $T$ be a simplex and let $V$ be the collection of vertices of $T$. Let $(\lambda_v)_{v\in V}$ be the barycentric coordinates on $T$. Then for all $z,z'\in T$
\be
\label{eqNormLambda}
\|z-z'\|^2_{\cH_T} = \frac {d+1} d \sum_{v\in V} |\lambda_v(z)-\lambda_v(z')|^2.
\ee
\end{lemma}

\proof
We first assume that $T$ is the reference equilateral simplex $\TEq$. 
For any two vertices $v,v'$ in the collection $\VEq$ of vertices of $\TEq$ one has 
\be
\label{eqScalTEq}
\<v,v'\> = 
\left\{
\begin{array}{cl}
1 &\text{ if } v=v,'\\
-1/d &\text{ if } v \neq v'.
\end{array}
\right.
\ee
Indeed $|v|=1$ by construction for any vertex $v$ of $\TEq$, and for any $v'\in \VEq$ distinct from $v$ the scalar product $\<v,v'\>$ has a value $\alpha$ independent of $v$ and $v'$ by symmetry. Since $\TEq$ is centered at the origin we obtain $0 = \<v,\sum_{v'\in V} v'\>  = |v|^2 + d \alpha$ which establishes \iref{eqScalTEq}.
We thus obtain for any $z,z'\in \TEq$ 
\begin{eqnarray*}
|z-z'|^2  &=& \left|\sum_{v\in \VEq} (\lambda_v(z) - \lambda_v(z'))v\right|^2 \\
&=& \sum_{v\in \VEq} |\lambda_v(z)-\lambda_v(z')|^2 -\frac 1 d \sum_{v\neq v'} (\lambda_v(z)-\lambda_v(z'))(\lambda_{v'}(z)-\lambda_{v'}(z')).
\end{eqnarray*}
Therefore 
$$
|z-z'|^2 =  \left(1+\frac 1 d\right) \sum_{v\in \VEq} |\lambda_v(z)-\lambda_v(z')|^2 -\frac 1 d \left(\sum_{v\in \VEq} \lambda_v(z)-\lambda_v(z')\right)^2,
$$
which concludes the proof in the case of $\TEq$ since $\cH_\TEq = \Id$ and since the barycentric coordinates of $z$ and of $z'$ sum to $1$.

We now consider an arbitrary simplex $T$ and an affine change change of coordinates $\Phi$ such that $\Phi(T) = \TEq$, $\Phi(z) = \phi z+ z_0$ where $\phi\in \GL_d$ and $z_0\in \R^d$. We have $\cH_T = \phi^\trans \phi$ according to \iref{invcH}, hence for all $z,z'\in T$
$$
 |\Phi(z) - \Phi(z')| = |\phi(z-z')| = \|z-z'\|_{\cH_T}.
$$
Furthermore the barycentric coordinates of $z$ and $z'$ in $T$ are the same as those of $\Phi(z)$ and $\Phi(z')$ in $\Phi(T) = \TEq$. This implies \iref{eqNormLambda} and concludes the proof of this lemma.
\sq

The next proposition immediately implies Point iii) of Proposition \ref{propMetGen}. Indeed for any mesh $\cT\in \bT_{g,C_0}$ the metric $d C_0^{10} H_\cT^g$ belongs to $\bH_g$ and is $C_0^6 \sqrt d$ equivalent to $\cT$.

\begin{prop}
\label{propHTG}
For any $\cT\in \bT_{g,C_0}$ the metric $H_\cT^g$ is $C_0$-adapted to $\cT$, and $d C_0^{10} H_\cT^g\in \bH_g$.
\end{prop}

\proof
We consider a simplex $T\in \cT$, and we recall that for any neighbor $T'$ of $T$ one has 
$$
C_0^{-2} \cH_T \leq \cH_{T'} \leq C_0^2 \cH_T.
$$
Averaging these matrices as in \iref{eqHcTg} we obtain for any vertex $v$ of $T$
\be
\label{eqHTHv}
C_0^{-2} \cH_T \leq H_\cT^g(v)\leq C_0^2 \cH_T.
\ee
Averaging with respect to the barycentric coordinates of a point $z\in T$ we obtain 
\be
\label{eqHTHz}
C_0^{-2} \cH_T \leq H_\cT^g(z)\leq C_0^2 \cH_T,
\ee
which establishes as announced that the metric $H_\cT^g$ is $C_0$-adapted to $\cT$.

We now establish the regularity property of $H_\cT^g$, and for that purpose we consider two points $z,z'\in T$ and a vector $u\in \R^d\sm\{0\}$. We thus have 
\begin{eqnarray*}
\frac {\|u\|^2_{H_\cT^g(z)}} {\|u\|^2_{H_\cT^g(z')}} &=& \frac {\sum_{v\in V} \lambda_v(z) \|u\|_{H_\cT^g(v)}^2}{\sum_{v\in V} \lambda_v(z') \|u\|_{H_\cT^g(v)}^2} \\
&=& 1+ \frac {\sum_{v\in V} (\lambda_v(z)-\lambda_v(z')) \|u\|_{H_\cT^g(v)}^2}{\sum_{v\in V} \lambda_v(z') \|u\|_{H_\cT^g(v)}^2}\\
&\leq& 1+ C_0^4 \frac{\sum_{v\in V} |\lambda_v(z)-\lambda_v(z')|}{\sum_{v\in V} \lambda_v(z')}\\
& \leq & 1+C_0^4 \sqrt {(d+1)\sum_{v\in V} |\lambda_v(z)-\lambda_v(z')|^2}\\
&=& 1+C_0^4 \sqrt d \| z-z'\|_{\cH_T},
\end{eqnarray*}
where we used \iref{eqHTHv} in the third line, the Cauchy-Schwartz inequality in the fourth line, and Lemma \ref{lemmaNormLambda} in the fifth line.
Proceeding similarly for  $\|u\|_{H_\cT^g(z')}/\|u\|_{H_\cT^g(z)}$ we obtain
\begin{eqnarray*}
d_\times (H_\cT^g(z), H_\cT^g(z')) &=& \max_{u\neq 0} \left|\ln \|u\|_{H_\cT^g(z)} - \ln \|u\|_{H_\cT^g(z')}\right|\\
&\leq& \frac 1 2\ln (1+C_0^4 \sqrt d \|z-z'\|_{\cH_T}) \\
&\leq& \frac {C_0^4 \sqrt d} 2 \|z-z'\|_{\cH_T}.
\end{eqnarray*}
Hence according to \iref{eqHTHz} and \iref{eqDilNorm} for all $z\in \interior (T)$
$$
\dil_z(H_\cT^g)_\times = \dil_z(H_\cT^g\ssep d_\times, \|\cdot \|_{H(z)}) \leq C_0 \dil_z(H_\cT^g\ssep d_\times, \|\cdot \|_{\cH_T})  \leq C_0^5 \sqrt d/2.
$$
If follows from Corollary \ref{corolLipNoGamma} that $\dil_z(H_\cT^g)_\times\leq C_0^5 \sqrt d/2$ for all $z\in \R^d$, and from Remark \ref{remHomog} that $(d C_0^{10}/4)  H \in \bH_g$, which concludes the proof.
\sq

We associate to any mesh $\cT\in \bT$ a metric $H_\cT^i$ defined as follows : for all $z\in \R^d$, 
$$
H_\cT^i(z) := \|H_\cT^g(z)\| \Id.
$$
The next result immediately implies Point i) of Proposition \ref{propMetGen}. Indeed for any mesh $\cT\in \bT_{i,C_0}$ the metric $d C_0^{10} H_\cT^i$ belongs to $\bH_i$ and is $C_0^7 \sqrt d$ equivalent to the mesh $\cT$.

\begin{corollary}
For any mesh $\cT\in \bT_{i,C_0}$ the metric $H_\cT^i$ is $C_0^2$-adapted to the mesh $\cT$, and $d C_0^{10} H_\cT^i \in \bH_i$.
\end{corollary}

\proof
We have for any simplex $T\in \cT$ 
\be
\label{eqcHNormcH}
C_0^{-2} \|\cH_T\| \Id\leq \frac{\|\cH_T\|} {\rho(T)^2} \Id = \|\cH_T^{-1}\|^{-1} \Id \leq \cH_T \leq \|\cH_T\|\Id.
\ee
For any point $z\in T$ we obtain since $C_0^{-2} H_\cT^g(z) \leq \cH_T \leq C_0^2 H_\cT^g(z)$
$$
C_0^{-2} H_\cT^i(z) = C_0^{-2} \|H_\cT^g(z)\| \Id \leq \|\cH_T\| \Id \leq C_0^2 \|H_\cT^g(z)\| =C_0^2 H_\cT^i(z).
$$
Combining this with \iref{eqcHNormcH} we obtain
$$
C_0^{-4} H_\cT^i(z) \leq \cH_T \leq C_0^2 H_\cT^i(z),
$$
which establishes as announced that the metric $H_\cT^i$ is $C_0^2$ equivalent to the mesh $\cT$.

According to Proposition \ref{propHTG} we have $d C_0^{10} H_\cT^g\in \bH_g$. Corollary \ref{corolLipNorm} thus implies that $d C_0^{10} H_\cT^i\in \bH_i$, which concludes the proof of this proposition.
\sq

Before turning to the case of quasi-acute metrics we prove a lemma on the neighborhood of a simplex in a graded mesh. For any mesh $\cT\in \bT$ and any closed set $E\subset \R^d$ the neighborhood $V_\cT(E)$ of $E$ in $\cT$ is  defined as the union of all simplices $T\in \cT$ which intersect $E$: 
$$
V_\cT(E) := \bigcup_{\substack{T\in \cT\\T \cap E \neq \emptyset}} T.
$$

\begin{prop}
\label{propNeighborT}
Let $\cT\in \bT_{g,C_0}$ and let $T \in \cT$. For all $x\in T$ and all $y \in \R^d \sm V_\cT(T)$ we have 
$$
\|x-y\|_{\cH_T} \geq (C_0\sqrt d)^{-1}.
$$
\end{prop}

\proof
We consider the function $\vp : \R^d \to \R$ which is piecewise linear on the mesh $\cT$, and which satisfies for any vertex $v$ of the mesh $\cT$
$$
\vp(v) = 
\left\{
\begin{array}{cl}
1 & \text{ if } v\in T,\\
0 & \text{ if } v \notin T.
\end{array}
\right.
$$
We denote by $V'$ the collection of vertices of a simplex $T'\in \cT$, and by $(\lambda_v(z))_{v\in V'}$ the barycentric coordinates of a point $z\in T'$.
We have for any $z,z'\in T'$
\begin{eqnarray*}
|\vp(z)- \vp(z')| &\leq& \sum_{v \in V'} |\lambda_v(z)-\lambda_v(z')| \vp(v)\\
& \leq &\sqrt{(d+1) \sum_{v\in V'} |\lambda_v(z)-\lambda_v(z')|^2} \\
&= &\sqrt d \|z-z'\|_{\cH_{T'}}.
\end{eqnarray*}
Since $\vp$ is identically $0$ on $T'$ if $T\cap T' = \emptyset$, and since $\cH_{T'} \leq C_0^2 \cH_T$ otherwise, we obtain for any $z,z'\in T'$
$$
|\vp(z)- \vp(z')| \leq C_0\sqrt d \|z-z'\|_{\cH_{T}}.
$$
Therefore $\dil_z(\vp, \|\cdot\|_{\cH_T})$ for any $T'\in \cT$ and any $z\in T'$, which implies that $\vp : (\R^d, \|\cdot \|_{\cH_T})  \to \R$ is $C_0\sqrt d$-Lipschitz according to Corollary \ref{corolLipNoGamma}.
For any $x\in T$ and any $y\in \R^d\sm V_\cT(T)$ we have $\vp(x) = 1$ and $\vp(y) = 0$, thus $1 = \vp(x)-\vp(y) \leq C_0\sqrt d \|x-y\|_{\cH_T}$ which concludes the proof of this proposition.
\sq

We assume from this point that the dimension is $d=2$. 
For each triangulation $\cT\in \bT$ we define a function $s_\cT : \R^2 \to \R_+^*$ and a metric $H_\cT^a\in \bH$ as follows : for all $z\in \R^2$
\begin{eqnarray}
s_\cT(z) &:=& \inf_{z'\in \sR^2} |z-z'| + \|H_\cT^g(z')^{-\frac 1 2}\|,\\
\label{defHTA}
H_\cT^a(z) &:=& \max\{s_\cT(z)^{-2} , \, H_\cT^g(z)\},
\end{eqnarray}
where the maximum of a real and of a symmetric matrix is defined at \iref{defMaxAlphaM}.
The next proposition immediately implies Point ii) of Proposition \ref{propMetGen}. Indeed the metric $32 C_0^{10} H_\cT^a$ belongs to $\bH_a$ and is $\sqrt{32} C C_0^{5}$ adapted to the mesh $\cT$.
\begin{prop}
\label{propHTA}
For any $C_0 \geq 1$ there exists $C=C(C_0)$ such that the following holds.
For any mesh $\cT\in \bT_{a,C_0}$ the metric $H_\cT^a$ is $C$-adapted to the mesh $\cT$, and $32 C_0^{10} H_\cT^a\in \bH_a$.
\end{prop}

The fact that $32 C_0^{10} H_\cT^a\in \bH_a$ directly follows from the fact that $2 C_0^{10} \bH_\cT^g\in \bH_g$ and from Corollary \ref{corolHgToHa}. 

We define $d_\cT(z) := \|H_\cT^g(z)^{-\frac 1 2}\|$ for each $z\in \R^2$. We establish below that there exists a constant $\eta = \eta(C_0)$, $0< \eta\leq 1$, such that for all $z,z'\in \R^2$
\be
\label{eqRegDT}
\text{ if } |z'-z| \leq \eta d_\cT(z) \stext{ then }  d_\cT(z') \geq \eta d_\cT(z).
\ee
Before turning to the proof of this property we show how it leads to the proof of Proposition \ref{propHTA}.
For all $z\in \R^d$ we have 
$$
s_\cT(z) = \inf_{z'\in \sR^2} |z-z'|+d_\cT(z')  \geq \inf_{z'\in \sR^2} \min \{|z-z'|, \, d_\cT(z')\} \geq \eta d_\cT(z'),
$$
indeed $|z-z'|\geq \eta d_\cT(z)$ or $d_\cT(z') \geq \eta d_\cT(z)$ for any $z'\in \R^d$ according to \iref{eqRegDT}.
Combining this with the definition \iref{defHTA} of $H_\cT^a$ we obtain for all $z\in \R^2$
$$
\eta^2 H_\cT^a(z) \leq H_\cT^g(z)\leq H_\cT^a(z). 
$$
Recalling that $H_\cT^g$ is $C_0$-adapted to the mesh $\cT$, see proposition \ref{propHTG}, we obtain for all $T\in \cT$ and all $z\in T$ 
$$
\eta^2 C_0^{-2} H_\cT^a(z) \leq C_0^{-2} H_\cT^g(z) \leq \cH_T \leq C_0^2 H_\cT^g(z) \leq C_0^2 H_\cT^a(z)
$$
which establishes that the metric $H_\cT^a$ is $C_0/\eta$ adapted to the mesh $\cT$, and concludes the proof of Proposition \ref{propHTA}.\\

For notational simplicity we denote from this point $H := H_\cT^g$ and $C_1 := \sqrt 2 C_0^5$, hence $C_1^2 H \in \bH_g$ and $\cT$ is $C_0$ adapted to $H$. 
It follows from \iref{eqDiamTHT2} that 
$$
C_0^{-1} \|H(z)^{-\frac 1 2}\| \leq \| \cH_T^{-\frac 1 2}\| \leq \diam(T) \leq 2\| \cH_T^{-\frac 1 2}\| \leq 2C_0 \|H(z)^{-\frac 1 2}\|,
$$
for any $T\in \cT$ and $z\in T$, which implies 
\be
\label{eqDiamDT}
\diam(T)/(2C_0)\leq d_\cT(z) \leq C_0\diam(T).
\ee
The quantity $d_\cT(z)$ should therefore be heuristically regarded as the diameter of the triangle $T\in \cT$ containing $z$. Property \iref{eqRegDT} heuristically states that the diameters are constrained to vary in a Lipschitz manner in a quasi-acute triangulation.

We now turn to the proof of \iref{eqRegDT}, and for that purpose we introduce some notations.
For each $u\in \R^2$, we denote by $u^\perp$ the vector obtained by rotating $u$ by $\cPi/2$ in direct trigonometric orientation.
For each segment $e=[\omega-u, \omega+u]\subset \R^2$ we define the diamond $\lozenge(e)$ of $e$ as follows
$$
\lozenge(e) := \{z\in \R^2 \sep |\<z-\omega,u\>|+ C_0 |\<z-\omega,u^\perp\>| < |u|^2\}.
$$
This set is illustrated on Figure \ref{figLozenge} (left).

\begin{figure}
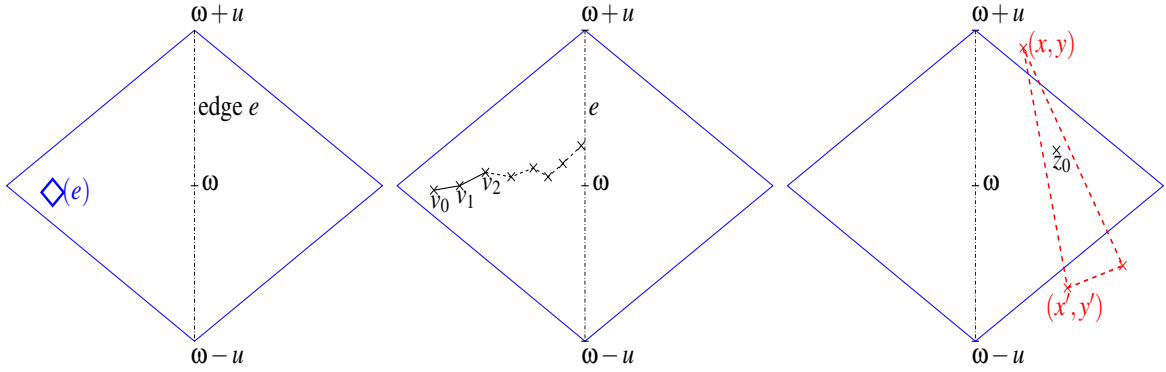

\centering
\includegraphics[width = 5cm, height=5cm]{\pathPic/Dessins/Dessins_JM/DessinJM_07.pdf} 
\includegraphics[width = 5cm, height=5cm]{\pathPic/Dessins/Dessins_JM/DessinJM_08.pdf} 
\includegraphics[width = 5cm, height=5cm]{\pathPic/Dessins/Dessins_JM/DessinJM_09.pdf}
\caption{An edge $e$ and the set $\lozenge(e)$ (left), sequence $(z_n)_{n \geq 0}$ (center), diameter of a triangle intersecting $\lozenge(e)$ (right). \label{figLozenge}}
\end{figure}

The next lemma considers a triangulation $\cT'$ on which the measure of sliverness is uniformly bounded by $C_0$, and compares the length an edge $e$ of $\cT'$ with the diameter of the triangles in $\cT'$ which intersect $\lozenge(e)$.

\begin{lemma}
\label{lemmaUUp}
Let $\cT'\in \bT$ be a triangulation such that $S(T) \leq C_0$ for all $T\in \cT'$. Let  $e=[\omega-u, \omega+u]$ be an edge of $\cT'$.
Then there is no vertex of $\cT'$ in the diamond $\lozenge(e)$, and for all $z\in \lozenge(e)$ and any triangle $T\in \cT'$ containing $z$, we have 
$$
\diam(T) \geq |u| -C_0 |\<z-\omega,u^\perp/ |u|\>|.
$$
\end{lemma}

\proof
We denote by $u_x := (1,0)$ and $u_y := (0,1)$ the canonical basis of $\R^2$, and we assume up to a translation and a rotation that $\omega = 0$ and $u = u_y$.
We recall that for any triangle $T$ with maximal angle $\theta$ one has 
$
S(T) = \max\{1, \tan(\theta/2)\}.
$

We assume for contradiction that there exists a vertex $v_0 = (x_0,y_0)\in \lozenge(e)$, and we assume without loss of generality that $x_0<0$. We consider a vertex $v_1 = (x_1,y_1)$ of $\cT'$ such that $\<v_1-v_0,  u_x\> \geq 0$ and such that the angle
$$
\theta_0 :=  \varangle (v_1-v_0, \, u_x) = \arccos\left(\frac{\<v_1-v_0,  u_x\>}{|v_1-v_0| |u_x|} \right)
$$
has the smallest possible value. By construction there exists a triangle $T\in \cT$ containing $v_0$ and $v_1$ which has an angle  larger or equal or equal than $2\theta_0$ at the vertex $v_0$. It follows that $\tan \theta_0 \leq S(T) \leq C_0$, and therefore  
$$
|y_1-y_0| \leq  (x_1 - x_0)\tan \theta_0\leq C_0 (x_1-x_0).
$$
Proceeding inductively, as illustrated on Figure \ref{figLozenge} (center), we define a sequence $v_n = (x_n, y_n)$, $n\geq 0$, of vertices of $\cT$ satisfying 
$
|y_{n+1} - y_n| \leq C_0 (x_{n+1} - x_n). 
$
Adding these inequalities together we obtain 
$$
|y_n-y_0|  \leq C_0 (x_n-x_0).
$$
It follows that one of the edges $[v_n, v_{n+1}]$, $n \geq 0$, of the mesh $\cT'$ intersects the edge $e$, which is a contradiction. As announced the diamond $\lozenge(e)$ therefore does not contain any vertex. \\

We now consider a point $z_0 = (x_0,y_0)\in \lozenge(e)$ and we assume without loss of generality that $x_0 \geq 0$. We consider a triangle $T$ containing $z_0$, and an edge $e'= [\omega'-u', \omega'+u']$ of $T$ which is on the left of $z_0$ (in the sense that the edge $e'$ and the segment $[(0,y_0),(x_0, y_0)]$ joining the point $z_0$ to its orthogonal projection on the edge $e$ intersect). The point $z_0$ and the triangle $T$ are illustrated on Figure \ref{figLozenge} (right).

Since there is no vertex of $\cT'$ in the diamond $\lozenge(e)$, the edge $e'$ intersects the boundary of $\lozenge(e)$ at two points. We denote these points by $(x,y)$ and $(x',y')$, and we assume without loss of generality that $y\geq 0 \geq y'$.
Observing that $0 \leq x \leq C_0^{-1}$, $0 \leq x'\leq C_0^{-1}$ and $\min\{x,x'\} \leq x_0$ we obtain  
$$
y-y' = (1-C_0 x) - (-1+C_0 x') \geq 1-C_0 x_0.
$$
Since $u = u_y = (0,1)$, $\omega = 0$ and $2 u' = \lambda (x-x', y-y')$ for some $|\lambda| \geq 1$, we obtain 
$$
\diam(T) = 2 |u'| \geq 2|\<u,u'\>| =
y-y' \geq 1-C_0 x_0 = |u|^2 - C_0 |\<z_0-\omega,u^\perp\>|,
$$
which concludes the proof of this proposition.
\sq

\begin{figure}
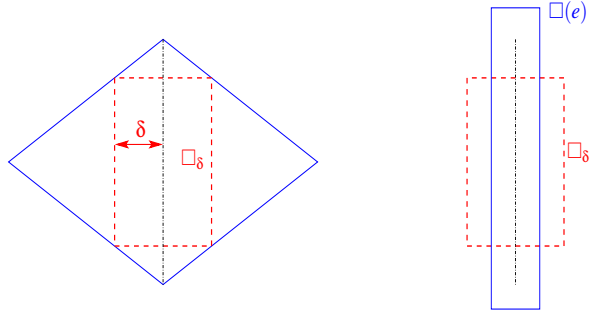

\centering
\includegraphics[width = 4.5cm, height = 4.5cm]{\pathPic/Dessins/Dessins_JM/DessinJM_05.pdf}
\includegraphics[width = 4.5cm, height = 4.5cm]{\pathPic/Dessins/Dessins_JM/DessinJM_06.pdf}
\caption{An edge $e$, the diamond $\lozenge(e)$ and the rectangle $\square_\delta$ involved in Lemma \ref{lemmaRect} (left), the rectangles $\square(e)$ and $\square_\delta$.
\label{figSquare}
}
\end{figure}

From this point we denote by $\cT'$ a $C_0$-refinement of the quasi-acute mesh $\cT\in \bT_{a,C_0}$, which satisfies $S(T) \leq C_0$ for all $T\in \cT'$. The next lemma involves a rectangular set $\square(e)$ associated to an edge $e$ of $\cT'$, illustrated on Figure \ref{figSquare} (right). Note that $e\subset \interior (\square(e))$, while $e\not \subset \lozenge(e)$.

It follows from \iref{eqDTimesDiffLog} and Proposition \ref{propEuclidRiemann} that for all $z,u\in \R^d$
\be
\label{eqRegDT2}
d_\cT(z+u) \geq (1-C_1 \|u\|_{H(z)} ) d_\cT(z).
\ee

\begin{lemma}
\label{lemmaRect}
There exists a constant $\delta_0 = \delta_0(C_0, d)>0$ such that the following holds.
For each edge $e= [\omega-u, \, \omega+u]$ of $\cT'$, let $\square(e)$ be the rectangle defined by 
$$
\square(e) = \{z\in \R^2 \sep |\<z-\omega, u\>|\leq (1+\delta_0)|u|^2 \text{ and } |\<z-\omega,u^\perp \>| \leq \delta_0 |u|^2\}.
$$
One has $d_\cT(z) \geq \delta_0|u|$ for all $z\in \square(e)$.
\end{lemma}

\proof
We denote by $u_x := (1,0)$ and $u_y := (0,1)$ the canonical basis of $\R^2$. 
Up to a rotation and a translation, we may assume that $\omega = 0$ and $u = u_y=(0,1)$.
We consider a triangle $T\in \cT$ containing the edge $e$ (note that $e\subset T$ but that $e$ may not be an edge of $T$, since $e$ is an edge in the triangulation $\cT'$ but may not be an edge in the triangulation $\cT$).
It follows from \iref{eqDiamDT} that 
$$
d_\cT(z) \geq \diam(T)/(2C_0) \geq \diam(e)/(2C_0) = C_0^{-1}
$$
for all $z\in T$, hence for all $z\in e$. 
Hence for all $z\in e$ and all $u\in \R^d$, we obtain using \iref{eqRegDT2}
$$
d_\cT(z+u) \geq  (1-C_1 \|u\|_{H(z)}) C_0^{-1}\geq (1-C_1C_0 \|u\|_{\cH_T})C_0^{-1}.
$$
We define $C_2 := 4C_0C_1$ and we observe that, since $\|u_y\|_{\cH_T} \leq 1$, any point $z'$ of the rectangle 
\be
\label{eqRectGrad}
\left[\frac {-1}{C_2 \|u_x\|_{\cH_T}}, \,\frac 1 {C_2 \|u_x\|_{\cH_T}}\right] \times \left[-(1+C_2^{-1}), \, 1+C_2^{-1}\right].
\ee
can be written under the form $z'=z+u$ where $z\in e$ and $\|u\|_{\cH_T} \leq 1/(2C_0C_1)$, and therefore 
$
d_\cT(z') \geq 1/(2C_0).
$

We introduce a small constant $\delta = \delta(C_0)>0$ which is specified later in the proof of this lemma. It follows from the above argument that the proof of this lemma is complete if $\|\cH_T\| \leq \delta^{-2}$, with $\delta_0 := \min \{\delta/C_2, \, 1/(2C_0)\}$.
We thus assume from this point that $\|\cH_T\|\geq \delta^{-2}$.\\

We assume that $\delta< C_0^{-1}$ and we define the rectangle
$$
\square_\delta := [-\delta, \delta] \times [1-C_0 \delta, -1+C_0 \delta],
$$
which is included in $\overline \lozenge(e)$ by construction.
Since $C_1^2 H\in \bH_g$, the function $z\mapsto \|H(z)\|^{-\frac 1 2}$ is $C_1$-Lipschitz according to Corollary \ref{corolLipNorm}. We thus have for all $z = (x,y)\in \square_\delta$
\be
\label{defr0}
\|H(z)\|^{-\frac 1 2} \leq \|H(0,y)\|^{-\frac 1 2}+C_1|x|\leq C_0 \|\cH_T\|^{-\frac 1 2} + C_1 |x| \leq r_0(\delta) := (C_0+C_1)\delta.
\ee
On the other hand it follows from Lemma \ref{lemmaUUp} that any $z = (x,y)\in \square_\delta$ is contained in a triangle $T'$ such that $\diam(T') \geq (1-C_0 \delta)$,
hence according to \iref{eqDiamDT}
\be
\label{defR0}
d_\cT(z) \geq \diam(T')/(2C_0) \geq  R_0(\delta) := (1-C_0 \delta)/(2C_0). 
\ee
For all $z\in \square_\delta$, we denote by $\theta(z) = (\cos \Theta(z), \sin \Theta(z))$, $\Theta(z)\in [0, \cPi]$, the direction of the eigenvector associated to the small eigenvalue of $H(z)$, in such way that 
\be
\label{eqNormTheta}
\|\theta(z)\|_{H(z)} = \|H(z)^{-\frac 1 2}\|^{-1} = d_\cT(z)^{-1}.
\ee
According to \iref{eqRegDT2} we have for all $z\in \square_\delta$ and all $r\in \R$
$$
d_\cT(z+r \theta(z)) \geq  d_\cT(z) (1- C_1 |r|/d_\cT(z)).
$$
For all $z\in \square_\delta$ and all $r\in \R$ such that $|r|\leq R_1(\delta) := R_0(\delta)/(2C_1)$ we thus obtain
\be
\label{DegDTR0}
d_\cT(z+r\theta(z)) \geq R_0(\delta)(1-C_1|r|/R_0(\delta))\geq R_0(\delta)/2.
\ee

We have for any $(0,y)\in e$
$$
C_0 \geq C_0 \|u_y\|_{\cH_T} \geq \|u_y\|_{H(0,y)} \geq |\cos \Theta(0,y)| \sqrt{\|H(0,y)\|} \geq |\cos \Theta(0,y)| /(\delta C_0) ,
$$
hence 
$
|\cos \Theta(0,y)| \leq  C_0^2 \delta.
$
For all $z\in \square_\delta$ we have according to Corollary \ref{corolDilThetaG} 
$$
(R_0(\delta)-r_0(\delta)) \dil_z(\Theta) \leq  (\|H(z)^{-\frac 1 2}\| - \|H(z)\|^{-\frac 1 2}) \dil_z(\Theta) \leq 2C_1.
$$
Hence we obtain for all $z= (x,y)\in \square_\delta$  
$$
(R_0(\delta)-r_0(\delta)) |\Theta(z) - \Theta(0,y)| \leq 2C_1 |x|\leq 2C_1 \delta,
$$
hence since $t \mapsto \cos t$ is a Lipschitz function 
\be
\label{defTheta0}
|\cos \Theta(z)| \leq \cos \Theta_0(\delta) := C_0^2 \delta + \frac{2C_1\delta}{R_0(\delta)-r_0(\delta)}.
\ee
%
It follows that $d_\cT(z')\geq R_0(\delta)/2$ for all $z'$ in the set 
$$
\{z+\ti r\theta(z) \sep z\in \square_\delta,\, |\ti r|\leq R_1(\delta)\} \supset [-\delta_0,\delta_0]\times [-(1+\mu_0),1+\mu_0]
$$
where 
\be
\label{eqD0MU0}
\delta_0 := \delta - r \cos \Theta_0(\delta) \stext{ and }  \mu_0 := -C_0 \delta+ r \sin \Theta_0(\delta)
\ee
and where the constant $0 \leq r\leq R_1(\delta)$ can be freely chosen.
The rest of the proof consists in choosing appropriate constants $\delta$ and $r$.\\

As $\delta \to 0$ we have $r_0(\delta)\to 0$, $R_0(\delta) \to 1/(2C_0)$ and $\cos\Theta_0(\delta) \to 0$. We thus choose $\delta = \delta(C_0)>0$ sufficiently small in such way that 
\be
\label{eqRCS}
R_0 \geq 1/(4C_0)
, \quad
\cos \Theta_0 \leq 1/(8C_0)
, \quad
\sin \Theta_0 \geq 1/2
.
\ee
and such that $r := 4C_0 \delta$, is smaller than $R_1(\delta)= R_0(\delta)/(2C_1)$. Injecting this choice of $r$ and \iref{eqRCS} in \iref{eqD0MU0} we thus obtain 
$$
\delta_0 \geq \delta/2 \stext{ and }\mu_0 \geq C_0 \delta.
$$
This establishes under the hypothesis $\|\cH_T\|\geq \delta^{-2}$ that  $d_\cT(z) \geq 1/(8C_0)$ for all $z$ in the rectangle 
$$
[-\delta/2, \delta/2]\times [-(1+C_0\delta), 1+C_0 \delta],
$$
which concludes the proof of the proposition.
\sq

Our next intermediate result establishes a technical property on the covering of a segment by smaller ones.

\begin{lemma}
\label{lemmaCoverSeg}
Let $x_1, \cdots , x_k\in [0,1]$ and let $h_1, \cdots, h_k >0$ be such that 
\be
\label{eqDefXiHi}
[0,1]\subset \bigcup_{1\leq i \leq k} [x_i-h_i, \, x_i+h_i].
\ee
Let $\delta>0$ and let $h_* = h_*(k,\delta) := \exp(1-2k/\delta)/\delta$. Then 
\be
\label{eqInclSeg}
[0,1]\subset \bigcup_{\substack{1\leq i \leq k\\ \text{s.t. } h_i \geq h_*}} [x_i-(1+\delta)h_i, \, x_i+(1+\delta)h_i].
\ee
\end{lemma}

\proof
For each $y\in [0,1]$ we denote by $h(y)$ the value $h_i$ associated to the segment $[x_i-h_i, x_i+h_i]$ of smallest length containing $y$.
We thus have 
$$
\int_{[0,1]} \frac {dy} {h(y)} \leq \sum_{1\leq i \leq k} \int_{x_i-h_i}^{x_i+h_i} \frac 1 {h_i} = 2k.
$$
We assume for contradiction that \iref{eqInclSeg} does not hold, and we consider a point $x\in [0,1]$ which does not belong to the right hand side of $[0,1]$. For each $y\in [0,1]$ there exists $i$, $1\leq i \leq k$, $h(y) = h_i$, such that $y\in [x_i-h_i, x_i+h_i]$ and $x\notin  [x_i-(1+\delta)h_i, \, x_i+(1+\delta)h_i]$. Therefore one of the two following possibilities holds
$$
\delta h(y) < |x-y| \stext{ or } h(y) < h_*.
$$
Therefore $h(y) \leq \max \{|x-y|/\delta, h_*\}$, which implies 
$$
\int_0^1 \frac {dy} {h(y)}  \geq \int_0^1 \frac {dy}{\max \{|x-y|/\delta, h_*\}} \geq \int_0^{\delta h_*} \frac 1 h_* + \int_{\delta h_*}^1 \frac \delta ydy = \delta(1-\ln(\delta h_*)).
$$
Note that the second inequality is an equality if $x=0$ or $1$.
Thus $2k>\delta(1-\ln(\delta h_*))$, which contradicts the definition of $h_*$. This concludes the proof of this lemma.
\sq

Observe that under the hypotheses of Lemma \ref{lemmaCoverSeg} we have
\be
\label{eqInclSegDouble}
[-h_*\delta,\, 1+h_*\delta]\subset \bigcup_{\substack{1\leq i \leq k\\ \text{s.t. } h_i \geq h_*}} [x_i-(1+2\delta)h_i, \, x_i+(1+2\delta)h_i].
\ee

We now conclude the proof of Proposition \ref{propHTA}.
We consider a point $z\in \R^2$ contained in a triangle $T\in \cT$, and we denote by $e$ the longest edge of $T$. Up to a rotation and a translation, we may assume that $e=[0, u_x]$ joins the origin of $\R^2$ to $u_x := (1,0)$.
Since the triangulation $\cT'$ is a $C_0$-refinement of $\cT$, the edge $e$ of $T$ is built of $k \leq C_0$ edges $(e_i)_{1 \leq i \leq k}$ of $\cT'$,
$$
e_i = [(x_i-h_i)u_x,\,  (x_i +h_i) u_x].
$$
The triangle $T$ and the edges $(e_i)_{1\leq i\leq k}$ are illustrated on Figure \ref{figTFlat} (left).

\begin{figure}
\centering
\includegraphics[width = 10cm, height = 2cm]{\pathPic/Dessins/Dessins_JM/DessinJM_10.pdf}
\includegraphics[width = 11cm, height = 3cm]{\pathPic/Dessins/Dessins_JM/DessinJM_11.pdf}
\caption{
A triangle $T\in \cT$ and the edges $(e_i)_{1\leq i \leq k}$ of $\cT'$ contained in $T$ (left), the set \iref{eqRectHoDo} (right). 
\label{figTFlat}}
\end{figure}

We define $h_* := h_*(C_0, \delta_0/2) = 2\exp(1-4C_0/\delta_0)/\delta_0$, where $\delta_0$ is the constant from Lemma \ref{lemmaRect}, and we obtain using \iref{eqInclSegDouble} that 
\be
\label{eqRectHoDo}
\bigcup_{h_* \leq h_i} \square(e_i) \supset \left[-\frac{h_* \delta_0} 2, 1+ \frac{h_* \delta_0} 2\right] \times [-h_* \delta_0, h_* \delta_0],
\ee
This set is illustrated on Figure \ref{figTFlat} (right).
Since $e$ is the longest edge of $T$, the angles of $T$ at the extremities of $e$ are acute.
Furthermore the height $h$ of $T$ orthogonal to the edge $e = [0,u_x]$ satisfies 
$$
h/2 = |T|  = \frac{|\TEq|} { \sqrt {\det \cH_T}} = \frac{|\TEq|}{\sqrt{\|\cH_T\| \|\cH_T^{-1}\|^{-1}}}  = \frac{|\TEq| \|\cH_T^{-1}\|}{\rho(T)} \geq \frac{|\TEq|}{4 \rho(T)},
$$
where we used that $\|\cH_T^{-\frac 1 2}\| \geq \diam(T)/2 \geq 1/2$ according to \iref{eqDiamTHT2}.
Hence we obtain with $c := |\TEq|/8$ the inclusion
\be
\label{eqTRectRho}
T \subset [0,1] \times [-c/\rho(T), c/\rho(T)]
\ee
We distinguish two cases depending on the value of $\rho(T)$. If $\rho(T) \geq \rho_0 := 2 c/(h_* \delta_0)$ we obtain combining \iref{eqRectHoDo} and \iref{eqTRectRho} that $d_\cT\geq h_*\delta_0$ on the set 
$$
z+ \left[-\frac {h_*\delta_0} 2, \frac {h_*\delta_0} 2\right]^2.
$$
Since $d_\cT(z) \leq C_0 \diam(T) = C_0$ this implies \iref{eqRegDT} with $\eta = h_* \delta_0 /(2C_0)$.
One the other hand if $\rho(T) \leq  \rho_0$, then 
$$
\|H(z)\|^\frac 1 2 \leq C_0 \|\cH_T\|^\frac 1 2 \leq \rho_0 C_0 \|\cH_T^{-\frac 1 2}\|^{-1} \leq \rho_0 C_0^2 \|H(z)^{-\frac 1 2}\|^{-1} \leq \frac{C_0^2 \rho_0}{d_\cT(z)}
$$
which implies for any vector $u\in \R^2$
$$
\frac{d_\cT(z+u)}{d_\cT(z)} \geq  1-C_1 \|u\|_{H(z)} \geq 1-\frac{C_* |u|}{d_\cT(z)}, 
$$
where $C_* := C_1 \rho_0 C_0^2$. We thus obtain \iref{eqRegDT} with the constant $\eta := 1/(2C_*)$, which concludes the proof of Proposition \ref{propHTA}.

\section{From metric to mesh}
\label{secMetToMesh}
\label{secMeshGen}

This section is devoted to the proof of the following result, on the generation of a mesh from a metric.
\begin{prop}
\label{propMeshGen}
There exists $C_0 = C_0(d)$ such that the following holds: 
\begin{enumerate}[i)]
\item For any $H\in \bH_i$ there exists $\cT \in \bT_{i,C_0}$ which is $C_0$-adapted to $H$.
\item If $d=2$, then for any $H \in \bH_a$ there exists $\cT\in \bT_a$ which is $C_0$-adapted to $\cT$.
\item If $d=2$, then for any $H \in \bH_g$ there exists $\cT\in \bT_g$ which is $C_0$-adapted to $\cT$.\end{enumerate}
\end{prop}

Point {\it i)} is established in \S\ref{secMeshGenI}, Point {\it iii)} in \S \ref{secMeshGenG} and eventually Point {\it ii)} in \S \ref{secMeshGenA}.


\begin{remark}
\label{remcC}
Let $\star \in \{i,a,g\}$ be a symbol, and let $1\geq c>0$ and $C\geq 1$ be numerical constants. 
Assume that for any metric $H$ such that $c^2 H\in \bH_\star$ there exists a mesh $\cT\in \bT_{\star, C}$ which is $C$ adapted to $H$.
Then for any metric $H'\in \bH_\star$ there exists a mesh $\cT' \in \bT_{\star, C}$ which is $C/c$-adapted to $H'$.
(Indeed, choose a mesh $\cT'$ which is $C$-equivalent to the metric $H'/c^2$, and observe that $\cT'$ is $C/c$-equivalent to $H'$.)
\end{remark}

\subsection{Isotropic Mesh generation}
\label{secMeshGenI}

\begin{figure}
\centering
\includegraphics[width = 5cm, height = 5cm]{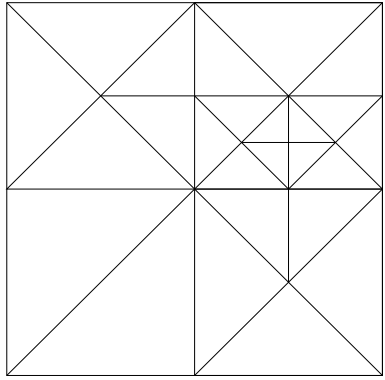}
\caption{A bidimensional triangulation generated by the isotropic refinement procedure exposed in \cite{Nochetto}.
\label{figHalfSquares}}
\end{figure}

The survey paper \cite{Nochetto} describes a the construction of a hierarchical family of meshes of the rectangular domain 
$
[0,1]^d.
$ 
All these meshes are the refinement of a ``fundamental'' mesh $\cT_0 = \{T_\sigma \sep \sigma \in \Sigma_d\}$, where $\Sigma_d$ stands for the collection of permutations of the set $\{1,\cdots,d\}$, and $T_\sigma$ for the Kuhn simplex
$$
T_\sigma := \{z\in [0,1]^d \sep z_{\sigma(1)} \leq \cdots \leq z_{\sigma(d)}\}.
$$
Starting from $\cT_0$, the algorithm described in \cite{Nochetto} produces, guided by the user, a family $(\cT_n)_{n \geq 0}$ of conforming meshes of the cube $[0,1]^d$ proceeding as follows. 
Assume that the mesh $\cT_n$ has already been generated.
\begin{itemize}
\item (Marking)
The user selects a subset $\cM_n \subset \cT_n$.
\item (Refinement)
The algorithm generates a $2$-refinement $\cT_{n+1}$ of $\cT_n$,  
$$
\cT_{n+1} := \Refine(\cT_n, \cM_n).
$$
which does not contain any of the marked simplices : $\cT_{n+1} \cap \cM_n = \emptyset.$
\end{itemize}
\noindent

The elements of the meshes $(\cT_n)_{n \geq 0}$ have a very specific structure, which main features are described below.
\begin{property}
\label{propertyIso}
Any simplex $T$ produced by this algorithm satisfies the following properties.
\begin{enumerate}[a)]
\item (Discrete Volumes) The volume of $T$ has the form
\be
\label{eqVolTG}
|T| = \frac {2^{-g(T)}}{d!},
\ee
where $g(T)\geq 0$ is an integer called the ``generation'' of $T$.
\item
(Finite number of classes)
The simplex 
$$
2^{\frac{g(T)} d} (T-z_T),
$$
belongs to a finite family $\gT_0$ of simplices.
\end{enumerate}
\end{property}
The survey paper \cite{Nochetto} also contains a result, Proposition 4.1, which establishes that $\cT_{n+1} := \Refine(\cT_n, \cM_n)$ is a \emph{local} refinement of $\cT_n$ in the following sense: for all $T'\in \cT_n \sm \cT_{n+1}$ there exists $T\in \cM_n$ such that 
\be
\label{eqDG}
g(T)\geq g(T') \stext{ and } d(T,T') \leq D \, 2^{-g(T')/d},
\ee
where $D = D(d)$ is a fixed constant and where 
$$
d(T,T') := \min\{|z-z'| \sep z\in T, \, z'\in T'\}.\\
$$

The next lemma shows that this refinement algorithm can be used to produce a conforming isotropic mesh $\cT$ of the unit cube $[0,1]^d$, such that the diameters $\diam(T)$, $T\in \cT$, are specified by a given Lipschitz function.

\begin{lemma}
\label{lemmaSCube}
Let $s:[0,1]^d\to  (0,\sqrt d]$ be a $c$-lipschitz function, where $0<c\leq 1$. Let $(\cT_n)_{n \geq 0}$ be the sequence of meshes produced by the algorithm when the collection $\cM_n\subset \cT_n$ of marked simplices is defined as follows:
\be
\label{defMn}
\cM_n := \left\{T\in \cT_n \sep \min_T s < \diam(T)\right\}.
\ee
Then the sequence of meshes stabilizes : there exists $N\geq 0$ such that $\cT_n = \cT_N$ for all $n \geq N$. Furthermore we have for all $T\in \cT_N$ and all $z\in T$
\be
\label{eqEquivCube}
\diam(T) \leq s(z) \leq K_c \diam(T),
\ee
where 
\be
\label{defKc}
K_c := 2^\frac 1 d \frac {D_+ (1+c) + c D}{D_-}
\ee
in which $D_+ := \max\{ \diam(T) \sep T \in \gT_0\}$ and $D_- := \min\{ \diam(T) \sep T \in \gT_0\}$.
\end{lemma}

\proof
We first show that for any $n \geq 0$ any $T_*\in \cT_n$, and any $z\in T_*$ we have $s(z) \leq K_c \diam(T_*)$.
Indeed if $T_*\in \cT_0$ then 
$$
s(z) \leq \sqrt d =  \diam(T_*) \leq K_c \diam(T_*)
$$
since $K_c\geq 1$ and since $s$ is uniformly bounded by $\sqrt d$.

Otherwise $T_*$ has a father $T'$ which was bisected in the refinement process, hence there exists a simplex $T$ such that \iref{eqDG} holds and which was selected for bisection at some stage of the algorithm.

We thus obtain for any $z\in T$, since $g(T_*)-1 = g(T') \geq g(T)$ and $\min_T s < \diam(T)$
\begin{eqnarray*}
s(z) &\leq& \min_T s+ c(\diam(T') + d(T,T'))\\
& \leq & \diam(T)+c(\diam(T') + d(T,T'))\\
& \leq & D_+ 2^{-g(T)/d} + c(D_+2^{-g(T')/d} + D 2^{-g(T')/d})\\
& \leq &  2^{-(g(T_*)-1)/d} (D_+(1+c)+cD)\\
&=& K_c D_- 2^{-g(T_*)/d}\\
& \leq & K_c \diam(T_*),
\end{eqnarray*}
which concludes the proof of the right part of \iref{eqEquivCube}.

From the fourth line of this inequality we also obtain an uniform lower bound on the volume of the simplices generated by the refinement procedure
$$
|T_*| = \frac{2^{-g(T_*)}}{d!} \geq c_0 := \frac 1 {d!} \left(\frac{\min\{ s(z) \sep z\in [0,1]^d\}}{2^\frac 1 d  (D_+(1+c)+cD)}\right)^d.
$$
It follows that $\#(\cT_n)$ is uniformly bounded by $1/c_0$, which immediately implies that the sequence $(\cT_n)_{n \geq 0}$ of meshes stabilizes.

Let $N$ be such that $\cT_n = \cT_N$ for all $n \geq N$. We thus have $\cM_N = \emptyset$ and therefore $\diam(T) \leq \min_T s$ for all $T\in \cT_N$, which implies the left part of \iref{eqEquivCube} and concludes the proof of this lemma.
\sq

The rest of this section is devoted to the proof the point i) of Proposition \ref{propMeshGen}.
We thus consider a metric $H\in \bH_i$, and we recall that there exists a Lipschitz function $s: \R^d \to \R_+^*$ such that for all $z\in \R^d$
$$
H(z) = \frac \Id {s(z)^2}.
$$
For each $n\geq 0$ we define a function $s_n : [0,1]^d \to \R_+^*$, which is clearly Lipschitz, as follows
\be
\label{defSn}
s_n(z) := 2^{-n} s(2^n (z-z_0))
\ee
where $z_0 = (1/2, \cdots , 1/2)$ is the barycenter of the cube $[0,1]^d$.
For all $z\in [0,1]^d$ and all $n \geq 0$ we have 
$$
s_n(z) \leq s_n(z_0) + |z-z_0| \leq 2^{-n} s(0) + \frac {\sqrt d} 2.
$$
The function $s_n$ is therefore uniformly bounded by $\sqrt d$ on $[0,1]^d$ when $n$ is sufficiently large.
We denote by $\cT^n$ the mesh of the cube $[0,1]^d$ described by Lemma \ref{lemmaSCube} for the function $s_n$. We thus have 
\be
\label{eqDiamSn}
\diam(T) \leq s_n(z) \leq C \diam(T) \stext{ for all } T \in \cT^n, \, z\in T.
\ee
where $C=C(d)$ is the constant from Lemma \ref{lemmaSCube}.
We denote by $\cT_n$ the mesh of the cube $[-2^{n-1}, 2^{n-1}]^d$ obtained by translating the mesh $\cT^n$ by $-z_0$ and dilating it by $2^n$. In mathematical terms: 
$$
\cT_n := \{2^n(T-z_0)\sep T \in \cT^n\}.
$$
In view of \iref{defSn} and \iref{eqDiamSn} we thus have 
\be
\label{eqMinSDiamT}
\diam(T) \leq s(z) \leq K_1 \diam(T) \stext{ for all } T \in \cT_n,\, z\in T,
\ee
We denote by $C_0$ the smallest constant such that for all $T$ in the finite set $\gT_0$ one has 
$$
 \frac{C_0^{-2} \Id}{\diam(T)^2} \leq \cH_T \leq \frac{C_0^2\Id}{K_1^2\diam(T)^2 },
$$
and we thus obtain for all $T\in \cT_n$ and all $z\in T$
\be
\label{eqHHTIso}
C_0^{-2} H(z) \leq \cH_T \leq C_0^2 H(z).
\ee
Unfortunately the mesh $\cT_n$ does not cover $\R^d$ but only the cube $[-2^{n-1}, 2^{n-1}]^d$. 
The next lemma shows that a global mesh $\cT$ of the infinite domain $\R^d$ can be extracted from the sequence of meshes $(\cT_n)_{n\geq 0}$.

\begin{lemma}
\label{lemmaExtractMesh}
Let $(\Omega_n)_{n \geq 0}$ be an increasing sequence of polygonal domains of $\R^d$, which exhausts $\R^d$. In other words
$$
\bigcup_{n \geq 0} \Omega_n = \R^d \stext{ and } \Omega_n \subset \Omega_{n+1} \text{ for all } n\geq 0.
$$
Let $H \in \bH_g$, let $(\cT_n)_{n \geq 0}$ be a sequence of conforming simplicial meshes of the domains $\Omega_n$. Assume that for all $n \geq 0$, all $T\in \cT_n$ and all $z\in T$ one has 
\be
\label{eqHHTConv}
C_0^{-2} H(z) \leq \cH_T \leq C_0^2 H(z).
\ee
Then there exists a mesh $\cT\in \bT$ such that 
$
C_0^{-2} H(z) \leq \cH_T \leq C_0^2 H(z)
$
for all $T\in \cT$ and all $z\in T$.
\end{lemma}

\proof
The heuristic idea of this proof is to ``extract a converging subsequence'' from the sequence $(\cT_n)_{n \geq 0}$ of meshes. Unfortunately meshes are discrete objects of combinatorial nature, and the convergence of meshes has no meaning a priori. We therefore use Lipschitz functions as an intermediate object, because their convergence is well defined and they benefit from a compactness property.\\

We associate a function $\phi_T : \R^d \to [0,1]$ to each simplex $T$, which is defined as follows. The function $\phi_T$ is supported on $T$, and for all $z\in T$ we have in terms of the barycentric coordinates $\lambda_v(z)$ of $z$ with respect to the vertices $v\in V$ of $T$:
$$
\phi_T(z) := \min_{v\in V} \lambda_v(z).
$$
According to Lemma \ref{lemmaNormLambda}, we have for any $z,z'\in T$ and any $v\in V$
$$
|\lambda_v(z) - \lambda_v(z')| \leq \|z-z'\|_{\cH_T},
$$
hence 
\be
\label{eqPhiLip}
|\phi_T(z) - \phi_T(z')| \leq \max_{v\in V} |\lambda_v(z) - \lambda_v(z')| \leq \|z-z'\|_{\cH_T}.
\ee
For all $n \geq 0$ we define the function $\phi_n : \R^d \to [0,1]$ as follows
$$
\phi_n := \sum_{T \in \cT_n} \phi_T.
$$
It follows from \iref{eqHHTConv} and \iref{eqPhiLip}, and \iref{eqDilNorm}, that for all $T\in \cT_n$ and all $z\in \interior(T)$
$$
\dil_z(\phi_n\ssep d_H) \leq C_0.
$$
Since in addition $\phi_n = 0$ on $\R^d \sm \Omega_n$, Corollary \ref{corolLipNoGamma} implies that $\phi_n : (\R^d,d_H) \to \R$ is $C_0$-Lipschitz. It follows from Ascoli's compactness theorem (for instance the version recalled in Theorem \ref{thAscoli} in the next section) that there exists 
subsequence $(\phi_{n_k})_{k \geq 0}$ which converges uniformly on all compact sets of $\R^d$ to a $C_0$-Lipschitz function $\phi$. 

We therefore assume, up to such an extraction, that 
$
\phi_n \to \phi
$
uniformly on all compact sets of $\R^d$.
We denote by $\cT$ the collection of closures of connected components of the set $\{z\in \R^d \sep \phi(z)>0\}$:
$$
\cT := \left\{ \overline E \sep E \text{ is a connected component of } \{z\in \R^d \sep \phi(z)>0\} \right\}.
$$
For each $z\in \R^d$ and each $n \geq 0$, we denote by $T_n(z)$ an element of $\cT_n$ containing $z$ if it exists, and $T_n(z) = \emptyset$ otherwise. Likewise we denote by $T(z)$ an element of $\cT$ containing $z$ if it exists, and $\emptyset$ otherwise.
If $\phi(z)>0$, then $z\in \interior (T_n(z))$ for all $n$ sufficiently large, and one easily checks that 
\be
\label{eqTTnz}
T_n(z)\text{ converges in Hausdorff distance to }T(z)\in \cT.
\ee
Hence $\cT$ is a collection of simplices, and the convergence $T_n(z)\to T(z)$ in Haussdorf distance, combined with the hypothesis $C_0^{-2} H(z) \leq T_n(z) \leq C_0^2 H(z)$, implies that  
\be
\label{eqHHTConv2}
C_0^{-2} H(z) \leq T(z) \leq C_0^2 H(z).
\ee
In order to conclude this proof we need to show that the collection $\cT$ of simplices is a conforming mesh of $\R^d$. We thus have to check the hypotheses of Definition \ref{defSimplMesh}.
The interiors $\interior(T)$ of the simplices $T\in \cT$ are pairwise disjoint since they are the connected components of the set $\{z\in \R^d \sep \phi(z)>0\}$.
We claim furthermore that $\cT$ covers $\R^d$. Indeed let $z\in \R^d$ be arbitrary and let $n\geq 0$ be such that $z\in \Omega_n$, hence $T_n(z) \neq \emptyset$. Since $C_0^{-2} H(z) \leq \cH_{T_n(z)}\leq C_0^2 H(z)$ there exists a subsequence $(T_{n_{\vp(k)}}(z))_{k\geq 0}$ which converges in the Haussdorf distance to a simplex $T$ containing $z$. We clearly have $T\in \cT$ which establishes as announced that $\cT$ covers $\R^d$.
The fact that $\cT$ is locally finite easily follows from \iref{eqHHTConv2}.\\

We now turn to the proof of the conformity property, and for that purpose we consider an arbitrary Lipschitz function $f : (\R^d, d_H) \to \R$ with compact support. 
Let $T$ be a simplex such that $C_0^{-2} H(z) \leq \cH_T \leq C_0^2 H(z)$. We have for any vertex $v$ of $T$
$$
|f(v) - f(z_T)| \leq d_H(z_T,v) \leq C_0 \|z_T - v\|_{\cH_T} = C_0.
$$
Denoting as before by $V$ the collection of vertices of $T$, and by $(\lambda_v)_{v\in V}$ the barycentric coordinates on $T$, we thus obtain for any $z,z'\in T$
$$
\interp_T f(z) - \interp_T f(z') = \sum_{v\in V} (\lambda_v(z) - \lambda_v(z')) (f(v) - f(z_T)),
$$
where $\interp_T$ denotes the $\P_1$ (piecewise affine) interpolation on a simplex $T$.
Hence using Lemma \ref{lemmaNormLambda}
\begin{eqnarray*}
|\interp_T f(z) - \interp_T f(z')| &\leq& C_0 \sum_{v\in V} |\lambda_v(z) - \lambda_v(z')|\\
 &\leq& C_0 \sqrt{(d+1)  \sum_{v\in V} |\lambda_v(z) - \lambda_v(z')|^2} \\
 &=& C_0 \sqrt d \|z-z'\|_{\cH_T}.
\end{eqnarray*}
The last inequality implies that 
$$
\dil_z(\interp_{\cT_n} f \ssep d_H)\leq C_0^2 \sqrt d
$$
for all $T\in \cT_n$ and all $z\in \interior(T)$, where $\interp_{\cT'}$ denotes the $\P_1$ interpolation on a mesh $\cT'\in \bT$.
Let $n_0 \geq 0$ be such that $\supp(f)\subset \Omega_{n_0}$, which implies that $\interp_{\cT_n} f$ is continuous. It follows from Corollary \ref{corolLipNoGamma} that $\interp_{\cT_n} f : (\R^d, d_H) \to \R$ is $C_0^2\sqrt d$-Lipschitz.

We define 
$
D := \cup_{T \in \cT} \interior(T),
$
and we (abusively) denote by $\interp_\cT f : D \to \R$ the function which coincides with $\interp_T f$ on the \emph{interior} of each $T\in \cT$. 
For any $z\in D$ we obtain using \iref{eqTTnz} that $\interp_{\cT_n} f(z) \to \interp_\cT f(z)$, hence for any $z,z'\in D$
$$
|\interp_\cT f(z) - \interp_\cT f(z')| = \lim_{n \to \infty} |\interp_{\cT_n} f(z)-\interp_{\cT_n} f(z')| \leq \sqrt d C_0^2 \, d_H(z,z').
$$
The piecewise interpolant $\interp_\cT f : D \to \R$ therefore extends to a $\sqrt d C_0^2$-Lipschitz function $\interp_\cT f : (\R^d, d_H) \to \R$. It easily follows that $\interp_\cT f$ is continuous for any compactly supported $C^1$ function $f$. This property characterizes conforming meshes, hence $\cT\in \bT$ which concludes the proof.
\sq

\begin{corollary}(Bidimensional triangulations)
\label{corolBidimIso}
Assume that the dimension is $d=2$. There exists $c>0$ such that the following holds: if $s: \R^2 \to \R_+^*$ is $c$-Lipschitz, then there exist a mesh $\cT$ of $\R^2$ built of half squares, and such that for all $T\in \cT$ and $z\in T$
\be
\label{eqDT32}
\diam(T) \leq s(z) \leq 3 \diam(T)/2.
\ee
\end{corollary}

\proof
When the dimension is $d=2$, the collection $\gT_0$ of triangles contains only half squares, of area $1/2$, hence  $D_+ = D_- = \sqrt 2$. 
As $c \to 0$ the constant $K_c$ defined by \iref{defKc} tends to $\sqrt 2$. Hence $K_c\leq 3/2$ when $c$ is sufficiently small.

From this point we construct as before a sequence of meshes $(\cT_n)_{n \geq n_0}$ of the squares $[-2^{n-1}, 2^{n-1}]$ which satisfies \iref{eqDiamSn} hence \iref{eqDT32}, and we extract a global mesh $\cT$ of $\R^2$ using Lemma \ref{lemmaExtractMesh}. The triangles $T\in \cT$ still satisfy the inequality \iref{eqDT32} since they are limits in the Haussdorff distance of triangles from the triangulations $\cT_n$.
\sq

\subsection{Graded mesh generation}
\label{secMeshGenG}

We prove in this section a result of bi-dimensional mesh generation, Theorem \ref{thEquiMesh}, which implies Point (iii) of Proposition \ref{propMeshGen}. The key ingredients of this section come from the paper \cite{Shew}.
We assume throughout this section that the dimension is $d=2$.

We first introduce the notion of equispaced points in an abstract metric space.
 Let $(X,d_X)$ be a metric space and let $\delta>0$, we say that a subset $\cV \subset X$ is $\delta$-equispaced if the two following properties hold.
\begin{enumerate}[a)]
\item (Covering) The distance from an arbitrary point $x\in X$ to $\cV$ is bounded by $1$ : for all $x\in X$
\be
\label{eqCover}
d_X(x, \cV) := \inf\{d_X(x,v) \sep v \in \cV\} \leq 1.
\ee
\item (Separation) The pairwise distances between the points of $\cV$ are larger than $\delta$ : for all $v, v'\in \cV$ such that $v \neq v'$
\be
\label{eqSep}
d_X(v,v') \geq \delta.
\ee
\end{enumerate}

The next lemma establishes the existence of such a set under certain assumptions. \begin{lemma}
\label{lemmaEquiSpaced}
Let $(X,d_X)$ be a metric space in which the closed balls $B'(x,r) := \{y\in X \sep d_X(x,y) \leq r\}$ are compact for all $x\in X$ and all $r\geq 0$. 
Let $\cV_0\subset X$ be such that $d_X(v,v') \geq 1$ for all $v,v'\in \cV_0$ such that $v\neq v'$. Then there exists a $1$-equispaced set $\cV$ of the metric space $(X,d_X)$, containing $\cV_0$.
\end{lemma}

\proof
We may assume that $\cV_0 \neq \emptyset$, up to including an arbitrary point of $X$. We choose an  arbitrary point  $x_0\in \cV_0$, and 
we define inductively a sequence $(x_n)_{0 \leq n< N}$ of points and $(\cV_n)_{0 \leq n < N}$ of subsets of $X$, where $N\in \N \cup \{\infty\}$ is specified later, proceeding as follows.
Assuming that $\cV_n$ is already defined we consider the closed set 
\be
\label{eqDefXn}
X_n :=  
\{x\in X \sep d_X(x,\cV_n) \geq 1\}.
\ee
If the set $X_n$ is empty, then the sequence ends and we define $N := n+1$. Otherwise we choose $x_{n+1}$ arbitrarily among the minimisers of $d_X(x,x_0)$ on $X^n$:
\be
\label{defxn}
x_{n+1} \in \underset {x\in X_n} \argmin \,d_X(x,x_0),
\ee
and we define $\cV_{n+1} := \cV_n \cup\{x_{n+1}\}$.
Such a minimizer exists since $X_n$ is closed and since the closed balls in $X$ are compact.
We define $N := \infty$ if none of the sets $(X_n)_{n \geq 0}$ is empty. 

We define $\cV := \cup_{0 \leq n <N} \cV_n$.
%
By construction, we thus have $d_X(v, v') \geq 1$ for all $v,v'\in \cV$
which establishes the separation property \iref{eqSep}. 

We now assume for contradiction that the covering property \iref{eqCover} does not hold, which implies that $N= \infty$. We choose a point $x_*\in X$ such that $d_X(x_*, \cV)>1$, and we remark that $x_*\in X_n$ for all $n \geq 0$. The definition \iref{defxn} of $x_{n+1}$ implies that $d_X(x_0, x_n) \leq d_X(x_0, x_*)$ for all $n \geq 0$, hence the collection of points
$
\{x_n \sep n \geq 0\}
$
is included in the closed ball $\{x\in X \sep d_X(x_0, x) \leq d_X(x_0, x_*)\}$ which is compact.
The sequence $(x_n)_{n \geq 0}$ therefore admits a converging subsequence, but this contradict the fact that the pairwise distances between these points are larger than $1$, which concludes the proof of this proposition.
\sq

%
%
%
%

We consider a fixed metric $H \in \bH$. For any $x,y\in \R^2$ we define 
$$
d_x(y) := \|x-y\|_{H(x)}.
$$
We also consider a fixed discrete collection of vertices, or ``sites'', $\cV \subset \R^2$, and we define the (anisotropic) Voronoi cell of a site $v\in V$ as follows 
\be
\label{defVorv}
\Vor(v) := \{z\in \R^2 \sep d_v(z) \leq d_w(z) \text{ for all } w \in V\}.
\ee
The (anisotropic) Voronoi cell $\Vor(v)$ is thus the collection of points $z$ which are \emph{closer to the vertex $v$ than to any other site}. The Voronoi diagram is the collection $\{\Vor(v) \sep v\in \cV\}$ of all Voronoi cells.
The classical isotropic Voronoi is obtained by choosing $H = \Id$ identically on $\R^d$.
An instance of an (anisotropic) Voronoi diagram is illustrated on Figure \ref{figShew}.

For each $z\in \R^2$ we denote by $\cV_z\subset \cV$ the collection of vertices which contain $z$ in their Voronoi region:
$$
\cV_z := \{v \in \cV \sep z \in \Vor(v)\}.
$$
The geometric dual of the Voronoi diagram is the graph $\cG(\cV,H)$ defined as follows:
\be
\label{defGVH}
\cG(\cV,H) := \{[v,v'] \sep \text{there exists } z\in \R^2 \text{ such that } \cV_z = \{v,v'\}\}.
\ee

\begin{figure}
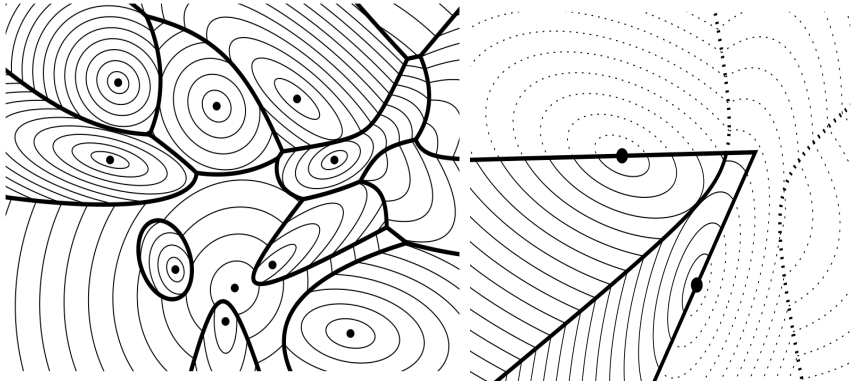

\centering
\includegraphics[width = 6cm, height=5cm]{\pathPic/Dessins/Shew2.png}
\vspace{1cm}
\includegraphics[width = 5cm, height=5cm]{\pathPic/Dessins/Shew1.png}
\caption{An (anisotropic) Voronoi Diagram (left), a wedge (right). (credit J.R.Shewchuck \cite{Shew})\label{figShew}}
\end{figure}
The next theorem establishes that $\cG(\cV,H)$ is, under certain conditions on the metric $H$ and the set $\cV$, the graph of a triangulation adapted to $H$.

\begin{theorem}
\label{thEquiMesh}
For each $\delta>0$ there exists $c=c(\delta)>0$ and $C=C(\delta)\geq 1$ such that the following holds.
Let $H$ be a metric such that $c^2 H\in \bH_g$, and let $\cV$ be a $\delta$-equispaced subset of $(\R^2, d_H)$.
Then the $\cG(\cV,H)$ is a planar graph which defines a partition of $\R^2$ into strictly convex polygons. Arbitrarily triangulating each polygon yields a triangulation $\cT\in \bT_{g,C}$ which is $C$-adapted to $H$, and called an anisotropic Delaunay triangulation. Furthermore if $\cV$ is in general position then $\cG(\cV,H)$ is already the graph of a triangulation.
\end{theorem}

Denote by $c$ and $C$ the constants attached to $\delta=1$ in this theorem.
It follows that for any metric $H$ such that $c^2 H \in \bH_g$, there exists a mesh $\cT\in \bT_{g,C}$ which is $C$-adapted to $H$, obtained as the Delaunay triangulation of a $1$-equispaced set of points $\cV\subset (\R^2, d_H)$. Such a set exists according to Lemma \ref{lemmaEquiSpaced}, since the balls of the metric space $(\R^2, d_H)$ are compact according to Proposition \ref{propEuclidRiemann}.
This immediately implies Point iii) of Proposition \ref{propMeshGen} according to Remark \ref{remcC}.

The rest of this section is devoted to the proof of Theorem \ref{thEquiMesh}, which is based on the methods and results presented in \cite{Shew}. 
For that purpose we recall some definitions of this paper. 
\begin{definition}
\label{defWedge}
\begin{itemize}
\item (Wedge of two points)
We define the wedge $\wed(v,w)$ of two distinct sites $v,w\in V$ as follows
$$
\wed(v,w) := \{z\in \R^2 \sep (z-v)^\trans H(v) (w-v) >0 \text{ and } (z-w)^\trans H(w) (v-w) >0\}
$$ 
\item (Wedge property)
We say that the Voronoi diagram of $V$ is wedged if for any two distinct $v,w \in V$ we have 
$$
\Vor(v) \cap \Vor(w) \subset \wed(v,w).
$$
\end{itemize}
\end{definition}

The wedge property expresses that the collection $\Vor(v) \cap \Vor(w)$ of points shared by the Voronoi regions of the sites $v$ and $w$ lie in a cone defined by two linear inequalities.
The wedge of two points is illustrated on Figure \ref{figShew} (right). Note that the wedge property is not satisfied on this figure.

For any metric $H$ and any constant $c>0$ such that $c^2H\in \bH_g$, it easily follows from Proposition \ref{propEuclidRiemann} that for all $x,y\in \R^2$
\be
\label{eqDHLog}
\ln (1+ c\, d_x(y)) \leq c \, d_H(x,y) \leq - \ln(1- c\, d_x(y)), 
\ee
which is equivalent to 
\be
\label{eqDHExp}
1-\exp(-c\, d_H(x,y)) \leq c \, d_x(y) \leq \exp(c \, d_H(x,y))-1.
\ee

\begin{lemma}
\label{lemmaDiamVor}
The following holds for any $\delta>0$ and any $c>0$. Let $H$ be a metric such that $c^2 H\in \bH_g$ and let $\cV$ be a $\delta$-equispaced subset of $(\R^2, d_H)$. Then
for any $v\in \cV$ and any $q\in \Vor(v)$ one has
\be
\label{eqDiamVor}
d_v(q) := \|q-v\|_{H(v)} \leq r_0(c) := (e^c-1)/c.
\ee
If $v,w \in \cV$ satisfy $\Vor(v) \cap \Vor(w) \neq \emptyset$, then 
\be
\label{eqVordvw}
d_H(v,w) \leq r_1(c) := -2\ln(2-e^c)/c.
\ee
Note that $r_0(c) \to 1$ and $r_1(c) \to 2$ as $c\to 0$. Note also that $d_H(v,w)\geq \delta$.
\end{lemma}

\proof
For any $q\in \cV$ there exists $v_* \in \cV$ such that 
$
d_H(q,v_*) \leq 1,
$
and therefore 
$$
\|q-v_*\|_{H(v_*)} \leq r_0(c) := (e^c-1)/c,
$$
which implies \iref{eqDiamVor} using the definition \iref{defVorv} of the Voronoi diagram.
We thus turn to the proof of \iref{eqVordvw} and for that purpose we consider a point $q\in \Vor (v)\cap\Vor(w)$. We thus have
$
\|q-v\|_{H(v)} = \|q-w\|_{H(w)} \leq  r_0(c)
$
hence using \iref{eqDHLog}
$$
d_H(v,w) \leq d_H(v,q)+d_H(q,w) \leq r_1(c) :=-2\ln(1-c r_0(c))/c,
$$
which concludes the proof.
\sq

The next proposition, which is partly is based on Lemma 5 from \cite{Shew},
 shows that the wedge property is automatically satisfied in our context.

\begin{prop}
\label{propWedge}
For each $\delta$ there exists $c=c(\delta)>0$ such that the following holds.\\
If a metric $H$ satisfies $c^2H\in \bH_g$ and if $\cV$ is a $\delta$-equispaced subset of $(\R^2,d_H)$, then the wedge property is satisfied.
\end{prop}

\proof
We assume for contradiction that there exists two distinct vertices $v,w \in V$, and a point $q\in \Vor(V) \cap \Vor(W)$ such that $q\notin \wed(v,w)$. We may assume, up to exchanging the roles of $v$ and $w$, that 
$$
(q-w)^\trans H(w) (v-w) \leq 0.
$$
Defining the ``relative distortion''
$
\tau := \exp d_\times (H(v), H(w))
$
we obtain using Lemma 5 from \cite{Shew} 
\be
\label{eqShewWedge}
 d_w(v) \leq d_w(q) \sqrt{\tau^2-1}.
\ee
We now estimate the quantities $d_w(v)$, $d_w(q)$ and $\tau$ appearing in this estimate, in order to obtain a contradiction when $c$ is sufficiently small.

We have 
\be
\label{eqTauWedge}
\ln \tau = d_\times (H(v), H(w)) \leq c d_H(v,w) \leq c r_1(c).
\ee
Since $\cV$ is $\delta$-equispaced we have 
$
d_H(v,w) \geq \delta
$
and therefore 
\be
\label{eqdvwWedge}
c \,d_w(v) \geq 1-\exp(-c d_H(v,w)) \geq 1-\exp(-c\delta).
\ee
Injecting \iref{eqDiamVor}, \iref{eqTauWedge} and \iref{eqdvwWedge} in \iref{eqShewWedge} we obtain
\be
\label{eqFinalWedge}
\frac{1-\exp(-c\delta )} c \leq \frac{\exp(c)-1} c \sqrt{\exp(2c r_1(c))-1}.
\ee
The left hand side tends to $\delta$ as $c\to 0$, while the right hand side tends to $0$ (and only depends on $c$). We thus obtain a contradiction when $c$ is sufficiently small which concludes the proof of this proposition.
\sq

\begin{theorem}[Dual triangulation theorem. Labelle, Shewchuck, \cite{Shew}]
\label{thShew}
Let $H\in \bH$ and let $\cV\subset \R^2$ be a discrete point set. If the Voronoi diagram of $\cV$ is wedged, then the geometric dual of this diagram is a polygonalization of $\R^2$ with strictly convex polygons, and is a triangulation if $\cV$ is in general position.
\end{theorem}

Under the hypotheses of this theorem the vertices $v_1, \cdots, v_k$ of a polygon in the dual Voronoi diagram satisfy by construction $\Vor(v_1)\cap \cdots \cap \Vor(v_k) \neq \emptyset$. Let $\cT$ be a triangulation obtained by arbitrarily triangulating these polygons. The vertices $v_1, v_2, v_3$ of a triangle $T\in \cT$ satisfy $\Vor(v_1) \cap \Vor(v_2) \cap \Vor(v_3)\neq \emptyset$.

The next theorem, originally stated in \cite{Shew} as Corollary 10 to Theorem 9 of the same paper, allows to control the angles of a triangle in the anisotropic Delaunay triangulation.

\begin{theorem}[Labelle, Shewchuk, \cite{Shew}]
\label{thAngle}
Let $v_1,v_2,v_3\in \cV$ and let $q\in \Vor(v_1) \cap \Vor(v_2) \cap \Vor(v_3)$. Let 
$$
r = d_{v_1}(q) = d_{v_2}(q) = d_{v_3}(q) \stext{ and }
l = \min\{d_{v_i}(v_j) \sep i,j \in \{1,2,3\}, i \neq j\}.
$$ 
Define $\beta := \max\{r/l,\, 1/\sqrt 2\}$, 
$$
\gamma := \exp \max \{d_\times (H(v_i), H(v_j)) 1\leq i<j \leq 3\},
$$
and
\be
\label{defChi}
\chi := \frac 1 {2 \beta} - \frac{(\gamma^2-1)\beta} 2.
\ee
If $\gamma \leq \sqrt 2$ and $\chi>0$ then any angle $\theta$ of the triangle of vertices $(\sqrt{H(v_1)} v_k)_{1 \leq k \leq 3}$ satisfies
$$
\arcsin ( \chi/ \gamma^2) \leq \theta \leq 2 \arccos(\chi/\gamma^2).
$$
\end{theorem}

\begin{corollary}
\label{corolAdaptG}
For each $\delta$, $0<\delta\leq 1$ there exists $c=c(\delta)>0$ and $C=C(\delta) \geq 1$ such that the following holds.
Let $H$ be a metric such that $c^2H\in \bH_g$ and let $\cV$ be a $\delta$-equispaced subset of $(\R^2, d_H)$. If $v_1,v_2,v_3\in \cV$ and if $\Vor(v_1) \cap \Vor(v_2) \cap \Vor(v_3)\neq \emptyset$, then denoting by $T$ the triangle of vertices $v_1,v_2,v_3$ we have for all $z\in T$
$$
C^{-2} H(z) \leq \cH_T \leq C^2 H(z).
$$
\end{corollary}
\proof
We first obtain obtain some explicit bounds on the quantities $r$, $l$, $\beta$, $\gamma$ and $\chi$ appearing in Theorem \ref{thAngle}.
It follows from Lemma \ref{lemmaDiamVor} that 
$$
r\leq r_0(c) \stext{ and } \delta \leq d_H(v_i,v_j) \leq r_1(c).
$$
Therefore 
\begin{eqnarray*}
l := \min_{i\neq j} d_{v_i}(v_j) &\geq& 
r_2(c) := (1-\exp(-c\delta))/ c,\\
L := \max_{i\neq j} d_{v_i}(v_j) &\leq& (\exp(c r_1(c))-1)/c,
\end{eqnarray*}
and 
$$
\gamma \leq \max_{i\neq j} \exp(c \,d_H(v_i,v_j)) \leq \exp(c \,r_1(c)).
$$
It follows that 
$
\beta \leq r_3(c) := \min\{r_0(c)/r_2(c),\, 1/\sqrt 2\}. 
$
Injecting this into \iref{defChi} yields 
$$
\chi\geq r_4(c) := \frac 1 {2r_3(c)} - \frac{ (\exp(2cr_1(c))-1) r_3(c)} 2,
$$ 
and therefore $\chi/\gamma^2\geq r_5(c) := r_4(c) \exp(-2cr_1(c))$.

For any fixed $\delta$, $0< \delta \leq 1$ we obtain as $c\to 0$ the limits $r_0(c) \to 1$, $r_1(c)\to 2 $, $r_2(c) \to \delta$, $r_3(c) \to \min\{\delta, 1/\sqrt 2\}$, $r_4(c) \to 1/(2\min\{\delta, 1/\sqrt 2\})$ and $r_5(c)$ has the same limit.

We may therefore choose $c$ sufficiently small in such way that 
$$
l \geq l_0 := \delta/2, \ L  \leq L_0 := 4 \text{ and }\chi/\gamma^2 \geq \Chi_0 := 1/(4\delta).  
$$
We denote by $T$ the triangle of vertices $v_1, v_2, v_3$. The length of any edge of $T' := \sqrt{H(v_1)}(T)$ is bounded above by $L_0$ and below by $l_0$ respectively, and the angles of $T'$ satisfy $\arcsin (\Chi_0) \leq \theta \leq 2 \arccos(\Chi_0)$. The collection of triangles centered at the origin which satisfy these inequalities is compact, thus there exists a constant $C_1 = C_1(\delta)$ such that 
$$
C_1^{-2} \Id \leq \cH_{T'} \leq C_1^2 \Id,
$$
hence $C_1^{-2} H(v_1) \leq \cH_{T'} \leq C_1^2 H(v_1)$. For any $z\in T$ we obtain, remarking that $\|z-v_1\|_{\cH_T} \leq 2$ for all $z\in T$ and using \iref{eqLocalNorm}, 
$$
C_1^{-2} (1-2cC_1)^2 H(z) \leq \cH_{T'} \leq C_1^2(1-2cC_1)^{-2} H(z).
$$
We may assume that $c\leq 1/(4C_1)$, which concludes the proof of this proposition with $C = 2C_1$.
\sq

We now prove of Theorem \ref{thEquiMesh}.
According to the above results for any $\delta>0$ there exists $c = c(\delta)$ and $C=C(\delta)$ such that the following holds. 
For any metric $H$ such that $c^2 H\in \bH_g$, and any set $\cV \subset (\R^2, d_H)$ which is $\delta$-equispaced, the Voronoi diagram associated to $\cV$ and $H$ is wedged according to Proposition \ref{propWedge}. The dual of this diagram is therefore a polygonalization of $\R^2$ with strictly convex polygons according to Theorem \ref{thShew}.
Randomly triangulating these cells yields a triangulation $\cT\in \bT$ which is $C$-adapted to $H$ according to Corollary \ref{corolAdaptG}. 

Furthermore consider to simplices $T,T'\in \cT$ such that $T \cap T'\neq \emptyset$, and a point $z\in T \cap T'$.
Then 
$$
C^{-2} \cH_T \leq H(z) \leq C^2 \cH_{T'},
$$
which shows that $\cT\in \bT_{g,C^2}$, and concludes the proof of Theorem \ref{thEquiMesh}.

\subsection{Quasi-Acute mesh generation}
\label{secMeshGenA}

This section is devoted to the proof of Point ii) of Proposition \ref{propMeshGen}, which is obtained by combining the  constructions of isotropic and graded triangulations presented in \S\ref{secMeshGenI} and \S\ref{secMeshGenG} to produce a quasi-acute triangulation. 
The key idea of this proof if that the anisotropic Delaunay triangulation described in the previous subsection contains only \emph{few} strongly obtuse angles if the set $\cV$ of vertices is structured along lines which are \emph{transverse to the direction of anisotropy}, as illustrated on Figure \ref{figAcuteGen} (in the introduction of this chapter). 

We consider a fixed metric $H$ such that $c^2 H\in \bH_a$, where $c$ is an absolute constant specified in the proof and such that $0<c\leq 1$. This subsection is devoted to the construction of a mesh $\cT\in \bT_{a,C}$ which is $C$-adapted to $H$, where $C$ is an absolute constant also specified in the proof.
This construction immediately implies Proposition \ref{propMeshGen} according to Remark \ref{remcC}.

For each $z\in \R^2$ we define 
$$
\lambda(z) := \|H(z)^{-\frac 1 2}\| \stext{ and } \mu(z) = \|H(z)\|^{-\frac 1 2},
$$
and we observe that $\lambda(z) \geq \mu(z)$.
The functions $\lambda$ and $\mu$ are $c$-Lipschitz according to Corollary \ref{corolLipNorm}.
If $c$ is sufficiently small, then there exists according to Corollary \ref{corolBidimIso} a mesh $\cT_0$ such that for all $T\in \cT_0$ and all $z\in T$
\be
\label{ratioDiamLambda}
\frac 2 3  \lambda(z) \leq \diam(T)\leq \lambda(z). 
\ee
The mesh $\cT_0$ is built of half squares by construction (as illustrated on Figure \ref{figHalfSquares}). The diameter of any triangle in $\cT_0$ is a power of $\sqrt 2$, hence for any two triangles $T,T'\in \cT_0$ sharing a vertex $v$ we obtain using \iref{ratioDiamLambda}
\be
\label{eqDiamRatio}
\diam(T)/\sqrt 2 \leq \diam(T') \leq \sqrt 2 \diam(T).
\ee
For any edge $[v,w]$ of any triangle $T\in \cT_0$ we have $\diam(T)/\sqrt 2 \leq |v-w|\leq \diam(T)$, hence for any two edges $[v,w]$ and $[v,w']$ (sharing the vertex $v$) of the triangulation $\cT_0$ 
\be
\label{eqEdgeRatio}
|v-w|/2 \leq |v-w'| \leq 2|v-w|.
\ee
Inequality \iref{ratioDiamLambda} also implies that for any edge $[v,w]$ of $\cT_0$
\be
\label{edgeLambda}
\frac{\lambda(z)\sqrt 2}3 \leq \frac{\diam(T)}{\sqrt 2}\leq |v-w| \leq \diam(T)\leq \lambda(z)
\ee
We denote by $\cV_0$ the collection of vertices of $\cT_0$. 
Consider a fixed $v\in \cV_0$ and any vertex $w\in \cV_0$ such that $|v-w|$ is minimal. Then $[v,w]$ is clearly an edge of $\cT_0$, since this triangulation is conforming and built of half squares. This implies for any distinct $v,w\in \cV_0$
\be
\label{pairwiseV0}
\lambda(v) \frac {\sqrt 2} 3 \leq |v-w| .
\ee

For each $z\in \R^2$ we define 
$$
\rho(z) := \lambda(z)/\mu(z)  = \sqrt{\|H(z)\| \|H(z)^{-1}\|} \in [1, \infty).
$$
Consider $z,z'\in \R^2$ such that $|z-z'| \leq r\lambda(z)$, where $r\geq 0$ is a constant. We obtain since $\lambda$ and $\mu$ are $c$-Lipschitz
\be
\label{eqrhozzp}
\rho(z') = \frac{\lambda(z')}{\mu(z')} \geq \frac{\lambda(z)(1-r c)}{\mu(z)+ r c \lambda(z)} \geq \frac {1-r c}{\rho(z)^{-1}+r c}.
\ee
For each $\rho_0\geq 1$, $r_0>0$ we denote by $c_\rho(\rho_0, r_0)>0$ the constant such that $\rho_0 = \frac {1-r_0 c_\rho(\rho_0, r_0)}{(\rho_0+1)^{-1}+r_0 c_\rho(\rho_0, c_0)}$.
Hence for all $z,z'\in \R^2$ 
$$
\left\{
\begin{array}{c}
|z-z'| \leq r_0\lambda(z)\\ 
\rho(z) \geq \rho_0+1\\
c\leq c_\rho(\rho_0, r_0)
\end{array}
\right.
\stext{ implies } \rho(z') \geq \rho_0
$$

For each $z\in \R^2$ such that $\rho(z)>1$ we denote by $\theta(z)\in \bbS := \{u\in \R^2\sep |u| = 1\}$ a unit vector such that
\be
\label{defAngle}
H(z) = \lambda(z)^{-2} \theta(z) \theta(z)^\trans+ \mu(z)^{-2} (\Id-\theta(z) \theta(z)^\trans).
\ee
Note that there exists two, opposite, choices for the vector $\theta(z)$.
For each $u,u'\in \R^2 \sm\{0\}$ we define the (unoriented) angle
$$
\varangle (u,u') := \arccos \left(\frac{\<u,u'\>}{|u| |u'|}\right) \in [0,\cPi].
$$
and 
$$
\lhd(u,u') := \min \{\varangle(u,u'),\, \varangle(u,-u')\} = \arccos\left(\frac {|\<u,u'\>|}{|u| |u'|}\right) \in [0, \cPi/2].
$$
\begin{lemma}
\label{lemmaAngle12}
If $c$ is sufficiently small, then the following holds.
For any $z,z'\in \R^2$ 
$$
\left\{
\begin{array}{c}
\rho(z) \geq 2\\
|z'-z|\leq 10 \lambda(z) 
\end{array}
\right.
\stext{ implies }
\lhd(\theta(z), \theta(z')) \leq \cPi/12.
$$
Likewise, $\rho(z) \geq 4$ and $|z-z'| \leq \lambda(z)$ implies $\lhd(\theta(z), \theta(z')) \leq 1/(10 \sqrt 2)$.
\end{lemma}

\proof
We consider a fixed point $z$ such that $\rho(z) \geq 2$ and we denote $B:=\{z'\in \R^2\sep |z'-z| \leq 10 \lambda(z)\}$.
For each $z'\in B$ we have
\begin{eqnarray*}
\lambda(z') - \mu(z') &\geq& (\lambda(z)-c|z-z'|) - (\mu(z)+ c |z-z'|) \\
&=& \lambda(z) -\mu(z) -2 c |z-z'| \\
&\geq& \lambda(z) -\lambda(z)/ 2 - 20 c\lambda(z)\\
&=& \lambda(z)(1/2 -20 c),
\end{eqnarray*}
hence $\lambda(z') - \mu(z')\geq \lambda(z)/3$ if $c \leq 1/120$, which we assume from this point.

Furthermore since $B$ is simply connected there exists a continuous function $\theta_* : B \to \bbS$ such that $\theta_*(z')  = \pm \theta(z')$ for all $z'\in B$.
It follows from Theorem \ref{th2Spaces} that for all $z'\in B$
$$
\lambda(z)\dil_{z'}(\theta_*) \leq 3(\lambda(z')-\mu(z')) \dil_{z'}(\theta_*) \leq 3c.
$$
Therefore, since the distance $d_\bbS$ on $\bbS$ involved in Theorem \ref{th2Spaces} coincides with the angle $\varangle$ we obtain
$$
\lhd(\theta(z), \theta(z'))\leq \varangle(\theta_*(z), \theta_*(z')) \leq |z-z'| \max_{\ti z \in [z,z']} \dil_{\ti z}(\theta_*) \leq 30c 
$$
which concludes the proof if $30c \leq \cPi/12$. The second estimate is obtained similarly.
\sq

For each vertex $v\in \cV_0$ we define a set $\Gamma_v\subset \R^2$ as follows : $\Gamma_v := \{v\}$ if $\rho(v) <3$, and 
\be
\label{defGammav}
\Gamma_v := \left[v, \, \frac v 3+ \frac{2w_-} 3\right] \cup \left[v,\frac v 3 +\frac{2w_+} 3\right] 
\ee
if $\rho(v)\geq 3$, where $[v,w_-]$ and $[v,w_+]$ are two edges of $\cT_0$ which respectively minimize and maximize the quantity $\<\theta(v)^\perp, w-v\>$ (among the edges containing $v$).
In other words for any edge of $\cT_0$ of the form $[v,w]$ we have
\be
\label{defWmWp}
\<\theta(v)^\perp, w_- -v\> \leq \<\theta(v)^\perp, w-v\> \leq  \<\theta(v)^\perp, w_+ - v\>. 
\ee
This property is illustrated on Figure \ref{figWmWp} (right).

\begin{figure}
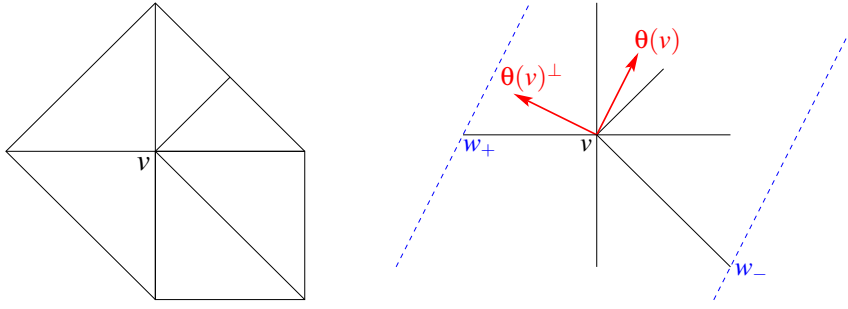

\centering
\includegraphics[width=4cm, height=4cm]{\pathPic/Dessins/Dessins_JM/DessinJM_02.pdf}\hspace{1cm}
\includegraphics[width=6cm, height=4cm]{\pathPic/Dessins/Dessins_JM/DessinJM_03.pdf}
\caption{Neighborhood $V_{\cT_0}(v)$ of a simplex $v$ in the triangulation $\cT_0$ (left), the neighbors $w_+$ and $w_-$ of a vertex $v$ in $\cT_0$ (right).
\label{figWmWp}}
\end{figure}

\begin{lemma}
\label{lemmaAngleWmWp}
let $v\in \cV_0$ be such that $\rho(v) \geq 3$, and let $[v,w_+]$, $[v,w_-]$ be two edges of $\cT_0$ satisfying \iref{defWmWp}.
We have 
\be
\label{eqAngleWmWp}
\varangle(\theta(v)^\perp, w_+ -v) \leq \cPi/4,
\stext{ and }
\varangle(\theta(v)^\perp, v-w_-) \leq \cPi/4.
\ee
\end{lemma}

\proof
We restrict without loss of generality our attention to the proof of the first inequality, and we assume for contradiction that $\vp := \varangle(\theta(v)^\perp, w_+-v)>\cPi/4$. 
Among the two vertices $w\in \cV_0$ such that the triangle of vertices $v,w,w_+$ belongs to $\cT_0$, we choose the one such that $\varangle(\theta^\perp, w-v)$ is minimal.
We denote $\psi := \varangle(w-v, w_+-v) \in \{\cPi/4, \cPi/2\}$ and we obtain 
$$
\<\theta(v)^\perp, w_+-v\> = |w_+-v|\cos \vp \stext{ and } \<\theta(v)^\perp, w-v\> = |w-v|\cos (\vp-\psi).
$$
If $\psi = \cPi/2$ then $|w_+-v| = |w-v|$ since the triangles $T\in \cT_0$ are half squares, and $\cos (\vp-\psi) = \sin \vp$. If follows from \iref{defWmWp} that $\cos \vp \geq \sin \vp$ which contradicts our assumption that $\vp>\cPi/4$.
If $\psi = \cPi/4$ then $|w_+-v| = \sqrt 2 |w-v|$ or $|w_+-v| =  |w-v|/\sqrt 2$. We thus have $\sqrt 2 \cos \vp \geq \cos (\vp-\cPi/4)$ which again contradicts our assumption that $\vp > \cPi/4$.
\sq

We define $\Gamma := \cup_{v\in \cV_0} \Gamma_v$, and for each $z\in \Gamma\sm\cV_0$ we define 
\be
\label{defGammaz}
\Gamma_z := \bigcup_{\substack{v\in \cV_0 \\ z \in \Gamma_v}} \Gamma_v.
\ee
Such a set $\Gamma_z$ is illustrated on Figure \ref{figDistZG} (right). For each $z\in \Gamma$ there exists $v\in \cV_0$ such that $\Gamma_z = \Gamma_v$, or there exists $v,v'\in \cV_0$, and $w,w'\in \cV_0$, such that $z\in [v,v']$ and 
$$
\Gamma_z= \Gamma_v \cup \Gamma_v' = \left[\frac{v+2w} 3,v\right]\cup [v,v'] \cup \left[v', \frac{v'+2w'} 3\right].
$$
Note that for any vertex $v\in \cV_0$ one has $\Gamma_v\cap \cV_0 = \{v\}$, hence \iref{defGammaz} agrees with \iref{defGammav} for any $z\in \cV_0$. Note also that for any $z,z'\in \Gamma$ the following are equivalent 
\begin{itemize}
\item $z'\in \Gamma_z$,
\item There exists $v\in \cV_0$ such that $\{z,z'\} \subset \Gamma_v$,
\item $z\in \Gamma_{z'}$.
\end{itemize}

The next proposition gives a robust estimate of the orientation of the set $\Gamma$.

\begin{lemma}
\label{lemmaPi3}
The following holds if $c$ is sufficiently small. Let $z\in \Gamma$, let $p,q\in \Gamma_z$ be two distinct points, and let $z'$ be such that $|z-z'| \leq 5 \lambda(z)$. Then
$$
\lhd(\theta(z')^\perp , p-q) \leq \cPi/3.
$$
\end{lemma}

\proof
We first remark that $\Gamma_z$ is not a singleton since the points $p$ and $q$ are distinct.
Therefore $\Gamma_z = \Gamma_v \cup \Gamma_{v'}$, or $\Gamma_z = \Gamma_v$, for some $v,v' \in \cV_0$ satisfying $\rho(v) \geq 3$ and $\rho(v') \geq 3$. We focus on the case of the union $\Gamma_v \cup \Gamma_{v'}$.

Since $|v-z|\leq \min \{\lambda(z), \lambda(v)\}$ we obtain
$$
|v-z'| \leq |v-z| + |z-z'| \leq \lambda(v)+ 5 \lambda(z) \leq \lambda(v)+5(\lambda(v)+ c|v-z|) \leq  \lambda(v) (6+5 c).
$$
If $c$ is sufficiently small, we thus have $|v-z'| \leq 10 \lambda(v)$ and therefore $\lhd(\theta(v), \theta(z') )\leq \cPi/12$ according to Lemma \ref{lemmaAngle12}.
Any of the two tangents $\bt$ to the polygonal line $\Gamma_v$ thus satisfies 
$$
\lhd(\theta(z')^\perp, \bt) \leq \lhd (\theta(v)^\perp, \bt)+ \lhd(\theta(v), \theta(z)) \leq 
 \cPi/4 + \cPi/12 = \cPi/3.
$$
Proceeding likewise we obtain that $\lhd(\theta(z')^\perp, \bt)\leq \cPi/3$ for any of the (three) tangents to the continuous polygonal line $\Gamma_z = \Gamma_v \cup \Gamma_{v'}$, and therefore $\lhd(\theta(z')^\perp, p-q)\leq\cPi/3$ for any $p,q$ on this line which concludes the proof.
\sq

The next proposition gives some upper and lower bounds on the distance from a point to the set $\Gamma$.
\begin{prop}
\label{propDistUpLow}
If $c$ is sufficiently small, then the following holds.
For any $z\in \Gamma$,
\be
\label{eqDistzGH0}
 \min\{|z-e| \sep  e\in \Gamma\sm \Gamma_z\} \geq \frac {10}{91} \lambda(z) 
\ee
and for any $z\in \R^2$
\be
\label{eqDistzGH}
\min\{\|z-e\|_{H(z)} \sep e \in \Gamma\} \leq 2 
\ee
\end{prop}

\begin{figure}
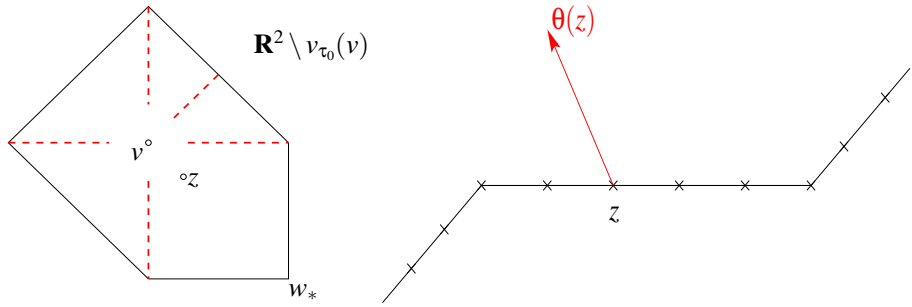

\centering
\includegraphics[width = 4.8cm, height = 4cm]{\pathPic/Dessins/Dessins_JM/DessinJM_04.pdf}
\includegraphics[width = 7cm, height = 4cm]{\pathPic/Dessins/Dessins_JM/DessinJM_12.pdf}
\caption{Illustration of the proof of Proposition \ref{propDistUpLow} \label{figDistZG}}
\end{figure}

\proof
We first establish \iref{eqDistzGH0}. 
Let $z\in \Gamma$ and let $v\in \cV_0$ be such that $z\in \Gamma_v$. There exists an edge $[v,w_*]$ of $\cT_0$ and a parameter $\alpha\in [0,2/3]$ such that $z= (1-\alpha) v+\alpha w_*$.
We denote by $V_{\cT_0}(v)$ the union of all triangles $T\in \cT_0$ containing $v$, and by $\cV_0(v)$ the collection of vertices $w$ such that $[v,w]$ is an edge of $\cT$. We thus have
$$
\Gamma\sm \Gamma_z \subset \overline {\R^2 \sm V_{\cT_0}(v)} \bigcup_{w\in \cV_0(v)\sm\{w_*\}} \left[\frac{2v+w}3, w\right].
$$
The set appearing in the right hand side is illustrated on Figure \ref{figDistZG} (left).
Denoting $d(z,E) := \inf\{ |z-e|\sep e\in E\}$ we obtain 
\begin{eqnarray*}
d(z,\Gamma\sm \Gamma_z) &\geq & \min \left\{ d(z, \R^2 \sm V_{\cT_0}(v)), \ \min_{w\in \cV_0(v) \sm \{w_*\} } d\left(z,\left[\frac{2v+w}3, w\right]\right)\right\}\\
& \geq & \min \left\{ \frac{|w_*-z|}{\sqrt 2}, \ \frac 1 {3\sqrt 2} \min_{w\in \cV_0(v)\sm \{w_*\}} |v-w|\right\}\\
& \geq & \frac 1 {3\sqrt 2} \min_{w\in \cV_0(v)} |v-w|\\
& \geq & \frac{1}{3\sqrt 2} \lambda(v)\frac {\sqrt 2} 3 \\
&\geq&  \frac{\lambda(z)-c|z-v|}{9} \geq \frac{\lambda(z)(1-c)}{9},
\end{eqnarray*}
where we used \iref{edgeLambda} in the fourth line, and the fact that $\lambda$ is $c$-Lipschitz in the last line. This establishes \iref{eqDistzGH0} if $c$ is sufficiently small.\\

We now turn to the proof of \iref{eqDistzGH}. 
Let $T \in \cT_0$ and let $z\in T$. Since $T$ is a half square there exists a vertex $v$ of $T$ such that $|z-v| \leq \diam(T)/2$, hence 
$$
\|z-v\|_{H(z)} \leq  \frac{|z-v|}{\mu(z)}   \leq \frac{\diam(T)}{2 \mu(z)} \leq \frac{\lambda(z)}{2 \mu(z)} = \frac{\rho(z)} 2,
$$
which concludes the proof if $\rho(z) \leq 4$. If $\rho(z) \geq 4$ and $c\leq c_\rho(3,1)$ then $\rho(v) \geq 3$ for any vertex $v$ of the triangle $T$.

We choose two vertices $v,w$ of $T$ such that the line $(z+ r \theta(z))_{r\in \sR}$ intersects $\partial T$ on the segment $[v,(v+w)/2]$.
We denote by $w_-,w_+\in \cV_0$ the neighbors of $v$ in $\cT_0$ which satisfy \iref{defWmWp}.
Our first step is to establish that the line $(z+ r \theta(z))_{r\in \sR}$ intersects the set $\Gamma_v$, which is  heuristically due to the fact that the angle $\lhd(\theta(z), \theta(v))$ is small and that the vertices $w_-$ and $w_+$ are the \emph{furthest away} from $v$ in the direction of $\theta(v)^\perp$ as illustrated on Figure \ref{figWmWp} (right).

It follows from \iref{eqEdgeRatio} an Lemma \ref{lemmaAngleWmWp} that
\be
\label{eqScalLower}
\<\theta(v)^\perp,w_+ - v\> \geq |w_+-v|\cos(\cPi/4) \geq \frac{|w-v|} 2 \frac 1 {\sqrt 2}. 
\ee
We assume without loss of generality that 
$\<\theta(z), \theta(v)\> \geq 0$ which implies according to Lemma \ref{lemmaAngle12}
$$
|\theta(z)-\theta(v)| \leq \varangle(\theta(z), \theta(v)) = \lhd(\theta(z), \theta(v)) \leq \delta := 1/(10 \sqrt 2).
$$
It follows that 
\begin{eqnarray*}
\<\theta(z)^\perp, w-v\> &=& \<\theta(z)^\perp-\theta(v)^\perp, w-v\>+\<\theta(v)^\perp, w-v\>\\
& \leq & |\theta(z)-\theta(v)| |w-v| + \<\theta(v)^\perp, w_+-v\>\\
& \leq & (1+2\delta \sqrt 2 )\<\theta(v)^\perp, w_+-v\>,
\end{eqnarray*}
where we injected \iref{defWmWp} in the second line, and \iref{eqScalLower} in the third line.
Denoting $\Theta := \varangle(\theta(v)^\perp, w_+-v)\leq \cPi/4$ and $\Phi := \varangle(\theta(z)^\perp, w_+-v)$ we obtain 
$$
\frac{\<\theta(z)^\perp, w_+-v\>}{\<\theta(v)^\perp, w_+-v\>} = \frac{\cos(\Phi)}{\cos(\Theta)} \geq \frac{\cos(\Theta) - |\Phi-\Theta|}{\cos(\Theta)} =1-\frac {|\Phi-\Theta|}{\cos\Theta}\geq 1-\delta \sqrt 2, 
$$
hence 
$$
\<\theta(z)^\perp, w-v\>  \leq \frac{1+2\delta \sqrt 2}{1-\delta \sqrt 2}\<\theta(z)^\perp, w_+-v\> = \frac 4 3 \<\theta(z)^\perp, w_+-v\>.
$$
Proceeding likewise for $w_-$ we obtain 
%
$$
\frac 4 3\<\theta(z)^\perp, w_- -v\> \leq \<\theta(z)^\perp, w-v\> \leq  \frac 4 3\<\theta(z)^\perp, w_+ - v\>,
$$
and therefore since $z+r\theta(z)\in [v,(v+w)/2]$ for some $r\in \R$,
$$
\frac 2 3\<\theta(z)^\perp, w_- -v\> \leq \<\theta(z)^\perp, z-v\> \leq  \frac 2 3\<\theta(z)^\perp, w_+ - v\>.
$$
It follows that the line $(z+r\theta(z))_{r\in \sR}$ intersects the set $\Gamma_v := [v, (v+2w_+)/3] \cup  [v, (v+2w_-)/3]$ at some point $z'$.
The point $z'$ belongs to a triangle $T'$ containing $v$, hence
\begin{eqnarray*}
\lambda(z) \|z'-z\|_{H(z)} &=& \lambda(z)\|r\theta(z)\|_{H(z)}\\
&=& |r|\\
&=& |z-z'|\\
&\leq& |z-v|+|v-z'|\\
&\leq& \diam(T)+ \frac 2 3\diam(T')\\
&\leq & \left(1+\frac 2 3 \sqrt 2\right) \diam(T)
 \leq  \left(1+\frac 2 3 \sqrt 2\right) \lambda(z),
\end{eqnarray*}
where we used \iref{eqDiamRatio} in the last line. This concludes the proof of this proposition.
\sq

We construct in the next lemma a collection $\cV$ of sites by sampling the metric space $(\Gamma, d_H)$, which is the collection of vertices of our future quasi-acute triangulation $\cT$. This set is illustrated by small crosses on $\Gamma_z$ on Figure \ref{figDistZG} (right).

\begin{lemma}
\label{lemmaConstV}
The following holds if $c$ is sufficiently small.
There exists a discrete set $\cV \subset\Gamma$ containing $\cV_0$ and satisfying the following:
\begin{enumerate}[i)]
\item (Separation) For all $v,w\in \cV$, $d_H(v,w) \geq 1/10$.
\item (Distance from $\Gamma$) For all $z\in \Gamma$, $d_H(z, \cV) \leq 1/10$.
\item (Distance from $\R^2$) For all $z\in \R^2$, $d_H(z,\cV) \leq 2+1/4$.
\end{enumerate}
\end{lemma}
\proof
For any distinct $v,w\in \cV_0$ we obtain using \iref{pairwiseV0} and Proposition \ref{propEuclidRiemann}
\begin{eqnarray*}
d_H(v,w) &\geq& \ln(1+c\|v-w\|_{H(v)})/c\\
&\geq& \ln(1+c|v-w|/\lambda(v))/c\\
&\geq& \ln(1+c\sqrt 2 /3)/c 
\end{eqnarray*}
which is larger than $1/10$ if $c$ is sufficiently small. 

We denote by $\cV$ a $1$-equispaced subset of the metric space $(\Gamma, 10\:d_H)$ containing $\cV_0$, which exists according to Lemma \ref{lemmaEquiSpaced}. The two first announced properties, separation and distance from $\Gamma$, follow directly from this construction.

The third property is obtained as follows: for any $z\in \R^2$ there exists according to Proposition \ref{propDistUpLow} a point $e\in \Gamma$ such that $\|z-e\|_{H(z)}\leq 2$. Furthermore there exists a vertex $v\in \cV$ such that $d_H(e,v) \leq 1/10$, hence
$$
d_H(z,\cV) \leq d_H(z,e)+ d_H(e,v) \leq -\ln(1-c\|z-e\|_{H(z)})/c + 1/10 \leq -\ln(1-2c)/c + 1/10.
$$
which is smaller than $2+1/4$ if $c$ is sufficiently small.
\sq

The next lemma introduces the anisotropic Delaunay triangulation $\cT$ associated to the metric $H$ and the set of points $\cV$. 

\begin{lemma}
\label{lemmaRecallG}
There exists an absolute constant $C_0\geq 1$ such that the following holds if $c$ is sufficiently small.
The set $\cV$ is wedged with respect to the metric $H$, and the anisotropic Delaunay triangulation $\cT$ obtained by arbitrarily triangulating the cells of the graph $\cG(\cV,H)$ is $C_0$-adapted to $H$ and belongs to $\bT_{g,C_0}$.
\end{lemma}

\proof
We denote $\alpha := (2+1/4)^{-1}$ throughout the proof of this lemma.
It follows from to Lemma \ref{lemmaConstV} that the set $\cV$ is a $\alpha/10$-equispaced subset of the metric space $(\R^2, \alpha d_H)$. 

Theorem \ref{thEquiMesh} applied to the metric $\alpha^2 H$, therefore implies that, if $c$ is sufficiently small, the set 
$\cV$ is wedged for the metric $\alpha^2 H$ and that arbitrarily triangulating the convex cells of the graph $\cG(\cV,\alpha^2 H)$, which are convex yields a triangulation $\cT$ which is $C$-adapted to the metric $\alpha^2H$ (hence $C_0 := C/\alpha$-adapted to $H$), where $C$ is an absolute constant.

It immediately follows from the definition \iref{defGVH} that the graphs $\cG(\cV,\alpha^2 H)$ and $\cG(\cV,H)$ are equal, and Definition \ref{defWedge} shows that the wedge property holds for $H$ if and only if it holds for $\alpha^2 H$. 
\sq

The key ingredient of the construction of the anisotropic Delaunay triangulation is the anisotropic Voronoi diagram introduced in \iref{defVorv}.
Our first lemma compares the Voronoi regions with some balls. We recall that $B_H(z,r) := \{z'\in \R^2 \sep \|z'-z\|\leq r\}$, and $B(z,r)$ denotes the usual euclidean ball of radius $r$ centered at $z$.

\begin{lemma}
\label{lemmaVorBall}
If $c$ is sufficiently small, then the following holds.
For any $v\in \cV$
\begin{enumerate}
\item $\overline B_H(v,1/25)  \subset \Vor(v)$
\item $B_H(v,2+1/3)  \supset \Vor(v)$
\item $B(v, 4\mu(v)/19) \supset  \Vor(v) \cap \Gamma = \Vor(v) \cap \Gamma_v$
\end{enumerate}
\end{lemma}

\proof
We first establish Point 1., and for that purpose we consider $z\in \partial \Vor(v)$, hence $z\in \Vor(v) \cap \Vor(w)$ for some vertex $w\in \cV$.
We define $r:=\|z-v\|_{H(v)} = \|z-w\|_{H(w)}$ and we obtain using Point i) of Lemma \ref{lemmaConstV}
\begin{eqnarray*}
1/10 &\leq& d_H(v,w)\\
&\leq& d_H(z,v)+d_H(z,w)\\
&\leq& -\ln (1-c\|z-v\|_{H(v)})/c - \ln (1-c\|z-v\|_{H(v)})/c\\
&=& -2 \ln (1-cr)/c.
\end{eqnarray*}
It follows that 
$$
r\geq (1-\exp(-c/20))/c,
$$
and therefore $r\geq 1/25$ if $c$ is sufficiently small, 
which concludes the proof of Point 1.

We now turn to Point 2, and for that purpose we remark that for any $z\in \R^2$ there exists according to Point iii) of Lemma \ref{lemmaConstV} a vertex $w\in \cV$ such that $d_H(z,w)\leq 2+1/4$. It follows that 
$$
\|z-w\|_{H(w)} \leq (\exp(c d_H(z,w))-1)/c \leq (\exp(c(2+1/4))-1)/c.
$$
Hence $\|z-w\|_{H(w)} \leq 2+1/3$ if $c$ is sufficiently small. It thus follows from the definition \iref{defVorv} of the Voronoi diagram that $\|z-v\|_{H(v)} \leq 2+1/3$ for all $z\in \Vor(v)$, which concludes the proof of Point 2.

We now turn to Point 3, and for that purpose we remark that for all $z\in \Gamma$ there exists according to Point ii) of Lemma \ref{lemmaConstV} a vertex $w\in \cV$ such that $d_H(z,w) \leq 1/10$. It follows that 
$$
\|z-w\|_{H(w)} \leq (\exp(c d_H(z,w))-1)/c \leq (\exp(c/10)-1)/c.
$$
Hence $\|z-w\|_{H(w)} \leq 2/19$ if $c$ is sufficiently small. It thus follows from the definition of the Voronoi diagram that $\|z-v\|_{H(v)} \leq 2/19$ for all $z\in \Vor(v) \cap \Gamma$. Therefore  
$$
\frac {|z-v|}{\lambda(v)} \leq  \|z-v\|_{H(v)}  \leq \frac 2 {19} < \frac{10}{91} 
$$
which implies that $z\in \Gamma_v$ according to Proposition \ref{propDistUpLow}.
We thus have $\lhd(\theta(v)^\perp , z-v) \leq \cPi/3$, according to Lemma \ref{lemmaPi3}, and therefore 
$$
\frac{|z-v|}{\mu(v)} \cos(\cPi/3) \leq  \|z-v\|_{H(v)}  \leq \frac 2 {19}
$$
which concludes the proof of this proposition.
\sq

We recall that, by definition of the anisotropic Delaunay triangulation, the Voronoi regions $\Vor(v)$ and $\Vor(w)$ intersect for any edge $[v,w]$ of $\cT$.
We say that an edge $[v,w]$ of $\cT$ is transverse if $\Vor(v) \cap \Vor(w) \cap \Gamma \neq \emptyset$, and we say that $[v,w]$ is aligned otherwise.

Our next intermediate result estimates the length of transverse and aligned edges.

\begin{lemma}
\label{lemmaLengthEdges}
Let $[v,w]$ be an edge of $\cT$ and let $z\in [v,w]$.
\begin{enumerate}
\item If $[v,w]$ is transverse, then $|v-w|\leq \mu(z)/2$. More precisely $\max\{|v-x|,|w-x|\} \leq \mu(z)/4$ for all $x\in \Vor(v) \cap \Vor(w) \cap \Gamma$.
\item In any case $\|v-w\|_{H(z)} \leq 5$.
\end{enumerate}
\end{lemma}

\proof
We first establish Point 1., hence we assume that $[v,w]$ is transverse and we consider a point $x\in \Vor(v) \cap \Vor(w) \cap \Gamma$.
For notational simplicity we denote $\beta := 4/19$.
It follows from Lemma \ref{lemmaVorBall} that 
$
|v-x| \leq \beta  \mu(v) \leq \beta (\mu(x)+ c|v-x|),
$
hence 
$
|v-x| \leq \beta (1-c\beta)^{-1}\mu(x).
$
Likewise $|w-x|\leq \beta (1-c\beta)^{-1}\mu(x)$. Furthermore $|z-x| \leq \max\{|v-x|, |w-x|\}$ since $z\in [v,w]$, hence
$$
\mu(x) \leq \mu(z) + c|z-x| \leq \mu(z) + c\beta  (1-c \beta)^{-1}\mu(x)
$$
which implies 
$
\mu(x) \leq  (1-c\beta (1-c \beta)^{-1})^{-1} \mu(z).
$
Therefore 
$$
|v-x| \leq \beta (1-c\beta)^{-1} \mu(x) \leq \beta (1-c\beta)^{-1} (1-c\beta (1-c \beta)^{-1})^{-1} \mu(z).
$$
The term in front of $\mu(z)$ in the right hand side tends to $\beta = 4/19<1/4$ as $c\to 0$. If $c$ is sufficiently small we therefore obtain as announced $|v-x| \leq \mu(z)/4$. Likewise $|w-x| \leq \mu(z)/4$ and therefore $|v-w|\leq \mu(z)/2$.

We now turn to the proof of the second point. 
We have according to Lemma \ref{lemmaDiamVor} 
$$
\alpha d_H(v,w) = d_{\alpha^2 H} (v,w) \leq -2\ln (2-e^{c'})/c'
$$
where $c' := c/\alpha$. Hence 
$$
d_H(v,w) \leq r_0(c) := -2 \ln (2-e^{c/\alpha})/c.
$$
If follows that 
$
\|v-w\|_{H(v)} \leq r_1(c) := (e^{c r_0(c)}-1)/c,
$
which implies using  \iref{eqLocalNorm} that
$$
\|v-w\|_{H(z)} \leq \|v-w\|_{H(v)} (1-c \|v-z\|_{H(v)})^{-1} \leq r_2(c) := r_1(c) (1-c r_1(c))^{-1}
$$
for any $z\in [v,w]$.
As $c\to 0$ the functions $r_0$, $r_1$ and $r_2$ all tend to $2 (1+1/4) <5$. Choosing $c$ sufficiently small we  conclude the proof of this proposition.
\sq

The next lemma describes the orientations of the transverse and aligned edges of $\cT$, in regions of sufficient anisotropy.
We introduce the constant 
$$
\rho_0 := \max\{50/\sin(\cPi/12), \, 11C_0\},
$$
where $C_0$  is the constant from Lemma \ref{lemmaRecallG}  (we only use in the sequel that $\cT$ is $C_0$ equivalent to the metric $H$).

For any $v,w\in \cV$, we say that $[v,w]$ is an edge of $\Gamma$ if 
$$
(v,w) \subset \Gamma\sm \cV,
$$
where $(v,w)$ denotes the relative interior of the segment $[v,w]$.
The next lemma characterizes the transverse edges of $\cT$ in regions where the anisotropy is sufficiently pronounced.

For any $x\in \Gamma$
$$
\Gamma_x^+ := \{v'\in \Gamma_v \sep \<\theta(x)^\perp ,v'-x\> > 0\},
$$
and we define $\Gamma_x^-$ similarly.
If $[v,w]$ is an edge of $\Gamma$, then one easily checks that $w$ is the closest element to $v$ in $\Gamma_v^+\cap \cV$ or $\Gamma_v^-\cap \cV$.

\begin{lemma}
\label{lemmaEdgeGamma}
The following holds if $c$ is sufficiently small.
Let $v,w\in \cV$ be a transverse edge, and assume that there exists $z\in [v,w]$ such that $\rho(z) \geq \rho_0$.
Then $[v,w]$ is and edge of $\Gamma$, $\Vor(v) \cap \Vor(w) \cap \Gamma$ is a singleton $\{x\}$, and $x\in [v,w]$.

\end{lemma}

\proof
We first consider an transverse edge $[v,w]$, a point $z\in [v,w]$ such that $\rho(z) \geq \rho_0$, and a point $x\in \Vor(v) \cap \Vor(w) \cap \Gamma$.
According to Lemma \ref{lemmaVorBall} we have $x\in \Gamma_v \cap \Gamma_w$, hence $v,w\in \Gamma_x$.

We may assume without loss of generality that $v\in \Gamma_x^+$. We denote by $v'$ the point of $\Gamma_x^+$ which is the closest to $x$, and our first objective is to to establish that $v=v'$. We assume for contradiction that this is not the case.

It follows from Lemma \ref{lemmaPi3} that the segment $[x,v]$ and the line $v'+ \theta(v')\R$ intersect. Hence there exists $\alpha\in [0,1]$ and $r\in \R$ such that 
$$
x+\alpha(v-x) = v'+ r \theta(v')
$$
Let $u$ be a unit vector orthogonal to $v-x$. We have $\<x,u\> = \<v',u\> + r \<\theta(v'),u\>$, and therefore using Lemma \ref{lemmaPi3},
\begin{eqnarray*}
|r|/2 & \leq & |r| \cos \lhd(\theta(v')^\perp , v-x)\\
& = & |r \<\theta(v'),u\>|\\
&=& |\<x-v', u\>|\\
& \leq & |x-v'| \leq |x-v| \leq \mu(z)/4
\end{eqnarray*}
where we used and Point 2. of Lemma \ref{lemmaVorBall} in the last.

According to Lemma 3 (Visibility Lemma) in \cite{Shew} the Voronoi regions of a wedged Voronoi diagram are star-shaped. 
It thus follows from Lemma \ref{lemmaVorBall} that the segment $[x,v]$ does not intersect $B_H(v',1/25)$, hence $|r|\geq \lambda(v')/25$ which yields
$$
\lambda(z) - c|z-v'| \leq \lambda(v') \leq \mu(z)25/2.
$$
It follows from Lemma \ref{lemmaLengthEdges} that 
$$
|z-v'| \leq |z-x|+|x-v'| \leq \max\{|v-x|, |w-x|\}+ |x-v|\leq \mu(z)/2.
$$
We thus obtain $25/2+ c/2 \geq \rho(z)\geq \rho_0$ which is a contradiction since $c\leq 1$. This contradiction is illustrated on Figure \ref{figContra} (left).

Hence $v=v'$ is the point of $\Gamma_x^+$ the closest from $x$.
If $w\in \Gamma_x^+$, then we obtain $w=v' = v$ which is a contradiction. Hence $w\in \Gamma_x^-$, and reasoning similarly we find that $w$ is the point of $\Gamma_x^-$ the closest to $x$.
This implies that $[v,w]$ is the unique edge of $\Gamma$ containing $x$ (note that the only points of $\Gamma$ which belong to two edges of $\Gamma$ are some vertices $v\in \cV$. But $x\notin \cV$ since any vertex $v\in \cV$ belongs to the interior of its own Voronoi region). 

Assume for contradiction that there exists another point $x'\in \Vor(v) \cap \Vor(w) \cap \Gamma$. The above argument shows that $[v,w]$ is the edge of $\Gamma$ containing $x'$. Therefore $x,x'\in [v,w]$, but this implies $x=x'$ since the Voronoi regions are star shaped.
\sq

\begin{figure}
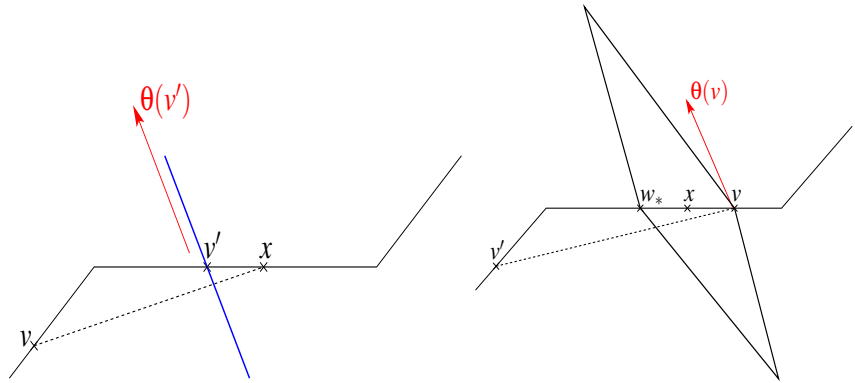

\centering
\includegraphics[width=6cm,height=4cm]{\pathPic/Dessins/Dessins_JM/DessinJM_13.pdf}
\includegraphics[width=5cm,height=5cm]{\pathPic/Dessins/Dessins_JM/DessinJM_14.pdf}
\caption{Illustration of the proof of Lemmas \ref{lemmaEdgeGamma} (left) and \ref{lemmaAligned} (right)\label{figContra}} 
\end{figure}

Our next purpose is to estimate the orientations of the edges of $\cT$.
\begin{lemma}
\label{lemmaNotInGammav}
The following holds if $c$ is sufficiently small.
Let $[v,w]$ be an edge of $\cT$ such that $w\notin \Gamma_v$ and assume that $\rho(z) \geq \rho_0$ for some $z\in [v,w]$. Then 
$$
\lhd(\theta(z), w-v) \leq \cPi/12.
$$
\end{lemma} 

\proof
It follows from Lemma \ref{lemmaRecallG} that 
$$
 \frac{|w-v|}{\mu(z)} \sin \lhd(\theta(z) , w-v) \leq \|w-v\|_{H(z)} \leq 5.
$$
On the other hand Proposition \ref{propDistUpLow} yields since $w\notin \Gamma_v$
$$
|v-w| \geq \lambda(v) 10/91 \geq (\lambda(z) - c|z-w|) 10/91,
$$
therefore $|v-w|/\lambda(z) \geq 10/(91+10c/91)$ which implies that $|v-w|\geq \lambda(z)/10$ if $c$ is sufficiently small.
It follows that 
$$
\sin \lhd(\theta(z) , w-v) \leq \frac{50}{\rho(z)},
$$
which concludes the proof since $\rho_0 \geq 50 / \sin (\cPi/12)$.
\sq

The next lemma establishes that the aligned edges of $\cT$ are, as their name indicates, aligned with the direction $\theta$ of anisotropy in regions where this anisotropy is sufficiently pronounced.

\begin{lemma}
\label{lemmaAligned}
The following holds if $c$ is sufficiently small.
\begin{enumerate}
\item
For any vertex $v\in \cV$ such that $\rho(v) \geq \rho_0$, there exists two aligned edges $[v,w]$, $[v,w']$ of $\cT$ such that 
\be
\label{eqDisjAngles}
\varangle(\theta(v), w-v) \leq \cPi/12 \stext{ and } \varangle(\theta(v), v-w') \leq \cPi/12
\ee
\item
Consider an aligned edge $[v,v']$ of $\cT$ and a point $z\in [v,v']$ such that $\rho(z) \geq \rho_0+1$.
Then $\lhd(\theta(z),w-v) \leq \cPi/12$.
\end{enumerate}
\end{lemma}

\proof
We first establish Point 1., and as an intermediate objective we prove that $\Gamma_v\sm \Vor(v) \neq \emptyset$. It follows from the definition \iref{defGammaz} of $\Gamma_v$ that there exists an edge $[v_0,w_0]$ of $\cT_0$ such that $v\in [v_0, (v_0+2w_0)/3] \subset \Gamma_v$.
Therefore according to \iref{edgeLambda}
$$
\max \{|v-v_0|, |v-(v_0+2w_0)/3| \} \geq |v_0-w_0|/3 \geq \lambda(v)\sqrt 2/9.
$$
On the other hand $|z-v| \leq 4\mu(v)/19$ for all $z\in \Vor(v)\cap \Gamma_v$ according to Lemma \ref{lemmaVorBall}. Since $\rho(v) \geq \rho_0 > 36/(19\sqrt 2)$ we find that one of the points $v_0$ or $w_0$ does not belong to $\Vor(v)$.

Consider a point $x\in \Gamma_v \cap \partial \Vor(v)$. There exists vertex $w_*\in \cV$ such that $x\in \Vor(w_*)$ and $[v,w_*]$ is an edge of $\cT$. 
Therefore $[v,w_*]$ is a transverse edge of $\cT$ such that $x\in \Vor(v) \cap \Vor(w_*) \cap \Gamma$ and $\rho(v) \geq \rho_0$. This implies according to Lemma \ref{lemmaEdgeGamma} that $[v,w_*]$ is an edge of $\Gamma$ containing $x$.

Let $T,T'\in \cT$ be the two triangles containing the edge $[v,w_*]$, and let $w,w'\in \cV$ be the third vertex of these triangles respectively.
Since $\cT$ is $C_0$-adapted to the metric $H$, we obtain using \iref{eqDiamTHT2}
\be
\label{eqDiamPlus}
|v-w_*|+|v-w| > \diam(T) \geq \|\cH_T^{-\frac 1 2}\| \geq \|H(v)^{-\frac 1 2}\|/C_0 = \lambda(v)/C_0.
\ee
If $w\in\Gamma_v$, then we obtain using Lemma \ref{lemmaPi3} and Lemma \ref{lemmaLengthEdges} that
$$
|v-w| /(2 \mu(v)) \leq |v-w| \cos(\lhd(\theta(v)^\perp,w-v))/\mu(v) \leq \|v-w\|_{H(v)} \leq 5,
$$
and $|v-w_*| \leq \mu(v)/2$. Injecting this in \iref{eqDiamPlus} we obtain 
$
\mu(v) (10+1/2) \leq \lambda(v)/C_0
$
hence 
$
\rho(v) \leq (10+1/2) C_0
$
which is a contradiction.
Therefore $w\notin \Gamma_v$, which implies $\lhd(\theta(v), w-v) \leq \cPi/12$ according to Lemma \ref{lemmaNotInGammav}. Hence
$$
\varangle(\theta(v), w-v)\in [0,\cPi/12]\cup[11\cPi/12, \cPi].
$$
and likewise $\varangle(\theta(v), w'-v)\in [0,\cPi/12]\cup[11\cPi/12, \cPi]$.
On the other hand we have $\lhd (\theta(v)^\perp, w_*-v) \leq \cPi/3$ according to Lemma \ref{lemmaPi3}, and therefore 
$\varangle (\theta(v), w_*-v) \in [\cPi/6, 5\cPi/6]$. Since the triangles $T$ and $T'$ are on opposite sides of the edge $[v,w_*]$ we obtain  \iref{eqDisjAngles} which concludes the proof of Point 1.\:

We now turn to the proof of Point 2.
We have $|v-z|/\lambda(v) \leq \|v-w\|_{H(z)} \leq 5$. Assuming that $c\leq c_\rho(\rho_0,5)$ we therefore obtain $\rho(v)\geq \rho_0$.
If $v'\notin \Gamma_v$, then Lemma  \ref{lemmaNotInGammav} gives the announced result.
We therefore assume for contradiction that $v'\in \Gamma_v$, and we may assume without loss of generality that $v'\in \Gamma_v^+$. It follows that $\Gamma_v^+\sm \Vor(v) \neq \emptyset$. Repeating our previous argument we consider a point $x\in \Gamma_v^+ \cap \partial \Vor(v)$, and we find that there exists a vertex $w_* \in \Gamma_v^+$ such that $[v,w_*]$ is an edge common to $\cT$ and $\Gamma$ containing $x$. We denote by $T,T'$  the two triangles containing this edge and by $w,w'$ the third vertex of these triangles.

Since $v'\in \Gamma_v^+$ we have $\lhd(\theta(v)^\perp, v-v') \leq \cPi/3$ according to Lemma \ref{lemmaPi3}, hence  $\varangle (\theta(v), v'-v) \in [\cPi/6, 5\cPi/6]$. Since $w,w'$ satisfy \iref{eqDisjAngles} we obtain that $[v,v']$ intersects one of the triangles $T$ or $T'$, which is a contradiction, as illustrated on Figure \ref{figContra} (right). This concludes the proof of this lemma.
\sq


For each $v\in \cV_0$ such that $\rho(v) \geq 2+\rho_0$ we define $\Gamma'_v := [v,w_-] \cup[v,w_+]$ where $w_-,w_+\in \cV_0$ are such that $\Gamma_v := [v,(2w_-+v)/3] \cup[v,(2w_++v)/3]$.
We denote by $\Gamma'$ the union of these sets
$$
\Gamma' := \bigcup_{\substack{v\in \cV_0\\ \rho(v) \geq 2+\rho_0}} \Gamma'_v,
$$
which is a union of segments joining some vertices in $\cV_0$, hence in $\cV$ according to Lemma \ref{lemmaConstV}.
We denote by $\cQ$ the partition into convex polygons of $\R^2$ obtained by bisecting the triangles $T\in \cT$ by these segments. Such a partition is generally referred to as the overlay  $\cT$ and $\Gamma'$.

\begin{prop}
The following holds if $c$ is sufficiently small.
Denote by $\cT'$ the triangulation obtained by arbitrarily triangulating the convex cells of $\cQ$.
Then $\cT'$ is a $300$-refinement of $\cT$ and $S(T') \leq \max\{C_0^2(3+\rho_0), \tan (11\cPi/24)\}$ for all $T'\in \cT'$.
\end{prop}

\proof
Our first objective is to give a uniform bound on the measure of sliverness in $\cT'$, and for that purpose we denote by $\cV'$ the collection of vertices of $\cT'$.
We have by construction $\cV' \subset \Gamma \cup \Gamma'$.

Assuming that $c\leq c_\rho(\rho_0+1, 1)$, we obtain that $\rho(z) \geq \rho_0+1$ for all $z\in \Gamma'$.

Consider a vertex $z\in \cV'$, and assume in a first time that $z\in \Gamma'$.
It follows from Lemma \ref{lemmaPi3} that there exists two edges $[z,v]$, $[z,v']$ of $\cT'$, contained in $\Gamma'$, and such that 
$$
\varangle(\theta(z)^\perp , v-z) \leq \cPi/3 \stext{ and } \varangle(\theta(z)^\perp , z-v') \leq \cPi/3.
$$
Furthermore, it follows from Lemma \ref{lemmaAligned} that there exists two edges $[z,w]$, $[z,w']$, of $\cT'$ such that 
$$
\varangle(\theta(v), w-z) \leq \cPi/12 \stext{ and } \varangle(\theta(v), z-w') \leq \cPi/12.
$$
In more detail two cases are possible : either $z$ is a vertex of the original triangulation $\cT$, and Point 1. of Lemma \ref{lemmaAligned} applies, or $z$ is a vertex of $\cT'$ created by the overlay of $\Gamma$ and $\cT$. In that case Point 2. of Lemma \ref{lemmaAligned} applies and the edges $[z,w]$, $[z,w']$, of $\cT'$ are contained in the same aligned edge of $\cT$.

We thus have 
\begin{eqnarray*}
\varangle(w-z,v-z) &\leq& \varangle(w-z,\theta(z)) + \varangle (\theta(z), \theta(z)^\perp)+ \varangle(\theta(z)^\perp, v-z)\\
& \leq & \cPi/12+ \cPi/2+ \cPi/3\\
& = & 11 \cPi/12.
\end{eqnarray*}
Likewise 
$$
\varangle(v-z, w'-z)  \leq 11 \cPi/12, \quad \varangle(w'-z, v'-z) \leq 11 \cPi/12 \stext{ and } \varangle(v'-z, w-z) \leq 11 \cPi/12,
$$
which immediately implies that any angle at the vertex $z$ is bounded by $11 \cPi/12$.

In a second time we consider a vertex $z\in \cV'\sm \Gamma'$. Since $\cV'\subset \Gamma \cup \Gamma'$ there exists $v\in \cV_0$ such that $z\in \Gamma_v$. If $\rho(z) \geq \rho_0+3$ and $c \leq c_\rho(\rho_0+3,1)$ then $\rho(v) \geq \rho_0+2$ and therefore $z\in \Gamma'$, which is a contradiction. Therefore $\rho(z) \leq \rho_0+3$, which implies that the triangles $T\in \cT$ containing $z$ satisfy $S(T) \leq \rho(T) \leq C_0^2 \rho(z) \leq C_0^2 (\rho_0+3)$. Any angle $\theta$ at the vertex $z$ in the triangulation $\cT$, hence also in the triangulation $\cT'$, thus satisfies $\tan (\theta/2) \leq C_0^2 (\rho_0+3)$, which concludes the proof of the upper bound of the measure of sliverness $S$ on $\cT'$.

We now prove that $\cT'$ is a bounded refinement of $\cT$, and for that purpose we consider an edge $e=[v,w]$ of $\cT$, and we denote $z := (v+w)/2$. 
We have by construction $|v-w|/\lambda(z) \leq \|v-w\|_{H(z)} \leq 5$.
If a triangle $T\in \cT_0$ intersects $[v,w]$, then for any $z'\in [v,w] \cap T$
$$
\lambda(z)(1-5c/2)\leq \lambda(z)- c|z-z'| \leq \lambda (z') \leq \lambda(z)+ c|z-z'| \leq  \lambda(z)(1+5c/2).
$$
Recalling that  $2\lambda(z')/3 \leq \diam(T)\leq \lambda(z')$ we obtain
$$
T \subset  B(z, 5 \lambda(z)/2 + \lambda(z')) \subset B(z, \lambda(z) (5/2+1+c5/2)), 
$$
and in the other hand since $T$ is a half square
$$
\sqrt{2 |T|} = \diam(T) \geq  2\lambda(z')/3 \geq 2\lambda(z)(1-5c/2)/3.
$$
Comparing the areas we obtain that the number of triangles $T\in \cT_0$ intersecting $[v,w]$ is bounded by 
$$
\cPi \left(\frac{ 5/2+1+5c/2}{ 2(1-5c/2)/3 }\right)^2
$$
which is smaller than $100$ if $c$ is sufficiently small.

It follows that any edge $[v,w]$ of $\cT$ is cut in at most $100$ parts by the set $\Gamma$.
For any $T\in \cT$ we thus obtain $\#(\partial T \cap \cV') \leq 3 \times 100 = 300$.
Since any conforming triangulation of the convex envelope of $n+2$ points uses $n$ triangles, we obtain that the triangulation $\cT'$ is a $298$ refinement of $\cT$ which concludes the proof of this proposition.
\sq


\chapter{Approximation theory based on metrics}
\minitoc
\label{chapApproxMet}
\section{Introduction} 
\label{secPer}

The main purpose of adaptive mesh generation, compared uniform mesh generation, is to reduce the number of simplices required to achieve a given task. Here are two relevant examples which 
are further discussed in this chapter:
\begin{enumerate}
\item
Given two domains of interest, generate the mesh with smallest possible cardinality that separates these
domains with a layer of simplices.
\item
Given a function $f$, generate the mesh of a prescribed cardinality $N$ that yields the best 
finite element approximation to $f$ of a given order $m-1$ in some prescribed norm $L^p$ or $W^{1,p}$.
\end{enumerate}
Both examples can be viewed as optimization problems posed on the set of meshes.
The second problem was dealt with in Chapter 2 and 3 in the case where $f$ is smooth, 
and an optimal mesh was constructed in the asymptotic regime $N\to +\infty$. The results
of the previous chapter have revealed that certain relevant classes of meshes can be equivalently
described in terms of riemannian metrics, and a natural objective is therefore
to translate the above problems as {\it optimization problems posed on the set of metrics}.
Optimization problems posed on metrics may indeed be easier to solve practically  
due to the continuous nature of these objects.
The goal of this chapter is to explain how this objective can be met for the
above two examples. The dimension $d\geq 2$ is fixed throughout this chapter, we denote by $\bT$ the collection of conforming simplicial meshes of $\R^d$ and by $\bH := C^0(\R^d, S_d^+)$ the collection of metrics.

In the previous chapter, we have introduced the classes of meshes
$$
\bT_{i,C} \subset \bT_{a,C} \subset \bT_{g,C},
$$
and characterized the respectively equivalent classes of metrics
$$
\bH_i\subset \bH_a\subset \bH_g.
$$
These classes of meshes and metrics are defined on the entire unbounded domain $\R^d$.
The above optimization problems do not make sense in this setting since the 
number of elements in such meshes is infinite. In order to circumvent this difficulty
while still avoiding the difficulties related to the boundary of a domain, we 
work in a periodized setting that we describe below.

\subsection{Periodic meshes and metrics}

We denote by $\bT_\per$ the collection of meshes $T\in \bT$ which satisfy the following properties:
\begin{enumerate}[i)]
\item (Translation invariance) For any $T\in \cT$ and any $u \in \Z^d$, one has $u+T \in \Z^d$.
\item (Constrained diameters) For any $T\in \cT$ one has $\cH_T\geq \Id$.
\end{enumerate}

Accordingly, we denote by $\bH_\per$ the collection of metrics $H \in \bH$ which satisfy the following properties:
\begin{enumerate}[i)]
\item (Translation invariance) For any $z\in \R^d$ and any $u \in \Z^d$, one has $H(z) = H(z+u)$.
\item (Constrained positivity) For any $z\in \R^d$ one has $H(z)\geq \Id$.
\end{enumerate}

A mesh $\cT\in \bT_\per$ can be regarded, thanks to the translation invariance property, 
as a mesh of the compact periodic space 
\be
\label{perSpace}
\Pi^d := (\R/ \ZZ)^d,
\ee
which is also called the $d$-dimensional torus. 
The diameter constraint property ensures that the simplices $T$ of a periodic mesh $\cT \in \bT_\per$ do not have a significantly larger size than the fundamental cell $[0,1]^d$ of the standard $\Z^d$ periodic tiling of $\R^d$. 

Strictly speaking, a periodic mesh $\cT\in \bT_\per$ has an infinite cardinality. However it only contains a finite number of simplices up to translation by an element of $\Z^d$.  
For any $\cT\in \bT_\per$, we therefore denote by $\#(\cT)$ the number of equivalence classes in $\cT$ for the relation
$$
T \sim T' \text{ if and only if there exists } u \in \Z^d \text{ such that } \ T' = T+u.
$$
The number $\#(\cT)$ is finite, and is the the cardinality of $\cT$ seen as a mesh of the compact periodic space \iref{perSpace}.
We define the {\it mass} $m(H)$ of a periodic metric $H \in \bH_\per$ as follows
$$
m(H) :=  \int_{[0,1]^d} \sqrt{\det H}.
$$
The next proposition shows that the mass of a metric $H$ and the cardinality of a mesh $\cT$ are close when
$H$ and $\cT$ are equivalent in the sense studied in the previous chapter.
We denote by $\TEq$ a fixed equilateral simplex centered at the origin of $\R^d$ and having its vertices on the unit sphere.
\begin{prop}
\label{propCardMass}
The cardinality of a mesh $\cT\in \bT_\per$ and the mass of a metric $H\in \bH_\per$ satisfy the following properties:
\begin{itemize}
\item If $\cH_T \leq C^2 H(z)$ for all $T\in \cT$ and all $z\in T$, then 
$
|\TEq|\#(\cT) \leq C^d  m(H).
$
\item 
If $H(z) \leq C^2 \cH_T$ for all $T\in \cT$ and all $z\in T$, then 
$
 m(H) \leq C^d |\TEq|\#(\cT).
$
\end{itemize}
In particular, if $\cT$ is $C$-adapted to $H$ in the sense of Definition \ref{defEqHT}, then 
\be
\label{eqCardMass}
C^{-d} m(H) \leq |\TEq| \,\#(\cT) \leq C^d m(H) .
\ee
\end{prop}

\proof
For any $z\in \R^d$ we denote by $T(z)$ a simplex $T\in \cT$ containing $z$.
Recalling (see Proposition \ref{propHT} in the previous chapter) that for any simplex $T$
$$
|T| \sqrt{\det \cH_T} = |\TEq|,
$$
we obtain
\be
\label{eqCardT}
\#(\cT) = \int_{[0,1]^d} \frac {dz}{|T(z)|} = \frac 1 {|\TEq|} \int_{[0,1]^d} \sqrt{\det \cH_{T(z)}} \, dz.
\ee
Furthermore we have $\det S \leq \det S'$ for any $S,S'\in S_d^+$ satisfying $S \leq S'$.
Applying this remark to the matrices $\cH_{T(z)}$ and $C^2 H(z)$ or $C^{-2} H(z)$ we obtain the announced inequalities, which concludes the proof of the proposition.
\sq

The equivalence between classes of meshes and metrics of $\R^d$ established in
Theorem \ref{thEquiv} of the previous chapter has an immediate generalization to periodic meshes and metrics.

\begin{theorem}
\label{thPer}
There exists a constant $C_0=C_0(d)$ such that for all $C \geq C_0$, 
\begin{enumerate}
\item The collection of meshes $\bT_{i,C} \cap \bT_\per$ is equivalent to the collection of metrics $\bH_i \cap \bH_\per$.
\item If $d=2$, then the collection of meshes $\bT_{a,C} \cap \bT_\per$ is equivalent to the collection of metrics $\bH_a \cap \bH_\per$.
\item If $d=2$, then the collection of meshes $\bT_{g,C} \cap \bT_\per$ is equivalent to the collection of metrics $\bH_g \cap \bH_\per$.
\end{enumerate}
\end{theorem}

\proof
The proof of Theorem \ref{thEquiv} presented in the previous chapter can be adapted to periodic metrics and triangulations with only slight changes, which are left to the reader.
\sq

\subsection{Compactness of metrics}
\label{subsecIntroCompact}

In order to ensure the existence of an ``optimal metric'' for a number of optimization problems 
posed on periodic and isotropic, quasi-acute or graded metrics, we need a property of 
compactness for such metrics.

We equip the set $\bH_\per$ of continuous and $\Z^d$ periodic metrics with a distance $d_\per$ which is defined as follows: for all $H,H'\in \bH_\per$
$$
d_\per(H,H') := \sup_{z\in [0,1]^d} d_\times (H(z), H(z)) = \sup_{z\in \R^d} d_\times (H(z) , H'(z)).
$$
Consider $H, H'\in \bH_\per$ and define $\delta = d_\per(H,H')$. We thus have for all $z\in \R^d$,
\be
\label{pointwiseCont}
e^{-2 \delta} H(z)\leq H'(z) \leq e^{2 \delta} H(z). 
\ee
It immediately follows that 
\be
\label{massCont}
e^{-d\delta} m(H) \leq m(H') \leq e^{d\delta} m(H),
\ee
and 
\be
\label{distCont}
e^{-\delta} d_{H}(x,y) \leq d_{H'}(x,y) \leq e^\delta d_{H}(x,y).
\ee
which shows that the mass $m \mapsto m(H)$ and the distance between to fixed points $H \mapsto d_H(x,y)$ depend continuously on $H \in \bH_\per$.

We establish in \S \ref{secCompact} the following compactness property : for any $M\geq 1$ and any symbol $\star \in \{i,a,g\}$ the collection of metrics 
$$
\{H \in \bH_\per \cap \bH_\star \sep m(H) \leq M\}
$$
is a compact subset of $\bH_\per$.

\subsection{Separation of two domains}

Isotropic and anisotropic mesh adaptation is often used to separate some geometric sets, in the sense given by the following definition. For any closed subset $E \subset \R^d$ and any mesh $\cT\in \bT$ we define the neighborhood $V_\cT(E)$ of $E$ in $\cT$ as follows :
$$
V_\cT(E) := \bigcup_{\substack{T \in \cT,\\ T \cap E\neq\emptyset}} T.
$$
We say that a mesh $\cT\in \bT$ separates two closed sets $X,Y\subset \R^d$ if and only if 
$$
V_\cT(X) \cap V_\cT(Y) = \emptyset,
$$
in other words if $X$ and $Y$ have disjoint neighborhoods in the mesh $\cT$. 

A natural objective is to build a mesh of minimal cardinality which separates two given regions. 
Consider two closed, disjoint and $\Z^d$-periodic sets $X,Y\subset \R^d$. We consider in \S \ref{secOptGeom} the optimization problem posed on meshes 
$$
m_{\star,C_0}(X,Y) := \min\{ \#(\cT) \sep \cT\in \bT_\per \cap \bT_{\star,C_0} \text{ and }  \cT \text{ separates } X \text{ and } Y\}
$$
where $\star \in \{i,a,g\}$ is a symbol and $C_0\geq 1$ is a constant.
We establish that this problem is equivalent to the following
optimization problem posed on metrics, if the dimension is $d=2$:
$$
m_\star(X,Y) := \min\{ m(H) \sep H \in \bH_\per\cap \bH_\star \text{ and } H \text{ separates } X \text{ and } Y\}, 
$$
where we say that a metric $H$ separates $X$ and $Y$ if 
$$
d_H(x,y) \geq 1 \text{ for all } (x,y) \in X\times Y.
$$
Precisely, we show that there exists a constant $C = C(C_0)\geq 1$, independent of the sets $X,Y$,
and such that 
$$
C^{-1} m_\star(X,Y)\leq m_{\star,C_0}(X,Y) \leq C m_\star(X,Y).
$$

Imposing more constraints on the triangulation typically raises the number of triangles required to achieve a given task, and likewise assuming more regularity of the metric raises the mass required for a given task.
We illustrate this remark in the context of the separation of domains as follows.
We fix the parameter $r_0 := 1/4$ and we define the periodic set 
$$
X := \{x\in \R^d \sep d(x,\Z^d) \leq r_0\},
$$
where $d(x,E) := \min \{|x-e|\sep e\in E\}$,
and for any $0<\delta\leq r_0$ 
$$
Y_\delta := \{y\in \R^d \sep d(y,\Z^d)\geq r_0 + \delta\}.
$$
The sets $X$ and $Y_\delta$ are illustrated on Figure \ref{figXYDelta}.
The following equivalences are established as $\delta\to 0$
\begin{eqnarray*}
m_i(X,Y_\delta) &\simeq& \delta^{-(d-1)},\\
m_a(X,Y_\delta) &\simeq& \delta^{-\frac{d-1} 2 } |\ln \delta|,\\
m_g(X,Y_\delta) &\simeq& \delta^{-\frac{d-1} 2}.
\end{eqnarray*}

A construction closely related to the separation of the sets $X$ and $Y_\delta$ is presented on Figure \ref{figDisk} in the main introduction of this thesis.
\begin{figure}
\centering
\includegraphics[width=4cm,height=4cm]{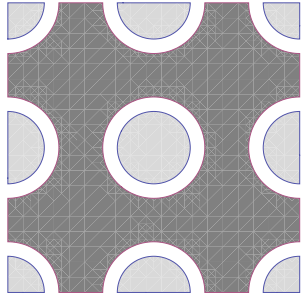}
\caption{The sets $X$ (light gray) and $Y_\delta$ (dark gray), with $d=2$ and $\delta = 0.1$. \label{figXYDelta}}
\end{figure}

\subsection{Approximation in a given norm}

One of the key purposes of anisotropic mesh generation is the adaptive approximation of functions,
as discussed extensively in Chapters 2 and 3. In order to study this problem from the point of view of metrics, we introduce an error $e_H(f)_p$ associated to a function and a metric, and we compare this quantity with the approximation error of $f$ on a mesh.

For that purpose we denote by $B_H(z)$ the open ball of radius $1$ around $z$ for the norm $\|\cdot \|_{H(z)}$
(therefore an ellipse):
$$
B_H(z) := \{z+u \sep \|u\|_{H(z)} < 1\}.
$$
For any exponent $p$, $1\leq p \leq \infty$, and any $f\in L^p_\loc(\R^d)$ we denote by $e_H(f\ssep z)_p$ the error of best approximation of $f$ on $B_H(z)$:
\be
\label{defEH}
e_H(f\ssep z)_p := \inf_{\mu\in \sP_{m-1} } \|f-\mu\|_{L^p(B_H(z))}.
\ee
We denote by $L^p_\per$ the collection of functions $f\in L^p_\loc(\R^d)$ which are $\Z^d$-periodic.
For all $f\in L^p_\per$ and all $H\in \bH_\per$ we define the approximation error $e_H(f)_p$ of $f$ associated to the metric $H$ as follows
\be
\label{defEHGlobal}
e_H(f)_p^p := \int_{[0,1]^d} \sqrt {\det H(z)} \, e_H(f\ssep z)_p^p \, dz.
\ee

Consider a metric $H\in \bH_\per \cap \bH_g$ and a mesh $\cT\in \bT_\per$ such that 
\be
\label{eqTH4}
\cH_T \geq 4^2 H(z), \;\;  \,z\in T,\; T\in \cT,
\ee
which heuristically means that that the mesh $\cT$ is ``more refined'' than the metric $H$.
We establish that for any $f\in L^p_\per$  
\be
\label{eqBestLP}
\inf_{g\in V_{m-1}(\cT)} \|f-g\|_p \leq C e_H(f)_p
\ee
where $V_k(\cT)$ stands for the space of finite elements of degree $k$ on $\cT$, and $C$ is an absolute constant. Furthermore we give an explicit expression of a finite element approximation $g\in V_k(\cT)$ which satisfies \iref{eqBestLP}, namely
\be
\label{defGOpt}
g := \interp_\cT^{m-1} f_{H'}, 
\ee
where $H' := 4^2 H$. We denote by $f_H$ the convolution \emph{distorted by the metric $H$} of a function $f$  with a fixed compactly supported mollifier $\vp$ (which satisfies a moments condition) 
$$
f_H(z) := \int_{\R^d} f(z+H(z)^{-\frac 1 2} u) \vp(|u|)du.
$$

The extension of the approximation result \iref{eqBestLP} to the $W^{1,p}$ semi-norm is not trivial due to the fact that interpolation on simplices of large measure of sliverness is not a stable procedure. 

For any vector field $v\in L^p_\loc(\R^d, \R^d)$ and any metric $H\in \bH$ we introduce the two quantities
\be
\label{defEHA}
e_H^a(v\ssep z)_p := \min_{\nu\in \sP_{m-2}^d} \|v-\nu\|_{L^p(B_H(z))}
\ee
and 
\be
\label{defEHG}
e_H^g(v\ssep z)_p := \|H(z)\|^{\frac 1 2} \min_{\nu\in \sP_{m-2}^d} \|H(z)^{-\frac 1 2} (v-\nu)\|_{L^p(B_H(z))}.
\ee
Note that 
\be
\label{ineqEHAG}
e_H^a(v\ssep z)_p \leq e_H^g(v\ssep z)_p
\ee
with equality if $H(z)$ is proportionnal to $\Id$.
If $v$ is $\Z^d$ periodic and if $H\in \bH_\per$, then we define $e_H^a(v)_p$ and $e_H^g(v)_p$ similarly to \iref{defEHGlobal}.

Consider a metric $H\in \bH_\per \cap \bH_a$ and mesh $\cT\in \bT_\per$ which satisfy \iref{eqTH4}. We establish that for any $f\in W^{1,p}_\per$  
\be
\label{eqBestW1PA}
\inf_{g\in V_{m-1}(\cT)} \|\nabla(f-g)\|_p \leq C \,e_H^a(\nabla f)_p \,\max_{T\in \cT} S(T),
\ee
where $C=C(d)$.

This estimate should be compared with the following. Consider a metric $H\in \bH_\per \cap \bH_g$ and a mesh $\cT\in \bT$ which satisfy \iref{eqTH4} and in addition $C_0^2 H(z) \geq \cH_T$ for all $T\in \cT$ and $z\in T$. 
We establish that 
\be
\label{eqBestW1PG}
\inf_{g\in V_{m-1}(\cT)} \|\nabla(f-g)\|_p \leq C \, e_H^g(\nabla f)_p,
\ee
where $C = C(C_0,d)$.

The function $g$ given by \iref{defGOpt} satisfies both estimates \iref{eqBestW1PA} and \iref{eqBestW1PG}.
These results put in light the price to pay to obtain the best error estimate \iref{eqBestW1PA} in the $W^{1,p}$ norm: the measure of sliverness should be uniformly bounded on $\cT$, and the metric $H$ should belong to $\bH_a$. If these conditions are not satisfied then we can only prove the estimate \iref{eqBestW1PG}, which penalizes the anisotropy of the metric $H$ as observed in \iref{ineqEHAG}.
Note that according to Theorem \ref{thPer}, if the dimension is $d=2$ then for any metric $H\in \bH_a\cap \bH_\per$ there exists a mesh $\cT_0\in \bT_{a,C}\cap \bT_\per$ which is $C'$-equivalent to $H$, where $C$ and $C'$ are absolute constants. By definition of quasi-acute meshes there exists a $C$-refinement $\cT$ of $\cT_0$ on which the measure of sliverness is uniformly bounded by $C$, which satisfies \iref{eqTH4}, and such that $\#(\cT) \leq C''m(H)$ where $C''$ is an absolute constant. 

\subsection{Asymptotic approximation and explicit metrics}
\label{subsecIntroAsympt}
We establish some counterparts for metrics of the asymptotic estimates presented in Chapters 2 and 3 for triangulations and more general meshes. 
Our estimates hold without restriction on the degree $m-1$ of interpolation, $m\geq 2$, or the dimension $d\geq 2$. They apply to a sufficiently smooth function and to metrics of asymptotically large mass.

We denote by $\H_m$ the vector space of all homogeneous polynomials of degree $m$ in $d$ variables, equipped with the norm
$$
\|\pi\| := \sup_{|u|\leq 1} |\pi(u)|.
$$
We define the \emph{shape function} $K$ as follows : for all $\pi \in \H_m$
\be
\label{defKIntro}
K(\pi) =\inf_{A\in \SL_d} \| \pi \circ A\|.
\ee
We establish that for any $\Z^d$ periodic and $C^m$ function
\be
\label{asymptLPMet}
\limsup_{M\to \infty} \left(M^{\frac m d} \inf_{\substack{H \in \bH_g\\  m(H)\leq M}} e_H(f)_p\right) \leq C \|K(d^m f)\|_\tau,
\ee
where the exponent $\tau$ is defined by $\frac 1 \tau := \frac m d+ \frac 1 p$.

The proof of this result would be greatly simplified if one could establish the existence of a smooth map $\pi\in \H_m \mapsto A(\pi)\in \SL_d$ such that $K(\pi) = \|\pi \circ A(\pi)\|$. This is unfortunately not the case, and the key ingredient of the proof of \iref{asymptLPMet} is the definition and the analysis of a family $K^{(\alpha)}$ of functions which tend to $K$ as $\alpha\to \infty$ and which are defined by a well posed optimization problem (for which there exists a unique minimizer depending continuously on the parameter $\pi$). We also prove in \S \ref{subsecKCont} that the shape function $K$ is uniformly equivalent to a continuous function on $\H_m$, as well as the shape functions $L_g$ and $L_g$ defined below.\\

In the case of $W^{1,p}$ norms we introduce two distinct \emph{shape functions}
\begin{eqnarray*}
L_a(\pi) &:=& \inf_{A\in \SL_d} \|(\nabla \pi)\circ A\|, \\
L_g(\pi) &:=& \inf_{A\in \SL_d} \|A^{-1}\| \|\nabla (\pi\circ A)\|, 
\end{eqnarray*}
where $\H_{m-1}^d$ is equipped with the norm
$$
\|\nu\| := \sup_{|u|=1} |\nu(u)|.
$$
We establish similar results to \iref{asymptLPMet}, for any $C^m$ and $\Z^d$ periodic function $f$
\be
\label{asymptA}
\limsup_{M\to \infty}\left(M^\frac {m-1} d \inf_{\substack{H \in \bH_a\\  m(H)\leq M}}  e_H^a(\nabla f)_p\right) \leq C \|L_a(d^m f)\|_\tau
\ee
and 
\be
\label{asymptG}
\limsup_{M\to \infty}\left(M^\frac {m-1} d \inf_{\substack{H \in \bH_g\\  m(H)\leq M}}  e_H^g(\nabla f)_p\right) \leq C \|L_g(d^m f)\|_\tau
\ee
where the exponent $\tau$ is defined by $\frac 1 \tau := \frac {m-1} d + \frac 1 p$.\\

For a quadratic or cubic (in two dimensions) homogeneous polynomial $\pi$, we give an explicit expression of near-minimizers of the optimization problems defining $L_a(\pi)$ and $L_g(\pi)$ in terms of the coefficients of $\pi$.
This expression can be used to compute efficiently in numerical applications a metric $H$ adapted to the approximated function $f$ using the available information on the derivatives of $f$.

We also use these minimizers to compute the ratio $L_a(\pi) / L_g(\pi)\leq 1$, and we discuss for which polynomials this ratio is significantly small. In other words for which types of anisotropic features of the approximated function $f$ the use of quasi-acute triangulations leads to a significantly smaller approximation error than the use of graded anisotropic triangulations (which is the case in all numerical software known to the author).

\section{Compactness of metrics of a given mass} 
\label{secCompact}

The purpose of this section is to establish the following compactness result.
We recall that the collection $\bH_\per$ of $\Z^d$ periodic metrics is equipped with the distance 
$$
d_\per(H,H') := \sup_{z\in \sR^d} d_\times (H(z), H'(z)).
$$
\begin{theorem}
\label{thCompact}
For any symbol $\star \in \{i,a,g\}$ and any $m_0\geq 1$ the collection of metrics 
\be
\label{defHStarM}
\{H \in \bH_\per \cap \bH_\star \sep m(H) \leq m_0\}
\ee
is a compact subset of $\bH_\per$ 
\end{theorem}

The rest of this section is devoted to the proof of this theorem. Our first intermediate result establishes that the set \iref{defHStarM} is closed in $\bH_\per$.

\begin{lemma}
\label{lemmaClosed}
The set \iref{defHStarM} is closed in $\bH_\per$, as well as the sets $\bH_i, \bH_a, \bH_g$.
\end{lemma}

\proof
As observed in \S \ref{subsecIntroCompact}, the pointwise evaluation $H \mapsto H(z)$ at a fixed point $z\in \R^d$, the distance $H \mapsto d_H(x,y)$ between two fixed points $x,y\in \R^d$ and the mass $H \mapsto m(H)$ define continuous maps on $\bH_\per$, see \iref{pointwiseCont},  \iref{massCont} and \iref{massCont}.
For any fixed $x,y\in \R^d$, the collection of all metrics $H \in \bH_\per$ which satisfy one of the following conditions
\begin{eqnarray*}
 d_\times (H(x),H(y)) &\leq& d_H(x,y),\\
 d_+ (H(x),H(y)) &\leq& |x-y|,\\
 \|H(x)\| \|H(x)^{-1}\| &=&1,\\
m(H) &\leq& m_0,
\end{eqnarray*}
is therefore a closed subset of $\bH_\per$.

A metric $H$ belongs to $\bH_g$ (resp. $\bH_a$, resp. $\bH_i$)
if and only if it satisfies the first of these inequalities for all $x,y\in \R^d$ (resp. the two first, resp. the three first).
The sets $\bH_i$, $\bH_a$ and $\bH_g$ are thus an intersection of closed subsets of $\bH_\per$, hence are closed.
The set \iref{defHStarM} is obtained by imposing the fourth inequality, and is therefore also closed.
\sq

We recall a classical result which states that the collection of Lipschitz functions between two compact metric spaces is itself a compact metric space. It is an immediate consequence of Ascoli theorem.
\begin{theorem}
\label{thAscoli}
If two metric spaces $(X,d_X)$ and $(Y,d_Y)$ are compact, then the collection $\Lip(X,Y)$  of Lipschitz functions from $X$ to $Y$ is also compact when
equipped with the distance 
$$
d(f,g) := \max_{x\in X} d_Y(f(x),g(x)).
$$
\end{theorem}

We use this theorem to establish that 
closed and bounded subsets of $\bH_\per \cap \bH_g$ are compact.

\begin{lemma}
\label{lemmaCompact}
Any closed and bounded subset $\bH_*$ of $\bH_\per \cap \bH_g$ is compact.
\end{lemma}

\proof
Since the collection $\bH_*$ of metrics is bounded, the following quantity $M$ is finite 
$$
M^2 := \sup_{H \in \bH_*} \sup_{z\in [0,1]^d} \|H(z)\|.
$$
For any $x,y\in \R^d$ and any $H \in \bH_*$ we thus have $d_H(x,y) \leq M|x-y|$.

We consider the compact metric space $X :=[0,1]^d$, equipped with the distance $d_X(x,y) := M|x-y|$, and the compact metric space 
$$
Y := \{M \in S_d^+ \sep \Id \leq M \leq C^2 \Id\}
$$
equipped with the distance $d_\times$.
According to Theorem \ref{thAscoli} the metric space $\Lip(X,Y)$ is thus compact.

Consider a sequence $(H_n)_{n\geq 0}$ of elements of $H_*$. For all $n \geq 0$, the restriction $H_{n|[0,1]^d}$ is an element of $\Lip(X,Y)$. Since this space is compact there exists a sub-sequence $H_{\vp(n)}$ which converges uniformly on $[0,1]^d$. Recalling that $H_n$ is periodic for all $n\geq 0$ we find that $H_{\vp(n)}$ converges uniformly to a periodic metric $H \in \bH_\per$. 

Since $\bH_*$ is closed in $\bH_\per \cap \bH_g$, and since $\bH_\per \cap \bH_g$ is closed in $\bH_\per$, we obtain that $\bH_*$ is closed in $\bH_\per$ and therefore $H \in \bH_*$ which concludes the proof.
\sq

The next lemma establishes that the set \iref{defHStarM} appearing in Theorem \ref{thCompact} is bounded. This set is closed according to Lemma \ref{lemmaClosed}, hence compact according to Lemma \ref{lemmaCompact}, which concludes the proof of Theorem \ref{thCompact}.

\begin{lemma}
\label{propCompact}
Let $m_0 \geq 1$ be a constant.
\begin{enumerate}
\item
There exists $C=C(m_0,d)$ such that for any $H \in \bH_g$ satisfying $H(z) \geq \Id \text{ for all } |z|\leq 1$ and 
$$
\int_{|z|\leq 1} \sqrt{\det H(z)} dz \leq m_0
$$
we have $\|H(0)\|\leq C^2$.
\item
The set 
$
\{H \in \bH_\per \cap \bH_g \sep m(H) \leq m_0\}
$
is bounded. 
\end{enumerate}
\end{lemma}

\proof
We first show how 2. can be obtained assuming 1., and for that purpose we consider $H \in \bH_\per \cap \bH_g$ and $z_* 
$. The ball $\{z\in \R^d \sep |z-z_*|\leq 1\}$ is covered by the $2^d$ cubes 
$z_*+u+[0,1]^d$ where $u \in \{-1,0\}^d$. Consequently 
$$
\int_{|z-z_*|\leq 1} \sqrt{\det H(z)} dz \leq \sum_{u\in \{-1,0\}^d } \int_{z_*+u+[0,1]^d} \sqrt{\det H} = 2^d m_0
$$
and we already know that $H(z) \geq \Id$ for all $z\in \R^d$. It follows that $\|H(z_*)\|\leq C'^2$ where $C'$ is the constant associated to $2^d m_0$ in Point 1., which concludes the proof of Point 2.\\

We now turn to the proof of Point 1. 
At each point $z\in \R^d$ we denote the eigenvalues of $H(z)$ by 
$$
\lambda_1^2(z) \geq \cdots \geq \lambda_d^2(z).
$$
For each $k \in \{1, \cdots ,d\}$ the function $\lambda_k: \R^d \to \R_+^*$ is continuous and $\lambda(z)\geq 1$ for all $|z|\leq 1$.

Our purpose is to obtain an upper bound for $\lambda_1(0) = \sqrt{\|H(0)\|}$ in terms of $m_0$.
We shall proceed as follows : we first give an upper bound for $\lambda_d(0)$ that depends only on $m_0$, and then for each $1\leq k \leq d-1$ we give an upper bound for $\lambda_k(0)$ in terms of $m_0$ and $\lambda_{k+1}(0)$.

\paragraph{The upper bound on $\lambda_d(0)$.}
We define the ellipse 
$$
\cE := \{z\in \R^d \sep \|z\|_{H(0)} \leq \lambda_d(0)\}
$$
and we observe that $|z|\leq 1$ for all $z\in \cE$.
We thus obtain using Corollary \ref{corolVol} that 
$$
c_d \ln \lambda_d(0) \leq \int_\cE \sqrt{\det H(z)} dz \leq  \int_{|z|\leq 1} \sqrt{\det H(z)} dz \leq m_0,
$$
where $c_d>0$ is a constant which depends only on the dimension $d$. This gives as announced an upper bound on $\lambda_d(0)$ in terms of $m_0$.

\paragraph{The upper bound on $\lambda_k(0)$.}
The vector space $\R^d$ is the orthogonal sum 
$$
\R^d = U\oplus V
$$
 where $U$ is the sum of the eigenspaces  associated to the eigenvalues $\lambda_1^2(0), \cdots, \lambda_k^2(0)$ of $H(0)$ and $V$ is associated to the other eigenvalues $\lambda_{k+1}^2(0), \cdots, \lambda_d^2(0)$ of $H(0)$.
We have for all $u\in U$ and all $v\in V$, 
$$
\|u\|_{H(0)} \geq \lambda_k(0) |u| \stext{ and } \|v\|_{H(0)} \leq \lambda_{k+1}(0) |v|.
$$

We consider an isometry $P\in \cO_{d,k}$ (in the sense that $P^\trans P$ is the $k\times k$ identity matrix) satisfying $P(\R^k) = U$. 
For each $v\in V$ we define a metric $H_v$ on $\R^k$ as follows: for all $u \in \R^k$
$$
H_v(u) := P^\trans H_v(v+ Pu) P.
$$
It follows from Proposition \ref{propHG} (in the previous chapter) that $H_v$ belongs to the collection $\bH_g(\R^k)$ of graded metrics on $\R^k$.
We define 
$$
B_U := \{u \in \R^k \sep |u|\leq 1/2\} \stext{ and } B_V := \{v\in V \sep \lambda_{k+1}(0) |v| \leq 1/2\},
$$
as we observe that 
\be
\label{eqFubiniSet}
|P u+v| \leq 1 \text{ for all }  (u,v)\in B_U \times  B_V.
\ee
Recalling that $\lambda_i(z)\geq 1$ for all $1\leq i \leq d$ and $|z| \leq 1$, we obtain for all $(u,v) \in B_U \times B_V$
\be
\label{eqDetHHV}
\begin{array}{rcl}
\det H_v(u) &\leq& \lambda_1^2(Pu +v) \cdots \lambda_k^2(Pu +v)\\
&\leq &  \lambda_1^2(Pu +v) \cdots \lambda_d^2(Pu +v)\\
& =& \det H(Pu+v).
\end{array}
\ee
For any $u\in \R^k$ and any $v\in V$ we have according to \iref{eqLocalNorm}
\be
\label{eqLowerHV}
\|u\|_{H_v(0)} = \|Pu\|_{H(v)} \geq (1-\|v\|_{H(0)}) \|Pu\|_{H(z_*)} \geq \lambda_k(0) |u|/2, 
\ee
and therefore 
$$
\{u\in \R^k \sep \|u\|_{H_v(0)} \leq \lambda_k(0)/4\} \subset B_U.
$$ 
Denoting by $c_k$ the constant from Corollary \ref{corolVol}, applied in dimension $k$,
we therefore obtain for any $v\in B_V$
\be
\label{eqHVBU}
\int_{B_U} \sqrt{\det H_v(u)} du \geq c_k \ln ( \lambda_k(0)/4).
\ee
We thus obtain successively 
\begin{eqnarray*}
\int_{|z|\leq 1} \sqrt{\det H(z)} dz &\geq& \int_{v\in B_V} \int_{u\in B_U} \sqrt{\det H(P u +v)} du dv\\
&\geq & \int_{v\in B_V} \int_{u \in B_U} \sqrt{\det H_v(u)} du dv\\
& \geq & \int_{v\in B_V} c_k \ln(\lambda_k(0)/4) dv\\
& = & |B_V| c_k \ln( \lambda_k(0)/4)\\
& = & \frac{\omega_{d-k}}{(2\lambda_{k+1}(0))^{d-k}} c_k \ln( \lambda_k(0)/4),
\end{eqnarray*}
where we used \iref{eqFubiniSet} and the Fubini integration formula in the first line, \iref{eqDetHHV} in the second line, and \iref{eqHVBU} in the third line. In the last line we denote by  $\omega_{d-k}$ the volume of the unit euclidean ball in $\R^{d-k}$.

This inequality yields an upper bound on $\lambda_k(0)$ in terms of $M$ and $\lambda_{k+1}(0)$
$$
\lambda_k(0) \leq 4 \exp\left(\frac {m_0 (2\lambda_{k+1}(0))^{d-k}}{\omega_{d-k}}\right),
$$
which concludes the proof of this theorem.
\sq

\section{Separation of two domains}
\label{secOptGeom}

We study in this section a geometrical problem: the separation of two sets using a mesh of minimal cardinality. Our first result, Proposition \ref{propEquivSep2}, establishes that this problem has an equivalent formulation in terms of metrics. 

\begin{proposition}
\label{propEquivSep2}
We assume that the collection of meshes $\bT_\per \cap \bT_{\star,C_0}$ and the collection of metrics $\bH_\per \cap \bH_\star$ are equivalent in dimension $d$, where $\star \in \{a,b,c\}$ is a symbol and $C_0\geq 1$ is a constant (According to Theorem \ref{thPer}, this condition holds at least if $d=2$ and $C_0$ is sufficiently large).\\

There exists a constant $C \geq 1$ such that the following holds:\\
Let $X,Y \subset \R^d$ be two closed, disjoint and periodic sets, in the sense that 
$$
x+u \in X \stext{  and }y+u \in Y, \ \text{ for any } x\in X, \, y\in Y \text{ and }u \in \Z^d.
$$
Consider the two optimization problems 
\begin{eqnarray*}
m_{\star,C_0}(X,Y) &:=& \min\{ \#(\cT) \sep \cT\in \bT_\per \cap \bT_{\star,C_0} \text{ and } \cT \text{ separates } X \text{ and } Y\},\\
m_\star(X,Y) &:=& \min\{ m(H) \sep H \in \bH_\per\cap \bH_\star \text{ and } H \text{ separates } Y \text{ and } Y\}.
\end{eqnarray*}
Then 
\be
\label{eqMetricMeshSep}
C^{-1} m_\star(X,Y)\leq m_{\star,C_0}(X,Y) \leq C m_\star(X,Y).
\ee
\end{proposition}

\proof 
We denote by $C_1\geq 1$ a constant such that for any $\cT\in \bT_\per \cap \bT_{\star, C_0}$ there exists $H\in \bH_\per\cap \bH_\star$ which is $C_1$-equivalent to $\cT$, and vice versa.

The proof of this proposition is a simple translation between the vocabulary of meshes and the vocabulary of metrics.
We first establish that $m_\star(X,Y) \leq C m_{\star,C_0}(X,Y)$, and then that $m_{\star,C_0}(X,Y) \leq C m_\star(X,Y)$.\\

Let $\cT\in \bT_\per \cap \bT_{\star,C_0}$ be a mesh which separates the sets $X$ and $Y$, and let $H \in \bH_\per \cap \bH_\star$ be a metric which is $C_1$ equivalent to $\cT$. We thus have for any $T\in \cT$ and any $z\in T$
$$
C_1^{-2}H(z)\leq \cH_T \leq C_1^2 H(z).
$$
Let $x\in X$ and let $y\in Y$. Let $T_0\in \cT$ be the simplex containing $x$ and let $V_\cT(T_0)$ be the neighborhood of $T_0$ in $\cT$. 
Since $x\in T_0$ and $y\notin V_\cT(T_0)$ we have $\|x-y\|_{\cH_{T_0}}\geq c := (C_0\sqrt d)^{-1}$ according to Proposition \ref{propNeighborT} (in the previous chapter). It follows that 
$$
d_H(x,y) \geq \ln (1+\|x-y\|_{H(x)}) \geq \ln(1+ \|x-y\|_{\cH_{T_0}}/C_1) \geq r:=\ln (1+c C_1^{-1}).
$$
 The metric $\ti H := r^{-2} H$ thus belongs to $\bH_\per \cap \bH_\star$ since $0<r\leq 1$, see Remark \ref{remHomog}. Furthermore this metric separates the sets $X$ and $Y$ and satisfies according to Proposition \ref{propCardMass}
$$
m(\ti H)  \leq (C_1/r)^d |\TEq|\#(\cT).
$$ 
Taking the infimum over all meshes $\cT\in \bT_{*,C_0}$ which separate $X$ and $Y$ we obtain the left part of \iref{eqMetricMeshSep} : $m_*(X,Y) \leq |\TEq| (C_1/r)^d \, m_{*,C_0}(X,Y)$.\\

We now turn to the proof of the right part of \iref{eqMetricMeshSep}, and for that purpose we consider a Riemannian metric $H\in \bH_*\cap \bH_\per$ which separates the sets $X$ and $Y$. We consider a parameter $\lambda \geq 1$, which value will be specified later, and we observe that $\ti H := \lambda^2 H \in \bH_*\cap \bH_\per$.  Therefore there exists a mesh $T\in \bT_{\star, C_0}\cap\bT_\per$ such that for all $T\in \cT$ and all $z\in \cT$ one has
$$
C_1^{-2} \ti H(z)\leq \cH_T \leq C_1^2 \ti H(z).
$$

Let us assume 
that $\cT$ does not separate the sets $X$ and $Y$, hence that $V_\cT(X) \cap V_\cT(Y) \neq \emptyset$. In that case there exists two simplices $T,T'\in \cT$ sharing a vertex $v$, and two points $x\in T \cap X$ and $y\in T'\cap Y$.
For any simplex $T$ and any $z,z'\in T$ we have $\|z-z'\|_{\cH_T} \leq 2$, see \iref{eqDiamTHT}.
Therefore 
$$
\lambda \leq d_{\ti H}(x,y) \leq d_{\ti H}(x,v) + d_{\ti H} (v,y) \leq C_1 (\|x-v\|_{\cH_T}+ \|v-y\|_{\cH_{T'}}) \leq 4 C_1.
$$
We now choose the particular the value $\lambda := 4C_1+1$, which contradicts the previous equation and shows that $\cT$ does separate $X$ and $Y$.
Furthermore
$$
|\TEq| \#(\cT) \leq  \, C_1^d \, m(\ti H) = (\lambda C_1)^d m(H).
$$
Taking the infimum among all metrics $H \in \bH_*$ which separate $X$ and $Y$ we obtain $|\TEq| m_{*,C_0}(X,Y) \leq  (\lambda C_1)^d m_*(X,Y)$, which 
concludes the proof of this proposition.
\sq

As a concrete example, we define $r_0 := 1/4$ and for any $0<\delta\leq r_0$ we consider in the following the sets 
\be
\label{defXYDelta}
X := \{x\in \R^d \sep d(x,\Z^d)\leq r_0\} \stext{ and } Y_\delta := \{y\in \R^d \sep d(y, \Z^d) \geq r_0+ \delta\},
\ee
which are illustrated on Figure \ref{figXYDelta}.
The next theorem gives a sharp estimate of the minimal mass of a metric which separates the sets $X$ and $Y_\delta$, as $\delta\to 0$.
The notation $\alpha_\delta \simeq \beta_\delta$ stands for the following : there exists a constant 
 $C\geq 1$ such that for all $0 <\delta \leq r_0$
$$
C^{-1}\alpha_\delta \leq \beta_\delta \leq C \alpha_\delta.
$$

\begin{theorem}
\label{thOptGeom}
We have the equivalences
\begin{eqnarray*}
m_i(X,Y_\delta) &\simeq& \delta^{-(d-1)}\\
m_a(X,Y_\delta) &\simeq& \delta^{-\frac{d-1} 2 } |\ln \delta|\\
m_g(X,Y_\delta) &\simeq& \delta^{-\frac{d-1} 2}.
\end{eqnarray*}
Furthermore for any $\star \in \{i,a,g\}$  and any $0< \delta \leq r_0$ the explicit metric $H_\delta^\star\in \bH_\star \cap \bH_\per$ defined below is a near minimizer of the problem $m_\star(X,Y_\delta)$. In other words $H_\delta^\star$ separates $X$ and $Y_\delta$ and satisfies $m(H_\delta^\star) \simeq m_\star(X,Y_\delta)$.
\end{theorem}

The proof is divided in two parts: we first give the explicit expression of the metrics $H_\delta^\star$ and we use it to obtain an upper bound on $m_\star(X,Y_\delta)$. We then prove a lower bound on $m_\star(X,Y_\delta)$.

\subsection{Explicit metrics and upper bound}

We define in this section the metrics $H_\delta^\star$, where $\star \in \{i,a,g\}$ and $0 < \delta \leq r_0 := 1/4$. We then establish in Lemma \ref{lemmaHInH} that $H_\delta^\star \in \bH_\star$, in Lemma \ref{lemmaHSeparates} that $H_\delta^\star$ separates $X$ and $Y_\delta$, and we prove in Lemma \ref{lemmaMHDelta} an upper estimate on the masses $m(H_\delta^\star)$ of these metrics. 

In order to define $H_\delta^\star$ we need to introduce some notations. For any $z\in \R^d \sm\{0\}$ we define
\be
\label{defsz}
s(z) := \min\{r_0, \left| |z|-r_0\right|\} \stext{ and } \theta(z) := \frac z {|z|}.
\ee
For any $0 < \delta \leq r_0$ and any $z\in \R^d\sm\{0\}$ we define three matrices $S_\delta^\star(z) \in S_d^+$, where $\star\in \{i,a,g\}$, as follows 
\be
\label{defSDeltaStar}
\begin{array}{rcl}
S_\delta^i(z) &:=& \max\{s(z), \delta\} \, \Id\vspace{0.2ex}\\
S_\delta^a(z) &:=& \max\{s(z), \delta\} \, \theta(z) \theta(z)^\trans + \max\left\{s(z), \sqrt \delta\right\} \, (\Id - \theta(z) \theta(z)^\trans)\vspace{0.2ex}\\
S_\delta^g(z) &:=& \max\{s(z), \delta\} \, \theta(z) \theta(z)^\trans + \sqrt{r_0\max\{s(z), \delta\}} \, (\Id - \theta(z) \theta(z)^\trans).
\end{array}
\ee
%
We also define $S_\delta^\star(0) = r_0\Id$ for any $\star\in \{i,a,g\}$. Eventually we define the metrics $H_\delta^*$ by the equality 
\be
\label{eqHS2}
2  H_\delta^\star(z+u)^{-\frac 1 2} := S_\delta^\star(z)
\ee
for all $z\in \left[-\frac 1 2, \frac 1 2\right]^d$ and all $u \in \Z^d$.
Observe that 
\be
\label{eqS2R0}
H_\delta^\star = 8^2 \Id \stext{ for all } z\in [-1/2,1/2]^d \text{ such that } |z| \geq 2  r_0 =  1 / 2,
\ee
since $S_\delta^\star(z) = r_0 \Id$ on this set. Figure \ref{figDisk} in the main introduction of this thesis illustrates some  (finite) triangulations which are respectively equivalent to the metrics $H_\delta^\star$, $\star\in \{i,a,g\}$.

\begin{lemma}
\label{lemmaHInH}
For all $0<\delta \leq  r_0$ and all $\star\in \{i,a,g\}$ one has $H_\delta^\star \in \bH_\per \cap \bH_\star$.
\end{lemma}
\proof
We first remark that $s = r_0$ on the boundary $\partial ([-\frac 1 2, \frac 1 2]^d)$, and therefore $H_\delta^\star = \Id$ on this set for any $\star \in \{i,a,g\}$. Using \iref{eqHS2} we thus obtain that $H_\delta^\star$ is continuous and $\Z^d$ periodic.

We denote  
$$
s_\per(z) := \min\{r_0, |d(z, \Z^d) -r_0|\}
$$
and we observe that $s_\per$ is a Lipschitz function.
Since 
$
H_\delta^i = 4\Id /\max\{s_\per , \delta\}^2,
$
we obtain that $H_\delta^i\in \bH_i$ (because $z\mapsto \min \{s(z), \delta\}/2$ is Lipschitz), which concludes the case $\star = i$.\\

We now consider the case $\star =g$, and for that purpose we use the notations of \S \ref{secEigen} (in the previous chapter). Our first step is to establish that the local dilatation $\dil_z(H)_\times$ defined in \iref{defDilH} is smaller than $1$ for all $z$ in the set 
$$
\Omega := \{z\in \R^d \sep 0 < |z| < r_0 - \delta\}. 
$$
We define the functions
$$
\lambda := s/2 \stext{ and } \mu := \sqrt s /4 
$$
and we observe that 
$
H_\delta^g = \lambda^{-2} \theta \theta^\trans + \mu^{-2} \theta \theta^\trans
$
on $\Omega$.
The functions $\lambda$ and $\mu$ are $C^\infty$ on $\Omega$ and we have
$$
\nabla \lambda = \frac{-\theta} 2 \stext{ and } \nabla \mu = \frac {-\theta}{8 \sqrt s}.
$$
Therefore for any $z\in \Omega$ 
$$
\dil_z(\ln \lambda\ssep d_H) = \frac{\dil_z(\lambda \ssep d_H)}{\lambda(z)} =  \frac{\|\nabla \lambda(z)\|_{H(z)^{-1}}}{ \lambda(z)}
= \frac{\|-\theta(z)/2\|_{H(z)^{-1}}}{ \lambda(z)} = 1/2,
$$
and 
$$
\dil_z(\ln \mu\ssep d_H) = \frac{\dil_z(\mu \ssep d_H)}{\mu(z)}
=  \frac{\|\nabla \mu(z)\|_{H(z)^{-1}} }{\mu(z)} 
= \frac{\lambda(z)}{8 \sqrt{s(z)} \mu(z)} = 1/4.
$$
On the other hand the derivative of $\theta$ at a point $z\in \R^d\sm \{0\}$ in the direction $u\in \R^d$ is the component of $u$ orthogonal to $\theta$ and divided by $|z|$ : 
$$
\frac{\partial \theta}{\partial u}(z) = \frac{u-\<u, \theta(z)\>\theta(z)}{|z|},
$$
therefore 
$$
\dil_z(\theta \ssep d_H) = \frac {\mu(z)} {|z|}.
$$
On the other hand for $z\in \Omega$
$$
\left|\frac{\lambda(z)}{\mu(z)} - \frac {\mu(z)}{\lambda(z)}\right| = \sqrt{\frac {r_0}{s(z)}} -  \sqrt{\frac {s(z)}{r_0}} \leq r_0-s(z) = |z|,
$$
where we used the fact that $\sqrt t- 1/\sqrt t\leq t$ for all $t\geq 1$ (which is easily checked by derivation).
Therefore 
$$
\frac 1 2 \left|\frac{\lambda(z)}{\mu(z)} - \frac {\mu(z)}{\lambda(z)}\right|\dil_z(\theta \ssep d_H) \leq \frac {\mu(z)} 2 \leq 1/8.
$$
It follows that the quantity $D_z(a)_\times$ defined in \iref{eqDZA} is smaller than $1/2$, and therefore that $\dil_z(H_\delta^g)_\times \leq 1$ on $\Omega$ according to Theorem \ref{th2Spaces}.
Proceeding likewise we also obtain that $\dil_z(H_\delta^g)_\times \leq 1$ on the annulus $\{z\in \R^d \sep r_0+\delta<|z|<2r_0\}$. Furthermore $H_\delta^g$ is constant on the annulus $\{z\in \R^d \sep r_0-\delta <|z|<r_0+\delta \}$ and on the domain $\{z\in (-1/2,1/2)^d  \sep |z|> 2r_0\}$. Defining 
$$
\Gamma_0 := \{z\in \R^d \sep |z| \in \{0,r_0-\delta, r_0+\delta, 2r_0\}\} \cup \partial ([-1/2,1/2]^d),
$$
and $\Gamma := \{z+u \sep z\in \Gamma_0, \, u\in \Z^d\}$ we have thus established that $\dil_z(H_\delta^g)_\times \leq 1$ on $\R^d\sm \Gamma$. Corollary \ref{corolLipNoGamma} thus implies that $H_\delta^g\in \bH_g$.

The proof that $H_\delta^a\in \bH_a$ is extremely similar to the proof that $H_\delta^g\in \bH_g$: 
the quantities $D_z(a)_\times$ and $D_z(a)_+$ defined in \iref{eqDZA} and \iref{eqDZAP} can be explicitely computed on $\R^d\sm\Gamma$ and are smaller that $1/2$. Theorem \ref{th2Spaces} thus implies that 
the local dilatations $\dil_z(H_\delta^a)_\times$ and $\dil_z(H_\delta^a)_+$ are bounded by $1$ on $\R^d\sm \Gamma$, and it follows from Corollary \ref{corolLipNoGamma} that $H_\delta^a\in \bH_a$. The details of the proof are left to the reader.
\sq

The next lemma establishes the that the metrics $H_\delta^\star$ do separate the sets of interest.

\begin{lemma}
\label{lemmaHSeparates}
For all $0<\delta\leq r_0$ and all $\star\in \{i,a,g\}$ the metric $H_\delta^\star$ separates the sets $X$ and $Y_\delta$.\end{lemma}

\proof
Let $\gamma\in C^1([0,1],\R^d)$ be such that $\gamma(0)\in X$ and $\gamma(1)\in Y_\delta$. Since the sets $X$, $Y_\delta$ are $\Z^d$-periodic, as well as the metrics $H_\delta^\star$, we may assume that $\gamma(0)\in \left[-\frac 1 2, \frac 1 2\right]^d$.
Recalling the definition \iref{defXYDelta} of $X$ and $Y_\delta$ we obtain
$$
|\gamma(0)|\leq r_0 \stext{ and } |\gamma(1)| \geq r_0+ \delta.
$$
We define 
$$
t_0 := \max \{t\in [0,1] \sep \gamma(t)\in X\} \stext{ and } t_1 := \min \{t\in [t_0, 1] \sep \gamma(t)\in Y_\delta\}.
$$
For any $t\in [t_0, t_1]$ we have 
$$
\gamma(t_0) = r_0  \leq |\gamma(t)| \leq r_0+\delta = \gamma(t_1).
$$
and therefore 
$$
\|\gamma'(t)\|_{H_\delta^\star(\gamma(t))} \geq \frac {2|\<\gamma'(t), \theta(\gamma(t))\>|}\delta = \frac{2 ||\gamma|'(t)|} \delta,
$$
where $|\gamma|$ stands for the function $z\mapsto |\gamma(z)|$.
Hence 
$$
l_{H_\delta^\star}(\gamma) \geq \int_0^1 \|\gamma'(t)\|_{H_\delta^\star(\gamma(t))} dt \geq \int_{t_0}^{t_1} \frac{2|\gamma|'(t)dt} \delta = 2\frac{|\gamma(t_1)|-|\gamma(t_0)|} \delta = 2.
$$
It follows that $d_{H_\delta^\star}(X,Y_\delta) \geq 2$ which concludes the proof.
\sq

The next lemma estimates the mass of the metrics $m(H_\delta^\star)$, which yields an upper estimate on the minimal mass $m_\star(X,Y_\delta)$ required to separate the sets $X$ and $Y_\delta$ using a metric in $\bH_\per \cap \bH_\star$.

\begin{lemma}
\label{lemmaMHDelta}
There exists a constant $C=C(d)>0$ such that for all $0<\delta \leq r_0$
\begin{eqnarray}
\label{eqHDeltaI}
m(H_\delta^i) &\leq& C\delta^{-(d-1)},\\
\label{eqHDeltaA}
m(H_\delta^a) &\leq& C\delta^{-\frac{d-1} 2 } |\ln \delta|,\\
\label{eqHDeltaG}
m(H_\delta^g) &\leq& C\delta^{-\frac{d-1} 2}.
\end{eqnarray}
\end{lemma}

\proof
Recalling the definition \iref{eqHS2} of $H_\delta^\star$ in terms of $S_\delta^\star$ we obtain 
$$
m(H_\delta^\star) = 8^d c_0 + 2^d \int_{|z|\leq 2r_0}  \frac {dz}{\det S_\delta^\star(z)},
$$
where $c_0$ stands for the volume of the following set, on which $H_\delta^\star = 8^2\Id$ according to \iref{eqS2R0},
$$
c_0 := \left|\left\{z\in \left[-\frac 1 2,\frac 1 2\right]^d \sep |z| \geq 2 r_0\right\}\right|.
$$
Furthermore 
$$
\int_{|z|\leq 2r_0} \frac {dz}{\det S_\delta^\star(z)} = c_1 \int_{r=0}^{2r_0} \frac {r^{d-1}dr}{s_\delta^\star(r)},
$$
where $c_1$ stands for the $d-1$-dimensional volume of the sphere $\{z\in \R^d \sep |z|=1\}$, and where
$$
s_\delta^\star(r) := \det S_\delta^\star(z)
$$
for any $z$ such that $|z| = r$.
Therefore there exists a constant $c_2>0$, independent of $\delta$, $0< \delta\leq r_0$, and of $\star\in \{i,a,g\}$, and such that 
$$
m(H_\delta^\star) \leq c_2 \int_{r=0}^{2r_0} \frac {dr}{s_\delta^\star(r)}.
$$
Recalling the definition \iref{defSDeltaStar} of $S_\delta^i(z)$ we obtain $s_\delta^i(r) = \max\{\delta,|r_0-r|\}^d$ for all $0 \leq r\leq 2r_0$, hence 
\begin{eqnarray*}
 \int_{r=0}^{2r_0} \frac {dr}{s_\delta^i(r)} &=& \int_0^{r_0 - \delta} \frac {dr}{(r_0-r)^d} + \int_{r_0 - \delta}^{r_0+\delta} \frac {dr}{\delta^d} + \int_{r_0+\delta}^{2r_0} \frac {dr}{(r-r_0)^d} \\
&=& 2 \delta^{-(d-1)} + 2\left(\frac{\delta^{-(d-1)}}{d-1} -\frac {r_0^{-(d-1)}} {d-1}\right)\\
&\leq & 4 \delta^{-(d-1)}
\end{eqnarray*}
which establishes \iref{eqHDeltaI}.
Similarly, we find 
$s_\delta^g(r) = r_0^{\frac {d-1} 2} \max\{\delta,|r_0-r|\}^{\frac {d+1} 2}$, hence 
\begin{eqnarray*}
r_0^{\frac {d-1} 2} \int_{r=0}^{2r_0} \frac {dr}{s_\delta^g(r)} &=& \int_0^{r_0 - \delta} \frac {dr}{(r_0-r)^{\frac{d+1} 2}} + \int_{r_0 - \delta}^{r_0+\delta} \frac {dr}{\delta^{\frac{d+1} 2}} + \int_{r_0+\delta}^{2r_0} \frac {dr}{(r-r_0)^{\frac{d+1} 2}} \\
&=& 2 \delta^{-\frac{d-1} 2} + 2\left( \frac{2\delta^{-\frac{d-1} 2} }{d-1} -\frac {2r_0^{-\frac{d-1} 2}} {(d-1)}\right)\\
&\leq & 6 \delta^{-\frac{d-1} 2}
\end{eqnarray*}
which establishes \iref{eqHDeltaG}.
Eventually, we have $s_\delta^a(r) = \max\{\delta,|r-r_0|\} \max\left\{\sqrt\delta, |r-r_0|\right\}^{d-1}$, and therefore 
\be
\label{eqThreeInt}
\int_{r=0}^{2r_0} \frac {dr}{s_\delta^a(r)} = \int_{I_\delta} \frac {dr} {|r_0-r|^d} + \int_{J_\delta} \frac {dr}{|r_0-r| \delta ^{\frac {d-1} 2}} + \int_{r_0-\delta}^{r_0+\delta} \frac {dr}{\delta^{\frac{d+1} 2}},
\ee
where 
$$
I_\delta = \left[0,r_0-\sqrt \delta\right] \cup \left[r_0+ \sqrt \delta, 2 r_0\right] \stext{ and }J_\delta = \left[r_0-\sqrt \delta, r_0-\delta\right] \cup \left[r_0+ \delta, r_0+ \sqrt \delta\right].
$$
Therefore 
\begin{eqnarray*}
\int_{r=0}^{2r_0} \frac {dr}{s_\delta^a(r)} &=& \frac 2 {d-1} \left(\left(\sqrt{\delta}\right)^{-(d-1)}-r_0^{-(d-1)} \right)\\
& & +  2 \delta^{-\frac {d-1} 2} \left(\ln \left(\sqrt \delta\right) - \ln \delta\right) \\
& & + 2 \delta^{-\frac {d-1} 2}\\
&\leq & \delta^{-\frac {d-1} 2} |\ln \delta| + 4 \delta^{-\frac {d-1} 2} 
\end{eqnarray*}
which establishes \iref{eqHDeltaA}.
\sq

\subsection{Lower bound}

We prove in this subsection the lower bound announced in Theorem \ref{thOptGeom}. 
\begin{lemma}
\label{lemmaMLower}
There exists a constant $c=c(d)>0$ such that for all $0<\delta\leq r_0$, 
\begin{eqnarray}
\label{eqHILower}
m_i(X,Y_\delta) &\geq& c\delta^{-(d-1)},\\
\label{eqHALower}
m_a(X,Y_\delta) &\geq& c\delta^{-\frac{d-1} 2 } |\ln \delta|,\\
\label{eqHGLower}
m_g(X,Y_\delta) &\geq& c\delta^{-\frac{d-1} 2}.
\end{eqnarray}
\end{lemma}
\proof 
We consider a metric $H\in \bH_\star$, where $\star\in \{i,a,g\}$, which separates the sets $X$ and $Y_\delta$. In particular $H \in \bH_g$.
We denote by $\bbS := \{v\in \R^d\sep |v| = 1\}$ the euclidean unit sphere of $\R^d$, and by $\omega$ the $d-1$-dimensional volume of $\bbS$. We define for all $v\in \bbS$
$$
n(v) := \|v\|_{H(r_0v)} \stext{ and } m(v) := \inf_{|w|=1} \|w\|_{H(r_0v)} =\|H(r_0v)^{-\frac 1 2}\|^{-1} .
$$
Proposition \ref{propEuclidRiemann} (in the previous chapter), states that for any $z,u\in \R^d$ we have 
\be
\label{eqER}
d_H(z,z+u) \leq -\ln(1-\|u\|_{H(z)}),
\ee
where the right hand side equals $\infty$ by convention if $\|u\|_{H(z)} \geq 1$. 


We now consider a point $v\in \bbS$ and we observe that $r_0 v\in X$ and $(r_0+\delta)v\in Y_\delta$. Therefore 
$$
1 \leq d_H(r_0v, \, (r_0+ \delta) v) \leq -\ln\left(1-\delta \|v\|_{H(r_0 v)}\right),
$$
hence $\lambda \delta n(v) \geq 1$
where
$
\lambda := (1-e^{-1})^{-1}.
$
Let $v,w\in \bbS$, and let $\sigma\in \{-1,1\}$ be such that $\sigma \<v,w\>\geq 0$.  Remarking that $2r_0+ \delta \leq 3/4 \leq 1$ we obtain 
$$
\left|r_0 v+\sigma \sqrt {\delta} w\right|^2 = r_0^2 +\delta + 2  r_0 \sqrt \delta \sigma\<v,w\> \geq r_0^2 +\delta \geq r_0^2 + (2r_0 +\delta) \delta = (r_0+\delta)^2.
$$
Hence $r_0 v \in X$ and $r_0 v+\sigma \sqrt {\delta} w\in Y_\delta$. Therefore 
$$
1 \leq d_H(r_0 v, \, r_0 v+ \sigma \sqrt \delta w) \leq -\ln\left(1-\sqrt \delta \|w\|_{H(r_0 v)}\right),
$$
which implies $\lambda \sqrt \delta m(v) \geq 1$.
We we now distinguish between the three cases $\star = i$, $a$ or $g$.
\begin{itemize}
\item If $\star = i$, then $H(z) = \Id/s(z)^2$ where $s$ is a Lipschitz function, and $s(r_0 v) = n(v)^{-1} \leq \lambda\delta$.
It follows that 
\begin{eqnarray*}
m(H) &\geq& \int_\bbS \int_0^{2r_0} \frac {r^{d-1} dr dv}{s(rv)^d} \\
&\geq& \int_\bbS\int_{r_0}^{r_0+\delta} \frac{r_0^{d-1} dr dv}{(n(v)^{-1}+\delta)^d}\\
&\geq& \frac{r_0^{d-1}\delta\omega }{(\lambda\delta+ \delta)^d}\\
&=& \frac{r_0^{d-1} \omega}{(\lambda+1)^d} \delta^{-(d-1)},
\end{eqnarray*}
where $\omega$ denotes the $d-1$ dimensional measure of $\bbS$.
Taking the infimum among all metrics $H \in \bH_i$ which separate the sets $X$ and $Y_\delta$, we obtain as announced
$$
m_i(X,Y_\delta) \geq c \delta^{-(d-1)} \stext{ where } c= \frac{r_0^{d-1} \omega}{(\lambda+1)^d}.
$$
\item If $\star =a$, then the functions $z\mapsto \|H(z)\|^{-\frac 1 2}$ as well as $z\mapsto \|H(z)^{-\frac 1 2}\|$ are Lipschitz. Hence for any $v\in \bbS$ and any $\mu\in \R$ we have 
\begin{eqnarray*}
\|H((r_0+\mu) v)\|^{-\frac 1 2} &\leq n(v)^{-1}+ |\mu| \leq& \lambda \delta + |\mu|,\\
\|H((r_0+\mu) v)^{-\frac 1 2}\| &\leq m(v)^{-1} + |\mu| \leq& \lambda \sqrt \delta + |\mu|.
\end{eqnarray*}
Recalling that $\det M \geq \|M\| \|M^{-1}\|^{-(d-1)}$ for any $M\in S_d^+$ we obtain 
$$
(\lambda\delta+|\mu|)(\lambda\sqrt\delta+|\mu|)^{d-1}\sqrt{\det H((r_0+\mu) v)} \geq 1.
$$
It follows that  
\begin{eqnarray*}
m(H) &\geq& \omega \int_{-r_0}^{r_0} \frac {(r_0+\mu)^{d-1}d\mu}{(\lambda\delta+|\mu|)(\lambda\sqrt\delta+|\mu|)^{d-1}}\\
&\geq& \omega \int_0^{\sqrt \delta} \frac {r_0^{d-1}d\mu}{(\lambda \delta +\mu)(\lambda\sqrt{\delta} + \sqrt \delta)^{d-1}}\\ 
&\geq& \frac{\omega r_0^{d-1} }{((\lambda+1)\sqrt \delta)^{d-1}} (\ln \sqrt \delta -\ln (\lambda\delta)).
\end{eqnarray*}
Hence 
$$
m_a(X,Y_\delta) \geq c\delta^{-\frac{d-1} 2} (-\ln\delta-2\ln \lambda) \stext{ where } c=\frac{\omega}{2(\lambda+1)^{d-1}}
$$
which concludes the proof of \iref{eqHALower}.
\item We now consider the case $\star = g$. 
We have for any $v\in \bbS$
$$
\sqrt{\det H(r_0 v)} \geq \frac 1 {n(v)m(v)^{d-1}} \geq \frac  1 {n(v) (\lambda \sqrt \delta)^{d-1}}.
$$
Using \iref{eqLocalDet} (in the previous chapter) we therefore obtain for any $\mu\in \R$
$$
\sqrt{\det H((r_0+\mu) v)} \geq (1-|\mu| n(v))^d\sqrt{\det H(r_0v)} \geq \frac{(1-|\mu| n(v))^d } {n(v)(\lambda\sqrt \delta)^{d-1} }
$$
hence 
\begin{eqnarray*}
m(H) &\geq &\int_\bbS \int_{-r_0}^{r_0} \frac{(r_0+\mu)^{d-1} (1-|\mu| n(v))^d d\mu dv}{n(v)(\lambda\sqrt \delta)^{d-1}}\\
&\geq &\int_\bbS \int_0^{\frac 1{2n(v)}} \frac{r_0^{d-1} (1/2)^d d\mu dv}{n(v)(\lambda\sqrt \delta)^{d-1}}\\
&=& \frac{\omega r_0^{d-1}}{2^d (\lambda\sqrt \delta)^{d-1}}.
\end{eqnarray*}
Therefore 
$$
m_g(X,Y_\delta) \geq c \delta^{-\frac {d-1} 2} \stext{ where } c=\frac{\omega r_0^{d-1}}{2^d \lambda^{d-1}},
$$
which concludes the proof of \iref{eqHGLower} and of this lemma. \sq
\end{itemize}

\begin{figure}
\centering
\includegraphics[width = 4cm, height=4cm]{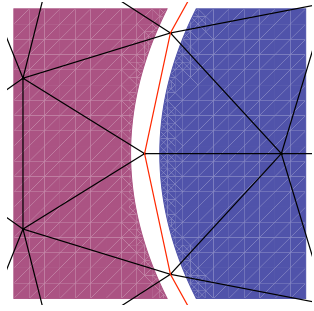}
\caption{Weak separation of two sets by a triangulation. \label{figWeakSeparation}}
\end{figure}

\begin{remark}[Constrained triangulations] 
A variation of the problem considered in this section is often considered in the litterature. 
We say that a mesh $\cT\in \bT$ weakly separates two regions $X, Y \subset \R^d$ if 
$$
\interior(V_\cT(X)) \cap \interior(V_\cT(Y)) = \emptyset,
$$
where $\interior(E)$ refers to the collection of interior points of $E$. 

A common approach for weakly separating two closed and disjoint sets $X,Y\subset \R^d$, is based on \emph{constrained mesh} generation.
The first step is to build a polygonal hypersurface $\Gamma$ such that for any path $\gamma \in C^0([0,1], \R^d)$ satisfying $\gamma(0) \in X$ and $\gamma(1)\in Y$, there exists $s\in [0,1]$ such that $\gamma(s) \in \Gamma$.
The second step is to build a mesh $\cT$ of the domain $\Omega$ which is \emph{constrained} to contain the polygonal surface $\Gamma$ in its skeleton $\cup_{T\in\cT} \partial T$.

The main advantage of weak separation is that it can be achieved efficiently without requiring anisotropic simplices, see Figure \ref{figWeakSeparation}. Heuristically the same number of simplices is required to separate the sets $X$ and $Y_\delta$ ``weakly'' with an isotropic triangulation, or ``strongly'' with a graded triangulation.
We refer the interested reader to \cite{Shew3} for a survey on constrained mesh generation.
\end{remark}

\section{Approximation in a given norm}

We compare in this section the approximation error, measured in the $L^p$ norm or the $W^{1,p}$ semi norm, of a function $f$ by finite elements on a mesh $\cT$ with the error associated to $f$ and a metric $H$. 
In this section, as well as in the next one, we fix an integer $m\geq 2$ 
and we consider finite elements of degree $m-1$.

The greek letter $\alpha$ always refers 
to a $d$-plet of non-negative integers
$$
\alpha = (\alpha_1, \cdots, \alpha_d) \in \Z_+^d,
$$
and we denote $|\alpha| = \alpha_1+\cdots+\alpha_d$.
We define the monomial
$$
Z^\alpha := \prod_{1\leq i\leq d} Z_i^{\alpha_i},
$$
where the variable is $Z=(Z_1,...,Z_d)$, and 
we denote by $\P_k$ the collection of all polynomials of degree at most $k$.
$$
\P_k := \Vect\{Z^\alpha \sep |\alpha|\leq k\}.
$$

We introduce a fixed function $\vp\in L^\infty(\R_+, \R)$ supported in $[0,1]$ and satisfying the following moments property: for all $\mu\in \P_{m-1}$ one has 
\be
\label{eqPhiMoments}
\int_{\sR^d} \mu(z) \vp(|z|) dz = \mu(0).
\ee
For any $f\in L^1_\loc(\R^d)$ and any $H\in \bH$ we denote by $f_H$ the convolution of $f$ with $\vp(|\cdot|)$ distorted by the metric $H$, also referred to as the distorted mollification of $f$ by $H$, and which is defined as follows : for all $z\in \R^d$
\begin{eqnarray}
\label{defConvH}
f_H(z) &:=& \int_{\sR^d} f(z+H(z)^{-\frac 1 2} u) \vp(|u|) du\\
\label{convHDet}
&=& \sqrt{\det H(z)}\int_{\sR^d} f(z+u) \vp(\|u\|_{H(z)}) du.
\end{eqnarray}
Note that $\mu_H = \mu$ for any polynomial $\mu\in \P_{m-1}$. Indeed we obtain using the moments property \iref{eqPhiMoments} that for all $z\in \R^d$ 
$$
\mu_H(z) = \mu(z+H(z)^{-\frac 1 2} 0) = \mu(z).
$$

The next proposition establishes that the distorted mollification $f\mapsto f_H$ is bounded in some functional spaces. 
For any metric $H\in \bH$, any point $z\in \R^d$ and any radius $r>0$ we define
$$
B_H(z,r) := \{z+u \sep \|u\|_{H(z)} < r\},
$$
as well as 
$$
B_H(z) := B_H(z,1) \stext{ and } B(z,r) := B_{\Id}(z,r)
$$
(the latter being the standard open euclidean ball).
Note that 
\be
\label{eqVolBall}
|B_H(z,r)|\sqrt {\det H(z)} = \omega \, r^d \stext{ where } \omega := |B(0,1)|.
\ee

\begin{lemma}
\begin{enumerate}
\label{lemmaBoundedConv}
\item (Boundedness from $L^1_\loc$ to $C^0$) For any $H\in \bH$ and any $f\in L^1_\loc$ one has $f_H\in C^0(\R^d)$. Furthermore for any $z\in \R^d$ 
\be
\label{eqFHBounded}
|f_H(z)|\leq \sqrt{\det H(z)} \, \|\vp\|_{L^\infty} \|f\|_{L^1(B_H(z))},
\ee
where $B_H(z) := \{u\in \R^d \sep \|u\|_{H(z)} <1\}$.
\item (Boundedness from $W^{1,1}_\loc$ to $W^{1,\infty}_\loc$) For any $H\in \bH_a$ and $f\in W^{1,1}_\loc$ one has $f_H\in W^{1,\infty}_\loc$. Furthermore for any $z\in \R^d$ 
$$
\dil_z (f_H)  \leq 2 \sqrt{\det H(z)} \|\vp\|_{L^\infty} \|\nabla f\|_{L^1(B_H(z))}
$$
where $\dil_z(f_H) := \lim_{\ve \to 0} \|\nabla f_H\|_{L^\infty(B(z,\ve))}$ is the local dilatation of $f$ at $z$.
\end{enumerate}
\end{lemma}

\proof
We first establish 1. The continuity of $f_H$ follows from the expression \iref{defConvH} of the distorted mollification, and from general results on parameter dependent integrals. Inequality \iref{eqFHBounded} follows from the fact that $\vp$ is supported on $[0,1]$ and from \iref{convHDet}.

We now turn to the proof of 2. Metrics $H\in \bH_a$ satisfy $\|H(z)^{-\frac 1 2}-H(z')^{-\frac 1 2}\|\leq |z-z'|$ for any $z,z'\in \R^d$, see Definition \iref{defHa} (in the previous chapter). Hence for any fixed $u\in B(0,1)$ the map 
$$
z\mapsto z+H(z)^{-\frac 1 2} u
$$
is $2$-Lipschitz.
Derivating \iref{defConvH} under the integral sign we therefore obtain 
\begin{eqnarray*}
\dil_z (f_H) &\leq& \int_{\sR^d} 2|\nabla f(z+H(z)^{-\frac 1 2} u)| \, |\vp(|u|)| du\\
&\leq& 2 \|\vp\|_{L^\infty} \int_{B_H(z)} |\nabla f(z+H(z)^{-\frac 1 2} u)|  du\\
&=& 2 \sqrt{\det H(z)} \|\vp\|_{L^\infty} \|\nabla f\|_{L^1(B_H(z))}
\end{eqnarray*}
which concludes the proof of this proposition.
\sq

\subsection{The $L^p$ error}
\label{subsecErrMeshMetLP}

This subsection is devoted to the proof of the following proposition, 
which states that the finite element approximation error measured in the $L^p$ norm of a function $f$ on a mesh $\cT$, is controlled by the error $e_H(f)_p$ associated to a metric $H$ when the mesh satisfies the size condition \iref{eqHT4}. 

\begin{prop}
\label{propErrMeshMet}
There exists a constant $C=C(m,d, \vp)$ such that the following holds. Let $H\in \bH_\per \cap \bH_g$ and  
$\cT\in \bT_\per$ be such that 
\be
\label{eqHT4}
4^2H(z) \leq \cH_T,\;\; T\in\cT, \; z\in T.
\ee 
Then for all $1\leq p \leq \infty$ and all $f\in L^p_\per$ we have 
$$
\|f-\interp_\cT^{m-1} f_{H'}\|_{L^p([0,1]^d)} \leq C e_H(f)_p,
$$
where we denoted $H' := 4^2 H$.
\end{prop}

Our first intermediate lemma studies the variations of the balls $B_H(z,r)$ as $z\in \R^d$ and $r>0$ vary and when $H\in \bH_g$ is a graded metric. 
We recall, see \iref{eqLocalNorm} and \iref{eqLocalDet} (in the previous chapter), that for any $H \in \bH_g$, and for any $z,u,v\in \R^d$ satisfying $\|u\|_{H(z)}<1$ one has
\be
\label{eqLocalNormDet}
\begin{array}{ccccc}
 (1-\|u\|_{H(z)} ) \|v\|_{H(z)} &\leq& \|v\|_{H(z+u)} &\leq&(1-\|u\|_{H(z)})^{-1}  \|v\|_{H(z)},\vspace{1mm}\\
(1-\|u\|_{H(z)} )^d \sqrt{\det H(z)}  &\leq& \sqrt{\det H(z+u)} &\leq& (1-\|u\|_{H(z)} )^{-d} \sqrt{\det H(z)} .
\end{array}
\ee

\begin{lemma}
\label{lemmaEllIncl}
Let $H\in \bH_g$ and let $z\in \R^d$. 
\begin{enumerate}
\item  For any $z'\in B_H(z,1/2)$ one has $B_H(z',1/4)\subset B_H(z)$. 
\item For any $z'\in B_H(z, 1/6)$ one has $B_H(z',1/2)\supset B_H(z,1/4)$. 
\end{enumerate}
\end{lemma}

\proof
We first prove Point 1., and for that purpose we observe that for any $q\in \R^d$
\begin{eqnarray*}
\|q-z'\|_{H(z')}
&\geq & (1-\|z'-z\|_{H(z)}) \|q-z'\|_{H(z)}\\
& \geq & (1-\|z'-z\|_{H(z)}) (\|q-z\|_{H(z)}-\|z'-z\|_{H(z)})\\
& \geq & (1-1/2) (\|q-z\|_{H(z)}-1/2)
\end{eqnarray*}
If $\|q-z\|_{H(z)}\geq 1$ then the previous equation implies $\|q-z'\|_{H(z')}\geq  1/4$, which establishes the first inclusion.

We now turn to the proof of Point 2.
We obtain for any $q\in \R^d$ 
\begin{eqnarray*}
\|q-z'\|_{H(z')} &\leq & (1-\|z'-z\|_{H(z)})^{-1} \|q-z'\|_{H(z)}\\
& \geq & (1-\|z-z'\|_{H(z)}) (\|q-z\|_{H(z)}+\|z'-z\|_{H(z)})\\
& \geq & (1-1/6)^{-1} (\|q-z\|_{H(z)}+1/6). 
\end{eqnarray*}
If $\|q-z\|_{H(z)} < 1/4$ then $\|q-z'\|_{H(z')} < \frac 6 5 (\frac 1 4+ \frac 1 6) = 1/2$. Hence $B_H(z',1/2)\supset B_H(z, 1/4)$ as announced which concludes the proof of this lemma.\sq

The next lemma estimates the finite element approximation error \emph{locally}. 

\begin{lemma}
There exists $C=C(m,d, \vp)$ such that the following holds. Let $H\in \bH_g$, $\cT\in \bT$, $z\in \R^d$ and let 
\be
\label{defTz}
\cT_z := \{T\in \cT \sep T \subset B_{H}(z,1/2)\}.
\ee
Then 
\be
\label{eqLocalMeshMet}
\|f-\interp_\cT^{m-1} f_{H'}\|_{L^p(\Omega_z)} \leq C \, e_H(f\ssep z)_p,
\ee
where $H':=4^2 H$ and $\Omega_z$ denotes the union of all triangles in $\cT_z$.
\end{lemma}

\proof
We denote by $C_{\interp}$ the norm of the interpolation operator 
$$
\interp_\TEq^{m-1} : C^0(\TEq) \to C^0(\TEq),
$$
where $\TEq$ denotes the reference equilateral simplex.
We obtain using a change of variables that for any triangle $T$ and any $g\in C^0(T)$ one has 
$$
\|\interp_T^{m-1} g\|_{L^\infty(T)}\leq C_{\interp} \|g\|_{L^\infty(T)}.
$$

Let $\mu\in \P_{m-1}$ be a polynomial satisfying $\|f-\mu\|_{L^p(B_H(z))} = e_H(f\ssep z)_p$.
Let $g := f-\mu$ and let $T\in \cT_z$. 
We have
\begin{eqnarray*}
\|\interp_T^{m-1} g_{H'}\|_{L^p(T)} &\leq & |T|^\frac 1 p \|\interp_T^{m-1} g_{H'}\|_{L^\infty(T)}\\
&\leq & C_{\interp} |T|^\frac 1 p \|g_{H'}\|_{L^\infty(T)}\\
&\leq & C_{\interp} |T|^\frac 1 p \|g_{H'}\|_{L^\infty(B_H(z,1/2))}.
\end{eqnarray*}
It follows from \iref{eqFHBounded} that 
\begin{eqnarray*}
\|g_{H'}\|_{L^\infty(B_H(z,1/2))} &\leq & \|\vp\|_{L^\infty} \max \{ \|g\|_{L^1(B_{H'}(z'))} \sqrt{\det H'(z')}\sep z'\in B_H(z,1/2)\}\\
&\leq & \|\vp\|_{L^\infty}   \|g\|_{L^1(B_{H}(z))} 4^d (1-1/2)^{-d} \sqrt {\det H(z)},
\end{eqnarray*}
where we used the inclusion stated in Point 1. of Lemma \ref{lemmaEllIncl} and \iref{eqLocalNormDet} in the second line.
Defining $C_* := 8^d C_{\interp}  \|\vp\|_{L^\infty}$ we thus obtain
\begin{eqnarray*}
\|\interp_\cT^{m-1} g_{H'}\|_{L^p(\Omega_z)} &=& \left(\sum_{T\in \cT_z} \|\interp_T^{m-1} g_{H'}\|_{L^p(T)}^p \right)^{\frac 1 p} \\
&\leq & C_* \|g\|_{L^1(B_{H}(z))} \sqrt{\det H(z)} \left(\sum_{T\in \cT_z} |T|\right)^\frac 1 p\\
& \leq & C_* \|g\|_{L^1(B_{H}(z))}\frac{\omega}{|B_H(z)|} |B_{H}(z)|^\frac 1 p\\
& \leq & C_*\omega \|g\|_{L^p(B_{H}(z))},
\end{eqnarray*}
where we used in the third line the expression \iref{eqVolBall} of the volume $|B_H(z)|$, and the inclusion \iref{defTz}.
Since the distorted convolution reproduces polynomials in $\P_{m-1}$ , as well as the interpolation operator $\interp_\cT^{m-1}$, we have
$$
g-\interp_\cT^{m-1} g_H = (f-\mu)-\interp_\cT^{m-1} (f-\mu)_H = f-\interp_\cT^{m-1} f_H.
$$ 
Therefore  
\begin{eqnarray*}
\|f-\interp_\cT^{m-1} f_H\|_{L^p(\Omega_z)} &=& \|g-\interp_\cT^{m-1} g_H\|_{L^p(\Omega_z)} \\
&\leq &\|g\|_{L^p(\Omega_z)}+ \|\interp_\cT^{m-1} g_H\|_{L^p(\Omega_z)}\\
& \leq & (1+C_*\omega) \, \|g\|_{L^p(B_{H}(z))}\\
&=& (1+C_*\omega) \, e_H(f\ssep z)_p
\end{eqnarray*}
which concludes the proof.
\sq

We now prove Proposition \ref{propErrMeshMet}.
For each triangle $T\in \cT$ we define 
$$
\Omega_T := \{z\in \R^d \sep T \subset B_H(z,1/2)\} \stext{ and } n(T) := \int_{\Omega_T} \sqrt{\det H(z)} dz,
$$
and we observe that 
$$
z\in \Omega_T \stext{ if and only if } T \in \cT_z,
$$
where $\cT_z$ is defined in \iref{defTz}.
Note also that $\Omega_T\subset B(z_T, 1/2)$ since $H(z) \geq \Id$ for all $z\in \R^d$, hence the sets $u+ \Omega_T$, $u\in \Z^d$, are pairwise disjoint.
We denote by $\cT_*$ a system of representatives of the set $\cT$ for the relation of equivalence  
$$
T \sim T' \stext{ if } T =T'+u \text{ for some } u \in \Z^d.
$$
It follows from the local error estimate \iref{eqLocalMeshMet} that 
\begin{eqnarray*}
C^p e_H(f)_p^p &=& C^p\int_{[0,1]^d} \sqrt{\det H(z)} \, e_H(f\ssep z)_p^p \, dz\\
&\geq&  \int_{[0,1]^d} \sqrt{\det H(z)}\, \|f-\interp_\cT^{m-1} f_H\|_{L^p(\Omega_z)}^p \, dz \\
&=&  \int_{[0,1]^d} \sqrt{\det H(z)}\, \sum_{T\in \cT_z} \|f-\interp_\cT^{m-1} f_H\|_{L^p(T)}^p \, dz \\
&=& \sum_{T\in \cT_*} \|f-\interp_\cT^{m-1} f_H\|_{L^p(T)}^p \int_{[0,1]^d} \sqrt{\det H(z)} \sum_{u\in \Z^d} \chi_{\Omega_T} (z-u) dz\\
&=&\sum_{T\in \cT_*} \|f-\interp_\cT^{m-1} f_H\|_{L^p(T)}^p \int_{\R^d} \sqrt{\det H(z)} \chi_{\Omega_T} (z) dz\\
&=& \sum_{T\in \cT_*}  \|f-\interp_\cT^{m-1} f_H\|_{L^p(T)}^p n_T \\
& \geq &  \|f-\interp_\cT^{m-1} f_H\|_{L^p([0,1]^d)}^p \min_{T \in \cT} n_T.
\end{eqnarray*}
We now establish that $n_T$ is uniformly bounded below, which concludes the proof of Proposition \ref{propErrMeshMet}.
It follows from  \iref{eqHT4} that $T\subset B_H(z_T,1/4)$, where $z_T$ denotes the barycenter of $T$.
The inclusion property established in Point 2. of Lemma \ref{lemmaEllIncl} thus implies that $B_H(z_T,1/6) \subset \Omega_T$. We thus obtain the lower bound 
\begin{eqnarray*}
n(T) &\geq& |B_H(z_T,1/6)| \min_{z\in B_H(z_T,1/6)}  \sqrt{\det H(z)}\\
&\geq& \frac{\omega 6^{-d}}{\sqrt{\det H(z_T)}}  \sqrt{\det H(z_T)} (1-1/6)^d\\
&=& \omega (5/36)^d,
\end{eqnarray*}
where we used the explicit expression \iref{eqVolBall} of the volume of $|B_H(z,r)|$ and the estimate \iref{eqLocalNormDet} on the variations of $\sqrt{\det H(z)}$.

\subsection{The $W^{1,p}$ error, when the measure of sliverness is uniformly bounded}
We establish in this subsection the following proposition 
which shows that the finite element approximation error, measured in the $W^{1,p}$ semi norm on a mesh with bounded measure of sliverness, is controlled by $e_H^a(\nabla f)_p$.

\begin{prop}
\label{propErrMeshMetA}
There exists a constant $C=C(m,d, \vp)$ such that the following holds. Let $H\in \bH_\per \cap \bH_a$ and  
$\cT\in \bT_\per$ be such that 
\be
\label{eqHT4A}
4^2H(z) \leq \cH_T,\;\; T\in\cT, \; z\in T.
\ee 
Then for all $1\leq p \leq \infty$ and all $f\in W^{1,p}_\per$ we have 
$$
\|\nabla (f-\interp_\cT^{m-1} f_{H'})\|_{L^p([0,1]^d)} \leq C e_H^a(f)_p \max_{T\in \cT} S(T),
$$
where we denoted $H' := 4^2 H$.
\end{prop}

The main difficulty of this proof is contained in the following lemma, which states that a near-best approximation of the gradient of a function, on an ellipse, has the form of the gradient of a polynomial. In this proposition we denote by $u\cdot v$ the scalar product of two vectors $u,v\in \R^d$.

\begin{lemma}
\label{lemmaProj}
\begin{enumerate}
\item (The orthogonal projection is compatible with gradients)
We equip the space $L^2(\bbB, \R^d)$, where $\bbB := B(0,1)$, with the scalar product 
$$
\<v,w\> := \int_\bbB v(z)\cdot w(z) dz.
$$
For any $f\in H^1(\bbB)$ there exists $\mu\in \P_{m-1}$ such that the orthogonal projection of $\nabla f$ onto $\P_{m-2}^d$ is $\nabla \mu$.
\item (Optimization among gradients)
There exists $C_\proj=C_\proj(m,d)$ such that for any ellipsoid $\cE\subset \R^d$, any $1\leq p \leq \infty$ and any $f\in W^{1,p}(\cE)$ one has 
\be
\label{eqGradApp}
\inf_{\mu \in \sP_{m-1}} \|\nabla f-\nabla \mu\|_{L^p(\cE)}
\leq 
C_\proj \inf_{\nu \in \sP_{m-2}^d} \|\nabla f-\nu\|_{L^p(\cE)}
\ee
\end{enumerate}
\end{lemma}

\proof
We denote $k:=m-1$, and 
we first establish Point 1. We define $\bbS := \partial \bbB$ and we denote by $D_k$ the collection of polynomial vector fields $\nu\in  \P_{k-1}^d$ satisfying 
$$
\left\{
\begin{array}{cl}
\diver \nu (z)=0 & \text{ for all } z\in \bbB\\
\nu(z) \cdot z = 0 &\text{ for all } z\in \bbS
\end{array}
\right.
$$
For any $g\in H^1(\bbB)$ and any $\nu \in D_k$ we thus have   
\be
\label{eqGradDiv}
\int_\bbB \nabla g \cdot \nu = -\int_\bbB g(z) \diver \nu(z)dz + \int_\bbS g(z)  \, \nu(z) \cdot z \, dz =0.
\ee
Our first objective is to establish the orthogonal decomposition
$$
\P_{k-1}^d = D_k \oplus \nabla \P_k.
$$
The spaces $\nabla \P_k$ and $D_k$ are orthogonal according to \iref{eqGradDiv}, and we know that $\dim \nabla \P_k = \dim \P_k -1$. Therefore we only have to show that 
$\dim D_k \geq \dim (\P_k^{d-1})  - \dim \P_k +1$.

The vector space $D_k$ is the kernel of the map $d_k$ defined as follows
$$
\begin{array}{ccc}
d_k :\P_{k-1}^d &\to& \P_{k-2} \oplus \P_k(\bbS)\\
\nu  &\mapsto & (\diver \nu, \ (\nu (z) \cdot z)_{|\bbS}),
\end{array}
$$
where  
$
\P_k(\bbS)
$
denotes the collection of restrictions to $\bbS$ of elements of $\P_k$. 
The kernel of the map
$$
\begin{array}{ccc}
\gamma_k : \P_k &\to & \P_k(\bbS)\\
\mu &\mapsto &\mu_{|\bbS}
\end{array}
$$
is 
$
\Ker(\gamma_k) = (1-|z|^2) \P_{k-2}.
$
Hence 
$
\dim(\P_k(\bbS)) = \dim \P_k - \dim \P_{k-2}.
$

The map $d_k$ is not surjective, since we have the linear relation 
$$
\int_{\bbB} \diver \nu(z)\, dz = \int_\bbS \nu(z)\cdot z \, dz.
$$
Therefore 
\begin{eqnarray*}
\dim D_k &=& \dim \Ker d_k \\
&=& \dim (\P_k^{d-1}) - \dim \Image d_k \\
&\geq &  \dim (\P_k^{d-1}) - (\dim \P_{k-2} + \dim \P_k(\bbS) -1) \\
&=& \dim (\P_k^{d-1})  - \dim \P_k +1,
\end{eqnarray*}
which implies the orthogonal decomposition $\P_{k-1}^d = D_k \oplus \nabla \P_k$ as announced.

Finally for any $f\in H^1(\bbB)$ the vector field $\nabla f$ is orthogonal to $D_k$ according to \iref{eqGradDiv}. The orthogonal projection of $\nabla f$ onto $\P_{k-1}^d$ is therefore also orthogonal to $D_k$, hence is an element of $\nabla \P_k$ which concludes the proof Point 1.\\

We now turn to the proof of Point 2.
We denote by $Q : L^2(\bbB,\R^d) \to \P_{k-1}^d$ the $L^2(\bbB)$ orthogonal projection.
We recall that for any $v\in L^2(\bbB,\R^d)$ the projection $Q(v)\in \P_{k-1}^d$ is defined by the collection of linear conditions
$$
\int_{\bbB} (v-Q(v))\cdot \nu  = 0 \text{ for all } \nu \in \P_{k-1}^d.
$$
Hence $Q(v)$ depends linearly on a finite number of moments of $v$, and therefore continuously extends to all $v\in L^1(\bbB, \R^d)$. 
Therefore there exists a constant $C_0$ such that for all $1\leq p \leq \infty$ and all $v\in L^p(\bbB,\R^d)$
$$
\|Q(v)\|_{L^p(\bbB)} \leq C_0 \|v\|_{L^p(\bbB)}.
$$
For any $A\in \GL_d$ and any $v\in L^2(\bbB,\R^d)$ we claim that $Q(A(v)) = A(Q(v))$. Indeed for any $\nu\in \P_{k-1}^d$ one has 
$$
\int_{\bbB} (A(v(z))-A(Q(v(z))))\cdot \nu(z) dz= \int_{\bbB} \left(v(z)-Q(v(z))\right)\cdot A^\trans(\nu(z)) dz = 0.
$$
Consider  $1\leq p \leq \infty$, $f\in W^{1,p}(\bbB)$ and $A\in \GL_d$. Let $\mu\in \P_{m-1}$ be such that $Q(\nabla f) = \nabla \mu$, and let $\nu \in \P_{m-1}^d$, we obtain 
\begin{eqnarray*}
\|A(\nabla f-\nabla \mu)\|_{L^p(\bbB)} &=& \|A(\nabla f)-A(Q(\nabla f))\|_{L^p(\bbB)}\\
&=& \|A(\nabla f)-Q(A(\nabla f))\|_{L^p(\bbB)}\\
& = & \|A(\nabla f-\nu)-Q(A(\nabla f)-\nu)\|_{L^p(\bbB)}\\
&\leq& (1+C_0)  \|A(\nabla f-\nu)\|_{L^p(\bbB)}.
\end{eqnarray*}
Therefore 
$$
\inf_{\mu\in \sP_{m-1}} \|A(\nabla f-\nabla \mu)\|_{L^p(\bbB)} \leq (1+C_0) \inf_{\nu\in \sP_{m-2}^d} \|A(\nabla f-\nu)\|_{L^p(\bbB)}.
$$
Using a linear change of variables $\cE \to \bbB$ we obtain the announced result \iref{eqGradApp}, with the constant $C_\proj := C_0+1$. This concludes the proof of this lemma.
\sq

The rest of the proof of Proposition \ref{propErrMeshMetA} is extremely similar to the proof of Proposition \ref{propErrMeshMet}.
The estimates 
\begin{eqnarray*}
|f_H(z)| &\leq& \sqrt{\det H(z)} \, \|\vp\|_{L^\infty} \|f\|_{L^1(B_H(z))}\\
\| \interp_T^{m-1} g\|_{L^\infty(T)} &\leq& C_I \|\nabla g\|_{L^\infty(T)},
\end{eqnarray*}
used in the proof of Proposition \ref{propErrMeshMet} have the counterparts
\begin{eqnarray*}
\dil_z (f_H)  &\leq& 2 \sqrt{\det H(z)} \|\vp\|_{L^\infty} \|\nabla f\|_{L^1(B_H(z))}\\
\| \nabla \interp_T^{m-1} g\|_{L^\infty(T)} &\leq& C_I' \|\nabla g\|_{L^\infty(T)} S(T),
\end{eqnarray*}
which are proved respectively in Lemma \ref{lemmaBoundedConv} and Lemma \ref{lemmaSD} (in Chapter \ref{chapW1P}).
From this point the adaptation of the proof is straightforward and is therefore left to the reader.

\subsection{The $W^{1,p}$ error, on a general mesh} 
\label{subsubsecErrMeshMetW1PG}

This subsection is devoted to the proof of Proposition \ref{propErrMeshMetG}, which estimates the $W^{1,p}$ approximation error of a function $f$ on a mesh with arbitrary measure of sliverness, in terms of the approximation error $e_H^g(\nabla f)_p$ with respect to a metric. This approach does not allow to recover the optimal anisotropic error estimates for $W^{1,p}$ norms developed in Chapter \ref{chapW1P}. However it is perhaps more taylored to current anisotropic mesh generation software, which generally do not guarantee any bound on the measure of sliverness (or the maximal angle in dimension $2$, since $S(T) = \max\{1, \tan( \theta/2)\}$ for any triangle $T$ of maximal angle $\theta$). 

\begin{prop}
\label{propErrMeshMetG}
For all $C_0 \geq 1$ there exists a constant $C=C(C_0, m,d, \vp)$ such that the following holds. Let $H\in \bH_\per \cap \bH_g$ and  let
$\cT\in \bT_\per$ be such that 
\be
\label{eqHT4G}
4^2H(z) \leq \cH_T\leq (4C_0)^2 H(z),\;\; T\in\cT, \; z\in T.
\ee 
Then for all $1\leq p \leq \infty$ and all $f\in W^{1,p}_\per$ we have 
$$
\|\nabla (f-\interp_\cT^{m-1} f_{H'})\|_{L^p([0,1]^d)} \leq C e_H^g(f)_p,
$$
where we denoted $H' := 4^2 H$.
\end{prop}

Our first step in the proof of this proposition is the next lemma which estimates the \emph{local} approximation error.

\begin{lemma}
For all $C_0\geq 1$ there exists $C=C(C_0,m,d)$ such that the following holds. 
Let $H\in \bH_\per$ and $\cT \in \bT_\per$ be such that \iref{eqHT4G} holds. Let $T\in \cT$ and $z\in T$. Then for all $1\leq p \leq \infty$ and all $f\in W^{1,p}_\per$ we have 
\be
\label{eqLocGW1P}
\|\nabla (f- \interp_\cT^{m-1} f_{H'})\|_{L^p(T)} \leq C e_H^g(\nabla f\ssep z)_p,
\ee
where we denoted $H' := 4^2 H$.
\end{lemma}

\proof 
We assume in a first time that $z=0$ and that $H(z) = \Id$, in such way that $B_H(z) = \bbB := B(0,1)$ is the euclidean unit ball.
We consider a polynomial $\mu\in \P_{m-1}$ such that 
$
\|\nabla (f-\mu)\|_{L^p(B)}
$
is minimal, and such that $\int_B (f-\mu) = 0$. The latter point can be ensured by adding an appropriate constant to $\mu$. Defining $g:=f-\mu$, we obtain using Lemma \ref{lemmaProj}
$$
\|\nabla g\|_{L^p(\bbB)} \leq C_\proj \,e_H^a(f\ssep z)_p.
$$
Since  $g\in W^{1,1}(B)$ and $\int_B g=0$ there exists according to Sobolev's injection theorem a constant $C_\sob=C_\sob(m,d)$ such that 
$$
\|g\|_{L^1(\bbB)} \leq C_\sob \|\nabla g\|_{L^1(\bbB)}.
$$
The function $g_{H'}$ is continuous on $T$ and satisfies according to Lemma \ref{lemmaBoundedConv}
\begin{eqnarray*}
\|g_{H'}\|_{L^\infty(T)} &\leq& \|\vp\|_{\infty} \max \{ \|g\|_{L^1(B_{H'}(z'))} \sqrt{\det H'(z')} \sep z'\in T\}\\
 &\leq& \|\vp\|_\infty \|g\|_{L^1(\bbB)} 4^d (1-1/2)^d\\ 
 &=& C_1\|g\|_{L^1(\bbB)},
\end{eqnarray*}
where we used \iref{eqLocalNormDet} and the fact that $\|z-z'\|_{H(z)} \leq 1/2$ for all $z'\in T$.

We denote by $C_{\interp}$ the norm of the operator $\nabla \interp_\TEq^{m-1} : C^0(\TEq)\to \P_{m-2}^d$, where $\TEq$ denotes the reference equilateral simplex : for any $h\in C^0(\TEq)$
$$
\|\nabla \interp_\TEq^{m-1} h\|_{L^\infty(\TEq)} \leq C_{\interp} \|h\|_{L^\infty(\TEq)}. 
$$
We recall that $\cH_T^{\frac 1 2}(T)$ is the image $\TEq$ by a translation and a rotation, see Proposition \ref{propHT} in the previous chapter. Choosing the function $h:=g_{H'}\circ \cH_T^{-\frac 1 2}\in C^0(\TEq)$, we obtain after a change of variables
$$
\|\cH_T^{-\frac 1 2}(\nabla \interp_T^{m-1} g_{H'})\|_{L^\infty(T)} \leq C_{\interp} \|g_{H'}\|_{L^\infty(T)},
$$
which implies 
$$
\|\nabla \interp_T^{m-1} g_{H'}\|_{L^\infty(T)} \leq C_{\interp}\|\cH_T^\frac 1 2\|  \|g_{H'}\|_{L^\infty(T)}.
$$
We thus obtain, since $\|\cH_T^\frac 1 2\|\leq C_0 \|H'(z)^\frac 1 2\| = 4 C_0$ and since $|T|\leq |\bbB|$,
\begin{eqnarray*}
\|\nabla \interp_T^{m-1} g_{H'}\|_{L^p(T)} &\leq& |\bbB|^\frac 1 p \|\nabla \interp_T^{m-1} g_{H'}\|_{L^\infty(T)}\\
& \leq & 4C_0 C_{\interp} |\bbB|^\frac 1 p   \|g_{H'}\|_{L^\infty(T)}\\
& \leq & 4C_0 C_1 C_{\interp} |\bbB|^\frac 1 p   \|g\|_{L^1(\bbB)}\\
& \leq & 4C_0 C_1 C_{\interp} C_\sob |\bbB|^\frac 1 p   \|\nabla g\|_{L^1(\bbB)}\\
&\leq &  4C_0 C_1 C_{\interp} C_\sob |\bbB|  \|\nabla g\|_{L^p(\bbB)}.
\end{eqnarray*}
Defining $C_* := 4C_0 C_1 C_{\interp} C_\sob |\bbB|$ we obtain 
\begin{eqnarray*}
\|\nabla (f-\interp_T^{m-1} f_{H'})\|_{L^p(T)} &=& \|\nabla (g-\interp_T^{m-1} g_{H'})\|_{L^p(T)}\\
&\leq& (1+C_*) \|\nabla g\|_{L^p(\bbB)}\\
&\leq& C\, e_H^a(\nabla f\ssep z)_p,
\end{eqnarray*}
where $C=(1+C_*) C_\proj$.

We now turn to the general case, and we do not any more assume that $z=0$ and $H(z) =\Id$.
We define the affine change of coordinates $\Phi(z') := H(z)^{-\frac 1 2} z'+z$, and we apply our previous reasonning to the metric $H_\Phi$ defined by \iref{defHPhi} (in the previous chapter), 
$$
H_\Phi(z') := H(z)^{-\frac 1 2} H(\Phi(z')) H(z)^{-\frac 1 2},
$$
which satisfies $H_\Phi(0) = \Id$ and $H_\Phi\in \bH_g$ according to Proposition \ref{propHG} (in the previous chapter).
We also define the triangle $T_\Phi := \Phi^{-1}(T)$, the function $f^\Phi := f\circ \Phi$ and the metric $H'_\Phi = 4^2 H_\Phi$.
The first part of this proof implies that 
$$
\|\nabla (f^\Phi-\interp_{T_\Phi}^{m-1} f^\Phi_{H_\Phi'})\|_{L^p(T_\Phi)} \leq C e_{H_\Phi}^a(\nabla f^\Phi \ssep 0)_p.
$$ 
We thus obtain
\begin{eqnarray*}
\|\nabla (f-\interp_T^{m-1} f_{H'})\|_{L^p(T)} &\leq & \|H(z)\|^{\frac 1 2} \|H(z)^{-\frac 1 2}(\nabla (f-\interp_T^{m-1} f_{H'}))\|_{L^p(T)}\\
&= &  \|H(z)\|^{\frac 1 2} (\det H(z))^{-\frac 1 {2p}}  \|\nabla (f^\Phi-\interp_{T_\Phi}^{m-1} f^\Phi_{H_\Phi'})\|_{L^p(T_\Phi)}\\
& \leq & C\|H(z)\|^{\frac 1 2} (\det H(z))^{-\frac 1 {2p}} e_{H_\Phi}^a(\nabla f^\Phi \ssep 0)_p\\
& = & C\|H(z)\|^{\frac 1 2} e_H^a(H(z)^{-\frac 1 2}(\nabla f) \ssep z)_p\\
&=&e_H^g(f\ssep z)_p,
\end{eqnarray*}
where we used the change of variables $\Phi^{-1} : T \to T_\Phi$ in the second and fourth lines. This concludes the proof of this lemma.\sq

We now conclude the proof of Proposition \ref{propErrMeshMetG}.
For any $T\in \cT$ and any $z\in T$,  recalling that $\cH_T\leq (4C_0)^2 H(z)$ we obtain 
$$
|\TEq|/|T| =  \sqrt{\det \cH_T} \leq (4C_0)^d  \sqrt{\det H(z)}.
$$
Averaging  \iref{eqLocGW1P} over $T$, it follows that 
$$
\|\nabla (f-\interp_T^{m-1} f)\|_{L^p(T)}^p \leq \frac C {|T|} \int_T e_H^g(f\ssep z)_p^p dz \leq C'  \int_T \sqrt{\det H(z)}\,  e_H^g(f\ssep z)_p^p dz
$$
where $C' := (4C_0)^d C/|\TEq|$. 

We denote by $\cT_*$ a system of representatives of the set $\cT$ for the relation of equivalence  
$$
T \sim T' \stext{ if } T =T'+u \text{ for some } u \in \Z^d.
$$
We thus obtain
\begin{eqnarray*}
\|\nabla (f-\interp_T^{m-1} f_{H'})\|_{L^p([0,1]^d)}^p &=& \sum_{T\in \cT_*} \|\nabla (f-\interp_T^{m-1} f_{H'})\|_{L^p(T)}^p \\
&\leq& C'  \sum_{T\in \cT_*} \int_T \sqrt{\det H(z)}\,  e_H^g(f\ssep z)_p^p dz\\
&=& C'  \int_{[0,1]^d} \sqrt{\det H(z)}\,  e_H^g(f\ssep z)_p^p dz
\end{eqnarray*}
which concludes the proof of this proposition.

\section{Asymptotic approximation and explicit metrics}

This section is devoted to the explicit construction of a metric $H$ adapted to a given function $f$, in the setting where the function $f$ is sufficiently smooth and the mass of metric is asymptotically large. 

Our main result is the following theorem, which is the counterpart for metrics of the results developed in Chapters 2 and 3. It involves the shape functions $K$, $L_a$ and $L_g$ which are defined in the introduction of this chapter, \S \ref{subsecIntroAsympt}. We do not establish that these estimates are optimal, but this is suggested by the lower bounds for the finite element interpolation error developed Chapters 2 and 3 for admissible sequences of triangulations.
\begin{theorem}
\label{thSmooth}
There exists a constant $C=C(m,d)$ such that for each $f\in C^m_\per$
$$
\limsup_{M\to \infty} \left(M^{\frac m d} \inf_{\substack{H \in \bH_g\\  m(H)\leq M}} e_H(f)_p\right) \leq C \|K(d^m f)\|_\tau, \\
$$
where $\frac 1 \tau := \frac m d+ \frac 1 p$.
Furthermore 
\begin{eqnarray*}
\limsup_{M\to \infty}\left(M^\frac {m-1} d \inf_{\substack{H \in \bH_a\\  m(H)\leq M}}  e_H^a(\nabla f)_p\right) &\leq& C \|L_a(d^m f)\|_\tau\\
\limsup_{M\to \infty}\left(M^\frac {m-1} d \inf_{\substack{H \in \bH_g\\  m(H)\leq M}}  e_H^g(\nabla f)_p\right) &\leq& C \|L_g(d^m f)\|_\tau\\
\end{eqnarray*}
where $\frac 1 \tau := \frac {m-1} d+ \frac 1 p$.
\end{theorem}

Our first lemma shows that the optimization problems on $\SL_d$ defining the shape functions $K$, $L_a$ and $L_g$ can be reformulated into optimization problems posed on the collection $S_d^+$ of symmetric positive definite matrices.

\begin{lemma}
\label{lemmaAM}
Let $\psi \in C^0( \GL_d, \R_+^*)$ be homogeneous of degree $r>0$ and such that $\psi(A) = \psi(AO)$ for any $A\in \GL_d$ and any $O\in \cO_d$. Then 
$$
\inf \{(\det M)^\frac r {2d} \sep M\in S_d^+ \text{ and } \psi(M^{-\frac 1 2}) \leq 1\} = \inf \{\psi(A) \sep A\in \SL_d\}.
$$
\end{lemma}

\proof
To each $M\in S_d^+$ we associate $A := (\det M)^\frac 1 {2d} M^{-\frac 1 2}\in \GL_d$, in such way that $\det(A) = 1$ and $\psi(A) = (\det M)^{\frac r {2d}} \psi(M^{-\frac 1 2})\leq (\det M)^{\frac r {2d}}$. Conversely to  each $A\in \GL_d$ we associate $M := \psi(A)^\frac 2 r (A A^\trans)^{-1}\in S_d^+$, in such way that $(\det M)^\frac r{2d} = \psi(A)$ and $\psi(M^{-\frac 1 2}) = \psi(\psi(A)^{-\frac 1 r} (A A^\trans)^\frac 1 2) = \psi(A)^{-1}\psi(AO) = 1$ where $O := A^{-1} (A A^\trans)^\frac 1 2$ is orthogonal. 
\sq

It follows from this lemma that for any $\pi \in \H_m$
\begin{eqnarray}
\label{defKMAT}
K(\pi) &=& \inf\{(\det M)^{\frac m {2d}} \sep \|\pi\circ M^{-\frac 1 2}\|\leq 1\}\\
\label{defLAM}
L_a(\pi) &=& \inf \{ (\det M)^{\frac {m-1} {2d}} \sep \|(\nabla \pi)\circ M^{-\frac 1 2}\| \leq 1\}\vspace{2mm}\\
\label{defLGM}
L_g(\pi) &=& \inf \{ (\det M)^{\frac {m-1} {2d}} \sep \sqrt{\|M\|} \|\nabla (\pi\circ M^{-\frac 1 2})\| \leq 1\}.
\end{eqnarray}

\begin{remark}
The optimizations among symmetric matrices appearing in these expressions have a geometrical interpretation: finding the ellipsoid of maximal volume included in a set of algebraic boundary.
Indeed one easily checks that the three following properties are equivalent for all $\pi \in \H_m$ and all $M\in S_d^+$
\begin{enumerate}[i)]
\item $\|\pi\circ M^{-\frac 1 2}\|\leq 1$. (resp. $\|(\nabla \pi)\circ M^{-\frac 1 2}\|\leq 1$, resp. $\sqrt{\|M\|} \|\nabla (\pi\circ M^{-\frac 1 2})\| \leq 1$)
\item $|\pi(z)|\leq \|z\|_M^m$ for all $z\in \R^d$. (resp. $|\nabla \pi(z)|\leq \|z\|_M^{m-1}$, resp. $\|M\|^\frac 1 2 |M^{-\frac 1 2}(\nabla \pi(z))| \leq \|z\|_M^{m-1}$)
\item The ellipsoid defined by the inequality $\|z\|_M\leq 1$, $z\in \R^d$, is included in the set of algebraic boundary defined by the inequality $|\pi(z)|\leq 1$. (resp. $|\nabla \pi(z)|\leq 1$, resp. $\|M\|^\frac 1 2 |M^{-\frac 1 2}(\nabla \pi(z))|\leq 1$ which depends on $M$)
\end{enumerate}
\end{remark}

We consider a fixed function $f\in C^m$ and for each $z\in \R^d$ we define the homogeneous polynomial $\pi_z\in \H_m$ as follows
\be
\label{defPiZAT}
\pi_z(Z) := \sum_{|\alpha| = m} \frac{\partial^\alpha f(z)}{(\partial z)^\alpha} \frac{Z^\alpha}{\alpha!} .
\ee
where $\alpha! = \alpha_1! \cdots \alpha_d!$ for $\alpha = (\alpha_1, \cdots, \alpha_d) \in \Z_+^d$.

The heuristic guideline of the proof of Theorem \ref{thSmooth} is the following (here in the case of approximation in the $L^p$ norm). 
Consider a $C^\infty$ periodic map $M : \R^d \to S_d^+$ such that $\|\pi_z \circ M(z)^{-\frac 1 2}\|\leq 1$ for all $z\in \R^d$, and define 
$$
H(z) := (\det M(z))^{\frac {-1}{mp+d}} M(z).
$$
Then is is not difficult to show that for all $\lambda$ \emph{sufficiently large} one has $\lambda H \in \bH_a$ and 
$$
m(\lambda H)^{\frac m d} e_{\lambda H}(f)_p \leq C \left\|(\det M)^\frac m {2d}\right\|_{L^\tau([0,1]^d)}.
$$
This leads to the estimate announced  in Theorem \ref{thSmooth} if $\left\|(\det M)^\frac m {2d}\right\|_{L^\tau([0,1]^d)}$ is comparable to $\left\|K(d^m f)\right\|_{L^\tau([0,1]^d)}$, in other words if $(\det M(z))^\frac m {2d}$ is comparable to $K(\pi_z)$ for all $z\in \R^d$, which means that $M(z)$ is a close to be a minimizer of the infimum \iref{defKMAT} defining $K(\pi_z)$. The same principles apply to the two other estimates stated in Theorem \ref{thSmooth}.

The main ingredient of the proof of Theorem \ref{thSmooth} is the construction of a continuous map $\pi \mapsto \cM(\pi)$ such that  $\|\pi\circ \cM(\pi)^{-\frac 1 2}\|\leq 1$ and such that $(\det \cM(\pi))^\frac m {2d}$ is sufficiently close to $K(\pi)$ (in the case of the approximation in the $L^p$ norm). Up to a few technicalities, which are adressed in \S\ref{subsecExplicitMetricLP}, defining $M(z) := \cM(\pi_z)$ and reasoning as above then concludes the proof. We give in \S \ref{subsecExplicit} an explicit expression of such a map in terms of the coefficients of $\pi$, in the case of piecewise linear and piecewise quadratic approximation in two dimensions. Such explicit expressions are valuable for numerical implementation in finite element software which use anisotropic meshes, and have been implemented by the author in \cite{FreeFem}.

We introduce in \S \ref{subsubsecGeomCvx} the property of ``geometric convexity'' for functions in $C^0(S_d^+, \R_+^*)$, which is a variant of the classical notion of convexity. This is used in \S \ref{subsecWellPosed} to define some variants $K^{(\alpha)}$, $\alpha\in (0,\infty)$ of the shape function $K$ which are defined by well posed optimization problems on $S_d^+$, in arbitrary dimension $d$ and degree $m$. In contrast for some choices of $\pi\in \H_m$ there exists no minimizer to the optimization problem appearing in the expression \iref{defKMAT} of $K(\pi)$, and other choices there exists an infinity of minimizers. 
We conclude the proof of Theorem \ref{thSmooth} in \S\ref{subsecExplicitMetricLP}. Eventually \S\ref{subsecKCont} is devoted to the proof of an additional property of the shape functions $K$, $L_a$ and $L_g$ : they are uniformly equivalent on $\H_m$ to some continuous functions.

\subsection{Explicit Metrics}
\label{subsecExplicit}

We give in this section some explicit minimizers, or near minimizers (in a sense defined below) of the optimization problems among symmetric matrices appearing in the expressions \iref{defKMAT}, \iref{defLAM} and \iref{defLGM} of the shape functions $K$, $L_a$ and $L_g$. Our results are so far limited to the case $m=2$, which corresponds to piecewise linear approximation, and to $m=3$ in dimension $d=2$, which corresponds to piecewise quadratic approximation.

These minimizers or near-minimizers have already been presented in Chapter 2 for the shape function $K$, and Chapter $3$ for the shape function $L_a$. The latter are recalled in the next proposition, and completed with their counterparts for the shape function $L_g$.

For any homogeneous quadratic polynomial $\pi$ we denote by $[\pi]\in S_d$ the associated symmetric matrix, which satisfies $z^\trans [\pi] z = \pi(z)$ for all $z\in \R^d$. Note that $\|\pi\| = \|[\pi]\|$.

\begin{prop}
\label{propLm}
\begin{enumerate}[i)]
\item (Piecewise linear approximation, in any dimension)
If $m=2$ and $d\geq 2$, then for any $\pi\in \H_2$, one has 
$$
L_a(\pi) = 2|\det \pi|^\frac 1 d \stext{ and } L_g(\pi) =2 \sqrt{\|\pi\|} \sqrt[2d]{|\det \pi|}.
$$
The minimizing matrices $M$ in the expressions \iref{defLAM} and \iref{defLGM} of $L_a(\pi)$ and $L_g(\pi)$ are  respectively
$$
\cM_a(\pi) = 4[\pi]^2 \stext{ and } \cM_g(\pi) = \|\pi\| \, |[\pi]|.
$$
\item (Piecewise quadratic approximation, in two dimensions)
If $m=3$ and $d=2$, let $\pi =  a x^3+ 3 b x^2 y+ 3 c x y^2+ d y^3\in \H_3$,
let 
$$
\cM_a(\pi) :=  
\sqrt{
\left(\begin{array}{cc}
a & b\\
b & c
\end{array}\right)^2
+
\left(\begin{array}{cc}
b & c\\
c & d
\end{array}\right)^2
}
$$
and let 
$$
\cM_g(\pi) := \cM_a(\pi) + \left(\frac{-\disc\pi}{\|\pi\|}\right)_+^{\frac 1 3}\Id
$$
where $\disc \pi$ is the discriminant of the cubic polynomial $\pi$ :  
\be
\label{eqDiscDet}
 \disc \pi= 4(ac-b^2)(bd-c^2) - (ad-bc)^2.
\ee

Then $\cM_a(\pi)$ and $\cM_g(\pi)$ are near minimizers of optimization problems defining $L_a(\pi)$ and $L_g(\pi)$ respectively, in the following sense: there exists a constant $C$ such that for all $\pi \in \H_3$
$$
\sqrt{\det \cM_a(\pi)}\leq C L_a (\pi) \stext{ and }  \sqrt{\det \cM_g(\pi)} \leq  L_g(\pi).
$$
Furthermore if $\pi$ is not univariate then $\cM_a(\pi)$ and $\cM_g(\pi)$ are non degenerate and we have
\be
\label{eqConstraint}
\|(\nabla \pi)\circ \cM_a(\pi)^{-\frac 1 2}\|\leq C \stext{ and } \|\cM_g(\pi)\|^\frac 1 2 \|\nabla (\pi\circ \cM_g(\pi)^{-\frac 1 2})\|\leq C.
\ee
\end{enumerate}
\end{prop}

\proof
We recall that $\nabla (\pi\circ A) = A^\trans ((\nabla \pi) \circ A)$ for any $\pi \in \H_m$ and $A\in \GL_d$.

We begin with the case i) of a quadratic polynomial $\pi\in \H_2$, and we observe that $\nabla \pi(z) = 2[\pi] z$. Hence $|\nabla \pi(z)|^2 = 4 z^\trans [\pi]^2 z$ and the announced result for $L_a$ easily follows.
For $L_g$, we are looking for a matrix $M\in S_d^+$ such that 
\be
\label{eqMPi}
\|M^\frac 1 2\| \|\nabla( \pi \circ M^{-\frac 1 2})\| = 2\|M^{\frac 1 2}\| \, \|M^{-\frac 1 2} [\pi] M^{-\frac 1 2}\|\leq 1,
\ee
and of minimal determinant.
Recalling that $\|AB\|\geq \|A\| \|B^{-1}\|^{-1}$ for all $A,B \in \GL_d$, therefore $\|M^{-\frac 1 2} [\pi] M^{-\frac 1 2}\|\geq \|\pi\| /\|M\|$, we obtain from \iref{eqMPi} that $2\|\pi \| \leq \|M^\frac 1 2\|$.
Hence 
\be
\label{eqMPiWeak}
4\|\pi\| \| M^{-\frac 1 2} [\pi] M^{-\frac 1 2}\|\leq 1.
\ee
A matrix of minimal determinant which achieves this inequality is clearly $M = 4\|\pi\| \ |[\pi]|$, where $|S|$ denotes the absolute value of a symmetric matrix $S\in S_d$. Observing that this matrix also satisfies the condition \iref{eqMPi} we conclude the proof of Point i).\\

We now turn to the case ii) of cubic polynomial in $\H_3$.
The case of $L_a$ and $\cM_a$ is exposed in Chapter \ref{chapW1P}. We therefore turn to the case of $L_g$. 
Combining the homogeneity and the continuity of the functions involved, one easily shows that there exists a constant $C_0$ such that for all $\pi\in \H_3$
\be
\label{eqContHomog}
\begin{array}{ccccc}
C_0^{-1} \|\pi \| &\leq& \|\nabla \pi \| &\leq& C_0 \| \pi\|\\
C_0^{-1} \|\pi \| &\leq&  \lambda(\pi) :=\|\cM_a(\pi) \| &\leq& C_0 \| \pi\|,
\end{array}
\ee
and 
\be
\label{eqContHomog2}
\begin{array}{ccc}
\mu(\pi) := \|\cM_a(\pi)^{-1} \|^{-1} & \leq& C_0 \|\pi\|\\
\nu(\pi) := \left(\frac{-\disc\pi}{\|\pi\|}\right)_+^\frac 1 3 &\leq&  C_0 \|\pi\|.
\end{array}
\ee
We consider below a fixed polynomial $\pi\in \H_3$, and we recall that 
$$
L_g(\pi) := \inf \{\sqrt{\det M} \sep \sqrt{\|M\|} \|\nabla (\pi \circ M^{-\frac 1 2})\| \leq 1\}.
$$

We first prove the lower bound $L_g(\pi)\geq c \sqrt{\det \cM_g(\pi)}$, where $c>0$ is an absolute constant.
For that purpose we consider a matrix $M\in S_d^+$ satisfying the constraint $\|M^{\frac 1 2}\| \| \nabla (\pi \circ M^{-\frac 1 2}) \|\leq 1$. Combining this with \iref{eqContHomog} we obtain
$$
1 \geq C_0^{-1}  \|M^{\frac 1 2}\| \| \pi \circ M^{-\frac 1 2} \| \geq C_0^{-1} \|M^{\frac 1 2}\| \frac{\| \pi \|}{\| M^{\frac 1 2}\|^3} = \frac{\|\pi\|}{C_0 \|M\|},
$$
hence $\|\pi\| \leq C_0 \|M\|$. It follows that 
$$
1 \geq C_0^{-1} (C_0^{-1} \|\pi\|)^{\frac 1 2} \| \pi\circ M^{-\frac 1 2}\| = \|\ti \pi \circ M^{-\frac 1 2}\|,
$$
where $\ti \pi := C_0^{-\frac 3 2} \| \pi\|^{\frac 1 2} \pi$. 
The solution $\ti M$ of the minimization problem
$$
\inf\{ \det M \sep \|\ti \pi \circ M^{-\frac 1 2}\| \leq 1\}
$$
is known exactly in the case of cubic polynomials in two dimensions, see Proposition \ref{propEllipseMax} in Chapter \ref{chapOptAniso}, and satisfies 
$$
(\det \ti M)^3 = \kappa |\disc \ti \pi| = \kappa (C_0^{-\frac 3 2} \| \pi\|^{\frac 1 2})^4 |\disc \pi| = \kappa C_0^{-6} \| \pi\|^2 |\disc \pi|.
$$
where $\kappa$ is a positive constant which depends only on the sign of $\disc \ti \pi$.
Defining $c_0:=\ti \kappa^\frac 1 6 C_0^{-1}$, where $\ti \kappa$ stands for the minimial value of $\kappa$ for the two possible signs, we obtain 
\begin{eqnarray*}
L_g(\pi) &\geq& \sqrt{\det \ti M} \\
&= & \kappa^\frac 1 6 C_0^{-1} \|\pi\|^{\frac 1 3} |\disc \pi|^{\frac 1 6} \\
&\geq& c_0 \sqrt{\|\pi\| \nu(\pi)} 
\end{eqnarray*}
On the other hand, there exists $c_1>0$ such that 
$$
L_g(\pi) \geq L_a(\pi) \geq c_1 \sqrt{\det \cM_a(\pi)} = c_1\sqrt{\lambda(\pi) \mu(\pi)}.
$$
Furthermore we have 
\begin{eqnarray*}
\det \cM_g(\pi) &=& (\lambda(\pi)+\|\pi\|)(\mu(\pi)+\nu(\pi)) \\
&\leq& (1+C_0)\min \{\lambda(\pi),\, \|\pi\|\} \left( 2\max\{\mu(\pi),\, \nu(\pi)\}\right)\\
&\leq & 2(1+C_0) \max \{\lambda(\pi) \mu(\pi), \, \|\pi\| \nu(\pi)\}.
\end{eqnarray*}
We thus obtain
$$
\sqrt{\det \cM_g(\pi)}\leq C_* L_g(\pi) \stext{ where } C_*:= 2(1+C_0) \max \{c_0^{-1}, \, c_1^{-1}\}.
$$
which concludes the proof of our lower estimate for $L_g$.\\

We now turn to the proof of the property \iref{eqConstraint}, which is equivalent to the following 
$$
\|\cM_g(\pi)\|^\frac 1 2 \|\cM_g(\pi)^{-\frac 1 2}( \nabla \pi(x,y))\| \leq C \stext{ when } \|(x,y)\|_{\cM_g(\pi)} \leq 1.
$$
For any rotation $U\in \cO_2$ we have $\cM_a(\pi \circ U) = U^\trans \cM_a(\pi) U$, see Proposition \ref{propInvMRot} in Chapter \ref{chapW1P}, and $\disc (\pi \circ U) = (\disc \pi) (\det U)^6 = \disc \pi$, see Chapter \ref{chapOptAniso}. We may therefore assume, up to a rotation, that $\cM_a(\pi)$ is a diagonal matrix and that the first diagonal  coefficient is larger than the second one.
Therefore 
$$
\lambda(\pi)^2 = a^2+2b^2+c^2 \stext{ and } \mu(\pi)^2 = b^2+2c^2+d^2,
$$
and 
$\cM_a(\pi) = 
\left(
\begin{array}{cc}
\lambda(\pi) & 0 \\
0 & \mu(\pi)
\end{array}
\right)
$.
Note that this matrix is degenerate if and only if $\mu(\pi) = 0$, which means that $\pi$ is the univariate polynomial $a x^3$.
In order to avoid notational clutter, and since the polynomial $\pi\in \H_3$ is fixed, we denote below by $\lambda$, $\mu$ and $\nu$ the quantities $\lambda(\pi)$, $\mu(\pi)$ and $\nu(\pi)$.
Our purpose is to establish an upper bound on the quantity
\be
\label{eqMGGradPi}
\begin{array}{cl}
& \displaystyle\|\cM_g(\pi)\|^\frac 1 2 \|\cM_g(\pi)^{-\frac 1 2}( \nabla \pi(x,y))\| \\
=&  \displaystyle(a x^2+ 2 b xy + c y^2)^2 + \left(\frac {\lambda+\nu} {\mu+\nu}\right) (b x^2+ 2 c xy+ d y^2)^2.
\end{array}
\ee
under the hypothesis 
$$
\|(x,y)\|_{\cM_g(\pi)}^2 = (\lambda+\nu) x^2+ (\mu+\nu) y^2 \leq 1,
$$
which implies in particular $\lambda x^2\leq 1$ and $\max\{\mu,\nu\} y^2\leq 1$.
Injecting this in \iref{eqMGGradPi}, and observing using \iref{eqContHomog} and \iref{eqContHomog2} that $\nu\leq C_0^2 \lambda$, we see that it is sufficient to bound the following quantities 
\be
\label{eqabcLl}
\frac {a^2} {\lambda^2}, \, \frac {b^2} {\lambda \mu}, \, \frac {c^2} {\mu^2}
\stext{ and }
\frac \lambda \mu \frac {b^2} {\lambda^2}, \ \frac \lambda \mu \frac {c^2} {\mu \lambda}, \ \frac \lambda {\max\{\mu,\nu\}} \frac {d^2} {\max\{\mu, \nu\}^2}.
\ee
Observing that 
\be
\label{eqLlMax}
\mu \geq \max\{b,c,d\}
\stext{ and }
\lambda \geq \max\{a,b,c,\mu\} \geq \max\{a,b,c,d\}
\ee
we find that all the quantities appearing in \iref{eqabcLl} are smaller than one, except perhaps the last one : 
%
\be
\label{eqd}
\frac {\lambda d^2}{\max\{\mu, \nu\}^3}
\ee
Using the expression \iref{eqDiscDet} of the discriminant we find that 
$$
-\disc \pi = (a^2+b^2+c^2) d^2 + a Q(b,c,d)+R(b,c,d)
$$
where $Q$ and $R$ are homogenous polynomials of degree $3$ and $4$ respectively. Hence there exists a constant $C_1$, independent of $\pi$, such that 
\begin{eqnarray*}
\lambda^2 d^2 &=& -\disc \pi - (a Q(b,c,d)+R(b,c,d)) \\
&\leq &\|\pi\| \nu^3 +C_1 \lambda \mu^3 \\
&\leq &\left(C_0 + C_1\right) \lambda \max \{\mu, \nu\}^3.
\end{eqnarray*}
The previous inequality yields a uniform bound on \iref{eqd}, which concludes the proof of this proposition.
\sq

Proposition \iref{propLm} immediately gives an explicit expression, up to a fixed multiplicative constant, of the shape functions $L_a$ and $L_g$ in the case $m=3$ and $d=2$ : $L_a(\pi) \simeq \sqrt {\cM_a(\pi)}$ and $L_g(\pi) \simeq \sqrt {\cM_g(\pi)}$. It follows that
$$
L_g(a x^3+ d y^3) \simeq \sqrt[3]{|ad| \max\{|a|,\, |d|\}} \stext{ and } L_a(a x^3+ d y^3) \simeq \sqrt{|ad|},
$$
and 
$$
L_g(ax^3+3cx y^2) \simeq L_a(a x^3+3c xy^2) \simeq \sqrt{\max\{|a|,\, |c|\} |c|},
$$
where $f(a,b,c,d) \simeq g(a,b,c,d)$ means that there exists a constant $C\geq 1$ such that $C^{-1} f(a,b,c,d) \leq g(a,b,c,d) \leq C f(a,b,c,d)$ for all $a,b,c,d\in \R$.

For a function which has anisotropic features of the type $a x^3+ d y^3$, where $a$ and $d$ have different orders of magnitude, the use of quasi-acute meshes and metrics may therefore lead to a substantial improvement of the approximation error compared to the use of graded meshes and metrics. This is not the case in contrast if all the anisotropic features of the approximated function are of the type $ax^3+3cx y^2$.

\subsection{Geometric convexity} 
\label{subsubsecGeomCvx}
We introduce in this section the geometric average of two symmetric positive definite matrices, and the property of geometric convexity for functions defined on $S_d^+$. This notion is a variant of the classical property of convexity, which is appropriate for the study of the minimization problems appearing in the expressions \iref{defKMAT}, \iref{defLAM} and \iref{defLGM} of the shape functions $K$, $L_a$ and $L_g$.

Let $M, M'\in S_d^+$, and let $S:=M^{-\frac 1 2}$. We define the geometric average $\Avg(M,M')\in S_d^+$ of $M$ and $M'$ as follows
\be
\label{defAvg}
S \Avg(M,M') S := \sqrt{S M' S}
\ee
For instance $\Avg(M,M^{-1}) = \Id$, and for any $k,k'>0$
$$
\Avg(k M , k'M') = \sqrt{k k'} \Avg(M,M').
$$

The following proposition gives an alternative characterization \iref{eqAvgAD} of $\Avg(M,M')$. This proposition also shows that $\Avg(M,M')$ is a natural midpoint between $M$ and $M'$ for the distance $d_\times$, and establishes two basic properties of the geometric average of matrices.

\begin{prop}
For any $M,M'\in S_d^+$ the following holds.
\begin{enumerate}
\item A matrix $\ti M\in S_d^+$ satisfies $\tilde M = \Avg(M,M')$ if and only if there exists $A \in \GL_d$ and a diagonal matrix $D$ such that 
\be
\label{eqAvgAD}
\tilde M = A^\trans A, \quad M = A^\trans e^{2D} A, \quad M' = A^\trans e^{-2D} A.
\ee
\item The geometric average $\Avg(M,M')$ is a midpoint between $M$ and $M'$ for the distance $d_\times$ on $S_d^+$:
$$
d_\times (M, \Avg(M,M')) = d_\times (\Avg(M,M'), M') = \frac 1 2 d_\times (M,M').
$$
\item The geometric average is commutative and compatible with the inversion of matrices:
\begin{eqnarray}
\label{eqAvgCommut}
\Avg(M,M') &=& \Avg(M',M)\\
\label{eqAvgInv}
\Avg(M,M')^{-1} &=& \Avg(M^{-1}, M'^{-1})
\end{eqnarray}
\end{enumerate}
\end{prop}

\proof
We first establish point 1.
Let us assume in a first time that $M,M'$ and $\tilde M$ have the form \iref{eqAvgAD}. We easily obtain 
$$
\tilde M M^{-1} \tilde M = M'
$$
Defining $S := M^{-\frac 1 2}$ we thus have 
$$
(S\tilde M S)^2 = S M' S.
$$
Taking the square root of the previous equation we obtain \iref{defAvg} which shows that $\tilde M = \Avg(M,M')$.

In order to establish the converse of this property we only need to show that for any $M,M'\in S_d^+$ there exists $A\in \GL_d$ and a diagonal matrix $D$ such that $M = A^\trans e^{2D} A$ and $M' = A^\trans e^{-2D} A$. The expression of $\Avg(M,M')$ then automatically follows from the previous argument.
Since the matrix $M'^{-\frac 1 2} M M'^{-\frac 1 2}$ is symmetric positive definite, there exists $U\in \cO_d$ and a diagonal matrix $D$ such that $M'^{-\frac 1 2} M M'^{-\frac 1 2}=U^\trans \exp (4D) U$. Choosing $A := e^D U M'^{\frac 1 2}$ we obtain the announced result.\\

We now turn to the proof of Point 2. Let $A\in \GL_d$ and let $D,D'$ be diagonal matrices, of diagonal coefficients $D_1, \cdots ,D_k$ and  $D'_1, \cdots , D'_d$ respectively. We define
$$
\delta := \max_{1\leq i \leq d} |D_i - D'_i|.
$$ 
We claim that the matrices $N := A^\trans e^{2D} A$ and $N':=A^\trans e^{2D'} A$ satisfy $d_\times(N,N') = \delta$.
Indeed
\begin{eqnarray*}
d_\times (N,\, N') &:= &\sup_{u\neq 0} \left|\ln \|u\|_N -\ln \|u\|_{N'}\right| \\
&=& \sup_{u\neq 0} \left|\ln |e^D Au| - \ln |e^{D'} Au|\right|\\
&=& \sup_{v\neq 0}  \left|\ln |e^D v| - \ln |e^{D'} v|\right|.
\end{eqnarray*}
For any $v=(v_1, \cdots ,v_d)\in \R^d$ we have 
$$
|e^D v|^2 = \sum_{1\leq i \leq d} e^{2D_i} v_i^2 \leq e^{2\delta}  \sum_{1\leq i \leq d} e^{2D'_i} v_i^2 = \left(e^\delta |e^{D'}v|\right)^2,
$$
and similarly $|e^{D'} v|\leq e^\delta |e^D v|$. Therefore $\left|\ln |e^D v| - \ln |e^{D'} v|\right|\leq \delta$, and taking the supremum among all $v\neq 0$ we obtain $d_\times(N,N') \leq \delta$. Furthermore denote by $i_0$, $1\leq i_0\leq d$, the position such that $|D_{i_0} - D'_{i_0}| = \delta$. Choosing $v = (0,\cdots ,1, \cdots ,0)$, with the nonzero coordinate at the position $i_0$, we obtain $\delta = \left|\ln |e^D v| - \ln |e^{D'} v|\right| \leq d_\times(N,N')$. We thus conclude that 
$$
d_\times (N,N') =  d_\times (A^\trans e^{2D} A, \, A^\trans e^{2D'} A) =\delta = \|D-D'\|.
$$

Recalling that $M,M'$ and $\ti M := \Avg(M,M')$ have the form \iref{eqAvgAD} we thus obtain 
$$
d_\times (M, \Avg(M,M')) = d_\times (\Avg(M,M'), M')= \|D\| \stext{ and } d_\times (M,M') = 2 \|D\|.
$$
which concludes the proof of this point.\\

Finally we establish Point 3.
Using the first point we obtain that $M,M'$ and $\tilde M := \Avg(M,M')$ have the form
$$
\tilde M = A^\trans A, \quad M = A^\trans e^{2D} A, \quad M' = A^\trans e^{-2D} A
$$
It follows that 
\be
\label{charCommut}
\tilde M = A^\trans A, \quad M' = A^\trans e^{-2D} A, \quad M = A^\trans e^{2D} A,
\ee
and 
\be
\label{charInv}
\tilde M^{-1} = (A^{-1})^\trans A^{-1}, \quad M^{-1} = (A^{-1})^\trans e^{-2D} A^{-1}, \quad M'^{-1} = (A^{-1})^\trans e^{2D} A^{-1}.
\ee
Equation  \iref{eqAvgCommut} (resp. \iref{eqAvgInv}) is the consequence of \iref{charCommut} (resp. \iref{charInv}) and of the characterization of the geometric average given in the first point.\sq

We recall that a function $\psi\in C^0(S_d, \R)$ is convex if and only if for all $M,M'\in S_d$ one has 
\be
\label{defConv}
\psi\left(\frac{M+M'} 2\right) \leq \frac{\psi(M)+ \psi(M')} 2.
\ee
The next proposition introduces the family of geometrically convex functions on $S_d^+$, which is defined by a property similar to \iref{defConv} but in which the arithmetic averages are replaced with geometric averages.
\begin{definition}
We say that a function $\vp \in C^0(S_d^+ , \R_+^*)$ is Geometrically Convex and Homogenous (GCH), if it satisfies the following properties.
\begin{enumerate}[i)]
\item (Sub multiplicativity) For any $M,M'\in S_d^+$ one has
\be
\label{eqSubMult}
\vp(\Avg(M,M')) \leq \sqrt{\vp(M) \vp(M')}.
\ee
\item (Homogenous) There exists $\alpha \in \R$ , called the degree of $\vp$,
such that for all $M\in S_d^+$ and all $\lambda\in \R_+^*$ one has 
$$
\vp(\lambda M) = \lambda^\alpha \vp(M).
$$
\end{enumerate}
We say that the function $\vp$ is Strictly Geometrically Convex and Homogeneous (SGCH), if $\vp$ is GCH and if the equality in \iref{eqSubMult} only holds if $M$ and $M'$ are proportional.
\end{definition}

For instance the function $\det$ is GCH of degree $d$, since 
$$
\det \Avg(M,M') = \sqrt{\det M \det M'} \stext{ and } \det (\lambda M) = \lambda^d M.
$$
We insist on point that GCH functions are by assumption continuous and strictly positive on $S_d^+$.
The next lemma enumerates some simple properties of geometrically convex functions.
\begin{lemma}
\label{lemmaAlgGCH}
\begin{enumerate}
\item The product of two GCH functions is GCH, and the product of a GCH function with a SGCH function is SGCH. Furthermore the degrees add.
\item The elevation to a positive power $\alpha\in \R_+^*$ of a GCH (resp. SGCH) function is also GCH (resp. SGCH). Furthermore the degree is multiplied by $\alpha$.
\item If a function $\vp$ is GCH (resp. SGCH) then $M \mapsto \vp(M^{-1})$ is also GCH (resp. SGCH), and has the opposite degree.
\end{enumerate}
\end{lemma}
\proof
The first two points are immediate, and the last point is a direct consequence of the compatibility \iref{eqAvgInv} of the geometric average with the inversion : if $\vp$ is GCH and $M,M'\in S_d^+$ then 
$$
\vp(\Avg(M,M')^{-1}) = \vp(\Avg(M^{-1}, M'^{-1})) \leq \sqrt{\vp(M^{-1}) \vp (M'^{-1})},
$$
which concludes the proof.
\sq

Aside from the determinant, a number of functions on $S_d^+$ are geometrically convex. The next proposition, which is based on an argument of complex analysis, enumerates some of them.
For all $\pi\in \H_m$ we define 
$$
\|\pi\|_\sC := \sup\{ |\pi(z)| \sep z\in \C^d, \, |z|\leq 1\}.
$$

\begin{prop}
\label{propExGCH}
\begin{enumerate}
\item The trace map $M \mapsto \Tr M$ is SGCH.
\item For any fixed polynomial $\pi\in \H_m\sm\{0\}$ the map $M \mapsto \|\pi \circ M^{\frac 1 2}\|_\sC$ is GCH, of degree $m/2$.
\item The norm map $M \mapsto \|M\|$ is GCH. 
\end{enumerate}

\end{prop}

\proof
We first establish Point 1.
Let $M,M'\in S_d^+$ and let $\ti M := \Avg(M)$. Let $A\in \GL_d$ and let $D$ be a diagonal matrix such that  $M = A^\trans e^{2D} A$, $M' = A^\trans e^{-2D} A$ and $\ti M = A^\trans A$.
We denote by $(A_{ij})_{1\leq i, j\leq d}$ the coefficients of $A$, and we define for all $1\leq i \leq d$
$$
\cA_i^2 := \sum_{1\leq j \leq d} A_{ij}^2.
$$
Note that $\cA_i>0$ for all $1\leq i \leq d$, since otherwise a full line of $A$ would be zero and $A$ would be degenerate.
Denoting by $D_i , \cdots , D_d$ the diagonal coefficients of $D$ we obtain  
$$
\Tr(M) = \sum_{1\leq i \leq d} \cA_i^2 e^{2D_i}, \quad \Tr(M') = \sum_{1\leq i \leq d} \cA_i^2 e^{-2D_i}, \quad \Tr(\ti M) = \sum_{1\leq i \leq d} \cA_i^2.
$$
It thus follows from the Cauchy Schwartz inequality that $\Tr(\ti M)^2 \leq \Tr(M) \Tr(M')$. Furthermore equality occurs if and only if the vectors 
$$
(\cA_i e^{D_i})_{1\leq i \leq d} \stext{ and }(\cA_i e^{-D_i})_{1\leq i \leq d}
$$
are proportional, which implies that $D_i =D_j$ for all $1\leq i,j\leq d$ and therefore that $M$ and $M'$ are proportional. \\

We now turn to the proof of Point 2., and for that purpose we recall a classical result of complex analysis, called Hadamard's three lines theorem. 
Let $\Omega:=\{\alpha\in \C \sep |\Re(\alpha)|<1\}$, where $\Re$ denotes the real part, and let $g$ be a continuous and bounded function on $\overline \Omega$ which is holomorphic on $\Omega$. Then 
\be
\label{eqComplex}
|g(0)|^2 \leq \left(\sup_{\Re \alpha = -1} |g(\alpha)|\right) \left(\sup_{\Re z = \alpha} |g(\alpha)|\right).
\ee


We consider a fixed $z\in \C^d$ satisfying $|z|\leq 1$, and we define for all $\alpha\in \overline \Omega$ 
$$
g(\alpha) := \pi(A^\trans e^{\alpha D} z).
$$
The function $g$ is continuous, uniformly bounded on $\overline \Omega$ by $\|\pi\|_\sC (\|A\| e^{\|D\|})^d$, and holomorphic on $\C$, hence on $\Omega$. Therefore Hadamard's three lines theorem applies.
For any $\sigma \in \{-1,1\}$ one has 
$$
\sup_{\Re(\alpha) = \sigma} |\pi(A^\trans e^{\alpha D} z)| \leq \|\pi\circ(A^\trans e^{\sigma D})\|_\sC,
$$
since $|e^{\b i\Im(\alpha) D}z| = |z|=1$, where $\Im$ denotes the imaginary part and $\b i$ the imaginary unit. 
Applying \iref{eqComplex} to the function $g$ we thus obtain 
$$
|\pi(A^\trans z)|^2 \leq \|\pi\circ(A^\trans e^D)\|_\sC \|\pi\circ(A^\trans e^{-D})\|_\sC.
$$
Taking the supremum among all $z\in \C^d$ such that $|z|\leq 1$ we conclude that 
$$
\|\pi\circ A^\trans\|_\sC \leq \|\pi\circ(A^\trans e^D)\|_\sC \|\pi\circ(A^\trans e^{-D})\|_\sC.
$$
One easily checks that the matrices $O, O', \tilde O\in \GL_d$ defined by  
$$
OM^\frac 1 2 =  e^D A, \quad O'M'{^\frac 1 2} = e^{-D} A \stext{ and } \tilde O \tilde M^\frac 1 2 = A,
$$
are orthogonal.
Since orthogonal matrices are also unitary matrices, in the sense that $|Oz| = |z|$ for all $z\in \C$, we obtain
\begin{eqnarray*}
\|\pi \circ A^Te^{-D}\|_\sC &=& \|\pi \circ M^{\frac 1 2}\|_\sC,\\
\|\pi \circ A^T e^{-D}\|_\sC &=& \|\pi \circ  M'^{\frac 1 2}\|_\sC\\
\|\pi \circ A^T\|_\sC &=& \|\pi \circ \ti M^{\frac 1 2}\|_\sC,
\end{eqnarray*}
thus $\|\pi\circ M^\frac 1 2\|_\sC$ satisfies the sub-multiplicativity property \iref{eqSubMult}. Furthermore the function $M\in S_d^+ \mapsto \|\pi\circ M^\frac 1 2\|_\sC$ is clearly continuous, strictly positive, and homogeneous of degree $m/2$ which concludes the proof of this point.\\

We now establish Point 3, and for that purpose we denote by $\pi\in \H_2$ the canonical quadratic form 
$$
\pi(z) := \sum_{1\leq j \leq d} z_j^2.
$$
Our purpose is to establish that for any $M\in S_d^+$
\be
\label{eqPiMM}
\|\pi\circ M^{\frac 1 2}\|_\sC = \|M\|,
\ee
which implies that $\|M\|$ is GCH in view of the previous point.
For any $z\in \C$ such that $|z|\leq 1$ one has 
$$
|\pi(M^\frac 1 2 z)| = |z^\trans M^\frac 1 2 M^\frac 1 2 z| \leq |M^\frac 1 2 z|^2 \leq \|M^\frac 1 2\|^2  = \|M\|.
$$
Furthermore choosing a normalized eigenvector $z\in \R^d$ associated to the maximal eigenvalue of $M$ we obtain $|\pi(M^\frac 1 2 z)|=\|M\|$, which establishes \iref{eqPiMM} and concludes the proof.
\sq

For each $M\in S_d^+$ we define
$$
\kappa(M) :=  \frac 1 d\sqrt{\Tr(M) \Tr(M^{-1})}.
$$
The function $\kappa$ is SGCH according to Lemma \ref{lemmaAlgGCH} and since $M \mapsto \Tr M$ is SGCH.
In terms of the eigenvalues $\lambda_1, \cdots, \lambda_d$ of $M$ we have 
$$
d^2 \kappa(M)^2 = \left(\sum_{1 \leq i \leq d} \lambda_i\right)\left(\sum_{1 \leq i \leq d} \lambda_i^{-1}\right),
$$ 
therefore $\kappa(M)\geq 1$, using the Cauchy Schwartz inequality, with equality if and only if $M = m\Id$ for some $m>0$.
Note also that  
$$
\kappa(M)\leq \sqrt{\|M\|\|M^{-1}\|} \leq d \, \kappa(M).
$$
We regard the function $\kappa$ as a close variation on the condition number $\sqrt{\|M\|\|M^{-1}\|}$ of the matrix $M^\frac 1 2$. 

\subsection{A well posed variant of the shape function}
\label{subsecWellPosed}
We define in this section three variants of the shape functions $K$, $L_a$ and $L_g$ : for all $\pi \in \H_m$
\begin{eqnarray*}
K^\sC(\pi) &:=& \inf_{A\in \SL_d} \|\pi\circ A\|_\sC\\
L_a^\sC(\pi) &:=& \inf_{A\in \SL_d} \sqrt{\|G(\pi)\circ A\|_\sC}\\
L_g^\sC(\pi) &:=& \inf_{A\in \SL_d} \|A^{-1}\| \|\pi\circ A\|_\sC,
\end{eqnarray*}
where $G(\pi) := |\nabla \pi|^2 \in \H_{2m-2}$.
For each $\alpha\in (0,\infty)$ we define three new variants of the shape functions
\begin{eqnarray*}
K^{(\alpha)}(\pi) &:=& \inf\{\kappa(M)^{\frac 1\alpha} (\det M)^{\frac m {2d}} \sep \|\pi\circ M^{-\frac 1 2}\|_\sC\leq 1\text{ and } \kappa(M)\leq e^\alpha\},\\
L_a^{(\alpha)}(\pi) &:=& \inf\{\kappa(M)^{\frac 1\alpha} (\det M)^{\frac {m-1} {2d}} \sep \|G(\pi)\circ M^{-\frac 1 2}\|_\sC\leq 1\text{ and } \kappa(M)\leq e^\alpha\},\\
L_g^{(\alpha)}(\pi) &:=& \inf\{\kappa(M)^{\frac 1\alpha} (\det M)^{\frac {m-1} {2d}} \sep \|M\|^\frac 1 2 \|\pi\circ M^{-\frac 1 2}\|_\sC\leq 1\text{ and } \kappa(M)\leq e^\alpha\}.
\end{eqnarray*}
This section is devoted to the proof of the following theorem. 
\begin{theorem}
\label{thWellPosed}
\begin{enumerate}
\item
The shape function $K$ (resp. $L_a$, resp. $L_g$) is uniformly equivalent to $K^\sC$ on $\H_m$ (resp. $L_a^\sC$, resp. $L_g^\sC$).
\item
For each $\pi\in \H_m$ we have the \emph{decreasing} convergence
$$
\lim_{\alpha \to \infty} K^{(\alpha)}(\pi) = K^\sC(\pi).
$$ 
(resp. $L_a^{(\alpha)}(\pi) \to L_a^\sC(\pi)$ decreasingly, resp. $L_g^{(\alpha)}(\pi) \to L_g^\sC(\pi)$ decreasingly, as $\alpha \to \infty$).
\item
For each $\alpha \in (0,\infty)$ and each $\pi \in \H_m\sm\{0\}$ there exists a unique minimizer $\cM^{(\alpha)}(\pi)\in S_d^+$ to the optimization problem defining $K^{(\alpha)}(\pi)$,
and the map 
$$
\pi\in \H_m\sm\{0\} \to \cM^{(\alpha)}(\pi)\in S_d^+
$$
is continuous.
(resp. likewise for $L^{(\alpha)}(\pi)$ and $\cM^{(\alpha)}_a(\pi)$, resp. likewise for $L_g^{(\alpha)}(\pi)$ and $\cM^{(\alpha)}_g(\pi)$)
\end{enumerate}
\end{theorem}

We first establish Point 1. of this theorem.
Since the vector spaces $\H_m$ and $\H_{2m-2}$ have finite dimension, there exists three constants $C_K$, $C_a$ and $C_g$ such that 
\be
\label{defCKAG}
\begin{array}{rcccll}
\|\pi\| &\leq & \|\pi\|_\sC&\leq &C_K \|\pi\|, & \pi \in \H_m,\\
\|\pi'\| &\leq & \|\pi'\|_\sC&\leq &C_a \|\pi'\|, & \pi \in \H_{2m-2},\\
C_g^{-1} \|\nabla \pi\| &\leq & \|\pi\|_\sC&\leq& C_g \|\nabla \pi\|, & \pi \in \H_m.
\end{array}
\ee
Hence it follows from the definition of the original shape functions, given in \S\ref{subsecIntroAsympt}, that $K\leq K^\sC\leq C_K K$, $L_a \leq L_a^\sC \leq C_a L_a$ and $C_g ^{-1}L_g  \leq L_g^\sC \leq C_g L_g$ on $\H_m$.

Point 2. of this theorem is an immediate consequence of the fact that $\kappa(M)\in [1, \infty)$ for all $M\in S_d^+$.

The following proposition illustrates the use of $\kappa$ as a regularization term, and immediately implies Point 3. of Theorem \ref{thWellPosed}.
\begin{prop}
\label{propContLift}
Let $\vp$ be a GSH function of degree $r>0$. 
Let $(X,d_X)$ be a metric space, and let $\psi\in C^0(S_d^+ \times X, \R_+^*)$ be such that $\psi(\cdot,x)$ is GCH of degree $-1$ for any fixed $x\in X$.
We define for all $x\in X$ and all $\alpha \in (0, \infty)$
$$
L^{(\alpha)}(x) := \inf\left\{\kappa(M)^{\frac 1 \alpha} \vp(M) \sep \psi(M,x) \leq 1 \text{ and } \kappa(M) \leq e^\alpha\right\}.
$$
For any $x\in X$ there exists a unique minimizer $\cM^{(\alpha)}(x)$ to the optimization problem defining $L^{(\alpha)}(x)$. Furthermore the map $\cM^{(\alpha)} : X \to S_d^+$ is continuous. 
\end{prop}

\proof
We consider a fixed parameter $\alpha\in (0, \infty)$.
We first establish the existence of a minimizer to the optimization problem defining $L^{(\alpha)}(x)$, and we thus consider a fixed $x\in X$.
For any $\beta \geq 0$ the following collection of matrices is clearly a compact subset of $S_d^+$
\be
\label{defKBeta}
K_\beta := \{M\in S_d^+ \sep |\ln(\det M)|\leq \beta \text{ and } \kappa(M) \leq e^\alpha\}.
\ee
We define
$$
\Lambda := \sup_{M\in K_0} \vp(M) \stext{ and } \lambda := \inf_{M\in K_0} \psi(M,x),
$$
in such way that for any $M\in S_d^+$ satisfying $\kappa(M) \leq e^\alpha$
\be
\label{eqPhiPsi}
\vp(M) \leq \Lambda (\det M)^{\frac r d} \stext{ and } \psi(M,x) \geq \lambda (\det M)^{-\frac 1 d}.
\ee
Let $(M_n)_{n \geq 0}$, $M_n\in S_d^+$, be a minimizing sequence for the optimization problem defining  $L^{(\alpha)}(x)$. It follows from \iref{eqPhiPsi} that $\det (M_n)$ uniformly bounded above and below: hence there exists $\beta \geq 0$ such that $|\ln \det M_n|\leq \beta$ for all $n \geq 0$. 
Furthermore $\kappa(M_n) \leq e^\alpha$ for each $n\geq 0$, hence $M_n \in K_\beta$. Since the set $K_\beta$ is compact there exists a converging subsequence $M_{\sigma(n)} \to M_\infty\in K_\beta$. By continuity of the functions $\vp$, $\psi(\cdot ,x)$ and $\kappa$ on $S_d^+$, the matrix $M_\infty$ is a minimizer of the optimization problem $L^{(\alpha)}(x)$.\\

We now consider two minimizers $M,M'$ to the optimization problem defining $L^{(\alpha)}(x)$, and we intend to show that $M=M'$. We denote $\ti M := \Avg(M,M')$. 
Since the functions $\psi(\cdot, x)$ and $\kappa$ are GCH, the matrix $\ti M$ satisfies the inequalities 
$$
\psi(\ti M,x) \leq 1 \stext{ and } \kappa(\ti M) \leq e^\alpha.
$$ 
We define for each $M\in S_d^+$ the quantity
\be
\label{defPhiSGCH}
\Phi(M) := \kappa(M)^\frac 1 \alpha \vp(M),
\ee
and we observe that $\Phi : S_d^+ \to \R_+^*$ is a SGCH function of degree $r>0$ according to Lemma \ref{lemmaAlgGCH}.
If $M$ is not proportionnal to $M'$ we thus have
$$
\Phi(\ti M) < \sqrt{\Phi( M) \Phi( M')} = L^{(\alpha)}(x),
$$
 which is a contradicts the fact that $M$ and $M'$ are minimizers of the optimization problem defining $L^{(\alpha)}(x)$. Thus $M = m M'$ for some $m>0$, but this implies 
 $$
L^{(\alpha)}(x) = \Phi( M) = \Phi( M') = m^r  \Phi( M),
 $$
 therefore $m=1$ and $M=M'$, which establishes as announced that the minimizer is unique.\\

Finally we establish that the map $x \mapsto \cM^{(\alpha)}(x)$ is continuous.
For all $x,x'\in X$, we define
$$
\omega(x,x') := \sup_{M \in K_0} |\ln \psi(M,x) - \ln\psi(M,x')|,
$$
where the compact set $K_0$ is defined by \iref{defKBeta}, and we observe that $\omega(x,x') \to 0$ as $x'\to x$.
(Note that $\omega(x,x')$ is the modulus of continuity of the continuous function $X \to C^0(K_0, \R)$, defined by $x\mapsto \ln \psi (\cdot,x)$.)

Let $x,x'\in X$ and let $M := e^{\omega(x,x')} \cM^{(\alpha)}(x)$ and $M' := e^{\omega(x,x')} \cM^{(\alpha)}(x')$. We thus have 
$$
\kappa(M) = \kappa(\cM^{(\alpha)}(x)) \leq e^\alpha,
$$
and likewise $\kappa(M') \leq e^\alpha$.
Since $\psi(\cdot ,x')$ is homogeneous of degree $-1$ we obtain
$$
\psi(M, x')\leq e^{\omega(x,x_n)} \psi(M, x) = e^{\omega(x,x')} \psi(e^{\omega(x,x')} \cM^{(\alpha)}(x),x) \leq 1,
$$
and likewise $\psi(M',x)\leq 1$.
Furthermore since the function $\Phi$ defined by \iref{defPhiSGCH} is homogeneous of degree $r$ we obtain
$$
L^{(\alpha)}(x') \leq \Phi(M) = e^{r\omega(x,x')} L^{(\alpha)}(x),
$$
which implies
$$
\Phi(M') = \Phi(e^{\omega(x,x')} \cM^{(\alpha)}(x'))  = e^{r\omega(x,x')} L^{(\alpha)}(x') \leq e^{2r\omega(x,x')} L^{(\alpha)}(x).
$$

Consider a converging sequence in $X$, $x_n \to x$, and define
$$
M_n := e^{\omega(x,x_n)} \cM^{(\alpha)} (x_n).
$$
The arguments above show that $M_n$ is a minimizing sequence for the optimization problem $L^{(\alpha)}(x)$. The arguments developed for the existence of a minimizer show that $\{M_n\sep n \geq 0\}$ is contained in the compact set $K_\beta$ if the constant $\beta\geq 0$ is sufficiently large, and that any converging subsequence of $(M_n)_{n \geq 0}$ tends to a minimizer of the optimization problem $L^{(\alpha)}(x)$, thus to $\cM^{(\alpha)}(x)$ by uniqueness.

Since the sequence $M_n$ takes its values in the compact set $K_\beta$ and has the single adherence value $\cM^{(\alpha)}(x)$, we have the convergence $M_n \to \cM^{(\alpha)}(x)$ as $n \to \infty$. Since $\omega(x,x_n) \to 0$ we obtain as announced 
$$
\cM^{(\alpha)}(x_n)\to \cM^{(\alpha)}(x),
$$
which concludes the proof of this proposition.
\sq

\begin{remark}
\label{remHomog2}
One easily checks that for each $\pi\in \H_m$, each $\alpha\in (0,\infty)$ and each $\lambda\in \R_+^*$ we have $K^{(\alpha)}(\lambda \pi) = \lambda L^{(\alpha)}(\pi)$ and 
$$
\cM^{(\alpha)}(\lambda \pi) = \lambda^\frac 2 m \cM^{(\alpha)}(\pi).
$$
It follows that $\cM^{(\alpha)}(\pi) \to 0$ as $\pi\to 0$, and therefore that the map 
from $\H_m$ to $S_d$ defined by 
$$
\pi \mapsto 
\left\{
\begin{array}{cc}
\cM^{(\alpha)}(\pi) &  \text{ if } \pi \neq 0\\
0 & \text{ if } \pi = 0
\end{array}
\right.
$$
is continuous.
Likewise $L^{(\alpha)}_a$ and $L^{(\alpha)}$ are $1$-homogeneous on $\H_m$, and $\cM^{(\alpha)}_a$ and $\cM^{(\alpha)}_g$ are $\frac 2 {m-1}$ homogeneous on $\H_m$ which implies that they tend to $0$ as $\pi \to 0$, and are continuous on $\H_m$.
\end{remark}

\subsection{Construction of the metric}
\label{subsecExplicitMetricLP}

This subsection is devoted to the construction of a metric $H$ which is \emph{adapted} to a given smooth function $f$ in the sense that the error $e_H(f)_p$ (resp. $e_H^a(\nabla f)_p$, resp. $e_H^g(\nabla f)_p$) is \emph{optimally small} among all metrics of the \emph{same mass} $m(H)$ up to a fixed multiplicative constant. Our construction only applies to an interpolation large mass of the metric $H$. We do not establish that it is optimal, but this is suggested by the lower error estimates given in Chapters 2 and 3 for the finite element interpolation error on admissible sequences of meshes.
 Producing a such a metric of a given mass independent of $f$ is so far an open question.

Our first lemma shows that the regularity constraints defining by the classes $\bH_a$ and $\bH_g$ of metrics ``disappear'' when one multiplies a metric with a large constant. 
\begin{lemma}
For any $H\in\bH_\per$ which is $C^1$, there exists $\lambda_0 = \lambda_0(H)$ such that $\lambda^2 H\in \bH_a$ for all $\lambda\geq \lambda_0$.
\end{lemma}
\proof
The following maps are $C^1$ since they are the composition of $C^1$ maps
$$
\begin{array}{ccl}
z\in \R^d &\mapsto& H(z)^{-\frac 1 2} \in S_d\\
z\in \R^d &\mapsto& \vp_\times (H(z)) := (u \mapsto \ln \|u\|_{H(z)}) \in C^0(\bbS, \R)
\end{array}
$$
where $\bbS := \{u\in \R^d \sep |u|=1\}$ denotes the unit euclidean sphere (see \S\ref{secEigen} in the previous Chapter for a more in depth discussion on the map $\vp_\times$).
Since $H$ is $C^1$ and $\Z^d$-periodic these maps are uniformly Lipschitz, hence there exists two constants $C_+, C_\times$ such that for all $z,z'\in \R^d$
$$
\begin{array}{ccccc}
d_+(H(z), H(z')) &:=& \|H(z)^{-\frac 1 2}-H(z')^{-\frac 1 2}\| &\leq& C_+ |z-z'|,\\
d_\times(H(z), H(z')) &:=& \|\vp_\times(H(z))-\vp_\times(H(z'))\|_{L^\infty(\bbS)} &\leq& C_\times |z-z'|.
\end{array}
$$
Furthermore, since $H$ is continuous and periodic, there exists $\ve>0$ such that $H(z) \geq \ve^2 \Id$ for all $z\in \R^d$. Hence
$$
d_H(z,z') \geq \ve |z-z'|,
$$
for all $z,z'\in \R^d$.
From this point, it follows from Remark \ref{remHomog} in the previous chapter that $\lambda^2 H\in \bH_a$ for all $\lambda \geq \lambda_0 := \max \{C_+, C_\times /\ve\}$.
\sq

Given a function $f\in C^m_\per$, an exponent $1\leq p\leq \infty$, and a triplet  $s=(\alpha, \rho, \delta)\in (\R_+^*)^3$ of parameters we construct a $C^\infty$ metric $H$ as follows.
For each $z\in \R^d$ we define 
$$
M(z) := \cM^{(\alpha)}(\pi_z) + \rho\Id,
$$
where the homogeneous polynomial $\pi_z\in \H_m$ is defined by \iref{defPiZAT}. Note that $M$ is continuous according to Remark \ref{remHomog2}.
We consider a fixed mollifier $\psi$, namely a radial non-negative compactly supported $C^\infty$ function of integral one, and we denote $\psi_\delta := \psi(\cdot/\delta)/\delta^d$. We then define 
\be
\label{defHs}
H(z) := (\det M'(z))^\frac{-1}{mp+d} M'(z) \stext{ where } M' := \psi_\delta*M.
\ee

Given a symbol $\star\in \{a,g\}$ we define similarly $M_\star := \cM^{(\alpha)}_\star(\pi_z) + \rho\Id$, $M'_\star := \psi_\delta*M_\star$ and 
$$
H_\star (z) := (\det M'_\star(z))^\frac{-1}{(m-1)p+d} M'_\star(z).
$$

The rest of this subsection is devoted to the proof of the following theorem, which immediately implies the asymptotic estimates \iref{asymptLPMet}, \iref{asymptA} and \iref{asymptG} announced in the introduction. 
\begin{theorem}
There exists a constant $C=C(m,d)$ such that the following holds for any $f\in C^m_\per$, any $1\leq p\leq \infty$ and any $\ve>0$. If the triplet of parameters $s=(\alpha, \rho, \delta)$ is well chosen then 
\begin{enumerate}
\item Defining $\tau$ by $\frac 1 \tau := \frac m d+ \frac 1p$ we have 
\be
\label{asymptK}
\lim_{\lambda \to \infty} m(\lambda^2 H)^\frac m d e_{\lambda^2 H}(f)_p \leq C\|K(d^m f)\|_{L^\tau([0,1]^d)}+\ve.
\ee
\item Defining $\tau$ by $\frac 1 \tau := \frac {m-1} d+ \frac 1p$ we have
$$
\lim_{\lambda \to \infty} m(\lambda^2 H_a)^\frac {m-1} d e_{\lambda^2 H_a}^a(\nabla f)_p \leq C\|L_a(d^m f)\|_{L^\tau([0,1]^d)}+\ve
$$
and
$$
\lim_{\lambda \to \infty} m(\lambda^2 H_g)^\frac {m-1} d e_{\lambda^2 H_g}^g(f)_p \leq C\|L_g(d^m f)\|_{L^\tau([0,1]^d)}+\ve.
$$
\end{enumerate}
\end{theorem}

The proof of the three estimates are completely similar. We therefore focus on the proof of \iref{asymptK}, and the details of the two other estimates are left to the reader.
For each $z\in \R^d$ we identify the polynomial $\pi_z$ with the collection $d^m f(z)/m!$ of $m$-th order derivatives of $f$ at $z$, and we observe that 
\begin{eqnarray*}
\|\det(\cM^{(\alpha)}(d^m f))^\frac m {2d}\|_{L^\tau([0,1]^d)} &=& \|K^{(\alpha)}(d^m f)\|_{L^\tau([0,1]^d)}\\ 
&\underset {\alpha\to \infty} \longrightarrow &\|K^\sC(d^m f)\|_{L^\tau([0,1])}\\
&\leq & C_K \|K(d^m f)\|_{L^\tau([0,1])}
\end{eqnarray*}
where the constant $C_K$ is defined in \iref{defCKAG}.
In the second line we used the fact that $K^{(\alpha)}(d^m f)$ converges pointwise \emph{decreasingly} to $K^\sC(d^m f)$, which implies the convergence of the integrals.
We may therefore choose $\alpha$ sufficiently large and $\rho$ sufficiently small in such way that 
$$
\|(\det M)^\frac m {2d}\|_{L^\tau([0,1]^d)} < C \|K(d^m f)\|_{L^\tau([0,1])} + \ve,
$$
where $C=C_K/m!$ (recall that $\cM^{(\alpha)}$ and $K^{(\alpha)}$ are homogeneous functions, as observed in Remark \ref{remHomog2}).

For each $z\in \R^d$ such that $\pi_z\neq 0$ we have 
\begin{eqnarray*}
\|\pi_z\circ M(z)^{-\frac 1 2}\| &=& \sup_{|u|=1} \frac{|\pi_z(u)|}{\|u\|_{M(z)}^m}\\
&=& \sup_{|u|=1} \frac {|\pi_z(u)|}{(\|u\|^2_{\cM(\pi_z)} +\rho)^\frac m 2}\\
&<& \sup_{|u|=1} \frac {|\pi_z(u)|}{\|u\|_{\cM(\pi_z)}^m}\\
&=& \|\pi_z\circ \cM(\pi_z)^{-\frac 1 2}\| \\
&\leq& 1.
\end{eqnarray*}
We may therefore choose $\delta>0$ sufficiently small in such way that 
$
\|\pi_z\circ M'(z)^{-\frac 1 2}\| <1
$
for all $z\in [0,1]^d$ 
and such that 
\be
\label{eqMK}
\|(\det M'(z))^\frac m {2d}\|_{L^\tau([0,1]^d)} < C \|K(d^m f)\|_{L^\tau([0,1])} + \ve.
\ee
Note that according to \iref{defHs}
$$
m(H)^\frac 1 \tau =\|(\det M'(z))^\frac m {2d}\|_{L^\tau([0,1]^d)}.
$$
We define for each $z\in \R^d$
\be
\label{defMuzLpMet}
\mu_z := \sum_{|\alpha| < m} \frac{\partial^{|\alpha|} f(z)}{(\partial z)^\alpha} \frac{Z^\alpha} {\alpha!}.
\ee
and for each $\lambda>0$
$$
f_{z,\lambda}(h) := \lambda^m(f(z+h/\lambda) - \mu_z(h/\lambda)).
$$
Using the isotropic change of variables $h \mapsto z+h/\lambda$, we obtain
$$
e_{\lambda^2 H} (f \ssep z)_p = \inf_{\mu\in \sP_{m-1}} \|f - \mu\|_{L^p(B_H(z,1/\lambda))}=  \lambda^{-(m+\frac d p)}\inf_{\mu\in \sP_{m-1}} \|f_{z,\lambda} - \mu\|_{L^p(B_H(z))}.
$$
Consider $R>0$ such that $H \geq \Id/R^2$ uniformly on $[0,1]^d$, hence on $\R^d$. The function $f_{z,\lambda}$ converges uniformly to $\pi_z$ on $B(0,R)$: 
\be
\label{eqConvfPi}
\lim_{\lambda\to \infty}\|f_{z,\lambda} -\pi_z\|_{L^\infty(B(0,R))} = 0,
\ee
since $\pi_z$ is the $m$-th term in the Taylor development of $f$ at $z$. Since $B_H(z) \subset B(0,R)$ we therefore obtain  
\begin{eqnarray*}
\lim_{\lambda \to \infty}  \lambda^{m+\frac d p}e_{\lambda^2 H} (f \ssep z)_p &=& 
\lim_{\lambda\to \infty}  \lambda^{m+\frac d p} \inf_{\mu\in \sP_{m-1}} \|f_{z,\lambda} - \mu\|_{L^p(B_H(z))}\\
&=&\inf_{\mu\in \sP_{m-1}} \|\pi_z - \mu\|_{L^p(B_H(z))}\\
& \leq & \|\pi_z\|_{L^p(B_H(z))}\\
& \leq & |B_H(z)|^\frac 1 p \|\pi_z\|_{L^\infty(B_H(z))}\\
&=& \omega \, (\det H(z))^{-\frac 1 {2p}} \|\pi_z \circ H(z)^{-\frac 1 2}\|\\
&=& \omega \|\pi_z \circ M(z)^{-\frac 1 2}\|\\
& \leq & \omega,
\end{eqnarray*}
where $\omega := |B(0,1)|$ denotes the volume of the unit euclidean ball. We used in the third line the proportionality relation $H(z)= (\det M(z))^{-\frac 1 {mp+d}}M(z)$ and the  fact that $\pi_z$ is homogeneous of degree $m$, which imply 
$$
\det H(z) = \det M(z)^{\frac {mp}{mp+d}} \stext{ and } \|\pi_z \circ H(z)^{-\frac 1 2}\| = (\det M(z))^{\frac m{2(mp+d)}}\|\pi_z \circ M(z)^{-\frac 1 2}\|.
$$

The convergence \iref{eqConvfPi} is uniform over all $z\in [0,1]^d$, and therefore 
\begin{eqnarray*}
\lim_{\lambda \to \infty} \lambda^{mp} e_{\lambda^2 H} (f)_p^p &=& \lim_{\lambda \to \infty} \lambda^{mp} \int_{[0,1]^d} \sqrt{\det (\lambda^2 H(z))}  e_{\lambda^2 H}(f\ssep z)_p^p dz\\
&=& \int_{[0,1]^d} \sqrt{\det  H(z)}  \left( \lim_{\lambda \to \infty} \lambda^{mp+d} e_{\lambda^2 H}(f\ssep z)_p^p\right) dz\\
& \leq &\omega^p \int_{[0,1]^d} \sqrt{\det H(z)}   dz.\\
&=&\omega^p m(H).
\end{eqnarray*}
Since $m(\lambda^2 H) = \lambda^d m(H)$ we therefore obtain injecting \iref{eqMK}
\begin{eqnarray*}
\lim_{\lambda \to \infty} m(\lambda^2 H)^\frac m d e_{\lambda^2 H}(f)_p & \leq & \omega \, m(H)^{\frac m d+ \frac 1 p}\\
& = & 
C\|K(d^m f)\|_{L^\tau([0,1]^d)}+\ve,
\end{eqnarray*}
which concludes the proof.

\subsection{The shape function is equivalent to a continuous function}
\label{subsecKCont}
The shape functions $K$, $L_a$ and $L_g$ are defined in \S\ref{subsecIntroAsympt} by an infimum, which may or may not be attained, over the collection $\SL_d$ of all matrices of unit determinant. 
One can infer from this property that the shape functions are upper-semi continuous, but they are not continuous in general as illustrated in Remark \ref{overfitting} (in Chapter 2) for the shape function $K$ when $m=4$ and $d=2$.

The purpose of this subsection is to establish the following theorem, which states that the shape functions are nevertheless uniformly equivalent to a continuous function on $\H_m$.

\begin{theorem}
\label{thContShape}
There exists a constant $C=C(m,d)$ and three continuous functions $K^c$, $L^c_a$ and $L^c_g$ on $\H_m$ such that for all $\pi \in \H_m$
$$
\begin{array}{ccccc}
C^{-1} K^c(\pi) &\leq& K(\pi) &\leq& C K^c(\pi)\\
C^{-1} L^c_a(\pi) &\leq& L_a(\pi) &\leq& C L^c_a(\pi)\\
C^{-1} L^c_g(\pi) &\leq& L_g(\pi) &\leq& C L^c_g(\pi).
\end{array}
$$
\end{theorem}

Note that in the case of the shape functions $K$ and $L_a$, in dimension $d=2$, this theorem is a consequence of the explicit algebraic equivalents of these shape functions given in Chapters 2 and 3.

The analysis presented below applies to all three shape functions, without restriction on the dimension $d$.
As observed in \S \ref{subsecWellPosed} the shape functions $L_a$ and $L_g$ are uniformly equivalent on $\H_m$ to the functions $L_a^*$ and $L_g^*$ respectively, defined for all $\pi\in \H_m$ by
\begin{eqnarray}
L_a^*(\pi) &:=& \inf_{A\in \SL_d} \sqrt{\|G(\pi)\circ A\|},\\
\label{defLGS}
L_g^*(\pi) &:=& \inf_{A\in \SL_d} \|A^{-1}\| \|\pi\circ A\|,
\end{eqnarray}
where $G(\pi) := |\nabla \pi|^2 \in \H_{2m-2}$. Let us also recall that 
\be
\label{defKCont}
K(\pi) := \inf_{A\in \SL_d} \|\pi\circ A\|.\\
\ee
The fact that $L_a^*$ is uniformly equivalent to a continuous function therefore directly follows from the same property for the shape function $K$ (in degree $2m-2$ instead of $m$), by composition with the continuous functions $G : \H_m \to \H_{2m-2}$ and $\sqrt{\cdot} : \R_+ \to \R_+$.

The next lemma gives a criterion on the lower semi-continuous envelope of a function, which guarantees that it is equivalent to a continuous function. 
\begin{lemma}
\label{lemmaToCont}
Let $(X,d_X)$ be an arbitrary metric space, and let $f : X \to \R_+$. We define the lower semi continuous envelope $\underline f$ of $f$ as follows : for all $x\in X$
$$
\underline f(x) := \lim_{\ve \to 0} \inf_{y\in B(x,\ve)} f(y).
$$
Assume that $f$ is upper semi-continuous and that there exists a constant $C$ such that 
\be
\label{eqffLow}
f(x) < C \underline f(x) \stext{ or } f(x) = 0
\ee
for all $x\in X$. Then there exists a continuous function $f^c \in C^0(X,\R_+)$ such that 
\be
\label{eqffc}
C^{-1} f^c \leq f \leq C f^c.
\ee
\end{lemma}

\proof
For each $n\in \Z$ we define the set 
\be
\label{defEn}
E_n := \{x\in X \sep f(x) \geq C^n\},
\ee
which is closed since $f$ is upper semi-continuous. 
For any $x\in \overline{X\sm E_n}$ there exists a sequence $(x_k)_{k\geq 0}$, $x_k\in X\sm E_n$ (hence $f(x_k) < C$), converging to $x$. Therefore
$$
\underline f(x)\leq \liminf_{k \to \infty} f(x_k)\leq C.
$$
Combining this with \iref{eqffLow} we obtain that $x\notin E_{n+1}$, hence the closed sets $E_{n+1}$ and $\overline{X\sm E_n}$ are disjoint.
We denote by $(r_n)_{n \in \Z}$,  $r_n\in C^0(X, [0,1])$, a family of continuous functions such that
$$
r_{n|E_{n+1}} = 1\stext{ and } r_{n|\overline{X\sm E_n}} = 0.
$$ 
The simplest construction of such functions is the following:
$$
r_n(x) := \frac{d_X(x, X \sm E_n)}{d_X(x, X\sm E_n)+ d_X(x,E_{n+1})},
$$
where 
$d_X(x,E) := \inf \{d_X(x,e)\sep e\in E\}$.
We define a function $r : X \to [-\infty, \infty)$ (observe that $-\infty$ is included) by the sum
$$
r(x) := \sum_{n <0} (r_n(x)-1) + \sum_{n \geq 0} r_n(x).
$$
For each $n\in \Z$ the function $r$ is equal to $n+r_n+r_{n+1}$ on the open set $\interior(E_n\sm E_{n+2})$. 
Furthermore 
$$
E_{n+1} \sm E_{n+2} \subset \{x\in X \sep f(x) < C^{n+2} \text{ and } \underline f(x) > C^n\} \subset \interior(E_n\sm E_{n+2}),
$$
therefore $r$ is continuous on the reunion
$$
\bigcup_{n \in \Z} \interior(E_n\sm E_{n+2}) \supset \bigcup_{n \in \Z} E_{n+1}\sm E_{n+2} = \{x\in X\sep f(x) >0\}.
$$
Consider $x\in X$ such that $f(x) = 0$, hence $r(x) = -\infty$. Since $f$ is upper semi-continuous there exists for each $n\in \Z$ a constant $\delta>0$ such that $f < C^n$ on $B(x,\delta)$. It follows that $r\leq n$ on $B(x,\delta)$, hence that $r$ is continuous at $x$, in the topology of $[-\infty, \infty)$.

The function $r$ is therefore continuous on the whole space $X$. Since $C\geq 1$ the function $f^c : X \to \R_+$ defined by 
$$
f^c(x) := C^{\,r(x)},
$$
is therefore also continuous. We have $f^c(x) = 0$ if and only if $r(x)=-\infty$, which is equivalent to $f(x) = 0$.
Furthermore for each $n\in \Z$ and each $x\in E_n\sm E_{n+1}$ we have $C^n \leq f(x) < C^{n+1}$ and $C^n \leq f^c(x) \leq C^{n+1}$ by construction, which establishes \iref{eqffc} and concludes the proof of this lemma.
\sq

The proofs that the functions $K$ and $L_g^*$ are uniformly equivalent on $\H_m$ to a continuous function are extremely similar. The case of $L_g^*$ is nevertheless slightly more involved due to the term $\|A^{-1}\|$ appearing in \iref{defLGS} which is not present in \iref{defKCont}. We therefore focus our attention on $L_g^*$, and we leave to the reader the details of the adaptation of the proof to the shape function $K$.

We define an auxiliary function $\gL$ on $M_d \times \H_m$ as follows : for all $(B, \pi)\in M_d\times \H_m$
\be
\label{defgL}
\gL(B,\pi) := \inf_{A\in \SL_d} \|A^{-1}B\| \|\pi\circ A\|,
\ee
hence $\gL(\Id, \pi) = L_g(\pi)$ for all $\pi \in \H_m$.

The function $\gL$ is upper semi continuous since it is defined as the infimum of a family of continuous functions. Hence for any converging sequence $(B_n, \pi_n) \to (B,\pi)$ in $M_d \times \H_m$ we have 
$$
\gL(B, \pi) \geq \limsup_{n \to \infty} \gL(B_n, \pi_n).
$$
The next lemma shows that this inequality becomes an equality if the limit is zero.
\begin{lemma}
\label{lemmaVanishGL}
For any converging sequence $(B_n, \pi_n) \to (B,\pi)$ in $M_d \times \H_m$ the following holds:
$$
\text{ if } \lim_{n \to \infty} \gL(B_n, \pi_n)=0 \stext{ then } \gL(B, \pi) = 0.
$$
\end{lemma}

\proof
If $\gL(B_n, \pi_n) \to 0$ then there exists a sequence $(A_n)_{n \geq 0}$, $A_n \in \SL_d$, such that 
\be
\label{eqAnBnPin}
\lim_{n \to \infty} \|A_n^{-1} B_n\| \|\pi_n\circ A_n\| =0.
\ee
For each $n$ we consider the principal value decomposition $A_n = U_n D_n V_n$, where $U_n, V_n \in \cO_d$ and $D_n$ is a diagonal matrix satisfying $\det D_n = 1$, with positive 
diagonal coefficients $t_n = (t_{n,1}, \cdots, t_{n,d})$. 
We have for any $n \geq 0$, since $V_n$ is orthogonal,
\begin{eqnarray*}
\|A_n^{-1} B_n \| \|\pi\circ A_n\| &=& \| V_n^\trans D_n^{-1} U_n^\trans B_n \| \|\pi\circ (U_n D_n V_n)\|\\
&=& \|D_n^{-1} U_n^\trans B_n \| \|\pi\circ (U_n D_n)\|.
\end{eqnarray*}
Since the collection $\cO_d$ of orthogonal matrices is compact, the sequence of orthogonal matrices $U_n$ admits a converging sub-sequence $U_{\vp(n)} \to U$. We define $B' := U^\trans B$, $\pi' := \pi \circ U$, and we observe that since $|\det U| = 1$
\be
\label{eqgLBPip}
\cL(B, \pi) = \cL(B', \pi').
\ee
We also define for each $n$
$$
B'_n := U_{\vp(n)}^\trans B_{\vp(n)}, \quad  \pi'_n := \pi_{\vp(n)} \circ U_{\vp(n)}, \stext{ and } D'_n := D_{\vp(n)},
$$
in such way that 
\be
\label{eqBpPip}
\lim_{n \to \infty} (B'_n,\pi'_n) = (B',\pi') \stext{ and } \lim_{n \to \infty} \|D'^{-1}_n B'_n\| \|\pi'_n\circ D'_n\| = 0.
\ee
We thus recognize our starting point \iref{eqAnBnPin} except that the matrices $A_n$ are replaced with diagonal matrices $D'_n$, of coefficients 
$$
t'_n := (t'_{n,1}, \cdots , t'_{n,d}) =  (t_{\vp(n),1}, \cdots , t_{\vp(n),d}).
$$
In order to avoid notational clutter we do not keep track of the sub-sequence extraction, and until the end of this proof we write $\pi$, $\pi_n$, $B$, $B_n$, $D_n$, $t_n$ and $t_{n,i}$ for the variables $\pi'$, $\pi_n'$, $B'$, $B_n'$, $D'_n$, $t_n'$ and $t_{n,i}'$.

We denote by $I\subset \{1, \cdots, d\}$ the collection of indices $i$ such that the associated line of $B$ is nonzero. There exists a constant $c>0$ such that for 
all $n$ sufficiently large
\be
\label{eqDnB}
\|D_n^{-1} B_n\| \geq c \max_{i\in I} t_{n,i}^{-1}.
\ee
We denote by $\Lambda$ the collection of exponents
$$
\Lambda := \{\alpha\in \Z_+^d \sep |\alpha| \leq m\}
$$
and by $\pi_\alpha$, $\pi_{n,\alpha}$ the coefficients of $\pi$ and $\pi_n$ respectively : 
$$
\pi = \sum_{\alpha \in \Lambda} \pi_\alpha Z^\alpha \stext{ and } \pi_n = \sum_{\alpha\in \Lambda} \pi_{n, \alpha} Z^\alpha.
$$
It follows from \iref{eqBpPip} and \iref{eqDnB} that for any $\alpha\in \Lambda$ we have 
\be
\label{eqTnPin}
\lim_{n \to \infty} \left(\max_{i\in I} t_{n,i}^{-1}\right) \pi_{n, \alpha} t_n^\alpha = 0.
\ee
For each $t = (t_1, \cdots , t_d)\in \R_+^d$ we define 
$$
\Lambda(t) := \left\{\alpha\in \Lambda \sep\left(\max_{i\in I} t_i^{-1}\right) t^{\alpha} \geq 1\right\}.
$$
The set $\Lambda$ has only a finite number of subsets, since it is finite. Hence there exists an extraction $\psi$ such that $\Lambda_* := \Lambda(t_{\psi(n)})$ does not depend on $n$. In view of \iref{eqTnPin} we thus have for all $\alpha\in \Lambda_*$
$$
\pi_\alpha = \lim_{n \to \infty} \pi_{\psi(n), \alpha} = 0. 
$$
Defining $D := D_{\psi(0)}$, and denoting by $t=(t_1, \cdots ,t_d)$, its diagonal coefficients we thus obtain 
$$
\lim_{n \to \infty} \|D^{-n} B\| \|\pi\circ D^n\| = 0. 
$$
Indeed $\|D^{-n} B\|\leq C \max_{i\in I}t_i^{-n}$ and 
$$
\left( \max_{i\in I}t_i^{-n}\right) \pi_\alpha t^{n\alpha} =  \pi_\alpha \left(t^{\alpha} \max_{i\in I}t_i^{-1}\right)^n,
$$
which equals $0$ if $\alpha\in \Lambda_*$, and tends to zero if $\alpha\in \Lambda\sm\Lambda_*$. It follows that $\gL(B,\pi) = 0$ as announced.
\sq

The next lemma compares the function $\gL$ with its lower semi-continuous envelope.
\begin{lemma}
\label{lemmagLLower}
There exists a constant $C$ such that $\gL \leq C \underline \gL$ on $M_d \times \H_m$.
\end{lemma}
\proof
We denote by $\cA$ the following compact subset of $M_d\times \H_m$
$$
\cA := \{(B,\pi) \in M_d \times \H_m \sep \|B\|=\|\pi\| = 1 \text{ and }  \|A^{-1}B\| \|\pi\circ A\|\geq 1 \text{ for all } A\in \SL_d \}.
$$
It follows from the definition \iref{defgL} of $\gL$ that $\gL(B,\pi) = 1$ for all $(B, \pi)\in \cA$. Therefore  $\underline \gL$ does not vanish on $\cA$ according to Lemma \ref{lemmaVanishGL}. Since $\underline \gL$ is lower semi continuous, it attains its minimum on $\cA$, we thus define 
$$
C^{-1} := \inf_\cA \underline \gL.
$$
Let $(B, \pi) \in M_d \times \H_m$, and let $(A_n)_{n \geq 0}$, $A_n \in \SL_d$, be a sequence such that 
\be
\label{eqAngLMin}
\lim_{n \to \infty} \|A_n^{-1} B \| \|\pi \circ A_n\| = \gL(B, \pi) := \inf_{A\in \SL_d} \|A^{-1} B \| \|\pi \circ A\|.
\ee
If $\gL(B, \pi) = 0$, then $\underline \gL(B, \pi)=0$ and there is nothing to prove, we may therefore assume that $\gL(B, \pi)>0$. 
We have for each $n\geq 0$
$$
\underline \gL\left(\frac{A_n^{-1} B}{\|A_n^{-1} B\|},\, \frac{\pi\circ A_n}{\|\pi\circ A_n\|}\right) = \frac{\underline \gL(B,\pi) }{\|A_n^{-1} B\|\|\pi\circ A_n\|}, 
$$
which tends to $\underline \gL(B,\pi)/\gL(B,\pi)$ as $n \to \infty$.
Consider an extraction $\vp$ such that the following quantities converge
$$
\lim_{n \to \infty} \frac{A_{\vp(n)}^{-1} B}{\|A_{\vp(n)}^{-1} B\|} = B_* \stext{ and } \lim_{n \to \infty} \frac{\pi\circ A_{\vp(n)}}{\|\pi\circ A_{\vp(n)}\|} = \pi_*.
$$

We now prove that $(B_*, \pi_*)\in \cA$.
We first remark that $\|B_*\| = \|\pi_*\|=1$ by construction. Assume for contradiction that there exists $A\in \GL_d$ such that $\|A^{-1} B_*\|\|\pi_*\circ A\| = \delta <1$. We obtain 
\begin{eqnarray*}
& &\lim_{n \to \infty}  \|(A_{\vp(n)} A)^{-1} B \| \|\pi \circ (A_{\vp(n)} A)\| \\
&=& \lim_{n \to \infty} \frac {\|A^{-1}(A_{\vp(n)}^{-1} B) \|}{\|A_{\vp(n)}^{-1} B \|} \frac { \|(\pi\circ A_{\vp(n)})\circ A\|}{\|\pi\circ A_{\vp(n)}\|} \|A_{\vp(n)}^{-1} B \| \|\pi\circ A_{\vp(n)}\| \\
&=& \delta \gL(B, \pi)
\end{eqnarray*}
which contradicts the definition \iref{defgL} of $\gL$, hence $(B_*, \pi_*)\in \cA$ as announced.
We thus obtain, since $\underline \gL$ is lower semi continuous, 
$$
\lim_{n \to \infty} \underline \gL\left(\frac{A_n^{-1} B}{\|A_n^{-1} B\|},\, \frac{\pi\circ A_n}{\|\pi\circ A_n\|}\right) \geq  \underline \gL(B_*, \pi_*) \geq C^{-1}. 
$$
If follows that $C\underline \gL(B,\pi) \geq \gL(B, \pi)$, which concludes the proof of this lemma.
\sq

We now conclude the proof that $L_g^*$ is equivalent to a continuous function. 
Applying Lemma \ref{lemmaToCont} to the function $\gL$ and the constant $C+1$, where $C$ is the constant from Lemma \ref{lemmagLLower}, we obtain that that $\gL$ is equivalent to a continuous function on $M_d \times \H_m$. Recalling that $L_g^*(\pi) = \gL(\Id,\pi)$ we obtain that $L_g^*$ is also equivalent to a continuous function, on $\H_m$, which concludes the proof of Theorem \ref{thContShape}.



\part{Hierarchical refinement algorithms}
\label{partGreedy}

\chapter[Adaptive and anisotropic multiresolution analysis]{Adaptive multiresolution analysis based on anisotropic triangulations}
\minitoc
\label{chapCDHM}
\section{Introduction}
\label{intro}

Approximation by piecewise polynomial functions
is a standard procedure which occurs in various applications. 
In some of them such as 
terrain data simplification or image compression,
the function to be approximated might be fully known,
while it might be only partially known or fully
unknown in other applications such as denoising, statistical learning or 
in the finite element discretization of PDE's.

In all these applications, one usually makes the 
distinction between {\it uniform} and {\it adaptive} approximation.
In the uniform case, the domain of interest is decomposed into
a partition where all elements have comparable shape and size,
while these attributes are allowed to vary strongly in the adaptive case.
In the context of adaptive triangulations, another important distinction
is between {\it isotropic} and {\it anisotropic} triangulations.
In the first case the triangles
satisfy a condition which guarantees that 
they do not differ too much from equilateral triangles. This can
either be stated in terms of a minimal value $\theta_0>0$ for
every angle, or by a uniform bound on the
aspect ratio
$$
\rho_T:=\frac {h_T}{r_T}
$$
of each triangle $T$ where $h_T$ and $r_T$ respectively denote
the diameter of $T$ and of its largest inscribed disc.
In the second case, which is in the scope of 
the present chapter, the aspect ratio is allowed
to be arbitrarily large, i.e. long and thin triangles are allowed.
In summary, adaptive and anisotropic triangulations mean that
we do not fix any constraint on the size and shape of the triangles.

Given a function $f$ and a norm $\|\cdot\|_X$ of interest, 
we can formulate the problem
of finding the {\it optimal triangulation} for $f$ in the $X$-norm 
in two related forms:
\begin{itemize}
\item
For a given $N$ find a triangulation $\cT_N$ with $N$ triangles
and a piecewise polynomial function $f_N$ (of some fixed degree) on $\cT_N$
such that $\|f-f_N\|_X$ is minimized.
\item
For a given $\e >0$ find a triangulation $\cT_N$ with minimal
number of triangles $N$ and a piecewise polynomial function $f_N$ 
such that $\|f-f_N\|_X\leq \e$.
\end{itemize}
In this chapter $X$ will be the $L^p$ norm for some arbitrary $1\leq p\leq \infty$.
The exact solution to such problems is usually out of reach both
analytically and algorithmically: even when restricting the search 
of the vertices to a finite grid, the number of possible triangulations
has combinatorial complexity and an exhaustive search is therefore
prohibitive. 

Concrete mesh generation algorithms have been
developed in order to generate in reasonable time triangulations
which are ``close'' to the above described optimal trade-off 
between error and complexity. They are typically governed
by two intuitively desirable features:
\begin{enumerate}
\item
The triangulation should {\it equidistribute} the local approximation error between 
each triangle. This rationale is typically used in local mesh refinement
algorithms for numerical PDE's
\cite{Ve}: a triangle is refined when the local approximation 
error (estimated by an a-posteriori error indicator) is large.
\item
In the case of anisotropic meshes, the local aspect ratio
should in addition be optimally adapted to 
the approximated function $f$. In the case of piecewise linear
approximation, this is achieved by imposing
that the triangles are isotropic with respect to a distorted
metric induced by the Hessian $d^2f$. 
We refer in particular to \cite{BFGLS}
where this task is executed using 
Delaunay mesh generation techniques.
\end{enumerate}

While these last algorithms fastly produce 
anisotropic meshes which are naturally
adapted to the approximated function,
they suffer from two intrinsic limitations:
\begin{enumerate}
\item
They are based on the evaluation of the Hessian $d^2f$,
and therefore do not in principle apply 
to arbitrary functions $f\in L^p(\Omega)$ for $1\leq p\leq \infty$ 
or to noisy data. 
\item
They are non-hierarchical: for $N>M$, the triangulation
$\cT_N$ is not a refinement of $\cT_M$. 
\end{enumerate}

One way to circumvent the first limitation is
to regularize the function $f$, either by projection onto a 
finite element space or by convolution by a mollifier. However
this raises the additional problem of appropriately tuning
the amount of smoothing, in particular depending on the
noise level in $f$.

The need for hierarchical triangulations is critical
in the construction of wavelet bases, which play
an important role in applications
to image and terrain data processing,
in particular data compression \cite{CDDD1}.
In such applications, the multilevel structure
is also of key use for the fast encoding
of the information.
Hierarchy is also useful in the 
design of optimally converging adaptive
methods for PDE's \cite{Dor,MNS,BDD,St}.
However, all these developments are so 
far mostly restricted to isotropic refinement methods.
Let us mention that hierarchical and anisotropic
triangulations have been investigated in \cite{KP},
yet in this work the triangulations are {\it fixed in advance}
and therefore generally not adapted to the approximated function.
\nl
\nl
{\it A natural objective is therefore to design
adaptive algorithmic techniques that combine
hierarchy and anisotropy, and that apply
to any function $f\in L^p(\Omega)$, without any need for regularization.}
\nl
\nl
In this chapter we propose and study a
simple {\it greedy refinement procedure}
that achieves this goal: starting from an initial 
triangulation $\cD_{0}$, the procedure
bisects every triangle from one of its vertices to the mid-point of 
the opposite segment. The choice of the vertex
is typically the one which minimizes the new approximation
error after bisection among the three options.

Surprisingly, it turns out that - in the case of piecewise
linear approximation - this elementary
strategy tends to generate anisotropic triangles 
with optimal aspect ratio. This fact is rigorously
proved in Chapter \ref{chapBisecOpt} which establishes optimal 
error estimates for the approximation
of {\it smooth and convex functions $f\in C^2$},
by adaptive triangulations $\cT_N$ with $N$ triangles.
These triangulations are obtained
by consecutively applying the refinement procedure 
to the triangle of maximal error. The estimates in Chapter \ref{chapBisecOpt}
are of the form
\be
\|f-f_N\|_{L^p}  \leq CN^{-1}\|\sqrt{|{\rm det}(d^2f)|}\|_{L^\tau},\;\; \frac 1 \tau=\frac 1 p + 1,
\label{aniserCDHM}
\ee
and were already established in \cite{CSX,BBLS} for functions which
are not necessarily convex, however based on triangulations which are non-hierarchical and based on the
evaluation of $d^2f$. Note that \iref{aniserCDHM} improves on the estimate
\be
\|f-f_N\|_{L^p}  \leq CN^{-1}\|d^2f\|_{L^\tau},\;\; \frac 1 \tau=\frac 1 p + 1.
\label{isoer}
\ee
which can be established, see \S \ref{secCM2} in the next chapter, for adaptive triangulations with 
isotropic triangles, and which itselfs improves on the
classical estimate
\be
\|f-f_N\|_{L^p}  \leq CN^{-1}\|d^2f\|_{L^p},
\label{unier}
\ee
which is known to hold for uniform triangulations.

The main objective of the present chapter is to introduce
the refinement procedure as well as several approximation
methods based on it, and
to study their convergence for {\it an arbitrary function $f\in L^p$}.
In \S \ref{secCDHM2}, we introduce notation that serves to
describe the refinement procedure and define
the anisotropic hierarchy of triangulations $(\cD_j)_{j\geq 0}$. 
We show how this general framework can be used 
to derive adaptive approximations of $f$ either by 
triangulations based on greedy or optimal trees,
or by wavelet thresholding.
In \S \ref{secCDHM3}, we show that as defined, the approximations
produced by the refinement
procedure may fail to converge for 
certain $f \in L^p$ and show how to modify 
the procedure so that convergence
holds for any arbitrary $f\in L^p$.
We finally present  in \S \ref{secCDHM4} some numerical tests which illustrate 
the optimal mesh adaptation, in the case of piecewise linear
elements, when the refinement procedure is applied 
either to synthetic functions
or to numerical images. 

\section{An adaptive and anisotropic multiresolution framework}
\label{secCDHM2}
\subsection{The refinement procedure}

Our refinement procedure is based on
a local approximation operator
$\cA_T$ 
acting from $L^p(T)$ onto $\P_m$ - the space of polynomials of
total degree less or equal to $m$. Here, the parameters $m\geq 0$ and $1\leq p\leq \infty$
are arbitrary but fixed. For a generic triangle
$T$, we denote by $(a,b,c)$
its edge vectors oriented in clockwise 
or anticlockwise direction so that
$$
a+b+c=0.
$$
We define the
local $L^p$ approximation error 
$$
e_T(f)_p:=\|f-\cA_Tf\|_{L^p(T)}.
$$
The most natural choice for $\cA_T$ is the operator $\cB_T$ of best $L^p(T)$ approximation
which is defined by
$$
\|f-\cB_Tf\|_{L^p(T)}=\min_{\pi\in\sP_m}\|f-\pi\|_{L^p(T)}.
$$
However this operator is non-linear and not easy to compute when $p\neq 2$. 
In practice, one prefers to use an operator 
which is easier to compute, yet nearly optimal in the sense
that
\be
\|f-\cA_Tf\|_{L^p(T)} \leq C \inf_{\pi\in\sP_m}\|f-\pi\|_{L^p(T)},
\label{unileb}
\ee
with $C$ a Lebesgue constant independent of $f$ and $T$. 
Two particularly simple admissible choices of approximation operators are the following:
\begin{enumerate}
\item
$\cA_T=P_T$, the $L^2(T)$-orthogonal projection onto $\P_m$, defined
by $P_Tf\in \P_m$ such that ${\int_T(f-P_T f)\pi}=0$ for all $\pi\in\P_m$. This operator 
has finite Lebesgue constant for all $p$, with $C=1$ when $p=2$
and $C\geq1$ otherwise.
\item
$\cA_T=I_T$, the local interpolation operator which is defined by  $I_Tf\in\P_m$ such that
$I_Tf(\gamma)=f(\gamma)$ for all $\gamma\in \Sigma:=\{ \sum \frac {k_i}m v_i\sep k_i\in \N,\; 
\sum k_i=m\}$
where $\{v_1,v_2,v_3\}$ are the vertices of $T$ (in the case $m=0$ we can take for $\Sigma$
the barycenter of $T$). This operator is only defined on
continuous functions and has Lebesgue constant $C>1$ in the $L^\infty$ norm. 
\end{enumerate}
All our results are simultaneously valid when $\cA_T$ is either $P_T$
or $I_T$ (in the case where $p=\infty$), or any linear operator that fulfills 
the continuity assumption \iref{unileb}.

Given a target function, our refinement procedure defines by induction a hierarchy of nested
triangulations  $(\cD_j)_{j\geq 0}$ with $\#(\cD_j)=2^j \#(\cD_0)$.
The procedure starts from the coarse triangulation $\cD_0$ of $\Omega$,
which is fixed independently of $f$. When 
$\Omega=[0,1]^2$ we may split it into
two symmetric triangles so that $\#(\cD_0)=2$. 
For every $T\in\cD_j$, we split $T$ into 
two sub-triangles of equal area by bisection
from one of its three vertices towards the
mid-point of the opposite edge $e\in \{a,b,c\}$. 
We denote by 
$T_e^1$ and $T_e^2$
the two resulting triangles. The choice of $e\in \{a,b,c\}$
is made according to a {\it refinement rule}
that selects this edge depending on the properties of $f$.
We denote by $\cR$ this refinement rule, which can therefore
be viewed as a mapping
$$
\cR : (f,T) \mapsto e.
$$
We thus obtain two {\it children} of $T$ corresponding to the choice $e$.
$\cD_{j+1}$ is the triangulation consisting of all such pairs 
corresponding to all $T\in \cD_j$.

In this chapter, we consider refinement rules
where the selected edge $e$ minimizes a {\it decision function} $e\mapsto d_T(e,f)$ 
among $\{a,b,c\}$. We refer to such rules
as {\it greedy refinement rules}. A more elaborate
type of refinement rule is also considered in \S \ref{secCDHM3p3}.

The role of the decision function is to drive the generation of 
anisotropic triangles according to the local properties of
$f$, in contrast to simpler procedures
such as {\it newest vertex bisection} (i.e. split $T$ from the most recently created vertex)
which is independent of $f$ and generates triangulations with isotropic
shape constraint. 

Therefore, the choice of $d_T(e,f)$ is critical in
order to obtain triangles with an optimal aspect ratio.
The most natural choice corresponds to the optimal split
\be
d_T(e,f)=e_{T_e^1}(f)_p^p+e_{T_e^2}(f)_p^p,
\ee
i.e. choose the edge that minimizes the resulting $L^p$ error after bisection. 
It is proved in Chapter \ref{chapBisecOpt} in the case of piecewise linear
approximation, that when $f$ is
a $C^2$ function which is strictly convex or concave
the refinement rule based on the decision function 
\be
d_T(e,f)=\|f-I_{T_e^1}f\|_{L^1(T_e^1)}+\|f-I_{T_e^2}f\|_{L^1(T_e^2)}.
\label{optil1CDHM}
\ee
generates triangles which tend to have have an optimal aspect ratio,
locally adapted to the Hessian $d^2f$. This aspect ratio
is independent of the $L^p$ norm in which one wants to minimize the error
between $f$ and its piecewise affine approximation.

\begin{remark}
If the minimizer $e$ is not unique, we may choose it among 
the multiple minimizers either randomly 
or according to some prescribed ordering of the edges (for example
the largest coordinate pair of the opposite vertex in lexicographical order).
\end{remark} 

\begin{remark}
The triangulations $\cD_j$ which are generated by
the greedy procedure are in general non-conforming, i.e.
exhibit hanging nodes. This is not problematic in the present setting since
we consider approximation in the $L^p$ norm which does
not require global continuity of the piecewise polynomial functions.
\end{remark}

The refinement rule $\cR$ defines a multiresolution
framework. For a given $f\in L^p(\Omega)$ and any 
triangle $T$ we denote by
$$
\cC(T):=\{T_1,T_2\},
$$
the {\it children} of $T$ which are 
the two triangles obtained by splitting $T$
based on the prescribed decision function $d_T(e,f)$. We also
say that $T$ is the parent of $T_1$ and $T_2$ and
write
$$
T=\cP(T_1) = \cP(T_2).
$$
Note that 
$$
\cD_j:=\cup_{T\in\cD_{j-1}}\cC(T).
$$
We also define
$$
\cD:=\cup_{j\geq 0}\cD_j,
$$
which has the structure of {\it an infinite binary tree}.
Note that $\cD_j$ depends on $f$ (except for $j=0$)
and on the refinement rule $\cR$,
and thus $\cD$ also depends on $f$ and $\cR$:
$$
\cD_j=\cD_j(f,\cR)\;\; {\rm and}\;\; \cD=\cD(f,\cR).
$$
For notational simplicity, we sometimes omit the dependence
in $f$ and $\cR$ when there is
no possible ambiguity.

\subsection{Adaptive tree-based triangulations}

A first application of the multiresolution framework
is the design of adaptive anisotropic triangulations $\cT_N$
for piecewise polynomial approximation, by a {\it greedy 
tree algorithm}. For any finite sub-tree $\cS\subset \cD$, we denote
by 
$$
\cL(\cS):=\{T\in \cS\;\; {\rm s.t.}\;\; \cC(T)\notin\cS\}
$$
its {\it leaves} which form a partition of $\Omega$. We also denote by
$$
\cI(\cS):=\cS\sm\cL(\cS),
$$
its {\it inner nodes}. Note that {\it any} finite partition 
of $\Omega$ by elements of $\cD$ 
is the set of leaves of a finite sub-tree. One easily checks that 
$$
\#(\cS)=2\#(\cL(\cS))-N_0.
$$
For each $N$, the greedy tree algorithm defines a finite sub-tree $\cS_N$
of $\cD$ which grows from $\cS_{N_0}:=\cD_0=\cT_{N_0}$, by adding
to $\cS_{N-1}$ the two children of the triangle $T^*_{N-1}$
{\it which maximizes the local $L^p$-error
$e_T(f)_p$ over all triangles in $\cT_{N-1}$}.

The adaptive partition $\cT_N$ associated
with the greedy algorithm is defined by
$$
\cT_N:=\cL(\cS_N).
$$
Similarly to $\cD$, the triangulation $\cT_N$ depends on $f$ and on
the refinement rule $\cR$, but also on $p$
and on the choice of the approximation operator $\cA_T$.
We denote by $f_N$ the 
piecewise polynomial approximation to $f$ which is defined as $\cA_T f$ on each $T\in \cT_N$. The global $L^p$ approximation
error is thus given by
$$
\|f-f_N\|_{L^p}=\|(e_T(f)_p)\|_{\ell^p(\cT_N)}.
$$
Stopping criterions for the algorithm
can be defined in various ways:
\begin{itemize}
\item
Number of triangles: stop once a prescribed
$N$ is attained.
\item
Local error: stop once $e_T(f)_p\leq \e$ for all $T\in\cT_N$, for some prescribed $\e>0$.
\item
Global error: stop once $\|f-f_N\|_{L^p}\leq \e$ for some prescribed $\e>0$.
\end{itemize}

\begin{remark}
\label{remDiffRoles}
The role of the
triangle selection based on the largest $e_T(f)_p$ is to
equidistribute the local $L^p$ error, a feature which is desirable
when we want to approximate $f$ in $L^p(\Omega)$
with the smallest number of triangles.
However, it should be well understood that the
refinement rule may still be chosen 
based on a decision function defined
by approximation errors in norms that differ from $L^p$.
In particular, as explained earlier, the decision 
function \iref{optil1CDHM} generates triangles 
which tend to have have an optimal aspect ratio,
locally adapted to the Hessian $d^2f$ when $f$ is strictly convex
or concave, and this aspect ratio
is independent of the $L^p$ norm in which one wants to minimize the error
between $f$ and its piecewise affine approximation.
\end{remark}

The greedy algorithm is one particular
way of deriving an adaptive triangulation
for $f$ within the multiresolution framework
defined by the infinite tree $\cD$. An interesting
alternative is 
to build
adaptive triangulations within $\cD$ which
offer an {\it optimal trade-off between 
error and complexity}. This can be done
when $1\leq p<\infty$ by 
solving the minimization problem
\be
\min_{\cS} \Bigl\{ \sum_{T\in\cL(\cS)} e_T(f)_p^p + \lambda \#(\cS) \Bigl\}
\label{mincart}
\ee
among all finite trees, for some fixed $\lambda>0$.
In this approach, we do not directly control the
number of triangles which depends on the penalty
parameter $\lambda$. However, it is immediate
to see that if $N=N(\lambda)$ is the cardinality of 
$\cT^*_N=\cL(\cS^*)$ where $\cS^*$ is the 
minimizing tree, then $\cT^*_N$
minimizes the $L^p$ approximation error
$$
\cT^*_N:=\underset {\#(\cT)\leq N}\argmin\sum_{T\in\cT}e_T(f)_p^p,
$$
where the minimum is taken among all partitions $\cT$ of $\Omega$ 
within $\cD$ of cardinality less than or equal to $N$.

Due to the additive structure of the error term, the
minimization problem \iref{mincart} can be performed
in fast computational time using an optimal pruning algorithm
of CART type, see \cite{BFOS,Do2}.
In the case $p=\infty$ the associated minimization problem
\be
\min_{\cS} \Bigl\{ \sup_{T\in\cL(\cS)} e_T(f)_\infty + \lambda \#(\cS)\Bigl\},
\label{mincartinf}
\ee
can also be solved by a similar fast algorithm.
It is obvious that this method improves over
the greedy tree algorithm: 
if $N$ is the cardinality of the
triangulation resulting from the minimization in \iref{mincart}
and $f_N^*$ the corresponding piecewise polynomial 
approximation of $f$ associated with this triangulation,
we have
$$
\|f-f_N^*\|_{L^p}\leq \|f-f_N\|_{L^p},
$$
where $f_N$ is built by the greedy tree algorithm.

\subsection{Anisotropic wavelets}

The multiresolution framework allows us to 
introduce the piecewise polynomial multiresolution spaces
$$
V_j=V_j(f,\cR):=\{g \;\; {\rm s.t.}\;\; g_{|T}\in \P_m,\; T\in \cD_j\},
$$
which depend on $f$ and on the refinement rule $\cR$. These spaces
are nested and we denote  by
$$
V=V(f,\cR)=\cup_{j\geq 0}V_j(f,\cR),
$$
their union. For notational simplicity, we sometimes omit the dependence
in $f$ and $\cR$ when there is
no possible ambiguity.

The $V_j$ spaces may be used
to construct wavelet bases, following the approach
introduced in \cite{Alp}
and that we describe in our present setting.

The space $V_j$ is equipped with an orthonormal
{\it scaling function basis}:
$$
\vp_T^{i},\;\; \; i=1,\cdots,\frac 1 2(m+1)(m+2),\;\; T\in\cD_j, 
$$
where the $\vp_T^i$ 
for $i=1,\cdots,\frac 1 2(m+1)(m+2)$
are supported in $T$ and constitute an orthonormal basis
of $\P_m$ in the sense of $L^2(T)$ for each $T\in\cD$. There
are several possible choices for such a basis. In the particular case where
$m=1$, a simple one is to take
for $T$ with vertices $(v_1,v_2,v_3)$,
$$
\vp_T^i(v_i)=|T|^{-1/2} \sqrt 3\;\;{\rm and}\;\; \vp_T^i(v_j)=-|T|^{-1/2} \sqrt 3,\; j\neq i.
$$
We denote by $P_j$ the orthogonal projection onto $V_j$:
$$
P_jg:=\sum_{T\in\cD_j} \sum_{i}\<g,\vp_T^i\>\vp_T^i.
$$
We next introduce for each $T\in\cD_j$
a set of {\it wavelets}
$$
\psi_T^{i},\;\;\; i=1,\cdots,\frac 1 2(m+1)(m+2),
$$
which constitutes an orthonormal basis
of the orthogonal complement of  $\P_m(T)$ into $\P_m(T')\oplus\P_m(T'')$ 
with $\{T',T''\}$ the children of $T$. 
In the particular case where
$m=1$, a simple choice for
such a basis is as follows: if $(v_1,v_2,v_3)$ and $(w_1,w_2,w_3)$
denote the vertices of $T'$ and $T''$, with the
convention that $v_1=w_1$ and $v_2=w_2$ denote
the common vertices, the second one being the midpoint of the segment $(v_3,w_3)$
(i.e. $T$ has vertices $(v_3,w_3,v_1)$), then
$$
\begin{array}{ll}
& \psi_T^1:=\frac{\vp^3_{T'}-\vp^3_{T''}}{\sqrt 2},\\
& \psi_T^2:=\frac{\vp^1_{T'}-\vp^2_{T'}-\vp^1_{T''}+\vp^2_{T''}} 2
,\\
&\psi_T^3:=\frac {\vp^1_{T'}-\vp^3_{T'}+\vp^1_{T''}-\vp^3_{T''}} 2.
\end{array}
$$
where $\vp^i_{T'}$ and $\vp^i_{T''}$ are the above defined scaling functions. 

The family 
$$
\psi_T^{i},\;\;\; i=1,\cdots,\frac 1 2(m+1)(m+2),\;\; T\in\cD_j
$$
constitutes an orthonormal basis of $W_j$, the $L^2$-orthogonal complement of
$V^j$ in $V^{j+1}$.  A multiscale orthonormal basis of $V_J$ is given by
$$
\{\vp_T\}_{T\in\cD_0}\cup \{\psi_T^i\}_{\substack{i=1,\cdots,\frac 1 2(m+1)(m+2),\\T\in\cD_j, \, j=0,\cdots,J-1}}.
$$
Letting $J$ go to $+\infty$ we thus obtain that 
$$
\{\vp_T\}_{T\in\cD_0}\cup \{\psi_T^i\}_{\substack{i=1,\cdots,\frac 1 2(m+1)(m+2),\\T\in\cD_j, \, j\geq 0}}
$$
is an orthonormal basis of the space
$$
V(f,\cR)_2:=\overline{V(f,\cR)}^{L^2(\Omega)}=\overline{\cup_{j\geq 0} V_j(f,\cR)}^{L^2(\Omega)}.
$$
For the sake of notational simplicity, we rewrite this basis as 
$$
(\psi_\lambda)_{\lambda\in\Lambda},
$$
Note that $V(f,\cR)$ is not necessarily dense
in $L^2(\Omega)$ and so $V(f,\cR)_2$ is not always equal to $L^2(\Omega)$.
Therefore, the expansion of an arbitrary function $g\in L^2(\Omega)$
in the above wavelet basis does not always converge towards $g$
in $L^2(\Omega)$. The same remark holds for the $L^p$ convergence
of the wavelet expansion of an arbitrary 
function $g\in L^p(\Omega)$ (or $\cC(\Omega)$
in the case $p=\infty$): $L^p$-convergence holds when the space
$$
V(f,\cR)_p:=\overline{V(f,\cR)}^{L^p(\Omega)}
$$
coincides with $L^p(\Omega)$ (or contains $\cC(\Omega)$ in the case $p=\infty$),
since we have
$$
\Big\|f-\sum_{|\lambda|<j}d_\lambda\psi_\lambda\Big\|_{L^p}=\|f-P_jf\|_{L^p} \leq C
\inf_{g\in V_j}\|f-g\|_{L^p},
$$
with $C$ the Lebesgue constant in \iref{unileb} for the orthogonal projector.

A sufficient condition for such a property to hold is obviously that the
size of all triangles goes to $0$ as the level $j$ increases, i.e.
$$
\lim_{j\to +\infty}\sup_{T\in\cD_j} {\rm diam}(T)=0.
$$
However, this condition might not hold for the
hierarchy $(\cD_j)_{j\geq 0}$ produced by the 
refinement procedure.
On the other hand, the multiresolution 
approximation being intrinsically adapted to $f$,
a more reasonable requirement is that the
expansion of $f$ converges towards $f$ in $L^p(\Omega)$
when $f\in L^p(\Omega)$ (or $\cC(\Omega)$
in the case $p=\infty$). This is equivalent to the property 
$$
f\in V(f,\cR)_p.
$$
We may then define an adaptive approximation of $f$
by thresholding its coefficients
at some level $\e>0$:
$$
f_\e:=\sum_{|f_\lambda|\geq \e} f_\lambda\psi_\lambda,
$$
where $f_\lambda:=\<f,\psi_\lambda\>$. 
When measuring the error in the $L^p$ norm, a more natural choice is
to perform thresholding on the component of the
expansion measured in this norm, defining therefore
$$
f_\e:=\sum_{\|f_\lambda\psi_\lambda\|_{L^p}\geq \e} f_\lambda\psi_\lambda.
$$
We shall next see that the condition $f\in V(f,\cR)_p$ also ensures
the convergence of the tree-based adaptive approximations $f_N$ and $f_N^*$ towards $f$
in $L^p(\Omega)$. We shall also see that this condition may not 
hold for certain functions $f$, but that this
difficulty can be circumvented by a modification
of the refinement procedure.

\section{Convergence analysis}
\label{secCDHM3}
\subsection{A convergence criterion}

The following result relates the convergence towards $f$ of
its approximations by projection onto the spaces $V_j(f,\cR)$, greedy and optimal tree algorithms, 
and wavelet thresholding. This result is valid for {\it any} refinement rule $\cR$.

\begin{theorem}
Let $\cR: (f,T)\mapsto e$ be an arbitrary refinement rule and let
$f\in L^p(\Omega)$. The following statements are equivalent:

{\rm (i)} $f\in V(f,\cR)_p$.

{\rm (ii)} The greedy tree approximation converges: $\lim_{N\to +\infty}\|f-f_N\|_{L^p}=0$.

{\rm (iii)} The optimal tree approximation converges:
$\lim_{N\to +\infty}\|f-f^*_N\|_{L^p}=0$.
\nl
In the case $p=2$, they are also equivalent to:

{\rm (iv)} The thresholding approximation converges:
$\lim_{\e\to 0}\|f-f_\e\|_{L^2}=0$.
\end{theorem}
\noindent
{\bf Proof:}  Clearly, (ii) implies (iii) since $\|f-f_N^*\|_{L^p}\leq \|f-f_N\|_{L^p}$.
Since the triangulation $\cD_N$ is a refinement of $\cT^*_N$,
we also find that  (iii) implies (i) as $\inf_{g\in V_N}\|f-g\|_{L^p}\leq \|f-f^*_N\|_{L^p}$.

We next show that (i) implies (ii). We first note that 
a consequence of (i) is that 
$$
\lim_{j\to+\infty}\sup_{T\in\cD_j} e_T(f)_p =0.
$$
It follows that for any $\eta>0$, there exists
$N(\eta)$ such that for $N>N(\eta)$, all triangles $T\in\cT_N$
satisfy
$$
e_T(f)_p\leq \eta.
$$
On the other hand, (i) means that for all $\e>0$, there exists
$J=J(\e)$ such that
$$
\inf_{g\in V_J}\|f-g\|_{L^p}\leq \e.
$$
For $N>N(\eta)$, we now split $\cT_N$ into $\cT_N^+\cup\cT_N^-$ where
$$
\cT_N^+:=\cT_N\cap (\cup_{j\geq J} \cD_j)\;\;{\rm and}
\;\; \cT_N^-:=\cT_N\cap (\cup_{j< J} \cD_j).
$$
We then estimate the error of the greedy algorithm by
\begin{eqnarray*}
\|f-f_N\|_{L^p}^p & = & \sum_{T\in\cT_N^+}e_T(f)_p^p+\sum_{T\in\cT_N^-}e_T(f)_p^p\\
& \leq & C^p\sum_{T\in\cT_N^+}\inf_{\pi\in \sP_m}\|f-\pi\|_{L^p(T)}^p+\eta^p \#(\cT_N^-) \\
& \leq & C^p\inf_{g\in V_J}\|f-g\|_{L^p}^p+\eta^p \#(\cT_N^-) \\
&\leq & C^p\e^p+2^JN_0\eta^p,
\end{eqnarray*}
%
where $C$ is the stability constant of \iref{unileb}.
This implies (ii) since for any $\delta>0$,
we can first choose $\e>0$ such that $C^p\e^p<\delta/2$, and then
choose $\eta>0$ such that $2^{J(\e)}N_0\eta^p<\delta/2$.
When $p=\infty$ the estimate is modified into
$$
\|f-f_N\|_{L^\infty}\leq \max\{C\e,\eta\},
$$
which also implies (ii) by a similar reasoning.

We finally prove the equivalence
between (i) and (iv) when $p=2$. Property (i) is  equivalent to
the $L^2$ convergence of the orthogonal projection 
$P_{j}f$ to $f$ as $j\to +\infty$, or 
equivalently of the partial sum 
$$
\sum_{|\lambda|< j} f_\lambda \psi_\lambda
$$
where $|\lambda|$ stands for the scale level of the wavelet $\psi_\lambda$.
Since $(\psi_\lambda)_{\lambda\in\Lambda}$ is an orthonormal basis of $V(f,\cR)_2$,
the summability and limit of $\sum_{\lambda\in\Lambda} f_\lambda \psi_\lambda$
do not depend on the order of the terms. Therefore (i) is equivalent
to the convergence of $f_\e$ to $f$.\hfill $\diamond$

\begin{remark}
The equivalence between statements (i) and (iv) can be extended to $1<p<\infty$
by showing that $(\psi_\lambda)_{\lambda\in\Lambda}$
is an $L^p$-unconditional system. Recall that this property means
that there exists an absolute constant $C>0$ such that
for any finitely supported sequences $(c_\lambda)$ and $(d_\lambda)$ such
$|c_\lambda|\leq |d_\lambda|$ for all $\lambda$, one has
$$
\| \sum c_\lambda\psi_\lambda\|_{L^p}\leq C\|\sum d_\lambda \psi_\lambda\|_{L^p}.
$$
A consequence of this property is that if $f\in L^p$
can be expressed as the $L^p$ limit 
$$
f=\lim_{j\to +\infty} \sum_ {|\lambda|<j} f_\lambda\psi_\lambda,
$$
then any rearrangement of the series $\sum f_\lambda \psi_\lambda$
converges towards $f$ in $L^p$. This easily implies the 
equivalence between (i) and (iv). The fact that 
$(\psi_\lambda)_{\lambda\in\Lambda}$ is an unconditional system
is well known for Haar systems \cite{KS} which correspond to the case $m=0$,
and can be extended to $m>0$ in a straightforward manner.
\end{remark}

\subsection{A case of non-convergence}
\label{secCDHM3p2}
We now show that if we use a greedy refinement
rule $\cR$ based on a decision function either based
on the interpolation or $L^2$ projection error
after bisection, there exists functions $f\in L^p(\Omega)$ such that
$f\notin V(f,\cR)_p$. Without loss of generality, it is enough to construct $f$
on the reference triangle $T_{\rm ref}$ of vertices $\{(0,0),(1,0),(1,1)\}$, since our construction
can be adapted to any triangle by an affine change of variables.

Consider first a decision function defined
from the interpolation error after bisection, such as \iref{optil1CDHM}, or
more generally 
$$
d_T(f,e):=\|f-I_{T_e^1}f\|_{L^p(T_e^1)}^p+\|f-I_{T_e^2}f\|_{L^p(T_e^2)}^p.
$$
Let $f$ be a continuous function which is not identically $0$
on $T_{\rm ref}$ and which vanishes
at all points $(x,y)$ such that $x=\frac k{2m}$ for $k=0,1,\cdots,2m$,
where $m$ is the degree of polynomial approximation.
For such an $f$, it is easy to see that
$I_Tf=0$ and that $I_{T'}f=0$ for all
subtriangles $T'$ obtained by one
bisection of $T_{\rm ref}$. This shows that there is 
no preferred bisection. Assuming
that we bisect from the vertex $(0,0)$ to the opposite mid-point $(1,\frac 1 2)$, we 
find that a similar situation occurs when splitting the two subtriangles.
Iterating this observation, we see that an admissible choice of bisections
leads after $j$ steps to a triangulation $\cD_j$ of $T_{\rm ref}$ consisting of the triangles $T_{j,k}$ with vertices
$\{(0,0),(1,2^{-j}k),(1,2^{-j}(k+1))\}$ with $k=0,\cdots,2^j-1$. On each of these triangles
$f$ is interpolated by the null function and therefore
by \iref{unileb} the best $L^\infty$ approximation in $V_j$ does not converge
to $f$ as $j\to +\infty$, i.e.  $f\notin V(f,\cR)_\infty$. It can also easily be
checked that $f\notin V(f,\cR)_p$.

Similar counter-examples can be constructed
when the decision function is defined
from the $L^2$ projection error, and has the form
$$
d_T(f,e):=\|f-P_{T_e^1}f\|_{L^p(T_e^1)}^p+\|f-P_{T_e^2}f\|_{L^p(T_e^2)}^p.
$$
Here, we describe such a construction in the case $m=1$. We define
$f$ on $\cR$ as a function of the first variable given by
$$
f(x,y)=u(x),\;\;{\rm if}\;\;x\in \left[0,\frac 1 2\right], \;\; f(x,y)=u\left(x-\frac 1 2\right),\;\;{\rm if}\;\;x\in \left(\frac 1 2,1\right],
$$
where $u$ is a non-trivial function in $L^2([0,\frac 1 2])$ such that
$$
\int_{0}^{\frac 1 2} u(x)dx=\int_{0}^{\frac 1 2} xu(x)dx=\int_{0}^{\frac 1 2} x^2u(x)dx=0.
$$
A possible choice is $u(x)=L_3(4x-1)=160x^3-120x^2+24x-1$ where $L_3$ is the Legendre polynomial
of degree $3$ defined on $[-1,1]$. With this choice, we have the following result.

\begin{lemma}
Let $T$ be any triangle such that its
vertices have $x$ coordinates either 
$(0,\frac 1 2,1)$ or $(0,1,1)$ or $(\frac 1 2,1,1)$.
Then $f$ is orthogonal to $\P_1$ in $L^2(T)$.
\end{lemma}
\noindent
{\bf Proof:} Define
$$
T_0:=\left\{(x,y)\in T,\; x\in \left[0,\frac 1 2\right]\right\}\;\;{\rm and}\;\; T_1:=\left\{(x,y)\in T,\; x\in \left(\frac 1 2,1\right]\right\}.
$$
Then, with $v(x,y)$ being either the function $1$ or $x$ or $y$, we have
$$
\begin{array}{ll}
\int_T f(x,y)v(x,y) dxdy&=\int_{T_0}u(x)v(x,y)dxdy +\int_{T_1}u(x)v(x,y)dxdy \\
& =\int_{0}^{\frac 1 2}u(x)q_0(x)dx +\int_{0}^{\frac 1 2}u(x)q_1(x)dx,
\end{array}
$$
with
$$
q_0(x):=\int_{T_{0,x}} v(x,y)dy,\;\; q_1(x):=\int_{T_{1,x}}v(x,y)dy,
$$
where $T_{i,x}=\{y\; :\; (x,y)\in T_i\}$ for $i=0,1$.
The functions $q_0$ and $q_1$ are polynomials of degree at most $2$
and we thus obtain from the properties of $u$ that  $\int_T fv=0$. \hfill $\diamond$
\nl
\nl
The above lemma shows that for any of the three possible choices
of bisection of $T_{\rm ref}$ based on the $L^2$ decision function, 
the error is left unchanged since the projection
of $f$ on all possible sub-triangle is $0$. There is 
therefore no preferred choice,
and assuming that we bisect from the vertex $(0,0)$ to the opposite mid-point $(1,\frac 1 2)$, then
we see that a similar situation occurs when splitting the two subtriangles.
The rest of the arguments showing that $f\notin V(f,\cR)_p$
are the same as in the previous counter-example.

The above two examples of non-convergence reflect the fact 
that when $f$ has some {\it oscillations}, the refinement procedure
cannot determine the most appropriate bisection. In order
to circumvent this difficulty one needs to modify
the refinement rule.

\subsection{A modified refinement rule}
\label{secCDHM3p3}
Our modification consists of bisecting from the most recently generated
vertex of $T$, in case the local error is not
reduced enough by all three bisections.  More precisely, 
we modify the choice of the bisection of any $T$ as follows:

Let $e$ be the edge which minimizes the decision function $d_T(e,f)$.
If 
$$
\left(e_{T_e^1}(f)_p^p+e_{T_e^2}(f)_p^p\right)^{1/p}\leq \theta e_{T}(f),
$$
we bisect $T$ towards the edge $e$ (greedy bisection).  Otherwise, we bisect $T$ 
from its most recently generated vertex (newest vertex bisection).
Here $\theta$ is a fixed number in $(0,1)$.
In the case $p=\infty$ we use the condition
$$
\max\{e_{T_e^1}(f),e_{T_e^2}(f)\} \leq \theta e_{T}(f).
$$
This new refinement rule benefits from the mesh size reduction properties
of newest vertex bisection. Indeed, a property illustrated in Figure \ref{fig1CDHM} 
is that a sequence $\{BNN\}$ of one arbitrary bisection ($B$) followed 
by two newest vertex bisections ($N$) produces triangles with diameters
bounded by half the diameter of the initial triangle. A more general
property - which proof is elementary yet tedious - is the following: 
a sequence of the type $\{BNB \cdots BN\}$ of 
length $k+2$ with a newest vertex bisection at iteration $2$ and $k+2$ 
produces triangles with diameter 
bounded by  $(1-2^{-k})$ times the diameter of the initial triangle,
the worst case being illustrated in Figure \ref{fig2CDHM} with $k=3$.

\begin{figure}
	\centering
		\includegraphics[width=9cm]{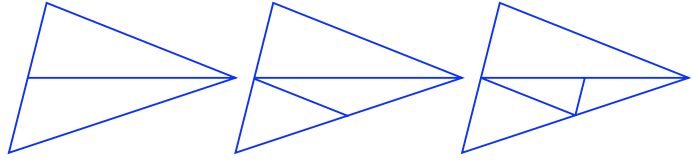} 
	\caption{\label{fig1CDHM}Diameter reduction by a $\{BNN\}$ sequence}
\end{figure}

\begin{figure}
	\centering
		\includegraphics[width=12cm]{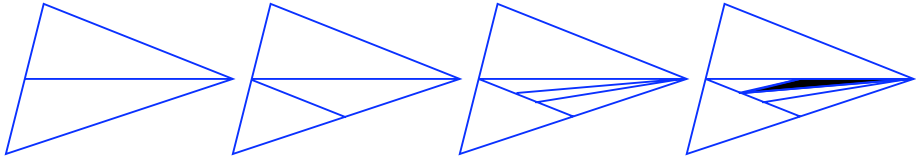} 
	\caption{\label{fig2CDHM}Diameter reduction by a $\{BNBBN\}$ sequence (the dark
triangle has diameter at most $7/8$ times the initial diameter)}
\end{figure}

Our next result shows that the modified algorithm now converges
for any $f\in L^p$.

\begin{theorem}
With $\cR$ defined as the modified bisection rule, we have
$$
f\in V(f,\cR)_p,
$$
for any $f\in L^p(\Omega)$ (or $\cC(\Omega)$ when $p=\infty$).
\end{theorem}
\noindent
{\bf Proof:} We first give the proof when $p<\infty$.
For each triangle $T\in \cD_j$ with $j\geq 1$, we introduce 
the two quantities 
$$
\alpha(T):=\frac {e_T(f)^p_p}{e_T(f)^p_p+e_{T'}(f)^p_p}\stext{ and }
\;\; \beta(T):=\frac {e_T(f)^p_p+e_{T'}(f)^p_p}{e_{\cP(T)}(f)^p_p},
$$
where $T'$ is the ``brother'' of $T$, i.e. $\cC(\cP(T))=\{T,T'\}$.
When a greedy bisection occurs in the split of $\cP(T)$, we have
$\beta(T)\leq \theta^p$. When a newest vertex bisection occurs, we have
$$
\beta(T)\leq  C^p \frac {\inf_{\pi\in\sP_m}\|f-\pi\|_{L^p(T)}^p+\inf_{\pi\in\sP_m}\|f-\pi\|_{L^p(T')}^p}
{\inf_{\pi\in\sP_m}\|f-\pi\|_{L^p(\cP(T))}^p}\leq C^p,
$$
where $C$ is the constant of \iref{unileb}.

We now consider a given level index $j>0$ 
and the triangulation $\cD_j$. For each $T\in \cD_j$,
we consider the chain of nested triangles
$(T_n)_{n=0}^j$ with $T_j=T$ and $T_{n-1}=\cP(T_n)$,
$n=j,j-1,\cdots,1$. 
We define
$$
\o\alpha(T)=\prod_{n=1}^j\alpha(T_n)
\quad \text{ and } \quad
\o\beta(T)=\prod_{n=1}^j\beta(T_n).
$$
It is easy to see that 
$$
\o\alpha(T)\o\beta(T)=\frac {e_T(f)^p_p}{e_{T_0}(f)^p_p}
$$
so that 
$$
e_T(f)^p_p\leq C_0\o\alpha(T)\o\beta(T),
$$
with $C_0:=\max_{T_0\in\cD_0}e_{T_0}(f)^p_p$.
It is also easy to check by induction on $j$ that
$$
\sum_{T\in\cD_j}\o\alpha(T)=\sum_{T\in\cD_{j-1}}\o\alpha(T)=\cdots=\sum_{T\in\cD_{1}}\o\alpha(T)=\#\cD_0.
$$
We denote by $f_j$ the approximation to $f$ in $V_j$ defined
by $f_j=\cA_T f$ on all $T\in\cT_j$ so that
$$
\|f-f_j\|_{L^p}^p=\sum_{T\in\cD_j}e_T(f)^p_p.
$$
In order to prove that $f_j$ converges to $f$
in $L^p$, it is sufficient to show that the sequence
$$
\e_j:=\max_{T\in\cD_j}\min\{\o\beta(T),{\rm diam}(T)\},
$$
tends to $0$ as $j$ grows. Indeed, if this holds, we
split $\cD_j$ into two sets $\cD_{j}^+$ and $\cD_j^-$ over
which $\o\beta(T)\leq \e_j$ and ${\rm diam}(T)\leq \e_j$
respectively. We can then write
\begin{eqnarray*}
\|f-f_j\|_{L^p}^p  & = & \sum_{T\in\cD_j^+}e_T(f)^p_p+\sum_{T\in\cD_j^-}e_T(f)^p_p \\
& \leq & C_0\e_j\sum_{T\in\cD_j^+}\o\alpha(T)+C^p\sum_{T\in\cD_j^-}\inf_{\pi\in\sP_m}\|f-\pi\|_{L^p(T)}^p \\
& \leq & C_0\#\cD_0\e_j+C^p\sum_{T\in\cD_j^-}\inf_{\pi\in\sP_m}\|f-\pi\|_{L^p(T)}^p,
\end{eqnarray*}
where $C$ is the constant of \iref{unileb}. Clearly the first term tends to $0$ and
so does the second term by standard properties of $L^p$
spaces since the diameter of the triangles in $\cD_j^-$ goes to $0$.
It thus remains to prove that
$$
\lim_{j\to +\infty}\e_j=0.
$$
Again we consider the chain $(T_n)_{n=0}^j$ which terminates at $T$,
and we associate to it a chain $(q_n)_{n=0}^{j-1}$ where
$q_n=1$ or $2$ if bisection of $T_n$ is greedy or newest vertex respectively.
If $r$ is the total number of $2$ in the chain $(q_n)$ we have
$$
\o\beta(T)\leq C^{pr}\theta^{p(j-r)},
$$
with $C$ the constant in \iref{unileb}.
Let a $k>0$ be a fixed number, large enough such that
$$
C \theta^{k-1} \leq 1.
$$
We thus have
\be
\o\beta(T)\leq (C^p \theta^{p(k-1)})^r\theta^{p(j-rk)}\leq \theta^{p(j-rk)}.
\label{controlbeta}
\ee
We now denote by $l$ the maximal number of disjoint 
sub-chains of the type 
$$
(\nu_1,2,\nu_2,\nu_3,\cdots,\nu_{q},2)
$$ 
with $\nu_j\in\{1,2\}$ and of length $q+2 \leq 2k+3$ which can be extracted
from $(q_n)_{n=0}^{j-1}$. From the remarks on the diameter reduction
properties of newest vertex bisection, we see that
$$
{\rm diam}(T) \leq B(1-2^{-2k})^l,
$$
with $B:=\max_{T_0\in\cD_0}{\rm diam}(T_0)$ a fixed constant. On the other
hand, it is not difficult to check that
\be
r\leq 3l+3+\frac{j-r}{2k}.
\label{lrjk}
\ee
Indeed let $\alpha_0$ be the total number of $1$
in the sequence $(q_n)$ which are not preceeded by a $2$,
and let $\alpha_i$ be the size of the series
of $1$ following the $i$-th occurence of $2$ in $(q_n)$ for
$i=1,\cdots,r$. Note that some $\alpha_i$ might be $0$.
Clearly we have
$$
j = \alpha_0+ \alpha_1+\cdots+\alpha_r+r.
$$
From the above equality, the number of $i$ such that $\alpha_i>2k$ is 
less than $\frac {j-r}{2k}$ and therefore there
is at least $m\geq r-\frac {j-r}{2k}$ indices
$\{i_0,\cdots,i_{m-1}\}$ such that $\alpha_i\leq 2k$.
Denoting by $\beta_i$ the position of the
$i$-th occurence of $2$ (so that $\beta_{i+1}=\beta_i+\alpha_i+1$),
we now consider the 
disjoint sequences of indices
$$
S_t=\{\beta_{i_{3t}},\cdots,\beta_{i_{3t+2}}\}, \;\; t=0,1,\cdots
$$
There is at least $\frac m 3-1$ such sequences within $\{1,\cdots,j\}$
and by construction each of them contains 
a sequence of the type $(\nu_1,2,\nu_2,\nu_3,\cdots,\nu_{q},2)$ 
with $\nu_j\in\{1,2\}$ and of length $q+2 \leq 2k+3$. Therefore
the maximal number of disjoint sequences of such type satisfies
$$
l\geq \frac 1 3(r-\frac {j-r}{2k})-1,
$$
which is equivalent to \iref{lrjk}. Therefore according to \iref{controlbeta}
$$
\o\beta(T)\leq \theta^{p(j-rk)} \leq \theta^{p(\frac j 2-3(l+1)k+\frac r 2)}
$$
If $3(l+1)k\leq \frac j 4$, we have 
$$
\o\beta(T)\leq \theta^{\frac {pj} 4}.
$$
On the other hand, if $3(l+1)k\geq \frac j 4$, we have
$$
{\rm diam}(T) \leq B(1-2^{-2k})^{\frac j{12 k}-1}.
$$
We therefore conclude that $\e_j$ goes to $0$ as $j$ grows, which proves the result
for $p<\infty$. 
\nl
\nl
We briefly sketch the proof for $p=\infty$, which is simpler. We now define $\beta(T)$ as
$$
\beta(T):=\frac {e_T(f)_\infty}{e_{\cP(T)}(f)_\infty},
$$
so that $\beta(T)\leq \theta$ if a greedy bisection occurs in the split of $\cP(T)$.
With the same definition of $\o\beta(T)$ we now have
$$
e_T(f)_\infty \leq C_0\o\beta(T),
$$
where $C_0:=\max_{T_0\in\cD_0} e_{T_0}(f)_\infty$. With the same definition of 
$\e_j$ and splitting
of $\cD_j$, we now reach
$$
\|f-f_j\|_{L^\infty}\leq \max\left \{C_0\e_j, C\max_{T\in\cD_j^-}\inf_{\pi\in\sP_m}\|f-\pi\|_{L^\infty(T)}\right \},
$$
which again tends to $0$ if $\e_j$ tends to $0$ and $f$ is continuous.
The proof that $\e_j$ tends to $0$ as $j$ grows is then similar to the
case $p<\infty$.
\hfill $\diamond$

\begin{remark}
The choice of the parameter $\theta<1$ deserves some attention: if it is chosen
too small, then most bisections are of type $N$ and we end up with an isotropic
triangulation. In the case $m=1$,
a proper choice can be found by observing that when 
$f$ has $C^2$ smoothness, it can be locally approximated
by a quadratic
polynomial $q\in \P_2$ with $e_T(f)_p\approx e_T(q)_p$
when the triangle $T$ is small enough.
For such quadratic functions $q\in\P_2$, one can
explicitely study the minimal error reduction which is always ensured by the
greedy refinement rule defined a given decision function.
In the particular case $p=2$ and with the choice $\cA_T=P_T$
which is considered in the numerical experiments of \S \ref{secCDHM4},
explicit formulas for the error $\|q-P_Tq\|_{L^2(T)}$ 
can be obtained by formal computing and can be used to prove a guaranteed
error reduction by a factor $\theta^*=\frac 3 5$.
It is therefore natural to choose $\theta$ such that $\theta^*<\theta <1$
(for example $\theta=\frac 2 3$)
which ensures that bisections of type $N$ only occur in the early steps of the
algorithm, when the function still exhibits too many oscillations on 
some triangles.
\end{remark}

\section{Numerical illustrations}
\label{secCDHM4}
The following numerical experiments were conducted with
piecewise linear approximation in the $L^2$ setting:
we use the $L^2$-based decision function
$$
d_T(f,e):=\|f-P_{T_e^1}f\|_{L^2(T_e^1)}^2+\|f-P_{T_e^2}f\|_{L^2(T_e^2)}^2,
$$
and we take for $\cA_T$ the $L^2(T)$-orthogonal projection $P_T$
onto $\P_1$. In these experiments, the function $f$
is either a quadratic polynomial or a function with a simple analytic expression
which allows us to compute the quantities $e_T(f)_2$ and $d_{T}(e,f)$
without any quadrature error,  or a numerical image in which case
the computation of these quantities is discretized on the pixel grid.

\subsection{Quadratic functions}

Our first goal is to illustrate numerically the optimal
adaptation properties of the refinement procedure 
in terms of triangle shape.
For this purpose, we take $f=\bq$ a quadratic form
i.e. an homogeneous polynomial of degree $2$. In this case,
all triangles should have the same aspect ratio 
since the Hessian is constant. In order to measure the
quality of the shape of a triangle $T$ in relation to $\bq$, we introduce
the following quantity: if $(a,b,c)$ are the edge vectors of $T$, we define
$$
\rho_\bq(T) := \frac{\max \{|\bq(a)|,|\bq(b)|,|\bq(c)|\}}{|T| \sqrt{|{\rm det}(\bq)|}},
$$
where ${\rm det}(\bq)$ is the determinant of the $2\times 2$ symmetric matrix $Q$
associated with $\bq$, i.e. such that 
$$
\bq(u)=\<Qu,u\>
$$
for all $u\in\R^2$. Using the reference triangle
and an affine change of variables, it is proved in \S \ref{secCM2} of the next chapter that 
$$
e_T(\bq)_p\sim |T|^{1+\frac 1 p}\rho_\bq(T)\sqrt{|{\rm det}(\bq)|},
$$
with equivalence constants independent of $\bq$ and $T$.
Therefore, if $T$ is a triangle of given area, its shape should be
designed in order to minimize $\rho_\bq(T)$.

In the case where $\bq$ is positive definite or negative
definite, $\rho_\bq(T)$ takes small values when $T$ is isotropic
{\it with respect to the metric} $|(x,y)|_{\bq}:=\sqrt{|\bq(x,y)|}$, the minimal value $\frac 4 {\sqrt{3}}$ being
attained for an equilateral triangle for this metric.
Specifically, we choose $\bq(x,y):=x^2+100 y^2$ and
display in Figure \ref{fig3CDHM} (left) the triangulation $\cD_{8}$ obtained after $j=8$ iterations
of the refinement procedure, starting with a triangle which is equilateral for 
the euclidean metric (and therefore not adapted to $\bq$). 
Triangles such that $\rho_\bq(T)\leq 4\sqrt 3$ (at most $3$ times the minimal value)
are displayed in white, 
others in grey. We observe that most triangles produced by
the refinement procedure are of the first type and therefore
have a good aspect ratio.

\begin{figure}
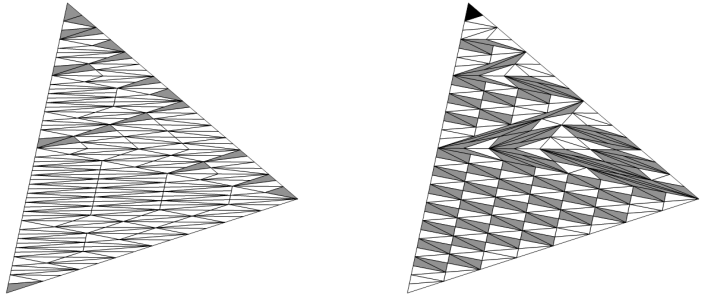

	\centering
		\includegraphics[width=4cm, height=4cm]{\pathPic/PaperCDHM/QuadPos2.pdf} 
		\hspace{1cm}
		\includegraphics[width=4cm, height=4cm]{\pathPic/PaperCDHM/QuadNeg.pdf} 
	\caption{\label{fig3CDHM} $\cD_8$ for $q(x,y):=x^2+100 y^2$ (left) and $q(x,y):=x^2-10 y^2$ (right).}
\end{figure}

The case of a
quadratic function of mixed signature is illustrated 
in Figure \ref{fig3CDHM} (right) with $\bq(x,y):=x^2-10 y^2$. 
For such quadratic functions, triangles which are isotropic with
respect to the metric $|\cdot|_{|\bq|}$ have a low value of $\rho_\bq$,
where $|\bq|$ denotes the positive quadratic form associated
to the absolute value $|Q|$ of the symmetric matrix $Q$ associated
to $\bq$.
Recall for any symmetric matrix $Q$ there exists $\lambda_1,\lambda_2\in \R$ and a rotation $R$ such that 
$$
Q = R^\trans 
\left(\begin{array}{cc}
\lambda_1 & 0\\
0 & \lambda_2
\end{array}\right)
R,
$$
and the absolute value $|Q|$ is defined as 
$$
|Q| = R^\trans
\left(\begin{array}{cc}
|\lambda_1| & 0\\
0 & |\lambda_2|
\end{array}\right)
R.
$$
In the present case, $R=I$ and $|\bq|(x,y)=x^2+10 y^2$.

But one can also check that 
$\rho_\bq$ is left invariant by any linear transformation
with eigenvalues $(t,\frac 1 t)$ for any $t>0$ and eigenvectors $(u,v)$ such that
$\bq(u)=\bq(v)=0$, i.e. belonging to 
the \emph{null cone} of $\b q$. More precisely, for any such transformation
$\psi$ and any triangle $T$, one has $\rho_{\bq}(\psi(T))=\rho_\bq(T)$.
In our example we have $u=(\sqrt {10},1)$ and $v=(\sqrt{10},-1)$).
Therefore long and thin triangles which are aligned with
these vectors also have a low value of $\rho_\bq$.
Triangles $T$
such that $\rho_{|\bq|}(T)\leq 4 \sqrt 3$ are displayed in white, those 
such that $\rho_{\bq}(T)\leq 4\sqrt 3$ while $\rho_{|\bq|}(T)> 4 \sqrt 3$ -
i.e. adapted to $\bq$ but not to $|\bq|$ - are displayed in grey,
and the others in dark. We observe that all the triangles triangles produced by
the refinement procedure except one are either of the first or second type
and therefore have a good aspect ratio. These empirical observations will be
rigorously justified in \S 9.1 of Chapter 9.

\subsection{Sharp transition}

We next study the adaptive triangulations produced by the greedy tree
algorithm for a function $f$ displaying a sharp transition along a curved edge. Specifically
we take
$$
f(x,y)=f_\delta(x,y) := g_\delta(\sqrt{x^2+y^2}),
$$ 
where $g_\delta$
 is defined by
$g_\delta(r)= \frac{5-r^2} 4$ for $0\leq r\leq 1$, $g_\delta(1+\delta+r)=-\frac{5-(1-r)^2} 4$ for $r\geq 0$, and $g_{\delta}$ is a polynomial of degree $5$ on $[1,1+\delta]$
which is determined by imposing that
$g_\delta$ is globally $C^2$.
The parameter $\delta$ therefore measures the sharpness of the transition as illustrated
in Figure \ref{fig4CDHM}.
It can be shown that the Hessian of $f_\delta$ is negative definite for $\sqrt{x^2+y^2}< 1+\delta/2$, and
of mixed type for $1+\delta/2 < \sqrt{x^2+y^2}\leq 2+\delta$.

\begin{figure}
	\centering
		\includegraphics[width=4cm, height=3cm]{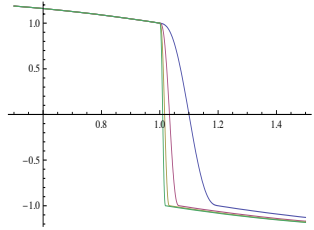} 
	\caption{\label{fig4CDHM}The function $g_\delta$, for $\delta=0.02,0.03,0.07,0.2$.}
\end{figure}

Figure \ref{fig5CDHM} displays the triangulation $\cT_{10000}$ obtained after $10000$ steps of
the algorithm for $\delta=0.2$. In particular, triangles $T$ such that $\rho_\bq(T)\leq 4$ 
- where $\bq$ is the quadratic form associated with $d^2f$ measured at the barycenter of $T$ - 
are displayed in white, others in grey. As expected, most triangles are of the
first type and therefore well adapted to $f$. We also display
on this figure the adaptive isotropic triangulation produced by the greedy tree algorithm
based on newest vertex bisection for the same number of triangles.

\begin{figure}
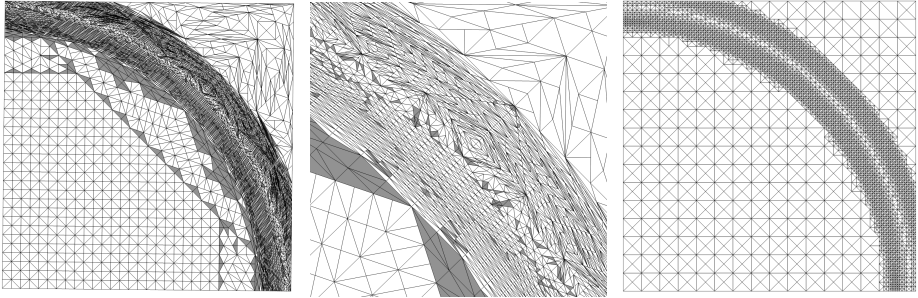

	\centering
		\includegraphics[width=4cm, height=4cm]{\pathPic/PaperCDHM/Rect42D5.pdf} 
		\includegraphics[width=3.9cm,height=3.9cm]{\pathPic/PaperCDHM/DetailR5.pdf}
		\includegraphics[width=4cm, height=4cm]{\pathPic/PaperCDHM/RectIsoD5.pdf} 
	\caption{\label{fig5CDHM}$\cT_{10000}$ (left),  detail (center), isotropic triangulation (right).}
\end{figure}

Since $f$ is a $C^2$ function, approximations by uniform, adaptive isotropic and adaptive anisotropic triangulations all yield the convergence rate $\cO(N^{-1})$.
However the constant 
$$
C:=\limsup_{N\to +\infty} N\|f-f_N\|_{L^2},
$$ 
strongly differs depending on the algorithm
and on the sharpness of the transition, as illustrated in the table below.
We denote by $C_U$, $C_I$ and $C_A$ the empirical constants
(estimated by $N \|f-f_N\|_2$ for $N=8192$)
in the uniform, adaptive isotropic and adaptive anisotropic case respectively, and by
$U(f):=\|d^2f\|_{L^{2}}$,
$I(f):=\|d^2f\|_{L^{2/3}}$ and
$A(f):=\|\sqrt{|{\rm det}(d^2f)|}\|_{L^{2/3}}$
the theoretical constants suggested by
\iref{unier}, \iref{isoer}
and \iref{aniserCDHM}. We observe that $C_U$
and $C_I$ grow in a similar way as $U(f)$ and $I(f)$
as $\delta\to 0$ (a detailed computation shows that $U(f)\approx 10.37 \, \delta^{-3/2}$
and $I(f)\approx 14.01 \, \delta^{-1/2}$).
In contrast $C_A$ and $A(f)$ remain uniformly bounded, a fact
which reflects the superiority of the anisotropic mesh
as the layer becomes thinner.

$$
\begin{array}{c|c|c|c|c|c|c|}
\delta & U(f)  & I(f) & A(f) & C_U & C_I   & C_A
\\
\hline
0.2	& 103  & 27 & 6.75 & 7.87	& 1.78 &  0.74\\
0.1 & 602 	& 60 & 8.50 & 23.7	& 2.98 & 0.92\\
0.05& 1705 & 82 & 8.48 & 65.5	&	4.13 & 0.92\\
0.02& 3670 & 105 & 8.47 & 200	&	6.60 & 0.92
\end{array}
$$

\subsection{Numerical images}

We finally apply the greedy tree
algorithm to numerical images. In this case
the data $f$ has the form of a discrete array of pixels,
and the $L^2(T)$-orthogonal projection is replaced
by the $\ell^2(S_T)$-orthogonal projection,
where $S_T$ is the set of pixels with centers contained in $T$.
The approximated $512\times 512$ image is displayed
in Figure \ref{fig6CDHM} which also shows its approximation $f_{N}$
by the greedy tree algorithm
based on newest vertex bisection with $N=2000$ triangles.
The systematic use of isotropic triangles results in
strong ringing artifacts near the edges.

\begin{figure}
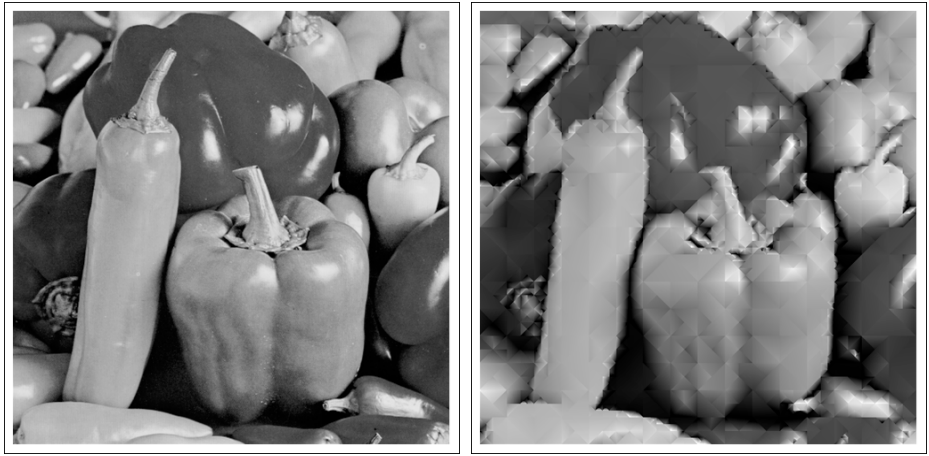

	\centering
		\includegraphics[width=6cm, height=6cm]{\pathPic/PaperCDHM/PeppersOrig2000G.pdf} 
		\includegraphics[width=6cm,height=6cm]{\pathPic/PaperCDHM/PeppersIso2000G.pdf}
	\caption{\label{fig6CDHM}The image ''peppers'' (left),  $f_{2000}$ with newest vertex (right).}
\end{figure}

\begin{figure}
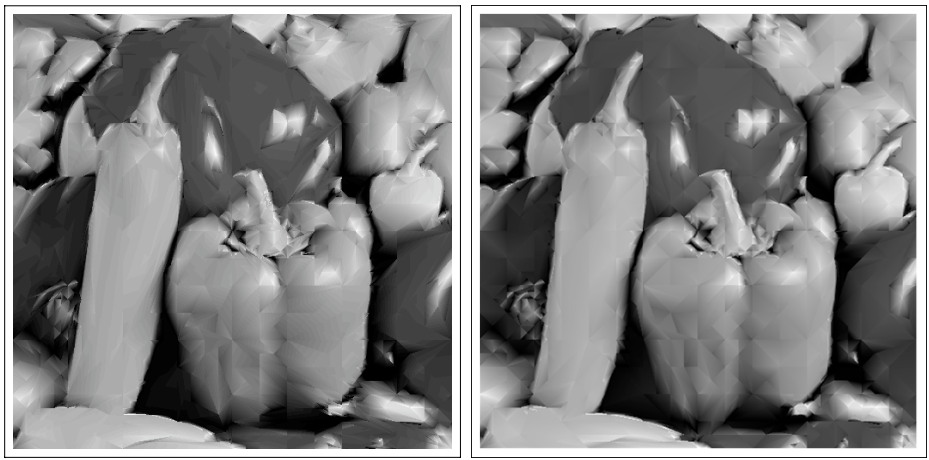

	\centering
		\includegraphics[width=6cm, height=6cm]{\pathPic/PaperCDHM/PeppersTri2000.pdf} 
		\includegraphics[width=6cm,height=6cm]{\pathPic/PaperCDHM/PeppersModif2000G.pdf}
	\caption{\label{fig7CDHM}$f_{2000}$ with greedy bisection (left),  modified procedure (right).}
\end{figure}

We display in Figure \ref{fig7CDHM} the result of the same algorithm 
now based on our greedy bisection procedure 
with the same number of triangles. As expected, the edges
are better approximated due to the presence of well oriented
anisotropic triangles. Yet artifacts persist on certain edges
due to oscillatory
features in the image which tend to mislead the algorithm
in its search for triangles with good aspect ratio, as explained
in \S \ref{secCDHM3p2}. These artifacts tend to disappear if we use the 
modified refinement rule proposed in \S \ref{secCDHM3p3} as also illustrated on
Figure \ref{fig7CDHM}. This modification is thus useful in the
practical application of the algorithm, in addition of being
necessary for proving convergence of the
approximations towards any $L^p$ function. Note that
encoding a triangulation resulting from $N$ iterations 
of the anisotropic refinement algorithm is more costly than for the newest
vertex rule: the algorithm encounters at most $2N$ triangles
and for each of them, one needs to encode one out of four options
(bisect towards edge $a$ or $b$ or $c$ or not bisect), therefore
resulting into $4N$ bits,
while only two options need to be encoded when using the newest vertex rule (bisect or not),
therefore resulting into $2N$ bits. In the perspective of
applications to image compression, another issue is the 
quantization and encoding of the piecewise affine function
as well as the treatment of the triangular visual artifacts that are inherent
to the use of discontinuous piecewise polynomials on 
triangulated domains. These issues will be discussed in 
a further work specifically dealing with image applications.

\section{Conclusions and perspectives}

In this chapter, we have studied a simple greedy 
refinement procedure which 
generates triangles that tend to 
have an {\it optimal aspect ratio}.
This fact is rigorously proved in Chapter \ref{chapBisecOpt},
together with the optimal 
convergence estimate \iref{aniserCDHM}
for the adaptive triangulations constructed by the
greedy tree algorithm in the case where 
the approximated function $f$ is
$C^2$ and convex. Our numerical results
illustrate these properties.

In the present chapter we also show
that for a general $f\in L^p$ the refinement procedure 
can be misled by 
oscillations in $f$, and that
this drawback may be circumvented by a 
simple modification of the refinement procedure.
This modification appears to be useful 
in image processing applications, as shown
by our numerical results.

Let us finally mention several perspectives that are raised
from our work, and that are the object of current investigation:
\begin{enumerate}
\item
Conforming triangulations: our algorithm inherently generates
hanging nodes, which might not be desirable in certain applications,
such as the numerical discretization of PDE's where anisotropic elements
are sometimes used \cite{Apel}.
When using the greedy tree algorithm, an obvious way of avoiding this phenomenon is to
bisect the chosen triangle together with an adjacent triangle in order
to preserve conformity. However, it is no more clear that
this strategy generates optimal triangulations. In fact, we observed
that many inappropriately oriented triangles can be generated by this
approach. An alternative strategy is to apply the non-conforming greedy tree algorithm
until a prescribed accuracy is met, followed by an additional
refinement procedure in order to remove hanging nodes.
\item
Discretization and encoding: 
our work is in part motivated by applications to image and terrain data
processing and compression. In such applications
the data to be approximated is usually given in discrete form
(pixels or point clouds) and the algorithm can be adapted
to such data, as shown in our numerical image examples. 
Key issues which need to be dealt with are then 
(i) the efficient encoding of the approximations and of the triangulations
using the tree structure in a similar spirit as in \cite{CDDD1} and
(ii) the removal of the triangular visual artifacts due to discontinuous piecewise polynomial
approximation by an appropriate post-processing step.
\item
Adaptation to curved edges: 
one of the motivation for the use of anisotropic triangulations is the 
approximation of functions with jump discontinuities along an edge.
For simple functions, such as characteristic functions of domains with smooth boundaries,
the $L^p$-error rate with an optimally adapted triangulation of $N$ elements
is known to be $\cO(N^{-\frac 2 p})$. This rate reflects an $\cO(1)$ error concentrated
on a strip of area $\cO(N^{-2})$ separating the curved edge
from a polygonal line. Our first investigations in this direction indicate that the
greedy tree algorithm based on our refinement procedure 
cannot achieve this rate, due to the fact that bisection does not offer enough
geometrical adaptation. This is in contrast with other splitting procedures,
such as in \cite{DLe} in which the direction of the new cutting edge is optimized within 
an infinite range of possible choices, or \cite{D} where the number of choices
grows together with the resolution level. An interesting question
is thus to understand if the optimal rate for edges can be achieved
by a splitting procedure with a small and fixed number 
of choices similar to our refinement procedure, which would be beneficial
from both a computational and encoding viewpoint. This question is addressed, and partially answered, in \S\ref{secAltBi}.
\end{enumerate}

\chapter{Greedy bisection generates optimally adapted triangulations} 
\minitoc
\label{chapBisecOpt}
\section{Introduction}

In finite element approximation,
a classical and important distinction is made between {\it uniform} 
and {\it adaptive} methods. In the first case
all the elements which constitute the mesh have comparable shape and size,
while these attributes are allowed to vary strongly in the second case.
An important feature of adaptive methods is the fact that 
the mesh is not fixed in advance
but rather tailored to the properties of the function $f$
to be approximated. Since the function approximating $f$ is not picked from a fixed linear space, 
adaptive finite elements can be considered as an 
instance of {\it non-linear approximation}. Other
instances include approximation by rational
functions, or by $N$-term linear combinations
of a basis or dictionary. We refer to \cite{De} for a general survey
on non-linear approximation.

In this chapter, we focus our interest
on {\it piecewise linear} finite element functions defined over triangulations
of a bidimensional polygonal domain $\Omega\subset \RR^2$.  Given a triangulation
$\cT$ we denote by $V_{\cT}:=\{v\;{\rm s.t.} \; v_{|T}\in\P_1,\; T\in\cT\}$
the associated finite element space. The norm in which we measure the approximation error is
the $L^p$ norm for $1\leq p\leq \infty$ and we therefore do not require 
that the triangulations are conforming and that the functions of $V_\cT$ are continuous
between triangles. For a given function $f$
we define 
$$
e_N(f)_{L^p}:=\inf_{\#(\cT)\leq N} \inf_{g\in V_{\cT}} \|f-g\|_{L^p},
$$
the best approximation error of $f$ when using at most $N$ elements.
In adaptive finite element approximation, critical questions are:
\begin{enumerate}
\item
Given a function $f$ and a number $N>0$, how can we characterize the {\it optimal mesh}
for $f$ with $N$ elements corresponding to the above defined best approximation error.
\item
What quantitative estimates
are available for the best approximation error $e_N(f)_{L^p}$ ? Such estimates should involve
the derivatives of $f$ in a different way than for non-adaptive meshes.
\item
Can we build by a simple algorithmic procedure
a mesh $\cT_N$ of cardinality $N$ and a finite element function
$f_N\in V_{\cT_N}$ such that $\|f-f_N\|_{L^p}$ 
is comparable to $e_N(f)_{L^p}$ ?
\end{enumerate}

While the optimal mesh
is usually difficult to characterize exactly, it should
satisfy two intuitively desirable features:
(i) the triangulation should {\it equidistribute} the local approximation error between 
each triangle and (ii) the aspect ratio of a triangle $T$ should be 
{\it isotropic} with respect to a distorted
metric induced by the local value of the hessian $d^2f$ on $T$
(and therefore anisotropic in the sense of the euclidean metric). 
Under such prescriptions on the mesh, 
quantitative error estimates have recently
been obtained in \cite{CSX,BBLS} when $f$ is 
a $C^2$ function. These estimates
are of the form
\be
\limsup_{N\to\infty} N e_N(f)_{L^p}  \leq C\|\sqrt{|\det(d^2f)|}\|_{L^\tau},\;\; \frac 1 \tau=\frac 1 p + 1,
\label{aniser}
\ee
where $\det(d^2f)$ is the determinant of the $2\times 2$ hessian matrix.
For a convex $C^2$ function $f$ this estimate has been proved to be {\it asymptotically optimal} in \cite{CSX}, in the following sense
\be
\liminf_{N\to +\infty} Ne_N(f)_{L^p} \geq c\|\sqrt{|\det(d^2f)|}\|_{L^\tau}.
\label{lowerboundasym}
\ee
The convexity assumption can actually be replaced by a mild assumption
on the sequence of triangulations which is used for the approximation of $f$: a
sequence $(\cT_N)_{N\geq N_0}$ is said to be admissible if $\#(\cT_N)\leq N$ and
$$
\sup_{N\geq N_0}\left( N^{1/2}\max_{T\in\cT_N} \diam(T) \right)< \infty.
$$
Then it is proved in Chapter \ref{chapOptAniso}, that for any admissible sequence and any $C^2$ function $f$,
one has
\be
\liminf_{N\to +\infty} N\inf_{g \in V_{\cT_N}}\|f-g\|_{L^p} \geq c\|\sqrt{|\det(d^2f)|}\|_{L^\tau}.
\label{lowerboundadmiss}
\ee
The admissibility assumption is not a severe limitation for an upper estimate
of the error since it is also proved that for all $\e>0$, there exist an admissible sequence
such that 
\be
\limsup_{N\to +\infty} N \inf_{g \in V_{\cT_N}} \|f-g\|_{L^p} \leq C\|\sqrt{|\det(d^2f)|}\|_{L^\tau}+\e.
\label{upperboundadmiss}
\ee
We also refer to Chapter \ref{chapOptAniso} for a generalization of such upper and lower estimates
to higher order elements. 

From the computational viewpoint, a commonly used strategy
for designing an optimal mesh consists therefore in
evaluating the hessian $d^2f$ and imposing that 
each triangle of the mesh is isotropic with respect
to a metric which is properly related to its local value.
We refer in particular to \cite{BFGLS}
where this program is executed using 
Delaunay mesh generation techniques.
While these algorithms fastly produce 
anisotropic meshes which are naturally
adapted to the approximated function,
they suffer from two intrinsic limitations:
\begin{enumerate}
\item
They use the data of $d^2f$,
and therefore do not apply 
to non-smooth or noisy functions. 
\item
They are non-hierarchical: for $N>M$, the triangulation
$\cT_N$ is not a refinement of $\cT_M$. 
\end{enumerate}

In Chapter \ref{chapCDHM}, an alternate strategy is proposed 
for the design of adaptive hierarchical meshes, based on a 
simple {\it greedy algorithm}: starting from an initial 
triangulation $\cT_{N_0}$, the algorithm
picks the triangle $T\in \cT_N$ with the largest local
$L^p$ error. This triangle is then bisected
from the mid-point of one of its edges to the opposite vertex. The choice of the edge
among the three options is 
the one that minimizes the new approximation
error after bisection.
The algorithm can be applied to any $L^p$ function,
smooth or not, in the context of 
piecewise polynomial approximation of any given order.
In the case of piecewise linear
approximation, numerical experiments 
in  Chapter \ref{chapCDHM} indicate that this elementary strategy
generates triangles with an optimal aspect ratio
and approximations $f_N\in V_{\cT_N}$ such that $\|f-f_N\|_{L^p}$
satisfies the same estimate as $e_N(f)_{L^p}$ in \iref{aniser}.

The goal of this chapter is to support these experimental observations
by a rigorous analysis. This chapter is organized as follows:

In \S \ref{secCM2}, we introduce 
notations which are used throughout the chapter
and collect some available approximation
theory results for piecewise linear finite elements,
making the distinction between (i) uniform,
(ii) adaptive isotropic and (iii) adaptive
anisotropic triangulations. In the last case, which
is in the scope of this chapter, we introduce a measure
of non-degeneracy of a triangle $T$ with respect
to a quadratic form. We show that the optimal
error estimate \iref{aniser} is met when each triangle
is non-degenerate in the sense of the above
measure with respect to the quadratic form
given by the local hessian $d^2f$.
We end by briefly recalling the greedy algorithm
which is introduced in Chapter \ref{chapCDHM}. 

In \S \ref{secCM3}, we study the behavior 
of the refinement procedure when applied to a
quadratic function $q$ such that its associated quadratic
form $\bq$ is of positive or negative sign.
A key observation is that the edge which is bisected
is the longest with respect to the metric induced by $\bq$.
This allows us to prove that the triangles generated 
by the refinement procedure adopt 
an optimal aspect ratio in the sense
of the non-degeneracy measure introduced in \S \ref{secCM2}.

In \S \ref{secCM4}, we study the behavior
of the algorithm when applied to a general $C^2$
function $f$ which is assumed to be
strictly convex (or strictly concave). We first establish a perturbation
result, which shows that when $f$ is locally
close to a quadratic function $q$ the algorithm
behaves in a similar manner as when applied
to $q$. We then prove that the diameters of
the triangles produced by the algorithm tend
to zero so that the perturbation result can be 
applied. This allows us to show that the
optimal convergence estimate 
\be
\limsup_{N \to \infty} N \|f-f_N\|_{L^p} \leq C \|\sqrt{|\det(d^2f)|}\|_{L^\tau}
\label{opticonvest}
\ee
is met by the sequence of approximations
$f_N\in V_{\cT_N}$ generated by the algorithm.

The extension of this result to an arbitrary $C^2$ function
$f$ remains an open problem. It is possible to
proceed to an analysis similar to \S \ref{secCM3} in the case
where the quadratic form $\bq$ is of mixed sign, also proving that
the triangles adopt an optimal aspect ratio as they get refined.
We describe this analysis in \S 9.1. However,
it seems difficult to extend the perturbation
analysis of \S \ref{secCM4} to this new setting. In particular
the diameters of the triangles are no more ensured to
tend to zero, and one can even exhibit examples of non-convex $C^2$
functions $f$ for which the approximation $f_N$ {\it fails} to converge towards $f$
due to this phenomenon. Such examples are
discussed in  Chapter \ref{chapCDHM} which also proposes
a modification of the algorithm for which convergence
is always ensured. However, we do not know if
the optimal convergence estimate \iref{opticonvest} holds for
any $f\in\cC^2$ with this modified algorithm, although this seems
plausible from the numerical experiments.

\section{Adaptive finite element approximation}
\label{secCM2}
\subsection{Notations}
\label{subsecAssumAT}
We shall make use of a linear approximation operator
$\cA_T$ that maps continuous functions defined on $T$ onto $\P_1$. 
For an arbitrary but fixed  $1\leq p\leq \infty$, we define the
local $L^p$ approximation error 
$$
e_T(f)_p:=\|f-\cA_Tf\|_{L^p(T)}.
$$
The critical assumptions in our analysis for the operator $\cA_T$ 
will be the following:
\begin{enumerate}
\item
$\cA_T$ is continuous in the $L^\infty$ norm. 
\item
$\cA_T$ commutes with affine changes of variables: $\cA_T(f) \circ \phi = \cA_{\phi^{-1}(T)} (f\circ \phi)$
for all affine $\phi$.
\item
$\cA_T$ reproduces $\P_1$:
$\cA_T(\pi)=\pi$,
for any $\pi\in\P_1$.
\end{enumerate}
Note that the commutation assumption implies that for any function $f$
and any affine transformation $\phi:x\mapsto x_0+Lx$ we have
\be
e_{\phi(T)}(f)_p=|\det(L)|^{1/p} e_T(f\circ \phi)_p,
\label{commuterror}
\ee
Two particularly simple admissible choices of approximation operators are the following:
\begin{itemize}
\item
$\cA_T=P_T$, the $L^2(T)$-orthogonal projection operator: $\int_T (f-P_T f)\pi=0$ for all $\pi\in \P_1$.
\item
$\cA_T=\interp_T$, the local interpolation operator: 
$\interp_T f(v_i)=f(v_i)$ with $\{v_0,v_1,v_2\}$ the vertices of $T$.
\end{itemize}
All our results are simultaneously valid when $\cA_T$ is either $P_T$
or $\interp_T$, or any linear operator that fulfills the three above assumptions.

Given a function $f$ and a triangulation 
$\cT_N$ with $N=\#(\cT_N)$, we can associate a finite element approximation $f_N$
defined on each $T\in\cT_N$ by
$f_N(x)=\cA_T f(x)$. The global approximation error is  given by
$$
\|f-f_N\|_{L^p}=\(\sum_{T\in\cT_N}e_T(f)_p^p\)^{\frac 1 p},
$$
with the usual modification when $p=\infty$. 

\begin{remark}
The operator $\cB_T$ of best $L^p(T)$ approximation
which is defined by
$$
\|f-\cB_Tf\|_{L^p(T)}=\min_{\pi\in\sP_m}\|f-\pi\|_{L^p(T)},
$$ 
does not fall in the above category of operators, since
it is non-linear (and not easy to compute) when $p\neq 2$. However, it is clear
that any estimate on $\|f-f_N\|_{L^p}$ with $f_N$ defined as $\cA_T f$ on each $T$
implies a similar estimate when $f_N$ is defined as $\cB_T f$ on each $T$.
\end{remark}

Here and throughout the chapter, when
$$
q(x,y)= a_{2,0} x^2+ 2 a_{1,1} xy+ a_{0,2} y^2+ a_{1,0} x+ a_{0,1} y + a_{0,0}
$$
we denote by $\bq$ the associated quadratic form : if $u=(x,y)$
$$
\bq(u)= a_{2,0} x^2+ 2 a_{1,1} xy+ a_{0,2} y^2.
$$
Note that $\bq(u)=\<Qu,u\>$ 
where 
$
Q=\left(\begin{array}{cc}
a_{2,0} & a_{1,1}\\
a_{1,1} & a_{0,2}
\end{array} \right).
$
We define
$$
\det (\bq) := \det (Q).
$$
If $\bq$ is a positive or negative quadratic form, we define the
{\it $\bq$-metric}
\be
|v|_\bq:=\sqrt{|\bq(v)|}
\label{qmetric}
\ee
which coincides with the euclidean norm when $\bq(v)=x^2+y^2$ for $v=(x,y)$.
If $\bq$ is a quadratic form of mixed sign, we define the associated
positive form $|\bq|$ which corresponds to the symmetric matrix $|Q|$
that has same eigenvectors as $Q$ with eigenvalues $(|\lambda|,|\mu|)$
if $(\lambda,\mu)$ are the eigenvalues of $Q$. Note that generally $|\bq|(u)\neq |\bq(u)|$
and that one always has $|\bq(u)|\leq |\bq|(u)$.

\begin{remark}
\label{rmCanonical}
If $\det Q>0$, then there exists a $2\times 2$ matrix $L$ and $\ve\in \{+1, -1\}$ such that 
$$
L^t Q L = \ve 
\left(
\begin{array}{cc}
1 & 0\\
0 & 1
\end{array}
\right).
$$
The linear change of coordinates $\phi(u) := Lu$, where $u = (x,y) \in \R^2$, therefore satisfies $\bq \circ \phi (u) = \ve (x^2+y^2)$. On the other hand, if $\det Q<0$ then there exists a $2 \times 2$ matrix $L$ such that 
$$
L^t Q L = 
\left(
\begin{array}{cc}
1 & 0\\
0 & -1
\end{array}
\right).
$$
Defining again $\phi(u) := Lu$ we obtain in this case $\bq \circ \phi (u) = x^2 - y^2$.
\end{remark}
\subsection{From uniform to adaptive isotropic triangulations}

A standard estimate in finite element approximation states that 
if $f\in W^{2,p}(\Omega)$ then
$$
\inf_{g\in V_h} \|f-g\|_{L^p} \leq C h^{2}\|d^2f\|_{L^{p}},
$$
where $V_h$ is the piecewise linear finite element space
associated with a triangulation $\cT_h$ of mesh size
$h:=\max_{T\in\cT_h} \diam(T)$. If we restrict our attention
to {\it uniform} triangulations, we have
$$
N:=\#(\cT_h) \sim  h^{-2}.
$$
Therefore, denoting by $e^{\rm unif}_N(f)_{L^p}$ the $L^p$ approximation error
by a uniform triangulation of cardinality $N$, we can re-express the above
estimate as
\be
\label{uniferror}
e^{\rm unif}_N(f)_{L^p} \leq CN^{-1}\|d^2f\|_{L^p}.
\ee
This estimate can be significantly improved when using
adaptive partitions. We give here some heuristic
arguments, which are based on the assumption that on each triangle $T$
the relative variation of $d^2f$ is small so that it can be considered
as constant over $T$ (which means that 
$f$ is replaced by a quadratic function on each $T$),
and we also indicate the available results which 
are proved more rigorously.

First consider {\it isotropic} triangulations,
i.e. such that all triangles satisfy a uniform estimate 
\be
\rho_T=\frac {h_T}{r_T}\leq A,
\label{regulartri}
\ee
where $h_T:=\diam(T)$ denotes the size of the longest edge of $T$, and $r_T$ is the radius of the largest disc contained in $T$. 
In such a case we start from the local approximation estimate on any $T$
$$
e_T(f)_p \leq Ch_T^2 \|d^2f\|_{L^p(T)},
$$
and notice that 
$$
h_T^2 \|d^2f\|_{L^p(T)}\sim |T|\, \|d^2f\|_{L^p(T)} = \|d^2f\|_{L^\tau(T)},
$$
with $\frac 1 \tau:=\frac 1 p +1$ and $|T|$ the area of $T$, where we have used
the isotropy assumption \iref{regulartri} in the equivalence and the fact that $d^2f$
is constant over $T$ in the equality. It follows that
$$
e_T(f)_p \leq C \|d^2f\|_{L^\tau(T)}, \;\; \frac 1 \tau:=\frac 1 p +1.
$$
Assume now that we can construct adaptive isotropic triangulations $\cT_N$ 
with $N:=\#(\cT_N)$ which {\it equidistributes} the local error in the
sense that for some prescribed $\e>0$
\be
c\e \leq e_T(f)_p  \leq \e,
\label{equil}
\ee
with $c>0$ a fixed constant independent of $T$ and $N$. Then
defining $f_N$ as $\cA_T(f)$ on each $T\in\cT_N$, we have on the one hand
$$
\|f-f_N\|_{L^p} \leq N^{1/p}\e,
$$
and on the other hand, with $\frac 1 \tau:=\frac 1 p +1$,
$$
N(c\e)^\tau \leq \sum_{T\in\cT_N}\|f-f_N\|_{L^p(T)}^\tau 
\leq C^\tau\sum_{T\in\cT_N} \|d^2f\|_{L^\tau(T)}^\tau\leq C^\tau \|d^2f\|_{L^\tau}^\tau.
$$
Combining both, one obtains for $e^{\rm iso}_N(f)_{L^p}:=\|f-f_N\|_{L^p}$
the estimate
\be
\label{isoerror}
e^{\rm iso}_N(f)_{L^p} \leq  CN^{-1}\|d^2f\|_{L^\tau}.
\ee
This estimate improves upon \iref{uniferror} since the rate $N^{-1}$ is now
obtained with the weaker smoothness condition $d^2f\in L^\tau$
and since, even for smooth $f$, the quantity $\|d^2f\|_{L^\tau}$ might
be significantly smaller than $\|d^2f\|_{L^p}$.
This type of result is classical in non-linear approximation and also occurs
when we consider best $N$-term approximation in
a wavelet basis.

The principle of error equidistribution suggests a
simple {\it greedy algorithm} to build an adaptive isotropic triangulation 
for a given $f$, similar to our algorithm but
where the bisection of the triangle $T$
that maximizes
the local error $e_T(f)_p$ is systematically done from
its {\it most recently created vertex} in order to preserve the estimate
\iref{regulartri}. Such an algorithm cannot exactly equilibrate the error
in the sense of \iref{equil} and therefore does not lead to the 
same the optimal estimate as in \iref{isoerror}.
However, it was proved in \cite{BDDP}  that it satisfies
$$
\|f-f_N\|_{L^p} \leq C |f|_{B^2_{\tau,\tau}}N^{-1},
$$
for all $\tau$ such that $\frac 1 \tau < \frac 1 p +1$, provided that the
local approximation operator $\cA_T$ is bounded in the $L^p$ norm. 
Here $B^2_{\tau,\tau}$ denotes the usual Besov space which is a 
natural substitute for $W^{2,\tau}$ when $\tau<1$. Therefore
this estimate is not far from \iref{isoerror}.

\subsection{Anisotropic triangulations: the optimal aspect ratio}

We now turn to anisotropic adaptive triangulations, and start by
discussing the optimal shape of a triangle $T$ for a given function $f$
at a given point. For this purpose,
we again replace $f$ by a quadratic function
assuming that $d^2f$ is constant over $T$. For such
a $q\in \P_2$ and its associated quadratic form $\bq$, we 
first derive an equivalent quantity for the local approximation error.
Here and as well as in \S \ref{secCM3} and \S \ref{secCM4}, we consider a triangle
$T$ and we denote by $(a,b,c)$
its edge vectors oriented in clockwise 
or anti-clockwise direction so that
$$
a+b+c=0.
$$
\begin{prop}
\label{propequiverrornondeg}
The local $L^p$-approximation error satisfies
$$
e_T(q)_p=e_T(\bq)_p\sim |T|^{\frac 1 p} \max\{|\bq(a)|,|\bq(b)|,|\bq(c)|\},
$$
where the constant in the equivalence is independent of $q$, $T$ and $p$.  
\end{prop}

\proof
The first equality is trivial since $q$ and $\bq$ differ by an affine function.
Let $\TEq$ be an equilateral triangle of area $|\TEq|=1$, and edges $a,b,c$. 
Let $E$ be the $3$-dimensional vector space of all quadratic forms.
Then the following quantities are norms on $E$, and thus equivalent: 
\be
\label{equivTEq}
e_\TEq(\bq)_p \sim \max\{|\bq(a)|,|\bq(b)|,|\bq(c)|\}.
\ee
Note that the constants in this equivalence are independent of $p$ since all $L^p(\TEq)$ norms
are uniformly equivalent on $E$. 

If $T$ is an arbitrary triangle, there exists an affine transform $\phi: x\mapsto x_0+Lx$ such that $T =\phi(\TEq)$.
For any quadratic function $q$, we thus obtain from \iref{commuterror}
$$
e_T(\bq)=e_T(q) = e_{\phi(\TEq)} (q) = |\det L|^{\frac 1 p} e_\TEq (q\circ \phi)= |\det L|^{\frac 1 p} e_\TEq (\bq\circ L)
$$
since $\bq \circ L$ is the homogeneous part of $q\circ \phi$.
By \iref{equivTEq}, we thus have
$$
e_T(\bq)
 \sim |\det L|^{\frac 1 p} \max\{|\bq(La)|,\ |\bq(Lb)|,\ |\bq(Lc)|\},
$$
where $\{a,b,c\}$ are again the edge vectors of $\TEq$.
Remarking that $|T| = |\det L|$ and that $\{La,Lb,Lc\}$ are the edge vectors of $T$, 
this concludes the proof of this proposition.
\sq
\nl
In order to describe the optimal shape of a triangle $T$ for the quadratic function
$q$, we fix the area of $|T|$ and try to minimize the error  $e_T(q)_p$
or equivalently $\max\{|\bq(a)|,|\bq(b)|,|\bq(c)|\}$. The solution to this 
problem can be found by introducing for any $\bq$ such that 
$\det(\bq)\neq 0$ the following measure of {\it non-degeneracy} for $T$:
\be
\label{rhodef}
\rho_\bq(T):=\frac {\max\{|\bq(a)|,|\bq(b)|,|\bq(c)|\}}{|T|\sqrt{|\det(\bq)|}}.
\ee
Let $\phi$ be a linear change of variables, $\bq$ a quadratic form and $T$ a triangle of edges $a,b,c$. Then $\det(\bq\circ \phi) = (\det \phi)^2 \det (\bq)$, the edges of $\phi(T)$ are $\phi(a), \phi(b), \phi(c)$ and $|\phi(T)| = |\det \phi| |T|$.
Hence we obtain
\be
\label{rhoinv}
\begin{array}{rcl}
\displaystyle \rho_{\bq\circ \phi}(T)
&=& \displaystyle \frac {\max\{|\bq\circ \phi (a)|,|\bq \circ \phi (b)|,|\bq \circ \phi (c)|\}}{|T|\sqrt{|\det(\bq \circ \phi)|}} \\
&=& \displaystyle \frac {\max\{|\bq(\phi(a))|,|\bq(\phi(b))|,|\bq(\phi(c))|\}}{|\det\phi| |T|\sqrt{|\det(\bq)|}} \\
&=& \displaystyle \rho_\bq(\phi(T)).
\end{array}
\ee


\noindent
The last equation, combined with Remark \ref{rmCanonical}, allows to reduce the study of $\rho_\bq(T)$ to two 
elementary cases by change of variable:
\begin{enumerate}
\item
The case where $\det(\bq)>0$ is reduced to
$\bq(x,y)=x^2+y^2$. 
Recall that for any triangle $T$ with edges $a,b,c$ we define $h_T := \diam(T) = \max\{|a|,|b|,|c|\}$,
with $|\cdot|$ the euclidean norm.
In this case we therefore have $\rho_\bq(T)=\frac {h_T^2}{|T|}$,
which corresponds to a standard measure of shape regularity in the sense that its boundedness
is equivalent to a property such as \iref{regulartri}. 
This quantity is minimized when the triangle $T$ is equilateral, with
minimal value $\frac 4{\sqrt 3}$ (in fact it was also proved in \cite{Chen1}
that the minimum of the interpolation error
$\|\bq-\interp_T\bq\|_{L^p(T)}$ among all triangles of area $|T|=1$ is attained 
when $T$ is equilateral). For a general quadratic form $\bq$ of positive
sign, we obtain by change of variable that the minimal value $\frac  4{\sqrt 3}$
is obtained for triangles which are equilateral with respect to the 
metric $|\cdot|_{\bq}$. More generally triangles with a good aspect ratio,
i.e. a small value of $\rho_{\bq}(T)$,
are those which are {\it isotropic with respect to this metric}. Of course, a similar
conclusion holds for a quadratic form of negative sign. 
\item
The case where $\det(\bq)<0$ is reduced to
$\bq(x,y)=x^2-y^2$. In this case, the analysis presented in \cite{Cao}
shows that the quantity $\rho_\bq(T)$ is minimized when
$T$ is a half of a square with sides parallel to the $x$ and $y$ axes,
with minimal value $2$. But using \iref{rhoinv} we also notice that 
$\rho_\bq(T) = \rho_\bq (L(T))$ for any linear transformation $L$
such that $\bq = \bq\circ L$. This holds if $L$ has eigenvalues $(\lambda,\frac 1 \lambda)$, 
where $\lambda \neq 0$, and eigenvectors $(1,1)$ and $(-1,1)$. 
Therefore, all images of the half square by such transformations
$L$ are also optimal triangles. Note that such triangles can be highly anisotropic.
For a general quadratic form $\bq$ of mixed sign, we notice that $\rho_\bq(T)\leq \rho_{|\bq|}(T)$, 
and therefore triangles which are equilateral
with respect to the metric $|\cdot|_{|\bq|}$ have 
a good aspect ratio, i.e. a small value of $\rho_\bq(T)$. 
In addition, by similar arguments, we find that all images
of such triangles by linear transforms $L$ with eigenvalues
$(\lambda,\frac 1 \lambda)$ and eigenvectors $(u,v)$ such that $\bq(u)=\bq(v)=0$
also have a good aspect ratio, since 
$\bq = \bq\circ L$ for such transforms.
\end{enumerate}
We leave aside the special case where $\det (\bq)=0$. In such a case, 
the triangles minimizing the error for a given area degenerate
in the sense that they should be infinitely long and thin,
aligned with the direction of the null eigenvalue of $\bq$.

Summing up, we find that triangles with a good aspect ratio
are characterized by the fact that $\rho_\bq(T)$ is small.
In addition, from Proposition \ref{propequiverrornondeg} and the definition
of $\rho_\bq(T)$, we have
\be
e_T(q)_p\sim |T|^{1+\frac 1 p}\sqrt{|\det(\bq)|}\rho_\bq(T) = \|\sqrt{|\det(\bq)|}\|_{L^\tau(T)}
\rho_\bq(T),\;\;  \frac 1 \tau:=\frac 1 p+1.
\label{errorrhoq}
\ee
We now return to a function $f$ such that $d^2f$ is assumed to be constant
on every $T\in\cT_N$. Assuming that all triangles have a good
aspect ratio in the sense that
$$
\rho_\bq(T) \leq C
$$
for some fixed constant $C$ and with $\bq$ the value of $d^2f$ over $T$, we find
up to a change in $C$ that 
\be
e_T(f)_p\leq C\|\sqrt{|\det(d^2f)|}\|_{L^\tau(T)}
\ee
By a similar reasoning as with isotropic triangulations, we now obtain
that if the triangulation equidistributes the error in the sense of \iref{equil}
\be
\label{aniserror2}
\|f-f_N\|_{L^p} \leq  CN^{-1}\|\sqrt{|\det(d^2f)|}\|_{L^\tau},
\ee
and therefore \iref{aniser} holds. This estimate improves
upon \iref{isoerror} since the quantity $\|d^2f\|_{L^\tau}$
might be significantly larger than $\|\sqrt{|\det(d^2f)|}\|_{L^\tau}$, in particular
when $f$ has some anisotropic features, such as sharp gradients
along curved edges.

The above derivation of \iref{aniser} is heuristic and non-rigorous.
Clearly, this estimate cannot be valid as such since 
$\det(d^2f)$ may vanish while the approximation error does not
(consider for instance $f$ depending only on a single variable). 
More rigorous versions were derived in 
\cite{CSX} and \cite{BBLS}. 
In these results $|d^2f|$ is typically replaced by a majorant $|d^2f|+\e I$,
avoiding that its determinant vanishes.
The estimate \iref{aniser} can then be rigorously proved but 
holds for $N\geq N(\e,f)$ large enough. This limitation
is unavoidable and reflects the
fact that enough resolution is needed so that the hessian
can be viewed as locally constant over each optimized triangle.
Another formulation,
which is rigorously proved in Chapter \ref{chapOptAniso}, reads as follows.

\begin{prop}
There exists an absolute constant $C>0$ such that for any polygonal domain $\Omega$ and any function $f\in C^2(\overline \Omega)$, one has
$$
\limsup_{N\to +\infty} Ne_N(f)_{L^p} \leq C\|\sqrt{|\det(d^2f)|}\|_{L^\tau}.
$$
 \end{prop}
\subsection{The greedy algorithm}

Given a target function $f$, our algorithm iteratively builds triangulations
$\cT_N$ with $N=\#(\cT_N)$ and finite element approximations $f_N$.
The starting point is
a coarse triangulation $\cT_{N_0}$. Given $\cT_N$, the algorithm
selects the triangle $T$ which maximizes the local error $e_T(f)_p$
among all triangles of $\cT_N$, and bisects it from the mid-point of one of its edges towards the opposite vertex.
This gives the new triangulation $\cT_{N+1}$.

The critical part of the algorithm lies in the choice of the edge
$e\in\{a,b,c\}$ from which $T$ is bisected. Denoting by
$T_e^1$ and $T_e^2$
the two resulting triangles, we choose $e$
as the minimizer of  a {\it decision function} $d_T(e,f)$, which role is
to drive the generated triangles towards an optimal aspect ratio.
While the most natural choice for $d_T(e,f)$ corresponds to the split
that minimizes the error after bisection, namely
$$
d_T(e,f)=e_{T_e^1}(f)_p^p+e_{T_e^2}(f)_p^p,
\label{optisplit}
$$
we shall instead focus our attention on a decision function which 
is defined as the $L^1$ norm of the interpolation error
\be
d_T(e,f)=\|f-\interp_{T_e^1}f\|_{L^1(T_e^1)}+\|f-\interp_{T_e^2}f\|_{L^1(T_e^2)}.
\label{optil1}
\ee
For this decision, the analysis of the 
algorithm is made simpler, due to the fact that
we can derive explicit expressions of $\|f-\interp_T f\|_{L^1(T)}$ when $f=q$ is a quadratic polynomial
with a positive homogeneous part $\bq$. We
prove in \S \ref{secCM3} that this choice leads to triangles with an optimal aspect ratio
in the sense of a small $\rho_\bq(T)$. This leads us in \S \ref{secCM4} to a proof that
the algorithm satisfies the optimal convergence estimate \iref{aniserror2}
in the case where $f$ is $C^2$ and strictly convex.

\begin{remark}
It should be well understood that while the decision function is based on the $L^1$ norm, the
selection of the triangle to be bisected is done by maximizing $e_T(f)_p$. The algorithm remains
therefore governed by the $L^p$ norm 
in which we wish to minimize the error $\|f-f_N\|_p$ for a given number of triangles. Intuitively, this
means that the $L^p$-norm influences the size of the triangles which have to equidistribute
the error, but not their optimal shape.
\end{remark}

\begin{remark}
 It was pointed out to us that
the $L^1$ norm of the interpolation error to a suitable convex function is 
also used to improve the mesh in the context of moving grid techniques, see \cite{Chen2}.
\end{remark}

\noindent
We define a variant of the decision function as follows 
$$
D_T(e,f) := \|f-\interp_Tf\|_{L^1(T)} - d_T(e,f).
$$
Note that $D_T(e,f)$ is the reduction of the $L^1$ interpolation error 
resulting from the bisection of the edge $e$, and that the selected edge that minimizes $d_T(\cdot,f)$ is also the one that maximizes $D(\cdot,f)$. 
The function $D_T$ has a simple expression in the case where 
$f$ is a convex function.
\begin{lemma}
Let $T$ be a triangle and let $f$ be a  convex function on $T$. Let $e$ be an edge of $T$ with endpoints $z_0$ and $z_1$. 
Then 
\be
\label{DDiff}
D_T(e,f) = \frac {|T|} 3  \left(\frac{f(z_0)+f(z_1)} 2-f\left(\frac {z_0+z_1} 2 \right)\right).
\ee
If in addition $f$ has $C^2$ smoothness, we also have
\be
\label{Dint01}
D_T(e,f) = \frac {|T|} 6 \int_0^1 \<d^2 f(z_t) e, e\> \min\{t,1-t\} dt, \text{ where } z_t := (1-t)z_0 + t z_1.
\ee

\end{lemma}

\proof
Since $f$ is convex, we have $\interp_T f\geq f$ on $T$, hence 
$$
 \|f-\interp_Tf\|_{L^1(T)}  = \int_T (\interp_T f-f).
$$
Similarly $\interp_{T_e^1} f\geq f$ on $T_e^1$ and $\interp_{T_e^2} f\geq f$ on $T_e^2$, hence 
$$
D_T(e,f) = \int_T \interp_T f - \int_{T_e^1} \interp_{T_e^1} f - \int_{T_e^2} \interp_{T_e^2} f.
$$
Let $z_2$ be the vertex of $T$ opposite the edge $e$. 
Since the function $f$ is convex, it follows the previous expression that 
$D_T(e,f)$ is the volume of the tetrahedron of vertices 
$$
\left(\frac{z_{0,x}+ z_{1,x}} 2 ,\frac{ z_{0,y}+ z_{1,y}} 2 , f\left(\frac{z_0+z_1} 2\right)\right) \text{ and } (z_{i,x},z_{i,y}, f(z_i)) \text{ for } i=0,1,2.
$$
where $(z_{i,x},z_{i,y})$ are the coordinates of $z_i$.
Let $u = z_0-z_2$ and $v=z_1-z_2$. 
We thus have $D_T(e,f) = \frac 1 6 |\det(M)|$ where 
$$
M:= \left(
\begin{array}{ccc}
u_x & v_x & \frac {u_x+v_x} 2\\
u_y & v_y & \frac {u_y+v_y} 2\\
f(z_0) - f(z_2) & f(z_1)-f(z_2) & f\left(\frac {z_0+z_1} 2 \right) - f(z_2)
\end{array} 
\right).
$$
Subtracting the half of the first two columns to the third one we find that $M$ has the same
determinant as
$$
\ti M:=\left(
\begin{array}{ccc}
u_x & v_x & 0\\
u_y & v_y & 0\\
f(z_0) - f(z_2) & f(z_1)-f(z_2) & f\left(\frac {z_0+z_1} 2 \right) - \frac{f(z_0)+f(z_1)} 2
\end{array} 
\right).
$$
Recalling that $2|T| = |\det(u,v)|$ we therefore obtain 
\iref{DDiff}. In order to
establish \iref{Dint01},
we observe that we have in the distribution sense $\partial_t^2 (\min\{t,1-t\}^+) = \delta_0-2 \delta_{1/2}+ \delta_1$, 
where $\delta_t$ is the one-dimensional Dirac function at a point $t$. Hence for any univariate function $h\in C^2([0,1])$, we have 
$$
\int_0^1 h''(t) \min\{t,1-t\} dt = h(0)-2h(1/2)+ h(1).
$$
Combining this result with \iref{DDiff} we obtain \iref{Dint01}.
\sq

\section{Positive quadratic functions}
\label{secCM3}
In this section, we study the algorithm 
when applied to a quadratic polynomial $q$ such that $\det(\bq)>0$. 
We shall assume without loss of generality that
$\bq$ is positive definite, since all our results
extend in a trivial manner to the negative definite case.

Our first observation is that 
the refinement procedure based on the 
decision function \iref{optil1} 
always selects for bisection the {\it longest edge} in the sense of the $\bq$-metric $|\cdot|_{\bq}$
defined by \iref{qmetric}.

\begin{lemma}
\label{lemmaLongestEdgeL1}
An edge $e$ of $T$ maximizes $D_T(e,q)$ among all edges of $T$
if and only if it maximizes $|e|_{\bq}$ among all edges of $T$.
\end{lemma}

\proof
The hessian $d^2q$ is constant and for all $e\in \R^2$ one has 
$$
\<d^2 q e,e\> = 2\bq(e).
$$
If $e$ is an edge of a triangle $T$, and if $q$ is a convex quadratic function, equation \iref{Dint01} 
therefore gives
\be
\label{eqDq}
D_T(e,q) = \frac {|T|} {3} \bq(e) \int_0^1\min\{t,1-t\}dt= \frac {|T|} {12} |e|_\bq^2.
\ee
This concludes the proof. 
\sq

It follows from this lemma that the longest edge of $T$ in the sense of the $\bq$-metric
is selected for bisection by the decision function. In the remainder of this
section, we use this fact to prove that the refinement procedure
produces triangles which tend to adopt an
optimal aspect ratio in the sense that $\rho_\bq(T)$ becomes small
in an average sense.

For this purpose, it is convenient to introduce a close variant to 
$\rho_{\bq}(T)$: if $T$ is a triangle with edges $a,b,c$, such that $|a|_\bq\geq |b|_\bq \geq |c|_\bq$, 
we define 
\be
\sigma_\bq(T):= \frac{\bq(b)+\bq(c)}{4 |T| \sqrt{\det \bq}}= \frac{|b|_\bq^2+|c|_\bq^2}{4 |T| \sqrt{\det \bq}}.
\label{sigmaq}
\ee
Using the inequalities $|b|^2_\bq+|c|^2_\bq\leq 2|a|^2_\bq$ and 
$|a|^2_\bq \leq 2(|b|^2_\bq+|c|^2_\bq)$, we obtain the equivalence
\be
 \frac {\rho_\bq(T)} 8 \leq \sigma_\bq(T)\leq \frac {\rho_\bq(T)} 2.
 \label{equivnondeg}
\ee
Similar to $\rho_\bq$, this quantity is invariant under a linear coordinate changes $\phi$,
in the sense that
$$
\sigma_{\bq\circ \phi}(T) = \sigma_\bq (\phi(T)),
$$
From \iref{errorrhoq} and \iref{equivnondeg}
we can relate $\sigma_\bq$ to the local
approximation error.
\begin{prop}
\label{propequiverrornondeg1}
There exists a constant $C_0$, which depends only on the choice of $\cA_T$, such that for any triangle $T$, quadratic function $q$ and exponent $1 \leq p \leq \infty$, 
the local $L^p$-approximation error satisfies
\be
\label{eqEQSigma}
C_0^{-1} e_T(q)_p \leq \sigma_\bq(T) \|\sqrt{\det \bq}\|_{L^\tau(T)} \leq C_0 e_T(q)_p.
\ee
where $\frac 1 \tau:=\frac 1 p+1$.
\end{prop}
Our next result shows that $\sigma_\bq(T)$ is always reduced by 
the refinement procedure.

\begin{prop} \label{PropSigmaDec}
If $T$ is a triangle with children $T_1$ and $T_2$ obtained by the
refinement procedure for the quadratic function $q$, then
$$
\max\{\sigma_\bq(T_1), \sigma_\bq(T_2)\}\leq \sigma_\bq(T).
$$
\end{prop}

\proof
Assuming that $|a|_\bq\geq |b|_\bq\geq |c|_\bq$, we know that the edge $a$ is cut
and that the children have area $|T|/2$ and edges $a/2,b,(c-b)/2$ and $a/2, (b-c)/2,c$ 
(recall that $a+b+c=0$). We then have
\begin{eqnarray}
	2 |T| \sqrt{\det \bq}\ \sigma_\bq(T_i) & \leq & \bq\left(\frac a 2\right)+\bq\left(\frac{b-c}{2}\right)\\
	& = & \bq\left(\frac{b+c}{2}\right)+\bq\left(\frac{b-c}{2}\right)\\
	& = & \frac{\bq(b)+\bq(c)}{2}\\
	& = & 2 |T|  \sqrt{\det \bq}\ \sigma_\bq(T).
	\end{eqnarray}
\sq

\begin{remark}
For any positive definite quadratic form $\bq$, the minimum of $\sigma_\bq$ is $1$,
as is easily seen using the inequality 
$$
2|\det (b,c)|\leq 2|b| |c| \leq |b|^2+ |c|^2, 
$$
in which we have equality if and only if $b,c\in \R^2$ are orthogonal vectors of the same norm.
When $\bq$ is the euclidean metric, the triangle that minimizes $\sigma_\bq$
is thus the half square. 
This is consistent with the above result since it
is the only triangle which is similar (i.e. identical 
up to a translation, a rotation and a dilation) to both of its children after one step of longest edge 
bisection. 
\end{remark}

\begin{remark}
A result of similar nature was already proved in \cite{Ri} : 
longest edge bisection
has the effect that the minimal angle in any triangle 
after an arbitrary number of refinements is at most twice
the minimal angle of the initial triangle.
\end{remark}

Our next objective is to show that as we iterate the refinement process,
the value of $\sigma_\bq(T)$ 
becomes bounded independently of $q$
for almost all generated triangles. For this purpose we introduce the following notation:
if $T$ is a triangle with edges such that $|a|_\bq \geq |b|_\bq \geq |c|_\bq$,
we denote by $\psi_\bq(T)$ the subtriangle of $T$ obtained after bisection of $a$
which contains the smallest edge $c$.
We first establish inequalities between the measures $\sigma_\bq$ and $\rho_\bq$
applied to $T$ and $\psi_\bq(T)$.

\begin{prop}
Let $T$ be a triangle, then
\begin{eqnarray}
\label{sigma58rho}
\sigma_\bq(\psi_\bq(T)) & \leq & \frac 5 8 \rho_\bq(T)\\
\label{rhodec}
\rho_\bq(\psi_\bq(T))   & \leq & \frac{\rho_\bq(T)} 2 \left(1+\frac{16}{\rho_\bq^2(T)}\right)
\end{eqnarray}
\end{prop}

\proof
We first prove \iref{sigma58rho}. Obviously,
$\psi_\bq(T)$ contains one edge $s\in \{a,b,c\}$ from $T$, 
and one half edge $t\in \{\frac a 2,\frac b 2,\frac c 2\}$ from $T$. Therefore
$$
\sigma_\bq(\psi_\bq(T)) \leq \frac{|s|_\bq^2+|t|_\bq^2}{4|\psi_\bq(T)|\sqrt{\det \bq}} \leq \frac{|a|_\bq^2+|\frac a 2|_\bq^2}{2|T|\sqrt{\det \bq}} = \frac 5 8 \rho_\bq(T).
$$
For the proof of \iref{rhodec}, we restrict our attention
to the case $\bq= x^2+y^2$, without loss of generality
thanks to the invariance formula \iref{rhoinv}.
Let $T$ be a triangle with edges $|a|\geq |b|\geq |c|$.
If $h$ is the width of $T$ in the direction perpendicular to $a$, then 
$$
h = \frac{2|T|}{|a|} = \frac{2 |a|}{\rho_\bq(T)}.
$$
The sub-triangle $\psi_\bq(T)$ of $T$ has edges $\frac a 2, c , d$ where $d=\frac{b-c} 2$, and the angles at the ends of $\frac a 2$ are acute. Indeed
$$
\<c,a/2\> = \frac 1 4 \left(|b|^2 - |a|^2-|c|^2\right)\leq 0 \stext{ and } \<d,a/ 2\> = \frac 1 4 \left(|c|^2-|b|^2\right)\leq 0. 
$$
By Pythagora's theorem we thus find
$$
\max\{\left|\frac a 2\right|^2, |c|^2,|d|^2\} \leq \left|\frac a 2\right|^2 +h^2 = \frac{|a|^2} 4 \left(1+\frac{16}{\rho_\bq^2(T)}\right).
$$ 
Dividing by the respective areas of $T$ and $\psi_\bq(T)$, we obtain the announced result.
\sq

Our next result shows that a significant reduction of $\sigma_\bq$ occurs at least for 
one of the triangles obtained by three successive refinements, unless it has reached a
small value of $\sigma_\bq$. We use the notation $\psi_\bq^2(T) := \psi_\bq(\psi_\bq(T))$ and $\psi_\bq^3(T) := \psi_\bq(\psi_\bq^2(T))$.

\begin{prop}
\label{PropFSigmaDec}
Let $T$ be a triangle such that $\sigma_\bq(\psi_\bq^3(T))\geq 5$. Then $\sigma_\bq(\psi_\bq^3(T))\leq 0.69 \sigma_\bq(T)$.
\end{prop}

\proof
The monotonicity of $\sigma_\bq$ established in Proposition \iref{PropSigmaDec} implies that 
$$
5 \leq \sigma_\bq(\psi_\bq^3(T)) \leq \sigma_\bq(\psi_\bq^2(T)) \leq \sigma_\bq(\psi_\bq(T)).
$$
Combining this with inequality \iref{sigma58rho} we obtain 
$$
8\leq \min\{\rho_\bq(\psi_\bq^2(T)), \ \rho_\bq(\psi_\bq(T)), \ \rho_\bq(T)\}.
$$
If a triangle $S$ obeys $\rho_\bq(S)\geq 4$, then 
$$
\frac 1 2 \left(1+\frac{16}{\rho_\bq^2(S)}\right)\leq 1
$$
and therefore $\rho_\bq(\psi_\bq(S)) \leq \rho_\bq(S)$ according to inequality \iref{rhodec}. We can apply
this to $S = \psi_\bq(T)$ and $S=T$, therefore obtaining
\be
\label{eqMonot}
\rho_\bq(\psi_\bq^2(T)) \leq \rho_\bq(\psi_\bq(T))  \leq \rho_\bq(T).
\ee
We now remark that inequality \iref{rhodec} is equivalent to 
$
(\rho_\bq(S)- \rho_\bq(\psi_\bq(S)))^2 \geq  \rho_\bq(\psi_\bq(S))^2 - 16,
$
hence
\be
\label{decSqrt}
\rho_\bq(S)\geq \rho_\bq(\psi_\bq(S))+ \sqrt{\rho_\bq(\psi_\bq(S))^2-16}
\ee
provided that $\rho_\bq(S) \geq \rho_\bq(\psi_\bq(S))$. 
Applying this to $S=\psi_\bq(T)$ and recalling that $\rho_\bq(\psi_\bq^2(T))\geq 8$ we obtain 
$$
\rho_\bq( \psi_\bq(T))\geq 8+ \sqrt{8^2-16} \geq 14.9.
$$
Applying again \iref{decSqrt} to $S = T$ we obtain 
$$
\rho_\bq(T)\geq 14.9+\sqrt{14.9^2-16} \geq  29.3.
$$
Using \iref{rhodec}, it follows that
$$
\frac{\rho(\psi_\bq^3(T))}{\rho(T)} \leq \frac 1 8 \left(1+\frac{16}{\rho_\bq^2(\psi_\bq^2(T))}\right)\left(1+\frac{16}{\rho_\bq^2(\psi_\bq(T))}\right)\left(1+\frac{16}{\rho_\bq^2(T)}\right)\leq 0.171.
$$
Eventually, the inequalities \iref{equivnondeg} imply that
$$
2 \sigma_\bq(\psi_\bq^3(T))\leq \rho_\bq(\psi_\bq^3(T)) \leq 0.171 \rho_\bq(T) \leq 0.171 (8\sigma_\bq(T))
$$
which concludes the proof.\sq
\nl
An immediate consequence of Propositions \ref{PropSigmaDec} and \ref{PropFSigmaDec} is the following.

\begin{corollary} \label{PropControlSigma2}
If $(T_i)_{i=1}^{8}$ are the eight children obtained from three successive refinement procedures 
from $T$ for the function $q$, then
\begin{itemize}
	\item for all $i$, $\sigma_\bq(T_i)\leq \sigma_\bq(T)$,
	\item there exists $i$ such that $\sigma_\bq(T_i)\leq 0.69\sigma_\bq(T)$ or $\sigma_\bq(T_i)\leq 5$.\\
\end{itemize}
\end{corollary}
We are now ready to prove that most triangles tend to adopt an optimal aspect ratio
as one iterates the refinement procedure.

\begin{theorem} \label{CorolStabilise}
Let $T$ be a triangle, and $\bq$ a positive definite 
quadratic function. Let  $k = \frac{\ln\sigma_\bq(T)-\ln 5}{-\ln(0.69)}$.
Then after $n$ applications of the refinement procedure starting from $T$,
at most $C n^k 7^{n/3}$ of the $2^n$ generated triangles satisfy $\sigma_\bq(S)\geq 5$, 
where $C$ is an absolute constant. Therefore the proportion 
of such triangles tends exponentially fast to $0$ as $n\to +\infty$.
\end{theorem}

\proof
If we prove the proposition for $n$ multiple of $3$, then it will hold for all $n$ (with a larger constant) since $\sigma_\bq$ decreases at each refinement step. We now assume that $n=3m$, and consider the octree with root $T$ obtained by only considering the 
triangles of generation $3 i$ for $i=0,\cdots,n$. 

According to 
Corollary \ref{PropControlSigma2}, for each node of this tree, 
one of its eight children either checks $\sigma_\bq\leq 5$
or has its non-degeneracy measure diminished by a factor $\theta:=0.69$.
We remark that if $\sigma_\bq$ is diminished at least
$k$ times on the path going from the root $T$ to a leaf $S$,
 then $\sigma_\bq(S)\leq 5$.
As a consequence, the number $N(m)$ of triangles $S$
which are such that $\sigma_\bq(S)> 5$ within the generation level $n=3m$ 
is bounded by the number of words in an eight letters 
alphabet $\{a_1,\cdots,a_8\}$ with length $m$ 
and that use the letter $a_8$ at most $k$ times, namely
$$
N(m)\leq \sum_{l=0}^k \binom m l 7^{m-l}
\leq  C m^k 7^m,
$$
which is the announced result.
\sq

The fact that most triangles tend to adopt an
optimal aspect ratio as one iterates the refinement procedure
is a first hint that the approximation
error in the greedy algorithm might satisfy the estimate
\iref{aniser} corresponding to an optimal triangulation.
The following result shows that this is indeed the case,
when this algorithm is applied
on a triangular domain $\Omega$ to a quadratic function $q$
with positive definite associated quadratic form $\bq$. 
The extension of this result to more general $C^2$ 
convex functions on polygonal domains requires a more involved analysis
based on local perturbation arguments
and is the object of the next section. 

\begin{corollary}
\label{corolQuadApprox}
Let $\Omega$ be a triangle, and let $q$ be a quadratic function with positive definite associated quadratic form $\bq$. 
Let $q_N$ be the approximant of $q$ on $\Omega$ obtained by the greedy algorithm 
for the $L^p$ metric, using the $L^1$ decision function \iref{optil1}.
Then 
$$
\limsup_{N\to \infty} N \|q-q_N\|_{L^p(\Omega)}  \leq C \|\sqrt{\det(\bq)}\|_{L^\tau(\Omega)},
$$
where $\frac 1 \tau=\frac 1 p +1$ and where the constant $C$ depends only on 
on the choice of the approximation operator $\cA_T$ used in the definition of the approximant.  
\end{corollary}

\proof
For any triangle $T$, quadratic function $q \in \P_2$, and exponent $p$, let 
$$
e'_T(q)_p := \inf_{\pi\in \sP_1} \|q-\pi\|_{L^p(T)}
$$
be the error of best approximation of $q$ on $T$.
Let $\TO$ be a fixed triangle of area $1$, then for any $q\in \P_2$ and $1\leq p\leq \infty$ one has 
$$
e'_\TO(q)_1 \leq e'_\TO(q)_p \leq e_\TO(q)_p \leq e_\TO(q)_\infty.
$$
Furthermore, $e'_\TO(\cdot)_1$ and $e_\TO(\cdot)_\infty$ are semi norms on the finite dimensional space $\P_2$ which vanish precisely on the same subspace of $\P_2$, namely $\P_1$. Hence these semi-norms are equivalent.
It follows that  
\be
\label{eqAOpt}
c_0 \, e_\TO(q)_p \leq e'_\TO(q)_p \leq e_\TO(q)_p
\ee
where $c_0$ is independent $q\in \P_2$ and of $p\geq 1$.
Using the invariance property \iref{commuterror}
we find that \iref{eqAOpt} holds for any triangle $T$ in place of $\TO$
with the same constant $c_0$. We also define for any triangulation $\cT$,
$$
e_\cT(f)_p^p := \sum_{T \in \cT} e_T(f)_p^p \quad \text{ and } \quad e'_\cT(f)_p^p := \sum_{T \in \cT} e'_T(f)_p^p,
$$
and we remark that 
$
c_0 \, e_\cT(q)_p \leq e'_\cT(q)_p \leq e_\cT(q)_p.
$
For each $n$, we denote by $\cT_n^u$ the triangulation of $\Omega$ produced by $n$ 
successive refinements based on the $L^1$ decision function \iref{optil1}
for the quadratic function $q$ of interest (note that $\#(\cT_n^u)=2^n$). We also define 
$
\cT_n^\sigma := \{T\in \cT_n^u\sep \sigma_{\bq}(T) > 5\}.
$
Therefore  $\sigma_{\bq}(T) \leq 5$ if $T\notin\cT_n^\sigma$, and
on the other hand we know from Proposition \ref{PropSigmaDec} that 
$\sigma_{\bq}(T) \leq \sigma_{\bq}(\Omega)$ for any $T\in \cT_n^u$.
It follows from Proposition \ref{propequiverrornondeg1} that 
$$
\begin{array}{ll}
e_{\cT_n^u}(q)_p & \leq C_0\(\sum_{T\in \cT_n^u} (\sigma_{\bq}(T)  |T|^{\frac 1 \tau}\sqrt {\det \bq} )^p\)^{\frac 1 p}\\
& \leq C_0 \( 5^p\times  2^n + \sigma_{\bq}(\Omega)^p\#(\cT_n^\sigma)\)^{\frac 1 p} \left(\frac{|\Omega|}{2^n}\right)^{\frac 1 \tau} \sqrt {\det \bq},
\end{array}
$$
where $C_0$ is the constant in \iref{eqEQSigma}. According to Theorem \ref{CorolStabilise},
we know that 
$$
\lim_{n\to +\infty} 2^{-n}\#(\cT_n^\sigma)=0.
$$
Hence 
$$
\limsup_{n \to \infty} 2^n e_{\cT_n^u}(q)_p \leq 5 C_0 \, |\Omega|^{\frac 1 \tau} \sqrt{\det \bq} = 5 C_0 \|\sqrt {\det \bq}\|_{L^\tau(\Omega)}.
$$
We now denote by $\cT^g_n$ the triangulation generated by the greedy procedure with stopping criterion based on the error $\eta_n := C_0^{-1}2^{-\frac n \tau}   \|\sqrt {\det \bq}\|_{L^\tau(\Omega)}$. 
It follows from \iref{eqEQSigma} that for all $T\in \cT_k^u$ with $k\leq n$, one has
$$
e_T(q)_p\geq C_0^{-1}\sigma_{\bq}(T) \|\sqrt {\det \bq}\|_{L^\tau(T)}\geq 
C_0^{-1}2^{-\frac k \tau}\|\sqrt {\det \bq} \|_{L^\tau(\Omega)}\geq \eta_n,
$$
where we used that $|T| = 2^{-k} |\Omega|$ and that the minimal value of $\sigma_\bq$ is $1$.
This shows that $\cT_g^n$ is a refinement of $\cT_n^u$. Furthermore any triangle  $T \in \cT_n^u$ has at most $2^{k(T)}$ children in $\cT_n^g$, where $k(T)$ is the smallest integer such that 
$$
\eta_n \geq  C_0 \, 2^{-\frac{n+k(T)} \tau } \, \sigma_{\bq}(T)  \|\sqrt {\det \bq}\|_{L^\tau(\Omega)}.
$$ 
Since $\frac 1 2 \leq \tau\leq 1$ we obtain $2^{k(T)} \leq 2^{\frac{k(T)} \tau} \leq 2^{\frac 1 \tau}C_0^2 \sigma_{q}(T)
\leq 4C_0^2 \sigma_{q}(T)$. Hence 
$$
\#(\cT_n^g) \leq 4  C_0^2 \sum_{T \in \cT_n^u} \sigma_{q}(T) \leq 4 C_0^2 \left(5 \times 2^n + \sigma_{\bq}(\Omega) \#(\cT_n^\sigma)\right) = C_1 2^n (1+\ve_n),
$$
where $C_1 = 20\,C_0^2$ and $\ve_n \to 0$ as $n \to \infty$.
If $\cT_N$ is the triangulation generated after $N$ steps
of the greedy algorithm, then there exists $n\geq 0$ such that $\cT_N$ is a refinement of $\cT_n^g$ (hence a refinement of $\cT_n^u$) and $\cT_{n+1}^g$ is a refinement of $\cT_N$. It follows that $\#(\cT_N) \leq \#(\cT_{n+1}^g) \leq C_1 2^{n+1}(1+ \ve_{n+1})$, and 
$$
c_0 \, e_{\cT_N}(q)_p \leq e'_{\cT_N}(q)_p  \leq e'_{\cT_n^u} (q)_p \leq e_{\cT_n^u} (q)_p,
$$
where we have used the fact that $e'_{\cT}(f)_p\leq e'_{\ti\cT}(f)_p$ whenever $\cT$ is a refinement
of $\ti\cT$. Eventually,
$$
\limsup_{N\to \infty} Ne_{\cT_N}(q)_p\leq \limsup_{n \to \infty} \frac {C_1}{c_0} 2^{n+1}(1+ \ve_{n+1}) e_{\cT_n^u}(q)\leq \frac {10C_0C_1}{c_0} \|\sqrt{\det q}\|_{L^\tau(\Omega)},
$$
which concludes the proof.
\sq

\section{The case of strictly convex functions}
\label{secCM4}
The goal of this section is to prove that the approximation
error in the greedy algorithm applied to a $C^2$ function $f$ satisfies the estimate
\iref{aniser} corresponding to an optimal triangulation.
Our main result is so far limited to 
the case where $f$ is strictly convex.

\begin{theorem}
\label{optitheo}
Let $f\in C^2(\overline \Omega)$ be such that 
$$
d^2f(x)\geq m I,\; \text{ for all } x\in \Omega
$$ 
for some arbitrary but fixed $m >0$ independent of $x$. Let $f_N$
be the approximant obtained by the greedy algorithm 
for the $L^p$ metric, using the $L^1$ decision function \iref{optil1}.
Then 
\be
\limsup_{N\to \infty} N \|f-f_N\|_{L^p}  \leq C \|\sqrt{\det(d^2f)}\|_{L^\tau},
\label{aniser1}
\ee
where $\frac 1 \tau=\frac 1 p +1$ and where $C$
is a constant independent of $p$, $f$ and $m$.  
\end{theorem}
Equation \iref{aniser1} can be rephrased as follows : there exists a sequence $\ve_N(f)$ such that $\ve_N(f)\to 0$ as $N \to \infty$ and 
$$
 \|f-f_N\|_{L^p}  \leq \(C \|\sqrt{\det(d^2f)}\|_{L^\tau}+\ve_N(f)\) N^{-1}. 
$$
Note also that since $\|\sqrt{\det(d^2f)}\|_{L^\tau}>0$, 
there exists $N_0(f)$ such that  $\|f-f_N\|_{L^p}  \leq 2C \|\sqrt{\det(d^2f)}\|_{L^\tau} N^{-1}$ for all $N\geq N_0(f)$.
It should be stressed hard that $N_0(f)$ can be arbitrarily large depending on the function $f$. Intuitively,
this means that when $f$ has very large hessian at certain point, it takes more iterations
for the algorithm to generate triangles with a good aspect ratio. 
The extension of this result to strictly concave functions is immediate
by a change of sign. Its extension to arbitrary $C^2$ functions
is so far incomplete, as it is explained in the end of 
the introduction.
The proof of Theorem \ref{optitheo} 
uses the fact that a strictly convex $C^2$ function is
{\it locally} close to a quadratic function with positive definite hessian,
which allows us to exploit the results obtained in \S \ref{secCM3} for 
these particular functions.

\subsection{A perturbation result}
\label{section:Perturb}

We consider a triangle $T$,  a function $f\in C^2(T)$, a convex quadratic function $q$ and $\mu>0$ such that on $T$
\be
\label{ineqfq}
d^2 q\leq d^2 f\leq (1+\mu) \ d^2 q.
\ee
It follows that
$
\det (d^2 q)\leq \det (d^2 f)\leq \det ((1+\mu) d^2 q ) = (1+\mu)^2 \det (d^2 q).
$
Since $\det (d^2 q) = 4 \det(\bq)$, we obtain
\be
\label{ineqdetfq}
2 \|\sqrt{\det \bq}\|_{L^\tau(T)}\leq \nsdf_{L^\tau(T)} \leq 2(1+\mu)  \|\sqrt{\det \bq}\|_{L^\tau(T)}.
\ee
The following Lemma shows how the local errors associated to $f$ and $q$ are close
\begin{prop}
\label{PropErC2}
The exists a constant $C_e>0$, depending only on the operator $\cA_T$ such that 
\be
\label{ineqefq}
(1-C_e\mu) e_T(q)_p \leq e_T(f)_p \leq (1+C_e\mu) e_T(q)_p.
\ee
\end{prop}

\proof
It follows from inequality \iref{ineqfq} that the functions $f-q$ and $(1+\mu)q -f$ are convex, hence
$$
\interp_T (f-q) - (f-q) \geq 0 \ \text{ and } \ \interp_T ((1+\mu)q-f) - ((1+\mu)q-f)\geq 0
$$
on the triangle $T$.
We therefore obtain
$$
0 \leq (\interp_T f-f) - (\interp_T q-q) \leq \mu (\interp_T q-q).
$$
There exists a constant $C_0>0$ depending only on $\cA_T$ such that for any $h\in C^0(T)$, 
$$
e_T(h)_p \leq C_0|T|^{\frac 1 p} \|h\|_{L^\infty(T)}.
$$
Furthermore according to Proposition \ref{propequiverrornondeg} there exists a constant $C_1>0$ depending only on $\cA_T$ such that 
$$
|T|^{1/p} \|q-\interp_T q\|_{L^\infty (T)} \leq C_1  e_T(q)_p.
$$
Hence 
\begin{eqnarray*}
|e_T(f)_p - e_T (q)_p| &\leq & e_T(f-q)_p\\
& = &e_T((\interp_T f-f) - (\interp_T q-q))_p \\
 & \leq & C_0 |T|^{\frac 1 p}\|(\interp_T f-f) - (\interp_T q-q)\|_{L^\infty (T)}\\
 & \leq & C_0|T|^{\frac 1 p} \|\mu (\interp_T q-q)\|_{L^\infty (T)}\\
 & \leq & C_0 C_1 \mu e_T(q)_p
\end{eqnarray*}
This concludes the proof of this Lemma, with $C_e=C_0C_1$.
\sq
\nl
Note that using Proposition {\rm \ref{propequiverrornondeg1}},
and assuming that $\mu\leq c_e:=\frac 1{2C_e}$, we have with
$\frac 1 \tau := 1+\frac 1 p$,
\be
\label{eqErrorFQ}
e_T(f)_p \sim e_T(q)_p \sim  \sigma_\bq(T)  \|\sqrt{|\det \bq|}\|_{L^\tau(T)}\sim \sigma_\bq(T) \nsdf_{L^\tau(T)},
\ee
with absolute constants in the equivalence.

We next study the behavior of the decision function $e\mapsto d_T(e,f)$.
For this purpose, we introduce
the following definition.

\begin{definition} \label{def:dg}
Let $T$ be a triangle with edges $a,b,c$. 
A $\delta$-near longest edge bisection with respect to 
the $\bq$-metric  is a 
bisection of any edge $e\in \{a,b,c\}$ such that 
$$
\bq(e)\geq (1-\delta)\max\{\bq(a), \bq(b), \bq(c)\}
$$
\end{definition}

\begin{proposition} 
Assume that $f$ and $q$ satisfy \iref{ineqfq}. Then, the bisection of $T$ prescribed by 
the decision function $e\mapsto d_T(e,f)$ is a $\mu$-near longest edge bisection
for the $\bq$-metric.
\end{proposition}

\proof
It follows directly from Equation \iref{Dint01} that for any edge $e$ of $T$,
$$
D_T(e,q)\leq D_T(e,f) \leq D_T(e,(1+\mu) q),
$$
hence we obtain using \iref{eqDq}
\be
\label{eqDfq}
\frac {|T|} {12} \bq(e) \leq D_T(e,f) \leq (1+\mu) \frac {|T|} {12} \bq(e).
\ee
Therefore the bisection of $T$ prescribed by 
the decision function $e\mapsto d_T(e,f)$ selects an $e$ such that
$$
(1+\mu) \bq(e) \geq \max\{\bq(a), \bq(b), \bq(c)\}.
$$
It is therefore a $\delta$-near longest edge bisection
for the $\bq$-metric with $\delta = \frac \mu {1+\mu}\leq \mu$ and therefore
also a $\mu$-near longest edge bisection. 
\sq

In the rest of this section,
we analyze the difference between a longest edge bisection in the $\bq$-metric and a $\delta$-near longest edge bisection. For that purpose we introduce a distance between triangles :
if $T_1, T_2$ are two triangles with edges $a_1,b_1,c_1$ and $a_2,b_2,c_2$ such that 
\be
\label{orderedEdges}
\bq(a_1)\geq \bq(b_1)\geq \bq(c_1) \ \text{ and } \ \bq(a_2)\geq \bq(b_2)\geq \bq(c_2),
\ee
we define
$$
\Delta_\bq(T_1,T_2) = \max\{|\bq(a_1)-\bq(a_2)|,|\bq(b_1)-\bq(b_2)|,|\bq(c_1)-\bq(c_2)|\}.
$$
Note that $\Delta_\bq$ is a distance up to rigid transformations.

\begin{lemma}
Let $T_1, T_2$ be two triangles, let $(R_1,U_1)$ and $(R_2,U_2)$ be the two pairs of children 
from the longest edge bisection of $T_1$ in the $\bq$-metric, and a $\delta$-near longest edge bisection of $T_2$ in the $\bq$-metric. 
Then, up to a permutation of the pair of triangles $(R_1,U_1)$,  
$$
\max \{\Delta_\bq(R_1,R_2),\Delta_\bq(U_1,U_2)\}\leq \frac 5 4 \Delta_\bq(T_1,T_2)+ \delta \bq(a_2).
$$
where $a_2$ is the longest edge of $T_2$ in the $\bq$-metric.
\end{lemma}

\proof
We assume that the edges of $T_1$ and $T_2$ are named and ordered as in \iref{orderedEdges}.
Up to a permutation,
$R_1$ and $U_1$ have edge vectors 
$b_1,a_1/2,(c_1-b_1)/2$ and $c_1,a_1/2,(b_1-c_1)/2$.
Two situations
might occur for the pair $(R_2,U_2)$:
\begin{itemize}
\item
$\bq(e) < (1-\delta) \bq(a_2)$ 
for $e=b_2$ and $c_2$. In such a case 
the triangle $T_2$ is bisected towards $a_2$, 
so that up to a permutation,
$R_2$ and $U_2$ have edge vectors 
$b_2,a_2/2,(c_2-b_2)/2$ and $c_2,a_2/2,(b_2-c_2)/2$.
Using that $\bq((c-b)/2)=\bq(c)/2+\bq(b)/2-\bq(a)/4$ when $a+b+c=0$, 
it clearly follows that
$$
\max \{\Delta_\bq(R_1,R_2),\Delta_\bq(U_1,U_2)\}\leq \frac 5 4 \Delta_\bq(T_1,T_2).
$$
\item
$\bq(e) \geq (1-\delta) \bq(a_2)$ 
for some $e=b_2$ or $c_2$. In such a case $T_2$ may be bisected say towards $b_2$, 
so that up to a permutation,
$R_2$ and $U_2$ have edge vectors 
$a_2,b_2/2,(c_2-a_2)/2$ and $c_2,b_2/2,(b_2-c_2)/2$.
But since $|\bq(b_2)-\bq(a_2)|\leq \delta \bq(a_2)$,
we obtain that
\be
\label{distPer}
\max \{\Delta_\bq(R_1,R_2),\Delta_\bq(U_1,U_2)\}\leq \frac 5 4 \Delta_\bq(T_1,T_2)+ \delta \bq(a_2).
\ee
\end{itemize}
\sq

We now introduce a perturbed version of the estimates 
describing the decay of the non-degeneracy measure
which were obtained in Proposition \ref{PropSigmaDec} and Corollary \ref{PropControlSigma2}.

\begin{prop}
\label{PropControlSigmaPer}
If $(T_i)_{i=1}^2$ are the two children obtained from a refinement of a triangle $T$ in which a $\delta$-near longest edge bisection in the $\bq$-metric is selected, then
\be
\label{eqSigmaDecPer}
\max\{\sigma_\bq(T_1),\sigma_\bq(T_2)\} \leq (1+4\delta) \sigma_\bq(T).
\ee
If $(T_i)_{i=1}^{8}$ are the eight children of a triangle $T$ 
obtained from three successive refinements in which a $\delta$-near longest edge bisection in the $\bq$-metric is selected, then
\begin{itemize}
       	\item for all $i$, $\sigma_\bq(T_i)\leq \sigma_\bq(T) (1+C_2\delta)$,
	\item there exists $i$ such that $\sigma_\bq(T_i)\leq 0.69 \, \sigma_\bq(T) (1+C_2\delta)$ or $\sigma_\bq(T_i)\leq M$,
	\end{itemize}
where $C_2 = \frac{61} 4$ and $M = 5(1+C_2 \delta)$. 
\end{prop}

\proof
We first prove \iref{eqSigmaDecPer}, and for that purpose we introduce the two children $T'_1, T'_2$  obtained by bisecting the longest edge of $T$ in the $\bq$-metric.
If follows from \iref{distPer} that, up to a permutation of the pair $(T'_1,T'_2)$, 
$$
\max \{\Delta_\bq(T_1,T'_1),\Delta_\bq(T_2,T'_2)\}\leq \delta \bq(a),
$$
where $a$ is the longest edge of $T$ in the $\bq$-metric.
Hence 
\be
\label{eqSigmaDelta}
|\sigma_\bq(T_i)-\sigma_\bq(T'_i)| \leq \frac {2 \Delta_\bq(T_i,T'_i)}{4|T_i|\sqrt{\det(\bq)}}
\leq 2 \delta \frac {\bq(a)}{4|T_i|\sqrt{\det(\bq)}}
\leq 4 \delta \sigma_\bq(T).
\ee
We know from Proposition \ref{PropSigmaDec} that $\max\{\sigma_\bq(T'_1), \sigma_\bq(T'_2)\}\leq \sigma_\bq(T)$. Combining this point with \iref{eqSigmaDelta} we conclude the proof of \iref{eqSigmaDecPer}.\\

We now turn to proof of the second part of the proposition and for that purpose we introduce the eight children $(T'_i)_{i=1}^8$ obtained from three successive refinements of $T$ in which the longest edge in the $\bq$-metric is selected. 
Iterating \iref{distPer}, we find that, up to a permutation of the triangles $(T'_i)_{i=1}^8$, one has 
$$
\max_{i=1,\cdots,8}\Delta_\bq(T_i,T'_i) \leq 
\left(1+\frac 5 4+\left(\frac 5 4\right)^2\right) \delta \bq(a)=\frac {61}{16} \delta \bq(a) = \frac {C_2 \delta} 4 q(a),
$$
where, again, $a$ is the longest edge of $T$ in the $\bq$-metric.
Repeating the argument \iref{eqSigmaDelta} we find that 
\be
\label{eqSigmaDelta2}
\max_{i=1,\cdots,8} |\sigma_\bq(T_i)-\sigma_\bq(T'_i)|\leq C_2\delta \sigma_\bq(T).
\ee
We know from Corollary \ref{PropControlSigma2}
that $\sigma_\bq(T_i')\leq \sigma_\bq(T)$ for all $i$
and that there exists $i$ such that either
$\sigma_\bq(T_i')\leq 0.69 \, \sigma_\bq(T)$ or $\sigma_\bq(T_i')\leq 5$.
Combining this point with \iref{eqSigmaDelta2} we conclude the proof of the proposition.
\sq

\subsection{Local optimality}
\label{section:OptQuad}

Our next step towards the proof of Theorem \ref{optitheo} is
to show that the triangulation produced by the greedy algorithm
is locally optimal in the following sense: 
if the refinement procedure for the function $f$ 
produces a triangle $T\in\cD$ on which $f$ is close
enough to a quadratic function $q$, then the triangles 
which are generated from the refinement of $T$ tend to adopt
an optimal aspect ratio in the $\bq$-metric, 
and a local version of the optimal
estimate \iref{aniser} holds on $T$. 

We first prove that most triangles adopt an optimal aspect ratio
as we iterate the refinement procedure. Our goal is thus to 
obtain a result similar to Theorem \ref{CorolStabilise} 
which was restricted to quadratic functions. However, 
due to the perturbations by $C_2\mu$
that appear in Proposition \ref{PropControlSigmaPer},
the formulation will be slightly different, yet sufficient for
our purposes: we shall prove that the measure of non-degeneracy
becomes bounded by an absolute constant in an average sense,
as we iterate the refinement procedure.

As in the previous section, we
assume that $f$ and $q$ satisfy \iref{ineqfq}. For any $T$, we define $\cT_n^u(T)$ the triangulation of $T$ 
which is built by iteratively applying the refinement procedure for the function $f$
to {\it all generated triangles} up to $3n$ generation levels.
Note that 
$$
\#(\cT_n^u(T))=2^{3n}\;\; {\rm and}\;\; |T'|=2^{-3n}|T|,\;\; T'\in \cT_n^u(T).
$$
For $r>0$, we define
the average $r$-th power of the measure of non-degeneracy of the $2^{3n}$
triangles obtained from $T$ after $3n$ iterations by
$$
\o {\sigma^r_\bq(n)} = \frac 1 {2^{3n}} \sum_{T'\in \cT_n^u(T) } \sigma^r_\bq(T').
$$
We also define
$$
\gamma(r,\mu):= \frac 1 8 \(0.69(1+C_2\mu)\)^r+\frac 7 8 (1+C_2\mu)^r,
$$
where $C_2$ is the constant in Proposition \ref{PropControlSigmaPer}. 
Note that for any $r>0$, the function $\gamma(r,\cdot)$ is continuous and increasing, and that $0<\gamma(r,0)<1$. Hence  for any $r>0$, there
exists $\mu(r)>0$
and $0<\gamma(r)<1$ such that $\gamma(r,\mu)\leq \gamma(r)$, if $0<\mu<\mu(r)$. 

\begin{prop} \label{CorolGamma}
Assume that $f$ and $q$ satisfy \iref{ineqfq} with $0<\mu\leq \mu(r)$. We then have
$$
\o {\sigma^r_\bq(n)} \leq \sigma^r_\bq(T) \gamma(r)^n+\frac {M^r} {8(1-\gamma(r))},
$$
where $M$ is the constant in Proposition \ref{PropControlSigmaPer}.
Therefore 
$$
\o {\sigma^r_\bq(n)}\leq C_3:= 1+\frac {M^r} {8(1-\gamma(r))},
$$ 
if 
$2^{3n}  \geq 8 \sigma_\bq(T)^{\lambda}$ with $\lambda := \frac{3 r \ln 2}{-\ln \gamma(r)}$.
\end{prop}

\proof
Let us use the notations $u= 0.69(1+C_2\mu)$ and $v=(1+C_2\mu)$.
According to Proposition \ref{PropControlSigmaPer}, we
have 
$$
\o {\sigma^r_\bq(n)}\leq  \mE(\sigma_n^r),
$$
where $\mE$ is the expectation operator and $\sigma_n$ is the Markov chain with value in $[1,+\infty[$ defined by
\begin{itemize}
	\item $\sigma_{n+1} = \max\{\sigma_n u, M\}$ with probability $\alpha:=\frac 1 8$,
	\item $\sigma_{n+1} = \sigma_n v$ with probability $\beta:=\frac 7 8$,
	\item $\sigma_0:=\sigma_\bq(T_0)$ with probability $1$.
\end{itemize}
Denoting by $\mu_n$ the probability distribution of $\sigma_n$, we have
\begin{eqnarray*}
\mE(\sigma_{n+1}^r) &=& \int_1^\infty \sigma^r d\mu_{n+1}(\sigma)  \\
&= &\int_1^\infty \left( \alpha (\max\{u\sigma,M\})^r+\beta (v\sigma)^r \right)d\mu_n(\sigma)\\
&=& \alpha M^r \int_1^{M/u} d\mu_n(\sigma) + \alpha u^r \int_{M/u}^\infty \sigma^r d\mu_n(\sigma)
+\beta v^r \int_1^{+\infty}\sigma^r d\mu_n(\sigma)\\
&\leq &\alpha M^r+(\alpha u^r + \beta v^r) \mE(\sigma_n^r)\\
&\leq &
\alpha M^r+\gamma(r) \mE(\sigma_n^r)
\end{eqnarray*}By iteration, it follows that
$$
\mE(\sigma_n^r)\leq \mE(\sigma_0^r)\gamma(r)^n+\frac{ \alpha M^r}{1-\gamma(r)},
$$ 
which gives the result.
\sq

Our next goal is to show that the greedy algorithm initialized from $T$
generates a triangulation which is a refinement 
of $\cT_n^u(T)$ and therefore more accurate, yet with a similar
amount of triangles. To this end, we apply the greedy 
algorithm with root $T$ and stopping criterion given
by the local error 
$$
\eta:=\min_{T'\in \cT_n^u(T)} e_{T'}(f)_p.
$$
Therefore $T'$ is splitted if and only if $e_{T'}(f)_p> \eta$.
We denote by $\cT_{N}(T)$ the resulting triangulation
where $N$ is its cardinality. From the definition of the 
stopping criterion, it is clear that $\cT_{N}(T)$ is a refinement 
of $\cT_n^u(T)$. 

\begin{prop} \label{PropQuadLP}
Assume that $f$ and $q$ satisfy \iref{ineqfq} with $\mu \leq \frac {1}{8}$,
and define
$r_0:=\frac{\ln 2}{\ln 4-\ln 3}>0$.
We then have
$$
N\leq C_42^{3n}\o {\sigma^{r_0}_\bq(n)},
$$
where $C_4$ is an absolute constant.
Assuming in addition that $\mu\leq \mu(r_0)$ as in Proposition \ref{CorolGamma},
we obtain that
$$
N\leq C_52^{3n},
$$
if $2^{3n} \geq 8 \sigma_\bq(T)^{\lambda}$ with $\lambda := \frac{3 r_0\ln 2}{-\ln \gamma(r_0)}$,
and where $C_5=C_3C_4$.
\end{prop}

\proof
Let $T_1$ be a triangle in $\cT_{n}^u(T)$ and $T_2$ a triangle in $\cT_{N}(T)$
such that $T_2\subset T_1$. We shall give a bound on the number of splits
$k$ which were applied between $T_1$ and $T_2$, i.e. such that
$|T_2|=2^{-k}|T_1|$. We first remark that according to Proposition 
\ref{propequiverrornondeg1} and \iref{eqErrorFQ}, we have
$$
\eta\geq c\min_{T'\in \cT_n^u(T)}|T'|^{1+\frac 1 p} \sigma_\bq(T')  \sqrt{\det \bq}
\geq c|T_1|^{1+\frac 1 p} \sqrt{\det \bq},
$$
where $c$ is an absolute constant.
On the other hand, using both
Proposition \ref{PropErC2} and Proposition \ref{PropControlSigmaPer}, we obtain
$$
\begin{array}{ll}
e_{T_2}(f)_\bq 
&\leq C|T_2|^{1+\frac 1 p} \sigma_\bq(T_2)  \sqrt{\det \bq} \\
&=C |T_1|^{1+\frac 1 p} 2^{-k(1+\frac 1 p)}\sigma_\bq(T_2)  \sqrt{\det \bq} \\
&\leq C|T_1|^{1+\frac 1 p} \sigma_\bq(T_1)\(2^{-(1+\frac 1 p)}(1 + 4\mu)\)^{k} \sqrt{\det \bq}.\\
&\leq \frac C c \sigma_\bq(T_1)\(\frac {1+4\mu} 2\)^{k} \eta \\
&\leq \frac C c \sigma_\bq(T_1)(\frac 3 4)^{k} \eta,
\end{array}
$$
where $C$ is an absolute constant.
Therefore we see that $k$ is at most the smallest integer 
such that $\frac C c \sigma_\bq(T_1)(\frac 3 4)^{k} \leq 1$.
It follows that the total number $n(T_1)$ of triangles $T_2\in \cT_{N}(T)$
which are contained in $T_1$
is bounded by
$$
n(T_1)\leq 2^k \leq 2 \left(\frac C c \sigma_\bq(T_1)\right)^{r_0},
$$
and therefore
$$
N=\sum_{T_1\in \cT^u_n(T)} n(T_1) \leq 2 \left(\frac C c\right)^{r_0} \sum_{T_1\in \cT^u_n(T)} 
\sigma_\bq(T_1)^{r_0} = C_42^{3n}\o {\sigma^{r_0}_\bq(n)},
$$
with $C_4= 2 \left(\frac C c\right)^{r_0}$. The fact that 
$N\leq C_52^{3n}$ when $2^{3n}  \geq 8 \sigma_\bq(T)^{\lambda}$ 
with $\lambda := \frac{3 r_0 \ln 2}{-\ln \gamma(r_0)}$ is an immediate
consequence of Proposition {\rm \ref{CorolGamma}}.
\sq

\subsection{Optimal convergence estimates}

Our last step towards the proof of Theorem
\ref{optitheo} consists in deriving local error estimates for the greedy algorithm.
For $\eta>0$, we 
denote by $f_{\eta}$ the approximant to $f$
obtained by the greedy algorithm with
stopping criterion given
by the local error $\eta$ : a triangle
$T$ is splitted if and only if $e_{T}(f)_p > \eta$.
The resulting triangulation is denoted by
$$
\cT_\eta=\cT_N, \text{ with } N=N(\eta)=\#(\cT_\eta).
$$
For this $N$, we thus have $f_\eta=f_N$.
For a given $T$ generated by the refinement procedure and
such that $\eta \leq e_T(f)_p$, we also define
$$
\cT_{\eta}(T)=\{T'\subset T\; ; \; T' \in \cT_\eta\}
$$
the triangles in $\cT_\eta$ which are contained in $T$
and
$$
N(T,\eta)=\#(\cT_\eta(T)).
$$
Our next result provides with estimates of the local
error $\|f-f_\eta\|_{L^p(T)}$ and of $N(T,\eta)$ in terms of 
$\eta$, provided that $\mu$ is small enough.

\begin{theorem} 
\label{localerror}
Assume that $f$ and $q$ satisfy \iref{ineqfq} with $\mu\leq c_2:= \min\{ \frac 1 8,\mu(r_0)\}$,
and that $\eta \leq \eta_0$, where
$$
\eta_0=\eta_0(T):=\left(\frac {|T|}{\sigma_\bq(T)^{\lambda}}\right)^{\frac 1 \tau}\sqrt{\det \bq},
$$
with $\lambda := \frac{3 r_0\ln 2}{-\ln \gamma(r_0)}$,
and $\frac 1 \tau=\frac 1 p+1$. Then
\be
\|f-f_\eta\|_{L^p(T)}\leq \eta N(T,\eta)^{\frac 1 p},
\label{localerreta}
\ee
and 
\be
N(T,\eta) \leq C_6\eta^{-\tau} \|\sqrt{\det(d^2f)}\|_{L^\tau(T)}^\tau,
\label{estimNeta}
\ee
where $C_6$ is an absolute constant.
\end{theorem}

\proof
The first estimate is trivial since
$$
\|f-f_\eta\|_{L^p(T)}=\(\sum_{T'\in \cT_{\eta}(T)} e_{T'}(f)_p^p\)^{\frac 1 p}
\leq \(\sum_{T'\in \cT_{\eta}(T)} \eta^p\)^{\frac 1 p} =\eta N(T,\eta)^{\frac 1 p}.
$$
In the case $p=\infty$, we trivially have
$$
\|f-f_\eta\|_{L^\infty(T)}\leq \eta.
$$
For the second estimate, we define $n_0=n_0(T)$ the smallest
positive integer such that 
$2^{3n_0(T)}\geq  8\sigma_\bq(T)^{\lambda}$ with 
$\lambda := \frac{3 r_0\ln 2}{-\ln \gamma(r_0)}$.
For any fixed $n\geq n_0$, we define
$$
\eta_n:=\min_{T'\in \cT_n^u(T)} e_{T'}(f)_p.
$$
We know from Proposition \ref{PropQuadLP} that with the choice $\eta=\eta_n$
\be
N(T,\eta_n)\leq C_5 2^{3n}.
\label{netan}
\ee
On the other hand, we know from Proposition \ref{CorolGamma}, 
that $\o{\sigma_\bq^{r_0}(n)}\leq C_3$,
from which it follows that 
$$
\min_{T'\in\cT_n^u(T)} \sigma_\bq(T') \leq C_3^{\frac 1 {r_0}}.
$$
According to Proposition 
\ref{PropErC2}, we also have 
$$
\eta_n
\leq C\min_{T'\in \cT_n^u(T)}|T'|^{1+\frac 1 p} \sigma_\bq(T')  \sqrt{\det \bq}
\leq C_3^{\frac 1 {r_0}} C  \(\frac {|T|}{2^{3n}}\)^{\frac 1 \tau}\sqrt{\det \bq},
$$
where $C$ is an absolute constant, which also reads
$$
2^{3n}\leq C_3^{\frac \tau {r_0}}C^\tau  \eta_n^{-\tau} |T|\sqrt{\det \bq}^{\,\tau}.
$$
Combining this with \iref{netan}, we have obtained the estimate
$$
N(T,\eta_n)\leq C_5C_3^{\frac \tau {r_0}}C^\tau  \eta_n^{-\tau} |T|\sqrt{\det \bq}^{\,\tau},
$$
which by Proposition \ref{PropErC2} is equivalent to \iref{estimNeta} 
with $\eta=\eta_n$. In order to obtain \iref{estimNeta} for all arbitrary
values of $\eta$, we write that 
$\eta_{n+1}< \eta \leq \eta_n$ for some $n\geq n_0$, then
$$
\begin{array}{ll}
N(T,\eta) &\leq N(T,\eta_{n+1}) \\
& \leq C_5 2^{3(n+1)}\\
& \leq  8C_5C_3^{\frac \tau {r_0}}C^\tau  \eta_n^{-\tau} |T|\sqrt{\det \bq}^{\,\tau}\\
&\leq 8C_5C_3^{\frac \tau {r_0}}C^\tau  \eta^{-\tau} |T|\sqrt{\det \bq}^{\,\tau},
\end{array}
$$
which by Proposition \ref{PropErC2} is equivalent to \iref{estimNeta}.
In the case where $\eta \geq \eta_{n_0}$, we simply write
$$
\begin{array}{ll}
N(T,\eta) & \leq  N(T,\eta_{n_0}) \\
& \leq C_5 2^{3n_0}\\
& \leq 64C_5\sigma_\bq(T)^{\lambda} \\
&=64C_5\eta_0^{-\tau} |T|\sqrt{\det \bq}^{\,\tau} \\
& \leq 64C_5\eta^{-\tau} |T|\sqrt{\det \bq}^{\,\tau},
\end{array}
$$
and we conclude in the same way.
\sq
\nl
We remark that combining the estimates
\iref{localerreta} and \iref{estimNeta} in the above theorem
yields the optimal local convergence estimate
$$
\|f-f_\eta\|_{L^p(T)}\leq C_6^{\frac 1 \tau} \|\sqrt{\det(d^2f)}\|_{L^\tau(T)}N(T,\eta)^{-1}.
$$
In order to obtain the global estimate of Theorem \ref{optitheo}, we 
need to be ensured that after sufficiently many steps of the greedy algorithm, 
the target $f$ can be well approximated by quadratic function $q=q(T)$ 
on each triangle $T$, so that our local results will apply on such triangles.
This is ensured due to the following key result.

\begin{prop} \label{PropDiamZero}
Let $f$ be a $C^2$ function such that 
$d^2 f(x)\geq m I $ 
for some arbitrary but fixed $m >0$ independent of $x$. 
Let $\cT_N$ be the triangulation generated by the greedy algorithm
applied to $f$ using the $L^1$ decision function given by \iref{optil1}.
Then
$$
\lim_{N\to +\infty} \max_{T\in\cT_N}\diam(T)=0,
$$
i.e. the diameter of all triangles tends to $0$.
\end{prop}
\proof
Let $T$ be a triangle with an angle $\theta$ at a vertex $z_0$. The other vertices of $T$ 
can be written as $z_1=z_0+\alpha u$ and $z_2=z_0+ \beta v$ where $\alpha, \beta\in \R_+$ and $u,v\in \R^2$ are unitary. We assume that $\alpha u$ is the longest edge of $T$, hence $\theta\leq \cPi/2$.
Observe that 
$$
\rho(T) := \frac {h_T^2} {|T|}=\frac{\alpha^2} {\frac 1 2\alpha \beta \sin \theta} = \frac {2 \alpha}{\beta \sin \theta},
$$
and
$$ 
|u-v| = 2 \sin\left( \frac \theta 2\right)=\frac {\sin \theta}{\cos(\frac \theta 2)}
=\frac {2\alpha}{\beta\rho(T)\cos(\frac \theta 2)}.
$$
Since $\frac {\sqrt 2} 2 \leq \cos\left(\frac \theta 2\right) $ we thus obtain 
$$
|u-v| \leq  \frac {2\sqrt 2\alpha}{\beta \rho(T)} \leq \frac {3\alpha}{\beta \rho(T)} .
$$
We now set $M := \|d^2 f\|_{L^\infty(\Omega)}$ and for all $\delta>0$ let 
$$
\omega(\delta) := \sup_{z,z'\in \Omega, \|z-z'\|\leq \delta} \|d^2f(z) - d^2f(z')\|.
$$
For $t\in\RR$, we define 
$$
H_t^u:= d^2 f_{z_0 + t u} \ \text{ and } \ H_t^v:= d^2 f_{z_0 + t v}.
$$
and notice that 
$
\|H_t^u - H_t^v\| \leq \omega(t|u-v|).
$
Hence, if $0\leq t\leq \beta$, we have
$$
\|H_t^u-H_t^v\|\leq \omega \left(\frac {3\alpha } {\rho(T)}\right).
$$ 
Furthermore, for all $t$ we have
$$
\begin{array}{ll}
|\<H_t^u \ u, u\> - \<H_t^u \ v, v\>| & 
=|\<H_t^u \ u, u\> - \<H_t^u \ u-(u-v), u-(u-v)\>| \\
& =|2\<H_t^u \ u, u-v\> -\<H_t^u \ (u-v), u-v\>|\\
&\leq 2M |u| |u-v|+ M |u-v|^2 \\
& \leq \frac M {\beta^2} \left( 2\frac {3\alpha \beta}{\rho(T)} + \left(\frac {3\alpha}{\rho(T)}\right)^2 \right).
\end{array}
$$
Applying the identity \iref{Dint01} to the edges $e=\alpha u$ and $\beta v$,
and using a change of variable, we can write
$$
D_T(\alpha u ,f)=\int_\R \min\{t,\alpha-t\}_+ \<H_t^u\, u,u\> dt
\;\; {\rm and}\;\; D_T(\beta v ,f)=\int_\R \min\{t,\beta-t\}_+ \<H_t^v\, v,v\> dt
$$
where we have used the notation $r_+ := \max\{r,0\}$.
Hence, noticing that 
$$
\int_\R \min\{t,\lambda-t\}_+ dt=\int_0^\lambda \min\{t,\lambda-t\}dt = \frac{(\lambda_+)^2} 4,
$$
and using the previous estimates
we obtain
\begin{eqnarray*}
D_T(\alpha u ,f) - D_T(\beta v,f) &=& \int_\R \left(\min\{t,\alpha-t\}_+ \<H_t^u\, u,u\> - \min\{t, \beta-t\}_+ \<H_t^v \, v,v\>\right) dt\\
&= &  \int_\R (\min\{t,\alpha-t\}_+ -  \min\{t, \beta-t\}_+) \<H_t^u \, u,u\> dt\\
& & -\int_\R  \min\{t, \beta-t\}_+ (\<H_t^v \, v,v\>- \<H_t^u\, u,u\>)dt \\
& \geq &m\int_\R (\min\{t,\alpha-t\}_+ -  \min\{t, \beta-t\}_+)dt \\
& & -\int_0^\beta \min\{t, \beta-t\} (|\<H_t^u \ u, u\> - \<H_t^u \ v, v\>| \\
& &\quad +| \<(H_t^u-H_t^v)v,v\>|)dt \\
&\geq &  m \frac {\alpha^2 -\beta^2} 4  -
\frac M 4 \left( 2\frac {3\alpha\beta}{\rho(T)} + \left(\frac {3\alpha}{\rho(T)}\right)^2 \right)
- \frac{\beta^2} 4 \omega\left(\frac { 3\alpha}{ \rho(T)}\right)\\
& \geq & m \frac {\alpha^2 -\beta^2} 4 - \frac {\alpha^2} 4 \left(
\frac 6 {\rho(T)}+\frac 9 {\rho(T)^2}+ \omega\left(\frac { 3\alpha}{ \rho(T)}\right)\right),
\end{eqnarray*}
where we have used the fact that $\alpha>\beta$ in the last line.
We can therefore write
\be
D_T(\alpha u ,f) - D_T(\beta v,f) \geq \frac m 4 \(\alpha^2\(1-\ti \omega\(\frac 1{\rho(T)}\)\)-\beta^2\),
\label{eqDT}
\ee
where we have set
$$ 
\tilde \omega(\delta) := \frac 1 m \( 6 \delta+9 \delta^2+\omega( 3\diam(\Omega) \delta)\).
$$ 
The inequality \iref{eqDT} shows that $d_T(\cdot,f)$ prescribes a
$\ti \omega\(\frac 1 {\rho(T)}\)$-near longest edge bisection in the euclidean metric for any triangle $T$. 
Indeed if the smaller edge $\beta v$ was selected, we would necessarily have
$$
|\beta v|^2=\beta^2 \geq \(1-\ti \omega\(\frac 1 {\rho(T)}\)\) \alpha^2 = \(1-\ti \omega\(\frac 1 {\rho(T)}\)\)
|\alpha u|^2.
$$
Notice that $\tilde \omega(\delta)\to 0$ as $\delta \to 0$.

Since $f$ is strictly convex, there does not exists any triangle $T\subset \Omega$ such that $e_T(f)_p = 0$. 
Let us assume for contradiction that the diameter of the triangles generated by the greedy
algorithm does not tend to zero. Then there 
exists a sequence $(T_i)_{i\geq 0}$ of triangles such that 
$T_{i+1}$ is one of the children of $T_i$, and $h_{T_i}\to d>0$ 
as $i \to \infty$, where $h_T$ denotes the diameter of a triangle $T$. Since $|T_i|\to 0$, this also implies that $\rho(T_i)\to +\infty$
as $i\to \infty$. We can therefore choose $i$ large enough such that
$h_{T_i}^2<\frac 4 3 d^2$ and 
$C_2\tilde \omega\(\frac 1 {\rho(T_j)}\)\leq \frac 1 2$ for all $j\geq i$,
where $C_2$ is the constant in Proposition \ref{PropControlSigmaPer}.
According to this Proposition, we have
$$
\sigma(T_{i+3}) \leq \frac 3 2 \sigma(T_i),
$$
where $\sigma$ stands for $\sigma_{\bq}$ in the euclidean case $\bq=x^2+y^2$.
On the other hand, we have for any triangle $T$,
$
\frac {h_T^2}{8|T|} \leq \sigma(T)\leq \frac {h_T^2}{2|T|},
$
from which it follows that 
$$
h_{T_{i+3}}^2 \leq 4 \frac {|T_{i+3}|\sigma(T_{i+3})}  {|T_{i}|\sigma(T_{i})}
h_{T_i}^2 \leq \frac 3 4
h_{T_i}^2.
$$
Therefore, $h_{T_{i+3}}<d$ 
which is a contradiction. This concludes the proof of
Proposition \ref{PropDiamZero}.
\sq

\noindent
{\bf Proof of Theorem \ref{optitheo}}
Since $f\in C^2$, an immediate consequence of Proposition
\ref{PropDiamZero} is that for
all $\mu>0$, there exists 
$$
N_1:=N_1(f,\mu),
$$
such that for all $T\in \cT_{N_1}$, there exists a quadratic function $q_T$
such that  
$$
d^2 q_T \leq d^2 f\leq (1+\mu) \, d^2 q_T.
$$
Therefore our local results apply on all $T\in \cT_{N_1}$. Specifically,
we choose
$$
N_1:=N_1(f,c_2),
$$ 
with $c_2$ the constant in Theorem \ref{localerror}. We then take 
$$
\eta \leq \eta_0:=\min_{T\in\cT_{N_1}}\left\{e_T(f)_p,\(\frac {|T|}{ \sigma_{\bq_T}(T)^{\lambda}}\)^{\frac 1 \tau}\sqrt{\det \bq_T}\right\}.
$$
We use the notations
$$
f_{\eta}=f_N,\;\; \cT_\eta=\cT_N,\;\, N=N(\eta)=\#(\cT_\eta)=\#(\cT_N),
$$
for the approximants and triangulation
obtained by the greedy algorithm with
stopping criterion given
by the local error $\eta$. Note that $\cT_{\eta}$ is a refinement of $\cT_{N_1}$,
since $\eta\leq \min_{T\in\cT_{N_1}}e_T(f)_p$, and therefore $N\geq N_1$.
We obviously have
$$
\|f-f_N\|_{L^p}\leq \eta N^{\frac 1 p}.
$$
Using Theorem \ref{localerror}, we also have
$$
N=\sum_{T\in \cT_{N_1}}N(T,\eta)\leq C_6\eta^{-\tau} \|\sqrt{\det(d^2f)}\|_{L^\tau(\Omega)}^\tau,
$$
and therefore 
$$
\|f-f_N\| \leq  C_6^{\frac 1 \tau} \|\sqrt{\det(d^2f)}\|_{L^\tau(\Omega)}N^{-1},
$$
which is the claimed estimate. Since we have assumed 
$\eta\leq \eta_0$, this estimate holds for 
$$
N> N_0,
$$
where $N_0$ is largest value of $N$ such that $e_T(f)_p\geq \eta_0$
for at least one $T\in\cT_N$.
\sq

\begin{remark}
In Chapter \ref{chapCDHM} a modification of the algorithm is proposed so that 
its convergence in the $L^p$ norm is ensured
for {\it any} function $f\in L^p(\Omega)$ (or $f\in C(\Omega)$ when $p=\infty$). 
However this modification is not needed 
in the proof of Theorem \ref{optitheo}, due to the assumption
that $f$ is convex. 
\end{remark}

\chapter{Variants of the greedy bisection algorithm}
\minitoc
\label{chapVariations}
\section{Introduction}

We study in this chapter several variants of the greedy algorithm that was
discussed in Chapter \ref{chapCDHM} and Chapter \ref{chapBisecOpt}. 

One of the key features of this algorithm is the decision function which governs the creation of anisotropy by selecting among the available directions of refinement of a triangle. It was proved in 
Chapter \ref{chapBisecOpt} that for piecewise linear approximation, 
the $L^1$ based decision function \iref{eqDecisionL1}
leads to optimally adapted triangles. We consider in \S \ref{secAltDec} two alternative
decision functions, based on the $L^2$ approximation error or the $L^\infty$ interpolation error respectively,
both in the context of piecewise linear approximation. 
The study of these decision functions is motivated by the following reasons. 
The $L^2$ based decision function can be computed at 
a significantly smaller numerical cost than the $L^1$ based or $L^\infty$ based decision functions, and is therefore the most suited method for numerical applications. 
The $L^\infty$ based decision function is computationally
more costly, but it leads to more general convergence results.

We then focus our attention in \S \ref{secAltBi} to the behaviour of the
greedy algorithm for piecewise linear approximation 
when applied to cartoon functions. 
Unfortunately, it turns out that the algorithm fails to achieve the
best convergence rate ($N^{-1}$ in the $L^2$ norm) expected for such functions. 
This is inherently due to the fact that the geometry of the bisections used in the
algorithm is too limited, regardless of the decision function which is being used.
We thus consider some modifications of the greedy algorithm 
which partially solve this problem using alternative bisection choices. 
Instead of bisecting a triangle from a vertex to the midpoint of the opposite edge, 
we offer different possibilities which lead to better directional selectivity. 
In turn we obtain optimal convergence rates for simple cartoon functions of the
form $f=\chi_P$ where $P$ is any half plane. 
The behavior of the greedy algorithm on general cartoon functions remains an open question.

Finally, we consider in \S \ref{secGreedyRect} a variant of the greedy algorithm which is in some sense much simpler. 
It is based on piecewise constant approximation, instead of piecewise linear, and it produces partitions of the original domain into rectangles aligned with the coordinate axes, instead of triangles of arbitrary direction. 
Thanks these simplifications we were able to analyze this algorithm in a rather general setting, and we establish in Theorem \ref{thRect} a convergence estimate which applies to any $C^1$ function
and which is in accordance with the convergence estimate established in Chapter 1 for optimized
rectangular partitions.

\section{Alternative decision functions}
\label{secAltDec}

Let us briefly recall the two steps of the greedy refinement algorithm studied in 
Chapter \ref{chapCDHM} and Chapter \ref{chapBisecOpt}:
\begin{enumerate}[a)]
\item A triangle $T\in \cT$ which maximizes the approximation error is selected.
$$
T = \underset {T'\in \cT} \argmax \ e_{T'}(f)_p.
$$
\item 
For each edge $e\in \{a,b,c\}$ of the triangle $T$, the bisection of $T$ from the midpoint of $e$ to the opposite vertex defines two children $T_e^1$ and $T_e^2$.
The edge $e$ which minimizes a given decision function $d_T(e,f)$ is selected
$$
e := \underset {e'\in \{a,b,c\}} \argmin \ d_T(e',f).
$$
and we define the new triangulation $\cT'$ by 
$$
\cT' := \cT - \{T\} + \{T_e^1, T_e^1\}.
$$
\end{enumerate}

Theorem \ref{optitheo} states that if the function $f$ to be approximated is $C^2$ and strictly convex, and if the decision function is given by 
\be
\label{eqDecisionL1}
d_T(e,f) := \|f-\interp_{T_e^1} f \|_{L^1(T_e^1)}+  \|f-\interp_{T_e^2} f \|_{L^1(T_e^2)},
\ee
then the sequence of triangulations generated by the greedy algorithm satisfies the optimal error estimate
\be
\label{eqOptiConv}
\limsup_{N\to \infty} N e_{\cT_N}(f)_p \leq C \left \| \sqrt{|\det(d^2 f)|} \right\|_{L^\tau(\Omega)}.
\ee

We now want to investigate two alternate choices for the decision function $d_T(e,f)$. 
The first of these choices is based on the $L^2$ approximation error 
\be
\label{eqDecisionL2}
d_T(e,f) := \|f-P_{T_e^1} f \|_{L^2(T_e^1)}^2+  \|f-P_{T_e^2} f \|_{L^2(T_e^2)}^2,
\ee
where, for any triangle $T$, we denote by $P_T$ the $L^2(T)$ orthogonal projection onto the space $\P_1$ of affine functions. The second choice is based on the $L^\infty$ interpolation error,
\be
\label{eqDecisionLInf}
d_T(e,f) := \|f-\interp_{T_e^1} f \|_{L^\infty(T_e^1)}+  \|f-\interp_{T_e^2} f \|_{L^\infty(T_e^2)}.
\ee

Our motivation for studying the $L^2$ based decision function
is computational. Indeed, from the point of view of computer run time, the biggest cost of the greedy algorithm comes from the evaluation of the decision function. As exposed below, a trick allows to evaluate the (discretized) $L^2$ based decision function \iref{eqDecisionL2} much faster than the $L^1$ based or $L^\infty$ based decision functions \iref{eqDecisionL1} and \iref{eqDecisionLInf}. 
The greedy algorithm has been implemented in C++ by Lihua Yang.
For an image of realistic resolution, say $512 \times 512$, the (discretized) greedy algorithm takes only 
two seconds on a standard laptop computer to generate $5000$ triangles using the $L^2$ based decision function if it is efficiently evaluated. In contrast it may take up to a minute to generate the same number of triangles based on a brute force evaluation of the $L^2$ based decision function, or the $L^1$ based or $L^\infty$ based decision functions. 
As a result the $L^2$ based decision function \iref{eqDecisionL2} 
is the preferred one in numerical experiments.
We shall prove that the greedy algorithm based on the $L^2$ decision function 
generates a sequence of triangulations satisfying the optimal estimate \iref{eqOptiConv} when it is applied to any quadratic function $f=q$ such that the quadratic form $\bq$ is non-degenerate.
In particular $q$ may be either stricly convex, concave or of hyperbolic type. 
The latter case, which corresponds to a mixed signature $(1,1)$ of the quadratic form $\bq$, is not treated in the  study in the previous Chapter \ref{chapBisecOpt} of the $L^1$ based decision function \iref{eqDecisionL1}.\\

The decision function based on the $L^\infty$ error  \iref{eqDecisionLInf}
is comparable to the $L^1$ based decision function in terms of computational cost,
but it leads to the most complete theoretical results. Indeed 
we shall prove that the optimal convergence estimate \iref{eqOptiConv} holds for any quadratic function $f=q$ such that the quadratic form $\bq$ is non degenerate, as well 
as for functions which are $C^2$ and stricly convex. The latter case was treated in 
Theorem \ref{optitheo} of Chapter \ref{chapBisecOpt} for the $L^1$ based decision function, 
but is not established for the $L^2$ based decision function.\\

Before turning to the detailed study of the different decision functions, 
we expose as announced why the (discretized) $L^2$ based decision function \iref{eqDecisionL2} is much less expensive to compute than the $L^1$ based or $L^\infty$ based decision functions, in terms of computer run time.

Let $T$ be a triangle, and let $M_T$ be the following $3\times 3$ symmetric positive definite matrix, the Grammian matrix of the basis $(1,x,y)$ of the subspace $\P_1$ of $L^2(T)$,
$$
M_T := \int_T \left( 
\begin{array}{ccc}
1 & x & y\\
x & x^2 & xy\\
y & xy & y^2
\end{array}
\right)
dx\, dy
= 
\int_T
 (1,x,y)(1,x,y)^\trans
dx\,dy.
$$
Let also $V_T(f)\in \R^3$ be the vector of the first order moments of the function $f$ on $T$, 
$$
V_T(f) := \int_T (1,x,y)^\trans f(x,y) \, dx\, dy.
$$
The $L^2(T)$ orthogonal projection of $f$ onto $\P_1$ has the expression 
$$
P_T(f) = (1,x,y)M_T^{-1} V_T.
$$
Therefore, denoting by $\<\cdot, \cdot\>$ the $L^2(T)$ scalar product,
\begin{eqnarray*}
\|f-P_T f\|_{L^2(T)}^2 &=& \|f\|_{L^2(T)}^2 - 2 \<f, P_T(f)\> + \|P_T(f)\|_{L^2(T)}^2\\
&=&  \|f\|_{L^2(T)}^2 - \|P_T(f)\|_{L^2(T)}^2\\
&=& \int_T f(x,y)^2 dx\, dy - V_T(f)^\trans M_T^{-1} V_T(f).
\end{eqnarray*}
Hence the $L^2(T)$ approximation error of $f$ on $T$ has an expression in terms of the integrals on $T$ of the functions $1,\, x,\, y,\, x^2, \, xy, \, y^2$, and $f, xf, yf, f^2$.
The Fubini formula transforms an integration on a bidimensional domain into two successive one dimensional  integrations. Indeed for any triangle $T$ and any $g\in L^1(T)$,
\be
\label{eqFubiniG}
\begin{array}{rcl}
\displaystyle \int_T g(x,y) \, dx \, dy &=& \displaystyle \int_{y_*(T)}^{y^*(T)} \left(\int_{x_*(T,y)}^{x^*(T,y)} g(x,y) dx \right) dy\\
& =& \displaystyle \int_{y_*(T)}^{y^*(T)} \(G(x^*(T,y),y) - G(x_*(T,y),y) \)dy.
\end{array}
\ee
Where we have used the notations 
\begin{eqnarray*}
G(x,y)    &:=& \int_{-\infty}^x g(x,y) dy,\\
y_*(T)    &:=& \min\{ y\in \R \sep \exists x \in \R \text{ such that } (x,y) \in T\},\\
x_*(T,y) &:=& \min\{ x\in \R \sep (x,y) \in T\},
\end{eqnarray*}
and where $y^*(T)$ and $x^*(T,y)$ are obtained by taking the maximum instead of the minimum in the definitions of $y_*(T)$ and $x_*(T,y)$ respectively. The important point in \iref{eqFubiniG} is that the function $G$ does not depend on the triangle $T$. If this function is known, then \iref{eqFubiniG} transforms the bidimensional integration of $g$ into a one dimensional integration.

This strategy extends to the discrete setting with the Lebesgue measure on $T$ 
replaced by the counting measure at pixels where the approximated function $f$ is sampled
and which have their {\it center} contained in $T$.
The integrals of the type $\int_T g(x,y) dx dy$ appearing in the previous equations are therefore replaced by discrete sums
\be
\label{eqIntDiscrete}
\sum_{(i,j)\in T} g(i,j),
\ee
where $(i,j)$ stands for the center of the pixel, and without loss of generality 
runs over $\{0,\cdots,N-1\}^2$ for a $N\times N$ image.
In particular the quantity $\|f-P_T f\|_2^2$ appearing in the definition \iref{eqDecisionL2} of the $L^2$ based decision function is replaced with 
\be
\label{eqErrorDiscrete}
\inf_{\pi \in \sP_1} \sum_{(i,j)\in T} |f(i,j) -\pi(i,j)|^2.
\ee
Again, this quantity can be expressed in terms of the (discrete) integrals of the functions $1,x,y, \cdots, xf,yf,f^2,$ which do not depend on $T$.
The Fubini formula \iref{eqFubiniG} is replaced with 
\be
\label{eqFubiniDiscrete}
\sum_{(i,j)\in T} g(i,j) = \sum_{j_*(T)\leq j \leq j^*(T)} G_{i^*(T,j), \, j} - G_{i_*(T,j), \, j}
\ee
where 
\begin{eqnarray*}
G_{i,j} &=& \sum_{0\leq i'\leq i} g(i',j),\\
j_*(T) &=& \lfloor y_*(T) \rfloor\\
i_*(T,j) &=&  \lfloor x_*(T,j) \rfloor,
\end{eqnarray*}
and similarly  $j^*(T) = \lfloor y^*(T) \rfloor$ and $i^*(T,j) =  \lfloor x^*(T,j) \rfloor$. The points $(i_*(T,j),\, j)$, for all $j_*(T)\leq j \leq j^*(T)$, are illustrated on the right of Figure \ref{figFastIntegration} and surrounded by grayed  squares. On the same figure the points $(i^*(T,j),\, j)$, for all $j_*(T)\leq j \leq j^*(T)$, are surrounded by white squares.

Consider a fixed function $g$, such as $1,x,y,x^2,xy, y^2$ or $f,xf,yf,f^2$ which are needed for the evaluation of \iref{eqErrorDiscrete}, for which the sum \iref{eqIntDiscrete} needs to be evaluated for a large number of distinct triangles. The matrix $(G_{i,j})$, $0\leq i \leq N$, $1\leq j \leq N$ is evaluated once, and for any triangle $T$ the sum \iref{eqIntDiscrete} is computed using \iref{eqFubiniG} which involves a much smaller collection of points, see Figure \ref{figFastIntegration}.
%
This method greatly accelerates the numerical implementation of the greedy algorithm, especially for the large and isotropic triangles that occur in its first steps.

\begin{figure}
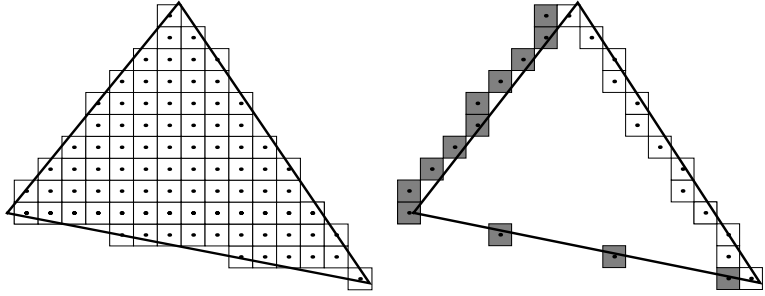

	\centering
	\includegraphics[height = 4cm, width = 5cm]{\pathPic/Triangles/AllPts.pdf}
	\includegraphics[height = 4cm, width = 5cm]{\pathPic/Triangles/FewPts.pdf}
	\caption{\label{figFastIntegration}The cost of evaluating the interpolation error of a function $f$ on a triangle $T$  in the $L^1$ or $L^\infty$ norm is proportional to the number of points with integer coordinates on $T$ (left). Fewer points need to be considered to compute the $L^2$ projection error (right).} 
\end{figure}

\subsection{Positive quadratic functions}
\label{secPosQuad}

In this section, we study the algorithm 
when applied to a quadratic polynomial $q$ such that $\det(\b q)>0$. 
We shall assume without loss of generality that
$\b q$ is positive definite, since all our results
extend in a trivial manner to the negative definite case.

We first establish that 
the refinement procedure, using either the $L^2$ based decision function \iref{eqDecisionL2} or the $L^\infty$ based \iref{eqDecisionLInf},
always selects for bisection the longest edge in the sense of the $\bq$-metric $|u|_{\bq} := \sqrt{\bq(u)}$ for $u\in \R^2$, as was already the case for the $L^1$ based decision function \iref{eqDecisionL1}.
This is used to prove that the refinement procedure
produces triangles which tend to adopt an
optimal aspect ratio. 

\subsubsection{The $L^\infty$-based split}
\noindent
Let us denote by
$$
\alpha_T(f) := \|f-I_Tf\|_{L^\infty(T)},
$$
the interpolation error in the sup norm.
The decision function \iref{eqDecisionLInf} can be re-expressed
as
\be
d_T(e,f) = \alpha_{T_e^1}(f) + \alpha_{T_e^2}(f).
\label{optilinfedge}
\ee

\begin{theorem} \label{longestedgelinf}
If $|a|_{\bq}> \max\{|b|_{\bq},|c|_{\bq}\}$,
then $d_T(a,q) <\min\{d_T(b,q),d_T(c,q)\}$. Therefore the refinement
procedure based on {\rm \iref{optilinfedge}} selects the longest edge in the sense of the $\bq$-metric.
\end{theorem}

In order to prove this result, we need to study the interpolation
error in detail.

\begin{prop} \label{PropAlgErLInf}
Let T be a triangle with edges $a,b,c$ such that $|a|_\bq\geq |b|_\bq\geq |c|_\bq$,
and let $w\in \R^2$ and $r>0$ be the center and radius of 
the circumscribed circle for the $\bq$-metric,
i.e. such that $|v-w|_\bq=r$ for all the vertices $v$ of $T$. 
Then
$$
\frac{|a|_\bq^2} 4 \leq \alpha_T(q) \leq r^2.
$$ 
Right equality holds if $T$ is acute, i.e. $\< Q b, c\>\leq 0$, where $Q$ is the symmetric matrix associated to the quadratic form $\bq$. 
Left equality holds if $T$ is obtuse, i.e. $\< Q b, c\>\geq 0$. 
\end{prop}

\proof At any point $u\in \RR^2$, we have 
$$
(q - I_T q)(u) = |u-w|_\bq^2 - r^2.
$$
Indeed, the difference $(q - I_T q)(u) - (|u-w|_\bq^2 - r^2)$ is an affine function of $u$, which vanishes at the three vertices of $T$. Hence this difference is zero.

The function $|u-w|_\bq^2 - r^2$ is negative on $T$ with maximal value $0$ at the
vertices. If $T$ is acute, then its minimal value on $T$ is $-r^2$
and is attained at $w\in T$. If $T$ is not acute, 
then the minimum is attained at $m_a$, the midpoint of $a$, 
and if we choose a vertex $v$ at one end of $a$, 
we obtain the value at the minimum by 
Pythagoras' identity which gives
$$
\begin{array}{ll}
(q - I_T q)(m_a) &= |m_a-w|_\bq^2 - r^2= |m_a-w|_\bq^2 -|v-w|_\bq^2\\
& =-|v-m_a|_\bq^2 =-|a|_\bq^2/4.
\end{array}
$$ \sq 

The dichotomy in the above result is illustrated in the case of the
euclidean metric on figure 4. Note that it would be
sufficient to establish the above proof in this particular case,
since we can perform an affine coordinate change $\phi=Q^{-\frac 1 2}$ 
such that $\b q\circ \phi$ is the standard euclidean form
and that the $L^\infty$ interpolation error is left invariant
by this coordinate change.

\begin{figure}[htbp]
\centerline{
\includegraphics[width=8cm]{\pathPic/PaperCM/maxpoint.pdf}
}
\centerline{Figure 4: maximum point for the $L^\infty$ interpolation error}
\label{figerrorlinf}
\end{figure}

\noindent
We now prove the following result which clearly implies Theorem
\ref{longestedgelinf}.

\begin{prop} \label{ErLInf} Let T be a triangle with edges $|a|_\bq\geq |b|_\bq\geq |c|_\bq$.
We then have :
\begin{eqnarray*}
d_T(b,q)-d_T(a,q)&\geq & \frac 1 4 (\b q(a)-\b q(b)) \\
d_T(c,q)-d_T(a,q)&\geq & \frac 1 4 \left(\b q(a)-\left(\frac{|b|_\bq+|c|_\bq} 2\right)^2\right).
\end{eqnarray*}
\end{prop}

\begin{figure}[htbp]
\centerline{
\includegraphics[width=15cm,height=4cm]{\pathPic/PaperCM/cutabc.pdf}
}
\centerline{Figure 5: Notations in the proof of Proposition \iref{ErLInf}}
	\label{fig:cuts}
\end{figure}

\proof  
We introduce sub-triangles $T^i_e$, $i=1,2$ and $e=a,b,c$,
as defined in Figure 5, which correspond 
to the three refinement scenarios. With such definitions, the
following inequalities are easily derived from
Proposition \iref{PropAlgErLInf}
\begin{eqnarray*}
4\alpha_{T^2_a}(q) &=& \b q(b) \text{ (since } T^2_a \text{ is obtuse)}\\
4\alpha_{T^1_b}(q) &\geq& \b q(a)\\
4\alpha_{T^1_c}(q) &\geq& \b q(a)\\
4\alpha_{T^2_c}(q) &\geq& \b q(b)
\end{eqnarray*}
On the other hand, we shall prove 
\be
\alpha_{T^1_a}(q) \leq \alpha_{T^2_b}(q),
\label{estim1}
\ee
and
\be
4\alpha_{T^1_a}(q) \leq \left(\frac{|b|_\bq+|c|_\bq} 2\right)^2.
\label{estim2}
\ee
The proof of \iref{estim1} follows from elementary geometric observations.
Let $L$ be a line which is parallel to $c$ but does not contain it, and for $x\in L$
denote by $T_x$ the triangle of vertices $x$ and the end points of $c$.
Denoting respectively by $u(x)$ and $v(x)$ the diameter of $T_x$ 
and of its circumscribed circle for the $\bq$-metric, we remark that these functions 
decrease monotonously as $x$ tends to a point $x_c$ which is the orthogonal
projection (also for the $\bq$-metric) of the mid-point of $c$ onto $L$. Since
the function $x\mapsto \alpha_{T_x}(q)$ is continuous in $x$ and
equal to $u(x)$ or $v(x)$ at all $x$, we conclude that this
function also decreases monotonously as $x$ tends to $x_c$.
Applying this observation to the line that contains $m_a$ and $m_b$
the mid-points of $a$ and $b$, and remarking that
$m_a$ is closer to $x_c$ than $m_b$, we conclude that \iref{estim1} holds.

From \iref{estim1} and the first set of inequalities, we obtain
the first statement of the theorem since
$$
\begin{array}{ll}
d_T(b,q)-d_T(a,q)&=\alpha_{T^1_b}(q)+\alpha_{T^2_b}(q)-\alpha_{T^2_a}(q)-\alpha_{T^1_a}(q)\\
& \geq \alpha_{T^1_b}(q)-\alpha_{T^2_a}(q)
\geq  \frac 1 4(\b q(a)-\b q(b)).
\end{array}
$$

The proof of \iref{estim2} also follows from elementary geometric observations.
In the case where $T^1_a$ is obtuse, one of its edges $e$ is such
that $4\alpha_{T^1_a}=\b q(e)$, and \iref{estim2} follows since
$|e|_\bq\leq \frac 1 2 (|b|_\bq+|c|_\bq)$ for all $e$, using triangle inequality.
Let $R$ be an acute triangle of vertices $u,v,w$ and let $m = (u+v)/2$. Remarking than the center of the (euclidean) circumscribed circle to $R$ lies inside $R$, one easily checks by convexity that the (euclidean) diameter of this circle is smaller than $|m-u|+ |m-w|$.
In the case where $T^1_a$ is acute, the diameter of its circumscribed circle is thus bounded by $\frac 1 2 (|b|_\bq+|c|_\bq)$, as illustrated
on Figure 6 when $\bq$ is the euclidean metric.

From \iref{estim2} and the first set of inequalities, we obtain
the second statement of the theorem since
$$
\begin{array}{ll}
d_T(c,q)-d_T(a,q)&=\alpha_{T^2_c}(q)+\alpha_{T^1_c}(q)-\alpha_{T^2_a}(q)-\alpha_{T^1_a}(q) \\
&\geq \frac 1 4\left(\b q(b)+\b q(a)-\b q(b)-\left(\frac{|b|_\bq+|c|_\bq} 2\right)^2\right)\\
&=\frac 1 4 \left(\b q(a)-\left(\frac{|b|_\bq+|c|_\bq} 2\right)^2\right).
\end{array}
$$
\sq

\begin{figure}[htbp]
\centerline{
\includegraphics[width=4cm]{\pathPic/PaperCM/acutecase.pdf}
}
\centerline{Figure 6: The case where $T^1_a$ is acute.}
	\label{fig:bcD}
\end{figure}

\subsubsection{The $L^2$-based split}

We now denote by
$$
\beta_T(f) := \|f-P_Tf\|_{L^2(T)},
$$
the orthogonal projection error in the $L^2$ norm.
The decision function \iref{eqDecisionL2} now writes
\be
d_T(e,f) =\beta_{T_e^1}(f) ^2+ \beta_{T_e^2}(f)^2.
\label{optil2edge}
\ee
We shall prove that the refinement procedure
based on \iref{optil2edge} behaves in a similar way
as \iref{optilinfedge}.

\begin{theorem} \label{longestedgel2}
If $d,e\in \{a,b,c\}$ are two edges such that $|d|_\bq<|e|_\bq$, then
$d_T(e,q) <d_T(d,q)$. Therefore the refinement
procedure based on {\rm \iref{optil2edge}} selects 
the longest edge in the sense of $|\cdot|_\bq$.
\end{theorem}

In order to prove this result, we first provide with an algebraic expression
of $\beta_T(q)$ which is valid for any quadratic function $q$.

\begin{prop} \label{PropAlgEr}
Let $T$ be a triangle with edges $a,b,c$ and area $|T|$, 
and let $q$ be a quadratic function.
Then
\be
\beta_T^2(q) = |T| \left(c_1 (\b q(a)+\b q(b)+\b q(c))^2 - c_2 \det(\b q) |T|^2\right).
\label{localprojqerror}
\ee
with constants $c_1 = \frac 1{1200}$
and $c_2 = c_1  \frac {64} 3  = \frac {4}{225}$.
\end{prop}

\proof 
We first prove \iref{localprojqerror}
on the triangle $R$ of vertices $\{(0,0),(0,1),(1,0)\}$.
It is easy to compute the integrals on $R$
of monomials $x^ky^l$, $k+l\leq 4$.
Using these quantities, we can derive the orthogonal 
projection of a quadratic function thanks to a formal computing program,
which gives us
$$
\bq = ux^2+vy^2+2wxy \quad \Ra \quad P_R \bq = -\frac{u+v+w}{10} + \frac{2x} 5 (2u+w)+\frac{2y} 5 (2v+w).
$$
This yields the following expression
for the $L^2$-squared error between $q$ and its projection  
$$
\int_{R} (q-P_R q)^2 = 
\int_{R} (\bq-P_R \bq)^2 = 
\frac 1 {300} \left(u^2 + \frac{2 u v} 3 + v^2 - 2 u w - 2 v w + \frac{7 w^2} 3\right),
$$
which is equivalent to \iref{localprojqerror}.

For an arbitrary triangle $T$, using an affine bijective transformation
$\phi$ from $R$ to $T$, we have
$$
\beta_T(q)^2=J_\phi \beta_R(\ti q)^2,
$$
where $\ti q=q\circ \phi$ and $J_\phi$ is the constant jacobian of $\phi$. 
Using the validity of \iref{localprojqerror} on $R$
and the fact that $|T|=J_\phi |R|$, we thus obtain 
$$
\beta_T(q)^2=|T| \left(c_1 (\b {\ti q}(\ti a)+\b {\ti q}(\ti b)+\b {\ti q}(\ti c))^2 - c_2 \det(\b {\ti q}) |R|^2\right),
$$
where $\b {\ti q}$ is the quadratic form associated to $\ti q$ and 
$\ti e$ denotes the edge segment of $R$ mapped onto $e$ by $\phi$.
Since $\b {\ti q}(\ti e)=\b q(e)$ and $\det(\b {\ti q})=J_\phi^2\det(\b q)$,
we obtain \iref{localprojqerror} for $T$.
\sq
\nl
We now prove the following result which clearly implies Theorem
\ref{longestedgel2}.

\begin{corollary} \label{PropErrorChildren} Let $T$ be a triangle with edges $a,b,c$ and area $|T|$, 
with $|a|_\bq\geq |b|_\bq$
and $|a|_\bq\geq  |c|_\bq$. Then
\begin{eqnarray}
\label{eqDecisionL2PosAB}
	d_T(b,q)-d_T(a,q) &\geq & \frac 5 4 c_1 |T| (\b q(a)^2-\b q(b)^2),\\
\label{eqDecisionL2PosAC}
	d_T(c,q)-d_T(a,q) &\geq & \frac 5 4 c_1 |T| (\b q(a)^2-\b q(c)^2).
\end{eqnarray}
\end{corollary}

\proof
The children triangles all have area $|T|/2$, and 
take their edges among $a,b,c$, 
$a/2,b/2,c/2$ and $\frac{a-b} 2,\frac{b-c} 2,\frac{c-a} 2$ (recall that $a+b+c=0$).
We use the identity
$$
\b q(u+v)+\b q(u-v) = 2\b q(u)+2\b q(v),
$$
which is valid for all quadratic forms, and implies
$$
\bq\left(\frac{b-c} 2\right) = \frac{\bq(b)+\bq(c)} 2-\frac{\bq(a)} 4.
$$
Using \iref{localprojqerror}, this allows us to compute the local projection 
errors for the children of $T$.
For example bisecting the edge $a$ creates two children $T'$ 
and $T''$ with edges $\frac a 2,b,\frac{c-b} 2$ and $\frac a 2,c,\frac{b-c} 2$, and therefore
\begin{eqnarray*}
\beta_{T'}^2(q) &=& |T'| \left(c_1 \left(\b q\left(\frac a 2\right)+\b q(b)+\b q\left(\frac{c-b} 2\right)\right)^2 - c_2 \det(\b q) |T'|^2\right)\\
&=& \frac{|T|} 2 \left(c_1 \left(\frac{3 \b q(b)+\b q(c)} 2\right)^2 - c_2 \det(\b q) \frac{|T|^2} 4\right)
\end{eqnarray*}
and similarly 
$$
\beta_{T''}^2(q) = \frac{|T|} 2 \left(c_1 \left(\frac{3 \b q(c)+\b q(b)} 2\right)^2 - c_2 \det(\b q) \frac{|T|^2} 4\right).
$$
Adding up, we thus obtain
$$
d_T(a,q) =\frac{|T|} 8 \left(c_1 (3 \b q(b)+\b q(c))^2 + c_1 (\b q(b)+3 \b q(c))^2 - 2 c_2 \det(\b q) |T|^2 \right).
$$
Subtracting this from the analogous expression for $d_T(b,q)$
we obtain 
\be
d_T(b,q) - d_T(a,q) = \frac{c_1 |T|} 4 \left(5 \left(\bq(a)^2 - \bq(b)^2\right) + 6 \bq(c) (\bq(a)-\bq(b))\right)
\label{dtaq}
\ee
which implies \iref{eqDecisionL2PosAB}.
Exchanging $b$ with $c$ we obtain \iref{eqDecisionL2PosAC} which concludes the proof.
\sq

\subsubsection{Convergence towards to optimal aspect ratio}
\label{subsecRatioConvex}

We have established in Theorems \ref{longestedgelinf} and \ref{longestedgel2} that for any triangle $T$ and any quadratic function $q\in \P_2$ such that the homogeneous component $\bq\in \H_2$ is positive definite, the decision functions \iref{eqDecisionL2} and \iref{eqDecisionLInf} based on the $L^2$ or $L^\infty$ error lead to the bisection of the same edge of $T$ : the longest edge of the $\bq$-metric.

Similarly it is established in Lemma \ref{lemmaLongestEdgeL1} in the previous chapter that the decision function \iref{eqDecisionL1} based on the $L^1$ interpolation error leads to the bisection of the same edge. As a result, the decision functions can be interchanged without changing the
result of Corollary \ref{corolQuadApprox} which extends as follows.

\begin{theorem}
\label{thOptiL2LInf}
Let $\Omega$ be a triangle, and let $q$ be a quadratic function with positive definite associated quadratic form $\bq$. 
Let $q_N$ be the approximant of $q$ on $\Omega$ obtained by the greedy algorithm 
for the $L^p$ metric, using the $L^1$, $L^2$ or $L^\infty$ decision function  \iref{eqDecisionL1},  \iref{eqDecisionL2}, \iref{eqDecisionLInf}.
Then 
\be
\limsup_{N\to \infty} N \|q-q_N\|_{L^p(\Omega)}  \leq C \|\sqrt{\det(\bq)}\|_{L^\tau(\Omega)},
\ee
where $\frac 1 \tau=\frac 1 p +1$ and where the constant $C$ depends only on 
on the choice of the approximation operator $\cA_T$ used in the definition of the approximant.  
\end{theorem}

The local perturbation analysis exposed in the previous chapter, \S \ref{section:OptQuad}, for the $L^1$ based decision function extends without difficulties to the $L^2$ and $L^\infty$ based decision functions. 
From this point, it is not difficult to extend Theorem \ref{thOptiL2LInf} 
to functions $f$ which are close enough to a positive definite quadratic function $q$, in the sense that 
\be
\label{eqFApproxQuad}
d^2 q \leq d^2 f\leq (1+\mu) d^2 q
\ee
for a constant $\mu>0$ small enough.

Let us consider a $C^2$ and strictly convex function $f$ defined on a polygonal domain $\Omega$. In order to establish an asymptotical error estimate, 
we need to be ensured that after sufficiently many steps of the greedy algorithm, 
the target function $f$ can be well approximated by a quadratic function $q=q_T$ 
on each triangle $T$, in the sense of \iref{eqFApproxQuad}, so that our local results will apply on such triangles.
This property follows from the next proposition, which is an extension of Proposition \ref{PropDiamZero} in which the $L^1$ based decision function \iref{eqDecisionL1} is replaced with the $L^\infty$ based decision function.
We postpone its proof to the appendix \S\ref{secPropDiamZero}.

\begin{prop} \label{propDiamZeroLInf}
Let $f$ be a $C^2$ function such that 
$d^2 f(x)\geq m I $ 
for some arbitrary but fixed $m >0$ independent of $x$. 
Let $\cT_N$ be the triangulation generated by the greedy algorithm
applied to $f$ using the $L^\infty$ decision function given by \iref{eqDecisionLInf}.
Then
$$
\lim_{N\to +\infty} \max_{T\in\cT_N}\diam(T)=0,
$$
i.e. the diameter of all triangles tends to $0$. 
\end{prop}

As a consequence, we may extend Theorem \ref{optitheo} to the $L^\infty$ based decision function.
\begin{theorem}
\label{optitheoLInf}
Let $\Omega$ be a polygonal domain and let $f\in C^2(\Omega)$ be such that 
$
d^2f(x)\geq m I
$ 
for all $x\in \Omega$,
for some arbitrary but fixed $m >0$ independent of $x$. Let $f_N$
be the approximant obtained by the greedy algorithm 
for the $L^p$ metric, using the $L^1$ decision function \iref{eqDecisionL1} or the $L^\infty$ decision function \iref{eqDecisionLInf}.
Then 
$$
\limsup_{N\to \infty} N \|f-f_N\|_{L^p}  \leq C \|\sqrt{\det(d^2f)}\|_{L^\tau},
$$
where $\frac 1 \tau=\frac 1 p +1$ and where $C$
is a constant independent of $p$, $f$ and $m$.  
\end{theorem}


\subsection{Quadratic functions of mixed sign}
\label{secNegQuad}

In this section, we study the algorithm 
when applied to a quadratic polynomial $q$ such that $\det(\b q)<0$. 
We shall follow the same steps, and reach similar conclusions, as in the positive definite case, using
a measure of non-degeneracy which is equivalent to $\rho_{\bq}(T)$.
If a triangle $T$ has edges $a,b,c$ such that $|\bq(a)|\geq |\bq(b)|\geq |\bq(c)|$, we will still refer to $a$ as the ``longest'' edge in the sense of $\bq$, although $\bq$ does not define a proper metric anymore. Recall that 
$$
\rho_{\bq}(T) := \frac{|\bq(a)|}{|T|\sqrt{|\det \bq|}}.
$$

The following inequalities that will be repeatedly used in this section can be derived when $\rho_\bq(T)$ is large enough. We postpone their proof to the appendix.
\begin{prop}
\label{rho48}
Let $T$ be a triangle such that $|\bq(a)|\geq |\bq(b)|\geq |\bq(c)|$, and define
$d= \frac{b-c} 2$.
\be 
\label{rho4}
\text{If }\rho_\bq(T) \geq 4, \text{ then } \bq(a)\bq(b)\geq 0, \ |\bq(a)|\geq 4 |\bq(c)| \text{ and } |\bq(a)|\geq 4 |\bq(d)|.
\ee
\be
\label{rho8}
\text{If }\rho_\bq(T) \geq 8, \text{ then } |\bq(a)|\leq \frac 3 8 |\bq(b)+\bq(c)|.
\ee
\end{prop}

\subsubsection{The $L^\infty$-based split}

\begin{theorem} 
\label{thLInfSplitNeg}
The refinement
procedure based on {\rm \iref{optilinfedge}} selects the longest edge in the sense of $\bq$:
if $|\bq(a)|> \max\{|\bq(b)|,|\bq(c)|\}$ and $\rho_\bq(T)\geq 4$,
then $$d_T(a,q) <\min\{d_T(b,q),d_T(c,q)\}.$$
\end{theorem}

This theorem is very similar to the one for positive quadratic functions. In order to prove it, we first study the interpolation error which has a simple form in this context.

\begin{prop}
Let $T$ be a triangle with edges $a,b,c$. Then
$$
\alpha_T(q) = \frac 1 4 \max\{|\bq(a)|,|\bq(b)|,|\bq(c)|\}.
$$
\end{prop}

\proof
Let $x_0$ be the point of $T$ at which the interpolation error is attained: $x_0 = \argmax_T |q-I_Tq|$. If $x_0$ is in the interior of $T$, then it must be a local extremum of $q-I_Tq$. However this function has only one critical point on $\R^2$, which is not an extremum since $\bq$ has mixed signature. Therefore $x_0$ must lie on an edge.
On each edge of $T$, the function $q-I_Tq$ is a one dimensional quadratic function vanishing at the endpoints. It follows that $x_0$ must lie in the middle of an edge
and the result follows.
\sq

\begin{prop}
Let $T$ be a triangle with edges $|\bq(a)|\geq |\bq(b)| \geq |\bq(c)|$ and such that $\rho_\bq(T)\geq 4$. Then
\begin{eqnarray*}
d_T(b,q)-d_T(a,q) &\geq& \frac{|\bq(a)|-|\bq(b)|} 8,\\
d_T(c,q)-d_T(a,q) &\geq& \frac{|\bq(a)|} 8.
\end{eqnarray*}
\end{prop}

\proof
The bisection through the edge $a$ creates two sub-triangles $T^1_a,T^2_a$ of edges respectively $\frac a 2, b,d$ and $\frac a 2, c, d$. Using the last two inequalities in \iref{rho4} we obtain that $4 \alpha_{T^1_a} = \max\{|\bq(a/2)|,|\bq(b)|\}$ and $4 \alpha_{T^2_a} = |\bq(a/2)|$.
Therefore 
$$
4 d_T(a,q) =\frac{|\bq(a)|} 4 +\max\left\{\frac{|\bq(a)|} 4,|\bq(b)|\right\}.
$$
On the other hand, the choice of bisecting the edge $b$ creates two subtriangles respectively containing  the edges $a$ and $\frac b 2$, and 
the choice of bisecting the edge $c$ creates two subtriangles respectively containing  the edges $a$ and $b$. This provides us with the lower bounds
\begin{eqnarray*}
4 d_T(b,q) &\geq & |\bq(a)|+\frac{|\bq(b)|} 4,\\
4 d_T(c,q) &\geq & |\bq(a)|+|\bq(b)|.
\end{eqnarray*}
The proposition follows easily, distinguishing between the two cases $|\bq(a)|\leq 4 |\bq(b)|$
and $|\bq(a)|\geq 4 |\bq(b)|$.
\sq

\subsubsection{The $L^2$-based split}

The same conclusions can be reached for the refinement procedure based on \iref{optil2edge}.

\begin{theorem} 
\label{thL2SplitNeg}
The refinement
procedure based on {\rm \iref{optil2edge}} selects the longest edge in the sense of $\bq$: if $|\bq(a)|> \max\{|\bq(b)|,|\bq(c)|\}$ and $\rho_\bq(T)\geq 4$,
then $d_T(a,q) <\min\{d_T(b,q),d_T(c,q)\}$.
\end{theorem}

\proof
The expression found in \iref{dtaq} remains valid when $\det(\bq)<0$.
Substituting $a$ by $b$ or $c$ and subtracting, we obtain
\begin{eqnarray*}
d_T(b,q)-d_T(a,q) &=& \frac{5 c_1} 2 |T| (\bq(a)-\bq(b))\left(s+\frac{\bq(b)} 5\right),\\ 
d_T(c,q)-d_T(a,q) &=& \frac{5 c_1} 2 |T| (\bq(a)-\bq(c))\left(s+\frac{\bq(c)} 5\right),
\end{eqnarray*}
where $s=\bq(a)+\bq(b)+\bq(c)$.
Using \iref{rho4}, we see that the quantities 
$$
s+\frac{\bq(b)} 5, \quad s+\frac{\bq(c)} 5, \quad \bq(a)-\bq(b) \stext{ and } \bq(a)-\bq(c)
$$
all have the same sign as $\bq(a)$
and are non-zero. It follows that 
$$
d_T(a,q) <\min\{d_T(b,q),d_T(c,q)\}
$$
which concludes the proof.
\sq

\subsubsection{Convergence toward the optimal aspect ratio.}

We have proved that the refinement procedure - either based on the
$L^\infty$ or $L^2$ decision function - systematically picks the longest
edge in the sense of $\bq$ if $\rho_\bq(T) \geq 4$. Similarly to the positive definite
case, we now study the iteration of several refinement steps and
show that the generated triangles tend to adopt an optimal ``aspect ratio''
in the sense of the measure of non-degeneracy $\rho_{\bq}(T)$
introduced in \S \ref{secCM2}. If $T$ is a triangle with edges $a,b,c$, let us recall that 
$$
\rho_\bq(T) := \frac{\max\{|\bq(a)|, \, |\bq(b)|, \, |\bq(c)| \}}{|T| \sqrt{|\det \bq|}}.
$$

As in \S \ref{secCM3}, we introduce a close variant to 
$\rho_{\bq}(T)$. If $T$ is a triangle with edges $a,b,c$, we define
\be
\sigma_\bq(T) := \frac{\min\left\{|\b q(a)+\b q(b)|,\, |\b q(b)+\b q(c)|,\, |\b q(c)+\b q(a)|\right\}}{4|T|\sqrt{|\det \b q|}}.
\ee
Note that if $\bq$ was a positive quadratic form, this definition is
consistent with \iref{sigmaq}. We define our measure of non-degeneracy $\kappa_\bq$ by
\be
\kappa_\bq(T) := \max\{\sigma_\bq(T), 5/ 2\}.
\ee
We first show that the quantities $\kappa_\bq$ and $\rho_\bq$ are equivalent.
\begin{prop}
\label{rhosigma}
For any triangle $T$, one has
\be
2\sigma_\bq(T) \leq \rho_\bq(T),
\label{sigmarho}
\ee
and
\be
\frac 4 5\kappa_\bq(T) \leq \rho_\bq(T) \leq \frac{32} 3 \kappa_\bq(T).
\label{rhokappa}
\ee
\end{prop} 

\proof
We denote by $a,b,c$ the edges of $T$, and we assume that $|\bq(a)|\geq |\bq(b)| \geq |\bq(c)|$.
The inequality \iref{sigmarho} follows directly from the triangle inequality:
$$
2 |T| \sqrt{|\det \bq|} \sigma_\bq(T) \leq \frac{|\bq(b)+\bq(c)|} 2 \leq |\bq(a)| \leq |T| \sqrt{|\det \bq|} \rho_\bq(T).
$$
As mentioned earlier, $\rho_\bq(T)$ is always larger than $2$
and therefore \iref{sigmarho} implies the left inequality in \iref{rhokappa}.

It remains to prove the right inequality in \iref{rhokappa}.
If $\rho_\bq(T)\leq 8$, it is immediate
since $\kappa_\bq(T)\geq \frac 5 2$ and $\frac{32} 3 \frac 5 2 \geq 8$.
If $\rho_\bq(T)\geq 8$ then $\bq(a)\bq(b)\geq 0$ according to \iref{rho4} and therefore 
$$
|\bq(b)+\bq(c)|\leq |\bq(a)+\bq(c)|\leq |\bq(a)+\bq(b)|.
$$
We obtain using \iref{rho8} that 
$$
\rho_\bq(T)=\frac{|\bq(a)|}{|T| \sqrt{|\det \bq|}} \leq \frac 8 3 \frac{|\bq(b)+\bq(c)|}{|T| \sqrt{|\det \bq|}}
= \frac{32} 3 \sigma_\bq(T)\leq \frac{32} 3 \kappa_\bq(T),
$$
which concludes the proof.
\sq

Similar to $\rho_\bq$, the quantity $\kappa_\bq$ is invariant by a linear coordinate changes $\phi$,
in the sense that
$$
\kappa_{\bq\circ \phi}(T) = \kappa_\bq (\phi(T)).
$$
Our next result shows that $\kappa_\bq(T)$ is always reduced by 
the refinement procedure.

\begin{prop} 
\label{PropSigmaDecmixed}
If $T$ is a triangle with children $T_1$ and $T_2$ obtained by the
refinement procedure for the quadratic function $q$, using the $L^2$ based or $L^\infty$ based decision function, then
$$
\max\{\kappa_\bq(T_1), \kappa_\bq(T_2)\}\leq \kappa_\bq(T).
$$
\end{prop}

\proof
Let us assume that $a$ is the longest edge in the sense of $\bq$.
In the case where $\rho_\bq(T)\geq 4$, we already noticed
in the proof of Proposition \ref{rhosigma} 
that 
$$
\sigma_\bq(T) = \frac{|\b q(b)+\b q(c)|}{4|T|\sqrt{|\det \b q|}}.
$$
Moreover, according Theorems \ref{thLInfSplitNeg} and \ref{thL2SplitNeg} 
the edge $a$ is selected by both decision functions.
It follows that the children $T_i$ have edges $a/2,b,(c-b)/2$ and $a/2, (b-c)/2,c$ 
(recall that $a+b+c=0$). We thus have
\begin{eqnarray*}
	2 |T| \sqrt{|\det \bq|}\ \sigma_\bq(T_i) & \leq & \left|\bq\left(\frac a 2 \right)+\bq\left(\frac{b-c}{2}\right)\right|\\
	& = & \left|\bq\left(\frac{b+c}{2}\right)+\bq\left(\frac{b-c}{2}\right)\right|\\
	& = & \frac{|\b q(b)+\b q(c)|}{2}\\
	& = & 2 |T|  \sqrt{|\det \b q|}\ \sigma_\bq(T).
\end{eqnarray*}
We have proved that $\sigma_\bq(T_i)\leq \sigma_\bq(T)$, and it readily follows that $\kappa_\bq(T_i)\leq \kappa_\bq(T)$.

In the case where $\rho_\bq(T)\leq 4$, we remark that $T_i$
contains at least one edge from $T$, say $s\in\{a,b,c\}$ 
and one half-edge $t\in \{\frac a 2,\frac b 2,\frac c 2\}$.
This provides an upper bound for $\sigma_\bq$ : 
\be
\label{rho4sigma}
\sigma_\bq(T_i) \leq \frac{|\bq(s)+ \bq(t)|}{2 |T| \sqrt{|\det \bq|}} \leq \frac{|\bq(a)|+ |\bq(\frac a 2)|}{2 |T| \sqrt{|\det \b q|}} = \frac 5 8 \rho_\bq(T) \leq \frac 5 2.
\ee
Therefore $\kappa_\bq(T_i) = \frac 5 2 \leq \kappa_\bq(T)$.
\sq

Our next objective is to show that as we iterate the refinement process,
the value of $\kappa_\bq(T)$ 
becomes bounded independently of $q$
for almost all generated triangles.
If $T$ is a triangle such that 
$|\bq(a)|\geq |\bq(b)|\geq |\bq(c)|$ and 
if the edge $a$ is cut (which is the case as soon as $\rho_\bq(T)\geq 4$, using the $L^2$ or $L^\infty$ based decision functions) 
we define $\psi_\bq(T)$ as the subtriangle containing the edge $c$.
We first prove a result which is analogous to Proposition \ref{PropFSigmaDec}.

\begin{prop} 
\label{PropFSigmaDecmixed} 
If $T$ is a triangle such that $\kappa_\bq(\psi_\bq^3(T))>\frac 5 2$, then 
$\kappa_\bq(\psi_\bq^3(T))\leq \frac 2 3 \kappa_\bq(T)$.
\end{prop}

\proof
Let $S$ be a triangle such that $\kappa_\bq(\psi_\bq(S))>\frac 5 2$. According to \iref{rho4sigma}, one must have $\rho_\bq(S)> 4$.
We assume that the edges of $S$ satisfy $|\bq(a)|\geq|\bq(b)|\geq|\bq(c)|$.
Since the three edges of $\psi_\bq(S)$ are $\frac a 2$, $c$ and $d=\frac{b-c} 2$, it follows 
from \iref{rho4} that the longest edge of $\psi_\bq(S)$ in the sense of 
$\bq$ is $\frac a 2$.

Since 
$$
\kappa_\bq(\psi_\bq(T)) \geq \kappa_\bq(\psi_\bq^2(T))\geq \kappa_\bq(\psi_\bq^3(T))>\frac 5 2,
$$
we can apply this observation to the triangles $T$, $\psi_\bq(T)$ and $\psi_\bq^2(T)$.
Therefore, denoting
by $a$ the longest edge of $T$ in the sense of $\bq$, 
we find that $\frac a 8$ is the longest edge of $\psi_\bq^3(T)$.
Since $|\psi_\bq^3(T)| = |T|/8$, we obtain that $\rho_\bq(\psi_\bq^3(T)) = \rho_\bq(T)/8$.
Using the results of Proposition \ref{rhosigma}, we thus have
$$
\kappa_\bq(\psi_\bq^3(T)) = \sigma_\bq(\psi_\bq^3(T))\leq \frac 1 2 \rho_\bq(\psi_\bq^3(T)) = \frac 1 {16} \rho_\bq(T) \leq \frac 2 3 \kappa_\bq(T),
$$
which concludes the proof.
\sq

An immediate consequence of Propositions \ref{PropSigmaDecmixed} and \ref{PropFSigmaDecmixed} is the following.

\begin{corollary} 
If $(T_i)_{i=1}^{8}$ are the 8 children obtained from 3 successive refinement procedures 
from $T$ for the function $q$, then
\begin{itemize}
	\item $\text{for all } i, \kappa_\bq(T_i)\leq \kappa_\bq(T)$,
	\item $\text{there exists } i \text{ such that }, \kappa_\bq(T_i)\leq \frac{2}{3}\kappa_\bq(T)$ or $\kappa_\bq(T_i)=\frac 5 2$.\\
\end{itemize}
\end{corollary}

\noindent
We eventually obtain
the two following results which proof are exactly similar to
the ones of Theorem \ref{CorolStabilise} and Corollary \ref{corolQuadApprox}.

\begin{theorem}
Let $T$ be a triangle, and let $\bq$ a 
quadratic function of mixed type. Let  $k = \frac{\ln (2\kappa_\bq(T)/5)}{\ln3-\ln 2}$.
Then after $n$ applications of the refinement procedure starting from $T$,
at most $C n^k 7^{n/3}$ of the $2^n$ generated triangles are such
that $\kappa_\bq(S)>\frac 5 2$ where $C$ is an absolute constant. Therefore the proportion 
of such triangles tends exponentially to $0$ as $n\to +\infty$.\\
\end{theorem}

\begin{theorem}
\label{thOptiL2LInfNeg}
Let $\Omega$ be a triangle, and let $q$ be a quadratic function such that $\det \bq<0$.
Let $q_N$ be the approximant of $q$ on $\Omega$ obtained by the greedy algorithm 
for the $L^p$ metric, using the $L^2$ or $L^\infty$ decision function \iref{eqDecisionL2}, \iref{eqDecisionLInf}.
Then 
$$
\limsup_{N\to \infty} N \|q-q_N\|_{L^p(\Omega)}  \leq C \|\sqrt{\det(\bq)}\|_{L^\tau(\Omega)},
$$
where $\frac 1 \tau=\frac 1 p +1$ and where the constant $C$ depends only on 
on the choice of the approximation operator $\cA_T$ used in the definition of the approximant.  
\end{theorem}

\section{Alternative bisection choices, and the approximation of cartoon functions}
\label{secAltBi}


A popular, although simplistic, model for images is the set of cartoon functions, see Definition \ref{defcartoon}. 
In summary a function $f$ defined on a domain $\Omega$ is a cartoon function if $f$ is $C^2$ except along a finite collection of $C^2$ curves where $f$ may have discontinuities.
For any cartoon function $f$ on the bi-dimensional domain $\Omega =[0,1]^2$, the approximations  $(f_N)_{N\geq 0}$ obtained by retaining the $N$ largest coefficients in the wavelet expansion of $f$ satisfy, as observed in 
\iref{wavN},
\be
\label{eqWav}
\|f-f_N\|_{L^2(\Omega)}\leq C N^{-\frac 1 2}. 
\ee
One of the purposes of anisotropic mesh adaptation is to improve on this estimate, 
and we therefore focus in this section on the $L^2$ approximation error.
For any triangle $T$ and any function $f\in L^2(T)$, we define 
\be
\label{defErrTL2}
e_T(f) := \inf_{\pi \in \sP_1} \|f-\pi\|_{L^2(T)} = \|f-P_T(f)\|_{L^2(T)},
\ee
where $P_T$ is the operator of $L^2(T)$ orthogonal projection onto the space $\P_1$ of affine functions.
For any triangulation $\cT$ of a polygonal domain $\Omega$, and any $f\in L^2(\Omega)$, we denote by $e_\cT(f)$ the error of approximation of $f$ in the $L^2(\Omega)$ norm by discontinuous piecewise affine functions on $\cT$. This quantity is defined by 
\be
\label{defErrCTL2}
e_\cT(f)^2 := \sum_{T\in \cT} e_T(f)^2 = \|f- P_\cT(f)\|_{L^2(\Omega)}^2,
\ee
where $P_\cT$ is the $L^2(\Omega)$ orthogonal projection onto the space of discontinuous piecewise affine functions on the triangulation $\cT$.
The heuristic analysis led in \iref{heurist} suggests that, if $f$ is a cartoon function on a polygonal domain $\Omega$, then  there exists a sequence $(\cT_N)_{N\geq N_0}$ of anisotropic triangulations of $\Omega$, satisfying $\#(\cT_N)\leq N$, and such that 
\be
\label{estimC}
e_{\cT_N}(f)\leq C N^{-1}. 
\ee

\begin{figure}
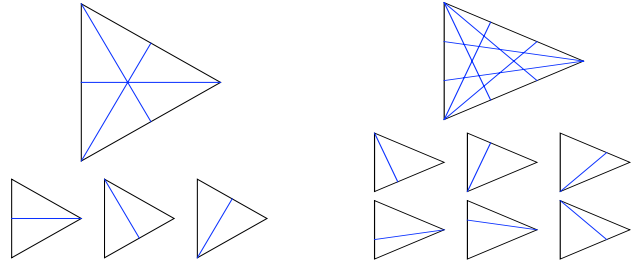

\centering
\includegraphics[width=3.5cm,height=3.5cm]{\pathPic/Subdivision/Bisec.pdf}
\hspace{1cm}
\includegraphics[width=3.5cm,height=3.5cm]{\pathPic/Subdivision/TriSec.pdf}
\caption{(left) Three bisection choices $S_T^3$, (right) Six bisection choices $S_T^6$.\label{figTriSec}}
\end{figure}

We make in this section an attempt, not totally conclusive, to answer the following question : is it possible to produce a sequence $(\cT_N)_{N \geq N_0}$ of triangulations satisfying \iref{estimC} using a hierarchical refinement algorithm? The analysis presented in \S\ref{secNegRes} shows that unfortunately the refinement algorithm studied in the previous chapters does \emph{not} satisfy \iref{estimC}. This failure is inherently due
to the limited geometric selectivity of the algorithm due to the choice of bisections.
For this reason, we study other types of bisections in \S \ref{secTriSec} and \S \ref{secConvexCut},
also based on splitting a selected element (which may not anymore be a triangle) by a line cut.
In this more general setting, the greedy algorithm has the following generic form:
\begin{enumerate}
\item An element $T$ which maximizes the approximation error is selected within the partition $\cT$.
$$
T = \underset {T'\in \cT} \argmax \ e_{T'}(f).
$$
\item The refinement procedure considers a finite set $S_T$ of segments $s$ along which the element $T$ may be split into two, thus creating two children elements $T_s^1$ and $T_s^2$. The segment $s\in S_T$ which minimizes a given decision function $d_T(s,f)$ is selected and we define 
$$
\cT' := \cT - \{T\} + \{T_s^1, T_s^1\}.
$$
\end{enumerate}
Two possible choices for the set $S_T$, denoted by $S_T^3$ and $S_T^6$,
are illustrated on Figure \ref{figTriSec}. 
The set $S_T^3$ contains three segments, which join one of the three vertices of $T$ to the midpoint of the opposite edge. It corresponds to the three bisection choices studied in the previous chapters.
The set $S_T^6$ contains a larger choice, namely six segments joining
 a vertex $v_1$ of $T$ to the barycenter $\frac 1 3 v_2+ \frac 2 3 v_3$ of the other vertices $v_2$ and $v_3$.
 Other bisection choices, that lead to non-triangular partitions,
are displayed on Figure \ref{figCvxCut}. 

We typically consider the $L^2$ based decision function : if an element $T$ is 
bisected along a segment $s$, creating two smaller elements $T_s^1$ and $T_s^2$, then 
\be
\label{eqDecisionL2Bi}
d_T(s,f) := e_{T_s^1}(f)^2+e_{T_s^2}(f)^2.
\ee

Starting from a triangulation $\cT_{N_0}$ of the polygonal domain $\Omega$ on which the function $f$ needs to be approximated, with $\#(\cT_{N_0}) = N_0$, the greedy algorithm creates step after step a sequence $\cT_{N_0+1}, \, \cT_{N_0+2}, \, \cdots$ of triangulations of $\Omega$ satisfying $\#(\cT_N) = N$.
At the time of writing, unfortunately, the author does not know how to choose the set of segments $S_T$ and the decision function $d_T(s,f)$ such that given any cartoon function $f$, the sequence $(\cT_N)_{N\geq N_0}$ of 
partitions produced by the greedy algorithm satisfies the desired estimate \iref{estimC}.

As an intermediate objective, we choose to study the behavior of the refinement procedure when $f$ is a particularly simple cartoon function, the characteristic function $f=\chi_P$ when $P$ is a half plane. Up to translation and
rotation, we may always assume that $P$ is of the form
$$
P := \{(x,y)\in \R^2 \sep y\geq 0\}.
$$
We denote by $D$ the line 
\be
\label{defD}
D = \partial P = \{(x,0) \sep x\in \R\},
\ee
and we observe that, denoting by $\mathring T$ the interior of a triangle $T$, one has 
\be
\label{eqErrorZero}
e_T(\chi_P) = 0 \quad \stext{ if and only if } \quad \mathring T \cap D = \emptyset.
\ee

We shall first prove that, for certain configurations between 
the line $D$ and the initial triangulation $\cT_0$,
the refinement procedure based on the three bisection choices $S_T^3$
produces a sequence of triangulations $(\cT_N)_{N\geq 1}$ that {\it does not} satisfy 
the optimal convergence estimate
\be
\label{eqEstimChi}
e_{\cT_N}(\chi_P) \leq CN^{-1}.
\ee
We then show that if the six bisection choices $S_T^6$ illustrated on Figure \iref{figTriSec} are allowed, and if an appropriate decision function is used, then \iref{eqEstimChi} is satisfied. Eventually, we consider
another variant of the refinement procedure, which only involves three bisection choices at each step, and for which \iref{eqEstimChi} is again satisfied.

\subsection{Three bisection choices : a negative result}
\label{secNegRes}


We show in this section that the greedy algorithm studied in the previous chapters beats the error estimate \iref{eqWav} associated to wavelet expansions when $f=\chi_P$. However it fails to achieve the desired error estimate \iref{eqEstimChi}.
More precisely, Proposition \ref{propSpeedBi} shows that, if the decision function $d_T(s,f)$ satisfies Assumptions \ref{assumBi},
then the greedy refinement procedure for the function $f=\chi_P$, starting from 
an arbitrary triangulation $\cT_{N_0}$, $\#(\cT_{N_0}) = N_0$,
produces a sequence $(\cT_N)_{N\geq N_0}$ of triangulations, $\#(\cT_N) = N$, such that 
\be
\label{eqUpperBi}
e_{\cT_N}(\chi_P) \leq C N^{-\lambda/2},
\ee
where $\lambda/2 = 0.732\cdots$. 
This estimate cannot be improved since we establish in Proposition \ref{propNegResult} that there exists a triangle $\TRef$ such that the sequence $(\cT_N)_{N\geq 1}$ produced by the greedy refinement procedure for the function $f=\chi_P$ starting from $\cT_1 = \{\TRef\}$  satisfies 
\be
\label{eqBelowBi}
e_{\cT_N}(\chi_P) \geq c N^{-\lambda/2}
\ee
where $c>0$.

We say that two triangles $T$, $T'$ are $D$-affine equivalent
if there exists an affine change of coordinates $\phi : \R^2\to \R^2$, that transforms $T$ into $T'$ and leaves the line $D$ defined at \iref{defD} invariant :  
$$
\phi(T) = T' \stext{ and } \phi(D) = D.
$$
If $T$ and $T'$ are $D$-affine equivalent, then it follows from a change of variables that 
\be
\label{eqErrorAffEq}
|T|^{-\frac 1 2} e_T(\chi_P) = |T'|^{-\frac 1 2}e_{T'}(\chi_P).
\ee

\begin{figure}
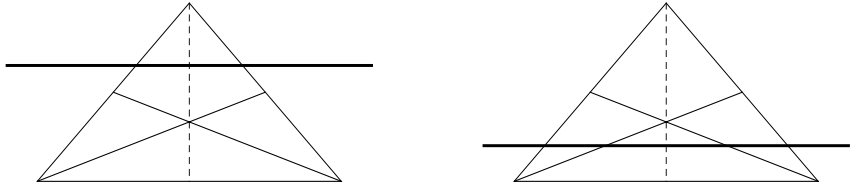

\centering
\includegraphics[height = 2.6cm, width = 5cm]{\pathPic/Subdivision/Cut2_4.pdf}
\hspace{1cm}
\includegraphics[height = 2.6cm, width = 5cm]{\pathPic/Subdivision/Cut2_1.pdf}
\caption{\label{figCut2} Illustrations of Assumptions \ref{assumBi}, a) left and b) right. The line $D$ (Thick), forbidden bisection choices (dashed), authorized bisection choices (full).}
\end{figure}

The following assumptions on the decision function are illustrated on Figure \iref{figCut2}. Their purpose is to create, through the refinement process, a thin layer of triangles along the line $D$.
\begin{assumptions}
\label{assumBi}
Let $T$ be a triangle such that $e_T(\chi_P) \neq 0$. 
In the following we work under the assumption that 
segment $s\in S_T^3$ which minimizes the decision function $d_T(s,\chi_P)$ satisfies the following properties.
\begin{enumerate}[a)]
\item If there is such a possibility in $S_T^3$, the bisection of $T$ along $s$ creates a child $T'$ such that $e_{T'}(\chi_P) = 0$. 
\item If a) is impossible then the decision function selects a segment $s\in S_T^3$ joining a vertex of the triangle $T$ to the midpoint of an edge which is intersected by the line $D$.
\end{enumerate}
\end{assumptions}

These two properties are illustrated in figure \ref{figCut2}. Numerical experiments strongly suggest that they are satisfied for the decision function based on the $L^2$ error which is defined by \iref{eqDecisionL2Bi}.

We first prove the lower error bound \iref{eqBelowBi}, and for that purpose we denote by $\TRef$ the triangle of vertices
$$
(0,-1), (3,-1), (0,2).
$$
$\TRef$ is the large triangle enclosing the others on Figure \ref{figSmallTree}.
The refinement procedure applied to the triangle $\TRef$ and the function $f=\chi_P$ creates an infinite master tree $\cPRef$ of triangles. 
If $\cP$ is a finite subtree of $\cPRef$, (in other words $\cP$ contains $\TRef$ and the nodes of $\cP$ have either two or no children), then the leaves of $\cP$ define a triangulation $\cT$ of $\TRef$. 

\begin{prop}
\label{propNegResult}
There exists a constant $c>0$ such that : 
for any triangulation $\cT$ associated to a finite subtree of $\cPRef$, we have 
\be
\label{eqErrorChiNeg}
e_\cT(\chi_P) \geq c \#(\cT)^{-\lambda/2}
\ee
where $c>0$ is a positive constant and $\lambda = 1.464\cdots$ is the solution of 
\be
\label{defLambdaNeg}
\left(\frac 1 4\right)^{\frac1{1+\lambda}} +  \left(\frac 1 8\right)^{\frac1{1+\lambda}}=1.
\ee
In particular let $(\cT_N)_{N\geq 1}$ be the sequence of triangulations  generated by the greedy algorithm applied to the function $f=\chi_P$ and starting from $\cT_1 = \{\TRef\}$. If the decision function satisfies Assumptions \ref{assumBi}, then for all $N \geq 1$ one has 
$$
e_{\cT_N}(\chi_P) \geq c N^{-\lambda/2}.
$$
\end{prop}

\begin{figure}
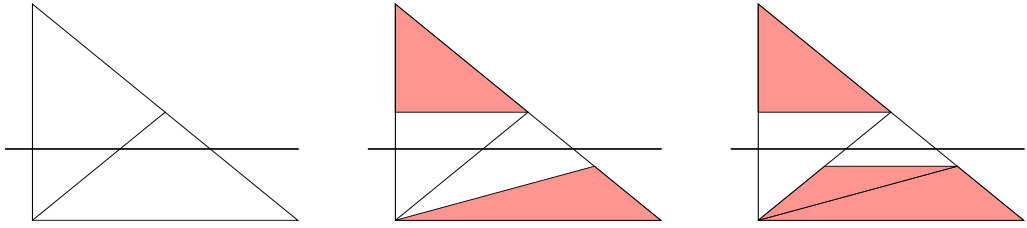

\centering
\includegraphics[height = 3cm, width = 4cm]{\pathPic/Subdivision/Tree1.pdf}
\hspace{0.5cm}
\includegraphics[height = 3cm, width = 4cm]{\pathPic/Subdivision/Tree2.pdf}
\hspace{0.5cm}
\includegraphics[height = 3cm, width = 4cm]{\pathPic/Subdivision/Cut2_3.pdf}
\caption{\label{figSmallTree}Subdivision of the triangle $\TRef$, some of its children and grand-children, with a decision function satisfying Assumptions \ref{assumBi}. The colored triangles satisfy $e_T(\chi_P) = 0$.}
\end{figure}

\begin{figure}
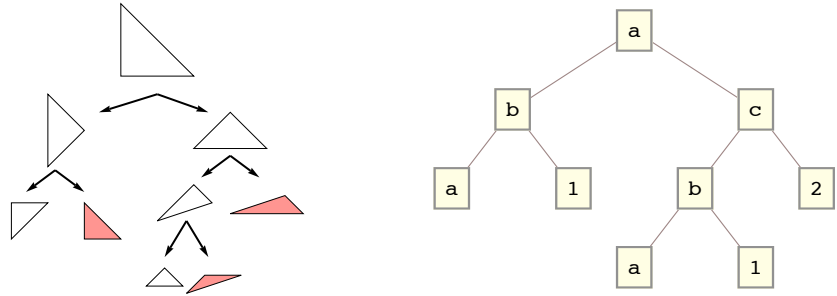

\centering
\includegraphics[height = 4cm, width = 4cm]{\pathPic/Subdivision/TreeTri.pdf}
\hspace{1cm}
\includegraphics[height = 4cm, width = 6cm]{\pathPic/Subdivision/TreeEq.pdf}
\caption{\label{figSmallTree2}Tree associated to the subdivision illustrated on Figure \ref{figSmallTree} (left), and equivalence classes of the elements of this tree for the relation of $D$-affine equivalence (right).}
\end{figure}

\proof
Figures \ref{figSmallTree} and \ref{figSmallTree2} reveal a surprising property of the master tree $\cPRef$. After a few steps, as illustrated on the right of Figure \ref{figSmallTree} we obtain two (white) triangles which are $D$-affine equivalent to the initial triangle $\TRef$, and three (colored) triangles on which the approximation error $e_T(\chi_P)$ is zero. It follows that, up to the relation of $D$-affine equivalence, the tree $\cPRef$ is \emph{self-similar} and is a repetition of the pattern presented on the right of Figure \ref{figSmallTree2}. 

We define an auxiliary tree $\cPRef'$ as follows. The nodes of $\cPRef'$ are the triangles $T_0\in \cPRef$ which are affine equivalent to $\TRef$ with respect to $D$. The children $T_1,T_2$ of a triangle $T_0$ in $\cPRef'$ are the grand-child and grand-grand-child of $T_0$ in $\cPRef$ which are affine equivalent to $\TRef$. Note that 
\be
\label{eqRatioArea}
|T_0| = 4 |T_1| = 8 |T_2|.
\ee
Let $\cT$ be a triangulation of $\TRef$ associated to a finite subtree $\cP$ of $\cPRef$. We denote by $\cP'$ the smallest subtree of $\cPRef'$ such that any leaf of $\cP'$ is contained in a leaf of $\cP$. 
We denote by $\cT'$ the collection of leaves of $\cP'$. For instance, if $\cT$ is any of the three triangulations illustrated on Figure \ref{figSmallTree}, then $\cT'$ consists of the two white triangles in the right of Figure \ref{figSmallTree}. A triangle $T\in \cT$ contains a single triangle $T'\in \cT'$ if $e_T(\chi_P)\neq 0$, and none otherwise. Therefore 
\be
\label{eqCardTNeg}
\#(\cT')\leq \#(\cT).
\ee
Furthermore, since all the elements of $\cT'$ are $D$-affine equivalent to $\TRef$, we have 
\be
\label{eqErrorTNeg}
e_\cT(\chi_P)^2 = \sum_{T\in \cT} e_T(\chi_P)^2 \geq \sum_{T\in \cT'} e_T(\chi_P)^2 = c_0 \sum_{T\in \cT'} |T|
\ee
where $c_0=e_{\TRef}(\chi_P)^2/|\TRef|$. 
For any triangle $T_0\in \cP'$, we denote by $\cE(T_0)$ the collection of the elements of $\cT'$ that it contains, 
$$
\cE(T_0) := \{T\in \cT' \sep T \subset T_0\},
$$
and we define 
$$
a(T_0) := \sum_{T\in \cE(T_0)} |T| \stext{ and } n(T_0) := \# (\cE(T_0)).
$$
We establish below that for any $T_0\in \cP'$, the following inequality holds
\be
\label{eqInducTree}
|T_0| \leq a(T_0) \, n(T_0)^\lambda.
\ee
Once this property is established, then choosing $T_0 = \TRef$ and combining this inequality 
with \iref{eqCardTNeg} and \iref{eqErrorTNeg} yields 
\begin{eqnarray*}
e_\cT(\chi_P)^2 &\geq&  c_0 \sum_{T\in \cT'} |T|\\
&=& c_0 \, a(\TRef)\\
&\geq & c_0 \, |\TRef| \, n(\TRef)^{-\lambda}\\
&=& c^2 \, \#(\cT')^{-\lambda}\\
&\geq & c^2 \, \#(\cT)^{-\lambda},
\end{eqnarray*}
where $c^2=c_0 \, |\TRef|=e_{\TRef}(\chi_P)^2$.
This establishes the announced result \iref{eqErrorChiNeg} and concludes the proof of this proposition.
In order to prove \iref{eqInducTree} we use an induction argument on the tree $\cP'$, and for that purpose we distinguish two cases.
\begin{enumerate}[i)]
\item If $T_0$ is a leaf of $\cP'$, then $T_0\in \cT'$ and therefore $a(T_0) = |T_0|$ and $n(T_0) =1$, hence \iref{eqInducTree} holds.
\item Otherwise $T_0$ has two children $T_1$ and $T_2$ in $\cP'$, and we may assume as an induction hypothesis that 
\be
\label{eqHR}
|T_1| \leq a(T_1) \, n(T_1)^\lambda \stext{ and } |T_2| \leq a(T_2) \, n(T_2)^\lambda.
\ee
Applying the holder inequality 
$$
u_1v_1+ u_2v_2 \leq (u_1^p+u_2^p)^{\frac 1 p} (v_1^q+ v_2^q)^{\frac 1 q}
$$ 
to the numbers 
$$
u_i = a(T_i)^{\frac 1 p}, \quad v_i = n(T_i)^{\frac 1 q}, \quad  p=1+\lambda, \quad q= \frac{1+\lambda}\lambda = \frac p \lambda,
$$
and elevating to the power $p$, we obtain 
\be
\label{eqHolderTree}
\left(\left(a(T_1) \, n(T_1)^{\lambda}\right)^{\frac 1 {p}} + \left(a(T_2) \, n(T_2)^{\lambda}\right)^{\frac 1 {p}} \right)^{p} \leq (a(T_1)+a(T_2)) \left(n(T_1)+ n(T_2)\right)^\lambda.
\ee
If follows that 
\begin{eqnarray*}
|T_0| &=& \left(\left(\frac {|T_0|} 4\right)^{\frac1{p}} +  \left(\frac {|T_0|} 8\right)^{\frac1{p}}\right)^p\\
&=& \left(|T_1|^{\frac 1 {p}}+ |T_2|^{\frac 1 {p}} \right)^{p}\\
&\leq & \left(\left(a(T_1) \, n(T_1)^{\lambda}\right)^{\frac 1 {p}} + \left(a(T_2) \, n(T_2)^{\lambda}\right)^{\frac 1 {p}} \right)^{p} \\
&\leq & (a(T_1)+a(T_2)) \left(n(T_1)+ n(T_2)\right)^\lambda.\\
&=& a(T_0) \, n(T_0)^\lambda,
\end{eqnarray*}
where we use successively \iref{defLambdaNeg}, \iref{eqRatioArea}, \iref{eqHR}, \iref{eqHolderTree} and the equalities $a(T_0) = a(T_1)+a(T_2)$ and $n(T_0)=n(T_1)+n(T_2)$. Again \iref{eqInducTree} holds, and by induction it holds for any element of $\cP'$.
\sq
\end{enumerate}

\begin{figure}
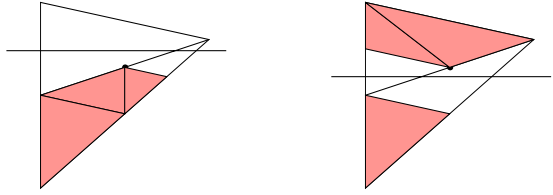

\centering
\includegraphics[height = 2.6cm, width = 3cm]{\pathPic/Subdivision/lemma2_1.pdf}
\hspace{1cm}
\includegraphics[height = 2.6cm, width = 3cm]{\pathPic/Subdivision/lemma2_2.pdf}
\caption{\label{figRefine2}Bisection of a triangle, some of its children, grand-children and grand-grand-children, with decision function satisfying Assumptions \ref{assumBi}}
\end{figure}

We now turn to the proof of the upper estimate \iref{eqUpperBi}, and for that purpose we begin with a geometrical lemma which uses the two assumptions \iref{assumBi} to analyse the first steps of the refinement procedure.

\begin{lemma}
\label{lemmaSmallRef2}
Let $T$ be a triangle such that $e_T(\chi_P)\neq 0$.
If the decision function satisfies \iref{assumBi}, then one of the following holds
\begin{enumerate}[i)]
\item A child $T'$ produced by the bisection of $T$ satisfies $e_{T'}(\chi_P) = 0$.
\item Applying iteratively the refinement procedure to $T$, some of its children, grand-children and grand-grand-children we obtain a triangulation $\cT$ of $T$ such that
\begin{itemize} 
\item $\#(\cT)\leq 5$, and $|T'| \geq \frac 1 {2^4} |T|$ for all $T'\in \cT$.
\item The approximation error $e_{T'}(\chi_P)$ is non-zero for at most two triangles $T'\in \cT$. Denoting these triangles by $T_1,T_2$, one of the following two possibilities holds 
$$
|T| = 2 |T_1| = 16|T_2| \stext{ or } |T| = 4 |T_1| = 8|T_2| .
$$
\end{itemize}
\end{enumerate}
\end{lemma}

\proof
Rather than reading a long discussion, the author invites the reader to check that, if i) is impossible, then one of the two possibilities illustrated on Figure \ref{figRefine2} occurs and satisfies ii).
\sq

\begin{prop}
\label{propSpeedBi}
Let $\cT_{N_0}$ be an arbitrary triangulation in $\R^2$ and let $(\cT_N)_{N\geq N_0}$ be the sequence of triangulations generated by the greedy algorithm applied to the function $f=\chi_P$ and starting from the triangulation $\cT_0$. If the decision function satisfies Assumptions \ref{assumBi}, then 
$$
e_{\cT_N} (\chi_P) \leq C N^{-\lambda/2},
$$
where $\lambda = 1.464\cdots$ is defined by \iref{defLambdaNeg}.
\end{prop}

\proof
A proof of this proposition can be obtained by a straightforward adaptation of Lemma \ref{lemmaTEpsMu} and Proposition \ref{propSpeedTri} in the next section, which is left to the reader.
\sq

\subsection{Six bisection choices}
\label{secTriSec}
The previous section contains a negative result, Proposition \ref{propNegResult}, showing that the greedy refinement procedure does not achieve the desired convergence rate \iref{eqEstimChi} when it is based on the three bisection choices $S_T^3$ and a decision function satisfying the natural assumptions \ref{assumBi}.

In this section the segment along which a triangle $T$ is bisected is picked among the six possibilities in $S_T^6$, illustrated on Figure \ref{figTriSec}, instead of the three possibilities in $S_T^3$. If the decision function satisfies Assumptions \ref{assumTri}, then the desired convergence rate \iref{eqEstimChi} is achieved, see Proposition \ref{propSpeedTri}. This result can be seen as a first step towards the construction of a hierarchical and anisotropic refinement algorithm well adapted to the approximation of cartoon functions.


The following assumptions on the decision function $d_T(s,f)$ are illustrated on Figure \ref{figDecision3}. Their purpose is to create, through the refinement process, a thin layer of triangles along the line $D$. Note that  for any triangle $T$ any for any segment $s\in S_T^6$, the two children of $T$ have areas $\frac 1 3 |T|$ and $\frac 2 3 |T|$.
\begin{assumptions}
\label{assumTri}
Let $T$ be a triangle such that $e_T(\chi_P) \neq 0$. 
In the following we work under the assumption that 
segment $s\in S_T^6$ 
which minimizes the decision function $d_T(s,\chi_P)$ satisfies the following properties.
%
\begin{enumerate}[a)]
\item If there is such a possibility in $S_T^6$, the bisection of $T$ along $s$ creates a child $T'$ such that $e_{T'}(\chi_P) = 0$ and $|T'| = \frac 2 3 |T|$.
\item If a) is impossible, then whenever there is such a possibility in $S_T^6$, the bisection of $T$ along $s$ creates a child $T'$ such that $e_{T'}(\chi_P) = 0$ and $|T'| = \frac 1 3 |T|$.
\item If a) and b) are impossible then the decision function selects a segment $s = [v_1,b] \in S_T^6$ joining a vertex $v_1$ of the triangle $T$ to the barycenter $b=\frac 1 3 v_2+ \frac 2 3 v_3$ of the other vertices. The choice of $s$ must be such that the line $D$ intersects the segment $[b,v_3]$.
\end{enumerate}
\end{assumptions}

\begin{figure}
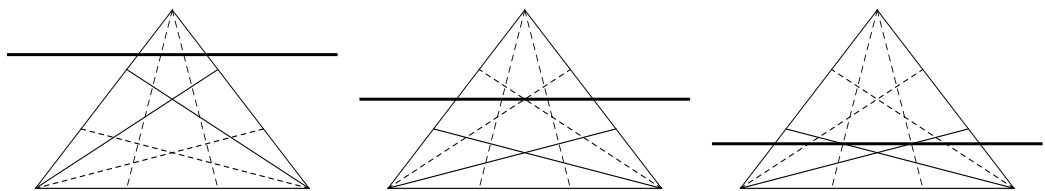

\centering
\includegraphics[height = 2.6cm, width = 4.5cm]{\pathPic/Subdivision/Cut3_6.pdf}
\includegraphics[height = 2.6cm, width = 4.5cm]{\pathPic/Subdivision/Cut3_4.pdf}
\includegraphics[height = 2.6cm, width = 4.5cm]{\pathPic/Subdivision/Cut3_5.pdf}
\caption{\label{figDecision3} Illustrations of Assumptions \ref{assumTri}, a) left, b) center and c) right. The line $D$ (Thick), forbidden bisection choices (dashed), authorized bisection choices (full).}
\end{figure}

The next lemma uses these three assumptions to analyse the first steps of the refinement procedure on a triangle $T$.

\begin{lemma}
\label{lemmaSmallRef3}
Let $T$ be a triangle such that $e_T(\chi_P)\neq 0$.
If the decision function satisfies \iref{assumTri}, then one of the following holds
\begin{enumerate}[i)]
\item A child $T'$ produced by the bisection of $T$ satisfies $e_{T'}(\chi_P) = 0$.
\item Applying iteratively the refinement procedure to $T$, some of its children, grand-children and grand-grand-children we obtain a triangulation $\cT$ of $T$ such that
\begin{itemize} 
\item $\#(\cT)\leq 7$, and $|T'| \geq \frac 1 {3^4} |T|$ for all $T'\in \cT$.
\item The approximation error $e_{T'}(\chi_P)$ is non-zero for at most two triangles $T'\in \cT$. Denoting these by $T_1,T_2$, there exists $i\in \{0,1,2\}$ such that 
\be
\label{eqRatioTChild}
\frac {|T_1|} {|T|} \leq \frac 1 3 \times \left(\frac {3-i} 3\right)^2 \stext{ and } \frac {|T_2|}{|T|} \leq \frac 2 9 \times \left(\frac {i+1} 3\right)^2.
\ee
\end{itemize}
\end{enumerate}
\end{lemma}

\proof
Rather than reading a long discussion, the author invites the reader to check that, if i) is impossible, then one of the three possibilities illustrated on Figure \ref{figRefine3} occurs and satisfies ii).
\sq

Let $T_r$ be a triangle such that $e_{T_r}(\chi_P)\neq 0$. The decision function applied to $T_r$ and its descendants defines a master tree $\cP_r$ of triangles.
If $\cP$ is a finite subtree of $\cP_r$, (in other words $\cP$ contains $T_r$, and the nodes of $\cP$ have either two or no children), then the leaves of $\cP$ define a triangulation $\cT$ of $T_r$. The next lemma describes some of these triangulations.

\begin{lemma}
\label{lemmaTEpsMu}
For any $0<\ve\leq 1$, there exists a triangulation $\cT_\ve$ of $T_r$ associated to a finite subtree of $\cP_r$, which satisfies
$
\#(\cT_\ve) \leq C \ve^{-\mu}
$
and 
\be
\label{eqTAreaEps}
\max \{ |T| \sep T \in \cT_\ve \text{ and } e_T(\chi_P) \neq \emptyset \} \ \leq  \ve |T_r|
\ee
where $C=30$ and $\mu = 0.329\cdots$ is the solution of 
\be
\label{defMuTri}
\left(\frac 4 {27}\right)^\mu + \left(\frac 8 {81}\right)^\mu = 1.
\ee
As a result, the master tree $\cP_r$ contains at most $C\ve^{-\mu}$ triangles $T$ such that $e_T(\chi_P)^2 \geq \ve |T_r|$.
%
\end{lemma}

\begin{figure}
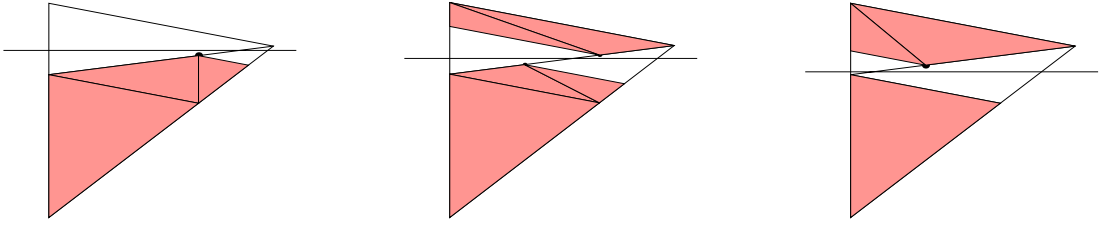

\centering
\includegraphics[height = 3cm, width = 4cm]{\pathPic/Subdivision/Cut3_1.pdf}
\hspace{1cm}
\includegraphics[height = 3cm, width = 4cm]{\pathPic/Subdivision/Cut3_2.pdf}
\hspace{1cm}
\includegraphics[height = 3cm, width = 4cm]{\pathPic/Subdivision/Cut3_3.pdf}
\caption{\label{figRefine3}Bisection of a triangle, some of its children, grand-children and grand-grand-children, with a decision function satisfying Assumptions \ref{assumTri}.}
\end{figure}

\proof
We consider a second tree $\cP_r'$ of triangles, with root $T_r$. 
A triangle $T\in \cP_r'$ such that $e_T(\chi_P) = 0$ has no children. A triangle $T\in \cP_r'$ such that  $e_T(\chi_P)\neq 0$ has at most $7$ children, those described in Lemma \ref{lemmaSmallRef3}. Two cases are possible : 
\begin{enumerate}[i)]
\item The triangle $T$ has two children and one of them, denoted by $T'$, satisfies $e_{T'}(\chi_P)=0$. The child $T'$ of $T$ is therefore a leaf of the tree $\cP_r'$. 
\item The triangle $T$ has at most $7$ children, and at most two of them $T_1,T_2$ are not leaves of the tree $\cP_r'$.  Furthermore the areas of $T_1,T_2$ satisfy \iref{eqRatioTChild} for some $i\in \{0,1,2\}$. 
\end{enumerate}

We denote by $\cP'_\ve$ the subtree of $\cP'_r$ created as follows : we start from the root $T_r$, and we include the children of a triangle $T$ if and only if $e_T(\chi_P) \neq 0$ and $|T|>\ve|T_r|$.
We define $\cT_\ve$ as the collection of leaves of $\cP'_\ve$ and we remark that \iref{eqTAreaEps} is satisfied. 
Note that the parent $T\in \cP'_\ve$ of any triangle $T'\in \cP'_\ve$ satisfies $|T| > \ve |T_r|$. Lemma \ref{lemmaSmallRef3} therefore implies that 
\be
\label{eqAreaBelow}
|T'|\geq \frac  {\ve |T_r|}  {3^4} \stext{ for any } T'\in \cP'_\ve.
\ee
We associate to each triangle $T_0\in \cP'_\ve$ the number 
$$
n(T_0) := \# \{T\in \cT_\ve \sep T\subset T_0\},
$$
and we intend to show that for any triangle $T_0$ in the tree $\cP'_\ve$ one has
\be
\label{eqRecTree3}
n(T_0) +5 \leq C \left(\frac {|T_0|}{\ve |T_r|}\right)^\mu
\ee
Once this point is established, choosing $T_0 =T_r$ concludes the proof of the first part of this proposition.
For that purpose, we proceed by induction on the tree $\cP'_\ve$. We first remark that 
for any leaf $T_0$ of $\cP'_\ve$ is follows from \iref{eqAreaBelow} that 
$$
C \left(\frac {|T_0|}{\ve |T_r|}\right)^\mu \geq C \left(\frac 1 {3^4}\right)^\mu \geq 6 = n(T_0)+5.
$$

If $T_0$ is not a leaf of $\cP'_\ve$, then 
we may assume as an induction hypothesis that \iref{eqRecTree3} holds for all the children of $T_0$. We now distinguish between the two types of nodes i) and ii) of the tree $\cP'_\ve$, and we obtain the following.
\begin{enumerate}[{Type} i)]
\item The triangle $T_0$ has two children, $T_1$, $T_2$, where $T_2$ is a leaf of the tree $\cP'_\ve$. Since $T_0$ is not a leaf of $\cP'_\ve$, it satisfies $|T_0|\geq \ve |T_r|$. Furthermore $|T_1| \leq \frac 2 3 |T_0|$, 
hence  
\begin{eqnarray*}
C \left(\frac {|T_0|}{\ve |T_r|}\right)^\mu &=& C \left(\frac 2 3 \times \frac { |T_0|}{\ve |T_r|}\right)^\mu + C\left(1-\left(\frac 2 3\right)^\mu\right) \left( \frac { |T_0|}{\ve |T_r|}\right)^\mu\\ 
&\geq & C  \left( \frac {|T_1|}{\ve |T_r|}\right)^\mu + C \left(1-\left(\frac 2 3\right)^\mu\right)\\
&\geq & n(T_1)+5 +1 \\
&=& n(T_0)+5.
\end{eqnarray*}
Therefore $T_0$ satisfies \iref{eqRecTree3}.
\item
The triangle $T_0$ has at most seven children, and at most two of them, denoted by $T_1$, $T_2$, are not leaves. Hence $n(T_0) \leq n(T_1)+n(T_2)+5$. 
Using the estimate \iref{eqRatioTChild} on the areas of $T_1$ and $T_2$, and the definition \iref{defMuTri} of $\mu$, we obtain
$$
\left(\frac{|T_1|}{|T_0|}\right)^\mu+ \left(\frac{|T_2|}{|T_0|}\right)^\mu \leq 1.
$$
Therefore 
\begin{eqnarray*}
C  \left(\frac {|T_0|}{\ve |T_r|}\right)^\mu &\geq & C  \left(\frac {|T_1|}{\ve |T_r|}\right)^\mu + C  \left(\frac {|T_2|}{\ve |T_r|}\right)^{\lambda}  \\
&\geq & n(T_1)+5+ n(T_2)+5\\
& \geq & n(T_0)+5,
\end{eqnarray*}
which concludes the proof of \iref{eqRecTree3}. 
\end{enumerate}
We now turn to the second part of the proposition, and for that purpose we denote by $\cP_\ve$ the subtree of $\cP_r$ associated to the triangulation $\cT_\ve$. Note that this tree is distinct from the subtree $\cP'_\ve$ of $\cP_r'$. In particular $\cP_\ve$ is a binary tree while $\cP'_\ve$ is not (the tree $\cP'_\ve$ can be obtained by grouping some of the elements of $\cP_\ve$ in a single node with some of their descendants).

Since $\cP_\ve$ is a binary tree, the cardinality of the set $\cP_\ve\sm \cT_\ve$ of its inner nodes has a simple expression
$$
\#(\cP_\ve \sm \cT_\ve) = \#(\cT_\ve)-1.
$$
Furthermore, any triangle $T \in \cP_r$ which does not belong to $\cP_\ve\sm \cT_\ve$ is an element of $\cT_\ve$ or one of its descendants, hence $T$ satisfies $e_T(\chi_P) = 0$ or $e_T(\chi_P)^2\leq |T|< \ve |T_r|$.
Therefore
$$
\#(\{T \in \cP_r \sep e_T(\chi_P)^2 \geq \ve |T_r|\}) \leq \#(\cP_\ve\sm \cT_\ve) = \#(\cT_\ve)-1\leq C \ve^{-\mu},
$$
which concludes the proof.
\sq

\begin{prop}
\label{propSpeedTri}
Let $\cT_{N_0}$ be an arbitrary triangulation in $\R^2$ and let $(\cT_N)_{N\geq N_0}$ be the sequence of triangulations generated by the greedy algorithm applied to the function $f=\chi_P$ and starting from the triangulation $\cT_0$. If the decision function satisfies the assumptions \ref{assumTri}, then 
\be
\label{eqErrorNu}
e_{\cT_N} (\chi_P) \leq C N^{-\nu},
\ee
where $\nu = \frac 1 {2 \mu} -\frac 1 2= 1.0188\cdots$. 
\end{prop}

\proof
Let $T_0$ be a triangle, let $f\in L^2(T)$ and let $T_1$, $T_2$ be the children produced by some bisection of $T$. Then
\begin{eqnarray*}
e_{T_1}(f)^2+e_{T_2}(f)^2 &=& \inf_{\pi_1, \pi_2\in \sP_1} \|f-\pi_1\|_{L^2(T_1)}^2+\|f-\pi_2\|_{L^2(T_2)}^2\\
&\leq & \inf_{\pi\in \sP_1} \|f-\pi\|_{L^2(T_1)}^2+\|f-\pi\|_{L^2(T_2)}^2\\
& = & e_{T_0}(f)^2.
\end{eqnarray*}
Hence 
$$
\max\{e_{T_1}(f),e_{T_2}(f)\} \leq e_{T_0}(f).
$$
The decision function applied to the triangles in $\cT_0$ and their descendants defines a collection $\cP_0$ of $N_0$ master trees. It follows from the previous inequality that the approximation error $e_T(\chi_P)$ is larger on a triangle $T$ than on any of its descendants.

At each step $N \geq N_0$, the greedy algorithm selects for bisection a triangle in $T_N \in \cT_N$ which maximizes the approximation error $e_T(\chi_P)$ among all $T \in \cT_N$.
As a result we have for 
any $N\geq N_0$,
$$
e_{T_{N_0}}(\chi_P) \geq e_{T_{N_0+1}}(\chi_P) \geq \cdots \geq e_{T_N}(\chi_P) \geq \max \left\{e_T(\chi_P) \sep T \in \cP_0 \sm \{T_{N_0}, \cdots, T_N\}\right\}.
$$
We thus have for any $N\geq N_0$, since $\#(\cT_N)=N$,
\be
\label{eqErrorGreedyDec}
e_{\cT_N}(\chi_P)^2 = \sum_{T\in \cT_N} e_T(\chi_P)^2 \leq N \ e_{T_N} (\chi_P)^2.
\ee
According to Lemma \ref{lemmaTEpsMu}, for any $0<\ve \leq 1$ the collection $\cP_0$ of $N_0$ trees contains at most $C_0 \ve^{-\mu}$ triangles $T$ such that $e_T(\chi_P)^2 \geq \ve$, where $C_0$ is independent of $\ve$. Choosing $\ve = C_0^{\frac 1 \mu} (N-N_0+1)^{-\frac 1 \mu}$ we obtain  
$$
e_{T_N}(\chi_P)^2 \leq C_0^{\frac 1 \mu} (N-N_0+1)^{-\frac 1 \mu}.
$$
Combining this with \iref{eqErrorGreedyDec} yields \iref{eqErrorNu} which concludes the proof.
\sq

\begin{remark}
\label{remarkExcessiveSpeed}
The convergence rate obtained in Proposition \ref{propSpeedTri} actually exceeds the desired error estimate \iref{eqEstimChi}. This point is due to the fact that  $f=\chi_P$ is a particularly simple cartoon function.
For a function cartoon $f$ which is discontinuous along curved edges, instead of a straight line, or which is not piecewise affine, the convergence rate \iref{eqEstimChi} is optimal.
\end{remark}

\subsection{Partitions of a domain into convex polygons} 
\label{secConvexCut}

We suggest in this section another modification of the refinement procedure, in an attempt to improve the result of Proposition \ref{propSpeedTri}. This proposition is indeed not completely satisfactory, in particular because of the following three reasons.
\begin{itemize}
\item The decision function $d_T(s,f)$ needs to satisfy Assumptions \ref{assumTri}. Unfortunately numerical results  suggest that the natural decision function \iref{eqDecisionL2Bi} does not satisfy these assumptions, and the author has not found a decision function with a simple expression which satisfies these assumptions.
\item The number of segments $s\in S_T^6$ along which a triangle $T$ may be bisected has increased from $3$ to $6$. 
From the point of view of data compression, the hierarchical bisection tree is therefore more expensive to encode. From the theoretical point of view, the algorithm could be regarded as less elegant. 
\item Since the cartoon function $f=\chi_P$ is particularly simple, see Remark \ref{remarkExcessiveSpeed}, a convergence rate faster that \iref{eqErrorNu} could be expected. For example an exponential convergence rate $e_{\cT_N}(\chi_P)\leq C e^{-\beta N}$ where $\beta>0$.
\end{itemize}

The modification of the algorithm described below solves these three points, at the following cost. The modified greedy algorithm does not produce triangulations any more, but partitions into convex polygons of the domain $\Omega$ on which the function $f$ to be approximated is defined. 

We denote by $\cC$ the collection of bounded closed and convex sets of nonempty interior included in $\R^2$. We denote by the capital letter $T$ the elements of $\cC$, whether they are triangles or not, and by $\cT$ any partition of a domain $\Omega$ into elements of $\cC$. 
The definitions \iref{defErrTL2} and \iref{defErrCTL2} of the approximation errors $e_T(f)$ and $e_\cT(f)$ are unchanged.
%
%
%


Let $T\in \cC$ be a bounded closed and convex set, and let $a,b,c\in \partial T$ be three distinct points on the boundary of $T$.
If none of the segments $[a,b]$, $[b,c]$ or $[c,a]$ is included in $\partial T$, then the complementary in $T$ of the triangle of vertices $a,b,c$ has three connected components.
We denote by $T_a$ the closure of the connected component which does not contain $a$, and  
we define $T_b$ and $T_c$ likewise, see Figure \ref{figCvxCut}. 
The bisection of $T$ along the segment $[b,c]$, $[c,a]$ or $[a,b]$ produces a pair of children $\{T_a, \overline {T\sm T_a}\}$,  $\{T_b, \overline {T\sm T_b}\}$ or $\{T_c, \overline {T\sm T_c}\}$ respectively which also belong to $\cC$. 
We define 
$$
\tau(T \ssep a,b,c) := \min \{|T_a|, |T_b|, |T_c|\}
$$
and we set the convention $\tau(T \ssep a,b,c) = 0$ for any $a,b,c\in \partial C$ such that the complementary in $T$ of the triangle of vertices $a,b,c$ has only two connected components. We also define
\be
\label{defTauT}
\tau(T) := \sup_{a,b,c\in \partial T} \tau(T \ssep a,b,c).
\ee
The next proposition gives a lower bound for the ratio $\tau(T)/|T|$.

\begin{figure}
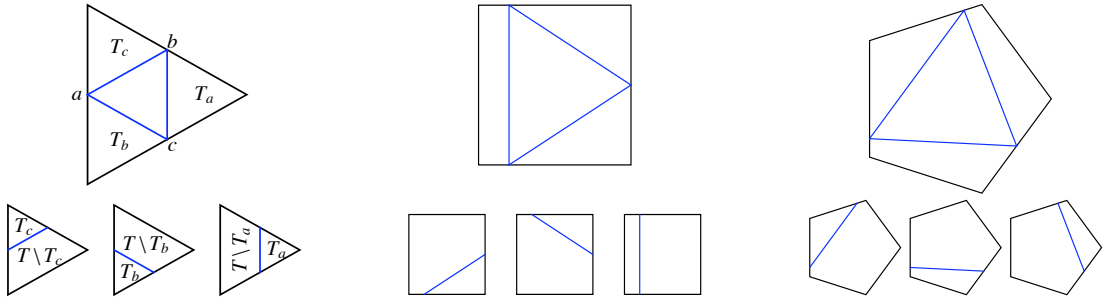

\centering
\includegraphics[width = 4cm, height=4cm]{\pathPic/Subdivision/Cvx1Names.pdf}
\hspace{1cm}
\includegraphics[width = 4cm, height=4cm]{\pathPic/Subdivision/Cvx2.pdf}
\hspace{1cm}
\includegraphics[width = 4cm, height=4cm]{\pathPic/Subdivision/Cvx3.pdf}
\caption{\label{figCvxCut} Bisection of a convex $T\in \cC$ along the set of segments $S_T^{3*}$ determined by three points $a,b,c\in \partial T$ which maximize \iref{defTauT}.}
\end{figure}

\begin{prop}
\label{propConvex}
The supremum \iref{defTauT} is attained and for any $a,b,c$ satisfying $\tau(T\ssep a,b,c) = \tau(T)$, we have $|T_a| = |T_b| = |T_c| = \tau(T)$.
Furthermore for any $T \in \cC$, one has 
\be
\label{eqDeltaTC}
|T| \leq 7\tau(T). 
\ee
\end{prop}

\proof
For any fixed $T\in \cC$, the boundary $\partial T$ is compact and the map
$$
(a,b,c)\in (\partial T)^3 \mapsto \tau(T\ssep a,b,c)
$$
is continuous. Hence the maximum \iref{defTauT} is attained. Furthermore let $a,b,c \in \partial T$ and let us assume that $|T_a| > |T_b|$, then moving the point $c\in \partial T$ in the direction of $T_a$ continuously decreases $|T_a|$ and increases $|T_b|$. As a result, if the triplet $(a,b,c)\in (\partial T)^3$ maximizes \iref{defTauT}, then one must have $|T_a| = |T_b|=|T_c|$.

For any affine change of coordinates $\Phi(z) := \phi z+ z_0$, where $\phi\in \GL_2$ and $z_0\in \R^2$, one easily checks that 
\be
\label{eqChgTau}
\frac{\tau(T)}{|T|} = \frac{\tau(\Phi(T))} {|\Phi(T)|}.
\ee
Let $T\in \cC$ and let $(a,b,c)\in (\partial T)^3$ be a triplet maximizing \iref{defTauT}. In order to establish \iref{eqDeltaTC} we may assume using \iref{eqChgTau}, up to an affine change of coordinates, that the points $a,b,c$ are the vertices of an equilateral triangle $\TEq$ centered at the origin $O = (0,0)\in \R^2$ and of edge-length $2 = |a-b| = |b-c| = |c-a|$. See Figure \ref{figChildCvx} for an illustration of this triangle and of the notations defined below.

\begin{figure}
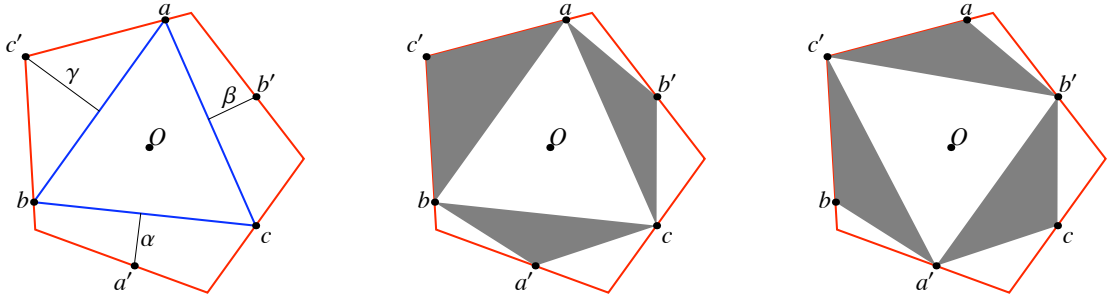
	
	\centering
		\includegraphics[width = 4cm, height=4cm]{\pathPic/Subdivision/ChildCvx1.pdf}
		\hspace{1cm}
		\includegraphics[width = 4cm, height=4cm]{\pathPic/Subdivision/ChildCvx2.pdf}
		\hspace{1cm}
		\includegraphics[width = 4cm, height=4cm]{\pathPic/Subdivision/ChildCvx3.pdf}
		\caption{\label{figChildCvx}
		Left : notations of Proposition \ref{propConvex}, Center : three triangles contained in the convex sets $T_a$, $T_b$ and $T_c$, Right : three triangles contained in $T_{a'}$, $T_{b'}$ and $T_{c'}$.}
\end{figure}

We denote by $a'$ the point on $\partial T$ such that $O\in [a,a']$, and we denote by $\alpha$ the distance from $a'$ to the segment $[b,c]$. 
Since the convex set $T_a$ contains the triangle of vertices $a',b,c$, we have 
\be
\label{eqLowerTA}
 |T_a|\geq \frac {|b-c| \,\alpha} 2 = \alpha.
\ee
We define $b',c'$ and $\beta, \gamma$ similarly. The triplet $(a',b',c')\in (\partial T)^3$ defines three convex sets $T_{a'}, T_{b'}, T_{c'}$. Recalling that 
\begin{eqnarray*}
\tau(T) &=& \tau(T\ssep a,b,c) =|T_a| = |T_b| = |T_c| \\
&\geq &\tau(T\ssep a', b',c') = \min\{|T_{a'}|, \, |T_{b'}|, \, |T_{c'}| \},
\end{eqnarray*}
we obtain 
\be
\label{eqLowerTauT}
\tau(T) \geq \max\{|T_a|,\, |T_b|, \, |T_c|, \ \min\{|T_{a'}|, \, |T_{b'}|, \, |T_{c'}| \} \}\\
\ee
Since $|T_{a'}|$ contains the triangle of vertices $a,b',c',$ an elementary computation shows that 
$$
|T_a'| \geq \frac  {|\det (b'-a, \,  c'-a)|} 2 = \frac {|\sqrt 3 (1-\beta \gamma) + \beta+ \gamma|} 4,
$$
hence 
$$
\tau(T) \geq 
\max\left\{ \alpha, \, \beta, \, \gamma, \, \frac {\sqrt 3} 4 + \frac 1 4 \min \{
\beta+\gamma-\sqrt 3\beta\gamma, \, 
\gamma+\alpha-\sqrt 3\gamma\alpha,\, 
\alpha +\beta-\sqrt 3\alpha\beta
\} \right\}.
$$
If the minimum appearing in the previous equation is negative, then $\alpha$, $\beta$ or $\gamma$ is larger than $1/\sqrt 3 > \sqrt 3/4$. It follows that $\tau(T)\geq \sqrt 3/4$, which concludes the proof since $|T| = |\TEq|+3\tau(T) = \sqrt 3+ 3\tau(T)$.
\sq

We denote by $S_T^{3*} := \{[a,b],[b,c],[c,a]\}$ the collection of segments joining the three points $a,b,c\in (\partial T)^3$ which maximize \iref{defTauT}. Some convex polygons and the associated optimal triplets $(a,b,c)\in (\partial T)^3$ are illustrated on Figure \ref{figCvxCut}. If there are several triplets $(a,b,c)$ realizing this maximum, for instance if $T$ is a disc or a regular polygon, then we consider the minimal triplet with respect to the lexicographic order on $(\R^2)^3$.

From an algorithmic point of view, for any convex polygon $T\in \cC$ the points $(a,b,c)$ which maximize \iref{defTauT} can be found by solving a finite number of third degree equations (one for each triplet of distinct faces of $T$). Note that Proposition \ref{propConvex} is not optimal. We conjecture that the minimal value $\tau(T)/|T|$ is attained when the convex set $T$ is the unit disc, and is equal to $(\cPi/3-\sqrt 3/4)/\cPi \simeq 0.195\cdots$, which is stricly larger than $1/7\simeq 0.142\cdots$.

\begin{assumptions}
\label{assumConv}
Let $T\in \cC$ be such that $e_T(\chi_P)\neq 0$. 
In the following we work under the assumption that the segment $s\in S_T^{3*}$ which minimizes the decision function $d_T(s,\chi_P)$ satisfies the following property. 
\begin{enumerate}[a)]
\item The bisection of $T$ along $s$ creates a child $T'$ such that $e_{T'}(\chi_P) = 0$. (This is always possible).
\end{enumerate}
\end{assumptions}

\begin{figure}
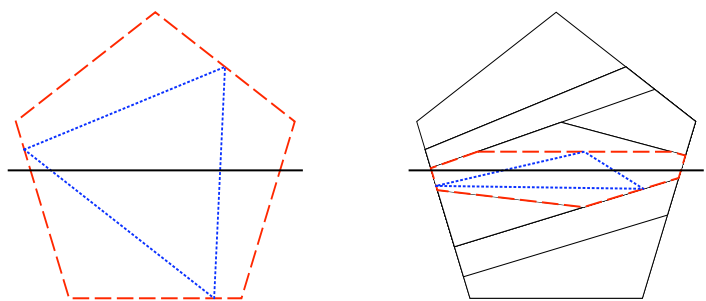
	
	\centering
		\includegraphics[width = 4cm, height=4cm]{\pathPic/Subdivision/Numerical/CutsLineHoriz.pdf}
		\hspace{1cm}
		\includegraphics[width = 4cm, height=4cm]{\pathPic/Subdivision/Numerical/CutsLine8Horiz.pdf}
		\caption{\label{figExpApprox}
		Left : A convex set (red, dashed), the segments $S_T^{3*}$ (blue, dotted), the line $D$ (black, thick).
		Right : After $8$ steps of the greedy algorithm, with a decision function satisfying Assumptions \ref{assumConv}.}
\end{figure}

The following corollary shows that, under the previous assumption on the decision function, the refinement procedure based on these bisection choices creates an exponentially thin layer around the discontinuity of the function $f=\chi_P$, as illustrated on Figure \ref{figExpApprox}.

\begin{corollary}
\label{corolExpApprox}
Let $\Omega\subset \R^2$ be a convex polygon, and let $(\cT_N)_{N\geq 1}$ be the sequence of partitions of $\Omega$ into convex polygons generated by the greedy algorithm applied to the function $f=\chi_P$ and starting from $\cT_1 = \{\Omega\}$. If the decision function satisfies Assumptions \ref{assumConv} then
$$
e_{\cT_N}(\chi_P)^2\leq |\Omega| \, (1-1/7)^{N-1}.
$$ 
\end{corollary}

\proof
It follows from the assumptions \ref{assumConv} that for all $N \geq 1$ at most a single element $T\in \cT_N$ satisfies $e_T(\chi_P)\neq 0$. We denote this element by $T_N$, and we remark that $T_{N+1}$ is a child of $T_N$.
It follows from Proposition \ref{propConvex} that  $|T_{N+1}|\leq (1-1/7) |T_N|$, for any $N \geq 1$, hence 
$$
e_{\cT_N}(\chi_P)^2 = e_{T_N} (\chi_P)^2 \leq |T_N| \leq  |\Omega| \, (1-1/7)^{N-1},
$$
which concludes the proof. See Figure \ref{figExpApprox} for an illustration.
\sq

\noindent
The bisection of $T\in \cC$ along a segment $s\in S_T^{3*}$ creates two children $T_s^1$ and $T_s^2$. In terms of the notations previously used in this section, the bisection along $s=[b,c]\in S_T^{3*}$ creates the two children $T_s^1 = T_a$ and $T_s^2 = \overline{T\sm T_a}$. 
From this point, we assume that the decision function $d_T(s,f)$ is defined by 
\be
\label{decisionConv}
d_T(s,f) := \min\left\{\frac {e_{T_s^1}(f)}{|T_s^1|^{\frac 3 2}}, \ \frac {e_{T_s^2}(f)}{|T_s^2|^{\frac 3 2}}\right\}.
\ee
If the approximation error of $f$ on $T_s^1$ or $T_s^2$ is zero, then $d_T(s,f)=0$. Therefore this decision function satisfies Assumptions \ref{assumConv}, and leads to an exponentially fast convergence if $f=\chi_P$ according to Corollary \ref{corolExpApprox}.

Let us consider a quadratic function $f(z) =\frac 1 2 z^\trans Q z$, where $Q$ is a $2\times 2$ symmetric matrix, to be approximated on a domain $\Omega$. As seen in the previous chapters, the elements of a partition $\cT$ of $\Omega$ should adopt, for an efficient approximation of $f$, an aspect ratio dictated by the hessian matrix $d^2 f=Q$ of $f$ (at least in an average sense).
The numerical experiments presented in Figure \ref{figCvxQuad} suggest that the greedy algorithm, based on the bisection choices $S_T^{3*}$ and the decision function \iref{decisionConv}, generates such partitions. Unfortunately, the author could not prove this property.

\begin{figure}
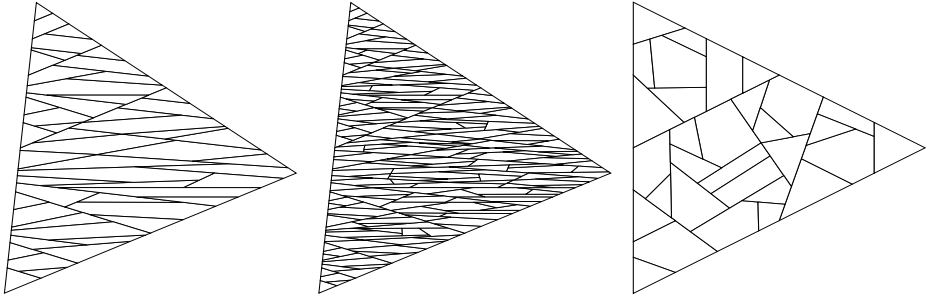

	\centering
		\includegraphics[width=4cm,height=4cm]{\pathPic/Subdivision/Numerical/CvxNum50_100.pdf}
		\includegraphics[width=4cm,height=4cm]{\pathPic/Subdivision/Numerical/CvxNum200_100.pdf}
		\includegraphics[width=4cm,height=4cm]{\pathPic/Subdivision/Numerical/CvxNum30_1.pdf}
		\caption{\label{figCvxQuad} Triangulations produced by the greedy algorithm applied to the function $f(x,y)=x^2+100 y^2$ with the decision function \iref{decisionConv} and 50 (i) and 200 (ii) elements. Likewise for the function $f(x,y) = x^2+y^2$ and 30 elements (iii).}
\end{figure}

Eventually, numerical experiments were also led for the characteristic function $f=\chi_D$ of the unit disc $D = \{(x,y)\in \R^2 \sep x^2+y^2\leq 1\}$. They seem to indicate that the partitions $(\cT_N)_{N \geq 1}$ of the triangle of vertices $(0,0), (2,0),(0,2)$ generated by the greedy algorithm do satisfy the desired error estimate
$$
e_{\cT_N}(\chi_D) \leq CN^{-1},
$$
see Figure \ref{figHorizonQuad}. However, again, the author could not establish this property.

For better numerical precision, the numerical experiments illustrated in Figures \ref{figExpApprox}, \ref{figCvxQuad} and \ref{figHorizonQuad} were conducted without discretizing function $f$ of interest 
but using its explicit form in order to have an exact expression of the $L^2$ approximation error.

\begin{figure}
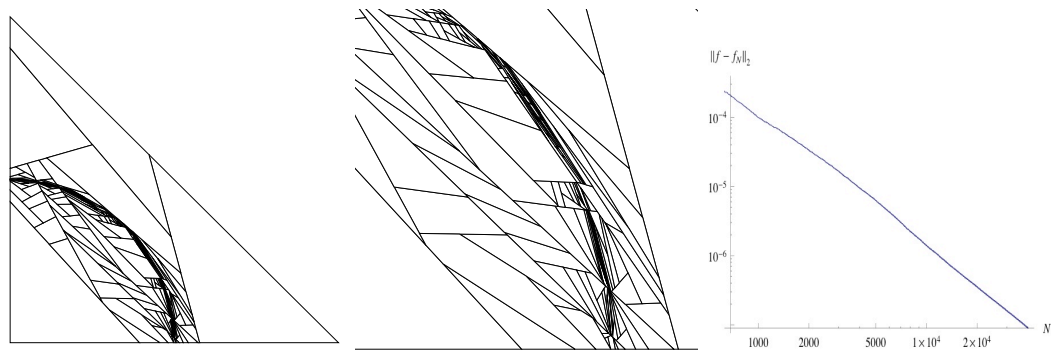

	\centering
		\includegraphics[width=4.5cm,height=4.5cm]{\pathPic/Subdivision/Numerical/Horizon200.pdf} 
		\includegraphics[width=4.5cm,height=4.5cm]{\pathPic/Subdivision/Numerical/Horizon200Detail.pdf}
		\includegraphics[width=4.5cm,height=4cm]{\pathPic/Subdivision/Numerical/GraphCvxHoriz.pdf} 		\caption{\label{figHorizonQuad}
		Triangulation produced by the greedy algorithm applied to the function $f=\chi_D$, with the decision function \iref{decisionConv} and $200$ elements (i). Detail (ii). Decay of the error $e_{\cT_N}(\chi_D)$, in logarithmic scale, for $1\leq N\leq 40000$.}
\end{figure}

\section{A greedy algorithm based on rectangles} 
\label{secGreedyRect}

We study in this section another variation of the greedy algorithm which is simpler than the previous ones.
This algorithm produces partitions $(\cT_N)_{N \geq 1}$ of the unit square $\Omega = [0,1]^2$ into rectangles aligned with the horizontal or vertical axes, as studied in Chapter 1 and illustrated on Figure \ref{figRectPart} in the main introduction of the thesis, instead of triangulations built of triangles of arbitrary direction. Here we work with piecewise constant approximants, and we measure the approximation 
error in the $L^\infty$ norm. These important simplifications allow us to establish 
an  error estimate which applies to any function $f\in C^1(\Omega)$ and is asymptotically optimal 
in accordance with the results in Chapter 1. We also
establish a second error estimate which gives some valuable 
(yet non optimal) information in a non-asymptotic sense, i.e. for any value of $N\geq 1$. 

In order to state our main result, we introduce some notation.
A rectangle $T\subset \R^2$ is a set of the form
$$
T = I \times J 
$$
where $I$ and $J$ are compact intervals of $\R$. 
For any rectangle $T$ and any $f$ in $C^0(T)$, we denote by $e_T(f)$ the difference between the supremum and the infimum of $f$ on $T$.
$$
e_T(f) := \sup_T f - \inf_T f.
$$
Note that $e_T(f)$ is also the double of the error of best approximation of $f$ on $T$ in the $L^\infty$ norm,
$$
e_T(f) = 2\inf_{c\in \R} \|f-c\|_{L^\infty(T)} .
$$
If a function $f\in C^1(T)$ attains its maximum on the rectangle $T$ at the point $z_M = (x_M,y_M)$, and its minimum at the point $z_m =(x_m,y_m)$ then observe that 
\begin{eqnarray*}
f(z_M) - f(z_m) &\leq& \left|\int_{x_m}^{x_M}  \partial_x f(x,y_m) dx\right| + \left|\int_{y_m}^{y_M}  \partial_y f(x_M,y) dy \right| \\
&\leq& |x_M-x_m| \|\partial_x f\|_{L^\infty(T)} + |y_M-y_m| \|\partial_y f\|_{L^\infty(T)}. 
\end{eqnarray*}
Therefore if $T = I \times J$, we obtain 
\be
\label{eqErrorBoundAff}
e_T(f) \leq |I| \|\partial_x f\|_{L^\infty(T)} + |J| \|\partial_y f\|_{L^\infty(T)}.
\ee
For any partition $\cT$ of a rectangle $T$ into smaller rectangles, and for any $f\in C^0(T)$, we denote by $e_\cT(f)$ the maximum of $e_T(f)$ among all $T \in \cT$:
$$
e_\cT(f) := \max_{T\in \cT} e_T(f).
$$
We clearly also have 
$$
e_\cT(f) = 2\min_{g\in V_\cT} \|f-g\|_{L^\infty(T)},
$$
where $V_\cT$ is the collection of functions on $T$ which are constant on the interior of each $T\in \cT$.

\begin{remark}
Please note that the definition of $e_T(f)$ in this section differs a lot from the definition of $e_T(f)$ in the previous section. Indeed, $e_T(f)$ stands here for the double of the best approximation of $f$ in $L^\infty(T)$ norm by a constant $c\in \R$, whereas $e_T(f)$ stands in the previous section for the  best approximation of $f$ in $L^2(T)$ norm by an affine function $\pi \in \P_1$.
\end{remark}


For any rectangle $T$ we denote by $(T_d,T_u)$ the down and up rectangles
which are obtained by a horizontal split of $T$ through its barycenter
and by $(T_l,T_r)$ the left and right rectangles
which are obtained by a vertical split of $T$ through its barycenter.
For any function $f\in C^0(T)$, we define 
$$
e_{T,h}(f) := \max\{e_{T_d}(f), e_{T_u}(f)\} \stext{ and } e_{T,v}(f) := \max\{e_{T_l}(f), e_{T_r}(f)\}.
$$

The greedy algorithm studied in this section takes as input a function $f\in C^0(\Omega)$, where $\Omega = [0,1]^2$ is the unit square, and a parameter $\rho\in [0,1]$. 
Given a partition $\cT$ of $\Omega$ into rectangles aligned with the coordinate axes, the greedy algorithm produces in one step a refinement $\cT'$ of $\cT$, satisfying $\#(\cT') = \#(\cT)+1$, proceeding as follows.
\begin{enumerate}
\item Select a rectangle $T\in \cT$ which maximizes $e_T(f)$.
\item (Greedy split) If 
\be
\label{eqCompTChild}
\min\{e_{T,h}(f), e_{T,v}(f)\}\leq \rho \,e_T(f),
\ee
then bisect the rectangle $T$ horizontally if $e_{T,h}(f)\leq e_{T,v}(f)$ and vertically otherwise, thus producing two children $T_1,T_2$ of $T$. Define 
$$
\cT' := \cT-\{T\}+ \{T_1,T_2\}.
$$
\item (Security split) If 
\be
\label{eqCompTChild2}
\min\{e_{T,h}(f), e_{T,v}(f)\} > \rho \,e_T(f),
\ee
then bisect the rectangle $T$ horizontally if $|I|\geq |J|$ and vertically otherwise where $I,J$ are the two compact intervals such that $T = I\times J$. This bisection produces two children $T_1,T_2$ of $T$ and we define 
$$
\cT' := \cT-\{T\}+ \{T_1,T_2\}.
$$
\end{enumerate}

Starting from the partition $\cT_1 = \{\Omega\}$ of the square $\Omega = [0,1]^2$, the greedy algorithm produces step after step a sequence $(\cT_N)_{N\geq 1}$ of partitions of $\Omega$ into dyadic rectangles, satisfying $\#(\cT_N) = N$.
The parameter $\rho\in [0,1]$ is chosen by the user. It should not be chosen too
small in order to avoid that all splits are of safety type which would then
lead to isotropic partitions. 
The choice $\rho=1$ is not indicated either. Indeed let $f(x,y) := \sin(4 \pi x)$, then for any rectangle $T$ of the form $T = [0,1] \times J$, where $J$ is a compact interval, one has 
$$
e_T(f) = e_{T,h}(f) = e_{T,v}(f) = 2.
$$
If a rectangle $T= [0,1] \times J$ is selected for bisection and if $\rho=1$, then $T$ is bisected horizontally and we are facing a similar situation for its children $T_u$ and $T_d$.
As a result, the greedy algorithm produces a sequence $(\cT_N)_{N\geq 1}$ of partitions of $\Omega$ consisting of rectangles of the form $[0,1]\times J$ where $J$ are dyadic intervals. This shows that the approximations produced by the algorithm fail to converge towards $f$, and the global error remains $e_{\cT_N}(f) = 1$ for all $N \geq 1$.
On the contrary, it is established in \cite{CMi2} that for any choice of $\rho\in [0,1[$ and any $f\in C^0(\Omega)$, the partitions $(\cT_N)_{N\geq 1}$ produced by the greedy algorithm satisfy $e_{\cT_N}(f) \to 0$ as $N \to \infty$.

We now prove that,
using the specific value $\rho=\frac 1 {\sqrt 2}$ and measuring the error in the $L^\infty$ norm,
the refinement algorithm has optimal convergence properties \iref{eqAsymptRect}
similar to Theorem \ref{thUpper}.
The authors suspect that this result can be extended to the $L^p$ norm, where $1\leq p \leq \infty$, and to an arbitrary parameter $\rho\in ]1/2,1[$,
at the price of more technicalities.

\begin{theorem}
\label{thRect}
There exists a constant $C>0$ such that for any $f\in C^1(\Omega)$, the partitions
$(\cT_N)_{N\geq 1}$ produced by the modified greedy refinement algorithm with parameter $\rho:=\frac 1 {\sqrt 2}$
satisfy the asymptotic convergence estimate
\be
\label{eqAsymptRect}
\limsup_{N\to +\infty} \, N^{\frac 1 2} e_{\cT_N}(f) \leq C \left\|\sqrt{|\partial_x f\ \partial_y f|}\right\|_{L^2}.
\ee
Furthermore for all $N \geq 1$ one has 
\be
\label{eqRobustRect}
e_{\cT_N}(f)\leq C \|\nabla f\|_{L^\infty(\Omega)}N^{-\frac 1 2}.
\ee
\end{theorem}

The asymptotical error estimate \iref{eqAsymptRect} is the best than one should hope to obtain. 
Indeed it follows from Theorem \ref{thLower} that there exists a constant $c>0$ independent of $f$ and such that 
\be
\label{eqLowerRect}
\liminf_{N\to +\infty} \, N^{\frac 1 2} e_{\cT'_N}(f) \geq c \left\|\sqrt{|\partial_x f\ \partial_y f|}\right\|_{L^2},
\ee
for any admissible sequence $(\cT'_N)_{N\geq N_0}$ of partitions of $\Omega$ into rectangles. (Admissibility is a minor technical condition, see Definition \ref{defAdmiRect}) 

Here and after, we use the $\ell^\infty$ norm on $\RR^2$ for
measuring the gradient: for $z=(x,y)\in\Omega$
$$
|\nabla f(z)|:=\max \{|\partial_x f(z)|,|\partial_y f(z)|\},
$$
and
$$
\|\nabla f\|_{L^\infty(T)}:= \sup_{z\in T}|\nabla f(z)|=\max\left\{\|\partial_x f\|_{L^\infty(T)},\|\partial_y f\|_{L^\infty(T)}\right\}.
$$
We consider a fixed function $f\in C^1(\Omega)$ and we denote by $(\cT_N)_{N\geq 1}$ the sequence of partitions of  the square $\Omega = [0,1]^2$ into dyadic rectangles generated by the greedy algorithm. Furthermore we define for each $N \geq 1$
$$
\ve(N) := e_{\cT_N}(f).
$$
Finally we sometimes use the notation $x(z)$ and $y(z)$ to denote the coordinates 
of a point $z\in\RR^2$.
The proof of Theorem \ref{thRect} requires a preliminary result.

\begin{lemma}
\label{CMlemmaRobust}
Let $T_0=I_0\times J_0 \in\cT_M$ be a dyadic rectangle obtained 
at some stage $M\geq 1$ of the refinement algorithm, and let $T=I\times J\in \cT_N$
be a dyadic rectangle obtained at some later stage $N>M$ and such
that $T\subset T_0$. We then have 
\be
\label{eqLemmaRobust}
|I|\geq \min\left \{|I_0|, \frac {\ve(N)} {4\|\nabla f \|_{L^\infty(T_0)}} \right \} \stext{ and } 
|J|\geq \min\left \{|J_0|, \frac {\ve(N)} {4\|\nabla f \|_{L^\infty(T_0)}}\right \}.
\ee
\end{lemma}
\proof
Since the coordinates $x$ and $y$ play symmetrical roles, it suffices to prove the first inequality.
We reason by contradiction. If the inequality does not hold, there exists 
a rectangle $T'=I'\times J'$ in the chain that led from $T_0$ to $T_1$
which is such that
\be
\label{eqIRobust}
|I'| < \frac  {\ve(N)} {2\|\nabla f \|_{L^\infty(T_0)}},
\ee
and such that $T'$ is split vertically by the algorithm. If this was a safety split, 
we would have that $|J'|\leq |I'|$
and therefore using \iref{eqErrorBoundAff}
$$
e_{T'}(f) \leq (|I'|+|J'|) \|\nabla f \|_{L^\infty(T)} \leq 2|I'|\|\nabla f \|_{L^\infty(T)} < \ve(N),
$$ 
which is a contradiction, since all ancestors of $T$ should satisfy $e_{T'}(f)\geq \ve(N)$. 
Hence this split was necessarily a greedy split.

Let $z_m:={\rm Argmin}_{z\in T'} f(z)$ and $z_M:={\rm Argmax}_{z\in T'} f(z)$, and let $T''$ be the child of $T'$ (after the vertical split) containing $z_M$. Then $T''$ also contains a point $z'_m$ such that $|x(z'_m) - x(z_m)|\leq |I'|/2$ and $y(z'_m) = y(z_m)$. It follows that
\begin{eqnarray*}
e_{T',v}(f) &\geq &e_{T''}(f) \\
& \geq& f(z_M) - f(z'_m)\\
 & \geq& f(z_M) - f(z_m)- \|\partial_x f\|_{L^\infty(T')} \frac {|I'|} 2\\
 & \geq & e_{T'}(f) - \frac 1 4 \ve(N)\\ 
&  \geq& \frac 3 4 e_{T'}(f) \\
& > & \rho e_{T'}(f).
\end{eqnarray*}
where we used \iref{eqIRobust} to obtain the fourth line.
The error was therefore insufficiently reduced which contradicts a greedy split.
\sq

\noindent
{\bf Proof of Inequality \iref{eqRobustRect}}
We first observe that for any $N \geq 1$,
$$
\ve(N) \leq \ve(1) = \sup_\Omega f - \inf_\Omega f \leq 2 \|\nabla f\|_{L^\infty(\Omega)}.
$$
Applying Lemma \ref{CMlemmaRobust} to the square $\Omega = [0,1] \times [0,1]$ we obtain 
$$
\min\{|I|, |J|\} \geq \frac {\ve(N)} {4 \|\nabla f \|_{L^\infty(\Omega)}}
$$
for any rectangle $T = I \times J \in \cT_N$.
Therefore 
$$
1= |\Omega| = \sum_{T = I \times J \in \cT_N} |I| \times |J| \geq N \left( \frac {\ve(N)} {4 \|\nabla f \|_{L^\infty(\Omega)}} \right)^2, 
$$
hence 
\be
\label{eqEpsToZero}
N^{\frac 1 2 } \ve(N) \leq 4 \|\nabla f \|_{L^\infty(\Omega)},
\ee
which concludes the proof of \iref{eqRobustRect}. 
\sq\\

\noindent
{\bf Proof of Inequality \iref{eqAsymptRect}} 
We consider a small but fixed $\delta >0$, and we define $h(\delta)$ as the maximal $h>0$ such that
\be
\label{eqUnifCont}
\forall z,z'\in \Omega,\ |z-z'|\leq 2h(\delta) \Rightarrow |\nabla f(z)-\nabla f(z')|\leq  \delta.
\ee
For any rectangle $T=I\times J\subset \Omega$, we thus have 
\be
\label{CMeqhdelta}
\begin{array}{ccc}
e_{T}(f) &\geq& (\|\partial_x f\|_{L^\infty(T)}-\delta) \min\{h(\delta),|I|\}, \\
e_{T}(f) &\geq& (\|\partial_y f\|_{L^\infty(T)}-\delta) \min\{h(\delta),|J|\}.
\end{array}
\ee
Let $\delta>0$ and let $M=M(f,\delta)$ be the smallest value of $N$ such that $\e(N)< 9\delta h(\delta)$.
For all $N\geq M$, and therefore $\e(N)< 9\delta h(\delta)$, we consider the 
partition $\cT_N$ which is a refinement of $\cT_{M}$. For any rectangle $T_0=I_0\times J_0 \in \cT_M$,
we denote by $\cT_N(T_0)$ the set of rectangles of $\cT_N$ that are contained $T_0$.
We thus have
$$
\cT_N:=\bigcup_{T_0\in\cT_M}\cT_N(T_0),
$$
and $\cT_N(T_0)$ is a partition of $T_0$. We shall next bound by below
the side length of $T=I\times J$ contained in $\cT_N(T_0)$,
distinguishing different cases depending on the behaviour of $f$ on $T_0$.
\nl
\nl
{\it Case 1.} If $T_0\in\cT_M$ is such that $\|\nabla f\|_{L^\infty(T_0)} \leq 10 \delta$,
then a direct application of Lemma \ref{CMlemmaRobust} shows that for all $T=I\times J\in \cT_N(T_0)$ we have 
\be
\label{CMeqFlat}
|I|\geq \min\left \{|I_0|, \frac {\ve(N)} {40\delta}\right \} \stext{ and } |J|\geq \min\left \{|J_0|, \frac {\ve(N)} {40\delta}\right \}
\ee
\nl
{\it Case 2.} If $T_0\in\cT_M$ is such that
$\|\partial_x f\|_{L^\infty(T_0)} \geq 10 \delta$ and $\|\partial_y f\|_{L^\infty(T_0)} \geq 10 \delta$,
we then claim that for all $T=I\times J\in \cT_N(T_0)$ we have
\be
\label{CMeqAffine}
|I| \geq \min\left\{|I_0|, \frac{\ve(N)} {20\|\partial_x f\|_{L^\infty(T_0)}}\right\} \stext{ and } 
|J| \geq \min\left\{|J_0|, \frac{\ve(N)} {20\|\partial_y f\|_{L^\infty(T_0)}}\right\},
\ee
and that furthermore 
\be
|T_0| \ \|\partial_x f\|_{L^\infty(T_0)} \|\partial_y f\|_{L^\infty(T_0)} \leq \left (\frac {10}9\right )^2 \int_{T_0} |\partial_x f\ \partial_y f| dx dy.
\label{CMintegralstatement}
\ee
This last statement easily follows by the following observation: 
recalling \iref{CMeqhdelta} we obtain
$$
9\delta h(\delta)> \ve(M) \geq e_{T_0}(f) \geq (\|\partial_y f\|_{L^\infty}(T_0)-\delta)\min\{h(\delta),|J_0|\}
\geq 9\delta\min\{h(\delta),|J_0|\},
$$
therefore $|J_0| \leq h(\delta)$. Similarly $|I_0|\leq h(\delta)$ and 
therefore  for all $z\in T_0$ we obtain 
\begin{eqnarray*}
|\partial_x f(z)| &\geq& \|\partial_x f\|_{L^\infty(T_0)} - \delta \geq \frac 9 {10}\|\partial_x f\|_{L^\infty(T_0)},\\
|\partial_y f(z)| &\geq& \|\partial_y f\|_{L^\infty(T_0)} - \delta \geq \frac 9 {10}\|\partial_y f\|_{L^\infty(T_0)}.
\end{eqnarray*}
Integrating over $T_0$ yields \iref{CMintegralstatement}.
Furthermore the previous equations show that $\partial_x f$ and $\partial_y f$ have a constant sign on $T_0$. Let us assume that both are positive, then for any rectangle $T = [a,b] \times [c,d]\subset T_0$ we obtain 
\begin{eqnarray*}
e_T(f) & = & \sup_T f-\inf_T f \\
& \geq & f(b,d) - f(a,c)\\
& \geq & \int_a^b \partial_x f(x,c) dx+ \int_c^d \partial_y f (b,y) dy\\
& \geq & \frac 9 {10} \left((b-a) \|\partial_x f\|_{L^\infty(T_0)} + (d-c) \|\partial_y f\|_{L^\infty(T_0)}\right)
\end{eqnarray*}
A similar reasoning can be applied if $\partial_x f$ or $\partial_y f$ is negative is negative on $T_0$. Recalling \iref{eqErrorBoundAff} we therefore obtain for any $T \subset T_0$
\be
\label{CMestimAffine}
\frac 9 {10} \leq \frac {e_{T}(f) }{\|\partial_x f\|_{L^\infty(T_0)}|I|+\|\partial_y f\|_{L^\infty(T_0)} |J|}\leq 1.
\ee
Clearly the two inequalities in \iref{CMeqAffine} are symmetrical, and it suffices to prove the first one. 
Similar to the proof of Lemma \ref{CMlemmaRobust}, we 
reason by contradiction, assuming that a rectangle $T'=I'\times J'$ with 
\be
\label{eqLenIRect}
|I'| \|\partial_x f\|_{L^\infty(T_0)}< \frac{\ve(N)}{10}
\ee
was split vertically by the algorithm in the chain leading from $T_0$ to $T$. 
A simple computation using \iref{CMestimAffine} shows that 
$$
\frac{e_{T',h}(f) }{e_{T'}(f)}\leq 
\frac{e_{T',h}(f)}{e_{T',v}(f)} \leq \frac 5 9 \times \frac{1+2\sigma}{1+\sigma/2}  \stext{ with } \sigma := \frac{\|\partial_x f\|_{L^\infty(T_0)} |I'|}{\|\partial_y f\|_{L^\infty(T_0)} |J'|}.
$$
In particular if $\sigma<0.2$ the we obtain that 
$$
e_{T',h}(f) \leq \rho e_{T'} (f) \stext{ and } e_{T',h}(f) \leq e_{T',v} (f),
$$
hence the algorithm performs a horizontal greedy split on $T'$, which contradicts our assumption.
Therefore $\sigma\geq 0.2$, but this also leads to a contradiction since, using \iref{eqLenIRect},
$$
\ve(N) \leq  e_{T'}(f) \leq \|\partial_x f\|_{L^\infty(T_0)} |I'|+ \|\partial_y f\|_{L^\infty(T_0)} |J'| 
\leq (1+\sigma^{-1}) \frac {\e(N)} {10} < \ve(N).
$$
\nl
\nl
{\it Case 3.} If $T_0\in\cT_M$ is such that $\|\partial_x f\|_{L^\infty(T_0)} \leq 10 \delta$ and $\|\partial_y f\|_{L^\infty(T_0)} 
\geq 10 \delta$, we then claim that for all $T=I\times J\in \cT_N(T_0)$ we have, \text{ with } $C_0:=200$,
\be
\label{CMeq1d}
|I| \geq \min\left \{|I_0|, \frac{\ve(N)} {C_0\delta}\right \} \stext{ and } |J| \geq \min\left \{|J_0|, \frac {\ve(N)} {4 \|\nabla f\|_{L^\infty}}\right \}
\ee
with a symmetrical result if $T_0$ is such that $\|\partial_x f\|_{L^\infty(T_0)} \geq 10 \delta$ and $\|\partial_y f\|_{L^\infty(T_0)} 
\leq 10 \delta$. The right part of \iref{CMeq1d} is a direct consequence of Lemma
\ref{CMlemmaRobust}, hence we focus on the left part.
Applying the second inequality of \iref{CMeqhdelta} to $T=T_0$, we obtain
$$
9\delta h(\delta)> e_{T_0}(f) \geq (\|\partial_y f\|_{L^\infty}(T_0)-\delta)\min\{h(\delta),|J_0|\}
\geq 9\delta\min\{h(\delta),|J_0|\},
$$
from which we infer that $|J_0| \leq h(\delta)$. If $z_1,z_2\in T_0$ and $x(z_1) = x(z_2)$ we therefore have 
$|\partial_y f(z_1)| \geq |\partial_y f(z_2)|-\delta$. It follows that for any rectangle $T=I \times J\subset T_0$ we have
\be
\label{CMerror1D}
(\|\partial_y f\|_{L^\infty(T)}-\delta) |J| \leq e_{T}(f) \leq \|\partial_y f\|_{L^\infty(T)} |J|+ 10\delta |I|.
\ee
We then again reason by contradiction, assuming that a rectangle $T'=I'\times J'$ with 
\be
\label{eqLengthR}
|I'| < \frac{2\ve(N)}{C_0\delta}
\ee
was split vertically by the algorithm in the chain leading from $T_0$ to $T$. 
If $\|\partial_y f\|_{L^\infty(T')} \leq 10\delta $, then $\|\nabla f\|_{L^\infty(T')}\leq 10\delta$ and 
Lemma \ref{CMlemmaRobust} shows that $T'$ should not have been split vertically, which is a contradiction.
Otherwise $\|\partial_y f\|_{L^\infty(T')} -\delta \geq \frac 9{10} \|\partial_y f\|_{L^\infty(T')}$, and we obtain using \iref{CMerror1D} and \iref{eqLengthR} that 
$$
\frac 9 {10} \|\partial_y f\|_{L^\infty(T')} |J'|\leq e_{T'}(f) \leq \|\partial_y f\|_{L^\infty(T')} |J'|+ 10 \delta \, \frac{2\ve(N)}{C_0\delta}.
$$
Recalling that $e_{T'}(f) \geq \ve(N)$ we obtain
\be
\label{CMestimDy}
\left(1-\frac{20} {C_0}\right)e_{T'}(f) \leq \|\partial_y f\|_{L^\infty(T')} |J'|\leq \frac {10} 9 e_{T'}(f).
\ee
We now consider the children $T'_v$ and $T'_h$ of $T'$
of maximal error after a horizontal and vertical split
respectively. We obtain 
\begin{eqnarray*}
e_{T',h}(f) &=&e_{T'_h}(f)\\
&  \leq  &\|\partial_y f\|_{L^\infty(T')} \frac{|J'|} 2+ 10\delta |I'| \\
&  \leq &\frac 5 9 e_{T'}(f)+ \frac{20\, \ve(N)}{C_0} \\
& \leq  &\left(\frac 5 9+\frac{20} {C_0}\right)  e_{T'}(f)= \frac{59}{90} e_{T'}(f),
\end{eqnarray*}
where we used \iref{CMerror1D} in the second line, and \iref{eqLengthR} and \iref{CMestimDy} in the third line.
On the other hand, using \iref{CMestimDy} in the fourth line
\begin{eqnarray*}
e_{T',v}(f) & =&e_{T'_h}(f)\\
 &\geq &(\|\partial_y f\|_{L^\infty(T')}-\delta) |J| \\
&  \geq &\frac 9 {10} \|\partial_y f\|_{L^\infty(T')} |J'|  \\
& \geq &\frac 9 {10} \left(1-\frac{20} {C_0}\right) e_{T'}(f)= \frac {81}{100}e_{T'}(f).
\end{eqnarray*}
Therefore $e_{T',v}(f) > e_{T',h}(f)$
which is a contradiction, since our decision rule would then select a horizontal split.\\

We now recall that $\ve(N)\to 0$ as $N \to \infty$, see \iref{eqEpsToZero}, and we choose $N$ large enough so that 
the minimum in \iref{CMeqFlat}, \iref{CMeqAffine} and \iref{CMeq1d} is
always equal to the second term. 
For all  $T\in \cT_N(T_0)$, we respectively find that 
$$
\frac{\ve(N)^2} {|T|}\leq  C^2
\left\{\begin{array}{ccl}
\delta^2 & \text{ if } & \|\nabla f\|_{L^\infty(T_0)}\leq 10\delta\\
\frac 1 {|T_0|} \int_{T_0} |\partial_x f\ \partial_y f| & \text{ if } &   \|\partial_x f\|_{L^\infty(T_0)} \geq 10 \delta \text{ and } \|\partial_y f\|_{L^\infty(T_0)} \geq 10 \delta\\
 \delta \|\nabla f\|_{L^\infty} & \text{ if } & \|\partial_x f\|_{L^\infty(T_0)} \leq 10 \delta \text{ and } \|\partial_y f\|_{L^\infty(T_0)} \geq 10 \delta \text{ (or reversed).}
\end{array}\right.
$$
with $C^2:= \max\{40^2,\ 20^2 (10/9)^2,\ 800\}=40^2$. For $z\in \Omega$, we set
$\psi_N(z) := \frac 1 {|T|}$ where $T\in \cT_N$ is such that $z\in T$, and we obtain
$$
N=\#(\cT_N)  = \int_{\Omega} \psi_N(z) dz \leq C^2\ve(N)^{-2} \left(  \int_{\Omega} |\partial_x f\ \partial_y f| dx dy+ \delta \|\nabla f\|_{L^\infty}+\delta^2\right).
$$
Therefore 
$$
\limsup_{N \to \infty} N^{\frac 1 2} \ve(N) \leq C \left(\int_{\Omega} |\partial_x f\ \partial_y f| dx dy+ \delta \|\nabla f\|_{L^\infty}+\delta^2\right)^\frac 1 2.
$$
Taking the limit as $\delta \to 0$, we obtain the announced result
$$
\limsup_{N\to \infty} \, N^{\frac 1 2} \ve(N)\leq C\left \| \sqrt {|\partial_x f\partial_y f|}\right \|_{L^2},
$$
which concludes the proof. \hfill $\Box$


\section{Appendix}
\subsection{Proof of Proposition \ref{rho48}}

Let $q$ be a quadratic function of mixed type, and let $T$ be a triangle with edges $a,b,c$ such that $|\bq(a)|\geq |\bq(b)|\geq |\bq(c)|$.
Up to a linear change of variables, we may assume that the quadratic part of $q$ is $\bq(x,y) = x^2-y^2$.
Up to a translation, linear rescaling and permutation between the $x$ and $y$
coordinates, we may assume that the edge $a$ is centered at $0$
and such that $\bq(a) = 4$. We write $\frac a 2 = (u,v)$, and assume that $u>0$ without loss of generality. Then $1 =  \bq(\frac a 2) = u^2-v^2$, and there must be $\theta\in \R$ such that $(u,v) = (\cosh(\theta),\sinh(\theta))$.

We next observe that the linear transformation $\phi$ of matrix 
$$
\left(\begin{array}{cc} 
\cosh(\theta) & -\sinh(\theta)\\
-\sinh(\theta) & \cosh(\theta)
\end{array}\right)
$$
leaves $\bq$ invariant, in the sense that $\bq\circ\phi = \bq$, and satisfies $\phi(\frac a 2) = (1,0)$. 
Up to such a transformation, we may therefore assume that $a=(2,0)$. Therefore
the two vertices of $T$ at the ends of the edge $a$ are $(-1,0)$ and $(1,0)$. 
We denote the third vertex by $(s,t)$.
There is no loss of generality, eventually, in assuming that $s\geq 0$ and $t> 0$. Note that $|T| = t$ and 
$\rho_\bq(T) = \frac 4 t$.

We now specialize to the case where $\rho_\bq(T) \geq 4$, which is equivalent to $t\leq 1$. Recall that the edges of $T$ are such that $|\bq(a)|\geq|\bq(b)|\geq|\bq(c)|$. 
The following lines show that $\bq(1+s,t)\geq |\bq(s-1,t)|$, which implies that 
$b=-(1+s,t)$ and $c=(s-1,t)$:
\begin{eqnarray*}
\bq(1+s,t)-\bq(s-1,t) &=& 4 s \geq 0,\\
\bq(1+s,t)+\bq(s-1,t) &=& 2(1+s^2 - t^2) \geq 0.
\end{eqnarray*}
In addition we see that $\bq(b)\geq 0$ and we thus have proved
that $\bq(a)\bq(b)\geq 0$. For the second and third inequality in \iref{rho4}
we remark that 
$\bq(c)=(s-1)^2-t^2$ and $\bq(d)=s^2-t^2$. Clearly
$$
-1\leq -t^2 \leq \min\{\bq(c),\bq(d)\}.
$$
If $0\leq s\leq 1$, we also clearly have 
$$
\max\{\bq(c),\bq(d)\} \leq 1.
$$
If $s\geq 1$ we have
$$
\bq(c)\leq \bq(d) = s^2-t^2 =(s+1)^2-t^2 - (2s+1) = \bq(b) -(2s+1) \leq \bq(a) -3 = 1.
$$
We thus always have 
\be
\max\{|\bq(c)|,|\bq(d)|\}\leq \frac{|\bq(a)|} 4 = 1,
\ee 
which concludes the proof of the inequalities \iref{rho4}.

Last, we specialize to the case where $\rho_\bq(T) = \frac 4 t \geq 8$, equivalently $t\leq \frac 1 2$, to prove \iref{rho8}:
$$
\bq(b)+\bq(c)= (s+1)^2-t^2 + (s-1)^2-t^2 = 2+2s^2-2t^2\geq \frac 3 2 = \frac 3 8 \bq(a).
$$

\subsection{Proof of Proposition \ref{propDiamZeroLInf}}
\label{secPropDiamZero}
For any triangle $T$ we denote by $T_x$ the interval defined as
the projection of $T$ on the $x$ axis, and by $|T_x|$ its length. 
We denote by $I_T$ and $I_{T_x}$ the
two dimensional and one dimensional local interpolations operators 
on $T$ and $T_x$ respectively.
It is clear that if $g$ is a $C^2$ convex function
of one variable and if $G(x,y):= g(x)$, then 
\be
\|G-I_T G\|_{L^\infty(T)} = \|g-I_{T_x}g\|_{L^\infty(T_x)}.
\label{inter1d2d}
\ee
The following lemma compares the 
interpolation error on an interval of $\R$ and on a sub-interval.

\begin{lemma} Let $g\in C^2(\R,\R)$ be such that $0<m\leq g'' \leq M$.
Let $x_1,x_2,x_3$ be real numbers satisfying $x_1<\frac{x_1+x_3} 2\leq x_2<x_3$, 
and denote $u := x_2-x_1\leq v := x_3-x_1$.
Then denoting by $I_u$ and $I_v$ the interpolation operators 
on the intervals $[x_1,x_2]$ and $[x_1,x_3]$ respectively, we have
$$
\|g-I_v g\|_{L^\infty([x_1,x_3])} \geq m v^2/8
$$
and
$$
\|g-I_u g\|_{L^\infty([x_1,x_2])}\leq \left(1-\alpha\frac{v-u}{v}\right) \|g-I_v g\|_{L^\infty([x_1,x_3])}
$$
with $\alpha :=\frac 1 4 \sqrt{m/M}$.
\end{lemma}

\proof
Let us define $h_v = I_v g-g$. Since $h_v'' = -g''$ and $h_v(x_1) = h_v(x_3) = 0$, 
$h_v$ can be represented as the integral
\be
h_v(x) = \int_{x_1}^{x_3} K_v(x,y) g''(y) dy.
\label{green}
\ee
where the Green kernel $K_v$ is given by
$$
K_v(x,y) := \frac 1 v
\left\{
\begin{array}{ll}
(x-x_1)(x_3-y) &\text{ if } x\leq y,\\
(x_3-x)(y-x_1) &\text{ if } x\geq y.
\end{array}
\right.
$$
Of course, we have a similar representation of $h_u=I_ug-g$ with a kernel $K_u$.\\
The first part of the proposition immediately follows from  
$$
h_v\left(\frac{x_1+x_3} 2\right)\geq m\int_{x_1}^{x_3} K_v\left(\frac{x_1+x_3} 2,y\right) dy \geq m v^2/8.
$$
In order to prove the second part, we shall compare
the Green Kernels $K_u$ and $K_v$. For this purpose, we define 
$$
\mu(x) := \frac{(x_2-x)/u}{(x_3-x)/ v}.
$$
For all $x,y\in [x_1,x_2]$, we thus have
$$
\frac{K_u(x,y)}{K_v(x,y)}=\min\{\mu(x),\mu(y)\}\leq \mu(x).
$$
Therefore, defining $x_u := \argmax_{t\in [x_1,x_2]}(I_u g-g)(t)$
and using \iref{green}, we obtain
$$
\begin{array}{ll}
\|g-I_u g\|_{L^\infty([x_1,x_2])}&=h_u(x_u) \\
& = 
\int_{x_1}^{x_2} K_u(x_u,y) g''(y)dy\\
& \leq \mu(x_u) \int_{x_1}^{x_2} K_v(x_u,y) g''(y)dy\\
& \leq \mu(x_u) \int_{x_1}^{x_3} K_v(x_u,y) g''(y)dy\\
&\leq \mu(x_u) \|g-I_v g\|_{L^\infty([x_1,x_3])}.
\end{array}
$$
In order to conclude, we need to estimate $\mu(x_u)$.
One easily checks by differentiation that $\mu$ is decreasing and concave 
on the interval $[x_1,x_2]$, and therefore for all $x\in [x_1,x_2]$
$$
\mu(x)\leq 1-(x-x_1)\frac{v-u}{v^2}.
$$
Differentiating the integral representation of $h_u$, we obtain
$$
u h_u'(x_u) = -\int_{x_1}^{x_u} (y-x_1)g''(y)dy+\int_{x_u}^{x_2} (x_2-y)g''(y)dy.
$$
Since $h'(x_u) = 0$ and $0<m\leq g''\leq M$, we obtain
$$
(x_2-x_u)^2m \leq (x_u-x_1)^2 M.
$$
Since $x_2\geq \frac{x_1+x_3} 2$, this gives $x_u-x_1\geq \sqrt{m/M} \frac v 4$.
Therefore,
$$
\mu(x_u)\leq 1-\frac 1 4 \sqrt{m/M} \frac{v-u} v,
$$
which concludes the proof.
\sq
\nl
The following corollary uses the above lemma to compare the values 
of the $L^\infty$-based decision functions for a convex function that depends only of one variable. 
For any vector $v\in \R^2$ we denote by $v_x$ the absolute value of its $x$ coordinate.

\begin{corollary} \label{CorolErCvx}
Let $T$ be a triangle with edges $a,b,c$, such that $a_x\geq b_x\geq c_x$.
Let $G(x,y) = g(x)$, where $g\in C^2$ and $0<m\leq g'' \leq M$.
Then, with $d_T$ defined by \iref{optilinfedge},
\begin{eqnarray*}
	d_T(b,G)-d_T(a,G) &\geq& C a_x (a_x-b_x),\\
	d_T(c,G)-d_T(a,G) &\geq& C a_x^2/2,
\end{eqnarray*}
with $C= \frac{m^{3/2}}{32\sqrt M}$.
\end{corollary}

\proof
We recall the notation $\alpha_T(f):=\|f-I_T f\|_{L^\infty}$. It follows from \iref{inter1d2d} that
$$
\alpha_T(G)=\|g-I_{T_x}g\|_{L^\infty(T_x)}.
$$ 
We label the extremities of $a$ by $i\in \{1,2\}$, and denote 
by $T^i_e$, $i\in \{1,2\}$,  $e\in\{a,b,c\}$ the child of $T$ 
resulting from the bisection through the edge $e$ in such a way 
that $T^i_e$ contains the extremity of $a$ of label $i$. We denote by $T^i_{e,x}$ the projection of $T^i_e$ onto the $x$-axis.
Then (up to exchanging the labels of the extremities of $a$),
$$
\begin{array}{llllll}
|T^1_{a, x}| &=& b_x,     		& |T^2_{a,x}| &=& a_x/2,\\
|T^1_{b,x}| &=& a_x,     		& |T^2_{b,x}| &=& c_x+b_x/2,\\
|T^1_{c,x}| &=& b_x+c_x/2,	& |T^2_{c,x}| &=& a_x.\\
\end{array}
$$
In particular, we have $T^2_{a,x}\subset T^2_{b,x}$ and $T^1_{a,x}\subset T^1_{b,x}$
and therefore
$$
\alpha_{T^2_a}(G)=\|g-I_{T^2_{a,x}}g\|_{L^\infty(T^2_{a,x})}
\leq \|g-I_{T^2_{b,x}}g\|_{L^\infty(T^2_{b,x})} = \alpha_{T^2_b}(G),
$$
and similarly $\alpha_{T^1_a}(G)\leq  \alpha_{T^1_c}(G)$. Moreover, 
we can apply the previous lemma with
$[x_1,x_2]=T^1_{a,x} \subset T^1_{b,x}=[x_1,x_3]$
or with $[x_1,x_2]=T^2_{a,x} \subset T^2_{c,x}=[x_1,x_3]$
which respectively leads to 
$$
\alpha_{T^1_b}(G)-\alpha_{T^1_a}(G) \geq \frac{m^{3/2}}{32\sqrt M} a_x (a_x-b_x),
$$
and 
$$
\alpha_{T^2_c}(G)-\alpha_{T^2_a}(G) \geq \frac{m^{3/2}}{32\sqrt M} a_x^2/2. 
$$
This allows us to conclude since 
$d_T(e,G) = \alpha_{T^1_e}+\alpha_{T^2_e}$, for $s\in \{a,b,c\}$.
\sq

Using the above result, we now prove that the decision function $d_T$
tends to prescribe a longest edge bisection with respect to the euclidean
metric when the triangle $T$ becomes too thin.

\begin{corollary} \label{CorolZero}
Let $f$ be a convex function such that $m\Id \leq d^2 f\leq M\Id$, 
and let $T$ be a triangle with measure of non-degeneracy 
$\sigma(T)$ for the euclidean metric 
and edges $a,b,c$, such that $|a|\geq |b|\geq |c|$. 
Then if $\sigma(T)> 2 K$ the bisection prescribed by the 
decision function $d_T(\cdot,f)$ is a $\delta$-near 
longest edge bisection with respect to the euclidean metric
(in the sense of definition \ref{def:dg}), with $K = 128 (\frac M m)^{3/2}$ and $\delta:=\frac K {\sigma(T)}$.
\end{corollary}

\proof 
We denote by $p(X)$ the
affine orthogonal projection of a point
$X\in \RR^2$ onto the line which includes
the edge $a$, and we denote by $O$ the midpoint of $a$. We then define
$$
\tilde f(X) = f(p(X)) + df_O(X-p(X)).
$$
Then, using the notation $X(z)=p(X)+z(X-p(X))$, we have
\begin{eqnarray*}
f(X) -\tilde f(X) &=& f(X) -f(p(X)) - df_O(X-p(X))\\
&=& \int_{0}^1 (df_{X(z)}-df_O)(X-p(X)) dz.
\end{eqnarray*}
But for $X \in T$, we have $X(z)\in T$ for all $z\in [0,1]$,
so that $|X(z)-O|\leq |a|$ and therefore
$$
\|df_{X(z)}-df_O\|\leq M |a|.
$$
We also have 
$$
|X-p(X)|\leq \frac{2 |T|}{|a|}\leq \frac{|a|}{\sigma(T)}.
$$
Therefore, $\|f-\tilde f\|_{L^\infty(T)}\leq M \frac{|a|^2}{\sigma(T)}$. 
We can apply the previous lemma to the function $\tilde f$, since it is the 
sum of a function of one variable and of an 
affine function which has no effect on 
the interpolation error. Assuming without loss of generality (up to a rotation) that
$a$ is parallel to the $x$ axis, this gives us
\begin{eqnarray*}
	d_T(b,\ti f)-d_T(a,\ti f) &\geq& C |a| (|a|-b_x) \geq C |a| (|a|-|b|)  ,\\
	d_T(c,\ti f)-d_T(a,\ti f) &\geq& C |a|^2/2.
\end{eqnarray*}
with $C= \frac{m^{3/2}}{32\sqrt M}$. Since the 
Lebesgue constant of the interpolation 
operator on any triangle is $2$, we have
$$
|d_T(e,f) - d_T(e,\tilde f)|\leq 4||f-\tilde f||_{L^\infty(T)}
$$
which implies the following inequalities
\begin{eqnarray*}
	d_T(b,f)-d_T(a,f) &\geq& C |a| (|a|-|b|) - 4M|a|^2/\sigma(T),\\
	d_T(c,f)-d_T(a,f) &\geq& |a|^2(C/2-4M/\sigma(T)).
\end{eqnarray*}
The second inequality shows that the edge $c$ cannot be cut if $C/2-4M/\sigma(T)>0$ which
is equivalent to $\sigma(T)>2K$. The first inequality shows that $b$ may be cut
provided that $C(|a|-|b|) - 4M|a|/\sigma(T)\leq 0$, i.e. $|b| \geq (1-\frac {4M}{C\sigma(T)})|a|$
which shows the property of $\delta$-near longest edge bisection 
with $\delta:=\frac K {\sigma(T)}$.
\sq

The proof of Proposition \ref{propDiamZeroLInf} directly follows from this last result, proceeding exactly as in the last paragraph of the proof of Proposition \ref{PropDiamZero}.

%


\end{document}